# SUPERBIMATRICES AND THEIR GENERALIZATIONS

**W. B. Vasantha Kandasamy**
**Florentin Smarandache**

**2009**

# SUPERBIMATRICES AND THEIR GENERALIZATIONS

**W. B. Vasantha Kandasamy**
e-mail: **vasanthakandasamy@gmail.com**
web: **http://mat.iitm.ac.in/~wbv**
**www.vasantha.in**

**Florentin Smarandache**
e-mail: **smarand@unm.edu**

**2009**

# CONTENTS

# PREFACE

The systematic study of supermatrices and super linear algebra has been carried out in 2008. These new algebraic structures find their applications in fuzzy models, Leontief economic models and data-storage in computers.

In this book the authors introduce the new notion of superbimatrices and generalize it to super trimatrices and super n-matrices. Study of these structures is not only interesting and innovative but is also best suited for the computerized world.

The main difference between simple bimatrices and super bimatrices is that in case of simple bimatrices we have only one type of product defined on them, whereas in case of superbimatrices we have different types products called minor and major defined using them.

This book has four chapters. Chapter one describes the basics concepts to make this book a self contained one. Superbimatrices, semi superbimatrices, symmetric superbimatrices are introduced in chapter two. Chapter three introduces the notion of super trimatrices and the products defined using them. Chapter four gives the most generalized form of superbimatrix, viz. super n-matrix.

This book has given several examples so as to make the reader understand this new concept. Further minor and major by product defined using these new concepts are illustrated by examples. These algebraic structures are best suited in data storage in computers.

They are also useful in constructing multi expert super models

Finally it is an immense pleasure to thank Dr.K.Kandasamy for proof-reading and Kama, Meena and Rahul without whose help the book would have been an impossibility.

We dedicate this book to the millions of Tamil children in Sri Lanka who have died or become disabled and displaced due to the recent Sri Lankan war.

W.B.VASANTHA KANDASAMY
FLORENTIN SMARANDACHE

Chapter One

# BASIC CONCEPTS

In this chapter we just recall the definition of supermatrix and some of its basic properties which comprises the section 1. In section two bimatrices and their generalizations are introduced.

## 1.1 Supermatrices

The general rectangular or square array of numbers such as

$$A = \begin{bmatrix} 2 & 3 & 1 & 4 \\ -5 & 0 & 7 & -8 \end{bmatrix}, \quad B = \begin{bmatrix} 1 & 2 & 3 \\ -4 & 5 & 6 \\ 7 & -8 & 11 \end{bmatrix},$$

$$C = [3, 1, 0, -1, -2] \text{ and } D = \begin{bmatrix} -7/2 \\ 0 \\ \sqrt{2} \\ 5 \\ -41 \end{bmatrix}$$

are known as matrices.

We shall call them as simple matrices [19]. By a simple matrix we mean a matrix each of whose elements are just an ordinary number or a letter that stands for a number. In other words, the elements of a simple matrix are scalars or scalar quantities.

A supermatrix on the other hand is one whose elements are themselves matrices with elements that can be either scalars or other matrices. In general the kind of supermatrices we shall deal with in this book, the matrix elements which have any scalar for their elements. Suppose we have the four matrices;

$$a_{11} = \begin{bmatrix} 2 & -4 \\ 0 & 1 \end{bmatrix}, \ a_{12} = \begin{bmatrix} 0 & 40 \\ 21 & -12 \end{bmatrix}$$

$$a_{21} = \begin{bmatrix} 3 & -1 \\ 5 & 7 \\ -2 & 9 \end{bmatrix} \text{ and } a_{22} = \begin{bmatrix} 4 & 12 \\ -17 & 6 \\ 3 & 11 \end{bmatrix}.$$

One can observe the change in notation $a_{ij}$ denotes a matrix and not a scalar of a matrix ($1 \leq i, j \leq 2$).

Let

$$a = \begin{bmatrix} a_{11} & a_{12} \\ a_{21} & a_{22} \end{bmatrix};$$

we can write out the matrix a in terms of the original matrix elements i.e.,

$$a = \left[\begin{array}{cc|cc} 2 & -4 & 0 & 40 \\ 0 & 1 & 21 & -12 \\ \hline 3 & -1 & 4 & 12 \\ 5 & 7 & -17 & 6 \\ -2 & 9 & 3 & 11 \end{array}\right].$$

Here the elements are divided vertically and horizontally by thin lines. If the lines were not used the matrix a would be read as a simple matrix.

Thus far we have referred to the elements in a supermatrix as matrices as elements. It is perhaps more usual to call the elements of a supermatrix as submatrices. We speak of the submatrices within a supermatrix. Now we proceed on to define the order of a supermatrix.

The order of a supermatrix is defined in the same way as that of a simple matrix. The height of a supermatrix is the number of rows of submatrices in it. The width of a supermatrix is the number of columns of submatrices in it.

All submatrices with in a given row must have the same number of rows. Likewise all submatrices with in a given column must have the same number of columns.

A diagrammatic representation is given by the following figure.

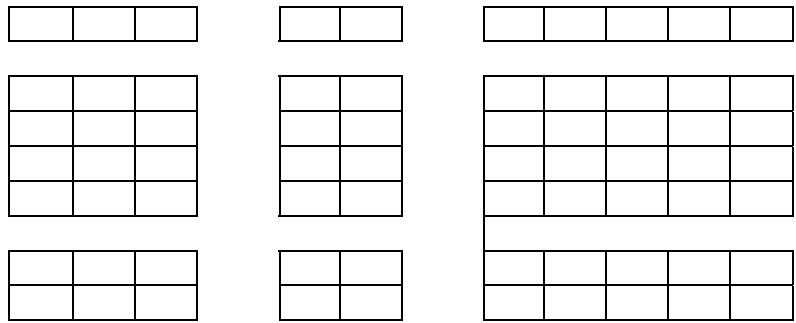

In the first row of rectangles we have one row of a square for each rectangle; in the second row of rectangles we have four rows of squares for each rectangle and in the third row of rectangles we have two rows of squares for each rectangle. Similarly for the first column of rectangles three columns of squares for each rectangle. For the second column of rectangles we have two column of squares for each rectangle, and for the third column of rectangles we have five columns of squares for each rectangle.

Thus we have for this supermatrix 3 rows and 3 columns.

One thing should now be clear from the definition of a supermatrix. The super order of a supermatrix tells us nothing about the simple order of the matrix from which it was obtained by partitioning. Furthermore, the order of supermatrix tells us nothing about the orders of the submatrices within that supermatrix.

Now we illustrate the number of rows and columns of a supermatrix.

***Example 1.1.1:*** Let

$$a = \left[\begin{array}{c|cccc} 3 & 3 & 0 & 1 & 4 \\ -1 & 2 & 1 & -1 & 6 \\ \hline 0 & 3 & 4 & 5 & 6 \\ 1 & 7 & 8 & -9 & 0 \\ 2 & 1 & 2 & 3 & -4 \end{array}\right].$$

a is a supermatrix with two rows and two columns.

Now we proceed on to define the notion of partitioned matrices. It is always possible to construct a supermatrix from any simple matrix that is not a scalar quantity.

The supermatrix can be constructed from a simple matrix this process of constructing supermatrix is called the partitioning.

A simple matrix can be partitioned by dividing or separating the matrix between certain specified rows, or the procedure may be reversed. The division may be made first between rows and then between columns.

We illustrate this by a simple example.

***Example 1.1.2:*** Let

$$A = \begin{bmatrix} 3 & 0 & 1 & 1 & 2 & 0 \\ 1 & 0 & 0 & 3 & 5 & 2 \\ 5 & -1 & 6 & 7 & 8 & 4 \\ 0 & 9 & 1 & 2 & 0 & -1 \\ 2 & 5 & 2 & 3 & 4 & 6 \\ 1 & 6 & 1 & 2 & 3 & 9 \end{bmatrix}$$

is a 6 × 6 simple matrix with real numbers as elements.

$$A_1 = \begin{bmatrix} 3 & 0 & | & 1 & 1 & 2 & 0 \\ 1 & 0 & | & 0 & 3 & 5 & 2 \\ 5 & -1 & | & 6 & 7 & 8 & 4 \\ 0 & 9 & | & 1 & 2 & 0 & -1 \\ 2 & 5 & | & 2 & 3 & 4 & 6 \\ 1 & 6 & | & 1 & 2 & 3 & 9 \end{bmatrix}.$$

Now let us draw a thin line between the 2$^{nd}$ and 3$^{rd}$ columns.

This gives us the matrix $A_1$. Actually $A_1$ may be regarded as a supermatrix with two matrix elements forming one row and two columns.

Now consider

$$A_2 = \begin{bmatrix} 3 & 0 & 1 & 1 & 2 & 0 \\ 1 & 0 & 0 & 3 & 5 & 2 \\ 5 & -1 & 6 & 7 & 8 & 4 \\ 0 & 9 & 1 & 2 & 0 & -1 \\ \hline 2 & 5 & 2 & 3 & 4 & 6 \\ 1 & 6 & 1 & 2 & 3 & 9 \end{bmatrix}$$

Draw a thin line between the rows 4 and 5 which gives us the new matrix $A_2$. $A_2$ is a supermatrix with two rows and one column.

Now consider the matrix

$$A_3 = \begin{bmatrix} 3 & 0 & | & 1 & 1 & 2 & 0 \\ 1 & 0 & | & 0 & 3 & 5 & 2 \\ 5 & -1 & | & 6 & 7 & 8 & 4 \\ 0 & 9 & | & 1 & 2 & 0 & -1 \\ \hline 2 & 5 & | & 2 & 3 & 4 & 6 \\ 1 & 6 & | & 1 & 2 & 3 & 9 \end{bmatrix},$$

$A_3$ is now a second order supermatrix with two rows and two columns. We can simply write $A_3$ as

$$\begin{bmatrix} a_{11} & a_{12} \\ a_{21} & a_{22} \end{bmatrix}$$

where

$$a_{11} = \begin{bmatrix} 3 & 0 \\ 1 & 0 \\ 5 & -1 \\ 0 & 9 \end{bmatrix},$$

$$a_{12} = \begin{bmatrix} 1 & 1 & 2 & 0 \\ 0 & 3 & 5 & 2 \\ 6 & 7 & 8 & 4 \\ 1 & 2 & 0 & -1 \end{bmatrix},$$

$$a_{21} = \begin{bmatrix} 2 & 5 \\ 1 & 6 \end{bmatrix} \text{ and } a_{22} = \begin{bmatrix} 2 & 3 & 4 & 6 \\ 1 & 2 & 3 & 9 \end{bmatrix}.$$

The elements now are the submatrices defined as $a_{11}$, $a_{12}$, $a_{21}$ and $a_{22}$ and therefore $A_3$ is in terms of letters.

According to the methods we have illustrated a simple matrix can be partitioned to obtain a supermatrix in any way that happens to suit our purposes.

The natural order of a supermatrix is usually determined by the natural order of the corresponding simple matrix. Further more we are not usually concerned with natural order of the submatrices within a supermatrix.

Now we proceed on to recall the notion of symmetric partition, for more information about these concepts please refer [19]. By a symmetric partitioning of a matrix we mean that the rows and columns are partitioned in exactly the same way. If the matrix is partitioned between the first and second column and between the third and fourth column, then to be symmetrically partitioning, it must also be partitioned between the first and second rows and third and fourth rows. According to this rule of symmetric partitioning only square simple matrix can be symmetrically partitioned. We give an example of a symmetrically partitioned matrix $a_s$,

*Example 1.1.3:* Let

$$a_s = \begin{bmatrix} 2 & 3 & 4 & 1 \\ 5 & 6 & 9 & 2 \\ \hline 0 & 6 & 1 & 9 \\ 5 & 1 & 1 & 5 \end{bmatrix}.$$

Here we see that the matrix has been partitioned between columns one and two and three and four. It has also been partitioned between rows one and two and rows three and four.

Now we just recall from [19] the method of symmetric partitioning of a symmetric simple matrix.

*Example 1.1.4:* Let us take a fourth order symmetric matrix and partition it between the second and third rows and also between the second and third columns.

$$a = \begin{bmatrix} 4 & 3 & 2 & 7 \\ 3 & 6 & 1 & 4 \\ \hline 2 & 1 & 5 & 2 \\ 7 & 4 & 2 & 7 \end{bmatrix}.$$

We can represent this matrix as a supermatrix with letter elements.

$$a_{11} = \begin{bmatrix} 4 & 3 \\ 3 & 6 \end{bmatrix}, a_{12} = \begin{bmatrix} 2 & 7 \\ 1 & 4 \end{bmatrix}$$

$$a_{21} = \begin{bmatrix} 2 & 1 \\ 7 & 4 \end{bmatrix} \text{ and } a_{22} = \begin{bmatrix} 5 & 2 \\ 2 & 7 \end{bmatrix},$$

so that

$$a = \begin{bmatrix} a_{11} & a_{12} \\ a_{21} & a_{22} \end{bmatrix}.$$

The diagonal elements of the supermatrix a are $a_{11}$ and $a_{22}$. We also observe the matrices $a_{11}$ and $a_{22}$ are also symmetric matrices.

The non diagonal elements of this supermatrix a are the matrices $a_{12}$ and $a_{21}$. Clearly $a_{21}$ is the transpose of $a_{12}$.

The simple rule about the matrix element of a symmetrically partitioned symmetric simple matrix are (1) The diagonal submatrices of the supermatrix are all symmetric matrices. (2) The matrix elements below the diagonal are the transposes of the corresponding elements above the diagonal.

The forth order supermatrix obtained from a symmetric partitioning of a symmetric simple matrix a is as follows.

$$a = \begin{bmatrix} a_{11} & a_{12} & a_{13} & a_{14} \\ a'_{12} & a_{22} & a_{23} & a_{24} \\ a'_{13} & a'_{23} & a_{33} & a_{34} \\ a'_{14} & a'_{24} & a'_{34} & a_{44} \end{bmatrix}.$$

How to express that a symmetric matrix has been symmetrically partitioned (i) $a_{11}$ and $a^t_{11}$ are equal. (ii) $a^t_{ij}$ ($i \neq j$); $a^t_{ij} = a_{ji}$ and $a^t_{ji} = a_{ij}$. Thus the general expression for a symmetrically partitioned symmetric matrix;

$$a = \begin{bmatrix} a_{11} & a_{12} & \ldots & a_{1n} \\ a'_{12} & a_{22} & \ldots & a_{2n} \\ \vdots & \vdots & & \vdots \\ a'_{1n} & a'_{2n} & \ldots & a_{nn} \end{bmatrix}.$$

If we want to indicate a symmetrically partitioned simple diagonal matrix we would write

$$D = \begin{bmatrix} D_1 & 0 & \cdots & 0 \\ 0' & D_2 & \cdots & 0 \\ & & & \\ 0' & 0' & \cdots & D_n \end{bmatrix}$$

0' only represents the order is reversed or transformed. We denote $a_{ij}^t = a'_{ij}$ just the ' means the transpose.

D will be referred to as the super diagonal matrix. The identity matrix

$$I = \begin{bmatrix} I_s & 0 & 0 \\ 0 & I_t & 0 \\ 0 & 0 & I_r \end{bmatrix}$$

s, t and r denote the number of rows and columns of the first second and third identity matrices respectively (zeros denote matrices with zero as all entries).

*Example 1.1.5:* We just illustrate a general super diagonal matrix d;

$$d = \left[\begin{array}{ccc|cc} 3 & 1 & 2 & 0 & 0 \\ 5 & 6 & 0 & 0 & 0 \\ \hline 0 & 0 & 0 & 2 & 5 \\ 0 & 0 & 0 & -1 & 3 \\ 0 & 0 & 0 & 9 & 10 \end{array}\right]$$

i.e., $d = \begin{bmatrix} m_1 & 0 \\ 0 & m_2 \end{bmatrix}$.

An example of a super diagonal matrix with vector elements is given, which can be useful in experimental designs.

*Example 1.1.6:* Let

$$\begin{bmatrix} 1 & 0 & 0 & 0 \\ 1 & 0 & 0 & 0 \\ 1 & 0 & 0 & 0 \\ \hline 0 & 1 & 0 & 0 \\ 0 & 1 & 0 & 0 \\ \hline 0 & 0 & 1 & 0 \\ 0 & 0 & 1 & 0 \\ 0 & 0 & 1 & 0 \\ 0 & 0 & 1 & 0 \\ \hline 0 & 0 & 0 & 1 \\ 0 & 0 & 0 & 1 \\ 0 & 0 & 0 & 1 \\ 0 & 0 & 0 & 1 \end{bmatrix}.$$

Here the diagonal elements are only column unit vectors. In case of supermatrix [19] has defined the notion of partial triangular matrix as a supermatrix.

*Example 1.1.7:* Let

$$u = \begin{bmatrix} 2 & 1 & 1 & 3 & 2 \\ 0 & 5 & 2 & 1 & 1 \\ 0 & 0 & 1 & 0 & 2 \end{bmatrix}$$

u is a partial upper triangular supermatrix.

*Example 1.1.8:* Let

$$L = \begin{bmatrix} 5 & 0 & 0 & 0 & 0 \\ 7 & 2 & 0 & 0 & 0 \\ 1 & 2 & 3 & 0 & 0 \\ 4 & 5 & 6 & 7 & 0 \\ 1 & 2 & 5 & 2 & 6 \\ \hline 1 & 2 & 3 & 4 & 5 \\ 0 & 1 & 0 & 1 & 0 \end{bmatrix};$$

L is partial upper triangular matrix partitioned as a supermatrix.

Thus $T = \begin{bmatrix} T \\ \overline{a'} \end{bmatrix}$ where T is the lower triangular submatrix, with

$$T = \begin{bmatrix} 5 & 0 & 0 & 0 & 0 \\ 7 & 2 & 0 & 0 & 0 \\ 1 & 2 & 3 & 0 & 0 \\ 4 & 5 & 6 & 7 & 0 \\ 1 & 2 & 5 & 2 & 6 \end{bmatrix} \text{ and } a' = \begin{bmatrix} 1 & 2 & 3 & 4 & 5 \\ 0 & 1 & 0 & 1 & 0 \end{bmatrix}.$$

We proceed on to define the notion of supervectors i.e., Type I column supervector. A simple vector is a vector each of whose elements is a scalar. It is nice to see the number of different types of supervectors given by [19].

*Example 1.1.9:* Let

$$v = \begin{bmatrix} 1 \\ 3 \\ 4 \\ \overline{5} \\ 7 \end{bmatrix}.$$

This is a type I i.e., type one column supervector.

$$v = \begin{bmatrix} v_1 \\ v_2 \\ \vdots \\ v_n \end{bmatrix}$$

where each $v_i$ is a column subvectors of the column vector v.

Type I row supervector is given by the following example.

**Example 1.1.10:** $v^1 = [2\ 3\ 1\ |\ 5\ 7\ 8\ 4]$ is a type I row supervector. i.e., $v' = [v'_1, v'_2, \ldots, v'_n]$ where each $v'_i$ is a row subvector; $1 \leq i \leq n$.

Next we recall the definition of type II supervectors.

Type II column supervectors.

**DEFINITION 1.1.1:** *Let*

$$a = \begin{bmatrix} a_{11} & a_{12} & \ldots & a_{1m} \\ a_{21} & a_{22} & \ldots & a_{2m} \\ \ldots & \ldots & \ldots & \ldots \\ a_{n1} & a_{n2} & \ldots & a_{nm} \end{bmatrix}$$

$$\begin{aligned} a_1^1 &= [a_{11} \ldots a_{1m}] \\ a_2^1 &= [a_{21} \ldots a_{2m}] \\ &\ldots \\ a_n^1 &= [a_{n1} \ldots a_{nm}] \end{aligned}$$

*i.e.,* $\quad a = \begin{bmatrix} a_1^1 \\ a_2^1 \\ \vdots \\ a_n^1 \end{bmatrix}_m$

*is defined to be the type II column supervector.*
*Similarly if*

$$a^1 = \begin{bmatrix} a_{11} \\ a_{21} \\ \vdots \\ a_{n1} \end{bmatrix},\ a^2 = \begin{bmatrix} a_{12} \\ a_{22} \\ \vdots \\ a_{n2} \end{bmatrix},\ \ldots,\ a^m = \begin{bmatrix} a_{1m} \\ a_{2m} \\ \vdots \\ a_{nm} \end{bmatrix}.$$

*Hence now* $a = [a^1\ a^2\ \ldots\ a^m]_n$ *is defined to be the type II row supervector.*

*Clearly*

$$a = \begin{bmatrix} a_1^1 \\ a_2^1 \\ \vdots \\ a_n^1 \end{bmatrix}_m = [a^1 \ a^2 \ ... \ a^m]_n$$

*the equality of supermatrices.*

***Example 1.1.11:*** Let

$$A = \begin{bmatrix} 3 & 6 & 0 & 4 & 5 \\ 2 & 1 & 6 & 3 & 0 \\ 1 & 1 & 1 & 2 & 1 \\ 0 & 1 & 0 & 1 & 0 \\ 2 & 0 & 1 & 2 & 1 \end{bmatrix}$$

be a simple matrix. Let a and b the supermatrix made from A.

$$a = \left[ \begin{array}{ccc|cc} 3 & 6 & 0 & 4 & 5 \\ 2 & 1 & 6 & 3 & 0 \\ 1 & 1 & 1 & 2 & 1 \\ \hline 0 & 1 & 0 & 1 & 0 \\ 2 & 0 & 1 & 2 & 1 \end{array} \right]$$

where

$$a_{11} = \begin{bmatrix} 3 & 6 & 0 \\ 2 & 1 & 6 \\ 1 & 1 & 1 \end{bmatrix}, a_{12} = \begin{bmatrix} 4 & 5 \\ 3 & 0 \\ 2 & 1 \end{bmatrix},$$

$$a_{21} = \begin{bmatrix} 0 & 1 & 0 \\ 2 & 0 & 1 \end{bmatrix} \text{ and } a_{22} = \begin{bmatrix} 1 & 0 \\ 2 & 1 \end{bmatrix}.$$

i.e., $$a = \begin{bmatrix} a_{11} & a_{12} \\ a_{21} & a_{22} \end{bmatrix}.$$

$$b = \begin{bmatrix} 3 & 6 & 0 & 4 & 5 \\ 2 & 1 & 6 & 3 & 0 \\ 1 & 1 & 1 & 2 & 1 \\ 0 & 1 & 0 & 1 & 0 \\ \hline 2 & 0 & 1 & 2 & 1 \end{bmatrix} = \begin{bmatrix} b_{11} & b_{12} \\ b_{21} & b_{22} \end{bmatrix}$$

where

$$b_{11} = \begin{bmatrix} 3 & 6 & 0 & 4 \\ 2 & 1 & 6 & 3 \\ 1 & 1 & 1 & 2 \\ 0 & 1 & 0 & 1 \end{bmatrix}, \quad b_{12} = \begin{bmatrix} 5 \\ 0 \\ 1 \\ 0 \end{bmatrix},$$

$b_{21} = [2\ 0\ 1\ 2\ ]$ and $b_{22} = [1]$.

$$a = \begin{bmatrix} 3 & 6 & 0 & 4 & 5 \\ 2 & 1 & 6 & 3 & 0 \\ 1 & 1 & 1 & 2 & 1 \\ \hline 0 & 1 & 0 & 1 & 0 \\ 2 & 0 & 1 & 2 & 1 \end{bmatrix}$$

and

$$b = \begin{bmatrix} 3 & 6 & 0 & 4 & 5 \\ 2 & 1 & 6 & 3 & 0 \\ 1 & 1 & 1 & 2 & 1 \\ 0 & 1 & 0 & 1 & 0 \\ \hline 2 & 0 & 1 & 2 & 1 \end{bmatrix}.$$

We see that the corresponding scalar elements for matrix a and matrix b are identical. Thus two supermatrices are equal if and only if their corresponding simple forms are equal.

Now we give examples of type III supervector for more refer [19].

*Example 1.1.12:*

$$a = \begin{bmatrix} 3 & 2 & 1 & 7 & 8 \\ 0 & 2 & 1 & 6 & 9 \\ 0 & 0 & 5 & 1 & 2 \end{bmatrix} = [T' \mid a']$$

and

$$b = \begin{bmatrix} 2 & 0 & 0 \\ 9 & 4 & 0 \\ 8 & 3 & 6 \\ \hline 5 & 2 & 9 \\ 4 & 7 & 3 \end{bmatrix} = \left[\frac{T}{b'}\right]$$

are type III supervectors.

One interesting and common example of a type III supervector is a prediction data matrix having both predictor and criterion attributes.

The next interesting notion about supermatrix is its transpose. First we illustrate this by an example before we give the general case.

*Example 1.1.13:* Let

$$a = \begin{bmatrix} 2 & 1 & 3 & 5 & 6 \\ 0 & 2 & 0 & 1 & 1 \\ 1 & 1 & 1 & 0 & 2 \\ \hline 2 & 2 & 0 & 1 & 1 \\ 5 & 6 & 1 & 0 & 1 \\ \hline 2 & 0 & 0 & 0 & 4 \\ 1 & 0 & 1 & 1 & 5 \end{bmatrix}$$

$$= \begin{bmatrix} a_{11} & a_{12} \\ a_{21} & a_{22} \\ a_{31} & a_{32} \end{bmatrix}$$

where

$$a_{11} = \begin{bmatrix} 2 & 1 & 3 \\ 0 & 2 & 0 \\ 1 & 1 & 1 \end{bmatrix}, a_{12} = \begin{bmatrix} 5 & 6 \\ 1 & 1 \\ 0 & 2 \end{bmatrix},$$

$$a_{21} = \begin{bmatrix} 2 & 2 & 0 \\ 5 & 6 & 1 \end{bmatrix}, a_{22} = \begin{bmatrix} 1 & 1 \\ 0 & 1 \end{bmatrix},$$

$$a_{31} = \begin{bmatrix} 2 & 0 & 0 \\ 1 & 0 & 1 \end{bmatrix} \text{ and } a_{32} = \begin{bmatrix} 0 & 4 \\ 1 & 5 \end{bmatrix}.$$

The transpose of a

$$a^t = a' = \left[\begin{array}{ccc|ccc|cc} 2 & 0 & 1 & 2 & 5 & 2 & 1 \\ 1 & 2 & 1 & 2 & 6 & 0 & 0 \\ 3 & 0 & 1 & 0 & 1 & 0 & 1 \\ \hline 5 & 1 & 0 & 1 & 0 & 0 & 1 \\ 6 & 1 & 2 & 1 & 1 & 4 & 5 \end{array}\right].$$

Let us consider the transposes of $a_{11}$, $a_{12}$, $a_{21}$, $a_{22}$, $a_{31}$ and $a_{32}$.

$$a'_{11} = a^t_{11} = \begin{bmatrix} 2 & 0 & 1 \\ 1 & 2 & 1 \\ 3 & 0 & 1 \end{bmatrix}$$

$$a'_{12} = a^t_{12} = \begin{bmatrix} 5 & 1 & 0 \\ 6 & 1 & 2 \end{bmatrix}$$

$$a'_{21} = a^t_{21} = \begin{bmatrix} 2 & 5 \\ 2 & 6 \\ 0 & 1 \end{bmatrix}$$

$$a'_{31} = a^t_{31} = \begin{bmatrix} 2 & 1 \\ 0 & 0 \\ 0 & 1 \end{bmatrix}$$

$$a'_{22} = a^t_{22} = \begin{bmatrix} 1 & 0 \\ 1 & 1 \end{bmatrix}$$

$$a'_{32} = a^t_{32} = \begin{bmatrix} 0 & 1 \\ 4 & 5 \end{bmatrix}.$$

$$a' = \begin{bmatrix} a'_{11} & a'_{21} & a'_{31} \\ a'_{12} & a'_{22} & a'_{32} \end{bmatrix}.$$

Now we describe the general case. Let

$$a = \begin{bmatrix} a_{11} & a_{12} & \cdots & a_{1m} \\ a_{21} & a_{22} & \cdots & a_{2m} \\ \vdots & \vdots & & \vdots \\ a_{n1} & a_{n2} & \cdots & a_{nm} \end{bmatrix}$$

be a n × m supermatrix. The transpose of the supermatrix a denoted by

$$a' = \begin{bmatrix} a'_{11} & a'_{21} & \cdots & a'_{n1} \\ a'_{12} & a'_{22} & \cdots & a'_{n2} \\ \vdots & \vdots & & \vdots \\ a'_{1m} & a'_{2m} & \cdots & a'_{nm} \end{bmatrix}.$$

a' is a m by n supermatrix obtained by taking the transpose of each element i.e., the submatrices of a.

Now we will find the transpose of a symmetrically partitioned symmetric simple matrix. Let a be the symmetrically partitioned symmetric simple matrix.

Let a be a m × m symmetric supermatrix i.e.,

$$a = \begin{bmatrix} a_{11} & a_{21} & \cdots & a_{m1} \\ a_{12} & a_{22} & \cdots & a_{m2} \\ \vdots & \vdots & & \vdots \\ a_{1m} & a_{2m} & \cdots & a_{mm} \end{bmatrix}$$

the transpose of the supermatrix is given by a'

$$a' = \begin{bmatrix} a'_{11} & (a'_{12})' & \cdots & (a'_{1m})' \\ a'_{12} & a'_{22} & \cdots & (a'_{2m})' \\ \vdots & \vdots & & \vdots \\ a'_{1m} & a'_{2m} & \cdots & a'_{mm} \end{bmatrix}$$

The diagonal matrix $a_{11}$ are symmetric matrices so are unaltered by transposition. Hence
$$a'_{11} = a_{11}, a'_{22} = a_{22}, \ldots, a'_{mm} = a_{mm}.$$

Recall also the transpose of a transpose is the original matrix. Therefore
$$(a'_{12})' = a_{12}, (a'_{13})' = a_{13}, \ldots, (a'_{ij})' = a_{ij}.$$

Thus the transpose of supermatrix constructed by symmetrically partitioned symmetric simple matrix a of a' is given by

$$a' = \begin{bmatrix} a_{11} & a_{12} & \cdots & a_{1m} \\ a'_{21} & a_{22} & \cdots & a_{2m} \\ \vdots & \vdots & & \vdots \\ a'_{1m} & a'_{2m} & \cdots & a_{mm} \end{bmatrix}.$$

Thus a = a'.

Similarly transpose of a symmetrically partitioned diagonal matrix is simply the original diagonal supermatrix itself;

i.e., if

$$D = \begin{bmatrix} d_1 & & & \\ & d_2 & & \\ & & \ddots & \\ & & & d_n \end{bmatrix}$$

$$D' = \begin{bmatrix} d'_1 & & & \\ & d'_2 & & \\ & & \ddots & \\ & & & d'_n \end{bmatrix}$$

$d'_1 = d_1$, $d'_2 = d_2$ etc. Thus $D = D'$.

Now we see the transpose of a type I supervector.

*Example 1.1.14:* Let

$$V = \begin{bmatrix} 3 \\ 1 \\ 2 \\ \hline 4 \\ 5 \\ 7 \\ \hline 5 \\ 1 \end{bmatrix}$$

The transpose of V denoted by V' or $V^t$ is

$$V' = [3\ 1\ 2\ |\ 4\ 5\ 7\ |\ 5\ 1].$$

If

$$V = \begin{bmatrix} v_1 \\ v_2 \\ v_3 \end{bmatrix}$$

where

$$v_1 = \begin{bmatrix} 3 \\ 1 \\ 2 \end{bmatrix}, v_2 = \begin{bmatrix} 4 \\ 5 \\ 7 \end{bmatrix} \text{ and } v_3 = \begin{bmatrix} 5 \\ 1 \end{bmatrix}$$

$$V' = [v'_1 \ v'_2 \ v'_3].$$

Thus if

$$V = \begin{bmatrix} v_1 \\ v_2 \\ \vdots \\ v_n \end{bmatrix}$$

then

$$V' = [v'_1 \ v'_2 \ \ldots \ v'_n].$$

*Example 1.1.15:* Let

$$t = \left[\begin{array}{cccc|cc} 3 & 0 & 1 & 1 & 5 & 2 \\ 4 & 2 & 0 & 1 & 3 & 5 \\ 1 & 0 & 1 & 0 & 1 & 6 \end{array}\right]$$

= [T | a ]. The transpose of t

$$\text{i.e., } t' = \left[\begin{array}{ccc} 3 & 4 & 1 \\ 0 & 2 & 0 \\ 1 & 0 & 1 \\ 1 & 1 & 0 \\ \hline 5 & 3 & 1 \\ 2 & 5 & 6 \end{array}\right] = \left[\begin{array}{c} T' \\ \hline a' \end{array}\right].$$

The addition of supermatrices may not be always be defined.

*Example 1.1.16:* For instance let

$$a = \begin{bmatrix} a_{11} & a_{12} \\ a_{21} & a_{22} \end{bmatrix}$$

and

$$b = \begin{bmatrix} b_{11} & b_{12} \\ b_{21} & b_{22} \end{bmatrix}$$

where

$$a_{11} = \begin{bmatrix} 3 & 0 \\ 1 & 2 \end{bmatrix}, \quad a_{12} = \begin{bmatrix} 1 \\ 7 \end{bmatrix}$$

$$a_{21} = [4 \; 3], \quad a_{22} = [6].$$

$$b_{11} = [2], \quad b_{12} = [1 \; 3]$$

$$b_{21} = \begin{bmatrix} 5 \\ 2 \end{bmatrix} \text{ and } b_{22} = \begin{bmatrix} 4 & 1 \\ 0 & 2 \end{bmatrix}.$$

It is clear both a and b are second order square supermatrices but here we cannot add together the corresponding matrix elements of a and b because the submatrices do not have the same order.

Now we proceed onto recall the definition of minor product of two supervectors.

Suppose

$$v_a = \begin{bmatrix} v_{a_1} \\ v_{a_2} \\ \vdots \\ v_{a_n} \end{bmatrix} \text{ and } v_b = \begin{bmatrix} v_{b_1} \\ v_{b_2} \\ \vdots \\ v_{b_n} \end{bmatrix}.$$

The minor product of these two supervectors $v_a$ and $v_b$ is given by

$$= v'_a v_b = \begin{bmatrix} v'_{a_1} & v'_{a_2} & \cdots & v'_{a_n} \end{bmatrix} \begin{bmatrix} v_{b_1} \\ v_{b_2} \\ \vdots \\ v_{b_n} \end{bmatrix}$$

$$= v'_{a_1} v_{b_1} + v'_{a_2} v_{b_2} + \cdots + v'_{a_n} v_{b_n}.$$

We illustrate this by the following example.

**Example 1.1.17:** Let $V_a$ and $V_b$ be two type I supervectors where

$$V_a = \begin{bmatrix} v_{a_1} \\ v_{a_2} \\ v_{a_3} \end{bmatrix}$$

with

$$v_{a_1} = \begin{bmatrix} 0 \\ 1 \\ 2 \end{bmatrix}, \quad v_{a_2} = \begin{bmatrix} 4 \\ 0 \\ 1 \\ -1 \end{bmatrix} \text{ and } v_{a_3} = \begin{bmatrix} 1 \\ 2 \end{bmatrix}.$$

Let

$$V_b = \begin{bmatrix} v_{b_1} \\ v_{b_2} \\ v_{b_3} \end{bmatrix}$$

where

$$v_{b_1} = \begin{bmatrix} 1 \\ -1 \\ 0 \end{bmatrix}, \quad v_{b_2} = \begin{bmatrix} -4 \\ 1 \\ 2 \\ 0 \end{bmatrix} \text{ and } v_{b_3} = \begin{bmatrix} -1 \\ 1 \end{bmatrix}.$$

$$V'_a V_b = \begin{bmatrix} v'_{a_1} & v'_{a_2} & v'_{a_3} \end{bmatrix} \begin{bmatrix} V_{b_1} \\ V_{b_2} \\ V_{b_3} \end{bmatrix}$$

$$= \quad v'_{a_1} v_{b_1} + v'_{a_2} v_{b_2} + v'_{a_3} v_{b_3}$$

$$= \quad \begin{bmatrix} 0 & 1 & 2 \end{bmatrix} \begin{bmatrix} 1 \\ -1 \\ 0 \end{bmatrix} + \begin{bmatrix} 4 & 0 & 1 & -1 \end{bmatrix} \begin{bmatrix} -4 \\ 1 \\ 2 \\ 0 \end{bmatrix} + \begin{bmatrix} 1 & 2 \end{bmatrix} \begin{bmatrix} -1 \\ 1 \end{bmatrix}$$

$$= \quad -1 + (-16+2) + (-1+2)$$
$$= \quad -1 - 16 + 2 - 1 + 2$$
$$= \quad -14.$$

It is easily proved $V'_a V_b = V'_b V_a$.

Now we proceed on to recall the definition of major product of type I supervectors.

Suppose

$$V_a = \begin{bmatrix} v_{a_1} \\ v_{a_2} \\ \vdots \\ v_{a_n} \end{bmatrix} \text{ and } V_b = \begin{bmatrix} v_{b_1} \\ v_{b_2} \\ \vdots \\ v_{b_m} \end{bmatrix}$$

be any two supervectors of type I. The major product is defined as

$$V_a V'_b = \begin{bmatrix} v_{a_1} \\ v_{a_2} \\ \vdots \\ v_{a_n} \end{bmatrix} \cdot \begin{bmatrix} v'_{b_1} & v'_{b_2} & \cdots & v'_{b_m} \end{bmatrix}$$

$$= \begin{bmatrix} v_{a_1} v'_{b_1} & v_{a_1} v'_{b_2} & \cdots & v_{a_1} v'_{b_m} \\ v_{a_2} v'_{b_1} & v_{a_2} v'_{b_2} & \cdots & v_{a_2} v'_{b_m} \\ \vdots & & & \\ v_{a_n} v'_{b_1} & v_{a_n} v'_{b_2} & \cdots & v_{a_n} v'_{b_m} \end{bmatrix}.$$

Now we illustrate this by the following example.

*Example 1.1.18:* Let

$$V_a = \begin{bmatrix} v_{a_1} \\ v_{a_2} \\ v_{a_3} \end{bmatrix} \text{ and } V_b = \begin{bmatrix} v_{b_1} \\ v_{b_2} \\ v_{b_3} \\ v_{b_4} \end{bmatrix}$$

where

$$v_{a_1} = [2], \ v_{a_2} = \begin{bmatrix} 1 \\ -1 \end{bmatrix} \text{ and } v_{a_3} = \begin{bmatrix} 1 \\ 2 \\ 0 \end{bmatrix}$$

and

$$v_{b_1} = \begin{bmatrix} 3 \\ 1 \\ 2 \end{bmatrix}, \ v_{b_2} = \begin{bmatrix} 1 \\ 2 \end{bmatrix}, \ v_{b_3} = \begin{bmatrix} 3 \\ 4 \\ -1 \\ 0 \end{bmatrix} \text{ and } v_{b_4} = [5].$$

$$V_a V'_b = \begin{bmatrix} 2 \\ 1 \\ -1 \\ 1 \\ 2 \\ 0 \end{bmatrix} \begin{bmatrix} 3 & 1 & 2 \mid 1 & 2 \mid 3 & 4 & -1 & 0 \mid 5 \end{bmatrix}$$

$$= \begin{bmatrix} [2][3\ 1\ 2] & [2][1\ 2] & [2][3\ 4\ -1\ 0] & [2][5] \\ \begin{bmatrix}1\\-1\end{bmatrix}[3\ 1\ 2] & \begin{bmatrix}1\\-1\end{bmatrix}[1\ 2] & \begin{bmatrix}1\\-1\end{bmatrix}[3\ 4\ -1\ 0] & \begin{bmatrix}1\\-1\end{bmatrix}[5] \\ \begin{bmatrix}1\\2\\0\end{bmatrix}[3\ 1\ 2] & \begin{bmatrix}1\\2\\0\end{bmatrix}[1\ 2] & \begin{bmatrix}1\\2\\0\end{bmatrix}[3\ 4\ -1\ 0] & \begin{bmatrix}1\\2\\0\end{bmatrix}[5] \end{bmatrix}$$

$$= \begin{bmatrix} 6 & 2 & 4 & 2 & 4 & 6 & 4 & -2 & 0 & 10 \\ 3 & 1 & 2 & 1 & 2 & 3 & 4 & -1 & 0 & 5 \\ -3 & -1 & -2 & -1 & -2 & -3 & -4 & 1 & 0 & -5 \\ 3 & 1 & 2 & 1 & 2 & 3 & 4 & -1 & 0 & 5 \\ 6 & 2 & 4 & 2 & 4 & 6 & 8 & -2 & 0 & 10 \\ 0 & 0 & 0 & 0 & 0 & 0 & 0 & 0 & 0 & 0 \end{bmatrix}.$$

We leave it for the reader to verify that $(V_a\ V'_b)' = V_b\ V'_a$.

*Example 1.1.19:* We just recall if

$$v = \begin{bmatrix} 3 \\ 4 \\ 7 \end{bmatrix}$$

is a column vector and v' the transpose of v is a row vector then we have

$$v'v = \begin{bmatrix} 3 & 4 & 7 \end{bmatrix} \begin{bmatrix} 3 \\ 4 \\ 7 \end{bmatrix}$$

$$= 3^2 + 4^2 + 7^2 = 74.$$

Thus if

$$V'_x = [x_1\ x_2\ \ldots\ x_n]$$

$$V'_x V_x = [x_1\ x_2\ \ldots\ x_n] \begin{bmatrix} x_1 \\ x_2 \\ \vdots \\ x_n \end{bmatrix}$$

$$= x_1^2 + x_2^2 + \ldots + x_n^2.$$

Also

$$[1\ 1\ \ldots\ 1\ ] \begin{bmatrix} x_1 \\ x_2 \\ \vdots \\ x_n \end{bmatrix} = [x_1 + x_2 + \ldots + x_n]$$

and

$$[x_1\ x_2\ \ldots\ x_n] \begin{bmatrix} 1 \\ 1 \\ \vdots \\ 1 \end{bmatrix} = [x_1 + x_2 + \ldots + x_n];$$

$$\text{i.e., } 1'v_x = v'_x 1 = \sum x_i$$

where

$$v_x = \begin{bmatrix} x_1 \\ x_2 \\ \vdots \\ x_n \end{bmatrix}$$

and

$$\sum x_i = x_1 + x_2 + \ldots + x_n.$$

We have the following types of products defined.

***Example 1.1.20:*** We have

$$[0\ 1\ 0\ 0]\begin{bmatrix}0\\1\\0\\0\end{bmatrix}=1,$$

$$[0\ 1\ 0\ 0]\begin{bmatrix}1\\0\\0\\0\end{bmatrix}=0,$$

$$[0\ 1\ 0\ 0]\begin{bmatrix}1\\1\\1\\1\end{bmatrix}=1$$

and

$$\begin{bmatrix}0\\1\\0\\0\end{bmatrix}[1\ 0\ 0]=\begin{bmatrix}0&0&0\\1&0&0\\0&0&0\\0&0&0\end{bmatrix}.$$

Recall

$$a=\begin{bmatrix}a_{11}&a_{12}&\cdots&a_{1m}\\a_{21}&a_{22}&\cdots&a_{2m}\\a_{n1}&a_{n2}&\cdots&a_{nm}\end{bmatrix}$$

we have

$$a=\begin{bmatrix}a_1^1\\a_2^1\\\vdots\\a_n^1\end{bmatrix}_m \quad (1)$$

and

$$a=[a^1\ a^2\ \ldots\ a^m]_n. \quad (2)$$

Now transpose of

$$a = \begin{bmatrix} a_1^1 \\ a_2^1 \\ \vdots \\ a_n^1 \end{bmatrix}_m$$

is given by the equation

$$a' = \begin{bmatrix} (a_1^1)' \ (a_2^1)' \cdots (a_n^1)' \end{bmatrix}_m$$

$$a' = \begin{bmatrix} (a^1)' \\ (a^2)' \\ \vdots \\ (a^m)' \end{bmatrix}_n .$$

The matrix

$$b = \begin{bmatrix} b_{11} & b_{12} & \cdots & b_{1s} \\ b_{21} & b_{22} & \cdots & b_{2s} \\ \vdots & \cdots & & \vdots \\ b_{t1} & b_{t2} & \cdots & b_{ts} \end{bmatrix}$$

row supervector of b is

$$b = [b_1 \ b_2 \ \ldots \ b_s]_t = [b^1 \ b^2 \ \ldots \ b^s]_t .$$

Column supervector of b is

$$b = \begin{bmatrix} b_1^1 \\ b_2^1 \\ \vdots \\ b_t^1 \end{bmatrix}_s .$$

Transpose of b;

$$b' = \begin{bmatrix} b_1^1 \\ b_2^1 \\ \vdots \\ b_s^1 \end{bmatrix}_t$$

$b' = [b_1 \; b_2 \; \ldots \; b_t]_s.$

The product of two matrices as a minor product of type II supervector.

$$ab = [a^1 \; a^2 \; \ldots \; a^m]_n \begin{bmatrix} b_1^1 \\ b_2^1 \\ \vdots \\ b_t^1 \end{bmatrix}_s$$

$$= \left[ a_1 b_1^1 + a_2 b_2^1 + \ldots + a_m b_t^1 \right]_{ns}.$$

How ever to make this point clear we give an example.

*Example 1.1.21:* Let

$$a = \begin{matrix} a_1^1 & a_2^1 \\ \begin{bmatrix} 2 & 1 \\ 3 & 5 \\ 6 & 1 \end{bmatrix} & \end{matrix}$$

and

$$b = \begin{bmatrix} 1 & 2 \\ 3 & 1 \end{bmatrix} \begin{matrix} b^1 \\ b^2 \end{matrix}.$$

$$ab = \begin{bmatrix} 2 \\ 3 \\ 6 \end{bmatrix} \begin{bmatrix} 1 & 2 \end{bmatrix} + \begin{bmatrix} 1 \\ 5 \\ 1 \end{bmatrix} \begin{bmatrix} 3 & 1 \end{bmatrix}$$

$$= \begin{bmatrix} 2 & 4 \\ 3 & 6 \\ 6 & 12 \end{bmatrix} + \begin{bmatrix} 3 & 1 \\ 15 & 5 \\ 3 & 1 \end{bmatrix}$$

$$= \begin{bmatrix} 5 & 5 \\ 18 & 11 \\ 9 & 13 \end{bmatrix}.$$

It is easily verified that if the major product of the type II supervector is computed between a and b, then the major product coincides with the minor product. From the above example.

$$ab = \begin{bmatrix} [2\ 1]\begin{bmatrix}1\\3\end{bmatrix} & [2\ 1]\begin{bmatrix}2\\1\end{bmatrix} \\ [3\ 5]\begin{bmatrix}1\\3\end{bmatrix} & [3\ 5]\begin{bmatrix}2\\1\end{bmatrix} \\ [6\ 1]\begin{bmatrix}1\\3\end{bmatrix} & [6\ 1]\begin{bmatrix}2\\1\end{bmatrix} \end{bmatrix}$$

$$= \begin{bmatrix} 2\times1+1\times3 & 2\times2+1\times1 \\ 3\times1+5\times3 & 3\times2+5\times1 \\ 6\times1+1\times3 & 6\times2+1\times1 \end{bmatrix}$$

$$= \begin{bmatrix} 5 & 5 \\ 18 & 11 \\ 9 & 13 \end{bmatrix}.$$

We can find the minor and major product of supervectors by reversing the order of the factors. Since the theory of multiplication of supermatrices involves lots of notations we have resolved to explain these concepts by working out these concepts with numerical illustrations, which we feel is easy for the grasp of the reader. Now we give the numerical illustration of the minor product of Type III vectors.

*Example 1.1.22:* Let

$$X = \begin{bmatrix} 2 & 3 & | & 4 & | & 2 & 2 & 2 \\ -1 & 1 & | & 1 & | & 1 & 0 & 1 \\ 0 & 0 & | & 2 & | & -4 & 0 & 0 \end{bmatrix}$$

and

$$Y = \begin{bmatrix} 2 & 0 \\ 1 & 1 \\ \hline 2 & 1 \\ \hline 5 & 3 \\ 1 & -1 \\ 0 & 2 \end{bmatrix}$$

be two type III supervectors. To find the product XY.

$$XY = \begin{bmatrix} 2 & 3 & | & 4 & | & 2 & 2 & 2 \\ -1 & 1 & | & 1 & | & 1 & 0 & 1 \\ 0 & 0 & | & 2 & | & -4 & 0 & 0 \end{bmatrix} \begin{bmatrix} 2 & 0 \\ 1 & 1 \\ \hline 2 & 1 \\ \hline 5 & 3 \\ 1 & -1 \\ 0 & 2 \end{bmatrix}$$

$$= \begin{bmatrix} 2 & 3 \\ -1 & 1 \\ 0 & 0 \end{bmatrix} \begin{bmatrix} 2 & 0 \\ 1 & 1 \end{bmatrix} + \begin{bmatrix} 4 \\ 1 \\ 2 \end{bmatrix} \begin{bmatrix} 2 & 1 \end{bmatrix} + \begin{bmatrix} 2 & 2 & 2 \\ 1 & 0 & 1 \\ -4 & 0 & 0 \end{bmatrix} \begin{bmatrix} 5 & 3 \\ 1 & -1 \\ 0 & 2 \end{bmatrix}$$

$$= \begin{bmatrix} 7 & 3 \\ -1 & 1 \\ 0 & 0 \end{bmatrix} + \begin{bmatrix} 8 & 4 \\ 2 & 1 \\ 4 & 2 \end{bmatrix} + \begin{bmatrix} 12 & 8 \\ 5 & 5 \\ -20 & -12 \end{bmatrix}$$

$$= \begin{bmatrix} 27 & 15 \\ 6 & 7 \\ -16 & -10 \end{bmatrix}.$$

$$Y^t X^t = \begin{bmatrix} 2 & 1 & | & 2 & | & 5 & 1 & 0 \\ 0 & 1 & | & 1 & | & 3 & -1 & 2 \end{bmatrix} \begin{bmatrix} 2 & -1 & 0 \\ 3 & 1 & 0 \\ \hline 4 & 1 & 2 \\ \hline 2 & 1 & -4 \\ 2 & 0 & 0 \\ 2 & 1 & 0 \end{bmatrix}$$

$$= \begin{bmatrix} 2 & 1 \\ 0 & 1 \end{bmatrix}\begin{bmatrix} 2 & -1 & 0 \\ 3 & 1 & 0 \end{bmatrix} + \begin{bmatrix} 2 \\ 1 \end{bmatrix}\begin{bmatrix} 4 & 1 & 2 \end{bmatrix} + \begin{bmatrix} 5 & 1 & 0 \\ 3 & -1 & 2 \end{bmatrix}\begin{bmatrix} 2 & 1 & -4 \\ 2 & 0 & 0 \\ 2 & 1 & 0 \end{bmatrix}$$

$$= \begin{bmatrix} 7 & -1 & 0 \\ 3 & 1 & 0 \end{bmatrix} + \begin{bmatrix} 8 & 2 & 4 \\ 4 & 1 & 2 \end{bmatrix} + \begin{bmatrix} 12 & 5 & -20 \\ 8 & 5 & -12 \end{bmatrix}$$

$$= \begin{bmatrix} 27 & 6 & -16 \\ 15 & 7 & -10 \end{bmatrix}.$$

From this example it is very clear.

$$(XY)^t = Y^t X^t.$$

Now we illustrate the minor product moment of type III row supervector by an example.

*Example 1.1.23:* Let

$$X = \begin{bmatrix} 2 & 3 & | & 4 & | & 3 & 4 & 5 & 0 \\ 1 & 4 & | & 1 & | & 1 & 1 & -1 & 6 \\ 2 & 1 & | & 2 & | & 0 & 2 & 1 & 1 \end{bmatrix}.$$

Consider

$$XX' = \begin{bmatrix} 2 & 3 & | & 4 & | & 3 & 4 & 5 & 0 \\ 1 & 4 & | & 1 & | & 1 & 1 & -1 & 6 \\ 2 & 1 & | & 2 & | & 0 & 2 & 1 & 1 \end{bmatrix} \begin{bmatrix} 2 & 1 & 2 \\ 3 & 4 & 1 \\ \hline 4 & 1 & 2 \\ \hline 3 & 1 & 0 \\ 4 & 1 & 2 \\ 5 & -1 & 1 \\ 0 & 6 & 1 \end{bmatrix}$$

$$= \begin{bmatrix} 2 & 3 \\ 1 & 4 \\ 2 & 1 \end{bmatrix} \begin{bmatrix} 2 & 1 & 2 \\ 3 & 4 & 1 \end{bmatrix} + \begin{bmatrix} 4 \\ 1 \\ 2 \end{bmatrix} \begin{bmatrix} 4 & 1 & 2 \end{bmatrix} +$$

$$\begin{bmatrix} 3 & 4 & 5 & 0 \\ 1 & 1 & -1 & 6 \\ 0 & 2 & 1 & 1 \end{bmatrix} \begin{bmatrix} 3 & 1 & 0 \\ 4 & 1 & 2 \\ 5 & -1 & 1 \\ 0 & 6 & 1 \end{bmatrix}$$

$$= \begin{bmatrix} 13 & 14 & 7 \\ 14 & 17 & 6 \\ 7 & 6 & 5 \end{bmatrix} + \begin{bmatrix} 16 & 4 & 8 \\ 4 & 1 & 2 \\ 8 & 2 & 4 \end{bmatrix} + \begin{bmatrix} 50 & 2 & 13 \\ 2 & 39 & 7 \\ 13 & 7 & 6 \end{bmatrix}$$

$$= \begin{bmatrix} 79 & 20 & 28 \\ 20 & 57 & 15 \\ 28 & 15 & 15 \end{bmatrix}.$$

Minor product of Type III column supervector is illustrated by the following example.

*Example 1.1.24:* Let

$$Y^t = \begin{bmatrix} 2 & 3 & 1 & | & 0 & | & 1 & 2 & 1 & 5 & 1 \\ 0 & 1 & 5 & | & 2 & | & 0 & 3 & 0 & 1 & 0 \end{bmatrix}$$

where Y is the column supervector

$$Y^tY = \begin{bmatrix} 2 & 3 & 1 & 0 & 1 & 2 & 1 & 5 & 1 \\ 0 & 1 & 5 & 2 & 0 & 3 & 0 & 1 & 0 \end{bmatrix} \begin{bmatrix} 2 & 0 \\ 3 & 1 \\ 1 & 5 \\ 0 & 2 \\ 1 & 0 \\ 2 & 3 \\ 1 & 0 \\ 5 & 1 \\ 1 & 0 \end{bmatrix}$$

$$= \begin{bmatrix} 2 & 3 & 1 \\ 0 & 1 & 5 \end{bmatrix} \begin{bmatrix} 2 & 0 \\ 3 & 1 \\ 1 & 5 \end{bmatrix} + \begin{bmatrix} 0 \\ 2 \end{bmatrix} \begin{bmatrix} 0 & 2 \end{bmatrix} + \begin{bmatrix} 1 & 2 & 1 & 5 & 1 \\ 0 & 3 & 0 & 1 & 0 \end{bmatrix} \begin{bmatrix} 1 & 0 \\ 2 & 3 \\ 1 & 0 \\ 5 & 1 \\ 1 & 0 \end{bmatrix}$$

$$= \begin{bmatrix} 14 & 8 \\ 8 & 26 \end{bmatrix} + \begin{bmatrix} 0 & 0 \\ 0 & 4 \end{bmatrix} + \begin{bmatrix} 32 & 11 \\ 11 & 10 \end{bmatrix} = \begin{bmatrix} 46 & 19 \\ 19 & 40 \end{bmatrix}.$$

Next we proceed on to illustrate the major product of Type III vectors.

*Example 1.1.25:* Let

$$X = \begin{bmatrix} 3 & 1 & 6 \\ 2 & 0 & -1 \\ 1 & 2 & 3 \\ 6 & 3 & 0 \\ 4 & 2 & 1 \\ 5 & 1 & -1 \end{bmatrix}$$

and

$$Y = \begin{bmatrix} 3 & | & 5 & 2 & 0 \\ 1 & | & 1 & 2 & 2 \\ 0 & | & 3 & 1 & -2 \end{bmatrix}.$$

$$XY = \begin{bmatrix} 3 & 1 & 6 \\ 2 & 0 & -1 \\ \hline 1 & 2 & 3 \\ \hline 6 & 3 & 0 \\ 4 & 2 & 1 \\ 5 & 1 & -1 \end{bmatrix} \begin{bmatrix} 3 & | & 5 & 2 & 0 \\ 1 & | & 1 & 2 & 2 \\ 0 & | & 3 & 1 & -2 \end{bmatrix}$$

$$= \begin{bmatrix} \begin{bmatrix} 3 & 1 & 6 \\ 2 & 0 & -1 \end{bmatrix} \begin{bmatrix} 3 \\ 1 \\ 0 \end{bmatrix} & \begin{bmatrix} 3 & 1 & 6 \\ 2 & 0 & -1 \end{bmatrix} \begin{bmatrix} 5 & 2 & 0 \\ 1 & 2 & 2 \\ 3 & 1 & -2 \end{bmatrix} \\ \hline [1\ 2\ 3] \begin{bmatrix} 3 \\ 1 \\ 0 \end{bmatrix} & [1\ 2\ 3] \begin{bmatrix} 5 & 2 & 0 \\ 1 & 2 & 2 \\ 3 & 1 & -2 \end{bmatrix} \\ \hline \begin{bmatrix} 6 & 3 & 0 \\ 4 & 2 & 1 \\ 5 & 1 & -1 \end{bmatrix} \begin{bmatrix} 3 \\ 1 \\ 0 \end{bmatrix} & \begin{bmatrix} 6 & 3 & 0 \\ 4 & 2 & 1 \\ 5 & 1 & -1 \end{bmatrix} \begin{bmatrix} 5 & 2 & 0 \\ 1 & 2 & 2 \\ 3 & 1 & -2 \end{bmatrix} \end{bmatrix}$$

$$= \begin{bmatrix} 10 & | & 34 & 14 & -10 \\ 6 & | & 7 & 3 & 2 \\ \hline 5 & | & 16 & 9 & -2 \\ \hline 21 & | & 33 & 18 & 6 \\ 14 & | & 25 & 13 & 2 \\ 16 & | & 23 & 11 & 4 \end{bmatrix}.$$

Now minor product of type IV vector is illustrated by the following example.

***Example 1.1.26:*** Let

$$X = \begin{bmatrix} 1 & 3 & 1 & 2 & 5 & 1 \\ 2 & 1 & 1 & 1 & 2 & 0 \\ \hline 1 & 5 & 1 & 1 & 1 & 2 \\ 4 & 1 & 0 & 2 & 2 & 1 \\ 3 & 2 & 1 & 0 & 1 & 1 \\ -1 & 0 & 1 & 1 & 0 & 1 \\ \hline 4 & 2 & 1 & 3 & 1 & 1 \end{bmatrix}$$

and

$$Y = \begin{bmatrix} 1 & 1 & 0 & 1 & 3 & 1 & 2 & 1 & 2 \\ 2 & 0 & 1 & 0 & 1 & 2 & 0 & 0 & 1 \\ 1 & 1 & 0 & 2 & 3 & 0 & 1 & 1 & 4 \\ 1 & 0 & 1 & 1 & 2 & 1 & 1 & 2 & 0 \\ \hline 1 & 2 & 0 & 1 & 1 & 0 & 1 & 1 & 2 \\ 0 & 1 & 1 & 0 & 1 & 1 & 0 & 2 & 1 \end{bmatrix}.$$

$$XY = \begin{bmatrix} 1 & 3 & 1 & 2 & 5 & 1 \\ 2 & 1 & 1 & 1 & 2 & 0 \\ \hline 1 & 5 & 1 & 1 & 1 & 2 \\ 4 & 1 & 0 & 2 & 2 & 1 \\ 3 & 2 & 1 & 0 & 1 & 1 \\ -1 & 0 & 1 & 1 & 0 & 1 \\ \hline 4 & 2 & 1 & 3 & 1 & 1 \end{bmatrix} \times$$

$$\begin{bmatrix} 1 & 1 & 0 & 1 & 3 & 1 & 2 & 1 & 2 \\ 2 & 0 & 1 & 0 & 1 & 2 & 0 & 0 & 1 \\ 1 & 1 & 0 & 2 & 3 & 0 & 1 & 1 & 4 \\ 1 & 0 & 1 & 1 & 2 & 1 & 1 & 2 & 0 \\ \hline 1 & 2 & 0 & 1 & 1 & 0 & 1 & 1 & 2 \\ 0 & 1 & 1 & 0 & 1 & 1 & 0 & 2 & 1 \end{bmatrix}$$

$$= \begin{bmatrix} 1 \\ 2 \\ \hline 1 \\ 4 \\ 3 \\ -1 \\ \hline 4 \end{bmatrix} [1\ 1\ 0\ 1\ |\ 3\ 1\ 2\ |\ 1\ 2] +$$

$$\begin{bmatrix} 3 & 1 & 2 \\ 1 & 1 & 1 \\ \hline 5 & 1 & 1 \\ 1 & 0 & 2 \\ 2 & 1 & 0 \\ 0 & 1 & 1 \\ \hline 2 & 1 & 3 \end{bmatrix} \begin{bmatrix} 2 & 0 & 1 & 0 & | & 1 & 2 & 0 & | & 0 & 1 \\ 1 & 1 & 0 & 2 & | & 3 & 0 & 1 & | & 1 & 4 \\ 1 & 0 & 1 & 1 & | & 2 & 1 & 1 & | & 2 & 0 \end{bmatrix}$$

$$+ \begin{bmatrix} 5 & 1 \\ 2 & 0 \\ \hline 1 & 2 \\ 2 & 1 \\ 1 & 1 \\ 0 & 1 \\ \hline 1 & 1 \end{bmatrix} \begin{bmatrix} 1 & 2 & 0 & 1 & | & 1 & 0 & 1 & | & 1 & 2 \\ 0 & 1 & 1 & 0 & | & 1 & 1 & 0 & | & 2 & 1 \end{bmatrix}$$

$$= \begin{bmatrix} \begin{bmatrix} 1 \\ 2 \end{bmatrix}[1\ 1\ 0\ 1] & \begin{bmatrix} 1 \\ 2 \end{bmatrix}[3\ 1\ 2] & \begin{bmatrix} 1 \\ 2 \end{bmatrix}[1\ 2] \\ \hline \begin{bmatrix} 1 \\ 4 \\ 3 \\ -1 \end{bmatrix}[1\ 1\ 0\ 1] & \begin{bmatrix} 1 \\ 4 \\ 3 \\ -1 \end{bmatrix}[3\ 1\ 2] & \begin{bmatrix} 1 \\ 4 \\ 3 \\ -1 \end{bmatrix}[1\ 2] \\ \hline [4][1\ 1\ 0\ 1] & [4][3\ 1\ 2] & [4][1\ 2] \end{bmatrix} +$$

$$\begin{bmatrix} \begin{bmatrix} 3 & 1 & 2 \\ 1 & 1 & 1 \end{bmatrix} \begin{bmatrix} 2 & 0 & 1 & 0 \\ 1 & 1 & 0 & 2 \\ 1 & 0 & 1 & 1 \end{bmatrix} & \begin{bmatrix} 3 & 1 & 2 \\ 1 & 1 & 1 \end{bmatrix} \begin{bmatrix} 1 & 2 & 0 \\ 3 & 0 & 1 \\ 2 & 1 & 1 \end{bmatrix} & \begin{bmatrix} 3 & 1 & 2 \\ 1 & 1 & 1 \end{bmatrix} \begin{bmatrix} 0 & 1 \\ 1 & 4 \\ 2 & 0 \end{bmatrix} \\ \begin{bmatrix} 5 & 1 & 1 \\ 1 & 0 & 2 \\ 2 & 1 & 0 \\ 0 & 1 & 1 \end{bmatrix} \begin{bmatrix} 2 & 0 & 1 & 0 \\ 1 & 1 & 0 & 2 \\ 1 & 0 & 1 & 1 \end{bmatrix} & \begin{bmatrix} 5 & 1 & 1 \\ 1 & 0 & 2 \\ 2 & 1 & 0 \\ 0 & 1 & 1 \end{bmatrix} \begin{bmatrix} 1 & 2 & 0 \\ 3 & 0 & 1 \\ 2 & 1 & 1 \end{bmatrix} & \begin{bmatrix} 5 & 1 & 1 \\ 1 & 0 & 2 \\ 2 & 1 & 0 \\ 0 & 1 & 1 \end{bmatrix} \begin{bmatrix} 0 & 1 \\ 1 & 4 \\ 2 & 0 \end{bmatrix} \\ \begin{bmatrix} 2 & 1 & 3 \end{bmatrix} \begin{bmatrix} 2 & 0 & 1 & 0 \\ 1 & 1 & 0 & 2 \\ 1 & 0 & 1 & 1 \end{bmatrix} & \begin{bmatrix} 2 & 1 & 3 \end{bmatrix} \begin{bmatrix} 1 & 2 & 0 \\ 3 & 0 & 1 \\ 2 & 1 & 1 \end{bmatrix} & \begin{bmatrix} 2 & 1 & 3 \end{bmatrix} \begin{bmatrix} 0 & 1 \\ 1 & 4 \\ 2 & 0 \end{bmatrix} \end{bmatrix}$$

$$+ \begin{bmatrix} \begin{bmatrix} 5 & 1 \\ 2 & 0 \end{bmatrix} \begin{bmatrix} 1 & 2 & 0 & 1 \\ 0 & 1 & 1 & 0 \end{bmatrix} & \begin{bmatrix} 5 & 1 \\ 2 & 0 \end{bmatrix} \begin{bmatrix} 1 & 0 & 1 \\ 1 & 1 & 0 \end{bmatrix} & \begin{bmatrix} 5 & 1 \\ 2 & 0 \end{bmatrix} \begin{bmatrix} 1 & 2 \\ 2 & 1 \end{bmatrix} \\ \begin{bmatrix} 1 & 2 \\ 2 & 1 \\ 1 & 1 \\ 0 & 1 \end{bmatrix} \begin{bmatrix} 1 & 2 & 0 & 1 \\ 0 & 1 & 1 & 0 \end{bmatrix} & \begin{bmatrix} 1 & 2 \\ 2 & 1 \\ 1 & 1 \\ 0 & 1 \end{bmatrix} \begin{bmatrix} 1 & 0 & 1 \\ 1 & 1 & 0 \end{bmatrix} & \begin{bmatrix} 1 & 2 \\ 2 & 1 \\ 1 & 1 \\ 0 & 1 \end{bmatrix} \begin{bmatrix} 1 & 2 \\ 2 & 1 \end{bmatrix} \\ \begin{bmatrix} 1 & 1 \end{bmatrix} \begin{bmatrix} 1 & 2 & 0 & 1 \\ 0 & 1 & 1 & 0 \end{bmatrix} & \begin{bmatrix} 1 & 1 \end{bmatrix} \begin{bmatrix} 1 & 0 & 1 \\ 1 & 1 & 0 \end{bmatrix} & \begin{bmatrix} 1 & 1 \end{bmatrix} \begin{bmatrix} 1 & 2 \\ 2 & 1 \end{bmatrix} \end{bmatrix}$$

$$= \begin{bmatrix} \begin{bmatrix} 1 & 1 & 0 & 1 \\ 2 & 2 & 0 & 2 \end{bmatrix} & \begin{bmatrix} 3 & 1 & 2 \\ 6 & 2 & 4 \end{bmatrix} & \begin{bmatrix} 1 & 2 \\ 2 & 4 \end{bmatrix} \\ \begin{bmatrix} 1 & 1 & 0 & 1 \\ 4 & 4 & 0 & 4 \\ 3 & 3 & 0 & 3 \\ -1 & -1 & 0 & -1 \end{bmatrix} & \begin{bmatrix} 3 & 1 & 2 \\ 12 & 4 & 8 \\ 9 & 3 & 6 \\ -3 & -1 & -2 \end{bmatrix} & \begin{bmatrix} 1 & 2 \\ 4 & 8 \\ 3 & 6 \\ -1 & -2 \end{bmatrix} \\ \begin{bmatrix} 4 & 4 & 0 & 4 \end{bmatrix} & \begin{bmatrix} 12 & 4 & 8 \end{bmatrix} & \begin{bmatrix} 4 & 8 \end{bmatrix} \end{bmatrix} +$$

$$\begin{bmatrix} 9 & 1 & 5 & 4 & 10 & 8 & 3 & 5 & 7 \\ 4 & 1 & 2 & 3 & 6 & 3 & 2 & 3 & 5 \\ \hline 12 & 1 & 6 & 3 & 10 & 11 & 2 & 3 & 9 \\ 4 & 0 & 3 & 2 & 5 & 4 & 2 & 4 & 1 \\ 5 & 1 & 2 & 2 & 5 & 4 & 1 & 1 & 6 \\ 2 & 1 & 1 & 3 & 5 & 1 & 2 & 3 & 4 \\ \hline 8 & 1 & 5 & 5 & 11 & 7 & 4 & 7 & 6 \end{bmatrix}$$

$$+ \begin{bmatrix} 5 & 11 & 1 & 5 & 6 & 1 & 5 & 7 & 11 \\ 2 & 4 & 0 & 2 & 2 & 0 & 2 & 2 & 4 \\ \hline 1 & 4 & 2 & 1 & 3 & 2 & 1 & 5 & 4 \\ 2 & 5 & 1 & 2 & 3 & 1 & 2 & 4 & 5 \\ 1 & 3 & 1 & 1 & 2 & 1 & 1 & 3 & 3 \\ 0 & 1 & 1 & 0 & 1 & 1 & 0 & 2 & 1 \\ \hline 1 & 3 & 1 & 1 & 2 & 1 & 1 & 3 & 3 \end{bmatrix}$$

$$= \begin{bmatrix} 15 & 13 & 6 & 10 & 19 & 10 & 10 & 13 & 20 \\ 8 & 7 & 2 & 7 & 14 & 5 & 8 & 7 & 13 \\ \hline 14 & 6 & 8 & 5 & 16 & 14 & 5 & 9 & 15 \\ 10 & 9 & 4 & 8 & 20 & 9 & 12 & 12 & 14 \\ 9 & 7 & 3 & 6 & 16 & 8 & 8 & 7 & 15 \\ 1 & 1 & 2 & 2 & 3 & 1 & 0 & 4 & 3 \\ \hline 13 & 8 & 6 & 10 & 25 & 12 & 13 & 14 & 17 \end{bmatrix}.$$

We now illustrate minor product moment of type IV row vector

*Example 1.1.27:* Let

$$X = \begin{bmatrix} 1 & 1 & 1 & 1 & 0 & 1 \\ 2 & 1 & 2 & 2 & 1 & 2 \\ \hline 0 & 1 & 1 & 3 & 1 & 1 \\ 1 & 0 & 1 & 1 & 3 & 2 \\ 5 & 1 & 0 & 2 & 1 & 3 \\ 1 & 1 & 0 & 1 & 2 & 4 \\ \hline 2 & 1 & 1 & 5 & 0 & 2 \end{bmatrix}.$$

$$XX^t = \begin{bmatrix} 1 & 1 & 1 & 1 & 0 & 1 \\ 2 & 1 & 2 & 2 & 1 & 2 \\ \hline 0 & 1 & 1 & 3 & 1 & 1 \\ 1 & 0 & 1 & 1 & 3 & 2 \\ 5 & 1 & 0 & 2 & 1 & 3 \\ 1 & 1 & 0 & 1 & 2 & 4 \\ \hline 2 & 1 & 1 & 5 & 0 & 2 \end{bmatrix} \begin{bmatrix} 1 & 2 & 0 & 1 & 5 & 1 & 2 \\ 1 & 1 & 1 & 0 & 1 & 1 & 1 \\ 1 & 2 & 1 & 1 & 0 & 0 & 1 \\ \hline 1 & 2 & 3 & 1 & 2 & 1 & 5 \\ 0 & 1 & 1 & 3 & 1 & 2 & 0 \\ \hline 1 & 2 & 1 & 2 & 3 & 4 & 2 \end{bmatrix}$$

$$= \begin{bmatrix} 1 & 1 & 1 \\ 2 & 1 & 2 \\ \hline 0 & 1 & 1 \\ 1 & 0 & 1 \\ 5 & 1 & 0 \\ 1 & 1 & 0 \\ \hline 2 & 1 & 1 \end{bmatrix} \begin{bmatrix} 1 & 2 & 0 & 1 & 5 & 1 & 2 \\ 1 & 1 & 1 & 0 & 1 & 1 & 1 \\ 1 & 2 & 1 & 1 & 0 & 0 & 1 \end{bmatrix} +$$

$$\begin{bmatrix} 1 & 0 \\ 2 & 1 \\ \hline 3 & 1 \\ 1 & 3 \\ 2 & 1 \\ 1 & 2 \\ \hline 5 & 0 \end{bmatrix} \begin{bmatrix} 1 & 2 & 3 & 1 & 2 & 1 & 5 \\ 0 & 1 & 1 & 3 & 1 & 2 & 0 \end{bmatrix} + \begin{bmatrix} 1 \\ 2 \\ \hline 1 \\ 2 \\ 3 \\ 4 \\ \hline 2 \end{bmatrix} [1\ 2\ |\ 1\ 2\ 3\ 4\ |\ 2] =$$

$$\begin{bmatrix} \begin{bmatrix} 1 & 1 & 1 \\ 2 & 1 & 2 \end{bmatrix}\begin{bmatrix} 1 & 2 \\ 1 & 1 \\ 1 & 2 \end{bmatrix} & \begin{bmatrix} 1 & 1 & 1 \\ 2 & 1 & 2 \end{bmatrix}\begin{bmatrix} 0 & 1 & 5 & 1 \\ 1 & 0 & 1 & 1 \\ 1 & 1 & 0 & 0 \end{bmatrix} & \begin{bmatrix} 1 & 1 & 1 \\ 2 & 1 & 2 \end{bmatrix}\begin{bmatrix} 2 \\ 1 \\ 1 \end{bmatrix} \\ \begin{bmatrix} 0 & 1 & 1 \\ 1 & 0 & 1 \\ 5 & 1 & 0 \\ 1 & 1 & 0 \end{bmatrix}\begin{bmatrix} 1 & 2 \\ 1 & 1 \\ 1 & 2 \end{bmatrix} & \begin{bmatrix} 0 & 1 & 1 \\ 1 & 0 & 1 \\ 5 & 1 & 0 \\ 1 & 1 & 0 \end{bmatrix}\begin{bmatrix} 0 & 1 & 5 & 1 \\ 1 & 0 & 1 & 1 \\ 1 & 1 & 0 & 0 \end{bmatrix} & \begin{bmatrix} 0 & 1 & 1 \\ 1 & 0 & 1 \\ 5 & 1 & 0 \\ 1 & 1 & 0 \end{bmatrix}\begin{bmatrix} 2 \\ 1 \\ 1 \end{bmatrix} \\ \begin{bmatrix} 2 & 1 & 1 \end{bmatrix}\begin{bmatrix} 1 & 2 \\ 1 & 1 \\ 1 & 2 \end{bmatrix} & \begin{bmatrix} 2 & 1 & 1 \end{bmatrix}\begin{bmatrix} 0 & 1 & 5 & 1 \\ 1 & 0 & 1 & 1 \\ 1 & 1 & 0 & 0 \end{bmatrix} & \begin{bmatrix} 2 & 1 & 1 \end{bmatrix}\begin{bmatrix} 2 \\ 1 \\ 1 \end{bmatrix} \end{bmatrix}$$

$$+ \begin{bmatrix} \begin{bmatrix} 1 & 0 \\ 2 & 1 \end{bmatrix}\begin{bmatrix} 1 & 2 \\ 0 & 1 \end{bmatrix} & \begin{bmatrix} 1 & 0 \\ 2 & 1 \end{bmatrix}\begin{bmatrix} 3 & 1 & 2 & 1 \\ 1 & 3 & 1 & 2 \end{bmatrix} & \begin{bmatrix} 1 & 0 \\ 2 & 1 \end{bmatrix}\begin{bmatrix} 5 \\ 0 \end{bmatrix} \\ \begin{bmatrix} 3 & 1 \\ 1 & 3 \\ 2 & 1 \\ 1 & 2 \end{bmatrix}\begin{bmatrix} 1 & 2 \\ 0 & 1 \end{bmatrix} & \begin{bmatrix} 3 & 1 \\ 1 & 3 \\ 2 & 1 \\ 1 & 2 \end{bmatrix}\begin{bmatrix} 3 & 1 & 2 & 1 \\ 1 & 3 & 1 & 2 \end{bmatrix} & \begin{bmatrix} 3 & 1 \\ 1 & 3 \\ 2 & 1 \\ 1 & 2 \end{bmatrix}\begin{bmatrix} 5 \\ 0 \end{bmatrix} \\ \begin{bmatrix} 5 & 0 \end{bmatrix}\begin{bmatrix} 1 & 2 \\ 0 & 1 \end{bmatrix} & \begin{bmatrix} 5 & 0 \end{bmatrix}\begin{bmatrix} 3 & 1 & 2 & 1 \\ 1 & 3 & 1 & 2 \end{bmatrix} & \begin{bmatrix} 5 & 0 \end{bmatrix}\begin{bmatrix} 5 \\ 0 \end{bmatrix} \end{bmatrix}$$

$$+ \begin{bmatrix} \begin{bmatrix} 1 \\ 2 \end{bmatrix}\begin{bmatrix} 1 & 2 \end{bmatrix} & \begin{bmatrix} 1 \\ 2 \end{bmatrix}\begin{bmatrix} 1 & 2 & 3 & 4 \end{bmatrix} & \begin{bmatrix} 1 \\ 2 \end{bmatrix}\begin{bmatrix} 2 \end{bmatrix} \\ \begin{bmatrix} 1 \\ 2 \\ 3 \\ 4 \end{bmatrix}\begin{bmatrix} 1 & 2 \end{bmatrix} & \begin{bmatrix} 1 \\ 2 \\ 3 \\ 4 \end{bmatrix}\begin{bmatrix} 1 & 2 & 3 & 4 \end{bmatrix} & \begin{bmatrix} 1 \\ 2 \\ 3 \\ 4 \end{bmatrix}\begin{bmatrix} 2 \end{bmatrix} \\ \begin{bmatrix} 2 \end{bmatrix}\begin{bmatrix} 1 & 2 \end{bmatrix} & \begin{bmatrix} 2 \end{bmatrix}\begin{bmatrix} 1 & 2 & 3 & 4 \end{bmatrix} & \begin{bmatrix} 2 \end{bmatrix}\begin{bmatrix} 2 \end{bmatrix} \end{bmatrix}$$

$$= \begin{bmatrix} 3 & 5 & 2 & 2 & 6 & 2 & 4 \\ 5 & 9 & 3 & 4 & 11 & 3 & 7 \\ \hline 2 & 3 & 2 & 1 & 1 & 1 & 2 \\ 2 & 4 & 1 & 2 & 5 & 1 & 3 \\ 6 & 11 & 1 & 5 & 26 & 6 & 11 \\ 2 & 3 & 1 & 1 & 6 & 2 & 3 \\ \hline 4 & 7 & 2 & 3 & 11 & 3 & 6 \end{bmatrix} +$$

$$\begin{bmatrix} 1 & 2 & 3 & 1 & 2 & 1 & 5 \\ 2 & 5 & 7 & 5 & 5 & 4 & 10 \\ \hline 3 & 7 & 10 & 6 & 7 & 5 & 15 \\ 1 & 5 & 6 & 10 & 5 & 7 & 5 \\ 2 & 5 & 7 & 5 & 5 & 4 & 10 \\ 1 & 4 & 5 & 7 & 4 & 5 & 5 \\ \hline 5 & 10 & 15 & 5 & 10 & 5 & 25 \end{bmatrix} + \begin{bmatrix} 1 & 2 & 1 & 2 & 3 & 4 & 2 \\ 2 & 4 & 2 & 4 & 6 & 8 & 4 \\ \hline 1 & 2 & 1 & 2 & 3 & 4 & 2 \\ 2 & 4 & 2 & 4 & 6 & 8 & 4 \\ 3 & 6 & 3 & 6 & 9 & 12 & 6 \\ 4 & 8 & 4 & 8 & 12 & 16 & 8 \\ \hline 2 & 4 & 2 & 4 & 6 & 8 & 4 \end{bmatrix}$$

$$= \begin{bmatrix} 5 & 9 & 6 & 5 & 11 & 7 & 11 \\ 9 & 18 & 12 & 13 & 22 & 15 & 21 \\ \hline 6 & 12 & 13 & 9 & 11 & 10 & 19 \\ 5 & 13 & 9 & 16 & 16 & 16 & 12 \\ 11 & 22 & 11 & 16 & 40 & 22 & 27 \\ 7 & 15 & 10 & 16 & 22 & 23 & 16 \\ \hline 11 & 21 & 19 & 12 & 27 & 16 & 35 \end{bmatrix}.$$

The minor product moment of type IV column vector is illustrated for the same X just given in case of row product.

*Example 1.1.28:* Let

$$X = \begin{bmatrix} 1 & 1 & 1 & 1 & 0 & 1 \\ 2 & 1 & 2 & 2 & 1 & 2 \\ \hline 0 & 1 & 1 & 3 & 1 & 1 \\ 1 & 0 & 1 & 1 & 3 & 2 \\ 5 & 1 & 0 & 2 & 1 & 3 \\ 1 & 1 & 0 & 1 & 2 & 4 \\ \hline 2 & 1 & 1 & 5 & 0 & 2 \end{bmatrix}.$$

$$X^t X = \begin{bmatrix} 1 & 2 & 0 & 1 & 5 & 1 & 2 \\ 1 & 1 & 1 & 0 & 1 & 1 & 1 \\ 1 & 2 & 1 & 1 & 0 & 0 & 1 \\ \hline 1 & 2 & 3 & 1 & 2 & 1 & 5 \\ 0 & 1 & 1 & 3 & 1 & 2 & 0 \\ \hline 1 & 2 & 1 & 2 & 3 & 4 & 2 \end{bmatrix} \begin{bmatrix} 1 & 1 & 1 & 1 & 0 & 1 \\ 2 & 1 & 2 & 2 & 1 & 2 \\ \hline 0 & 1 & 1 & 3 & 1 & 1 \\ 1 & 0 & 1 & 1 & 3 & 2 \\ 5 & 1 & 0 & 2 & 1 & 3 \\ 1 & 1 & 0 & 1 & 2 & 4 \\ \hline 2 & 1 & 1 & 5 & 0 & 2 \end{bmatrix}$$

$$= \begin{bmatrix} 1 & 2 \\ 1 & 1 \\ 1 & 2 \\ \hline 1 & 2 \\ 0 & 1 \\ \hline 1 & 2 \end{bmatrix} \begin{bmatrix} 1 & 1 & 1 & 1 & 0 & 1 \\ 2 & 1 & 2 & 2 & 1 & 2 \end{bmatrix}$$

$$+ \begin{bmatrix} 0 & 1 & 5 & 1 \\ 1 & 0 & 1 & 1 \\ 1 & 1 & 0 & 0 \\ \hline 3 & 1 & 2 & 1 \\ 1 & 3 & 1 & 2 \\ \hline 1 & 2 & 3 & 4 \end{bmatrix} \begin{bmatrix} 0 & 1 & 1 & 3 & 1 & 1 \\ 1 & 0 & 1 & 1 & 3 & 2 \\ 5 & 1 & 0 & 2 & 1 & 3 \\ 1 & 1 & 0 & 1 & 2 & 4 \end{bmatrix} + \begin{bmatrix} 2 \\ 1 \\ 1 \\ \hline 5 \\ 0 \\ \hline 2 \end{bmatrix} [2\ 1\ 1\ |\ 5\ 0\ |\ 2]$$

$$= \left[ \begin{array}{c|c|c} \begin{bmatrix} 1 & 2 \\ 1 & 1 \\ 1 & 2 \end{bmatrix} \begin{bmatrix} 1 & 1 & 1 \\ 2 & 1 & 2 \end{bmatrix} & \begin{bmatrix} 1 & 2 \\ 1 & 1 \\ 1 & 2 \end{bmatrix} \begin{bmatrix} 1 & 0 \\ 2 & 1 \end{bmatrix} & \begin{bmatrix} 1 & 2 \\ 1 & 1 \\ 1 & 2 \end{bmatrix} \begin{bmatrix} 1 \\ 2 \end{bmatrix} \\ \hline \begin{bmatrix} 1 & 2 \\ 0 & 1 \end{bmatrix} \begin{bmatrix} 1 & 1 & 1 \\ 2 & 1 & 2 \end{bmatrix} & \begin{bmatrix} 1 & 2 \\ 0 & 1 \end{bmatrix} \begin{bmatrix} 1 & 0 \\ 2 & 1 \end{bmatrix} & \begin{bmatrix} 1 & 2 \\ 0 & 1 \end{bmatrix} \begin{bmatrix} 1 \\ 2 \end{bmatrix} \\ \hline \begin{bmatrix} 1 & 2 \end{bmatrix} \begin{bmatrix} 1 & 1 & 1 \\ 2 & 1 & 2 \end{bmatrix} & \begin{bmatrix} 1 & 2 \end{bmatrix} \begin{bmatrix} 1 & 0 \\ 2 & 1 \end{bmatrix} & \begin{bmatrix} 1 & 2 \end{bmatrix} \begin{bmatrix} 1 \\ 2 \end{bmatrix} \end{array} \right] +$$

$$\left[ \begin{array}{c|c|c} \begin{bmatrix} 0 & 1 & 5 & 1 \\ 1 & 0 & 1 & 1 \\ 1 & 1 & 0 & 0 \end{bmatrix} \begin{bmatrix} 0 & 1 & 1 \\ 1 & 0 & 1 \\ 5 & 1 & 0 \\ 1 & 1 & 0 \end{bmatrix} & \begin{bmatrix} 0 & 1 & 5 & 1 \\ 1 & 0 & 1 & 1 \\ 1 & 1 & 0 & 0 \end{bmatrix} \begin{bmatrix} 3 & 1 \\ 1 & 3 \\ 2 & 1 \\ 1 & 2 \end{bmatrix} & \begin{bmatrix} 0 & 1 & 5 & 1 \\ 1 & 0 & 1 & 1 \\ 1 & 1 & 0 & 0 \end{bmatrix} \begin{bmatrix} 1 \\ 2 \\ 3 \\ 4 \end{bmatrix} \\ \hline \begin{bmatrix} 3 & 1 & 2 & 1 \\ 1 & 3 & 1 & 2 \end{bmatrix} \begin{bmatrix} 0 & 1 & 1 \\ 1 & 0 & 1 \\ 5 & 1 & 0 \\ 1 & 1 & 0 \end{bmatrix} & \begin{bmatrix} 3 & 1 & 2 & 1 \\ 1 & 3 & 1 & 2 \end{bmatrix} \begin{bmatrix} 3 & 1 \\ 1 & 3 \\ 2 & 1 \\ 1 & 2 \end{bmatrix} & \begin{bmatrix} 3 & 1 & 2 & 1 \\ 1 & 3 & 1 & 2 \end{bmatrix} \begin{bmatrix} 1 \\ 2 \\ 3 \\ 4 \end{bmatrix} \\ \hline \begin{bmatrix} 1 & 2 & 3 & 4 \end{bmatrix} \begin{bmatrix} 0 & 1 & 1 \\ 1 & 0 & 1 \\ 5 & 1 & 0 \\ 1 & 1 & 0 \end{bmatrix} & \begin{bmatrix} 1 & 2 & 3 & 4 \end{bmatrix} \begin{bmatrix} 3 & 1 \\ 1 & 3 \\ 2 & 1 \\ 1 & 2 \end{bmatrix} & \begin{bmatrix} 1 & 2 & 3 & 4 \end{bmatrix} \begin{bmatrix} 1 \\ 2 \\ 3 \\ 4 \end{bmatrix} \end{array} \right]$$

$$+ \left[ \begin{array}{c|c|c} \begin{bmatrix} 2 \\ 1 \\ 1 \end{bmatrix} \begin{bmatrix} 2 & 1 & 1 \end{bmatrix} & \begin{bmatrix} 2 \\ 1 \\ 1 \end{bmatrix} \begin{bmatrix} 5 & 0 \end{bmatrix} & \begin{bmatrix} 2 \\ 1 \\ 1 \end{bmatrix} \begin{bmatrix} 2 \end{bmatrix} \\ \hline \begin{bmatrix} 5 \\ 0 \end{bmatrix} \begin{bmatrix} 2 & 1 & 1 \end{bmatrix} & \begin{bmatrix} 5 \\ 0 \end{bmatrix} \begin{bmatrix} 5 & 0 \end{bmatrix} & \begin{bmatrix} 5 \\ 0 \end{bmatrix} \begin{bmatrix} 2 \end{bmatrix} \\ \hline \begin{bmatrix} 2 \end{bmatrix} \begin{bmatrix} 2 & 1 & 1 \end{bmatrix} & \begin{bmatrix} 2 \end{bmatrix} \begin{bmatrix} 5 & 0 \end{bmatrix} & \begin{bmatrix} 2 \end{bmatrix} \begin{bmatrix} 2 \end{bmatrix} \end{array} \right] =$$

$$\begin{bmatrix} 5 & 3 & 5 & 5 & 2 & 5 \\ 3 & 2 & 3 & 3 & 1 & 3 \\ 5 & 3 & 5 & 5 & 2 & 5 \\ \hline 5 & 3 & 5 & 5 & 2 & 5 \\ 2 & 1 & 2 & 2 & 1 & 2 \\ 5 & 3 & 5 & 5 & 2 & 5 \end{bmatrix} + \begin{bmatrix} 27 & 6 & 1 & 12 & 10 & 21 \\ 6 & 3 & 1 & 6 & 4 & 8 \\ 1 & 1 & 2 & 4 & 4 & 3 \\ \hline 12 & 6 & 4 & 15 & 10 & 15 \\ 10 & 4 & 4 & 10 & 15 & 18 \\ 21 & 8 & 3 & 15 & 18 & 30 \end{bmatrix} +$$

$$\begin{bmatrix} 4 & 2 & 2 & 10 & 0 & 4 \\ 2 & 1 & 1 & 5 & 0 & 2 \\ 2 & 1 & 1 & 5 & 0 & 2 \\ \hline 10 & 5 & 5 & 25 & 0 & 10 \\ 0 & 0 & 0 & 0 & 0 & 0 \\ 4 & 2 & 2 & 10 & 0 & 4 \end{bmatrix}$$

$$= \begin{bmatrix} 36 & 11 & 8 & 27 & 12 & 30 \\ 11 & 6 & 5 & 14 & 5 & 13 \\ 8 & 5 & 8 & 14 & 6 & 10 \\ \hline 27 & 14 & 14 & 45 & 12 & 30 \\ 12 & 5 & 6 & 12 & 16 & 20 \\ 30 & 13 & 10 & 30 & 20 & 39 \end{bmatrix}.$$

Now we proceed on to illustrate the major product of type IV vectors

*Example 1.1.29:* Let

$$X = \begin{bmatrix} \begin{bmatrix} 1 & 2 & 1 & 1 & 2 & 3 \\ 3 & 1 & 2 & 3 & 1 & 1 \end{bmatrix} \\ \begin{bmatrix} 1 & 1 & 3 & 1 & 1 & 1 \\ 2 & 3 & 1 & 2 & 0 & 1 \\ 3 & 4 & 2 & 0 & 1 & 0 \\ 4 & 2 & 4 & 1 & 0 & 0 \end{bmatrix} \\ \begin{bmatrix} 5 & 0 & 1 & 1 & 1 & 1 \end{bmatrix} \end{bmatrix}$$

and

$$Y = \begin{bmatrix} \begin{bmatrix} 1 & 1 & 2 & 1 \\ 1 & 0 & 2 & 4 \\ \hline 0 & 1 & 0 & 3 \\ 1 & 1 & 0 & 0 \\ \hline 1 & 0 & 1 & 1 \\ 0 & 1 & 0 & 1 \end{bmatrix} & \begin{bmatrix} 2 & 1 \\ 3 & 1 \\ \hline 1 & 0 \\ 2 & 1 \\ \hline 1 & 2 \\ 1 & 1 \end{bmatrix} & \begin{bmatrix} 3 & 1 & 0 \\ 4 & 1 & 1 \\ \hline 1 & 2 & 1 \\ 1 & 2 & 1 \\ \hline 2 & 1 & 2 \\ 1 & -1 & 0 \end{bmatrix} \end{bmatrix}.$$

Now we find the major product of XY. The product of the first row of X with first column of Y gives

$$\begin{bmatrix} 1 & 2 & 1 & 1 & 2 & 3 \\ 3 & 1 & 2 & 3 & 1 & 1 \end{bmatrix} \begin{bmatrix} 1 & 1 & 2 & 1 \\ 1 & 0 & 2 & 4 \\ 0 & 1 & 0 & 3 \\ 1 & 1 & 0 & 0 \\ 1 & 0 & 1 & 1 \\ 0 & 1 & 0 & 1 \end{bmatrix}$$

$$= \begin{bmatrix} 1 \\ 3 \end{bmatrix} \begin{bmatrix} 1 & 1 & 2 & 1 \end{bmatrix} + \begin{bmatrix} 2 & 1 \\ 1 & 2 \end{bmatrix} \begin{bmatrix} 1 & 0 & 2 & 4 \\ 0 & 1 & 0 & 3 \end{bmatrix} +$$

$$\begin{bmatrix} 1 & 2 & 3 \\ 3 & 1 & 1 \end{bmatrix} \begin{bmatrix} 1 & 1 & 0 & 0 \\ 1 & 0 & 1 & 1 \\ 0 & 1 & 0 & 1 \end{bmatrix}$$

$$= \begin{bmatrix} 1 & 1 & 2 & 1 \\ 3 & 3 & 6 & 3 \end{bmatrix} + \begin{bmatrix} 2 & 1 & 4 & 11 \\ 1 & 2 & 2 & 10 \end{bmatrix} + \begin{bmatrix} 3 & 4 & 2 & 5 \\ 4 & 4 & 1 & 2 \end{bmatrix}$$

$$= \begin{bmatrix} 6 & 6 & 8 & 17 \\ 8 & 9 & 9 & 15 \end{bmatrix}.$$

Now

$$\begin{bmatrix} 1 & 2 & 1 & 1 & 2 & 3 \\ 3 & 1 & 2 & 3 & 1 & 1 \end{bmatrix} \begin{bmatrix} 2 & 1 \\ 3 & 1 \\ 1 & 0 \\ 2 & 1 \\ 1 & 2 \\ 1 & 1 \end{bmatrix}$$

$$= \begin{bmatrix} 1 \\ 3 \end{bmatrix} \begin{bmatrix} 2 & 1 \end{bmatrix} + \begin{bmatrix} 2 & 1 \\ 1 & 2 \end{bmatrix} \begin{bmatrix} 3 & 1 \\ 1 & 0 \end{bmatrix} + \begin{bmatrix} 1 & 2 & 3 \\ 3 & 1 & 1 \end{bmatrix} \begin{bmatrix} 2 & 1 \\ 1 & 2 \\ 1 & 1 \end{bmatrix}$$

$$= \begin{bmatrix} 2 & 1 \\ 6 & 3 \end{bmatrix} + \begin{bmatrix} 7 & 2 \\ 5 & 1 \end{bmatrix} + \begin{bmatrix} 7 & 8 \\ 8 & 6 \end{bmatrix}$$

$$= \begin{bmatrix} 16 & 11 \\ 19 & 10 \end{bmatrix}.$$

Consider the product of first row with the 3$^{rd}$ column.

$$\begin{bmatrix} 1 & 2 & 1 & 1 & 2 & 3 \\ 3 & 1 & 2 & 3 & 1 & 1 \end{bmatrix} \begin{bmatrix} 3 & 1 & 0 \\ 4 & 1 & 1 \\ 1 & 2 & 1 \\ 1 & 2 & 1 \\ 2 & 1 & 2 \\ 1 & -1 & 0 \end{bmatrix}$$

$$= \begin{bmatrix} 1 \\ 3 \end{bmatrix} \begin{bmatrix} 3 & 1 & 0 \end{bmatrix} + \begin{bmatrix} 2 & 1 \\ 1 & 2 \end{bmatrix} \begin{bmatrix} 4 & 1 & 1 \\ 1 & 2 & 1 \end{bmatrix} + \begin{bmatrix} 1 & 2 & 3 \\ 3 & 1 & 1 \end{bmatrix} \begin{bmatrix} 1 & 2 & 1 \\ 2 & 1 & 2 \\ 1 & -1 & 0 \end{bmatrix}$$

$$= \begin{bmatrix} 3 & 1 & 0 \\ 9 & 3 & 0 \end{bmatrix} + \begin{bmatrix} 9 & 4 & 3 \\ 6 & 5 & 3 \end{bmatrix} + \begin{bmatrix} 8 & 1 & 5 \\ 6 & 6 & 5 \end{bmatrix}$$

$$= \begin{bmatrix} 20 & 6 & 8 \\ 21 & 14 & 8 \end{bmatrix}.$$

The product of $2^{nd}$ row of X with first column of Y gives

$$\begin{bmatrix} 1 & 1 & 3 & 1 & 1 & 1 \\ 2 & 3 & 1 & 2 & 0 & 1 \\ 3 & 4 & 2 & 0 & 1 & 0 \\ 4 & 2 & 4 & 1 & 0 & 0 \end{bmatrix} \begin{bmatrix} 1 & 1 & 2 & 1 \\ 1 & 0 & 2 & 4 \\ 0 & 1 & 0 & 3 \\ 1 & 1 & 0 & 0 \\ 1 & 0 & 1 & 1 \\ 0 & 1 & 0 & 1 \end{bmatrix} =$$

$$\begin{bmatrix} 1 \\ 2 \\ 3 \\ 4 \end{bmatrix} \begin{bmatrix} 1 & 1 & 2 & 1 \end{bmatrix} + \begin{bmatrix} 1 & 3 \\ 3 & 1 \\ 4 & 2 \\ 2 & 4 \end{bmatrix} \begin{bmatrix} 1 & 0 & 2 & 4 \\ 0 & 1 & 0 & 3 \end{bmatrix} + \begin{bmatrix} 1 & 1 & 1 \\ 2 & 0 & 1 \\ 0 & 1 & 0 \\ 1 & 0 & 0 \end{bmatrix} \begin{bmatrix} 1 & 1 & 0 & 0 \\ 1 & 0 & 1 & 1 \\ 0 & 1 & 0 & 1 \end{bmatrix}$$

$$= \begin{bmatrix} 1 & 1 & 2 & 1 \\ 2 & 2 & 4 & 2 \\ 3 & 3 & 6 & 3 \\ 4 & 4 & 8 & 4 \end{bmatrix} + \begin{bmatrix} 1 & 3 & 2 & 13 \\ 3 & 1 & 6 & 15 \\ 4 & 2 & 8 & 22 \\ 2 & 4 & 4 & 20 \end{bmatrix} + \begin{bmatrix} 2 & 2 & 1 & 2 \\ 2 & 3 & 0 & 1 \\ 1 & 0 & 1 & 1 \\ 1 & 1 & 0 & 0 \end{bmatrix}$$

$$= \begin{bmatrix} 4 & 6 & 5 & 16 \\ 7 & 6 & 10 & 18 \\ 8 & 5 & 15 & 26 \\ 7 & 9 & 12 & 24 \end{bmatrix}.$$

The product of $3^{rd}$ row of X with the $3^{rd}$ column of Y.

$$[5 \mid 0 \ 1 \mid 1 \ 1 \ 1] \begin{bmatrix} 3 & 1 & 0 \\ 4 & 1 & 1 \\ \hline 1 & 2 & 1 \\ 1 & 2 & 1 \\ \hline 2 & 1 & 2 \\ 1 & -1 & 0 \end{bmatrix}$$

$$= [5][3 \ 1 \ 0] + [0 \ 1]\begin{bmatrix} 4 & 1 & 1 \\ 1 & 2 & 1 \end{bmatrix} + [1 \ 1 \ 1]\begin{bmatrix} 1 & 2 & 1 \\ 2 & 1 & 2 \\ 1 & -1 & 0 \end{bmatrix}$$

$$= [15 \ 5 \ 0] + [1 \ 2 \ 1] + [4 \ 2 \ 3]$$
$$= [20 \ 9 \ 4].$$

The product of second row of X with second column of Y.

$$\begin{bmatrix} 1 & | & 1 & 3 & | & 1 & 1 & 1 \\ 2 & | & 3 & 1 & | & 2 & 0 & 1 \\ 3 & | & 4 & 2 & | & 0 & 1 & 0 \\ 4 & | & 2 & 4 & | & 1 & 0 & 0 \end{bmatrix} \begin{bmatrix} 2 & 1 \\ 3 & 1 \\ \hline 1 & 0 \\ 2 & 1 \\ \hline 1 & 2 \\ 1 & 1 \end{bmatrix}$$

$$= \begin{bmatrix} 1 \\ 2 \\ 3 \\ 4 \end{bmatrix}[2 \ 1] + \begin{bmatrix} 1 & 3 \\ 3 & 1 \\ 4 & 2 \\ 2 & 4 \end{bmatrix}\begin{bmatrix} 3 & 1 \\ 1 & 0 \end{bmatrix} + \begin{bmatrix} 1 & 1 & 1 \\ 2 & 0 & 1 \\ 0 & 1 & 0 \\ 1 & 0 & 0 \end{bmatrix}\begin{bmatrix} 2 & 1 \\ 1 & 2 \\ 1 & 1 \end{bmatrix}$$

$$= \begin{bmatrix} 2 & 1 \\ 4 & 2 \\ 6 & 3 \\ 8 & 4 \end{bmatrix} + \begin{bmatrix} 6 & 1 \\ 10 & 3 \\ 14 & 4 \\ 10 & 2 \end{bmatrix} + \begin{bmatrix} 4 & 4 \\ 5 & 3 \\ 1 & 2 \\ 2 & 1 \end{bmatrix}$$

$$= \begin{bmatrix} 12 & 6 \\ 19 & 8 \\ 21 & 9 \\ 20 & 7 \end{bmatrix}.$$

The product of the 2nd row with the last column of Y.

$$\begin{bmatrix} 1 & 1 & 3 & 1 & 1 & 1 \\ 2 & 3 & 1 & 2 & 0 & 1 \\ 3 & 4 & 2 & 0 & 1 & 0 \\ 4 & 2 & 4 & 1 & 0 & 0 \end{bmatrix} \begin{bmatrix} 3 & 1 & 0 \\ 4 & 1 & 1 \\ 1 & 2 & 1 \\ 1 & 2 & 1 \\ 2 & 1 & 2 \\ 1 & -1 & 0 \end{bmatrix}$$

$$= \begin{bmatrix} 1 \\ 2 \\ 3 \\ 4 \end{bmatrix} \begin{bmatrix} 3 & 1 & 0 \end{bmatrix} + \begin{bmatrix} 1 & 3 \\ 3 & 1 \\ 4 & 2 \\ 2 & 4 \end{bmatrix} \begin{bmatrix} 4 & 1 & 1 \\ 1 & 2 & 1 \end{bmatrix} + \begin{bmatrix} 1 & 1 & 1 \\ 2 & 0 & 1 \\ 0 & 1 & 0 \\ 1 & 0 & 0 \end{bmatrix} \begin{bmatrix} 1 & 2 & 1 \\ 2 & 1 & 2 \\ 1 & -1 & 0 \end{bmatrix}$$

$$= \begin{bmatrix} 3 & 1 & 0 \\ 6 & 2 & 0 \\ 9 & 3 & 0 \\ 12 & 4 & 0 \end{bmatrix} + \begin{bmatrix} 7 & 7 & 4 \\ 13 & 5 & 4 \\ 18 & 8 & 6 \\ 12 & 10 & 6 \end{bmatrix} + \begin{bmatrix} 4 & 2 & 3 \\ 3 & 3 & 2 \\ 2 & 1 & 2 \\ 1 & 2 & 1 \end{bmatrix}$$

$$= \begin{bmatrix} 14 & 10 & 7 \\ 22 & 10 & 6 \\ 29 & 12 & 8 \\ 25 & 16 & 7 \end{bmatrix}.$$

The product of 3$^{rd}$ row of X with 1$^{st}$ column of Y

$$[5 \mid 0 \ 1 \mid 1 \ 1 \ 1] \begin{bmatrix} 1 & 1 & 2 & 1 \\ 1 & 0 & 2 & 4 \\ \hline 0 & 1 & 0 & 3 \\ 1 & 1 & 0 & 0 \\ \hline 1 & 0 & 1 & 1 \\ 0 & 1 & 0 & 1 \end{bmatrix} =$$

$$[5][1 \ 1 \ 2 \ 1] + [0 \ 1]\begin{bmatrix} 1 & 0 & 2 & 4 \\ 0 & 1 & 0 & 3 \end{bmatrix} + [1 \ 1 \ 1]\begin{bmatrix} 1 & 1 & 0 & 0 \\ 1 & 0 & 1 & 1 \\ 0 & 1 & 0 & 1 \end{bmatrix}$$

$= [5 \ 5 \ 10 \ 5] + [0 \ 1 \ 0 \ 3] + [2 \ 2 \ 1 \ 2]$

$= [7 \ 8 \ 11 \ 10]$.

The product of 3$^{rd}$ row of X with 2$^{nd}$ column of Y.

$$[5 \mid 0 \ 1 \mid 1 \ 1 \ 1] \begin{bmatrix} 2 & 1 \\ 3 & 1 \\ \hline 1 & 0 \\ 2 & 1 \\ \hline 1 & 2 \\ 1 & 1 \end{bmatrix} =$$

$$[5] \ [2 \ 1] + [0 \ 1]\begin{bmatrix} 3 & 1 \\ 1 & 0 \end{bmatrix} + [1 \ 1 \ 1]\begin{bmatrix} 2 & 1 \\ 1 & 2 \\ 1 & 1 \end{bmatrix}$$

$= \ [10 \ 5] + [1 \ 0] + [4 \ 4]$
$= \ [15 \ 9]$.

$$XY = \begin{bmatrix} 6 & 6 & 8 & 17 & 16 & 11 & 20 & 6 & 8 \\ 8 & 8 & 9 & 15 & 19 & 10 & 21 & 14 & 8 \\ 4 & 6 & 5 & 16 & 12 & 6 & 14 & 10 & 7 \\ 7 & 6 & 10 & 18 & 19 & 8 & 22 & 10 & 6 \\ 8 & 5 & 15 & 26 & 21 & 9 & 29 & 12 & 8 \\ 7 & 9 & 12 & 24 & 20 & 7 & 25 & 16 & 7 \\ 7 & 8 & 11 & 10 & 15 & 9 & 20 & 9 & 4 \end{bmatrix}_{7 \times 9}$$

On similar lines we can find the transpose of major product of Type IV vectors.

Now we proceed on to just show the major product moment of a type IV vector.

*Example 1.1.30:* Suppose

$$X = \begin{bmatrix} 1 & 2 & 1 & 3 & 2 & 1 \\ 2 & 3 & 1 & 2 & 1 & 2 \\ 1 & 4 & 2 & 3 & 2 & 2 \\ 4 & 1 & 3 & 2 & 1 & 1 \\ 2 & 3 & 2 & 3 & 2 & 3 \\ 3 & 4 & 1 & 1 & 4 & 2 \\ 2 & 1 & 2 & 2 & 1 & 3 \end{bmatrix}$$

and

$$X^t = \begin{bmatrix} 1 & 2 & 1 & 4 & 2 & 3 & 2 \\ 2 & 3 & 4 & 1 & 3 & 4 & 1 \\ 1 & 1 & 2 & 3 & 2 & 1 & 2 \\ 3 & 2 & 3 & 2 & 3 & 1 & 2 \\ 2 & 1 & 2 & 1 & 2 & 4 & 1 \\ 1 & 2 & 2 & 1 & 3 & 2 & 3 \end{bmatrix}.$$

$$X^t X = \begin{bmatrix} 1 & 2 & 1 & 4 & 2 & 3 & 2 \\ 2 & 3 & 4 & 1 & 3 & 4 & 1 \\ 1 & 1 & 2 & 3 & 2 & 1 & 2 \\ \hline 3 & 2 & 3 & 2 & 3 & 1 & 2 \\ 2 & 1 & 2 & 1 & 2 & 4 & 1 \\ 1 & 2 & 2 & 1 & 3 & 2 & 3 \end{bmatrix} \times$$

$$\begin{bmatrix} 1 & 2 & 1 & 3 & 2 & 1 \\ 2 & 3 & 1 & 2 & 1 & 2 \\ \hline 1 & 4 & 2 & 3 & 2 & 2 \\ 4 & 1 & 3 & 2 & 1 & 1 \\ 2 & 3 & 2 & 3 & 2 & 3 \\ 3 & 4 & 1 & 1 & 4 & 2 \\ \hline 2 & 1 & 2 & 2 & 1 & 3 \end{bmatrix}.$$

Product of 1st row of $X^t$ with 1st column of X

$$\begin{bmatrix} 1 & 2 & | & 1 & 4 & 2 & 3 & | & 2 \end{bmatrix} \begin{bmatrix} 1 \\ 2 \\ \hline 1 \\ 4 \\ 2 \\ 3 \\ \hline 2 \end{bmatrix}$$

$$= \begin{bmatrix} 1 & 2 \end{bmatrix} \begin{bmatrix} 1 \\ 2 \end{bmatrix} + \begin{bmatrix} 1 & 4 & 2 & 3 \end{bmatrix} \begin{bmatrix} 1 \\ 4 \\ 2 \\ 3 \end{bmatrix} + \begin{bmatrix} 2 \end{bmatrix} \begin{bmatrix} 2 \end{bmatrix}$$

$$= 5 + 30 + 4$$
$$= 39.$$

Product of $1^{st}$ row of $X^t$ with $2^{nd}$ column of X.

$$[1 \quad 2 \mid 1 \quad 4 \quad 2 \quad 3 \mid 2] \begin{bmatrix} 2 & 1 \\ 3 & 1 \\ \hline 4 & 2 \\ 1 & 3 \\ 3 & 2 \\ 4 & 1 \\ \hline 1 & 2 \end{bmatrix}$$

$$= [1 \quad 2]\begin{bmatrix} 2 & 1 \\ 3 & 1 \end{bmatrix} + [1 \quad 4 \quad 2 \quad 3]\begin{bmatrix} 4 & 2 \\ 1 & 3 \\ 3 & 2 \\ 4 & 1 \end{bmatrix} + [2][1 \quad 2]$$

$= [8 \ 3] + [26 \ 21] + [2 \ 4]$
$= [36 \ 28]$.

The product of $1^{st}$ row of $X^t$ with $3^{rd}$ column of X.

$$[1 \quad 2 \mid 1 \quad 4 \quad 2 \quad 3 \mid 2] \begin{bmatrix} 3 & 2 & 1 \\ 2 & 1 & 2 \\ \hline 3 & 2 & 2 \\ 2 & 1 & 1 \\ 3 & 2 & 3 \\ 1 & 4 & 2 \\ \hline 2 & 1 & 3 \end{bmatrix}$$

$$= [1 \quad 2]\begin{bmatrix} 3 & 2 & 1 \\ 2 & 1 & 2 \end{bmatrix} + [1 \quad 4 \quad 2 \quad 3]\begin{bmatrix} 3 & 2 & 2 \\ 2 & 1 & 1 \\ 3 & 2 & 3 \\ 1 & 4 & 2 \end{bmatrix} + [2][2 \quad 1 \quad 3]$$

$$= [7\ 4\ 5] + [20\ 22\ 18] + [4\ 2\ 6]$$

$$= [31\ 28\ 29].$$

The product of 2$^{nd}$ row of X$^t$ with 1$^{st}$ column of X.

$$\begin{bmatrix} 2 & 3 & | & 4 & 1 & 3 & 4 & | & 1 \\ 1 & 1 & | & 2 & 3 & 2 & 1 & | & 2 \end{bmatrix} \begin{bmatrix} 1 \\ 2 \\ \hline 1 \\ 4 \\ 2 \\ 3 \\ \hline 2 \end{bmatrix}$$

$$= \begin{bmatrix} 2 & 3 \\ 1 & 1 \end{bmatrix} \begin{bmatrix} 1 \\ 2 \end{bmatrix} + \begin{bmatrix} 4 & 1 & 3 & 4 \\ 2 & 3 & 2 & 1 \end{bmatrix} \begin{bmatrix} 1 \\ 4 \\ 2 \\ 3 \end{bmatrix} + \begin{bmatrix} 1 \\ 2 \end{bmatrix} [2]$$

$$= \begin{bmatrix} 8 \\ 3 \end{bmatrix} + \begin{bmatrix} 26 \\ 21 \end{bmatrix} + \begin{bmatrix} 2 \\ 4 \end{bmatrix} = \begin{bmatrix} 36 \\ 28 \end{bmatrix}.$$

The product of 2$^{nd}$ row of X$^t$ with 2$^{nd}$ column of X.

$$\begin{bmatrix} 2 & 3 & | & 4 & 1 & 3 & 4 & | & 1 \\ 1 & 1 & | & 2 & 3 & 2 & 1 & | & 2 \end{bmatrix} \begin{bmatrix} 2 & 1 \\ 3 & 1 \\ \hline 4 & 2 \\ 1 & 3 \\ 3 & 2 \\ 4 & 1 \\ \hline 1 & 2 \end{bmatrix}$$

$$= \begin{bmatrix} 2 & 3 \\ 1 & 1 \end{bmatrix} \begin{bmatrix} 2 & 1 \\ 3 & 1 \end{bmatrix} + \begin{bmatrix} 4 & 1 & 3 & 4 \\ 2 & 3 & 2 & 1 \end{bmatrix} \begin{bmatrix} 4 & 2 \\ 1 & 3 \\ 3 & 2 \\ 4 & 1 \end{bmatrix} + \begin{bmatrix} 1 \\ 2 \end{bmatrix} \begin{bmatrix} 1 & 2 \end{bmatrix}$$

$$= \begin{bmatrix} 13 & 5 \\ 5 & 2 \end{bmatrix} + \begin{bmatrix} 42 & 21 \\ 21 & 18 \end{bmatrix} + \begin{bmatrix} 1 & 2 \\ 2 & 4 \end{bmatrix}$$

$$= \begin{bmatrix} 56 & 28 \\ 28 & 24 \end{bmatrix}.$$

The product of $2^{nd}$ row of $X^t$ with $3^{rd}$ column of X.

$$\begin{bmatrix} 2 & 3 & | & 4 & 1 & 3 & 4 & | & 1 \\ 1 & 1 & | & 2 & 3 & 2 & 1 & | & 2 \end{bmatrix} \begin{bmatrix} 3 & 2 & 1 \\ 2 & 1 & 2 \\ \hline 3 & 2 & 2 \\ 2 & 1 & 1 \\ 3 & 2 & 3 \\ 1 & 4 & 2 \\ \hline 2 & 1 & 3 \end{bmatrix}$$

$$= \begin{bmatrix} 2 & 3 \\ 1 & 1 \end{bmatrix} \begin{bmatrix} 3 & 2 & 1 \\ 2 & 1 & 2 \end{bmatrix} + \begin{bmatrix} 4 & 1 & 3 & 4 \\ 2 & 3 & 2 & 1 \end{bmatrix} \begin{bmatrix} 3 & 2 & 2 \\ 2 & 1 & 1 \\ 3 & 2 & 3 \\ 1 & 4 & 2 \end{bmatrix} + \begin{bmatrix} 1 \\ 2 \end{bmatrix} \begin{bmatrix} 2 & 1 & 3 \end{bmatrix}$$

$$= \begin{bmatrix} 12 & 7 & 8 \\ 5 & 3 & 3 \end{bmatrix} + \begin{bmatrix} 27 & 31 & 26 \\ 19 & 15 & 15 \end{bmatrix} + \begin{bmatrix} 2 & 1 & 3 \\ 4 & 2 & 6 \end{bmatrix}$$

$$= \begin{bmatrix} 41 & 39 & 37 \\ 28 & 20 & 24 \end{bmatrix}.$$

The product of $3^{rd}$ row of $X^t$ with $1^{st}$ column of X.

$$\begin{bmatrix} 3 & 2 & | & 3 & 2 & 3 & 1 & | & 2 \\ 2 & 1 & | & 2 & 1 & 2 & 4 & | & 1 \\ 1 & 2 & | & 2 & 1 & 3 & 2 & | & 3 \end{bmatrix} \begin{bmatrix} 1 \\ 2 \\ \hline 1 \\ 4 \\ 2 \\ 3 \\ \hline 2 \end{bmatrix}$$

$$= \begin{bmatrix} 7 \\ 4 \\ 5 \end{bmatrix} + \begin{bmatrix} 20 \\ 22 \\ 18 \end{bmatrix} + \begin{bmatrix} 4 \\ 2 \\ 6 \end{bmatrix} = \begin{bmatrix} 31 \\ 28 \\ 29 \end{bmatrix}.$$

The product of $3^{rd}$ row of $X^t$ with $2^{nd}$ column of X.

$$\begin{bmatrix} 3 & 2 & | & 3 & 2 & 3 & 1 & | & 2 \\ 2 & 1 & | & 2 & 1 & 2 & 4 & | & 1 \\ 1 & 2 & | & 2 & 1 & 3 & 2 & | & 3 \end{bmatrix} \begin{bmatrix} 2 & 1 \\ 3 & 1 \\ \hline 4 & 2 \\ 1 & 3 \\ 3 & 2 \\ 4 & 1 \\ \hline 1 & 2 \end{bmatrix}$$

$$= \begin{bmatrix} 12 & 5 \\ 7 & 3 \\ 8 & 3 \end{bmatrix} + \begin{bmatrix} 27 & 19 \\ 31 & 15 \\ 26 & 15 \end{bmatrix} + \begin{bmatrix} 2 & 4 \\ 1 & 2 \\ 3 & 6 \end{bmatrix}$$

$$= \begin{bmatrix} 41 & 28 \\ 39 & 20 \\ 37 & 24 \end{bmatrix}.$$

The product 3$^{rd}$ row of X$^t$ with 3$^{rd}$ column of X.

$$\begin{bmatrix} 3 & 2 \\ 2 & 1 \\ 1 & 2 \end{bmatrix} \begin{bmatrix} 3 & 2 & 3 & 1 \\ 2 & 1 & 2 & 4 \\ 2 & 1 & 3 & 2 \end{bmatrix} \begin{bmatrix} 2 \\ 1 \\ 3 \end{bmatrix} \begin{bmatrix} 3 & 2 & 1 \\ 2 & 1 & 2 \\ 3 & 2 & 2 \\ 2 & 1 & 1 \\ 3 & 2 & 3 \\ 1 & 4 & 2 \\ 2 & 1 & 3 \end{bmatrix}$$

$$= \begin{bmatrix} 13 & 8 & 7 \\ 8 & 5 & 4 \\ 7 & 4 & 5 \end{bmatrix} + \begin{bmatrix} 23 & 18 & 19 \\ 18 & 25 & 19 \\ 19 & 19 & 18 \end{bmatrix} + \begin{bmatrix} 4 & 2 & 6 \\ 2 & 1 & 3 \\ 6 & 3 & 9 \end{bmatrix}$$

$$= \begin{bmatrix} 40 & 28 & 32 \\ 28 & 31 & 26 \\ 32 & 26 & 32 \end{bmatrix}.$$

$$X^t X = \begin{bmatrix} 39 & 36 & 28 & 31 & 28 & 29 \\ 36 & 56 & 28 & 41 & 39 & 37 \\ 28 & 28 & 24 & 28 & 20 & 24 \\ 31 & 41 & 28 & 40 & 28 & 32 \\ 28 & 39 & 20 & 28 & 31 & 26 \\ 29 & 37 & 24 & 32 & 26 & 32 \end{bmatrix}.$$

On similar lines interested reader can find the major product moment of type IV column vector.

### 1.2 Bimatrices and their Generalizations

In this section we recall some of the basic properties of bimatrices and their generalizations which will be useful for us in the definition of linear bicodes and linear n-codes respectively.

In this section we recall the notion of bimatrix and illustrate them with examples and define some of basic operations on them.

**DEFINITION 1.2.1:** *A bimatrix $A_B$ is defined as the union of two rectangular array of numbers $A_1$ and $A_2$ arranged into rows and columns. It is written as follows $A_B = A_1 \cup A_2$ where $A_1 \neq A_2$ with*

$$A_1 = \begin{bmatrix} a_{11}^1 & a_{12}^1 & \cdots & a_{1n}^1 \\ a_{21}^1 & a_{22}^1 & \cdots & a_{2n}^1 \\ \vdots & & & \\ a_{m1}^1 & a_{m2}^1 & \cdots & a_{mn}^1 \end{bmatrix}$$

*and*

$$A_2 = \begin{bmatrix} a_{11}^2 & a_{12}^2 & \cdots & a_{1n}^2 \\ a_{21}^2 & a_{22}^2 & \cdots & a_{2n}^2 \\ \vdots & & & \\ a_{m1}^2 & a_{m2}^2 & \cdots & a_{mn}^2 \end{bmatrix}$$

*'$\cup$' is just the notational convenience (symbol) only.*

The above array is called a m by n bimatrix (written as B(m × n) since each of $A_i$ (i = 1, 2) has m rows and n columns). It is to be noted a bimatrix has no numerical value associated with it. It is only a convenient way of representing a pair of array of numbers.

*Note:* If $A_1 = A_2$ then $A_B = A_1 \cup A_2$ is not a bimatrix. A bimatrix $A_B$ is denoted by $\left(a_{ij}^1\right) \cup \left(a_{ij}^2\right)$. If both $A_1$ and $A_2$ are m × n matrices then the bimatrix $A_B$ is called the m × n rectangular bimatrix.

But we make an assumption the zero bimatrix is a union of two zero matrices even if $A_1$ and $A_2$ are one and the same; i.e., $A_1 = A_2 = (0)$.

*Example 1.2.1:* The following are bimatrices

i. $A_B = \begin{bmatrix} 3 & 0 & 1 \\ -1 & 2 & 1 \end{bmatrix} \cup \begin{bmatrix} 0 & 2 & -1 \\ 1 & 1 & 0 \end{bmatrix}$

   is a rectangular $2 \times 3$ bimatrix.

ii. $A'_B = \begin{bmatrix} 3 \\ 1 \\ 2 \end{bmatrix} \cup \begin{bmatrix} 0 \\ -1 \\ 0 \end{bmatrix}$

   is a column bimatrix.

iii. $A''_B = (3, -2, 0, 1, 1) \cup (1, 1, -1, 1, 2)$

   is a row bimatrix.

In a bimatrix $A_B = A_1 \cup A_2$ if both $A_1$ and $A_2$ are $m \times n$ rectangular matrices then the bimatrix $A_B$ is called the rectangular $m \times n$ bimatrix.

**DEFINITION 1.2.2:** *Let $A_B = A_1 \cup A_2$ be a bimatrix. If both $A_1$ and $A_2$ are square matrices then $A_B$ is called the square bimatrix.*

*If one of the matrices in the bimatrix $A_B = A_1 \cup A_2$ is a square matrix and other is a rectangular matrix or if both $A_1$ and $A_2$ are rectangular matrices say $m_1 \times n_1$ and $m_2 \times n_2$ with $m_1 \neq m_2$ or $n_1 \neq n_2$ then we say $A_B$ is a mixed bimatrix.*

The following are examples of a square bimatrix and the mixed bimatrix.

***Example 1.2.2:*** Given

$$A_B = \begin{bmatrix} 3 & 0 & 1 \\ 2 & 1 & 1 \\ -1 & 1 & 0 \end{bmatrix} \cup \begin{bmatrix} 4 & 1 & 1 \\ 2 & 1 & 0 \\ 0 & 0 & 1 \end{bmatrix}$$

is a 3 × 3 square bimatrix.

$$A'_B = \begin{bmatrix} 1 & 1 & 0 & 0 \\ 2 & 0 & 0 & 1 \\ 0 & 0 & 0 & 3 \\ 1 & 0 & 1 & 2 \end{bmatrix} \cup \begin{bmatrix} 2 & 0 & 0 & -1 \\ -1 & 0 & 1 & 0 \\ 0 & -1 & 0 & 3 \\ -3 & -2 & 0 & 0 \end{bmatrix}$$

is a 4 × 4 square bimatrix.

***Example 1.2.3:*** Let

$$A_B = \begin{bmatrix} 3 & 0 & 1 & 2 \\ 0 & 0 & 1 & 1 \\ 2 & 1 & 0 & 0 \\ 1 & 0 & 1 & 0 \end{bmatrix} \cup \begin{bmatrix} 1 & 1 & 2 \\ 0 & 2 & 1 \\ 0 & 0 & 4 \end{bmatrix}$$

then $A_B$ is a mixed square bimatrix.
Let

$$A'_B = \begin{bmatrix} 2 & 0 & 1 & 1 \\ 0 & 1 & 0 & 1 \\ -1 & 0 & 2 & 1 \end{bmatrix} \cup \begin{bmatrix} 2 & 0 \\ 4 & -3 \end{bmatrix},$$

$A'_B$ is a mixed bimatrix.

Now we proceed on to give the operations on bimatrices.

Let $A_B = A_1 \cup A_2$ and $C_B = C_1 \cup C_2$ be two bimatrices we say $A_B$ and $C_B$ are equal written as $A_B = C_B$ if and only if $A_1$ and $C_1$ are identical and $A_2$ and $C_2$ are identical i.e., $A_1 = C_1$ and $A_2 = C_2$.

If $A_B = A_1 \cup A_2$ and $C_B = C_1 \cup C_2$, we say $A_B$ is not equal to $C_B$, we write $A_B \neq C_B$ if and only if $A_1 \neq C_1$ or $A_2 \neq C_2$.

*Example 1.2.4:* Let

$$A_B = \begin{bmatrix} 3 & 2 & 0 \\ 2 & 1 & 1 \end{bmatrix} \cup \begin{bmatrix} 0 & -1 & 2 \\ 0 & 0 & 1 \end{bmatrix}$$

and

$$C_B = \begin{bmatrix} 1 & 1 & 1 \\ 0 & 0 & 0 \end{bmatrix} \cup \begin{bmatrix} 2 & 0 & 1 \\ 1 & 0 & 2 \end{bmatrix}$$

clearly $A_B \neq C_B$. Let

$$A_B = \begin{bmatrix} 0 & 0 & 1 \\ 1 & 1 & 2 \end{bmatrix} \cup \begin{bmatrix} 0 & 4 & -2 \\ -3 & 0 & 0 \end{bmatrix}$$

$$C_B = \begin{bmatrix} 0 & 0 & 1 \\ 1 & 1 & 2 \end{bmatrix} \cup \begin{bmatrix} 0 & 0 & 0 \\ 1 & 0 & 1 \end{bmatrix}$$

clearly $A_B \neq C_B$.

If $A_B = C_B$ then we have $C_B = A_B$.

We now proceed on to define multiplication by a scalar. Given a bimatrix $A_B = A_1 \cup B_1$ and a scalar $\lambda$, the product of $\lambda$ and $A_B$ written $\lambda A_B$ is defined to be

$$\lambda A_B = \begin{bmatrix} \lambda a_{11} & \cdots & \lambda a_{1n} \\ \vdots & & \vdots \\ \lambda a_{m1} & \cdots & \lambda a_{mn} \end{bmatrix} \cup \begin{bmatrix} \lambda b_{11} & \cdots & \lambda b_{1n} \\ \vdots & & \vdots \\ \lambda b_{m1} & \cdots & \lambda b_{mn} \end{bmatrix}$$

each element of $A_1$ and $B_1$ are multiplied by $\lambda$. The product $\lambda A_B$ is then another bimatrix having m rows and n columns if $A_B$ has m rows and n columns.

We write
$$\lambda A_B = \left[\lambda a_{ij}\right] \cup \left[\lambda b_{ij}\right]$$
$$= \left[a_{ij}\lambda\right] \cup \left[b_{ij}\lambda\right]$$
$$= A_B \lambda.$$

*Example 1.2.5:* Let
$$A_B = \begin{bmatrix} 2 & 0 & 1 \\ 3 & 3 & -1 \end{bmatrix} \cup \begin{bmatrix} 0 & 1 & -1 \\ 2 & 1 & 0 \end{bmatrix}$$

and $\lambda = 3$ then
$$3A_B = \begin{bmatrix} 6 & 0 & 3 \\ 9 & 9 & -3 \end{bmatrix} \cup \begin{bmatrix} 0 & 3 & -3 \\ 6 & 3 & 0 \end{bmatrix}.$$

If $\lambda = -2$ for
$$A_B = [3\ 1\ 2\ -4] \cup [0\ 1\ -1\ 0],$$
$$\lambda A_B = [-6\ -2\ -4\ 8] \cup [0\ -2\ 2\ 0].$$

Let $A_B = A_1 \cup B_1$ and $C_B = A_2 \cup B_2$ be any two $m \times n$ bimatrices. The sum $D_B$ of the bimatrices $A_B$ and $C_B$ is defined as $D_B = A_B + C_B = [A_1 \cup B_1] + [A_2 \cup B_2] = (A_1 + A_2) \cup [B_2 + B_2]$; where $A_1 + A_2$ and $B_1 + B_2$ are the usual addition of matrices i.e., if

$$A_B = \left(a_{ij}^1\right) \cup \left(b_{ij}^1\right)$$

and

$$C_B = \left(a_{ij}^2\right) \cup \left(b_{ij}^2\right)$$

then

$$A_B + C_B = D_B = \left(a_{ij}^1 + a_{ij}^2\right) \cup \left(b_{ij}^1 + b_{ij}^2\right)\ (\forall ij).$$

If we write in detail

$$A_B = \begin{bmatrix} a_{11}^1 & \cdots & a_{1n}^1 \\ \vdots & & \vdots \\ a_{m1}^1 & \cdots & a_{mn}^1 \end{bmatrix} \cup \begin{bmatrix} b_{11}^1 & \cdots & b_{1n}^1 \\ \vdots & & \vdots \\ b_{m1}^1 & \cdots & b_{mn}^1 \end{bmatrix}$$

$$C_B = \begin{bmatrix} a_{11}^2 & \cdots & a_{1n}^2 \\ \vdots & & \vdots \\ a_{m1}^2 & \cdots & a_{mn}^2 \end{bmatrix} \cup \begin{bmatrix} b_{11}^2 & \cdots & b_{1n}^2 \\ \vdots & & \vdots \\ b_{m1}^2 & \cdots & b_{mn}^2 \end{bmatrix}$$

$$A_B + C_B =$$

$$\begin{bmatrix} a_{11}^1 + a_{11}^2 & \cdots & a_{1n}^1 + a_{1n}^2 \\ \vdots & & \vdots \\ a_{m1}^1 + a_{m1}^2 & \cdots & a_{mn}^1 + a_{mn}^2 \end{bmatrix} \cup \begin{bmatrix} b_{11}^1 + b_{11}^2 & \cdots & b_{1n}^1 + b_{1n}^2 \\ \vdots & & \vdots \\ b_{m1}^1 + b_{m1}^2 & \cdots & b_{mn}^1 + b_{mn}^2 \end{bmatrix}.$$

The expression is abbreviated to

$$\begin{aligned} D_B &= A_B + C_B \\ &= (A_1 \cup B_1) + (A_2 \cup B_2) \\ &= (A_1 + A_2) \cup (B_1 + B_2). \end{aligned}$$

Thus two bimatrices are added by adding the corresponding elements only when compatibility of usual matrix addition exists.

*Note*: If $A_B = A^1 \cup A^2$ be a bimatrix we call $A^1$ and $A^2$ as the components of $A_B$ or component matrices of the bimatrix $A_B$.

*Example 1.2.6:*

(i) Let
$$A_B = \begin{bmatrix} 3 & 1 & 1 \\ -1 & 0 & 2 \end{bmatrix} \cup \begin{bmatrix} 4 & 0 & -1 \\ 0 & 1 & -2 \end{bmatrix}$$
and

$$C_B = \begin{bmatrix} -1 & 0 & 1 \\ 2 & 2 & -1 \end{bmatrix} \cup \begin{bmatrix} 3 & 3 & 1 \\ 0 & 2 & -1 \end{bmatrix},$$

then,

$$D_B = A_B + C_B$$

$$= \begin{bmatrix} 3 & 1 & 1 \\ -1 & 0 & 2 \end{bmatrix} + \begin{bmatrix} -1 & 0 & 1 \\ 2 & 2 & -1 \end{bmatrix} \cup$$

$$\begin{bmatrix} 4 & 0 & -1 \\ 0 & 1 & -2 \end{bmatrix} + \begin{bmatrix} 3 & 3 & 1 \\ 0 & 2 & -1 \end{bmatrix}$$

$$= \begin{bmatrix} 2 & 1 & 2 \\ 1 & 2 & 1 \end{bmatrix} \cup \begin{bmatrix} 7 & 3 & 0 \\ 0 & 3 & -3 \end{bmatrix}.$$

(ii) Let
$$A_B = (3\ 2\ -1\ 0\ 1) \cup (0\ 1\ 1\ 0\ -1)$$
and
$$C_B = (1\ 1\ 1\ 1\ 1) \cup (5\ -1\ 2\ 0\ 3),$$

$$A_B + C_B = (4\ 3\ 0\ 1\ 2) \cup (5\ 0\ 3\ 0\ 2).$$

*Example 1.2.7:* Let

$$A_B = \begin{bmatrix} 6 & -1 \\ 2 & 2 \\ 1 & -1 \end{bmatrix} \cup \begin{bmatrix} 3 & 1 \\ 0 & 2 \\ -1 & 3 \end{bmatrix}$$

and

$$C_B = \begin{bmatrix} 2 & -4 \\ 4 & -1 \\ 3 & 0 \end{bmatrix} \cup \begin{bmatrix} 1 & 4 \\ 2 & 1 \\ 3 & 1 \end{bmatrix}.$$

$$A_B + A_B = \begin{bmatrix} 12 & -2 \\ 4 & 4 \\ 2 & -2 \end{bmatrix} \cup \begin{bmatrix} 6 & 2 \\ 0 & 4 \\ -2 & 6 \end{bmatrix} = 2A_B$$

$$C_B + C_B = \begin{bmatrix} 4 & -8 \\ 8 & -2 \\ 6 & 0 \end{bmatrix} \cup \begin{bmatrix} 2 & 8 \\ 4 & 2 \\ 6 & 2 \end{bmatrix} = 2C_B.$$

Similarly we can add

$$A_B + A_B + A_B = 3A_B = \begin{bmatrix} 18 & -3 \\ 6 & 6 \\ 3 & -3 \end{bmatrix} \cup \begin{bmatrix} 9 & 3 \\ 0 & 6 \\ -3 & 9 \end{bmatrix}.$$

*Note:* Addition of bimatrices are defined if and only if both the bimatrices are m × n bimatrices.

Let

$$A_B = \begin{bmatrix} 3 & 0 & 1 \\ 1 & 2 & 0 \end{bmatrix} \cup \begin{bmatrix} 1 & 1 & 1 \\ 0 & 2 & -1 \end{bmatrix}$$

and

$$C_B = \begin{bmatrix} 3 & 1 \\ 2 & 1 \\ 0 & 0 \end{bmatrix} \cup \begin{bmatrix} 1 & 1 \\ 2 & -1 \\ 3 & 0 \end{bmatrix}.$$

The addition of $A_B$ with $C_B$ is not defined for $A_B$ is a 2 × 3 bimatrix where as $C_B$ is a 3 × 2 bimatrix.

Clearly $A_B + C_B = C_B + A_B$ when both $A_B$ and $C_B$ are m × n matrices.

Also if $A_B$, $C_B$, $D_B$ be any three m × n bimatrices then $A_B + (C_B + D_B) = (A_B + C_B) + D_B$.

Subtraction is defined in terms of operations already considered for if

$$A_B = A_1 \cup A_2$$

and

$$B_B = B_1 \cup B_2$$

then

$$\begin{aligned} A_B - B_B &= A_B + (-B_B) \\ &= (A_1 \cup A_2) + (-B_1 \cup -B_2) \\ &= (A_1 - B_1) \cup (A_2 - B_2) \\ &= [A_1 + (-B_1)] \cup [A_2 + (-B_2)]. \end{aligned}$$

*Example 1.2.8:*

i. Let

$$A_B = \begin{bmatrix} 3 & 1 \\ -1 & 2 \\ 0 & 3 \end{bmatrix} \cup \begin{bmatrix} 5 & -2 \\ 1 & 1 \\ 3 & -2 \end{bmatrix}$$

and

$$B_B = \begin{bmatrix} 8 & -1 \\ 4 & 2 \\ -1 & 3 \end{bmatrix} \cup \begin{bmatrix} 9 & 2 \\ 2 & 9 \\ -1 & 1 \end{bmatrix}$$

$$A_B - B_B = A_B + (-B_B).$$

$$= \left\{ \begin{bmatrix} 3 & 1 \\ -1 & 2 \\ 0 & 3 \end{bmatrix} \cup \begin{bmatrix} 5 & -2 \\ 1 & 1 \\ 3 & -2 \end{bmatrix} \right\} + \left\{ - \begin{bmatrix} 8 & -1 \\ 4 & 2 \\ -1 & 3 \end{bmatrix} \cup \begin{bmatrix} 9 & 2 \\ 2 & 9 \\ -1 & 1 \end{bmatrix} \right\}$$

$$= \left\{ \begin{bmatrix} 3 & 1 \\ -1 & 2 \\ 0 & 3 \end{bmatrix} - \begin{bmatrix} 8 & -1 \\ 4 & 2 \\ -1 & 3 \end{bmatrix} \right\} \cup \left\{ \begin{bmatrix} 5 & -2 \\ 1 & 1 \\ 3 & -2 \end{bmatrix} - \begin{bmatrix} 9 & 2 \\ 2 & 9 \\ -1 & 1 \end{bmatrix} \right\}$$

$$= \begin{bmatrix} -5 & 2 \\ -5 & 0 \\ 1 & 0 \end{bmatrix} \cup \begin{bmatrix} 4 & -4 \\ -1 & -8 \\ 4 & -3 \end{bmatrix}.$$

ii. Let

$$A_B \;=\; (1, 2, 3, -1, 2, 1) \cup (3, -1, 2, 0, 3, 1)$$

and

$$B_B \;=\; (-1, 1, 1, 1, 1, 0) \cup (2, 0, -2, 0, 3, 0)$$

then

$$A_B + (-B_B) \;=\; (2, 1, 2, -2, 1, 1) \cup (1, -1, 4, 0, 0, 1).$$

Now we have defined addition and subtraction of bimatrices. Unlike in matrices we cannot say if we add two bimatrices the sum will be a bimatrix.

Now we proceed onto define the notion of n-matrices.

**DEFINITION 1.2.3:** *A n matrix A is defined to be the union of n rectangular array of numbers $A_1, \ldots, A_n$ arranged into rows and columns. It is written as $A = A_1 \cup \ldots \cup A_n$ where $A_i \neq A_j$ with*

$$A_i = \begin{bmatrix} a_{11}^i & a_{12}^i & \ldots & a_{1p}^i \\ a_{21}^i & a_{22}^i & \ldots & a_{2p}^i \\ \vdots & \vdots & & \vdots \\ a_{m1}^i & a_{m2}^i & \ldots & a_{mp}^i \end{bmatrix}$$

$i = 1, 2, \ldots, n.$

'$\cup$' is just the notional convenience (symbol) only ($n \geq 3$).

**Note**: If $n = 2$ we get the bimatrix.

*Example 1.2.9:* Let

$$A = \begin{bmatrix} 3 & 1 & 0 & 1 \\ 0 & 0 & 1 & 1 \end{bmatrix} \cup \begin{bmatrix} 2 & 1 & 1 & 0 \\ 0 & 1 & 1 & 0 \end{bmatrix} \cup$$

$$\begin{bmatrix} 1 & 0 & 0 & 1 \\ 0 & 1 & 0 & 1 \end{bmatrix} \cup \begin{bmatrix} 5 & 1 & 0 & 2 \\ 7 & -1 & 0 & 3 \end{bmatrix}$$

A is a 4-matrix.

*Example 1.2.10:* Let

$$A = A_1 \cup A_2 \cup A_3 \cup A_4 \cup A_5$$

$$= [1\ 0\ 0] \cup \begin{bmatrix} 1 \\ 2 \\ -1 \\ 0 \\ 0 \end{bmatrix} \cup \begin{bmatrix} 3 & 1 & 2 \\ 0 & 1 & 1 \\ 9 & 7 & -8 \end{bmatrix}$$

$$\cup \begin{bmatrix} 2 & 1 & 3 & 5 \\ 0 & 1 & 0 & 2 \end{bmatrix} \cup \begin{bmatrix} 7 & 9 & 8 & 11 & 0 \\ 1 & -2 & 0 & 9 & 7 \\ 0 & 5 & 7 & -1 & 8 \\ -4 & -6 & 6 & 0 & 1 \end{bmatrix};$$

A is a 5-matrix. Infact A is a mixed 5-matrix.

*Example 1.2.11:* Consider the 7-matrix

$$A = \begin{bmatrix} 2 & 0 \\ 1 & 1 \end{bmatrix} \cup \begin{bmatrix} 1 & 1 \\ 0 & 1 \end{bmatrix} \cup \begin{bmatrix} 1 & -1 \\ 0 & 4 \end{bmatrix} \cup$$

$$\begin{bmatrix} 3 & 1 \\ 1 & 0 \\ 0 & 1 \\ 2 & 1 \end{bmatrix} \cup \begin{bmatrix} 5 \\ 6 \\ 7 \\ 8 \\ -1 \\ 3 \\ 2 \end{bmatrix} \cup [3 \ 7 \ 8 \ 1 \ 0] \cup \begin{bmatrix} 2 & 1 & 0 & 0 \\ 1 & 0 & 2 & -1 \\ -1 & 1 & 0 & 0 \\ 0 & 2 & 0 & 1 \\ 2 & 0 & 0 & 1 \end{bmatrix}.$$

$= A_1 \cup A_2 \cup \ldots \cup A_7$. A is a mixed 7-matrix.

Chapter Two

# SUPERBIMATRICES AND THEIR PROPERTIES

In this chapter we introduce the notion of superbimatrices and explain some of its properties. We also give the type of products defined on them. Also the notion of semi superbimatrices and symmetric semi superbimatrices are introduced.

**DEFINITION 2.1:** *Let $A_1$ and $A_2$ be any two supermatrices, we call $A = A_1 \cup A_2$ to be a superbimatrix; '$\cup$' is just the symbol.*

**Note:** Further if $A_1 = A_2$, as non partitioned matrices then they must have distinct partitions. If $A_1 = A_2$ and they have the same set of partitions then we don't call $A = A_1 \cup A_2$ to be a superbimatrix.

We first illustrate this by the following examples.

*Example 2.1:* Let

$$A_1 = \begin{bmatrix} 3 & 1 & 0 & 2 \\ 1 & 1 & 6 & 0 \\ 0 & 1 & 0 & -1 \end{bmatrix}$$

and

$$A_2 = \begin{bmatrix} 1 & 1 & 0 & 1 & 0 \\ 2 & 0 & 2 & 0 & 1 \\ \hline 5 & 2 & 1 & 0 & 5 \\ 1 & 1 & 0 & 1 & 2 \end{bmatrix}$$

be any two supermatrices $A = A_1 \cup A_2$ is a superbimatrix.

$$A = \begin{bmatrix} 3 & 1 & 0 & 2 \\ 1 & 1 & 6 & 0 \\ 0 & 1 & 0 & -1 \end{bmatrix} \cup \begin{bmatrix} 1 & 1 & 0 & 1 & 0 \\ 2 & 0 & 2 & 0 & 1 \\ \hline 5 & 2 & 1 & 0 & 5 \\ 1 & 1 & 0 & 1 & 2 \end{bmatrix}$$

is a superbimatrix.

*Example 2.2:* Let

$$A_1 = \begin{bmatrix} 3 & 0 & 1 & 2 \\ 0 & 1 & 0 & 3 \\ \hline 1 & 1 & 5 & 2 \\ 0 & 0 & 2 & -1 \end{bmatrix}$$

and

$$A_2 = \begin{bmatrix} 3 & 0 & 1 & 2 \\ 0 & 1 & 0 & 3 \\ 1 & 1 & 5 & 2 \\ 0 & 0 & 2 & -1 \end{bmatrix}$$

be two supermatrices. $A = A_1 \cup A_2$ is a superbimatrix. We see clearly $A_1$ and $A_2$ are identical but only the partition on $A_1$ and $A_2$ is different. Hence A is a superbimatrix.

*Example 2.3:* Let

$$A_1 = \begin{bmatrix} 3 & 0 & 1 \\ 2 & 1 & 1 \\ 5 & 2 & 0 \end{bmatrix}$$

and

$$A_2 = \begin{bmatrix} 3 & | & 0 & 1 \\ 2 & | & 1 & 1 \\ 5 & | & 2 & 0 \end{bmatrix} = A_1.$$

Clearly $A = A_1 \cup A_2$ is not a superbimatrix.

*Example 2.4:* Let $A = A_1 \cup A_2$ be a superbimatrix where

$$A_1 = [3\ 0\ 1\ 2\ |\ 1\ 1\ 2\ |\ 1\ 5]$$

and

$$A_2 = \begin{bmatrix} 2 & 1 \\ 0 & 5 \\ \hline 1 & 1 \\ 0 & 1 \\ 3 & 2 \\ 1 & 0 \end{bmatrix}.$$

*Example 2.5:* Let $A = A_1 \cup A_2$ where $A_1 = [3\ 1\ 0\ 1\ |\ 5\ 0\ 2\ 3\ 1]$ and $A_2 = [\ 3\ 0\ 1\ |\ 2\ 2\ 0\ 5\ 3\ 1\ |\ 0\ 1\ 1]$. A is a superbimatrix in which we see both $A_1$ and $A_2$ are row supermatrices.

*Example 2.6:* Let $A = A_1 \cup A_2$ where

$$A_1 = \begin{bmatrix} 3 \\ 0 \\ 2 \\ \hline 1 \\ 10 \\ -1 \\ 5 \\ 4 \end{bmatrix} \text{ and } A_2 = \begin{bmatrix} 1 \\ 2 \\ 3 \\ 4 \\ \hline 5 \\ 6 \\ 7 \end{bmatrix}$$

are two super column matrices. Then A is a superbimatrix.

***Example 2.7:*** Let $A = A_1 \cup A_2$ where $A_1 = [3\ 1 \mid 2\ 0\ 5]$ and

$$A_2 = \begin{bmatrix} 1 \\ 1 \\ \hline 0 \\ 5 \\ -1 \end{bmatrix}$$

be two supermatrices. Clearly $A = A_1 \cup A_2$ is not a super column bimatrix or a super row bimatrix.

Now we have seen several examples of superbimatrices and we see each of them is of a specific type, so now we proceed on to define them.

**DEFINITION 2.2:** *Let $A = A_1 \cup A_2$ where $A_1 = [a_{11} \mid a_{12} \mid ... \mid a_{1n}]$ and $A_2 = [a_{21}\ a_{22}\ a_{23} \mid a_{24} ... a_{2m}]$ are both distinct super row matrices. Then we define $A = A_1 \cup A_2$ to be a super row bimatrix.*

The superbimatrix given in example 2.5 is a super row bimatrix.
    Now we proceed onto define the notion of super column bimatrix.

**DEFINITION 2.3:** *Let $A = A_1 \cup A_2$ where both $A_1$ and $A_2$ distinct column supermatrices,*

$$A_1 = \begin{bmatrix} a_1^1 \\ a_2^1 \\ a_3^1 \\ a_4^1 \\ \hline a_5^1 \\ \vdots \\ a_m^1 \end{bmatrix}$$

*and*

$$A_2 = \begin{bmatrix} a_1^2 \\ \hline a_2^2 \\ a_3^2 \\ \vdots \\ a_n^2 \end{bmatrix}.$$

*A is a superbimatrix which we define as the super column bimatrix or a column superbimatrix.*

The superbimatrix given in example 2.6 is a super column bimatrix.

Now we proceed onto define square superbimatrix.

**DEFINITION 2.4:** *Let $A = A_1 \cup A_2$ be a superbimatrix. If both $A_1$ and $A_2$ are distinct $m \times m$ square supermatrices then we call $A = A_1 \cup A_2$ to be a square superbimatrix.*

The superbimatrix given in example 2.2 is a square superbimatrix or to be more specific A is a 4 × 4 is square superbimatrix.

**Note:** If in the square supermatrix $A = A_1 \cup A_2$ if we have $A_1$ to be a m × m square supermatrix and $A_2$ to be a n × n square supermatrix (m ≠ n) then we call $A = A_1 \cup A_2$ to be a mixed square superbimatrix.

We now illustrate this by a simple example.

*Example 2.8:* Let $A = A_1 \cup A_2$ where

$$A_1 = \begin{bmatrix} 1 & 0 & 1 & 1 & 0 & 1 \\ 2 & 1 & 1 & 0 & 1 & 0 \\ \hline 1 & 0 & 1 & 2 & 0 & 3 \\ 3 & 5 & 1 & 0 & 5 & 1 \\ 1 & 1 & 1 & 0 & 1 & 0 \\ 7 & 2 & 1 & 3 & 1 & 0 \end{bmatrix}$$

and

$$A_2 = \left[\begin{array}{ccc|c} 3 & 1 & 0 & 2 \\ 1 & 0 & 1 & 1 \\ 2 & 1 & 0 & 0 \\ \hline 5 & 1 & 1 & 2 \end{array}\right];$$

both $A_1$ and $A_2$ are square supermatrices of different order. Thus $A = A_1 \cup A_2$ is a mixed square superbimatrix.

Now we proceed onto define rectangular superbimatrix and mixed rectangular superbimatrix.

**DEFINITION 2.5:** *Let $A = A_1 \cup A_2$ if both $A_1$ and $A_2$ are distinct $m \times n$ rectangular superbimatrices then we define $A = A_1 \cup A_2$ to be a rectangular superbimatrix. If $A_1$ is a $m_1 \times n_1$ rectangular supermatrix and $A_2$ is a $m_2 \times n_2$ rectangular supermatrix with $m_1 \neq m_2$ (or $n_1 \neq n_2$) then we call $A = A_1 \cup A_2$ to be a mixed rectangular superbimatrix.*

The example 2.1 is a superbimatrix which is a mixed rectangular superbimatrix.

Now we proceed onto give an example of a rectangular superbimatrix.

*Example 2.9:* Let $A = A_1 \cup A_2$ where

$$A_1 = \left[\begin{array}{cccc|cc} 3 & 1 & 2 & 3 & 5 & 0 \\ 1 & 0 & 2 & 1 & 2 & 1 \\ \hline 1 & 0 & 1 & 0 & 1 & 0 \\ 2 & 1 & 1 & 1 & 0 & 1 \end{array}\right]$$

and

$$A_2 = \left[\begin{array}{cc|c|ccc} 1 & 1 & 1 & 0 & 0 & 0 \\ 1 & 0 & 0 & 1 & 0 & 1 \\ 2 & 5 & 0 & 0 & 1 & 0 \\ \hline 3 & 1 & 2 & 2 & 1 & 0 \end{array}\right]$$

be any two rectangular supermatrices of same order i.e., both $A_1$ and $A_2$ are $4 \times 6$ supermatrices. We define $A = A_1 \cup A_2$ to be a rectangular superbimatrix of $4 \times 6$ order.

Now lastly we proceed onto define the notion of mixed superbimatrix.

**DEFINITION 2.6:** *Let $A = A_1 \cup A_2$ where $A_1$ and $A_2$ are supermatrices if $A_1$ is a square supermatrix and $A_2$ is a rectangular supermatrix then we define $A = A_1 \cup A_2$ to be a mixed superbimatrix.*

The superbimatrix given in all the examples is not a mixed superbimatrix. So now we proceed onto give an example of the same.

*Example 2.10:* Let $A = A_1 \cup A_2$ where

$$A_1 = \begin{bmatrix} 3 & 1 & | & 0 & 2 & 5 \\ 1 & 0 & | & 1 & 1 & 2 \\ \hline 2 & 1 & | & 0 & 5 & 1 \\ 5 & 1 & | & 3 & 2 & 1 \\ 3 & 2 & | & 0 & 0 & 3 \end{bmatrix}$$

and

$$A_2 = \begin{bmatrix} 4 & 0 & | & 1 & 1 & 1 & 2 \\ 2 & 1 & | & 0 & 2 & 1 & 5 \\ 3 & 1 & | & 3 & 0 & 1 & 0 \\ \hline 1 & 2 & | & 5 & 2 & 1 & 0 \end{bmatrix}$$

are two supermatrices where $A_1$ is a $5 \times 5$ square supermatrix and $A_2$ is a $4 \times 6$ rectangular supermatrix. Then $A = A_1 \cup A_2$ is a mixed superbimatrix.

Now having seen examples and definitions of several types of superbimatrices now we proceed onto define operations on them and the conditions under which operations are defined on them.

We first illustrate by some examples before we abstractly define those concepts.

*Example 2.11:* Let $A = A_1 \cup A_2$ and $B = B_1 \cup B_2$ be two superbimatrices. Suppose

$$A = A_1 \cup A_2 = [3\ 1\ 2\ |\ 0\ 1\ 5\ 1] \cup \begin{bmatrix} 3 & 0 & | & 1 \\ 2 & 1 & | & 1 \end{bmatrix}$$

and

$$B = B_1 \cup B_2 = [0\ -1\ 0\ |\ 1\ 0\ -1\ 5] \cup \begin{bmatrix} 0 & 0 & | & 1 \\ -2 & 0 & | & 5 \end{bmatrix}.$$

Then we can define biaddition of the superbimatrices A and B.

$$\begin{aligned}
A + B &= (A_1 \cup A_2) + (B_1 \cup B_2) \\
&= (A_1 + B_1) \cup (A_2 + B_2) \\
&= [3\ 1\ 2\ |\ 0\ 1\ 5\ 1] + [0\ -1\ 0\ |\ 1\ 0\ -1\ 5] \\
&\quad \cup \begin{bmatrix} 3 & 0 & | & 1 \\ 2 & 1 & | & 1 \end{bmatrix} + \begin{bmatrix} 0 & 0 & | & 1 \\ -2 & 0 & | & 5 \end{bmatrix}. \\
&= [3\ 0\ 2\ |\ 1\ 1\ 4\ 6] \cup \begin{bmatrix} 3 & 0 & | & 2 \\ 0 & 1 & | & 6 \end{bmatrix}.
\end{aligned}$$

We see both A and B are mixed superbimatrices and A + B is also a mixed superbimatrix of the same type.

Now we give yet another example.

*Example 2.12:* Let A and B be any two superbimatrices where

$$A = A_1 \cup A_2 = [3\ 1\ 1\ |\ 2] \cup \begin{bmatrix} 0 & 1 \\ 5 & 2 \end{bmatrix}$$

and

$$B = B_1 \cup B_2 = [3\ |\ 1\ 1\ 2] \cup \begin{bmatrix} 0 & | & 1 \\ 5 & | & 2 \end{bmatrix}.$$

Clearly A + B cannot be defined for $A_1$ and $B_1$ though are the supermatrices they enjoy different partitions. Similarly we see $B_2$ and $A_2$ are supermatrices yet on them are defined different partitions so addition of them cannot be defined. Thus we see unlike matrices of same order can be added; in case of supermatrices for addition to be compatible we need the matrices should be of same order and also they should have the same partition defined on them.

Now we proceed on to define addition of superbimatrices.

**DEFINITION 2.7:** *Let $A = A_1 \cup A_2$ and $B = B_1 \cup B_2$ be any two superbimatrices. For their addition $A + B$ to be defined we demand the following conditions to be satisfied.*

1. *$A_1$ and $B_1$ should be supermatrices of same order and the partition on $A_1$ and $B_1$ must be the same or identical then alone $A_1 + B_1$ is defined.*

2. *$A_2$ and $B_2$ should be supermatrices of the same order and the partitions on $A_2$ and $B_2$ must be identical then alone the sum of $A_2$ and $B_2$ can be defined.*

Thus when $A = A_1 \cup A_2$ and $B = B_1 \cup B_2$ the sum of A and B is defined to be
$$\begin{aligned} A + B &= (A_1 \cup A_2) + (B_1 \cup B_2) \\ &= (A_1 + B_1) \cup (A_2 + B_2). \end{aligned}$$

Thus only when all the above condition are satisfied we have the sum or addition of two superbimatrices to be defined and existing.

**Note:** If $A = A_1 \cup A_2$ is a superbimatrix then we have
$$\begin{aligned} A + A &= (A_1 \cup A_2) + (A_1 \cup A_2) \\ &= (A_1 + A_1) \cup (A_2 + A_2) \end{aligned}$$
is always defined and
$$\begin{aligned} A + A &= (A_1 + A_1) \cup (A_2 + A_2) \\ &= 2A_1 \cup 2A_2. \end{aligned}$$

Thus $\underbrace{A + \ldots + A}_{n-\text{times}} = nA = nA_1 \cup nA_2$.

*Example 2.13:* Let

$$A = \begin{bmatrix} 3 & | & 2 \\ 1 & | & 0 \\ 0 & | & 5 \end{bmatrix} \cup \left[ \begin{array}{cc|ccc} 1 & 1 & 3 & 0 & 2 \\ 3 & 0 & 5 & 2 & -1 \\ \hline 1 & 1 & 3 & 2 & -5 \end{array} \right] = A_1 \cup A_2$$

be a mixed rectangular superbimatrix.

$$\begin{aligned}
A + A &= (A_1 \cup A_2) + (A_1 \cup A_2) \\
&= (A_1 + A_1) \cup (A_2 + A_2) \\
&= \begin{bmatrix} 3 & | & 2 \\ 1 & | & 0 \\ 0 & | & 5 \end{bmatrix} + \begin{bmatrix} 3 & | & 2 \\ 1 & | & 0 \\ 0 & | & 5 \end{bmatrix} \\
&\quad \cup \left[ \begin{array}{cc|ccc} 1 & 1 & 3 & 0 & 2 \\ 3 & 0 & 5 & 2 & -1 \\ \hline 1 & 1 & 3 & 2 & -5 \end{array} \right] + \left[ \begin{array}{cc|ccc} 1 & 1 & 3 & 0 & 2 \\ 3 & 0 & 5 & 2 & -1 \\ \hline 1 & 1 & 3 & 2 & -5 \end{array} \right] \\
&= \begin{bmatrix} 6 & | & 4 \\ 2 & | & 0 \\ 0 & | & 10 \end{bmatrix} \cup \left[ \begin{array}{cc|ccc} 2 & 2 & 6 & 0 & 4 \\ 6 & 0 & 10 & 4 & -2 \\ \hline 2 & 2 & 6 & 4 & -10 \end{array} \right] \\
&= 2A_1 \cup 2A_2.
\end{aligned}$$

Thus $8A = 8A_1 \cup 8A_2$

$$= \begin{bmatrix} 24 & | & 16 \\ 8 & | & 0 \\ 0 & | & 40 \end{bmatrix} \cup \left[ \begin{array}{cc|ccc} 8 & 8 & 24 & 0 & 16 \\ 24 & 0 & 40 & 16 & -8 \\ \hline 8 & 8 & 24 & 16 & -40 \end{array} \right].$$

Now addition of mixed square superbimatrices, mixed column superbimatrices etc., can be defined in a similar way, provided they enjoy the same order and identical partition.

*Example 2.14:* Let

$$A = A_1 \cup A_2$$

$$= \begin{bmatrix} 3 & 5 & 1 & 2 \\ 0 & 0 & 5 & -1 \\ 2 & 0 & 1 & 2 \\ 0 & 1 & 0 & 1 \end{bmatrix} \cup \begin{bmatrix} 3 & 1 & 0 & 0 & 5 \\ 1 & 0 & 1 & 1 & 2 \\ 3 & 0 & 1 & 1 & 0 \end{bmatrix}$$

and

$$B = B_1 \cup B_2$$

$$= \begin{bmatrix} 5 & 3 & -1 & 2 \\ 1 & 1 & 1 & 1 \\ 0 & 0 & 5 & -1 \\ -1 & 2 & 0 & 3 \end{bmatrix} \cup \begin{bmatrix} 0 & 1 & 1 & 1 & -2 \\ 1 & 1 & 1 & 0 & 3 \\ 2 & 0 & 1 & 0 & 1 \end{bmatrix}$$

be two mixed superbimatrices.

$$A + B = (A_1 \cup A_2) + (B_1 \cup B_2)$$

$$= (A_1 + B_1) \cup (A_2 + B_2)$$

$$= \begin{bmatrix} 3 & 5 & 1 & 2 \\ 0 & 0 & 5 & -1 \\ 2 & 0 & 1 & 2 \\ 0 & 1 & 0 & 1 \end{bmatrix} + \begin{bmatrix} 5 & 3 & -1 & 2 \\ 1 & 1 & 1 & 1 \\ 0 & 0 & 5 & -1 \\ -1 & 2 & 0 & 3 \end{bmatrix} \cup$$

$$\begin{bmatrix} 3 & 1 & 0 & 0 & 5 \\ 1 & 0 & 1 & 1 & 2 \\ 3 & 0 & 1 & 1 & 0 \end{bmatrix} + \begin{bmatrix} 0 & 1 & 1 & 1 & -2 \\ 1 & 1 & 1 & 0 & 3 \\ 2 & 0 & 1 & 0 & 1 \end{bmatrix}$$

$$= \begin{bmatrix} 8 & | & 8 & 0 & 4 \\ 1 & | & 1 & 6 & 0 \\ 2 & | & 0 & 6 & 1 \\ -1 & | & 3 & 0 & 4 \end{bmatrix} \cup \begin{bmatrix} 3 & 2 & | & 1 & 1 & 3 \\ 2 & 1 & | & 2 & 1 & 5 \\ 5 & 0 & | & 2 & 1 & 1 \end{bmatrix}.$$

Now we have to define the transpose of a superbimatrix $A = A_1 \cup A_2$.

**DEFINITION 2.8:** *Let $A = A_1 \cup A_2$ be any superbimatrix. The transpose of the superbimatrix A denoted by $A^T$ is defined to be $A^T = (A_1 \cup A_2)^T = A_1^T \cup A_2^T$. Clearly the transpose of a superbimatrix is again a superbimatrix. If A is a mixed rectangular superbimatrix then its transpose, $A^T$ will also be a mixed rectangular superbimatrix.*

It has become pertinent to mention here that however a column superbimatrix transpose would be a row superbimatrix and a row superbimatrix transpose would be a column superbimatrix.

Now we proceed onto illustrate them with examples.

*Example 2.15:* Let $A = A_1 \cup A_2 = [3\ 0\ 1\ 1\ |\ -1\ 5\ 2\ 3\ 1] \cup [1\ 0\ 1\ |\ 5\ 2\ 0\ |\ 1\ 1\ 1\ 0\ 2]$ be a row superbimatrix. The transpose of A denoted by

$$A^T = \begin{bmatrix} 3 \\ 0 \\ 1 \\ 1 \\ \overline{-1} \\ 5 \\ 2 \\ 3 \\ 1 \end{bmatrix} \cup \begin{bmatrix} 1 \\ 0 \\ 1 \\ \overline{5} \\ 2 \\ 0 \\ \overline{1} \\ 1 \\ 1 \\ 0 \\ 2 \end{bmatrix} = A_1^T \cup A_2^T.$$

Clearly $A^T$ is a column superbimatrix.

*Example 2.16:* Let

$$B = B_1 \cup B_2 = \begin{bmatrix} 3 \\ 1 \\ 1 \\ \frac{1}{2} \\ 2 \\ \overline{1} \\ 1 \\ 0 \end{bmatrix} \cup \begin{bmatrix} \frac{1}{2} \\ -1 \\ 0 \\ 1 \\ 1 \\ \overline{3} \\ 5 \\ 1 \\ 2 \\ 6 \end{bmatrix}$$

be a column superbimatrix. The transpose of B which is $B^T$ = [3 1 1 1 | 2 2 | 1 1 | 0] $\cup$ [1 | 2 –1 0 1 1 | 3 5 1 2 6] = $B_1^T \cup B_2^T$. Clearly $B^T$ is a row superbimatrix.

*Example 2.17:* Let $C = C_1 \cup C_2$ be any superbimatrix where

$$C_1 = \begin{bmatrix} 3 & 1 & 0 & 2 & 1 \\ 6 & 2 & 1 & 0 & 5 \\ -1 & 0 & 1 & -1 & 6 \end{bmatrix}$$

and

$$C_2 = \begin{bmatrix} 3 & 1 & 0 \\ 1 & 1 & 6 \\ 2 & 1 & 5 \\ -1 & 0 & 1 \\ 2 & 3 & 2 \\ 1 & 0 & 5 \end{bmatrix}.$$

C is a mixed rectangular superbimatrix.

$$C^T = (C_1 \cup C_2)^T = C_1^T \cup C_2^T$$

$$= \begin{bmatrix} 3 & 6 & -1 \\ 1 & 2 & 0 \\ \hline 0 & 1 & 1 \\ 2 & 0 & -1 \\ 1 & 5 & 6 \end{bmatrix} \cup \begin{bmatrix} 3 & 1 & 2 & -1 & 2 & 1 \\ 1 & 1 & 1 & 0 & 3 & 0 \\ 0 & 6 & 5 & 1 & 2 & 5 \end{bmatrix}.$$

$C^T$ is also a mixed rectangular superbimatrix.

***Example 2.18:*** Let $D = D_1 \cup D_2$ where D is a mixed square superbimatrix with

$$D_1 = \begin{bmatrix} 3 & 1 & 1 \\ 0 & 5 & 2 \\ \hline 1 & 0 & 1 \end{bmatrix}$$

and

$$D_2 = \begin{bmatrix} 3 & 0 & -1 & 2 \\ 1 & 2 & 0 & 1 \\ \hline 0 & 1 & 1 & 0 \\ -1 & 3 & 0 & 1 \end{bmatrix}.$$

Now

$$D^T = (D_1 \cup D_2)^T = D_1^T \cup D_2^T$$

$$= \begin{bmatrix} 3 & 0 & 1 \\ 1 & 5 & 0 \\ \hline 1 & 2 & 1 \end{bmatrix} \cup \begin{bmatrix} 3 & 1 & 0 & -1 \\ 0 & 2 & 1 & 3 \\ \hline -1 & 0 & 1 & 0 \\ 2 & 1 & 0 & 1 \end{bmatrix}.$$

We see $D^T$ is also a square mixed superbimatrix.

Now we are interested to know what to define or call the following type of bimatrices.

*Example 2.19:* Let $A = A_1 \cup A_2$ where

$$A_1 = \begin{bmatrix} 3 & 1 & 2 & 5 & 6 \\ 0 & 2 & 0 & 1 & 0 \\ 1 & 1 & 5 & 3 & -1 \end{bmatrix}$$

and

$$A_2 = \left[\begin{array}{c|ccc} 3 & 1 & 1 & 0 \\ 1 & 0 & 1 & 1 \\ \hline -1 & -1 & 0 & 0 \\ 0 & 2 & -2 & 1 \end{array}\right]$$

where $A_1$ is just a $3 \times 5$ matrix and $A_2$ is a square supermatrix. Since $A_1$ is not a supermatrix we cannot define $A$ to be a superbimatrix since $A_2$ is a supermatrix we cannot define $A$ to be a bimatrix.

So we define a new notion called semi superbimatrix in such cases.

**DEFINITION 2.9:** *Let $A = A_1 \cup A_2$ where $A_1$ is just a simple matrix and $A_2$ is a supermatrix then we define $A = A_1 \cup A_2$ to be a semi superbimatrix.*

*Example 2.20:* Let

$$A = A_1 \cup A_2 = \begin{bmatrix} 3 & 1 & 1 & 2 \\ 0 & 5 & 1 & 0 \end{bmatrix} \cup \left[\begin{array}{cc|cc} 3 & 1 & 2 & 0 \\ 5 & 1 & 1 & 1 \\ \hline 2 & 0 & 2 & 6 \\ 1 & 0 & 1 & 5 \end{array}\right]$$

where $A_1$ is just a $2 \times 4$ matrix and $A_2$ is a square supermatrix, $A$ is defined as the semi superbimatrix.

***Example 2.21:*** Let $A = A_1 \cup A_2$ where $A_1 = [3\ 2\ 3\ |\ 1\ 0\ 0\ 5]$ and $A_2 = [1\ 1\ 0\ 1\ 1\ 1]$. We see A is a semi superbimatrix called as a row semi superbimatrix.

***Example 2.22:*** Let $A = A_1 \cup A_2$ where

$$A_1 = \begin{bmatrix} 1 \\ 0 \\ 1 \\ \hline 2 \\ 2 \\ -1 \\ 6 \end{bmatrix}$$

and

$$A_2 = \begin{bmatrix} 3 \\ 1 \\ 4 \\ 5 \end{bmatrix}.$$

A is a semi superbimatrix which we define as a column semi superbimatrix.

***Example 2.23:*** Let $A = A_1 \cup A_2$ where

$$A_1 = \begin{bmatrix} 3 & 2 & 1 & 3 & 5 & 3 \\ 1 & 0 & 0 & -1 & 2 & 1 \end{bmatrix}$$

and

$$A_2 = \begin{bmatrix} 1 & 2 & 3 & 5 & 3 & 7 & 0 & 1 \\ 7 & 6 & 3 & 5 & 1 & 2 & -1 & 1 \\ 1 & 0 & 1 & 0 & 1 & 5 & 2 & 3 \end{bmatrix}.$$

A is a mixed rectangular semi superbimatrix.

*Example 2.24:* Let $A = A_1 \cup A_2$ where

$$A_1 = \begin{bmatrix} 3 & 1 & 2 & 0 \\ 1 & 0 & 1 & 0 \\ 3 & 1 & 2 & 1 \\ 5 & 0 & 1 & 1 \end{bmatrix}$$

and

$$A_2 = \begin{bmatrix} 3 & 1 & 2 \\ 0 & 0 & 1 \\ \hline 1 & 0 & 1 \end{bmatrix}.$$

A is a mixed square semi superbimatrix.

*Example 2.25:* Let $A = A_1 \cup A_2$ where

$$A_1 = \begin{bmatrix} 3 & 0 & 1 & 0 \\ 1 & 1 & -1 & 1 \\ 2 & 1 & 0 & 1 \\ 5 & -1 & -1 & 0 \end{bmatrix}$$

and

$$A_2 = \begin{bmatrix} 8 & 1 & 2 \\ 0 & 1 & 1 \\ \hline 1 & 0 & 1 \\ 1 & 1 & 0 \\ 1 & 0 & 0 \\ 0 & 1 & 0 \\ 2 & 1 & 8 \end{bmatrix}$$

A is a mixed semi superbimatrix for $A_1$ is just a square matrix and $A_2$ is a super rectangular $7 \times 3$ matrix. Thus we can define 7 types of semi superbimatrices viz. row semi superbimatrix, column semi superbimatrix, $n \times n$ square semi superbimatrix, m

× n rectangular semi superbimatrix, mixed square semi superbimatrix, mixed rectangular semi superbimatrix and mixed semi superbimatrix.

*Example 2.26:* Let

$$A = \begin{bmatrix} 3 & 1 & 1 \\ 1 & 0 & 1 \end{bmatrix} \cup \left[ \begin{array}{cc|c} 0 & 1 & 2 \\ 3 & 4 & 5 \end{array} \right]$$

$$= A_1 \cup A_2.$$

A is a 2 × 3 rectangular semi superbimatrix.

*Example 2.27:* Let

$$A = \begin{bmatrix} 1 & 0 & 3 & -1 \\ 1 & 0 & 1 & 1 \\ 1 & 0 & 5 & 2 \\ 1 & 2 & 0 & 3 \end{bmatrix} \cup \left[ \begin{array}{c|ccc} 0 & -1 & 2 & 3 \\ \hline 1 & 2 & 3 & 4 \\ 5 & 6 & 7 & 8 \\ 4 & 3 & 2 & 1 \end{array} \right]$$

$$= A_1 \cup A_2.$$

A is a 4 × 4 square semi superbimatrix.

We see as in case of superbimatrices the transpose of a semi superbimatrix is also a semi superbimatrix. The transpose of a row semi superbimatrix is a column semi superbimatrix and the transpose of a column semi superbimatrix is a row semi superbimatrix.

*Example 2.28:* Let

A  =  $A_1 \cup A_2$
   =  [3 1 0 0 0 3 5] ∪ [1 1 2 3 | 5 5 4 | 4 1 2]

be a row semi superbimatrix.
  $A^T$ =  $(A_1 \cup A_2)^T$

$$A_1^T \cup A_2^T = \begin{bmatrix} 3 \\ 1 \\ 0 \\ 0 \\ 0 \\ 3 \\ 5 \end{bmatrix} \cup \begin{bmatrix} 1 \\ 1 \\ 2 \\ \dfrac{3}{5} \\ 5 \\ \dfrac{4}{4} \\ 1 \\ 2 \end{bmatrix}$$

is a column semi superbimatrix.

*Example 2.29:* Let $A = A_1 \cup A_2$ where

$$A_1 = \begin{bmatrix} 3 \\ 1 \\ 1 \\ 1 \\ 0 \\ 2 \end{bmatrix} \text{ and } A_2 = \begin{bmatrix} 1 \\ 2 \\ \dfrac{3}{4} \\ 5 \\ 6 \\ 7 \end{bmatrix}.$$

$$\begin{aligned} A^T &= (A_1 \cup A_2)^T = A_1^T \cup A_2^T \\ &= [3\ 1\ 1\ 1\ 0\ 2] \cup [1\ 2\ 3\ |\ 4\ 5\ 6\ 7]. \end{aligned}$$

A is clearly a column semi superbimatrix but $A^T$ is a row semi superbimatrix.

*Example 2.30:* Let $A = A_1 \cup A_2$ be a mixed square semi superbimatrix where

$$A_1 = \begin{bmatrix} 3 & 0 & 1 \\ 1 & 1 & 1 \\ 0 & 1 & 0 \end{bmatrix}$$

and

$$A_2 = \begin{bmatrix} 2 & 1 & 0 & 5 \\ \hline 1 & 1 & 2 & 8 \\ 0 & 2 & 1 & 3 \\ 5 & 8 & 3 & 0 \end{bmatrix}.$$

$$A^T = (A_1 \cup A_2)^T = A_1^T \cup A_2^T$$

$$= \begin{bmatrix} 3 & 1 & 0 \\ 0 & 1 & 1 \\ 1 & 1 & 0 \end{bmatrix} \cup \begin{bmatrix} 2 & 1 & 0 & 5 \\ \hline 1 & 1 & 2 & 8 \\ \hline 0 & 2 & 1 & 3 \\ 5 & 8 & 3 & 0 \end{bmatrix}.$$

Clearly $A^T$ is also a mixed square semi superbimatrix.

*Example 2.31:* Let $A = A_1 \cup A_2$ be a mixed rectangular semi superbimatrix where

$$A_1 = \begin{bmatrix} 3 & 1 & 2 \\ 0 & 1 & 1 \\ 1 & 1 & 0 \\ 2 & 2 & 0 \\ 2 & 1 & 0 \\ 1 & 0 & 3 \end{bmatrix}$$

and

$$A_2 = \begin{bmatrix} 1 & 3 & 5 & 8 & 1 & 3 & 0 \\ 2 & 1 & 6 & 9 & 2 & 1 & 1 \end{bmatrix}.$$

The transpose of A

$$A^T = (A_1 \cup A_2)^T = A_1^T \cup A_2^T$$

where

$$A_1^T = \begin{bmatrix} 3 & 0 & 1 & 2 & 2 & 1 \\ 1 & 1 & 1 & 2 & 1 & 0 \\ 2 & 1 & 0 & 0 & 0 & 3 \end{bmatrix}$$

$$A_2^T = \begin{bmatrix} 1 & 2 \\ 3 & 1 \\ 5 & 6 \\ \hline 8 & 9 \\ 1 & 2 \\ 3 & 1 \\ 0 & 1 \end{bmatrix}.$$

$A^T$ is also a mixed rectangular semi superbimatrix.

***Example 2.32:*** Let $A = A_1 \cup A_2$ be a $7 \times 3$ rectangular semi superbimatrix where

$$A_1 = \begin{bmatrix} 1 & 0 & 1 \\ 0 & 1 & 1 \\ 1 & 0 & 0 \\ 0 & 0 & 1 \\ 0 & 1 & 0 \\ 1 & 1 & 1 \\ 5 & 7 & 8 \end{bmatrix}$$

and

$$A_2 = \begin{bmatrix} 2 & 1 & 2 \\ 1 & 0 & 1 \\ -1 & 1 & 0 \\ 1 & 1 & 5 \\ \hline 2 & 3 & 5 \\ 1 & 2 & 3 \\ 1 & 0 & 0 \end{bmatrix}.$$

$$A^T = A_1^T \cup A_2^T = \begin{bmatrix} 1 & 0 & 1 & 0 & 0 & 1 & 5 \\ 0 & 1 & 0 & 0 & 1 & 1 & 7 \\ 1 & 1 & 0 & 1 & 0 & 1 & 8 \end{bmatrix} \cup$$

$$\begin{bmatrix} 2 & 1 & -1 & 1 & 2 & 1 & 1 \\ 1 & 0 & 0 & 1 & 3 & 2 & 0 \\ 2 & 1 & 1 & 5 & 5 & 3 & 0 \end{bmatrix}$$

is also a 3 × 7 rectangular semi superbimatrix.

As in case of superbimatrices the sum of two semi superbimatrices can be added i.e., $A = A_1 \cup A_2$ and $B = B_1 \cup B_2$ can be added if and only if both $A_1$ and $B_1$ and same order matrices and $A_2$ and $B_2$ are same order supermatrices with same or identical partition on it.

*Example 2.33:* Let $A = A_1 \cup A_2$ and $B = B_1 \cup B_2$ be any two semi superbimatrices. Here

$$A_1 = \begin{bmatrix} 0 & 1 & 2 \\ 3 & 4 & 5 \\ 6 & 7 & 8 \end{bmatrix}$$

and

$$A_2 = \begin{bmatrix} 3 & 1 & 0 & 2 & 1 \\ 1 & 1 & 5 & 0 & 2 \\ \hline -1 & 3 & 1 & 0 & 1 \end{bmatrix}$$

and

$$B_1 = \begin{bmatrix} 0 & 1 & 0 \\ 1 & 1 & 1 \\ 2 & 1 & 0 \end{bmatrix}$$

and

$$B_2 = \begin{bmatrix} 1 & 1 & 0 & 1 & 0 \\ 0 & 1 & 2 & 3 & 5 \\ \hline 6 & 7 & 8 & 9 & 1 \end{bmatrix}.$$

$$\begin{aligned} A + B &= (A_1 \cup A_2) + (B_1 \cup B_2) \\ &= (A_1 + B_1) \cup (A_2 + B_2) \end{aligned}$$

$$= \begin{bmatrix} 0 & 1 & 2 \\ 3 & 4 & 5 \\ 6 & 7 & 8 \end{bmatrix} + \begin{bmatrix} 0 & 1 & 0 \\ 1 & 1 & 1 \\ 2 & 1 & 0 \end{bmatrix} \cup$$

$$\left[ \begin{array}{cc|ccc} 3 & 1 & 0 & 2 & 1 \\ 1 & 1 & 5 & 0 & 2 \\ \hline -1 & 3 & 1 & 0 & 1 \end{array} \right] + \left[ \begin{array}{cc|ccc} 1 & 1 & 0 & 1 & 0 \\ 0 & 1 & 2 & 3 & 5 \\ \hline 6 & 7 & 8 & 9 & 1 \end{array} \right]$$

$$= \begin{bmatrix} 0 & 2 & 2 \\ 4 & 5 & 6 \\ 8 & 8 & 8 \end{bmatrix} \cup \left[ \begin{array}{cc|ccc} 4 & 2 & 0 & 3 & 1 \\ 1 & 2 & 7 & 3 & 7 \\ \hline 5 & 10 & 9 & 9 & 2 \end{array} \right].$$

A+B is also a semi superbimatrix of the same type.

Suppose we have

$$A = A_1 \cup A_2 = \left[ \begin{array}{cc|c} 0 & 1 & 2 \\ 0 & 1 & -1 \\ 1 & 1 & 0 \end{array} \right] \cup \left[ \begin{array}{ccc|ccc} 3 & 1 & -1 & 1 & 1 & 0 \\ 4 & 2 & 0 & 1 & 0 & -1 \\ 1 & 0 & 1 & 0 & -1 & 2 \end{array} \right]$$

to be a mixed superbimatrix then we have always A + A is defined and is 2A. In fact A + A + … + A, n-times in nA. We see

$$A + A = \left[ \begin{array}{cc|c} 0 & 1 & 2 \\ 0 & 1 & -1 \\ 1 & 1 & 0 \end{array} \right] + \left[ \begin{array}{cc|c} 0 & 1 & 2 \\ 0 & 1 & -1 \\ 1 & 1 & 0 \end{array} \right] \cup$$

$$\left[ \begin{array}{ccc|ccc} 3 & 1 & -1 & 1 & 1 & 0 \\ 4 & 2 & 0 & 1 & 0 & -1 \\ 1 & 0 & 1 & 0 & -1 & 2 \end{array} \right] + \left[ \begin{array}{ccc|ccc} 3 & 1 & -1 & 1 & 1 & 0 \\ 4 & 2 & 0 & 1 & 0 & -1 \\ 1 & 0 & 1 & 0 & -1 & 2 \end{array} \right]$$

$$= \begin{bmatrix} 0 & 2 & 4 \\ 0 & 2 & -2 \\ 2 & 2 & 0 \end{bmatrix} \cup \begin{bmatrix} 6 & 2 & -2 & 2 & 2 & 0 \\ 8 & 4 & 0 & 2 & 0 & -2 \\ 2 & 0 & 2 & 0 & -2 & 4 \end{bmatrix} = 2A_1 \cup 2A_2.$$

Now $nA = nA_1 \cup nA_2$ for any $n > 1$. Also

$$\begin{aligned} A - A &= (A_1 \cup A_2) - (A_1 \cup A_2) \\ &= (A_1 - A_1) \cup (A_2 - A_2) \\ &= \begin{bmatrix} 0 & 0 & 0 \\ 0 & 0 & 0 \\ 0 & 0 & 0 \end{bmatrix} \cup \begin{bmatrix} 0 & 0 & 0 & 0 & 0 & 0 \\ 0 & 0 & 0 & 0 & 0 & 0 \\ 0 & 0 & 0 & 0 & 0 & 0 \end{bmatrix}. \end{aligned}$$

Thus we get the difference of $A - A$ to be a zero superbimatrix.
Now if

$$A = A_1 \cup A_2$$

$$= \begin{bmatrix} 1 & 2 & 3 \\ 0 & 1 & 0 \\ 1 & 0 & 1 \end{bmatrix} \cup \begin{bmatrix} 3 & 2 & 1 & 0 & 5 & 2 & 3 & 1 \\ 1 & 2 & 0 & 1 & -1 & 3 & 1 & 1 \end{bmatrix}$$

be a mixed semi superbimatrix. Then

$$\begin{aligned} A + A &= (A_1 \cup A_2) + (A_1 \cup A_2) \\ &= (A_1 + A_1) \cup (A_2 + A_2) \\ &= \begin{bmatrix} 1 & 2 & 3 \\ 0 & 1 & 0 \\ 1 & 0 & 1 \end{bmatrix} + \begin{bmatrix} 1 & 2 & 3 \\ 0 & 1 & 0 \\ 1 & 0 & 1 \end{bmatrix} \cup \\ &\quad \begin{bmatrix} 3 & 2 & 1 & 0 & 5 & 2 & 3 & 1 \\ 1 & 2 & 0 & 1 & -1 & 3 & 1 & 1 \end{bmatrix} + \end{aligned}$$

$$\begin{bmatrix} 3 & 2 & 1 & 0 & 5 & 2 & 3 & 1 \\ 1 & 2 & 0 & 1 & -1 & 3 & 1 & 1 \end{bmatrix}$$

$$= \begin{bmatrix} 2 & 4 & | & 6 \\ 0 & 2 & | & 0 \\ \hline 2 & 0 & | & 2 \end{bmatrix} \cup \begin{bmatrix} 6 & 4 & 2 & 0 & 10 & 4 & 6 & 2 \\ 2 & 4 & 0 & 2 & -2 & 6 & 2 & 2 \end{bmatrix}$$

$$= 2A_1 \cup 2A_2 = 2A.$$

Thus we see sum of a semi superbimatrix A with itself is 2A. On similar lines we can say if A is a semi superbimatrix then A + A + … + A, n-times is $nA = nA_1 \cup nA_2$.

Having defined transpose and sum of semi superbimatrix whenever it is defined we proceed onto define some product of these superbimatrices and semi superbimatrices.

As in case of supermatrices we in case of superbimatrices and semi superbimatrices first define the notion of the product of a superbimatrix with its transpose. We also for this need the simple definition of symmetric superbimatrices, symmetric semi superbimatrices, quasi symmetric superbimatrices and quasi symmetric semi superbimatrices.

**DEFINITION 2.10:** *Let $A = A_1 \cup A_2$ be a superbimatrix. We call A to be a symmetric superbimatrix if both $A_1$ and $A_2$ are symmetric supermatrices.*

*Example 2.34:* Let $A = A_1 \cup A_2$ be a superbimatrix where

$$A_1 = \begin{bmatrix} 0 & 1 & 2 & | & 3 & 4 \\ 1 & 2 & 1 & | & 2 & 3 \\ 2 & 1 & 4 & | & 2 & 5 \\ \hline 3 & 2 & 2 & | & 3 & 1 \\ 4 & 3 & 5 & | & 1 & 1 \end{bmatrix}$$

and

$$A_2 = \begin{bmatrix} 6 & 1 & 2 & 0 \\ 1 & 3 & 1 & 4 \\ 2 & 1 & 5 & 1 \\ 0 & 4 & 1 & 7 \end{bmatrix},$$

since both $A_1$ and $A_2$ are symmetric supermatrices we see $A = A_1 \cup A_2$ is a symmetric superbimatrix.

**DEFINITION 2.11:** *Let $A = A_1 \cup A_2$ be a superbimatrix, A is said to be quasi symmetric superbimatrix if and only if one of $A_1$ or $A_2$ is a symmetric supermatrix.*

The following result is obvious every symmetric superbimatrix is trivially a quasi symmetric superbimatrix. However a quasi symmetric superbimatrix is never a symmetric superbimatrix.

*Example 2.35:* Let $A = A_1 \cup A_2$ where

$$A_1 = \begin{bmatrix} 3 & 2 & 1 & 0 & 5 \\ 2 & 1 & 3 & 2 & 1 \\ 1 & 3 & 5 & 1 & 2 \\ 0 & 2 & 1 & 4 & 3 \\ 5 & 1 & 2 & 3 & 2 \end{bmatrix}$$

and

$$A_2 = \begin{bmatrix} 4 & 1 & 2 & 3 \\ 1 & 2 & 3 & 4 \\ 2 & 3 & 4 & 1 \\ 3 & 4 & 1 & 2 \end{bmatrix}$$

be a superbimatrix. Clearly $A_1$ is a symmetric supermatrix where as $A_2$ is only a square supermatrix. Hence A is only a quasi symmetric superbimatrix.

**DEFINITION 2.12:** *Let $A = A_1 \cup A_2$ be a semi superbimatrix. If both $A_1$ and $A_2$ are symmetric matrices and $A_1$ or $A_2$ is a super* *symmetric matrix then we all A to be a symmetric semi superbimatrix 'or' used in the mutually exclusive sense.*

**Example 2.36:** Let $A = A_1 \cup A_2$ where

$$A_1 = \begin{bmatrix} 5 & 4 & | & 1 & 2 & 3 \\ 4 & 1 & | & 2 & 3 & 4 \\ \hline 1 & 2 & | & 4 & 1 & 2 \\ 2 & 3 & | & 1 & 3 & 1 \\ 3 & 4 & | & 2 & 1 & 2 \end{bmatrix}$$

and

$$A_2 = \begin{bmatrix} 0 & 1 & 0 & 1 \\ 1 & 4 & 3 & 1 \\ 0 & 3 & 2 & 1 \\ 1 & 1 & 1 & 5 \end{bmatrix}.$$

A is a symmetric semi superbimatrix as $A_1$ is a symmetric supermatrix and $A_2$ is a symmetric matrix.

**DEFINITION 2.13:** *Let $A = A_1 \cup A_2$ be a semi superbimatrix. If only one of $A_1$ or $A_2$ is a symmetric matrix (or a symmetric supermatrix) then we call A to be a quasi symmetric semi superbimatrix.*

**Example 2.37:** Let $A = A_1 \cup A_2$ be a semi superbimatrix where

$$A_1 = \begin{bmatrix} 2 & 1 & 3 & 0 & 5 \\ 1 & 5 & 1 & 6 & 2 \\ 3 & 1 & 7 & 1 & 2 \\ 0 & 6 & 1 & 3 & 1 \\ 5 & 2 & 2 & 1 & 4 \end{bmatrix}$$

and

$$A_2 = \begin{bmatrix} 3 & 1 & 2 & 5 & 6 & 7 \\ 1 & 2 & 3 & 4 & 5 & 6 \\ \hline 2 & 3 & 4 & 5 & 6 & 1 \\ 3 & 4 & 5 & 6 & 1 & 2 \\ 4 & 5 & 6 & 1 & 2 & 3 \\ 5 & 6 & 1 & 2 & 3 & 4 \end{bmatrix}$$

is a quasi symmetric semi superbimatrix.

*Example 2.38:* Let $A = A_1 \cup A_2$ where

$$A_1 = \begin{bmatrix} 0 & 1 & 0 & 1 \\ 1 & 2 & 3 & 4 \\ 1 & 0 & 1 & 0 \\ 4 & 3 & 2 & 1 \end{bmatrix}$$

and

$$A_2 = \begin{bmatrix} 1 & 2 & 3 & 4 & 5 \\ 2 & 5 & 4 & 3 & 2 \\ 3 & 4 & 2 & 4 & 3 \\ 4 & 3 & 4 & 3 & 1 \\ \hline 5 & 2 & 3 & 1 & 6 \end{bmatrix}.$$

A is a semi superbimatrix which is also a quasi symmetric semi superbimatrix.
   We show later the byproduct which we define on superbimatrices we get a class of symmetric bimatrices.

   We first illustrate the product of two superbimatrices.

*Example 2.39:* Let $A = A_1 \cup A_2$ and $B = B_1 \cup B_2$ be any two superbimatrices. Here

$$A_1 = \begin{bmatrix} 2 & 0 & | & 3 & 0 & 1 & | & 4 \\ 1 & 1 & | & 1 & 1 & 0 & | & 1 \\ 1 & 2 & | & 0 & 1 & 1 & | & 0 \end{bmatrix}$$

and

$$A_2 = \begin{bmatrix} 3 & 1 & 0 & | & 3 & 3 & 0 & 1 \\ 4 & 5 & 1 & | & 1 & 0 & 1 & 1 \\ 3 & 4 & 1 & | & 0 & 1 & 0 & 1 \\ 1 & 2 & 2 & | & 4 & 2 & 5 & 6 \end{bmatrix}.$$

$$B_1 = \begin{bmatrix} 0 & 1 \\ 3 & 0 \\ 1 & 0 \\ \hline 1 & 1 \\ 2 & 0 \\ \hline 0 & 1 \end{bmatrix}$$

and

$$B_2 = \begin{bmatrix} 1 & 0 & 3 \\ 3 & 1 & 1 \\ 5 & 1 & 2 \\ \hline 1 & 1 & 0 \\ 0 & 1 & 1 \\ 1 & 0 & 0 \\ 0 & 1 & 0 \end{bmatrix}.$$

A and B are row superbimatrix and column superbimatrix respectively.

$$\begin{aligned} AB &= [A_1 \cup A_2][B_1 \cup B_2] \\ &= A_1B_1 \cup A_2B_2 \end{aligned}$$

where the product AB is defined as the minor product (refer 37-40 of chapter one).

$$AB = \begin{bmatrix} 2 & 0 & | & 3 & 0 & 1 & | & 4 \\ 1 & 1 & | & 1 & 1 & 0 & | & 1 \\ 1 & 2 & | & 0 & 1 & 1 & | & 0 \end{bmatrix} \begin{bmatrix} 0 & 1 \\ 3 & 0 \\ \hline 1 & 0 \\ 1 & 1 \\ 2 & 0 \\ \hline 0 & 1 \end{bmatrix} \cup$$

$$\begin{bmatrix} 3 & 1 & 0 & | & 3 & 3 & 0 & 1 \\ 4 & 5 & 1 & | & 1 & 0 & 1 & 1 \\ 3 & 4 & 1 & | & 0 & 1 & 0 & 1 \\ 1 & 2 & 2 & | & 4 & 2 & 5 & 6 \end{bmatrix} \begin{bmatrix} 1 & 0 & 3 \\ 3 & 1 & 1 \\ 5 & 1 & 2 \\ \hline 1 & 1 & 0 \\ 0 & 1 & 1 \\ 1 & 0 & 0 \\ 0 & 1 & 0 \end{bmatrix}$$

$$= \left\{ \begin{bmatrix} 2 & 0 \\ 1 & 1 \\ 1 & 2 \end{bmatrix} \begin{bmatrix} 0 & 1 \\ 3 & 0 \end{bmatrix} + \begin{bmatrix} 3 & 0 & 1 \\ 1 & 1 & 0 \\ 0 & 1 & 1 \end{bmatrix} \begin{bmatrix} 1 & 0 \\ 1 & 1 \\ 2 & 0 \end{bmatrix} + \begin{bmatrix} 4 \\ 1 \\ 0 \end{bmatrix} \begin{bmatrix} 0 & 1 \end{bmatrix} \right\} \cup$$

$$\left\{ \begin{bmatrix} 3 & 1 & 0 \\ 4 & 5 & 1 \\ 3 & 4 & 1 \\ 1 & 2 & 2 \end{bmatrix} \begin{bmatrix} 1 & 0 & 3 \\ 3 & 1 & 1 \\ 5 & 1 & 2 \end{bmatrix} + \begin{bmatrix} 3 & 3 & 0 & 1 \\ 1 & 0 & 1 & 1 \\ 0 & 1 & 0 & 1 \\ 4 & 2 & 5 & 6 \end{bmatrix} \begin{bmatrix} 1 & 1 & 0 \\ 0 & 1 & 1 \\ 1 & 0 & 0 \\ 0 & 1 & 0 \end{bmatrix} \right\}$$

$$= \left\{ \begin{bmatrix} 0 & 2 \\ 3 & 1 \\ 6 & 1 \end{bmatrix} + \begin{bmatrix} 5 & 0 \\ 2 & 1 \\ 3 & 1 \end{bmatrix} + \begin{bmatrix} 0 & 4 \\ 0 & 1 \\ 0 & 0 \end{bmatrix} \right\} \cup$$

$$\left\{ \begin{bmatrix} 6 & 1 & 10 \\ 24 & 6 & 19 \\ 20 & 5 & 15 \\ 17 & 4 & 9 \end{bmatrix} + \begin{bmatrix} 3 & 7 & 3 \\ 2 & 2 & 0 \\ 0 & 2 & 1 \\ 9 & 12 & 2 \end{bmatrix} \right\} = \begin{bmatrix} 5 & 6 \\ 5 & 3 \\ 9 & 2 \end{bmatrix} \cup \begin{bmatrix} 9 & 8 & 13 \\ 26 & 8 & 19 \\ 20 & 7 & 16 \\ 26 & 16 & 11 \end{bmatrix}.$$

is a bimatrix and is clearly not a superbimatrix. Thus this sort of product leads only to a bimatrix.

We give yet another example of the same type before we proceed onto define more complicated products.

*Example 2.40:* Let $A = A_1 \cup A_2$ and $B = B_1 \cup B_2$ be any two superbimatrices where

$$A_1 = \begin{bmatrix} 6 & 1 & 0 & 3 & 2 & 0 & 1 & -1 & 5 \\ 1 & 1 & 1 & 0 & 0 & 0 & 1 & 0 & 2 \end{bmatrix}$$

and

$$A_2 = \begin{bmatrix} 3 & 1 & 1 & 1 & 1 & 0 & 1 & 1 \\ 6 & 2 & 0 & 1 & 0 & 1 & 0 & 1 \\ 1 & 0 & 0 & 0 & 1 & 1 & 1 & 0 \end{bmatrix}.$$

$$B_1 = \begin{bmatrix} 3 & 1 & 0 & 1 \\ 1 & 0 & 1 & 0 \\ -1 & 1 & 0 & 1 \\ 2 & 0 & 0 & 0 \\ 0 & 0 & 1 & 1 \\ 1 & 0 & 0 & 1 \\ 0 & 0 & 1 & 0 \\ 2 & 1 & 2 & 0 \\ 1 & 0 & 1 & 1 \end{bmatrix}$$

and

$$B_2 = \begin{bmatrix} 1 & 0 \\ 0 & 1 \\ \hline 2 & 1 \\ 1 & 2 \\ 0 & 1 \\ 1 & 1 \\ 0 & 1 \\ \hline 1 & 1 \end{bmatrix}.$$

A is a row superbimatrix and B is a column superbimatrix. Now

$$
\begin{aligned}
A.B &= (A_1 \cup A_2) . (B_1 \cup B_2) \\
&= A_1 B_1 \cup A_2 B_2 \\
&= \begin{bmatrix} 6 & 1 & 0 & 3 & | & 2 & 0 & 1 & | & -1 & 5 \\ 1 & 1 & 1 & 0 & | & 0 & 0 & 1 & | & 0 & 2 \end{bmatrix} \begin{bmatrix} 3 & 1 & 0 & 1 \\ 1 & 0 & 1 & 0 \\ -1 & 1 & 0 & 1 \\ 2 & 0 & 0 & 0 \\ 0 & 0 & 1 & 1 \\ 1 & 0 & 0 & 1 \\ 0 & 0 & 1 & 0 \\ \hline 2 & 1 & 2 & 0 \\ 1 & 0 & 1 & 1 \end{bmatrix} \cup \\
&\quad \begin{bmatrix} 3 & 1 & | & 1 & 1 & 1 & 0 & 1 & | & 1 \\ 6 & 2 & | & 0 & 1 & 0 & 1 & 0 & | & 1 \\ 1 & 0 & | & 0 & 0 & 1 & 1 & 1 & | & 0 \end{bmatrix} \begin{bmatrix} 1 & 0 \\ 0 & 1 \\ \hline 2 & 1 \\ 1 & 2 \\ 0 & 1 \\ 1 & 1 \\ 0 & 1 \\ \hline 1 & 1 \end{bmatrix}
\end{aligned}
$$

$$= \left\{ \begin{bmatrix} 6 & 1 & 0 & 3 \\ 1 & 1 & 1 & 0 \end{bmatrix} \begin{bmatrix} 3 & 1 & 0 & 1 \\ 1 & 0 & 1 & 0 \\ -1 & 1 & 0 & 1 \\ 2 & 0 & 0 & 0 \end{bmatrix} + \right.$$

$$\begin{bmatrix} 2 & 0 & 1 \\ 0 & 0 & 1 \end{bmatrix} \begin{bmatrix} 0 & 0 & 1 & 1 \\ 1 & 0 & 0 & 1 \\ 0 & 0 & 1 & 0 \end{bmatrix} +$$

$$\left. \begin{bmatrix} -1 & 5 \\ 0 & 2 \end{bmatrix} \begin{bmatrix} 2 & 1 & 2 & 0 \\ 1 & 0 & 1 & 1 \end{bmatrix} \right\} \cup$$

$$\left\{ \begin{bmatrix} 3 & 1 \\ 6 & 2 \\ 1 & 0 \end{bmatrix} \begin{bmatrix} 1 & 0 \\ 0 & 1 \end{bmatrix} + \begin{bmatrix} 1 & 1 & 1 & 0 & 1 \\ 0 & 1 & 0 & 1 & 0 \\ 0 & 0 & 1 & 1 & 1 \end{bmatrix} \begin{bmatrix} 2 & 1 \\ 1 & 2 \\ 0 & 1 \\ 1 & 1 \\ 0 & 1 \end{bmatrix} + \begin{bmatrix} 1 \\ 1 \\ 0 \end{bmatrix} \begin{bmatrix} 1 & 1 \end{bmatrix} \right\}$$

$$= \left\{ \begin{bmatrix} 25 & 6 & 1 & 6 \\ 3 & 2 & 1 & 2 \end{bmatrix} + \begin{bmatrix} 0 & 0 & 3 & 2 \\ 0 & 0 & 1 & 0 \end{bmatrix} + \begin{bmatrix} 3 & -1 & 3 & 5 \\ 2 & 0 & 2 & 2 \end{bmatrix} \right\} \cup$$

$$\left\{ \begin{bmatrix} 3 & 1 \\ 6 & 2 \\ 1 & 0 \end{bmatrix} + \begin{bmatrix} 3 & 5 \\ 2 & 3 \\ 1 & 3 \end{bmatrix} + \begin{bmatrix} 1 & 1 \\ 1 & 1 \\ 0 & 0 \end{bmatrix} \right\}$$

$$= \begin{bmatrix} 28 & 5 & 7 & 13 \\ 5 & 2 & 4 & 4 \end{bmatrix} \cup \begin{bmatrix} 7 & 6 \\ 9 & 6 \\ 2 & 3 \end{bmatrix}.$$

We see the product of two superbimatrices is only a bimatrix.

To this end we define two new notions.

**DEFINITION 2.14:** *Let $A = A_1 \cup A_2$ be a mixed rectangular $m \times n$ superbimatrix with $m < n$. If in both $A_1$ and $A_2$ partition is only between the columns i.e., only vertical partition then we call A to be a row superbivector.*

*Example 2.41:* Let $A = A_1 \cup A_2$ be a mixed rectangular superbimatrix where

$$A_1 = \begin{bmatrix} 3 & 0 & 1 & 2 & 1 & 3 & 3 & 2 & 1 \\ 1 & 1 & 1 & 3 & 1 & 0 & 1 & 1 & 0 \\ 2 & 1 & 1 & 4 & 1 & 0 & 0 & 1 & 0 \end{bmatrix}$$

and

$$A_2 = \begin{bmatrix} 3 & 1 & 5 & 1 & 0 & 1 & 1 & 1 & 1 & 0 \\ 1 & 1 & 7 & 2 & 1 & 2 & 0 & 0 & 1 & 1 \\ 0 & 0 & 8 & 3 & 2 & 3 & 2 & 5 & 7 & 8 \\ 2 & 1 & 9 & 4 & 3 & 4 & 1 & 2 & 3 & 4 \end{bmatrix};$$

We call A to be a row superbivector. It is clear the partition of the matrices are only by the vertical lines i.e., between the columns and no partition is carried out between the rows.

**DEFINITION 2.15:** *Let $A = A_1 \cup A_2$ be a mixed rectangular $m \times n$ superbimatrix with $m > n$. If we have on both $A_1$ and $A_2$ only horizontal partitions then we call A to be a column superbivector. When we say horizontal partitions it takes place only between the rows.*

*Example 2.42:* Let $A = A_1 \cup A_2$ a rectangular superbimatrix, where

$$A_1 = \begin{bmatrix} 1 & 0 & 1 & 1 \\ 2 & 1 & 2 & 0 \\ 0 & 1 & 0 & 1 \\ \hline 3 & 1 & 2 & 5 \\ 1 & 2 & 3 & 4 \\ 5 & 6 & 7 & 8 \\ 0 & 1 & 2 & 3 \\ 1 & 1 & 0 & 5 \\ 2 & 5 & 7 & 1 \end{bmatrix}$$

and

$$A_2 = \begin{bmatrix} 1 & 1 & 0 \\ \hline 1 & 1 & 1 \\ 0 & 1 & 2 \\ 3 & 4 & 5 \\ 6 & 7 & 8 \\ 9 & 0 & 1 \\ \hline 0 & 1 & 1 \\ 1 & 0 & 1 \\ 1 & 1 & 0 \\ \hline 1 & 2 & 3 \\ 4 & 5 & 6 \end{bmatrix}$$

is a column superbivector.

Now we proceed onto define the notion of product of these type of superbivectors.

**DEFINITION 2.16:** *Let $A = A_1 \cup A_2$ a column superbivector and $B = B_1 \cup B_2$ be a row superbivector. Now how to define the product BA. BA is defined if the following conditions are satisfied.*

*If $A = A_1 \cup A_2$ where*

$$A_1 = \begin{bmatrix} A_1^1 \\ \vdots \\ A_{m_1}^1 \end{bmatrix}$$

and

$$A_2 = \begin{bmatrix} A_1^2 \\ \vdots \\ A_{m_2}^2 \end{bmatrix}$$

and in $B = B_1 \cup B_2$ we have $B_1 = \begin{bmatrix} B_1^1 & | & \cdots & | & B_{m_1}^1 \end{bmatrix}$ and $B_2 = \begin{bmatrix} B_1^2 & | & \cdots & | & B_{m_2}^2 \end{bmatrix}$ then first

$$\begin{aligned} B.A &= (B_1 \cup B_2)(A_1 \cup A_2) \\ &= B_1 A_1 \cup B_2 A_2, \end{aligned}$$

*$B_1 A_1$ and $B_2 A_2$ is the super vector product defined only if the number of columns in each of $B_i^1$, $1 \leq i \leq m_1$ is equal to the number of rows in each of $A_j^1$ ; $1 \leq j \leq m_1$ and the number of columns in each of $B_i^2$, $1 \leq i \leq m_2$ is equal to the number of rows in each of $A_j^2$, $1 \leq j \leq m$. Now*

$$\begin{aligned} BA &= B_1 A_1 \cup B_2 A_2 \\ &= \begin{bmatrix} B_1^1 & | & \cdots & | & B_{m_1}^1 \end{bmatrix} \begin{bmatrix} A_1^1 \\ \vdots \\ A_{m_1}^1 \end{bmatrix} \cup \begin{bmatrix} B_1^2 & | & \cdots & | & B_{m_2}^2 \end{bmatrix} \begin{bmatrix} A_1^2 \\ \vdots \\ A_{m_2}^2 \end{bmatrix} \\ &= \{B_1^1 A_1^1 + \ldots + B_{m_1}^1 A_{m_1}^1\} \cup \{B_1^2 A_1^2 + \ldots + B_{m_2}^2 A_{m_2}^2\} \\ &= C \cup D; \end{aligned}$$

*we see C and D are not super vectors they are only just matrices. Thus $C \cup D$ is only a bimatrix and not a superbimatrix.*

We have illustrated this type of product in example 2.39 and 2.40. The product defined using row superbivectors and column superbivectors results only in bivectors and this product will be known as the minor product of superbivectors.

Now we proceed onto define the notion of major product of these superbivectors. Before we define them abstractly we give some examples of them.

*Example 2.43:* Let $A = A_1 \cup A_2$ and $B = B_1 \cup B_2$ be two superbivectors where $A = A_1 \cup A_2$ with

$$A_1 = \begin{bmatrix} 1 & 2 & 3 \\ 0 & 1 & 1 \\ \hline 3 & 0 & 1 \\ 1 & -1 & 0 \\ 5 & 0 & 1 \end{bmatrix}$$

and

$$A_2 = \begin{bmatrix} 2 & 1 \\ 3 & 0 \\ \hline 1 & 1 \\ \hline 1 & 4 \\ 2 & 5 \\ 3 & 0 \end{bmatrix}$$

i.e., $A = A_1 \cup A_2$ is a column superbivector. Given $B = B_1 \cup B_2$ where

$$B_1 = \begin{bmatrix} 1 & 2 & 5 & 1 & 3 & 1 \\ 0 & 0 & 1 & 0 & -1 & 0 \\ 1 & 1 & 2 & 0 & 1 & 1 \end{bmatrix}$$

and

$$B_2 = \begin{bmatrix} 2 & 1 & 3 & 1 & 1 & 5 & 1 & 3 & 1 \\ 1 & 1 & 1 & 2 & 1 & 0 & 1 & 0 & 1 \end{bmatrix}$$

is the row superbivector. The major byproduct

$$
\begin{aligned}
AB &= (A_1 \cup A_2)(B_1 \cup B_2) \\
&= A_1 B_1 \cup A_2 B_2
\end{aligned}
$$

$$
= \begin{bmatrix} 1 & 2 & 3 \\ 0 & 1 & 1 \\ \hline 3 & 0 & 1 \\ 1 & -1 & 0 \\ 5 & 0 & 1 \end{bmatrix} \begin{bmatrix} 1 & 2 & 5 & 1 & 3 & 1 \\ 0 & 0 & 1 & 0 & -1 & 0 \\ 1 & 1 & 2 & 0 & 1 & 1 \end{bmatrix} \cup
$$

$$
\begin{bmatrix} 2 & 1 \\ 3 & 0 \\ \hline 1 & 1 \\ \hline 1 & 4 \\ 2 & 5 \\ 3 & 0 \end{bmatrix} \begin{bmatrix} 2 & 3 & 1 & 1 & 1 & 5 & 1 & 3 & 1 \\ 1 & 1 & 1 & 2 & 1 & 0 & 1 & 0 & 1 \end{bmatrix} =
$$

$$
\begin{bmatrix}
\begin{bmatrix} 1 & 2 & 3 \\ 0 & 1 & 1 \end{bmatrix}\begin{bmatrix} 1 \\ 0 \\ 1 \end{bmatrix} & \begin{pmatrix} 1 & 2 & 3 \\ 0 & 1 & 1 \end{pmatrix}\begin{pmatrix} 2 & 5 & 1 \\ 0 & 1 & 0 \\ 1 & 2 & 0 \end{pmatrix} & \begin{pmatrix} 1 & 2 & 3 \\ 0 & 1 & 1 \end{pmatrix}\begin{bmatrix} 3 & 1 \\ -1 & 0 \\ 1 & 1 \end{bmatrix} \\
\hline
\begin{bmatrix} 3 & 0 & 1 \\ 1 & -1 & 0 \\ 5 & 0 & 1 \end{bmatrix}\begin{bmatrix} 1 \\ 0 \\ 1 \end{bmatrix} & \begin{pmatrix} 3 & 0 & 1 \\ 1 & -1 & 0 \\ 5 & 0 & 1 \end{pmatrix}\begin{pmatrix} 2 & 5 & 1 \\ 0 & 1 & 0 \\ 1 & 2 & 0 \end{pmatrix} & \begin{pmatrix} 3 & 0 & 1 \\ 1 & -1 & 0 \\ 5 & 0 & 1 \end{pmatrix}\begin{bmatrix} 3 & 1 \\ -1 & 0 \\ 1 & 1 \end{bmatrix}
\end{bmatrix}
$$

$$
\cup \begin{bmatrix}
\begin{pmatrix} 2 & 1 \\ 3 & 0 \end{pmatrix}\begin{pmatrix} 2 & 3 & 1 \\ 1 & 1 & 1 \end{pmatrix} & \begin{pmatrix} 2 & 1 \\ 3 & 0 \end{pmatrix}\begin{pmatrix} 1 & 1 \\ 2 & 1 \end{pmatrix} & \begin{pmatrix} 2 & 1 \\ 3 & 0 \end{pmatrix}\begin{pmatrix} 5 & 1 & 3 & 1 \\ 0 & 1 & 0 & 1 \end{pmatrix} \\
\hline
(1 \; 1)\begin{pmatrix} 2 & 3 & 1 \\ 1 & 1 & 1 \end{pmatrix} & (1 \; 1)\begin{pmatrix} 1 & 1 \\ 2 & 1 \end{pmatrix} & (1 \; 1)\begin{pmatrix} 5 & 1 & 3 & 1 \\ 0 & 1 & 0 & 1 \end{pmatrix} \\
\hline
\begin{pmatrix} 1 & 4 \\ 2 & 5 \\ 3 & 0 \end{pmatrix}\begin{pmatrix} 2 & 3 & 1 \\ 1 & 1 & 1 \end{pmatrix} & \begin{pmatrix} 1 & 4 \\ 2 & 5 \\ 3 & 0 \end{pmatrix}\begin{pmatrix} 1 & 1 \\ 2 & 1 \end{pmatrix} & \begin{pmatrix} 1 & 4 \\ 2 & 5 \\ 3 & 0 \end{pmatrix}\begin{pmatrix} 5 & 1 & 3 & 1 \\ 0 & 1 & 0 & 1 \end{pmatrix}
\end{bmatrix}
$$

$$= \begin{bmatrix} 4 & 5 & 13 & 1 & 4 & 4 \\ 1 & 1 & 3 & 0 & 0 & 1 \\ \hline 4 & 7 & 17 & 3 & 10 & 4 \\ 1 & 2 & 4 & 1 & 4 & 1 \\ 6 & 11 & 27 & 5 & 16 & 6 \end{bmatrix} \cup$$

$$\begin{bmatrix} 5 & 7 & 3 & 4 & 3 & 10 & 3 & 6 & 3 \\ 6 & 9 & 3 & 3 & 3 & 15 & 3 & 9 & 3 \\ 3 & 4 & 2 & 3 & 2 & 5 & 2 & 3 & 2 \\ \hline 6 & 7 & 5 & 9 & 5 & 5 & 5 & 3 & 5 \\ 9 & 11 & 7 & 12 & 7 & 10 & 7 & 6 & 7 \\ 6 & 9 & 3 & 3 & 3 & 15 & 3 & 9 & 3 \end{bmatrix}$$

$= C \cup D$ where both C and D are superbimatrices and not superbivectors.

We illustrate the same with one more example.

***Example 2.44:*** Let $A = A_1 \cup A_2$ and $B = B_1 \cup B_2$ be two superbivectors where $A = A_1 \cup A_2$ is a column superbivector and $B = B_1 \cup B_2$ is a row superbivector with

$$A_1 = \begin{bmatrix} 2 & 0 & 1 \\ 1 & 1 & 1 \\ \hline 3 & 7 & 0 \\ 1 & 1 & 0 \\ 2 & 0 & 1 \\ 3 & 5 & 1 \\ \hline 1 & 1 & 1 \\ 0 & 1 & 0 \\ 1 & 0 & 1 \\ 0 & 0 & 1 \end{bmatrix}$$

and

$$A_2 = \begin{bmatrix} 1 & 2 \\ 3 & 4 \\ 5 & 6 \\ \hline 1 & 0 \\ 0 & 1 \\ \hline 1 & 1 \\ 1 & 0 \\ 1 & 1 \end{bmatrix}$$

is a column superbivector. Now $B = B_1 \cup B_2$ where

$$B_1 = \begin{bmatrix} 1 & 0 & 1 & 1 & 0 & 5 & 1 & 2 \\ 2 & 1 & 0 & 1 & 1 & 0 & 1 & 1 \\ 3 & 1 & 1 & 0 & 1 & 0 & 1 & 2 \end{bmatrix}$$

and

$$B_2 = \begin{bmatrix} 1 & 0 & 1 & 1 & 1 & 1 & 1 & 1 & 1 \\ 1 & 1 & 0 & 0 & 1 & 1 & 0 & 1 & 0 & 1 \end{bmatrix}$$

is the row superbivector. Now the major byproduct of AB is defined as $AB = (A_1 \cup A_2)(B_1 \cup B_2) = A_1 B_1 \cup A_2 B_2$

$$= \begin{bmatrix} 2 & 0 & 1 \\ 1 & 1 & 1 \\ 3 & 7 & 0 \\ \hline 1 & 1 & 0 \\ 2 & 0 & 1 \\ 3 & 5 & 1 \\ \hline 1 & 1 & 1 \\ 0 & 1 & 0 \\ 1 & 0 & 1 \\ 0 & 0 & 1 \end{bmatrix} \begin{bmatrix} 1 & 0 & 1 & 1 & 0 & 5 & 1 & 2 \\ 2 & 1 & 0 & 1 & 1 & 0 & 1 & 1 \\ 3 & 1 & 1 & 0 & 1 & 0 & 1 & 2 \end{bmatrix} \cup$$

$$\begin{bmatrix} 1 & 2 \\ 3 & 4 \\ 5 & 6 \\ \hline 1 & 0 \\ 0 & 1 \\ \hline 1 & 1 \\ 1 & 0 \\ 1 & 1 \end{bmatrix} \left[\begin{array}{cc|cccc|cccc} 1 & 0 & 1 & 1 & 1 & 1 & 1 & 1 & 1 & 1 \\ 1 & 1 & 0 & 0 & 1 & 1 & 0 & 1 & 0 & 1 \end{array}\right] =$$

$$\left[\begin{array}{c|c|c} \begin{bmatrix} 2 & 0 & 1 \\ 1 & 1 & 1 \end{bmatrix}\begin{bmatrix} 1 & 0 & 1 & 1 \\ 2 & 1 & 0 & 1 \\ 3 & 1 & 1 & 0 \end{bmatrix} & \begin{bmatrix} 2 & 0 & 1 \\ 1 & 1 & 1 \end{bmatrix}\begin{bmatrix} 0 \\ 1 \\ 1 \end{bmatrix} & \begin{bmatrix} 2 & 0 & 1 \\ 1 & 1 & 1 \end{bmatrix}\begin{bmatrix} 5 & 1 & 2 \\ 0 & 1 & 1 \\ 0 & 1 & 2 \end{bmatrix} \\ \hline \begin{bmatrix} 3 & 7 & 0 \\ 1 & 1 & 0 \\ 2 & 0 & 1 \\ 3 & 5 & 1 \end{bmatrix}\begin{bmatrix} 1 & 0 & 1 & 1 \\ 2 & 1 & 0 & 1 \\ 3 & 1 & 1 & 0 \end{bmatrix} & \begin{bmatrix} 3 & 7 & 0 \\ 1 & 1 & 0 \\ 2 & 0 & 1 \\ 3 & 5 & 1 \end{bmatrix}\begin{bmatrix} 0 \\ 1 \\ 1 \end{bmatrix} & \begin{bmatrix} 3 & 7 & 0 \\ 1 & 1 & 0 \\ 2 & 0 & 1 \\ 3 & 5 & 1 \end{bmatrix}\begin{bmatrix} 5 & 1 & 2 \\ 0 & 1 & 1 \\ 0 & 1 & 2 \end{bmatrix} \\ \hline \begin{bmatrix} 1 & 1 & 1 \\ 0 & 1 & 0 \\ 1 & 0 & 1 \\ 0 & 0 & 1 \end{bmatrix}\begin{bmatrix} 1 & 0 & 1 & 1 \\ 2 & 1 & 0 & 1 \\ 3 & 1 & 1 & 0 \end{bmatrix} & \begin{bmatrix} 1 & 1 & 1 \\ 0 & 1 & 0 \\ 1 & 0 & 1 \\ 0 & 0 & 1 \end{bmatrix}\begin{bmatrix} 0 \\ 1 \\ 1 \end{bmatrix} & \begin{bmatrix} 1 & 1 & 1 \\ 0 & 1 & 0 \\ 1 & 0 & 1 \\ 0 & 0 & 1 \end{bmatrix}\begin{bmatrix} 5 & 1 & 2 \\ 0 & 1 & 1 \\ 0 & 1 & 2 \end{bmatrix} \end{array}\right]$$

$$\subset \left[\begin{array}{c|c|c} \begin{pmatrix} 1 & 2 \\ 3 & 4 \\ 5 & 6 \end{pmatrix}\begin{bmatrix} 1 & 0 \\ 1 & 1 \end{bmatrix} & \begin{pmatrix} 1 & 2 \\ 3 & 4 \\ 5 & 6 \end{pmatrix}\begin{bmatrix} 1 & 1 & 1 & 1 \\ 0 & 0 & 1 & 1 \end{bmatrix} & \begin{pmatrix} 1 & 2 \\ 3 & 4 \\ 5 & 6 \end{pmatrix}\begin{bmatrix} 1 & 1 & 1 & 1 \\ 0 & 1 & 0 & 1 \end{bmatrix} \\ \hline \begin{pmatrix} 1 & 0 \\ 0 & 1 \end{pmatrix}\begin{bmatrix} 1 & 0 \\ 1 & 1 \end{bmatrix} & \begin{pmatrix} 1 & 0 \\ 0 & 1 \end{pmatrix}\begin{bmatrix} 1 & 1 & 1 & 1 \\ 0 & 0 & 1 & 1 \end{bmatrix} & \begin{pmatrix} 1 & 0 \\ 0 & 1 \end{pmatrix}\begin{bmatrix} 1 & 1 & 1 & 1 \\ 0 & 1 & 0 & 1 \end{bmatrix} \\ \hline \begin{pmatrix} 1 & 1 \\ 1 & 0 \\ 1 & 1 \end{pmatrix}\begin{bmatrix} 1 & 0 \\ 1 & 1 \end{bmatrix} & \begin{pmatrix} 1 & 1 \\ 1 & 0 \\ 1 & 1 \end{pmatrix}\begin{bmatrix} 1 & 1 & 1 & 1 \\ 0 & 0 & 1 & 1 \end{bmatrix} & \begin{pmatrix} 1 & 1 \\ 1 & 0 \\ 1 & 1 \end{pmatrix}\begin{bmatrix} 1 & 1 & 1 & 1 \\ 0 & 1 & 0 & 1 \end{bmatrix} \end{array}\right]$$

$$= \begin{bmatrix} 5 & 1 & 3 & 2 & 1 & 10 & 3 & 6 \\ 6 & 2 & 2 & 2 & 2 & 5 & 3 & 5 \\ \hline 17 & 7 & 3 & 10 & 7 & 15 & 10 & 13 \\ 3 & 1 & 1 & 2 & 1 & 5 & 2 & 3 \\ 5 & 1 & 3 & 2 & 1 & 10 & 3 & 6 \\ 16 & 6 & 4 & 8 & 6 & 15 & 9 & 13 \\ 6 & 2 & 2 & 2 & 2 & 5 & 3 & 5 \\ 2 & 1 & 0 & 1 & 1 & 0 & 1 & 1 \\ 4 & 1 & 2 & 1 & 1 & 5 & 2 & 4 \\ 3 & 1 & 1 & 0 & 1 & 0 & 1 & 2 \end{bmatrix} \cup$$

$$\begin{bmatrix} 3 & 2 & 1 & 1 & 3 & 3 & 1 & 3 & 1 & 3 \\ 7 & 4 & 3 & 3 & 7 & 7 & 3 & 7 & 3 & 7 \\ 11 & 6 & 5 & 5 & 11 & 11 & 5 & 11 & 5 & 11 \\ \hline 1 & 0 & 1 & 1 & 1 & 1 & 1 & 1 & 1 & 1 \\ 1 & 1 & 0 & 0 & 1 & 1 & 0 & 1 & 0 & 1 \\ \hline 2 & 1 & 1 & 1 & 2 & 2 & 1 & 2 & 1 & 2 \\ 1 & 0 & 1 & 1 & 1 & 1 & 1 & 1 & 1 & 1 \\ 2 & 1 & 1 & 1 & 2 & 2 & 1 & 2 & 1 & 2 \end{bmatrix}$$

$= C_1 \cup C_2$. Thus the resultant of a column superbivector with a row superbivector is a superbimatrix.

Now we proceed on to define the major byproduct of superbivectors.

**DEFINITION 2.17:** *Let $A = A_1 \cup A_2$ be a column superbivector and $B = B_1 \cup B_2$ be a row superbivector. The major byproduct of these two superbivectors is defined to be AB where $AB = (A_1 \cup A_2)(B_1 \cup B_2) = A_1 B_1 \cup A_2 B_2$ is compatible if and only if the number of columns in $A_1$ must be equal to the number of rows in $B_1$ and the number of columns in $A_2$ must be equal to the number of rows in $B_2$ respectively. Then the resultant bimatrix is* *always a superbimatrix. Thus the major byproduct yields a superbimatrix which is neither a row superbivector nor a column superbivector.*

Now the immediate application of these major byproduct of superbivectors is the product of a superbivector with its transpose. The example 2.43 and 2.44 are illustrations of major byproduct of superbivectors.

Now we proceed onto give an example of a superbivector with its transpose.

***Example 2.45:*** Let $A = A_1 \cup A_2$ be a column superbivector. Let $A^T$ be its transpose. The byproduct $AA^T$ gives a superbimatrix which is neither a column superbivector nor a row superbivector. Given $A = A_1 \cup A_2$ where

$$A_1 = \begin{bmatrix} 3 & 1 & 0 \\ -1 & 1 & 6 \\ 0 & 1 & 1 \\ \hline 2 & 1 & 0 \\ \hline 1 & 2 & 3 \\ 1 & 0 & 1 \\ 0 & 1 & 0 \\ 1 & 0 & 1 \end{bmatrix}$$

and

$$A_2 = \begin{bmatrix} 2 & 1 & 0 & 4 \\ 1 & 1 & 6 & 0 \\ \hline 0 & 0 & 1 & 1 \\ 1 & 0 & 1 & 1 \\ 0 & 5 & 2 & 3 \\ 1 & 1 & 0 & 1 \\ 2 & 0 & 2 & 1 \end{bmatrix}.$$

Now

$$A^T = (A_1 \cup A_2)^T = A_1^T \cup A_2^T.$$

$$A_1^T = \begin{bmatrix} 3 & -1 & 0 & | & 2 & | & 1 & 1 & 0 & 1 \\ 1 & 1 & 1 & | & 1 & | & 2 & 0 & 1 & 0 \\ 0 & 6 & 1 & | & 0 & | & 3 & 1 & 0 & 1 \end{bmatrix}$$

and

$$A_2^T = \begin{bmatrix} 2 & 1 & | & 0 & 1 & 0 & 1 & | & 2 \\ 1 & 1 & | & 0 & 0 & 5 & 1 & | & 0 \\ 0 & 6 & | & 1 & 1 & 2 & 0 & | & 2 \\ 4 & 0 & | & 1 & 1 & 3 & 1 & | & 1 \end{bmatrix}$$

$$\begin{aligned} AA^T &= (A_1 \cup A_2)(A_1 \cup A_2)^T \\ &= (A_1 \cup A_2)(A_1^T \cup A_2^T) \\ &= A_1 A_1^T \cup A_2 A_2^T. \end{aligned}$$

$$= \begin{bmatrix} 3 & 1 & 0 \\ -1 & 1 & 6 \\ 0 & 1 & 1 \\ \hline 2 & 1 & 0 \\ \hline 1 & 2 & 3 \\ 1 & 0 & 1 \\ 0 & 1 & 0 \\ 1 & 0 & 1 \end{bmatrix} \begin{bmatrix} 3 & -1 & 0 & | & 2 & | & 1 & 1 & 0 & 1 \\ 1 & 1 & 1 & | & 1 & | & 2 & 0 & 1 & 0 \\ 0 & 6 & 1 & | & 0 & | & 3 & 1 & 0 & 1 \end{bmatrix} \cup$$

$$\begin{bmatrix} 2 & 1 & 0 & 4 \\ 1 & 1 & 6 & 0 \\ \hline 0 & 0 & 1 & 1 \\ 1 & 0 & 1 & 1 \\ 0 & 5 & 2 & 3 \\ 1 & 1 & 0 & 1 \\ \hline 2 & 0 & 2 & 1 \end{bmatrix} \begin{bmatrix} 2 & 1 & | & 0 & 1 & 0 & 1 & | & 2 \\ 1 & 1 & | & 0 & 0 & 5 & 1 & | & 0 \\ 0 & 6 & | & 1 & 1 & 2 & 0 & | & 2 \\ 4 & 0 & | & 1 & 1 & 3 & 1 & | & 1 \end{bmatrix} =$$

$$\left[ \begin{array}{c|c|c} \begin{pmatrix} 3 & 1 & 0 \\ -1 & 1 & 6 \\ 0 & 1 & 1 \end{pmatrix} \begin{pmatrix} 3 & -1 & 0 \\ 1 & 1 & 1 \\ 0 & 6 & 1 \end{pmatrix} & \begin{pmatrix} 3 & 1 & 0 \\ -1 & 1 & 6 \\ 0 & 1 & 1 \end{pmatrix} \begin{pmatrix} 2 \\ 1 \\ 0 \end{pmatrix} & \begin{pmatrix} 3 & 1 & 0 \\ -1 & 1 & 6 \\ 0 & 1 & 1 \end{pmatrix} \begin{pmatrix} 1 & 1 & 0 & 1 \\ 2 & 0 & 1 & 0 \\ 3 & 1 & 0 & 1 \end{pmatrix} \\ \hline \begin{pmatrix} 2 & 1 & 0 \end{pmatrix} \begin{pmatrix} 3 & -1 & 0 \\ 1 & 1 & 1 \\ 0 & 6 & 1 \end{pmatrix} & \begin{pmatrix} 2 & 1 & 0 \end{pmatrix} \begin{pmatrix} 2 \\ 1 \\ 0 \end{pmatrix} & \begin{pmatrix} 2 & 1 & 0 \end{pmatrix} \begin{pmatrix} 1 & 1 & 0 & 1 \\ 2 & 0 & 1 & 0 \\ 3 & 1 & 0 & 1 \end{pmatrix} \\ \hline \begin{pmatrix} 1 & 2 & 3 \\ 1 & 0 & 1 \\ 0 & 1 & 0 \\ 1 & 0 & 1 \end{pmatrix} \begin{pmatrix} 3 & -1 & 0 \\ 1 & 1 & 1 \\ 0 & 6 & 1 \end{pmatrix} & \begin{pmatrix} 1 & 2 & 3 \\ 1 & 0 & 1 \\ 0 & 1 & 0 \\ 1 & 0 & 1 \end{pmatrix} \begin{pmatrix} 2 \\ 1 \\ 0 \end{pmatrix} & \begin{pmatrix} 1 & 2 & 3 \\ 1 & 0 & 1 \\ 0 & 1 & 0 \\ 1 & 0 & 1 \end{pmatrix} \begin{pmatrix} 1 & 1 & 0 & 1 \\ 2 & 0 & 1 & 0 \\ 3 & 1 & 0 & 1 \end{pmatrix} \end{array} \right]$$

$\cup$

$$\left[ \begin{array}{c|c|c} \begin{pmatrix} 2 & 1 & 0 & 4 \\ 1 & 1 & 6 & 0 \end{pmatrix} \begin{bmatrix} 2 & 1 \\ 1 & 1 \\ 0 & 6 \\ 4 & 0 \end{bmatrix} & \begin{pmatrix} 2 & 1 & 0 & 4 \\ 1 & 1 & 6 & 0 \end{pmatrix} \begin{bmatrix} 0 & 1 & 0 & 1 \\ 0 & 0 & 5 & 1 \\ 1 & 1 & 2 & 0 \\ 1 & 1 & 3 & 1 \end{bmatrix} & \begin{pmatrix} 2 & 1 & 0 & 4 \\ 1 & 1 & 6 & 0 \end{pmatrix} \begin{bmatrix} 2 \\ 0 \\ 2 \\ 1 \end{bmatrix} \\ \hline \begin{pmatrix} 0 & 0 & 1 & 1 \\ 1 & 0 & 1 & 1 \\ 0 & 5 & 2 & 3 \\ 1 & 1 & 0 & 1 \end{pmatrix} \begin{bmatrix} 2 & 1 \\ 1 & 1 \\ 0 & 6 \\ 4 & 0 \end{bmatrix} & \begin{pmatrix} 0 & 0 & 1 & 1 \\ 1 & 0 & 1 & 1 \\ 0 & 5 & 2 & 3 \\ 1 & 1 & 0 & 1 \end{pmatrix} \begin{bmatrix} 0 & 1 & 0 & 1 \\ 0 & 0 & 5 & 1 \\ 1 & 1 & 2 & 0 \\ 1 & 1 & 3 & 1 \end{bmatrix} & \begin{pmatrix} 0 & 0 & 1 & 1 \\ 1 & 0 & 1 & 1 \\ 0 & 5 & 2 & 3 \\ 1 & 1 & 0 & 1 \end{pmatrix} \begin{bmatrix} 2 \\ 0 \\ 2 \\ 1 \end{bmatrix} \\ \hline \begin{pmatrix} 2 & 0 & 2 & 1 \end{pmatrix} \begin{bmatrix} 2 & 1 \\ 1 & 1 \\ 0 & 6 \\ 4 & 0 \end{bmatrix} & \begin{pmatrix} 2 & 0 & 2 & 1 \end{pmatrix} \begin{bmatrix} 0 & 1 & 0 & 1 \\ 0 & 0 & 5 & 1 \\ 1 & 1 & 2 & 0 \\ 1 & 1 & 3 & 1 \end{bmatrix} & \begin{pmatrix} 2 & 0 & 2 & 1 \end{pmatrix} \begin{bmatrix} 2 \\ 0 \\ 2 \\ 1 \end{bmatrix} \end{array} \right]$$

$$= \begin{bmatrix} 10 & -2 & 1 & 7 & 5 & 3 & 1 & 3 \\ -2 & 38 & 7 & -1 & 19 & 5 & 1 & 5 \\ 1 & 7 & 2 & 1 & 5 & 1 & 1 & 1 \\ \hline 7 & -1 & 1 & 5 & 4 & 2 & 1 & 2 \\ \hline 5 & 19 & 5 & 4 & 14 & 4 & 2 & 4 \\ 3 & 5 & 1 & 2 & 4 & 2 & 0 & 2 \\ 1 & 1 & 1 & 1 & 2 & 0 & 1 & 0 \\ 3 & 5 & 1 & 2 & 4 & 2 & 0 & 2 \end{bmatrix} \cup$$

$$\begin{bmatrix} 21 & 3 & 4 & 6 & 17 & 7 & 8 \\ 3 & 38 & 6 & 7 & 17 & 2 & 14 \\ \hline 4 & 6 & 2 & 2 & 5 & 1 & 3 \\ 6 & 7 & 2 & 3 & 5 & 2 & 5 \\ 17 & 17 & 5 & 5 & 38 & 8 & 7 \\ 7 & 2 & 1 & 2 & 8 & 3 & 3 \\ \hline 8 & 14 & 3 & 5 & 7 & 3 & 9 \end{bmatrix}$$

$= S_1 \cup S_2 = S$ we see both $S_1$ and $S_2$ are symmetric supermatrices thus $AA^T = S$ is a symmetric superbimatrix. Thus this product of a column superbivector with its transpose

$$A_1^T = \begin{bmatrix} 1 & 0 & 1 \\ 2 & 3 & 4 \\ \hline 1 & 0 & 3 \\ \hline 1 & 1 & 0 \\ 5 & 1 & 1 \\ 2 & 0 & 0 \\ 3 & 1 & 1 \\ 1 & 0 & 0 \end{bmatrix}$$

and

$$A_2^T = \begin{bmatrix} 1 & 2 & 3 & 4 \\ 1 & 0 & 1 & 1 \\ 3 & 6 & 1 & 0 \\ 1 & 1 & 0 & 1 \\ \hline 2 & 1 & 3 & 5 \\ 0 & 0 & 1 & 2 \\ 3 & 1 & 2 & 1 \\ \hline 1 & 1 & 0 & 1 \end{bmatrix}.$$

$$\begin{aligned} A^T A &= (A_1 \cup A_2)^T (A_1 \cup A_2) \\ &= (A_1^T \cup A_2^T)(A_1 \cup A_2) \\ &= A_1^T A_1 \cup A_2^T A_2 \end{aligned}$$

Now $A_1^T A_1 \cup A_2^T A_2$ yields a superbimatrix which is always symmetric.

*Example 2.46:* Let $A = A_1 \cup A_2$ be a row superbivector. Now we find the product of $A^T$ with A. Given $A = A_1 \cup A_2$ where

$$A_1 = \begin{bmatrix} 1 & 2 & 1 & 1 & 5 & 2 & 3 & 1 \\ 0 & 3 & 0 & 1 & 1 & 0 & 1 & 0 \\ 1 & 4 & 3 & 0 & 1 & 0 & 1 & 0 \end{bmatrix}$$

and

$$A_2 = \begin{bmatrix} 1 & 1 & 3 & 1 & 2 & 0 & 3 & 1 \\ 2 & 0 & 6 & 1 & 1 & 0 & 1 & 1 \\ 3 & 1 & 1 & 0 & 3 & 1 & 2 & 0 \\ 4 & 1 & 0 & 1 & 5 & 2 & 1 & 1 \end{bmatrix}.$$

Now

$$A^T = (A_1 \cup A_2)^T = A_1^T \cup A_2^T$$

$$\begin{bmatrix} 1 & 0 & 1 \\ 2 & 3 & 4 \\ \hline 1 & 0 & 3 \\ \hline 1 & 1 & 0 \\ 5 & 1 & 1 \\ 2 & 0 & 0 \\ 3 & 1 & 1 \\ 1 & 0 & 0 \end{bmatrix} \cup \begin{bmatrix} 1 & 2 & 3 & 4 \\ 1 & 0 & 1 & 1 \\ 3 & 6 & 1 & 0 \\ \hline 1 & 1 & 0 & 1 \\ 2 & 1 & 3 & 5 \\ 0 & 0 & 1 & 2 \\ \hline 3 & 1 & 2 & 1 \\ 1 & 1 & 0 & 1 \end{bmatrix} \cdot \left[\begin{array}{cc|c|ccccc} 1 & 2 & 1 & 1 & 5 & 2 & 3 & 1 \\ 0 & 3 & 0 & 1 & 1 & 0 & 1 & 0 \\ 1 & 4 & 3 & 0 & 1 & 0 & 1 & 0 \end{array}\right] \times$$

$$\begin{bmatrix} 1 & 0 & 1 \\ 2 & 3 & 4 \\ \hline 1 & 0 & 3 \\ \hline 1 & 1 & 0 \\ 5 & 1 & 1 \\ 2 & 0 & 0 \\ 3 & 1 & 1 \\ 1 & 0 & 0 \end{bmatrix} \cup \begin{bmatrix} 1 & 2 & 3 & 4 \\ 1 & 0 & 1 & 1 \\ 3 & 6 & 1 & 0 \\ \hline 1 & 1 & 0 & 1 \\ 2 & 1 & 3 & 5 \\ 0 & 0 & 1 & 2 \\ \hline 3 & 1 & 2 & 1 \\ 1 & 1 & 0 & 1 \end{bmatrix} \cdot \left[\begin{array}{cccc|ccc|c} 1 & 1 & 3 & 1 & 2 & 0 & 3 & 1 \\ 2 & 0 & 6 & 1 & 1 & 0 & 1 & 1 \\ 3 & 1 & 1 & 0 & 3 & 1 & 2 & 0 \\ 4 & 1 & 0 & 1 & 5 & 2 & 1 & 1 \end{array}\right] =$$

$$\cdot \left[\begin{array}{c|c|c} \begin{pmatrix} 1 & 0 & 1 \\ 2 & 3 & 4 \end{pmatrix}\begin{pmatrix} 1 & 2 \\ 0 & 3 \\ 1 & 4 \end{pmatrix} & \begin{pmatrix} 1 & 0 & 1 \\ 2 & 3 & 4 \end{pmatrix}\begin{pmatrix} 1 \\ 0 \\ 3 \end{pmatrix} & \begin{pmatrix} 1 & 0 & 1 \\ 2 & 3 & 4 \end{pmatrix}\begin{bmatrix} 1 & 5 & 2 & 3 & 1 \\ 1 & 1 & 0 & 1 & 0 \\ 0 & 1 & 0 & 1 & 0 \end{bmatrix} \\ \hline \begin{pmatrix} 1 & 0 & 3 \end{pmatrix}\begin{bmatrix} 1 & 2 \\ 0 & 3 \\ 1 & 4 \end{bmatrix} & \begin{pmatrix} 1 & 0 & 3 \end{pmatrix}\begin{pmatrix} 1 \\ 0 \\ 3 \end{pmatrix} & \begin{pmatrix} 1 & 0 & 3 \end{pmatrix}\begin{bmatrix} 1 & 5 & 2 & 3 & 1 \\ 1 & 1 & 0 & 1 & 0 \\ 0 & 1 & 0 & 1 & 0 \end{bmatrix} \\ \hline \begin{bmatrix} 1 & 1 & 0 \\ 5 & 1 & 1 \\ 2 & 0 & 0 \\ 3 & 1 & 1 \\ 1 & 0 & 0 \end{bmatrix}\begin{bmatrix} 1 & 2 \\ 0 & 3 \\ 1 & 4 \end{bmatrix} & \begin{bmatrix} 1 & 1 & 0 \\ 5 & 1 & 1 \\ 2 & 0 & 0 \\ 3 & 1 & 1 \\ 1 & 0 & 0 \end{bmatrix}\begin{pmatrix} 1 \\ 0 \\ 3 \end{pmatrix} & \begin{bmatrix} 1 & 1 & 0 \\ 5 & 1 & 1 \\ 2 & 0 & 0 \\ 3 & 1 & 1 \\ 1 & 0 & 0 \end{bmatrix}\begin{bmatrix} 1 & 5 & 2 & 3 & 1 \\ 1 & 1 & 0 & 1 & 0 \\ 0 & 1 & 0 & 1 & 0 \end{bmatrix} \end{array}\right]$$

$$\left[\begin{array}{c|c|c}\begin{bmatrix}1&2&3&4\\1&0&1&1\\3&6&1&0\\1&1&0&1\end{bmatrix}\begin{bmatrix}1&1&3&1\\2&0&6&1\\3&1&1&0\\4&1&0&1\end{bmatrix} & \begin{bmatrix}1&2&3&4\\1&0&1&1\\3&6&1&0\\1&1&0&1\end{bmatrix}\begin{bmatrix}2&0&3\\1&0&1\\3&1&2\\5&2&1\end{bmatrix} & \begin{bmatrix}1&2&3&4\\1&0&1&1\\3&6&1&0\\1&1&0&1\end{bmatrix}\begin{bmatrix}1\\1\\0\\1\end{bmatrix}\\ \hline \begin{bmatrix}2&1&3&5\\0&0&1&2\\3&1&2&1\end{bmatrix}\begin{bmatrix}1&1&3&1\\2&0&6&1\\3&1&1&0\\4&1&0&1\end{bmatrix} & \begin{bmatrix}2&1&3&5\\0&0&1&2\\3&1&2&1\end{bmatrix}\begin{bmatrix}2&0&3\\1&0&1\\3&1&2\\5&2&1\end{bmatrix} & \begin{bmatrix}2&1&3&5\\0&0&1&2\\3&1&2&1\end{bmatrix}\begin{bmatrix}1\\1\\0\\1\end{bmatrix}\\ \hline \begin{bmatrix}1&1&0&1\end{bmatrix}\begin{bmatrix}1&1&3&1\\2&0&6&1\\3&1&1&0\\4&1&0&1\end{bmatrix} & (1\ 1\ 0\ 1)\begin{bmatrix}2&0&3\\1&0&1\\3&1&2\\5&2&1\end{bmatrix} & (1\ 1\ 0\ 1)\begin{bmatrix}1\\1\\0\\1\end{bmatrix}\end{array}\right]$$

$$=\left[\begin{array}{cc|c|cccc|c}2&6&4&1&6&2&4&1\\6&29&14&5&17&4&13&2\\\hline 4&14&10&1&8&2&6&1\\\hline 1&5&1&2&6&2&4&1\\6&14&8&6&27&10&17&5\\2&4&2&2&10&4&6&2\\4&13&6&4&17&6&11&3\\1&2&1&1&5&2&3&1\end{array}\right]\cup$$

$$\left[\begin{array}{cccc|ccc|c}30&8&18&7&33&11&15&7\\8&3&4&2&10&3&6&2\\18&4&46&9&15&1&17&9\\7&2&9&3&8&2&5&3\\\hline 33&10&15&8&39&13&18&8\\11&3&1&2&13&5&4&2\\15&6&17&5&18&4&15&5\\\hline 7&2&9&3&8&2&5&3\end{array}\right]$$

$$\begin{aligned}
&= \quad S_1 \cup S_2 \\
&= \quad A^T A = S.
\end{aligned}$$

This S is a symmetric superbimatrix.

**THEOREM 2.1:** *Let $A = A_1 \cup A_2$ be a column superbivector. Then $AA^T$ is a symmetric superbimatrix.*

Proof is left as an exercise for the reader.
Hint: Let

$$A = \begin{bmatrix} A_1^1 \\ \hline A_2^1 \\ \hline \vdots \\ \hline A_{n_1}^1 \end{bmatrix} \cup \begin{bmatrix} A_1^2 \\ \hline A_2^2 \\ \hline \vdots \\ \hline A_{n_2}^2 \end{bmatrix} = A_1 \cup A_2$$

be the given column superbivector. Now

$$\begin{aligned}
A^T &= (A_1 \cup A_2)^T \\
&= A_1^T \cup A_2^T \\
&= \begin{bmatrix} A_1^1 \mid A_2^1 \mid \cdots \mid A_{n_1}^1 \end{bmatrix} \cup \begin{bmatrix} A_1^2 \mid A_2^2 \mid \cdots \mid A_{n_2}^2 \end{bmatrix}
\end{aligned}$$

$$\begin{aligned}
AA^T &= (A_1 \cup A_2)(A_1 \cup A_2)^T \\
&= (A_1 \cup A_2)(A_1^T \cup A_2^T) \\
&= A_1 A_1^T \cup A_2 A_2^T \\
&= \begin{bmatrix} A_1^1 \\ \hline A_2^1 \\ \hline \vdots \\ \hline A_{n_1}^1 \end{bmatrix} \begin{bmatrix} A_1^{1T} \mid A_2^{1T} \mid A_3^{1T} \mid \cdots \mid A_{n_1}^{1T} \end{bmatrix} \cup \\
&\quad \begin{bmatrix} A_1^2 \\ \hline A_2^2 \\ \hline \vdots \\ \hline A_{n_2}^2 \end{bmatrix} \begin{bmatrix} A_1^{2T} \mid A_2^{2T} \mid A_3^{2T} \mid \cdots \mid A_{n_2}^{2T} \end{bmatrix}
\end{aligned}$$

$$= \begin{bmatrix} A_1^1 A_1^{tT} & A_1^1 A_2^{tT} & \cdots & A_1^1 A_{n_1}^{tT} \\ \hline A_2^1 A_2^{tT} & A_2^1 A_2^{tT} & \vdots & A_2^1 A_{n_1}^{tT} \\ \hline \vdots & \vdots & & \vdots \\ \hline A_{n_1}^1 A_1^{tT} & A_{n_1}^1 A_2^{tT} & \cdots & A_{n_1}^1 A_{n_1}^{tT} \end{bmatrix} \cup$$

$$\begin{bmatrix} A_1^2 A_1^{2T} & A_1^2 A_2^{2T} & \cdots & A_1^2 A_{n_2}^{2T} \\ \hline A_2^2 A_1^{2T} & A_2^2 A_2^{2T} & \vdots & A_2^2 A_{n_2}^{2T} \\ \hline \vdots & \vdots & & \vdots \\ \hline A_{n_2}^2 A_1^{2T} & A_{n_2}^2 A_2^{2T} & \cdots & A_{n_2}^2 A_{n_2}^{2T} \end{bmatrix}.$$

We see $\left(A_1^1 A_i^{1T}\right)^T = A_i^1 A_1^{1T}$; i = 1, 2, …, $n_1$. Also $\left(A_j^2 A_k^{2T}\right)^T = A_k^2 A_j^{2T}$, $1 \le k, j \le n_2$. Thus we easily see the product $AA^T$ gives a symmetric superbimatrix.

Now we proceed onto define the minor product of semi superbimatrix, to this end first we define semi superbivector.

**DEFINITION 2.18:** *Let $A = A_1 \cup A_2$ be any semi supermatrix. We say A is a semi superbivector if the supermatrix $A_1$ (or $A_2$) is just partitioned only vertically or horizontally 'or' in the mutually exclusive sense. The other matrix may be a square matrix or a rectangular matrix or a column vector or row vector.*

*Example 2.47:* Let $A = A_1 \cup A_2$ where

$$A_1 = \begin{bmatrix} 3 & 2 & 1 & 0 & 5 & 1 \\ 1 & 2 & 0 & 5 & 6 & 3 \end{bmatrix}$$

and

$$A_2 = \left[\begin{array}{c|cccc|ccc} 3 & 1 & 1 & 0 & -1 & 3 & 1 & 1 \\ 1 & 1 & 0 & 1 & 2 & 1 & 3 & 0 \\ 0 & 0 & 1 & 1 & 3 & 2 & 5 & 1 \end{array}\right].$$

A is a semi superbimatrix, because $A_2$ happens to be a row super vector we call A to be a row semi superbimatrix.

**Note:** Even if $A_1$ is not a rectangular matrix still we call A to be only a row semi superbimatrix.

***Example 2.48:*** Let $B = B_1 \cup B_2$ where

$$B_1 = \begin{bmatrix} 2 & 0 & 1 & 1 \\ 1 & 2 & 0 & 1 \\ 1 & 1 & 0 & 2 \\ 1 & 0 & 1 & 2 \end{bmatrix}$$

and

$$B_2 = \begin{bmatrix} 3 & 1 & 5 & 6 \\ 1 & 0 & 2 & 1 \\ \hline 7 & 6 & 5 & 4 \\ 3 & 2 & 1 & 0 \\ 1 & 1 & 2 & 1 \\ \hline 0 & 7 & 2 & 5 \end{bmatrix};$$

then we call B to be a column semi superbimatrix though $B_1$ is just a 4 × 4 square matrix.

Now we call these row semi superbimatrix and column semi superbimatrix as semi superbivectors, even if the non super component is a square matrix or a column vector or a row vector or a rectangular matrix.

Now we illustrate the minor product of semi superbivector.

***Example 2.49:*** Let $A = A_1 \cup A_2$ and $B = B_1 \cup B_2$ be two semi superbimatrices; where

$$A_1 = \begin{bmatrix} 3 & 0 & 1 & 2 \\ 1 & 1 & 0 & 1 \\ 5 & 0 & 1 & 3 \end{bmatrix}$$

and

$$A_2 = \begin{bmatrix} 3 & 1 & 2 & 5 & 3 & 1 & 4 & 3 \\ 1 & 0 & 1 & 1 & 0 & 1 & 1 & 1 \\ 0 & 1 & 0 & 2 & 1 & 0 & 0 & 2 \end{bmatrix}.$$

A is the row semi superbivector and $B = B_1 \cup B_2$ with

$$B_1 = \begin{bmatrix} 3 & 1 \\ 0 & 6 \\ 1 & 1 \\ 2 & 0 \end{bmatrix}$$

and

$$B_2 = \begin{bmatrix} 2 & 1 \\ 3 & 0 \\ 1 & 1 \\ \hline 1 & 2 \\ 2 & 1 \\ 1 & 0 \\ 0 & 1 \\ \hline 1 & 5 \end{bmatrix}$$

be the column semi superbivector. Now

AB  = $(A_1 \cup A_2)(B_1 \cup B_2)$
    = $A_1 B_1 \cup A_2 B_2$

be the minor product of semi superbivectors.

$$A_1 B_1 \cup A_2 B_2 = \begin{bmatrix} 3 & 0 & 1 & 2 \\ 1 & 1 & 0 & 1 \\ 5 & 0 & 1 & 3 \end{bmatrix} \begin{bmatrix} 3 & 1 \\ 0 & 6 \\ 1 & 1 \\ 2 & 0 \end{bmatrix} \cup$$

$$\begin{bmatrix} 3 & 1 & 2 & | & 5 & 3 & 1 & 4 & | & 3 \\ 1 & 0 & 1 & | & 1 & 0 & 1 & 1 & | & 1 \\ 0 & 1 & 0 & | & 2 & 1 & 0 & 0 & | & 2 \end{bmatrix} \begin{bmatrix} 2 & 1 \\ 3 & 0 \\ 1 & 1 \\ \hline 1 & 2 \\ 2 & 1 \\ 1 & 0 \\ 0 & 1 \\ \hline 1 & 5 \end{bmatrix}$$

$$= \begin{bmatrix} 14 & 4 \\ 5 & 7 \\ 22 & 6 \end{bmatrix} \cup$$

$$\left\{ \begin{bmatrix} 3 & 1 & 2 \\ 1 & 0 & 1 \\ 0 & 1 & 0 \end{bmatrix} \begin{bmatrix} 2 & 1 \\ 3 & 0 \\ 1 & 1 \end{bmatrix} + \begin{bmatrix} 5 & 3 & 1 & 4 \\ 1 & 0 & 1 & 1 \\ 2 & 1 & 0 & 0 \end{bmatrix} \begin{bmatrix} 1 & 2 \\ 2 & 1 \\ 1 & 0 \\ 0 & 1 \end{bmatrix} + \begin{bmatrix} 3 \\ 1 \\ 2 \end{bmatrix} \begin{bmatrix} 1 & 5 \end{bmatrix} \right\}$$

$$= \begin{bmatrix} 14 & 4 \\ 5 & 7 \\ 22 & 6 \end{bmatrix} \cup \left\{ \begin{bmatrix} 11 & 5 \\ 3 & 2 \\ 3 & 0 \end{bmatrix} + \begin{bmatrix} 12 & 17 \\ 2 & 2 \\ 4 & 5 \end{bmatrix} + \begin{bmatrix} 3 & 15 \\ 1 & 5 \\ 2 & 10 \end{bmatrix} \right\}.$$

$$= \begin{bmatrix} 14 & 4 \\ 5 & 7 \\ 22 & 6 \end{bmatrix} \cup \begin{bmatrix} 26 & 37 \\ 6 & 9 \\ 9 & 15 \end{bmatrix}.$$

*Example 2.50:* Let $A = A_1 \cup A_2$ where

$$A_1 = [1\ 2\ 3\ 0\ 1\ 5\ 6] \cup \begin{bmatrix} 4 & 2 & | & 3 & 5 & 3 & | & 5 & 0 & 7 & 1 \\ 5 & 3 & | & 2 & 6 & 3 & | & 1 & 0 & 1 & 1 \\ 1 & 0 & | & 1 & 1 & 4 & | & 0 & 1 & 0 & 1 \\ 0 & 1 & | & 1 & 0 & 5 & | & 2 & 0 & 1 & 0 \end{bmatrix}$$

be a semi superbivector and $B = B_1 \cup B_2$ with

$$B_1 = \begin{bmatrix} 0 \\ 1 \\ 2 \\ 3 \\ 0 \\ 1 \\ 5 \end{bmatrix}$$

and

$$B_2 = \begin{bmatrix} 0 & 5 & 1 & 2 \\ 1 & 2 & 0 & 2 \\ \hline 1 & 0 & 1 & 0 \\ 0 & 1 & 0 & 1 \\ 1 & 1 & 1 & 1 \\ \hline 0 & -1 & 0 & -1 \\ 1 & 0 & 1 & 0 \\ -1 & 1 & -1 & 1 \\ 0 & 1 & 1 & 0 \end{bmatrix}$$

be a column semi superbivector. Now

$$\begin{aligned}
AB &= (A_1 \cup A_2)(B_1 \cup B_2) \\
&= A_1 A_2 \cup B_1 B_2 \\
&= [1\ 2\ 3\ 0\ 1\ 5\ 6] \begin{bmatrix} 0 \\ 1 \\ 2 \\ 3 \\ 0 \\ 1 \\ 5 \end{bmatrix} \cup
\end{aligned}$$

$$\left\{ \begin{bmatrix} 4 & 2 & | & 3 & 5 & 3 & | & 5 & 0 & 7 & 1 \\ 5 & 3 & | & 2 & 6 & 3 & | & 1 & 0 & 1 & 1 \\ 1 & 0 & | & 1 & 1 & 4 & | & 0 & 1 & 0 & 1 \\ 0 & 1 & | & 1 & 0 & 5 & | & 2 & 0 & 1 & 0 \end{bmatrix} \times \begin{bmatrix} 0 & 5 & 1 & 2 \\ 1 & 2 & 0 & 2 \\ \hline 1 & 0 & 1 & 0 \\ 0 & 1 & 0 & 1 \\ 1 & 1 & 1 & 1 \\ \hline 0 & -1 & 0 & -1 \\ 1 & 0 & 1 & 0 \\ -1 & 0 & -1 & 1 \\ 0 & 1 & 1 & 0 \end{bmatrix} \right\}$$

$$= [4 \ 3] \cup \left\{ \begin{bmatrix} 4 & 2 \\ 5 & 3 \\ 1 & 0 \\ 0 & 1 \end{bmatrix} \begin{bmatrix} 0 & 5 & 1 & 2 \\ 1 & 2 & 0 & 2 \end{bmatrix} + \begin{bmatrix} 3 & 5 & 3 \\ 2 & 6 & 3 \\ 1 & 1 & 4 \\ 1 & 0 & 5 \end{bmatrix} \begin{bmatrix} 1 & 0 & 1 & 0 \\ 0 & 1 & 0 & 1 \\ 1 & 1 & 1 & 1 \end{bmatrix} \right.$$

$$\left. + \begin{bmatrix} 5 & 0 & 7 & 1 \\ 1 & 0 & 1 & 1 \\ 0 & 1 & 0 & 1 \\ 2 & 0 & 1 & 0 \end{bmatrix} \begin{bmatrix} 0 & -1 & 0 & -1 \\ 1 & 0 & 1 & 0 \\ -1 & 1 & -1 & 1 \\ 0 & 1 & 1 & 0 \end{bmatrix} \right\} =$$

$$[4 \ 3] \cup \left\{ \begin{bmatrix} 2 & 24 & 4 & 12 \\ 3 & 31 & 5 & 16 \\ 0 & 5 & 1 & 2 \\ 1 & 2 & 0 & 2 \end{bmatrix} + \begin{bmatrix} 6 & 8 & 6 & 8 \\ 5 & 9 & 5 & 9 \\ 5 & 5 & 5 & 5 \\ 6 & 5 & 6 & 5 \end{bmatrix} + \begin{bmatrix} -7 & -3 & -6 & 2 \\ -1 & 1 & 0 & 0 \\ 1 & 1 & 2 & 0 \\ -1 & -1 & -1 & -1 \end{bmatrix} \right\}$$

$$= [4 \ 3] \cup \begin{bmatrix} 1 & 29 & 4 & 22 \\ 7 & 39 & 10 & 25 \\ 6 & 11 & 8 & 7 \\ 6 & 6 & 5 & 6 \end{bmatrix}.$$

Clearly AB is only a usual bimatrix and is not a semi superbimatrix. Thus we see the minor product yields only a bimatrix, the semi super quality is lost by this product.

Now we proceed on to define the minor product of semi superbivectors.

**DEFINITION 2.19:** *Let $A = A_1 \cup A_2$ and $B = B_1 \cup B_2$ be two semi superbivectors. The minor product of AB is defined if and only if in the product $AB = (A_1 \cup A_2)(B_1 \cup B_2) = A_1B_1 \cup A_2B_2$, the usual matrix product $A_1B_1$ of the $A_1$ and $B_1$ is defined i.e., if $A_1$ is a $m \times n$ matrix then $B_1$ must be a $n \times t$ matrix and $A_2B_2$ is defined only if $A_2$ is a row super vector say*
$$A_2 = \left[ A_1^2 \mid A_2^2 \mid \cdots \mid A_{n_2}^2 \right]$$
*and $B_2$ is a column super vector such that if*
$$B_2 = \begin{bmatrix} B_1^2 \\ \hline B_2^2 \\ \hline \vdots \\ \hline B_{n_2}^2 \end{bmatrix}$$
*then each of the product of matrices $A_i^2 B_i^2$ is defined for $i = 1, 2, \ldots, n_2$ and all of them are of same order, i.e.,*
$$A_2B_2 = \left[ A_1^2 \mid A_2^2 \mid \cdots \mid A_{n_2}^2 \right] \begin{bmatrix} B_1^2 \\ \hline B_2^2 \\ \hline \vdots \\ \hline B_{n_2}^2 \end{bmatrix} = A_1^2 B_1^2 + A_2^2 B_2^2 + \cdots + A_{n_2}^2 B_{n_2}^2$$
*since each $A_i^2 B_i^2$ is only a ordinary simple matrix and not a supermatrix we see $A_2B_2 = \sum_{i=1}^{n_2} A_i^2 B_i^2$ is just a single $m_2 \times n_2$ matrix. This us see $AB = A_1B_1 \cup A_2B_2$ is only a bimatrix.*

Now we just define major product of semi superbivectors. We first illustrate it by examples before we go for the abstract definition.

*Example 2.51:* Let $A = A_1 \cup A_2$ and $B = B_1 \cup B_2$ be two semi superbivectors where $A = A_1 \cup A_2$ with

$$A_1 = \begin{bmatrix} 3 & 1 & 0 & 5 & 7 & 2 & 1 & 0 \\ 0 & 1 & 0 & 1 & 2 & 3 & 1 & 6 \\ -1 & 1 & 1 & 0 & 0 & 1 & 0 & 1 \end{bmatrix}$$

and

$$A_2 = \begin{bmatrix} 3 & 1 & 0 & 1 \\ 1 & 1 & 1 & 0 \\ \hline 1 & 1 & 0 & 1 \\ 2 & 1 & 2 & 0 \\ 1 & 0 & 1 & 0 \\ \hline 1 & 2 & 3 & 1 \\ 1 & 1 & 1 & 0 \\ 0 & 1 & 1 & 0 \\ 1 & 0 & 0 & 1 \end{bmatrix}$$

be a column semi superbivector and $B = B_1 \cup B_2$ with

$$B_1 = \begin{bmatrix} 1 & 1 & 0 \\ 2 & 0 & 1 \\ 3 & 1 & 0 \\ 4 & 0 & 1 \\ 5 & 1 & 0 \\ 6 & 0 & 1 \\ 7 & 1 & 0 \\ 8 & 0 & 1 \end{bmatrix}$$

and

$$B_2 = \begin{bmatrix} 1 & 2 & 1 & 4 & 3 & 0 & 1 & 1 & 1 & 1 \\ 0 & 3 & 0 & 1 & 1 & 1 & 0 & 1 & 0 & 1 \\ 1 & 1 & 1 & 2 & 3 & 0 & 1 & 0 & 1 & 0 \\ 2 & 0 & 1 & 0 & 1 & 1 & 0 & 0 & 1 & 1 \end{bmatrix}$$

$B_2$ be the row superbivector so that B is a row semi superbivector. Now AB the major product of the two semi superbivectors is defined as follows

$$AB = (A_1 \cup A_2)(B_1 \cup B_2)$$
$$= A_1 B_1 \cup A_2 B_2$$

$$\begin{bmatrix} 3 & 1 & 0 & 5 & 7 & 2 & 1 & 0 \\ 0 & 1 & 0 & 1 & 2 & 3 & 1 & 6 \\ -1 & 1 & 1 & 0 & 0 & 1 & 0 & 1 \end{bmatrix} \begin{bmatrix} 1 & 1 & 0 \\ 2 & 0 & 1 \\ 3 & 1 & 0 \\ 4 & 0 & 1 \\ 5 & 1 & 0 \\ 6 & 0 & 1 \\ 7 & 1 & 0 \\ 8 & 0 & 1 \end{bmatrix} \cup$$

$$\begin{bmatrix} 3 & 1 & 0 & 1 \\ 1 & 1 & 1 & 0 \\ 1 & 1 & 0 & 2 \\ 2 & 1 & 2 & 0 \\ 1 & 0 & 1 & 0 \\ 1 & 2 & 3 & 1 \\ 1 & 1 & 1 & 0 \\ 0 & 1 & 1 & 0 \\ 1 & 0 & 0 & 1 \end{bmatrix} \begin{bmatrix} 1 & 2 & 1 & 4 & 3 & 0 & 1 & 1 & 1 & 1 \\ 0 & 3 & 0 & 1 & 1 & 1 & 0 & 1 & 0 & 1 \\ 1 & 1 & 1 & 2 & 3 & 0 & 1 & 0 & 1 & 0 \\ 2 & 0 & 1 & 0 & 1 & 1 & 0 & 0 & 1 & 1 \end{bmatrix}$$

$$= \begin{bmatrix} 79 & 11 & 8 \\ 89 & 3 & 11 \\ 18 & 0 & 3 \end{bmatrix} \cup$$

$$\left[ \begin{array}{c|c|c} \begin{pmatrix} 3 & 1 & 0 & 1 \\ 1 & 1 & 1 & 0 \end{pmatrix} \begin{bmatrix} 1 \\ 0 \\ 1 \\ 2 \end{bmatrix} & \begin{pmatrix} 3 & 1 & 0 & 1 \\ 1 & 1 & 1 & 0 \end{pmatrix} \begin{bmatrix} 2 & 1 \\ 3 & 0 \\ 1 & 1 \\ 0 & 1 \end{bmatrix} & \begin{pmatrix} 3 & 1 & 0 & 1 \\ 1 & 1 & 1 & 0 \end{pmatrix} \begin{bmatrix} 4 & 3 & 0 \\ 1 & 1 & 1 \\ 2 & 3 & 0 \\ 0 & 1 & 1 \end{bmatrix} \\ \hline \begin{pmatrix} 1 & 1 & 0 & 2 \\ 2 & 1 & 2 & 0 \\ 1 & 0 & 1 & 0 \end{pmatrix} \begin{bmatrix} 1 \\ 0 \\ 1 \\ 2 \end{bmatrix} & \begin{pmatrix} 1 & 1 & 0 & 2 \\ 2 & 1 & 2 & 0 \\ 1 & 0 & 1 & 0 \end{pmatrix} \begin{bmatrix} 2 & 1 \\ 3 & 0 \\ 1 & 1 \\ 0 & 1 \end{bmatrix} & \begin{pmatrix} 1 & 1 & 0 & 2 \\ 2 & 1 & 2 & 0 \\ 1 & 0 & 1 & 0 \end{pmatrix} \begin{bmatrix} 4 & 3 & 0 \\ 1 & 1 & 1 \\ 2 & 3 & 0 \\ 0 & 1 & 1 \end{bmatrix} \\ \hline \begin{pmatrix} 1 & 2 & 3 & 1 \\ 1 & 1 & 1 & 0 \\ 0 & 1 & 1 & 0 \\ 1 & 0 & 0 & 1 \end{pmatrix} \begin{bmatrix} 1 \\ 0 \\ 1 \\ 2 \end{bmatrix} & \begin{pmatrix} 1 & 2 & 3 & 1 \\ 1 & 1 & 1 & 0 \\ 0 & 1 & 1 & 0 \\ 1 & 0 & 0 & 1 \end{pmatrix} \begin{bmatrix} 2 & 1 \\ 3 & 0 \\ 1 & 1 \\ 0 & 1 \end{bmatrix} & \begin{pmatrix} 1 & 2 & 3 & 1 \\ 1 & 1 & 1 & 0 \\ 0 & 1 & 1 & 0 \\ 1 & 0 & 0 & 1 \end{pmatrix} \begin{bmatrix} 4 & 3 & 0 \\ 1 & 1 & 1 \\ 2 & 3 & 0 \\ 0 & 1 & 1 \end{bmatrix} \end{array} \right]$$

$$\left[ \begin{array}{c} \begin{pmatrix} 3 & 1 & 0 & 1 \\ 1 & 1 & 1 & 0 \end{pmatrix} \begin{bmatrix} 1 & 1 & 1 & 1 \\ 0 & 1 & 0 & 1 \\ 1 & 0 & 1 & 0 \\ 0 & 0 & 1 & 1 \end{bmatrix} \\ \hline \begin{pmatrix} 1 & 1 & 0 & 2 \\ 2 & 1 & 2 & 0 \\ 1 & 0 & 1 & 0 \end{pmatrix} \begin{bmatrix} 1 & 1 & 1 & 1 \\ 0 & 1 & 0 & 1 \\ 1 & 0 & 1 & 0 \\ 0 & 0 & 1 & 1 \end{bmatrix} \\ \hline \begin{pmatrix} 1 & 2 & 3 & 1 \\ 1 & 1 & 1 & 0 \\ 0 & 1 & 1 & 0 \\ 1 & 0 & 0 & 1 \end{pmatrix} \begin{bmatrix} 1 & 1 & 1 & 1 \\ 0 & 1 & 0 & 1 \\ 1 & 0 & 1 & 0 \\ 0 & 0 & 1 & 1 \end{bmatrix} \end{array} \right]$$

$$= \begin{bmatrix} 79 & 11 & 8 \\ 89 & 3 & 11 \\ 18 & 0 & 3 \end{bmatrix} \cup$$

$$\begin{bmatrix} 5 & 9 & 4 & 13 & 11 & 2 & 3 & 4 & 4 & 5 \\ 2 & 6 & 2 & 7 & 7 & 1 & 2 & 2 & 2 & 2 \\ \hline 5 & 5 & 3 & 5 & 6 & 3 & 1 & 2 & 3 & 4 \\ 4 & 9 & 4 & 13 & 13 & 1 & 4 & 3 & 4 & 3 \\ 2 & 3 & 2 & 6 & 6 & 0 & 2 & 1 & 2 & 1 \\ \hline 6 & 11 & 4 & 12 & 14 & 3 & 4 & 3 & 5 & 4 \\ 2 & 6 & 2 & 7 & 7 & 1 & 2 & 2 & 2 & 2 \\ 1 & 4 & 1 & 3 & 4 & 1 & 1 & 1 & 1 & 1 \\ 3 & 2 & 2 & 4 & 4 & 1 & 1 & 1 & 2 & 2 \end{bmatrix}.$$

Clearly the major product yield a semi superbimatrix which is not a row or column semi superbivector.

Now we give yet another example before we proceed to give the abstract definition.

***Example 2.52:*** Let $A = A_1 \cup A_2$ and $B = B_1 \cup B_2$ be any two semi superbivectors given by

$$A_1 = \begin{bmatrix} 1 & 2 & 0 & 4 & 5 \\ 0 & 3 & 1 & 1 & 0 \\ 1 & 2 & 0 & 1 & 0 \\ 3 & 1 & 0 & 0 & 1 \end{bmatrix} \text{ and } A_2 = \begin{bmatrix} 2 & 3 & 1 \\ 1 & 2 & 0 \\ \hline 1 & 1 & 0 \\ 3 & 1 & 5 \\ 5 & 1 & 7 \\ \hline 1 & 1 & 0 \\ 0 & 1 & 1 \\ 1 & 0 & 1 \\ 0 & 1 & 0 \\ 1 & 0 & 0 \\ 0 & 0 & 1 \\ 1 & 1 & 1 \end{bmatrix}.$$

A be the column semi superbivector.
$$B = B_1 \cup B_2$$

$$= \begin{bmatrix} 1 & 0 \\ 2 & 1 \\ 3 & 2 \\ 4 & 1 \\ 5 & 0 \end{bmatrix} \cup \begin{bmatrix} 1 & 1 & 0 & | & 0 & 1 & 1 & 1 & -1 & | & 0 & 2 & 4 \\ 0 & 2 & 1 & | & 1 & 0 & 1 & 1 & 0 & | & 1 & 0 & 1 \\ 2 & 0 & 2 & | & 0 & 1 & 0 & 1 & 1 & | & 2 & 1 & 1 \end{bmatrix}$$

be the row semi superbivector. $AB = (A_1 \cup A_2)(B_1 \cup B_2) = A_1B_1 \cup A_2B_2$ where we define the major product of the semi superbivectors. $AB = A_1B_1 \cup A_2B_2$

$$= \begin{bmatrix} 1 & 2 & 0 & 4 & 5 \\ 0 & 3 & 1 & 1 & 0 \\ 1 & 2 & 0 & 1 & 0 \\ 3 & 1 & 0 & 0 & 1 \end{bmatrix} \begin{bmatrix} 1 & 0 \\ 2 & 1 \\ 3 & 2 \\ 4 & 1 \\ 5 & 0 \end{bmatrix} \cup$$

$$\begin{bmatrix} 2 & 3 & 1 \\ 1 & 2 & 0 \\ \hline 1 & 1 & 0 \\ 3 & 1 & 5 \\ 5 & 1 & 7 \\ \hline 1 & 1 & 0 \\ 0 & 1 & 1 \\ 1 & 0 & 1 \\ 0 & 1 & 0 \\ 1 & 0 & 0 \\ 0 & 0 & 1 \\ 1 & 1 & 1 \end{bmatrix} \begin{bmatrix} 1 & 1 & 0 & | & 0 & 1 & 1 & 1 & -1 & | & 0 & 2 & 4 \\ 0 & 2 & 1 & | & 1 & 0 & 1 & 1 & 0 & | & 1 & 0 & 1 \\ 2 & 0 & 2 & | & 0 & 1 & 0 & 1 & 1 & | & 2 & 1 & 1 \end{bmatrix}$$

$$= \begin{bmatrix} 46 & 6 \\ 13 & 6 \\ 9 & 3 \\ 10 & 1 \end{bmatrix} \cup$$

$$
\begin{bmatrix}
\begin{bmatrix} 2 & 3 & 1 \\ 1 & 2 & 0 \end{bmatrix}\begin{bmatrix} 1 & 1 & 0 \\ 0 & 2 & 1 \\ 2 & 0 & 2 \end{bmatrix} & \begin{bmatrix} 2 & 3 & 1 \\ 1 & 2 & 0 \end{bmatrix}\begin{bmatrix} 0 & 1 & 1 & 1 & -1 \\ 1 & 0 & 1 & 1 & 0 \\ 0 & 1 & 0 & 1 & 1 \end{bmatrix} & \begin{bmatrix} 2 & 3 & 1 \\ 1 & 2 & 0 \end{bmatrix}\begin{bmatrix} 0 & 2 & 4 \\ 1 & 0 & 1 \\ 2 & 1 & 1 \end{bmatrix} \\
\hline
\begin{bmatrix} 1 & 1 & 0 \\ 3 & 1 & 5 \\ 5 & 1 & 7 \end{bmatrix}\begin{bmatrix} 1 & 1 & 0 \\ 0 & 2 & 1 \\ 2 & 0 & 2 \end{bmatrix} & \begin{bmatrix} 1 & 1 & 0 \\ 3 & 1 & 5 \\ 5 & 1 & 7 \end{bmatrix}\begin{bmatrix} 0 & 1 & 1 & 1 & -1 \\ 1 & 0 & 1 & 1 & 0 \\ 0 & 1 & 0 & 1 & 1 \end{bmatrix} & \begin{bmatrix} 1 & 1 & 0 \\ 3 & 1 & 5 \\ 5 & 1 & 7 \end{bmatrix}\begin{bmatrix} 0 & 2 & 4 \\ 1 & 0 & 1 \\ 2 & 1 & 1 \end{bmatrix} \\
\hline
\begin{bmatrix} 1 & 1 & 0 \\ 0 & 1 & 1 \\ 1 & 0 & 1 \\ 0 & 1 & 0 \\ 1 & 0 & 0 \\ 0 & 0 & 1 \\ 1 & 1 & 1 \end{bmatrix}\begin{bmatrix} 1 & 1 & 0 \\ 0 & 2 & 1 \\ 2 & 0 & 2 \end{bmatrix} & \begin{bmatrix} 1 & 1 & 0 \\ 0 & 1 & 1 \\ 1 & 0 & 1 \\ 0 & 1 & 0 \\ 1 & 0 & 0 \\ 0 & 0 & 1 \\ 1 & 1 & 1 \end{bmatrix}\begin{bmatrix} 0 & 1 & 1 & 1 & -1 \\ 1 & 0 & 1 & 1 & 0 \\ 0 & 1 & 0 & 1 & 1 \end{bmatrix} & \begin{bmatrix} 1 & 1 & 0 \\ 0 & 1 & 1 \\ 1 & 0 & 1 \\ 0 & 1 & 0 \\ 1 & 0 & 0 \\ 0 & 0 & 1 \\ 1 & 1 & 1 \end{bmatrix}\begin{bmatrix} 0 & 2 & 4 \\ 1 & 0 & 1 \\ 2 & 1 & 1 \end{bmatrix}
\end{bmatrix}
$$

$$
= \begin{bmatrix} 46 & 6 \\ 13 & 6 \\ 9 & 3 \\ 10 & 1 \end{bmatrix} \cup \left[\begin{array}{ccc|cccc|ccc}
4 & 8 & 5 & 3 & 3 & 5 & 6 & -1 & 5 & 5 & 12 \\
1 & 5 & 2 & 2 & 1 & 3 & 3 & -1 & 2 & 2 & 6 \\
\hline
1 & 3 & 1 & 1 & 1 & 2 & 2 & -1 & 1 & 2 & 5 \\
13 & 5 & 11 & 1 & 8 & 4 & 9 & 2 & 11 & 11 & 18 \\
19 & 7 & 15 & 1 & 12 & 6 & 13 & 2 & 15 & 17 & 28 \\
\hline
1 & 3 & 1 & 1 & 1 & 2 & 2 & -1 & 1 & 2 & 5 \\
2 & 2 & 3 & 1 & 1 & 1 & 2 & 0 & 3 & 1 & 2 \\
3 & 1 & 2 & 0 & 2 & 1 & 2 & 0 & 2 & 3 & 5 \\
0 & 2 & 1 & 1 & 0 & 1 & 1 & 0 & 1 & 0 & 1 \\
1 & 1 & 0 & 0 & 1 & 1 & 1 & -1 & 0 & 2 & 4 \\
2 & 0 & 2 & 0 & 1 & 0 & 1 & 1 & 2 & 1 & 1 \\
3 & 3 & 3 & 1 & 2 & 2 & 3 & 0 & 3 & 3 & 6
\end{array}\right]
$$

is a semi superbimatrix.

This major product converts product of semi superbivectors in to semi superbimatrix where as minor product makes the product of semi superbivectors into just a bimatrix.

**DEFINITION 2.20:** *Let $A = A_1 \cup A_2$ and $B = B_1 \cup B_2$ be two semi superbivectors. The major product of the two semi superbivectors A and B is AB, is defined as $AB = (A_1 \cup A_2)(B_1 \cup B_2) = A_1B_1 \cup A_2B_2$ if*

1. *$A_1B_1$ must be compatible with respect to usual matrix product that is if $A_1$ is a m×n matrix then $B_1$ must be a n×t matrix.*

2. *$A_2B_2$ is defined only if $A_2$ is a super column vector and $B_2$ is a super row vector such that the number of columns in $A_2$ must be equal to the number of rows in $B_2$.*

Now we find the product of A with $A^T$ or $A^T$ with A which ever is compatible where A is a semi superbivector.

*Example 2.53:* Let $A = A_1 \cup A_2$ be a semi superbivector given by

$$A = \begin{bmatrix} 1 & 0 & 1 & 2 & 3 \\ 1 & 1 & 0 & 1 & 2 \\ 3 & 0 & 1 & 0 & 1 \end{bmatrix}$$

$$\cup \begin{bmatrix} 1 & 0 & 2 & 1 & 5 & 2 & 1 & 1 & 0 \\ 2 & 1 & 1 & 0 & 1 & 0 & 0 & 1 & 1 \\ 5 & 3 & 0 & 1 & 0 & 1 & 1 & 0 & 1 \end{bmatrix}.$$

Now

$$A^T = (A_1 \cup A_2)^T$$
$$= A_1^T \cup A_2^T$$

$$= \begin{bmatrix} 1 & 1 & 3 \\ 0 & 1 & 0 \\ 1 & 0 & 1 \\ 2 & 1 & 0 \\ 3 & 2 & 1 \end{bmatrix} \cup \begin{bmatrix} 1 & 2 & 5 \\ 0 & 1 & 3 \\ 2 & 1 & 0 \\ 1 & 0 & 1 \\ \hline 5 & 1 & 0 \\ 2 & 0 & 1 \\ \hline 1 & 0 & 1 \\ 1 & 1 & 0 \\ 0 & 1 & 1 \end{bmatrix}$$

$$\begin{aligned} A^T A &= (A_1 \cup A_2)^T (A_1 \cup A_2) \\ &= (A_1^T \cup A_2^T)(A_1 \cup A_2)^T \\ &= A_1^T A_1 \cup A_2^T A_2 \end{aligned}$$

$$= \begin{bmatrix} 1 & 1 & 3 \\ 0 & 1 & 0 \\ 1 & 0 & 1 \\ 2 & 1 & 0 \\ 3 & 2 & 1 \end{bmatrix} \begin{bmatrix} 1 & 0 & 1 & 2 & 3 \\ 1 & 1 & 0 & 1 & 2 \\ 3 & 0 & 1 & 0 & 1 \end{bmatrix} \cup$$

$$\begin{bmatrix} 1 & 2 & 5 \\ 0 & 1 & 3 \\ 2 & 1 & 0 \\ 1 & 0 & 1 \\ \hline 5 & 1 & 0 \\ 2 & 0 & 1 \\ \hline 1 & 0 & 1 \\ 1 & 1 & 0 \\ 0 & 1 & 1 \end{bmatrix} \begin{bmatrix} 1 & 0 & 2 & 1 & 5 & 2 & 1 & 1 & 0 \\ 2 & 1 & 1 & 0 & 1 & 0 & 0 & 1 & 1 \\ 5 & 3 & 0 & 1 & 0 & 1 & 1 & 0 & 1 \end{bmatrix}$$

$$= \begin{bmatrix} 11 & 1 & 4 & 3 & 8 \\ 1 & 1 & 0 & 1 & 2 \\ 4 & 0 & 2 & 2 & 4 \\ 3 & 1 & 2 & 5 & 8 \\ 8 & 2 & 4 & 8 & 14 \end{bmatrix} \cup$$

$$\begin{bmatrix} \begin{bmatrix} 1 & 2 & 5 \\ 0 & 1 & 3 \\ 2 & 1 & 0 \\ 1 & 0 & 1 \end{bmatrix} \begin{bmatrix} 1 & 0 & 2 & 1 \\ 2 & 1 & 1 & 0 \\ 5 & 3 & 0 & 1 \end{bmatrix} & \begin{bmatrix} 1 & 2 & 5 \\ 0 & 1 & 3 \\ 2 & 1 & 0 \\ 1 & 0 & 1 \end{bmatrix} \begin{bmatrix} 5 & 2 \\ 1 & 0 \\ 0 & 1 \end{bmatrix} & \begin{pmatrix} 1 & 2 & 5 \\ 0 & 1 & 3 \\ 2 & 1 & 0 \\ 1 & 0 & 1 \end{pmatrix} \begin{bmatrix} 1 & 1 & 0 \\ 0 & 1 & 1 \\ 1 & 0 & 1 \end{bmatrix} \\ \hline \begin{pmatrix} 5 & 1 & 0 \\ 2 & 0 & 1 \end{pmatrix} \begin{bmatrix} 1 & 0 & 2 & 1 \\ 2 & 1 & 1 & 0 \\ 5 & 3 & 0 & 1 \end{bmatrix} & \begin{pmatrix} 5 & 1 & 0 \\ 2 & 0 & 1 \end{pmatrix} \begin{pmatrix} 5 & 2 \\ 1 & 0 \\ 0 & 1 \end{pmatrix} & \begin{pmatrix} 5 & 1 & 0 \\ 2 & 0 & 1 \end{pmatrix} \begin{pmatrix} 1 & 1 & 0 \\ 0 & 1 & 1 \\ 1 & 0 & 1 \end{pmatrix} \\ \begin{pmatrix} 1 & 0 & 1 \\ 1 & 1 & 0 \\ 0 & 1 & 1 \end{pmatrix} \begin{pmatrix} 1 & 0 & 2 & 1 \\ 2 & 1 & 1 & 0 \\ 5 & 3 & 0 & 1 \end{pmatrix} & \begin{pmatrix} 1 & 0 & 1 \\ 1 & 1 & 0 \\ 0 & 1 & 1 \end{pmatrix} \begin{pmatrix} 5 & 2 \\ 1 & 0 \\ 0 & 1 \end{pmatrix} & \begin{pmatrix} 1 & 0 & 1 \\ 1 & 1 & 0 \\ 0 & 1 & 1 \end{pmatrix} \begin{pmatrix} 1 & 1 & 0 \\ 0 & 1 & 1 \\ 1 & 0 & 1 \end{pmatrix} \end{bmatrix}$$

$$= \begin{bmatrix} 11 & 1 & 4 & 3 & 8 \\ 1 & 1 & 0 & 1 & 2 \\ 4 & 0 & 2 & 2 & 4 \\ 3 & 1 & 2 & 5 & 8 \\ 8 & 2 & 4 & 8 & 14 \end{bmatrix} \cup \left[ \begin{array}{cccc|cc|ccc} 30 & 17 & 4 & 6 & 7 & 7 & 6 & 3 & 7 \\ 17 & 10 & 1 & 3 & 1 & 3 & 3 & 1 & 4 \\ 4 & 1 & 5 & 2 & 11 & 4 & 2 & 3 & 1 \\ 6 & 3 & 2 & 2 & 5 & 3 & 2 & 1 & 1 \\ \hline 7 & 1 & 11 & 5 & 26 & 10 & 5 & 6 & 1 \\ 7 & 3 & 4 & 3 & 10 & 5 & 3 & 2 & 1 \\ \hline 6 & 3 & 2 & 2 & 5 & 3 & 2 & 1 & 1 \\ 3 & 1 & 3 & 1 & 6 & 2 & 1 & 2 & 1 \\ 7 & 4 & 1 & 1 & 1 & 1 & 1 & 1 & 2 \end{array} \right]$$

$= S_1 \cup S_2 = A^T A,$

we see $A^TA$ is a symmetric semi superbimatrix. This product helps one to construct any number of symmetric semi superbimatrices.

***Example 2.54:*** Let $A = A_1 \cup A_2$ be column semi superbivector. Then we can find $AA^T$. Given

$$A = A_1 \cup A_2$$

$$= \begin{bmatrix} 3 & 0 & 1 & 2 \\ 1 & 1 & 0 & 1 \\ 5 & 2 & 0 & 1 \\ 0 & 1 & 1 & 0 \end{bmatrix} \cup \left[ \begin{array}{cccc} 0 & 1 & 2 & 3 \\ 1 & 1 & 0 & 1 \\ 5 & 3 & 1 & 2 \\ 0 & 1 & 0 & 1 \\ \hline 3 & 3 & 0 & 1 \\ 1 & 2 & 0 & 0 \\ 0 & 3 & 1 & 0 \\ 1 & 0 & 1 & 0 \\ 1 & 1 & 1 & 1 \\ 0 & 1 & 0 & 1 \\ 0 & 1 & 1 & 0 \\ 1 & 0 & 0 & 1 \end{array} \right]$$

to be a column semi superbivector.
Now

$$\begin{aligned} A^T &= (A_1 \cup A_2)^T \\ &= A_1^T \cup A_2^T \end{aligned}$$

$$= \begin{bmatrix} 3 & 1 & 5 & 0 \\ 0 & 1 & 2 & 1 \\ 1 & 0 & 0 & 1 \\ 2 & 1 & 1 & 0 \end{bmatrix} \cup$$

$$\begin{bmatrix} 0 & 1 & 5 & 0 & | & 3 & 1 & | & 0 & 1 & 1 & 0 & 0 & 1 \\ 1 & 1 & 3 & 1 & | & 3 & 2 & | & 3 & 0 & 1 & 1 & 1 & 0 \\ 2 & 0 & 1 & 0 & | & 0 & 0 & | & 1 & 1 & 1 & 0 & 1 & 0 \\ 3 & 1 & 2 & 1 & | & 1 & 0 & | & 0 & 0 & 1 & 1 & 0 & 1 \end{bmatrix}.$$

Now

$$\begin{aligned} AA^T &= (A_1 \cup A_2)(A_1 \cup A_2)^T \\ &= (A_1 \cup A_2)(A_1^T \cup A_2^T) \\ &= A_1 A_1^T \cup A_2 A_2^T \end{aligned}$$

$$= \begin{bmatrix} 3 & 0 & 1 & 2 \\ 1 & 1 & 0 & 1 \\ 5 & 2 & 0 & 1 \\ 0 & 1 & 1 & 0 \end{bmatrix} \begin{bmatrix} 3 & 1 & 5 & 0 \\ 0 & 1 & 2 & 1 \\ 1 & 0 & 0 & 1 \\ 2 & 1 & 1 & 0 \end{bmatrix} \cup \begin{bmatrix} 0 & 1 & 2 & 3 \\ 1 & 1 & 0 & 1 \\ 5 & 3 & 1 & 2 \\ 0 & 1 & 0 & 1 \\ \hline 3 & 3 & 0 & 1 \\ 1 & 2 & 6 & 0 \\ \hline 0 & 3 & 1 & 0 \\ 1 & 0 & 1 & 0 \\ 1 & 1 & 1 & 1 \\ 0 & 1 & 0 & 1 \\ 0 & 1 & 1 & 0 \\ 1 & 0 & 0 & 1 \end{bmatrix}$$

$$\begin{bmatrix} 0 & 1 & 5 & 0 & | & 3 & 1 & | & 0 & 1 & 1 & 0 & 0 & 1 \\ 1 & 1 & 3 & 1 & | & 3 & 2 & | & 3 & 0 & 1 & 1 & 1 & 0 \\ 2 & 0 & 1 & 0 & | & 0 & 6 & | & 1 & 1 & 1 & 0 & 1 & 0 \\ 3 & 1 & 2 & 1 & | & 1 & 0 & | & 0 & 0 & 1 & 1 & 0 & 1 \end{bmatrix}$$

$$= \begin{bmatrix} 14 & 5 & 17 & 1 \\ 5 & 3 & 8 & 1 \\ 17 & 8 & 30 & 2 \\ 1 & 1 & 2 & 2 \end{bmatrix} \cup$$

$$\left[\begin{pmatrix}0&1&2&3\\1&1&0&1\\5&3&1&2\\0&1&0&1\end{pmatrix}\begin{pmatrix}0&1&5&0\\1&1&3&1\\2&0&1&0\\3&1&2&1\end{pmatrix}\right|\left.\begin{pmatrix}0&1&2&3\\1&1&0&1\\5&3&1&2\\0&1&0&1\end{pmatrix}\begin{bmatrix}3&1\\3&2\\0&6\\1&0\end{bmatrix}\right.$$

$$\begin{pmatrix}3&3&0&1\\1&2&6&0\end{pmatrix}\begin{pmatrix}0&1&5&0\\1&1&3&1\\2&0&1&0\\3&1&2&1\end{pmatrix} \quad \begin{pmatrix}3&3&0&1\\1&2&6&0\end{pmatrix}\begin{bmatrix}3&1\\3&2\\0&6\\1&0\end{bmatrix}$$

$$\begin{pmatrix}0&3&1&0\\1&0&1&0\\1&1&1&1\\0&1&0&1\\0&1&1&0\\1&0&0&1\end{pmatrix}\begin{pmatrix}0&1&5&0\\1&1&3&1\\2&0&1&0\\3&1&2&1\end{pmatrix} \quad \begin{pmatrix}0&3&1&0\\1&0&1&0\\1&1&1&1\\0&1&0&1\\0&1&1&0\\1&0&0&1\end{pmatrix}\begin{bmatrix}3&1\\3&2\\0&6\\1&0\end{bmatrix}\left.\right]$$

$$\left[\begin{pmatrix}0&1&2&3\\1&1&0&1\\5&3&1&2\\0&1&0&1\end{pmatrix}\begin{bmatrix}0&1&1&0&0&1\\3&0&1&1&1&0\\1&1&1&0&1&0\\0&0&1&1&0&1\end{bmatrix}\right.$$

$$\begin{pmatrix}3&3&0&1\\1&2&6&0\end{pmatrix}\begin{bmatrix}0&1&1&0&0&1\\3&0&1&1&1&0\\1&1&1&0&1&0\\0&0&1&1&0&1\end{bmatrix}$$

$$\begin{pmatrix}0&3&1&0\\1&0&1&0\\1&1&1&1\\0&1&0&1\\0&1&1&0\\1&0&0&1\end{pmatrix}\begin{bmatrix}0&1&1&0&0&1\\3&0&1&1&1&0\\1&1&1&0&1&0\\0&0&1&1&0&1\end{bmatrix}\left.\right]$$

$$= \begin{bmatrix} 14 & 5 & 17 & 1 \\ 5 & 3 & 8 & 1 \\ 17 & 8 & 30 & 2 \\ 1 & 1 & 2 & 2 \end{bmatrix} \cup$$

$$\begin{bmatrix} 14 & 4 & 11 & 4 & 6 & 14 & 5 & 2 & 6 & 4 & 3 & 3 \\ 4 & 3 & 10 & 2 & 7 & 3 & 3 & 1 & 3 & 2 & 1 & 2 \\ 11 & 10 & 39 & 5 & 26 & 17 & 10 & 6 & 11 & 5 & 4 & 7 \\ 4 & 2 & 5 & 2 & 4 & 2 & 3 & 0 & 2 & 2 & 1 & 1 \\ \hline 6 & 7 & 26 & 4 & 19 & 9 & 9 & 3 & 7 & 4 & 3 & 4 \\ 14 & 3 & 17 & 2 & 9 & 41 & 12 & 7 & 9 & 2 & 8 & 1 \\ \hline 5 & 3 & 10 & 3 & 9 & 12 & 10 & 1 & 4 & 3 & 4 & 0 \\ 2 & 1 & 6 & 0 & 3 & 7 & 1 & 2 & 2 & 0 & 1 & 1 \\ 6 & 3 & 11 & 2 & 7 & 9 & 4 & 2 & 4 & 2 & 2 & 2 \\ 4 & 2 & 5 & 2 & 4 & 2 & 3 & 0 & 2 & 2 & 1 & 1 \\ 3 & 1 & 4 & 1 & 3 & 8 & 4 & 1 & 2 & 1 & 2 & 0 \\ 3 & 2 & 7 & 1 & 4 & 1 & 0 & 1 & 2 & 1 & 0 & 2 \end{bmatrix}.$$

$AA^T$ is a symmetric semi superbimatrix. Thus this major product when A is multiplied by its transpose where A is only semi bivector yields a symmetric semi superbimatrix which is not a bivector.

Like in case of superbimatrices we can easily prove the following theorem.

**THEOREM 2.2:** *(1) Let $A = A_1 \cup A_2$ be a column semi superbivector then $AA^T$ is a symmetric semi superbimatrix.*
*(2) If $A = A_1 \cup A_2$ is a row semi superbivector then $A^TA$ is a symmetric semi superbimatrix.*

The examples 2.53 and 2.54 substantiate the above theorem. Now we proceed onto define the product of superbimatrices which are not superbivectors and semi superbimatrices which are not semi superbivectors.

*Example 2.55:* Let $A = A_1 \cup A_2$ and $B = B_1 \cup B_2$ be any two superbimatrices we define the byproduct AB of A and B. Given

$$A = A_1 \cup A_2$$

$$= \begin{bmatrix} 0 & 3 & 1 & 1 & 1 & 2 \\ 1 & 1 & 0 & 2 & 1 & 1 \\ 2 & 1 & 1 & 0 & 1 & 1 \\ 0 & 0 & 1 & 1 & 0 & 1 \\ \hline 1 & 2 & 4 & 0 & 1 & 2 \\ 0 & 0 & 1 & 1 & 0 & 1 \\ \hline 1 & 1 & 1 & 3 & 0 & 1 \\ 0 & 0 & 1 & 1 & 0 & 3 \\ 1 & 1 & 0 & 0 & 1 & 0 \end{bmatrix} \cup \begin{bmatrix} 3 & 0 & 1 & 1 & 0 & 1 & 1 & 5 & 0 \\ 1 & 1 & 0 & 1 & 1 & 0 & 1 & 0 & 1 \\ 0 & 1 & 1 & 0 & 0 & 1 & 2 & 1 & 0 \\ \hline 2 & 1 & 1 & 1 & 0 & 1 & 1 & 0 & 0 \\ 1 & 2 & 1 & 0 & 1 & 1 & 0 & 1 & 0 \\ 0 & 1 & 0 & 1 & 1 & 1 & 0 & 0 & 1 \\ \hline 1 & 0 & 0 & 0 & 1 & 0 & 1 & 0 & 2 \\ 5 & 2 & 2 & 1 & 0 & 2 & 0 & 1 & 2 \end{bmatrix}$$

and

$$B = B_1 \cup B_2$$

$$= \begin{bmatrix} 1 & 0 & 1 & 1 & 2 & 4 & 3 \\ 0 & 1 & 1 & 0 & 1 & 2 & 0 \\ \hline 1 & 0 & 0 & 1 & 0 & 1 & 2 \\ \hline 3 & 1 & 0 & 1 & 0 & 1 & 0 \\ 1 & 1 & 1 & 0 & 1 & 1 & 1 \\ 0 & 1 & 0 & 1 & 1 & 0 & 5 \end{bmatrix} \cup \begin{bmatrix} 1 & 0 & 1 & 1 & 1 & 2 \\ 1 & 0 & 1 & 0 & 1 & 0 \\ \hline 1 & 1 & 0 & 1 & 0 & 1 \\ 0 & 1 & 2 & 5 & 1 & 2 \\ 1 & 0 & 1 & 1 & 1 & 0 \\ 1 & 0 & 0 & 1 & 0 & 1 \\ \hline 3 & 1 & 2 & 1 & 1 & 2 \\ 0 & 1 & 1 & 0 & 6 & 0 \\ 1 & 0 & 1 & 2 & 1 & 1 \end{bmatrix}$$

AB  =  $(A_1 \cup A_2)(B_1 \cup B_2)$
    =  $A_1 B_1 \cup A_2 B_2$

$$= \begin{bmatrix} 0 & 3 & 1 & 1 & 1 & 2 \\ 1 & 1 & 0 & 2 & 1 & 1 \\ 2 & 1 & 1 & 0 & 1 & 1 \\ 0 & 0 & 1 & 1 & 0 & 1 \\ \hline 1 & 2 & 4 & 0 & 1 & 2 \\ 0 & 0 & 1 & 1 & 0 & 1 \\ \hline 1 & 1 & 1 & 3 & 0 & 1 \\ 0 & 0 & 1 & 1 & 0 & 3 \\ 1 & 1 & 0 & 0 & 1 & 0 \end{bmatrix} \begin{bmatrix} 1 & 0 & 1 & 1 & 2 & 4 & 3 \\ \hline 0 & 1 & 1 & 0 & 1 & 2 & 0 \\ 1 & 0 & 0 & 1 & 0 & 1 & 2 \\ \hline 3 & 1 & 0 & 1 & 0 & 1 & 0 \\ 1 & 1 & 1 & 0 & 1 & 1 & 1 \\ 0 & 1 & 0 & 1 & 1 & 0 & 5 \end{bmatrix} \subset$$

$$\begin{bmatrix} 3 & 0 & 1 & 1 & 0 & 1 & 1 & 5 & 0 \\ 1 & 1 & 0 & 1 & 1 & 0 & 1 & 0 & 1 \\ 0 & 1 & 1 & 0 & 0 & 1 & 2 & 1 & 0 \\ \hline 2 & 1 & 1 & 1 & 0 & 1 & 1 & 0 & 0 \\ 1 & 2 & 1 & 0 & 1 & 1 & 0 & 1 & 0 \\ 0 & 1 & 0 & 1 & 1 & 1 & 0 & 0 & 1 \\ 1 & 0 & 0 & 0 & 1 & 0 & 1 & 0 & 2 \\ \hline 5 & 2 & 2 & 1 & 0 & 2 & 0 & 1 & 2 \end{bmatrix} \begin{bmatrix} 1 & 0 & 1 & 1 & 1 & 2 \\ 1 & 0 & 1 & 0 & 1 & 0 \\ 1 & 1 & 0 & 1 & 0 & 1 \\ 0 & 1 & 2 & 5 & 1 & 2 \\ 1 & 0 & 1 & 1 & 1 & 0 \\ 1 & 0 & 0 & 1 & 0 & 1 \\ \hline 3 & 1 & 2 & 1 & 1 & 2 \\ 0 & 1 & 1 & 0 & 6 & 0 \\ 1 & 0 & 1 & 2 & 1 & 1 \end{bmatrix}$$

$$= \begin{bmatrix} \begin{bmatrix} 0 \\ 1 \\ 2 \\ 0 \end{bmatrix} [1 \ 0] & \begin{bmatrix} 0 \\ 1 \\ 2 \\ 0 \end{bmatrix} [1 \ 1 \ 2] & \begin{bmatrix} 0 \\ 1 \\ 2 \\ 0 \end{bmatrix} [4 \ 3] \\ \hline \begin{bmatrix} 0 \\ 1 \end{bmatrix} [1 \ 0] & \begin{bmatrix} 0 \\ 1 \end{bmatrix} [1 \ 1 \ 2] & \begin{bmatrix} 0 \\ 1 \end{bmatrix} [4 \ 3] \\ \hline \begin{bmatrix} 1 \\ 0 \\ 1 \end{bmatrix} [1 \ 0] & \begin{bmatrix} 1 \\ 0 \\ 1 \end{bmatrix} [1 \ 1 \ 2] & \begin{bmatrix} 1 \\ 0 \\ 1 \end{bmatrix} [4 \ 3] \end{bmatrix} +$$

$$\begin{bmatrix} \begin{bmatrix} 3 & 1 \\ 1 & 0 \\ 1 & 1 \\ 0 & 1 \end{bmatrix}\begin{bmatrix} 0 & 1 \\ 1 & 0 \end{bmatrix} & \begin{bmatrix} 3 & 1 \\ 1 & 0 \\ 1 & 1 \\ 0 & 1 \end{bmatrix}\begin{bmatrix} 1 & 0 & 1 \\ 0 & 1 & 0 \end{bmatrix} & \begin{bmatrix} 3 & 1 \\ 1 & 0 \\ 1 & 1 \\ 0 & 1 \end{bmatrix}\begin{bmatrix} 2 & 0 \\ 1 & 2 \end{bmatrix} \\ \begin{bmatrix} 2 & 4 \\ 0 & 1 \end{bmatrix}\begin{bmatrix} 0 & 1 \\ 1 & 0 \end{bmatrix} & \begin{bmatrix} 2 & 4 \\ 0 & 1 \end{bmatrix}\begin{bmatrix} 1 & 0 & 1 \\ 0 & 1 & 0 \end{bmatrix} & \begin{bmatrix} 2 & 4 \\ 0 & 1 \end{bmatrix}\begin{bmatrix} 2 & 0 \\ 1 & 2 \end{bmatrix} \\ \begin{bmatrix} 1 & 1 \\ 0 & 1 \\ 1 & 0 \end{bmatrix}\begin{bmatrix} 0 & 1 \\ 1 & 0 \end{bmatrix} & \begin{bmatrix} 1 & 1 \\ 0 & 1 \\ 1 & 0 \end{bmatrix}\begin{bmatrix} 1 & 0 & 1 \\ 0 & 1 & 0 \end{bmatrix} & \begin{bmatrix} 1 & 1 \\ 0 & 1 \\ 1 & 0 \end{bmatrix}\begin{bmatrix} 2 & 0 \\ 1 & 2 \end{bmatrix} \end{bmatrix} =$$

$$\begin{bmatrix} \begin{bmatrix} 1 & 1 & 2 \\ 2 & 1 & 1 \\ 0 & 1 & 1 \\ 1 & 0 & 0 \end{bmatrix}\begin{bmatrix} 3 & 1 \\ 1 & 1 \\ 0 & 1 \end{bmatrix} & \begin{bmatrix} 1 & 1 & 2 \\ 2 & 1 & 1 \\ 0 & 1 & 1 \\ 1 & 0 & 0 \end{bmatrix}\begin{bmatrix} 0 & 1 & 0 \\ 1 & 0 & 1 \\ 0 & 1 & 1 \end{bmatrix} & \begin{bmatrix} 1 & 1 & 2 \\ 2 & 1 & 1 \\ 0 & 1 & 1 \\ 1 & 0 & 0 \end{bmatrix}\begin{bmatrix} 1 & 0 \\ 1 & 1 \\ 0 & 5 \end{bmatrix} \\ \begin{bmatrix} 0 & 1 & 2 \\ 1 & 0 & 1 \end{bmatrix}\begin{bmatrix} 3 & 1 \\ 1 & 1 \\ 0 & 1 \end{bmatrix} & \begin{bmatrix} 0 & 1 & 2 \\ 1 & 0 & 1 \end{bmatrix}\begin{bmatrix} 0 & 1 & 0 \\ 1 & 0 & 1 \\ 0 & 1 & 1 \end{bmatrix} & \begin{bmatrix} 0 & 1 & 2 \\ 1 & 0 & 1 \end{bmatrix}\begin{bmatrix} 1 & 0 \\ 1 & 1 \\ 0 & 5 \end{bmatrix} \\ \begin{bmatrix} 3 & 0 & 1 \\ 1 & 0 & 3 \\ 0 & 1 & 0 \end{bmatrix}\begin{bmatrix} 3 & 1 \\ 1 & 1 \\ 0 & 1 \end{bmatrix} & \begin{bmatrix} 3 & 0 & 1 \\ 1 & 0 & 3 \\ 0 & 1 & 0 \end{bmatrix}\begin{bmatrix} 0 & 1 & 0 \\ 1 & 0 & 1 \\ 0 & 1 & 1 \end{bmatrix} & \begin{bmatrix} 3 & 0 & 1 \\ 1 & 0 & 3 \\ 0 & 1 & 0 \end{bmatrix}\begin{bmatrix} 1 & 0 \\ 1 & 1 \\ 0 & 5 \end{bmatrix} \end{bmatrix}$$

$$\subset \begin{bmatrix} \begin{bmatrix} 3 & 0 \\ 1 & 1 \\ 0 & 1 \end{bmatrix}\begin{bmatrix} 1 & 0 & 1 \\ 1 & 0 & 1 \end{bmatrix} & \begin{bmatrix} 3 & 0 \\ 1 & 1 \\ 0 & 1 \end{bmatrix}\begin{bmatrix} 1 & 1 \\ 0 & 1 \end{bmatrix} & \begin{bmatrix} 3 & 0 \\ 1 & 1 \\ 0 & 1 \end{bmatrix}\begin{bmatrix} 2 \\ 0 \end{bmatrix} \\ \begin{bmatrix} 2 & 1 \\ 1 & 2 \\ 0 & 1 \\ 1 & 0 \end{bmatrix}\begin{bmatrix} 1 & 0 & 1 \\ 1 & 0 & 1 \end{bmatrix} & \begin{bmatrix} 2 & 1 \\ 1 & 2 \\ 0 & 1 \\ 1 & 0 \end{bmatrix}\begin{bmatrix} 1 & 1 \\ 0 & 1 \end{bmatrix} & \begin{bmatrix} 2 & 1 \\ 1 & 2 \\ 0 & 1 \\ 1 & 0 \end{bmatrix}\begin{bmatrix} 2 \\ 0 \end{bmatrix} \\ \begin{bmatrix} 5 & 2 \end{bmatrix}\begin{bmatrix} 1 & 0 & 1 \\ 1 & 0 & 1 \end{bmatrix} & \begin{bmatrix} 5 & 2 \end{bmatrix}\begin{bmatrix} 1 & 1 \\ 0 & 1 \end{bmatrix} & \begin{bmatrix} 5 & 2 \end{bmatrix}\begin{bmatrix} 2 \\ 0 \end{bmatrix} \end{bmatrix} +$$

$$\begin{bmatrix} \begin{bmatrix} 1 & 1 & 0 & 1 \\ 0 & 1 & 1 & 0 \\ 1 & 0 & 0 & 1 \end{bmatrix} \begin{bmatrix} 1 & 1 & 0 \\ 0 & 1 & 2 \\ 1 & 0 & 1 \\ 1 & 0 & 0 \end{bmatrix} & \begin{bmatrix} 1 & 1 & 0 & 1 \\ 0 & 1 & 1 & 0 \\ 1 & 0 & 0 & 1 \end{bmatrix} \begin{bmatrix} 1 & 0 \\ 5 & 1 \\ 1 & 1 \\ 1 & 0 \end{bmatrix} & \begin{bmatrix} 1 & 1 & 0 & 1 \\ 0 & 1 & 1 & 0 \\ 1 & 0 & 0 & 1 \end{bmatrix} \begin{bmatrix} 1 \\ 2 \\ 0 \\ 1 \end{bmatrix} \\ \begin{bmatrix} 1 & 1 & 0 & 1 \\ 1 & 0 & 1 & 1 \\ 0 & 1 & 1 & 1 \\ 0 & 0 & 1 & 0 \end{bmatrix} \begin{bmatrix} 1 & 1 & 0 \\ 0 & 1 & 2 \\ 1 & 0 & 1 \\ 1 & 0 & 0 \end{bmatrix} & \begin{bmatrix} 1 & 1 & 0 & 1 \\ 1 & 0 & 1 & 1 \\ 0 & 1 & 1 & 1 \\ 0 & 0 & 1 & 0 \end{bmatrix} \begin{bmatrix} 1 & 0 \\ 5 & 1 \\ 1 & 1 \\ 1 & 0 \end{bmatrix} & \begin{bmatrix} 1 & 1 & 0 & 1 \\ 1 & 0 & 1 & 1 \\ 0 & 1 & 1 & 1 \\ 0 & 0 & 1 & 0 \end{bmatrix} \begin{bmatrix} 1 \\ 2 \\ 0 \\ 1 \end{bmatrix} \\ \begin{bmatrix} 2 & 1 & 0 & 2 \end{bmatrix} \begin{bmatrix} 1 & 1 & 0 \\ 0 & 1 & 2 \\ 1 & 0 & 1 \\ 1 & 0 & 0 \end{bmatrix} & \begin{bmatrix} 2 & 1 & 0 & 2 \end{bmatrix} \begin{bmatrix} 1 & 0 \\ 5 & 1 \\ 1 & 1 \\ 1 & 0 \end{bmatrix} & \begin{bmatrix} 2 & 1 & 0 & 2 \end{bmatrix} \begin{bmatrix} 1 \\ 2 \\ 0 \\ 1 \end{bmatrix} \end{bmatrix}$$

$$+ \begin{bmatrix} \begin{bmatrix} 1 & 5 & 0 \\ 1 & 0 & 1 \\ 2 & 1 & 0 \end{bmatrix} \begin{bmatrix} 3 & 1 & 2 \\ 0 & 1 & 1 \\ 1 & 0 & 1 \end{bmatrix} & \begin{bmatrix} 1 & 5 & 0 \\ 1 & 0 & 1 \\ 2 & 1 & 0 \end{bmatrix} \begin{bmatrix} 1 & 1 \\ 0 & 6 \\ 2 & 1 \end{bmatrix} & \begin{bmatrix} 1 & 5 & 0 \\ 1 & 0 & 1 \\ 2 & 1 & 0 \end{bmatrix} \begin{bmatrix} 2 \\ 0 \\ 1 \end{bmatrix} \\ \begin{bmatrix} 1 & 0 & 0 \\ 0 & 1 & 0 \\ 0 & 0 & 1 \\ 1 & 0 & 2 \end{bmatrix} \begin{bmatrix} 3 & 1 & 2 \\ 0 & 1 & 1 \\ 1 & 0 & 1 \end{bmatrix} & \begin{bmatrix} 1 & 0 & 0 \\ 0 & 1 & 0 \\ 0 & 0 & 1 \\ 1 & 0 & 2 \end{bmatrix} \begin{bmatrix} 1 & 1 \\ 0 & 6 \\ 2 & 1 \end{bmatrix} & \begin{bmatrix} 1 & 0 & 0 \\ 0 & 1 & 0 \\ 0 & 0 & 1 \\ 1 & 0 & 2 \end{bmatrix} \begin{bmatrix} 2 \\ 0 \\ 1 \end{bmatrix} \\ \begin{bmatrix} 0 & 1 & 2 \end{bmatrix} \begin{bmatrix} 3 & 1 & 2 \\ 0 & 1 & 1 \\ 1 & 0 & 1 \end{bmatrix} & \begin{bmatrix} 0 & 1 & 2 \end{bmatrix} \begin{bmatrix} 1 & 1 \\ 0 & 6 \\ 2 & 1 \end{bmatrix} & \begin{bmatrix} 0 & 1 & 2 \end{bmatrix} \begin{bmatrix} 2 \\ 0 \\ 1 \end{bmatrix} \end{bmatrix}$$

$$= \left[\begin{array}{cc|ccc|cc} 0 & 0 & 0 & 0 & 0 & 0 & 0 \\ 1 & 0 & 1 & 1 & 2 & 4 & 3 \\ 2 & 0 & 2 & 2 & 4 & 8 & 6 \\ 0 & 0 & 0 & 0 & 0 & 0 & 0 \\ \hline 0 & 0 & 0 & 0 & 0 & 0 & 0 \\ 1 & 0 & 1 & 1 & 2 & 4 & 3 \\ \hline 1 & 0 & 1 & 1 & 2 & 4 & 3 \\ 0 & 0 & 0 & 0 & 0 & 0 & 0 \\ 1 & 0 & 1 & 1 & 2 & 4 & 3 \end{array}\right] +$$

$$\left[\begin{array}{cc|ccc|cc} 1 & 3 & 3 & 1 & 3 & 7 & 2 \\ 0 & 1 & 1 & 0 & 1 & 2 & 0 \\ 1 & 1 & 1 & 1 & 1 & 3 & 2 \\ 1 & 0 & 0 & 1 & 0 & 1 & 2 \\ \hline 4 & 2 & 2 & 4 & 2 & 8 & 8 \\ 1 & 0 & 0 & 1 & 0 & 1 & 2 \\ \hline 1 & 1 & 1 & 1 & 1 & 3 & 2 \\ 1 & 0 & 0 & 1 & 0 & 1 & 2 \\ 0 & 1 & 1 & 0 & 1 & 2 & 0 \end{array}\right] +$$

$$\left[\begin{array}{cc|ccc|cc} 4 & 4 & 1 & 3 & 3 & 2 & 11 \\ 7 & 4 & 1 & 3 & 2 & 3 & 6 \\ 1 & 2 & 1 & 1 & 2 & 1 & 6 \\ 3 & 1 & 0 & 2 & 1 & 1 & 0 \\ \hline 1 & 3 & 1 & 2 & 3 & 1 & 11 \\ 3 & 2 & 0 & 2 & 1 & 1 & 5 \\ \hline 9 & 4 & 0 & 4 & 1 & 3 & 5 \\ 3 & 4 & 0 & 4 & 3 & 1 & 15 \\ 1 & 1 & 1 & 0 & 1 & 1 & 1 \end{array}\right] \cup$$

$$\begin{bmatrix} 3 & 0 & 3 & 3 & 3 & 6 \\ 2 & 0 & 2 & 1 & 2 & 2 \\ 1 & 0 & 1 & 0 & 1 & 0 \\ \hline 3 & 0 & 3 & 2 & 3 & 4 \\ 3 & 0 & 3 & 1 & 3 & 1 \\ 1 & 0 & 1 & 0 & 1 & 0 \\ 1 & 0 & 1 & 1 & 1 & 2 \\ \hline 7 & 0 & 7 & 5 & 7 & 10 \end{bmatrix} +$$

$$\begin{bmatrix} 2 & 2 & 2 & 7 & 1 & 4 \\ 1 & 1 & 3 & 6 & 2 & 2 \\ 2 & 1 & 0 & 2 & 0 & 2 \\ \hline 2 & 2 & 2 & 7 & 1 & 4 \\ 3 & 1 & 1 & 3 & 1 & 2 \\ 2 & 1 & 3 & 7 & 2 & 3 \\ 1 & 0 & 1 & 1 & 1 & 0 \\ \hline 4 & 3 & 2 & 9 & 1 & 6 \end{bmatrix} +$$

$$\begin{bmatrix} 3 & 6 & 7 & 1 & 31 & 2 \\ 4 & 1 & 3 & 3 & 2 & 3 \\ 6 & 3 & 5 & 2 & 8 & 4 \\ \hline 3 & 1 & 2 & 1 & 1 & 2 \\ 0 & 1 & 1 & 0 & 6 & 0 \\ 1 & 0 & 1 & 2 & 1 & 1 \\ 5 & 1 & 4 & 5 & 3 & 4 \\ \hline 2 & 1 & 3 & 4 & 8 & 2 \end{bmatrix} =$$

$$\begin{bmatrix} 5 & 7 & 4 & 4 & 6 & 9 & 13 \\ 8 & 5 & 3 & 4 & 5 & 9 & 9 \\ 4 & 3 & 4 & 4 & 7 & 12 & 14 \\ 4 & 1 & 0 & 3 & 1 & 2 & 2 \\ 5 & 5 & 3 & 6 & 5 & 9 & 19 \\ 5 & 2 & 1 & 4 & 3 & 6 & 10 \\ 11 & 5 & 2 & 6 & 4 & 10 & 10 \\ 4 & 4 & 0 & 5 & 3 & 2 & 17 \\ 2 & 2 & 3 & 1 & 4 & 7 & 4 \end{bmatrix} \cup$$

$$\begin{bmatrix} 8 & 8 & 12 & 11 & 35 & 12 \\ 7 & 2 & 8 & 10 & 6 & 7 \\ 9 & 4 & 6 & 4 & 9 & 6 \\ 8 & 3 & 7 & 10 & 5 & 10 \\ 6 & 2 & 5 & 4 & 10 & 3 \\ 4 & 1 & 5 & 9 & 4 & 4 \\ 7 & 1 & 6 & 7 & 5 & 6 \\ 13 & 4 & 12 & 18 & 16 & 18 \end{bmatrix}.$$

Clearly the resultant under the minor product of two superbimatrices is a superbimatrix.

Thus we see the minor product of two superbimatrices results in a superbimatrix. We have observed from the example 2.55 that only for the product AB to be compatible we need the number of columns in A must equal number of rows of B but also the way the columns of A are partitioned must be identical with the way the rows of B are partitioned.

Then alone we have the product to be defined. We give yet another example of a minor product of two superbimatrices before we proceed on to define them.

***Example 2.56:*** Let $A = A_1 \cup A_2$ and $B = B_1 \cup B_2$ be any two superbimatrices. We will find the minor product AB.

Given A = A₁ ∪ A₂ where

$$A_1 = \left[\begin{array}{ccc|cc|cc} 0 & 1 & 2 & 3 & 1 & 0 & 1 \\ 2 & 1 & 1 & 0 & 3 & 5 & 4 \\ 1 & 0 & 1 & 1 & 1 & 0 & 2 \\ 0 & 1 & 0 & 0 & 1 & 4 & 1 \\ \hline 1 & 1 & 0 & 2 & 1 & 1 & 1 \\ 1 & 0 & 1 & 5 & 1 & 0 & 1 \\ \hline 1 & 2 & 1 & 0 & 1 & 1 & 1 \\ 2 & 0 & 2 & 1 & 0 & 0 & 1 \\ 3 & 1 & 0 & 0 & 0 & 1 & 0 \end{array}\right]$$

and

$$A_2 = \left[\begin{array}{cc|cccc} 3 & 1 & 1 & 1 & 1 & 0 \\ 1 & 0 & 0 & 1 & 0 & 1 \\ 2 & 1 & 1 & 1 & 1 & 0 \\ \hline 0 & 3 & 2 & 1 & 0 & 1 \\ 1 & 1 & 0 & 1 & 2 & 1 \\ 5 & 2 & 0 & 1 & 0 & 2 \end{array}\right].$$

B = B₁ ∪ B₂ with

$$B_1 = \left[\begin{array}{ccc|ccc} 1 & 0 & 1 & 2 & 1 & 5 \\ 0 & 2 & 1 & 6 & 1 & 2 \\ 1 & 1 & 1 & 1 & 0 & 1 \\ \hline 2 & 5 & 3 & 1 & 1 & 1 \\ 1 & 3 & 1 & 1 & 0 & 1 \\ \hline 4 & 2 & 1 & 3 & 1 & 2 \\ 0 & 0 & 6 & 2 & 1 & 3 \end{array}\right]$$

and

$$B_2 = \begin{bmatrix} 1 & 1 & 1 & 1 & 1 \\ 2 & 2 & 0 & 1 & 1 \\ \hline 3 & 4 & 0 & 1 & 0 \\ 4 & 1 & 1 & 0 & 1 \\ 0 & 6 & 2 & 0 & 5 \\ 1 & 0 & 1 & 0 & 3 \end{bmatrix}.$$

$$\begin{aligned} AB &= (A_1 \cup A_2)(B_1 \cup B_2) \\ &= A_1 B_1 \cup A_2 B_2 \end{aligned}$$

$$= \begin{bmatrix} 0 & 1 & 2 & 3 & 1 & 0 & 1 \\ 2 & 1 & 1 & 0 & 3 & 5 & 4 \\ 1 & 0 & 1 & 1 & 1 & 0 & 2 \\ 0 & 1 & 0 & 0 & 1 & 4 & 1 \\ \hline 1 & 1 & 0 & 2 & 1 & 1 & 1 \\ 1 & 0 & 1 & 5 & 1 & 0 & 1 \\ \hline 1 & 2 & 1 & 0 & 1 & 1 & 1 \\ 2 & 0 & 2 & 1 & 0 & 0 & 1 \\ 3 & 1 & 0 & 0 & 0 & 1 & 0 \end{bmatrix} \begin{bmatrix} 1 & 0 & 1 & 2 & 1 & 5 \\ 0 & 2 & 1 & 6 & 1 & 2 \\ 1 & 1 & 1 & 1 & 0 & 1 \\ \hline 2 & 5 & 3 & 1 & 1 & 1 \\ 1 & 3 & 1 & 1 & 0 & 1 \\ \hline 4 & 2 & 1 & 3 & 1 & 2 \\ 0 & 0 & 6 & 2 & 1 & 3 \end{bmatrix} \cup$$

$$\begin{bmatrix} 3 & 1 & 1 & 1 & 1 & 0 \\ 1 & 0 & 0 & 1 & 0 & 1 \\ 2 & 1 & 1 & 1 & 1 & 0 \\ \hline 0 & 3 & 2 & 1 & 0 & 1 \\ 1 & 1 & 0 & 1 & 2 & 1 \\ 5 & 2 & 0 & 1 & 0 & 2 \end{bmatrix} \begin{bmatrix} 1 & 1 & 1 & 1 & 1 \\ 2 & 2 & 0 & 1 & 1 \\ \hline 3 & 4 & 0 & 1 & 0 \\ 4 & 1 & 1 & 0 & 1 \\ 0 & 6 & 2 & 0 & 5 \\ 1 & 0 & 1 & 0 & 3 \end{bmatrix}$$

$$= \left[ \begin{array}{c|c} \begin{bmatrix} 0 & 1 & 2 \\ 2 & 1 & 1 \\ 1 & 0 & 1 \\ 0 & 1 & 0 \end{bmatrix} \begin{bmatrix} 1 & 0 & 1 \\ 0 & 2 & 1 \\ 1 & 1 & 1 \end{bmatrix} & \begin{bmatrix} 0 & 1 & 2 \\ 2 & 1 & 1 \\ 1 & 0 & 1 \\ 0 & 1 & 0 \end{bmatrix} \begin{bmatrix} 2 & 1 & 5 \\ 6 & 1 & 2 \\ 1 & 0 & 1 \end{bmatrix} \\ \hline \begin{bmatrix} 1 & 1 & 0 \\ 1 & 0 & 1 \end{bmatrix} \begin{bmatrix} 1 & 0 & 1 \\ 0 & 2 & 1 \\ 1 & 1 & 1 \end{bmatrix} & \begin{bmatrix} 1 & 1 & 0 \\ 1 & 0 & 1 \end{bmatrix} \begin{bmatrix} 2 & 1 & 5 \\ 6 & 1 & 2 \\ 1 & 0 & 1 \end{bmatrix} \\ \hline \begin{bmatrix} 1 & 2 & 1 \\ 2 & 0 & 2 \\ 3 & 1 & 0 \end{bmatrix} \begin{bmatrix} 1 & 0 & 1 \\ 0 & 2 & 1 \\ 1 & 1 & 1 \end{bmatrix} & \begin{bmatrix} 1 & 2 & 1 \\ 2 & 0 & 2 \\ 3 & 1 & 0 \end{bmatrix} \begin{bmatrix} 2 & 1 & 5 \\ 6 & 1 & 2 \\ 1 & 0 & 1 \end{bmatrix} \end{array} \right] +$$

$$\left[ \begin{array}{c|c} \begin{bmatrix} 3 & 1 \\ 0 & 3 \\ 1 & 1 \\ 0 & 1 \end{bmatrix} \begin{bmatrix} 2 & 5 & 3 \\ 1 & 3 & 1 \end{bmatrix} & \begin{bmatrix} 3 & 1 \\ 0 & 3 \\ 1 & 1 \\ 0 & 1 \end{bmatrix} \begin{bmatrix} 1 & 1 & 1 \\ 1 & 0 & 1 \end{bmatrix} \\ \hline \begin{bmatrix} 2 & 1 \\ 5 & 1 \end{bmatrix} \begin{bmatrix} 2 & 5 & 3 \\ 1 & 3 & 1 \end{bmatrix} & \begin{bmatrix} 2 & 1 \\ 5 & 1 \end{bmatrix} \begin{bmatrix} 1 & 1 & 1 \\ 1 & 0 & 1 \end{bmatrix} \\ \hline \begin{bmatrix} 0 & 1 \\ 1 & 0 \\ 0 & 0 \end{bmatrix} \begin{bmatrix} 2 & 5 & 3 \\ 1 & 3 & 1 \end{bmatrix} & \begin{bmatrix} 0 & 1 \\ 1 & 0 \\ 0 & 0 \end{bmatrix} \begin{bmatrix} 1 & 1 & 1 \\ 1 & 0 & 1 \end{bmatrix} \end{array} \right] +$$

$$\left[ \begin{array}{c|c} \begin{bmatrix} 0 & 1 \\ 5 & 4 \\ 0 & 2 \\ 4 & 1 \end{bmatrix} \begin{bmatrix} 4 & 2 & 1 \\ 0 & 0 & 6 \end{bmatrix} & \begin{bmatrix} 0 & 1 \\ 5 & 4 \\ 0 & 2 \\ 4 & 1 \end{bmatrix} \begin{bmatrix} 3 & 1 & 2 \\ 2 & 1 & 3 \end{bmatrix} \\ \hline \begin{bmatrix} 1 & 1 \\ 0 & 1 \end{bmatrix} \begin{bmatrix} 4 & 2 & 1 \\ 0 & 0 & 6 \end{bmatrix} & \begin{bmatrix} 1 & 1 \\ 0 & 1 \end{bmatrix} \begin{bmatrix} 3 & 1 & 2 \\ 2 & 1 & 3 \end{bmatrix} \\ \hline \begin{bmatrix} 1 & 1 \\ 0 & 1 \\ 1 & 0 \end{bmatrix} \begin{bmatrix} 4 & 2 & 1 \\ 0 & 0 & 6 \end{bmatrix} & \begin{bmatrix} 1 & 1 \\ 0 & 1 \\ 1 & 0 \end{bmatrix} \begin{bmatrix} 3 & 1 & 2 \\ 2 & 1 & 3 \end{bmatrix} \end{array} \right] \cup$$

$$\begin{bmatrix} \begin{bmatrix} 3 & 1 \\ 1 & 0 \\ 2 & 1 \end{bmatrix}\begin{bmatrix} 1 & 1 \\ 2 & 2 \end{bmatrix} & \begin{bmatrix} 3 & 1 \\ 1 & 0 \\ 2 & 1 \end{bmatrix}\begin{bmatrix} 1 & 1 & 1 \\ 0 & 1 & 1 \end{bmatrix} \\ \hline \begin{bmatrix} 0 & 3 \\ 1 & 1 \\ 5 & 2 \end{bmatrix}\begin{bmatrix} 1 & 1 \\ 2 & 2 \end{bmatrix} & \begin{bmatrix} 0 & 3 \\ 1 & 1 \\ 5 & 2 \end{bmatrix}\begin{bmatrix} 1 & 1 & 1 \\ 0 & 1 & 1 \end{bmatrix} \end{bmatrix} +$$

$$\begin{bmatrix} \begin{bmatrix} 1 & 1 & 1 & 0 \\ 0 & 1 & 0 & 1 \\ 1 & 1 & 1 & 0 \end{bmatrix}\begin{bmatrix} 3 & 4 \\ 4 & 1 \\ 0 & 6 \\ 1 & 0 \end{bmatrix} & \begin{bmatrix} 1 & 1 & 1 & 0 \\ 0 & 1 & 0 & 1 \\ 1 & 1 & 1 & 0 \end{bmatrix}\begin{bmatrix} 0 & 1 & 0 \\ 1 & 0 & 1 \\ 2 & 0 & 5 \\ 1 & 0 & 3 \end{bmatrix} \\ \hline \begin{bmatrix} 2 & 1 & 0 & 1 \\ 0 & 1 & 2 & 1 \\ 0 & 1 & 0 & 2 \end{bmatrix}\begin{bmatrix} 3 & 4 \\ 4 & 1 \\ 0 & 6 \\ 1 & 0 \end{bmatrix} & \begin{bmatrix} 2 & 1 & 0 & 1 \\ 0 & 1 & 2 & 1 \\ 0 & 1 & 0 & 2 \end{bmatrix}\begin{bmatrix} 0 & 1 & 0 \\ 1 & 0 & 1 \\ 2 & 0 & 5 \\ 1 & 0 & 3 \end{bmatrix} \end{bmatrix}$$

$$= \begin{bmatrix} 2 & 4 & 3 & 8 & 1 & 4 \\ 3 & 3 & 4 & 11 & 3 & 13 \\ 2 & 1 & 2 & 3 & 1 & 6 \\ 0 & 2 & 1 & 6 & 1 & 2 \\ \hline 1 & 2 & 2 & 8 & 2 & 7 \\ 2 & 1 & 2 & 3 & 1 & 6 \\ \hline 2 & 5 & 4 & 15 & 3 & 10 \\ 4 & 2 & 4 & 6 & 2 & 12 \\ 3 & 2 & 4 & 12 & 4 & 17 \end{bmatrix} +$$

$$\begin{bmatrix} 7 & 18 & 10 & 4 & 3 & 4 \\ 3 & 9 & 3 & 3 & 0 & 3 \\ 3 & 8 & 4 & 2 & 1 & 2 \\ 1 & 3 & 1 & 1 & 0 & 1 \\ \hline 5 & 13 & 7 & 3 & 2 & 3 \\ 11 & 28 & 16 & 6 & 5 & 6 \\ \hline 1 & 3 & 1 & 1 & 0 & 1 \\ 2 & 5 & 3 & 1 & 1 & 1 \\ 0 & 0 & 0 & 0 & 0 & 0 \end{bmatrix} +$$

$$\begin{bmatrix} 0 & 0 & 6 & 2 & 1 & 3 \\ 20 & 10 & 29 & 23 & 9 & 22 \\ 0 & 0 & 12 & 4 & 2 & 6 \\ 16 & 8 & 10 & 14 & 5 & 11 \\ \hline 4 & 2 & 7 & 5 & 2 & 5 \\ 0 & 0 & 6 & 2 & 1 & 3 \\ \hline 4 & 2 & 7 & 5 & 2 & 5 \\ 0 & 0 & 6 & 2 & 1 & 3 \\ 4 & 2 & 1 & 3 & 1 & 2 \end{bmatrix} \cup$$

$$\begin{bmatrix} 5 & 5 & 3 & 4 & 4 \\ 1 & 1 & 1 & 1 & 1 \\ 4 & 4 & 2 & 3 & 3 \\ \hline 6 & 6 & 0 & 3 & 3 \\ 3 & 3 & 1 & 2 & 2 \\ 9 & 9 & 5 & 7 & 7 \end{bmatrix} + \begin{bmatrix} 7 & 11 & 3 & 1 & 6 \\ 5 & 1 & 2 & 0 & 4 \\ 7 & 11 & 3 & 1 & 6 \\ \hline 11 & 9 & 2 & 2 & 4 \\ 5 & 13 & 6 & 0 & 14 \\ 6 & 1 & 3 & 0 & 7 \end{bmatrix} =$$

$$\begin{bmatrix} 9 & 22 & 19 & 14 & 5 & 11 \\ 26 & 22 & 36 & 37 & 12 & 38 \\ 5 & 9 & 18 & 9 & 4 & 14 \\ 17 & 13 & 12 & 21 & 6 & 14 \\ \hline 10 & 17 & 16 & 16 & 6 & 15 \\ 13 & 29 & 24 & 11 & 7 & 15 \\ \hline 7 & 10 & 12 & 21 & 5 & 16 \\ 6 & 7 & 13 & 9 & 4 & 16 \\ 7 & 4 & 5 & 15 & 5 & 19 \end{bmatrix} \cup \begin{bmatrix} 12 & 16 & 6 & 5 & 10 \\ 6 & 2 & 3 & 1 & 5 \\ 11 & 15 & 5 & 4 & 9 \\ \hline 17 & 15 & 2 & 5 & 7 \\ 8 & 16 & 7 & 2 & 16 \\ 15 & 10 & 8 & 7 & 14 \end{bmatrix}$$

$= S_1 \cup S_2 = AB$ is once again a superbimatrix.

**DEFINITION 2.21:** *Let $A = A_1 \cup A_2$ and $B = B_1 \cup B_2$ be two superbimatrices. The minor byproduct of the two superbimatrices $AB = (A_1 \cup A_2)(B_1 \cup B_2) = A_1 B_1 \cup A_2 B_2$ is defined if and only if the following conditions are satisfied.*

1. *The number of columns in $A_i$ is equal to the number of rows in $B_i$; $i = 1, 2$.*

2. *The partition of $A_i$ along the columns is equal or identical with the partition of $B_i$ along the rows $i = 1, 2$.*

We see the minor byproduct of AB when it exists for the superbimatrices A and B is again a superbimatrix. The examples 2.55 and 2.56 show explicitly how this minor byproduct of any two superbimatrices are defined resulting in a superbimatrix. Now we proceed onto first illustrate by examples how the product of the transpose of a superbimatrix with a superbimatrix is defined.

***Example 2.57:*** Let $A = A_1 \cup A_2$ be a superbimatrix where

$$A_1 = \begin{bmatrix} 1 & 0 & 1 & 0 & 1 & 0 \\ 2 & 1 & 0 & 2 & 1 & 2 \\ \hline 0 & 2 & 1 & 0 & 2 & 1 \\ 1 & 1 & 2 & 1 & 2 & 0 \\ 1 & 1 & 0 & 2 & 1 & 0 \\ 2 & 0 & 1 & 2 & 0 & 1 \\ \hline 5 & 1 & 2 & 1 & 0 & 1 \end{bmatrix}$$

and

$$A_2 = \begin{bmatrix} 1 & 1 & 3 & 1 & 0 & 1 & 4 \\ 5 & 1 & 0 & 2 & 0 & 1 & 0 \\ \hline 7 & 2 & 1 & 0 & 1 & 0 & 1 \\ 8 & 1 & 2 & 1 & 1 & 0 & 1 \\ \hline 2 & 0 & 1 & 0 & 2 & 2 & 1 \\ \hline 3 & 5 & 7 & 0 & 1 & 6 & 0 \end{bmatrix}.$$

Now
$$\begin{aligned} A^T &= (A_1 \cup A_2)^T \\ &= (A_1^T \cup A_2^T) \end{aligned}$$

$$= \begin{bmatrix} 1 & 2 & 0 & 1 & 1 & 2 & 5 \\ \hline 0 & 1 & 2 & 1 & 1 & 0 & 1 \\ 1 & 0 & 1 & 2 & 0 & 1 & 2 \\ \hline 0 & 2 & 0 & 1 & 2 & 2 & 1 \\ 1 & 1 & 2 & 2 & 1 & 0 & 0 \\ 0 & 2 & 1 & 0 & 0 & 1 & 1 \end{bmatrix} \cup \begin{bmatrix} 1 & 5 & 7 & 8 & 2 & 3 \\ 1 & 1 & 2 & 1 & 0 & 5 \\ \hline 3 & 0 & 1 & 2 & 1 & 7 \\ 1 & 2 & 0 & 1 & 0 & 0 \\ 0 & 0 & 1 & 1 & 2 & 1 \\ \hline 1 & 1 & 0 & 0 & 2 & 6 \\ 4 & 0 & 1 & 1 & 1 & 0 \end{bmatrix}.$$

Now

$$\begin{aligned} AA^T &= (A_1 \cup A_2)(A_1 \cup A_2)^T \\ &= (A_1 \cup A_2)(A_1^T \cup A_2^T) \\ &= A_1 A_1^T \cup A_2 A_2^T \end{aligned}$$

$$= \begin{bmatrix} 1 & 0 & 1 & 0 & 1 & 0 \\ 2 & 1 & 0 & 2 & 1 & 2 \\ \hline 0 & 2 & 1 & 0 & 2 & 1 \\ 1 & 1 & 2 & 1 & 2 & 0 \\ 1 & 1 & 0 & 2 & 1 & 0 \\ 2 & 0 & 1 & 2 & 0 & 1 \\ \hline 5 & 1 & 2 & 1 & 0 & 1 \end{bmatrix} \begin{bmatrix} 1 & 2 & 0 & 1 & 1 & 2 & 5 \\ 0 & 1 & 2 & 1 & 1 & 0 & 1 \\ \hline 1 & 0 & 1 & 2 & 0 & 1 & 2 \\ 0 & 2 & 0 & 1 & 2 & 2 & 1 \\ 1 & 1 & 2 & 2 & 1 & 0 & 0 \\ 0 & 2 & 1 & 0 & 0 & 1 & 1 \end{bmatrix} \cup$$

$$\begin{bmatrix} 1 & 1 & 3 & 1 & 0 & 1 & 4 \\ 5 & 1 & 0 & 2 & 0 & 1 & 0 \\ 7 & 2 & 1 & 0 & 1 & 0 & 1 \\ \hline 8 & 1 & 2 & 1 & 1 & 0 & 1 \\ 2 & 0 & 1 & 0 & 2 & 2 & 1 \\ \hline 3 & 5 & 7 & 0 & 1 & 6 & 0 \end{bmatrix} \begin{bmatrix} 1 & 5 & 7 & 8 & 2 & 3 \\ 1 & 1 & 2 & 1 & 0 & 5 \\ \hline 3 & 0 & 1 & 2 & 1 & 7 \\ 1 & 2 & 0 & 1 & 0 & 0 \\ 0 & 0 & 1 & 1 & 2 & 1 \\ 1 & 1 & 0 & 0 & 2 & 6 \\ 4 & 0 & 1 & 1 & 1 & 0 \end{bmatrix}$$

$$= \begin{bmatrix} \begin{bmatrix}1\\2\end{bmatrix}\begin{bmatrix}1 & 2\end{bmatrix} & \begin{bmatrix}1\\2\end{bmatrix}\begin{bmatrix}0 & 1 & 1 & 2\end{bmatrix} & \begin{bmatrix}1\\2\end{bmatrix}[5] \\ \hline \begin{bmatrix}0\\1\\1\\2\end{bmatrix}\begin{bmatrix}1 & 2\end{bmatrix} & \begin{bmatrix}0\\1\\1\\2\end{bmatrix}\begin{bmatrix}0 & 1 & 1 & 2\end{bmatrix} & \begin{bmatrix}0\\1\\1\\2\end{bmatrix}[5] \\ \hline [5]\begin{bmatrix}1 & 2\end{bmatrix} & [5]\begin{bmatrix}0 & 1 & 1 & 2\end{bmatrix} & [5][5] \end{bmatrix} +$$

$$\begin{bmatrix} \begin{pmatrix} 0 & 1 \\ 1 & 0 \end{pmatrix}\begin{pmatrix} 0 & 1 \\ 1 & 0 \end{pmatrix} & \begin{pmatrix} 0 & 1 \\ 1 & 0 \end{pmatrix}\begin{pmatrix} 2 & 1 & 1 & 0 \\ 1 & 2 & 0 & 1 \end{pmatrix} & \begin{pmatrix} 0 & 1 \\ 1 & 0 \end{pmatrix}\begin{pmatrix} 1 \\ 2 \end{pmatrix} \\ \hline \begin{pmatrix} 2 & 1 \\ 1 & 2 \\ 1 & 0 \\ 0 & 1 \end{pmatrix}\begin{pmatrix} 0 & 1 \\ 1 & 0 \end{pmatrix} & \begin{pmatrix} 2 & 1 \\ 1 & 2 \\ 1 & 0 \\ 0 & 1 \end{pmatrix}\begin{pmatrix} 2 & 1 & 1 & 0 \\ 1 & 2 & 0 & 1 \end{pmatrix} & \begin{pmatrix} 2 & 1 \\ 1 & 2 \\ 1 & 0 \\ 0 & 1 \end{pmatrix}\begin{pmatrix} 1 \\ 2 \end{pmatrix} \\ \hline \begin{pmatrix} 1 & 2 \end{pmatrix}\begin{pmatrix} 0 & 1 \\ 1 & 0 \end{pmatrix} & \begin{pmatrix} 1 & 2 \end{pmatrix}\begin{pmatrix} 2 & 1 & 1 & 0 \\ 1 & 2 & 0 & 1 \end{pmatrix} & \begin{pmatrix} 1 & 2 \end{pmatrix}\begin{pmatrix} 1 \\ 2 \end{pmatrix} \end{bmatrix} +$$

$$\begin{bmatrix} \begin{bmatrix} 0 & 1 & 0 \\ 2 & 1 & 2 \end{bmatrix}\begin{bmatrix} 0 & 2 \\ 1 & 1 \\ 0 & 2 \end{bmatrix} & \begin{bmatrix} 0 & 1 & 0 \\ 2 & 1 & 2 \end{bmatrix}\begin{bmatrix} 0 & 1 & 2 & 2 \\ 2 & 2 & 1 & 0 \\ 1 & 0 & 0 & 1 \end{bmatrix} & \begin{bmatrix} 0 & 1 & 0 \\ 2 & 1 & 2 \end{bmatrix}\begin{bmatrix} 1 \\ 0 \\ 1 \end{bmatrix} \\ \hline \begin{bmatrix} 0 & 2 & 1 \\ 1 & 2 & 0 \\ 2 & 1 & 0 \\ 2 & 0 & 1 \end{bmatrix}\begin{bmatrix} 0 & 2 \\ 1 & 1 \\ 0 & 2 \end{bmatrix} & \begin{bmatrix} 0 & 2 & 1 \\ 1 & 2 & 0 \\ 2 & 1 & 0 \\ 2 & 0 & 1 \end{bmatrix}\begin{bmatrix} 0 & 1 & 2 & 2 \\ 2 & 2 & 1 & 0 \\ 1 & 0 & 0 & 1 \end{bmatrix} & \begin{bmatrix} 0 & 2 & 1 \\ 1 & 2 & 0 \\ 2 & 1 & 0 \\ 2 & 0 & 1 \end{bmatrix}\begin{bmatrix} 1 \\ 0 \\ 1 \end{bmatrix} \\ \hline \begin{bmatrix} 1 & 0 & 1 \end{bmatrix}\begin{bmatrix} 0 & 2 \\ 1 & 1 \\ 0 & 2 \end{bmatrix} & \begin{bmatrix} 1 & 0 & 1 \end{bmatrix}\begin{bmatrix} 0 & 1 & 2 & 2 \\ 2 & 2 & 1 & 0 \\ 1 & 0 & 0 & 1 \end{bmatrix} & \begin{bmatrix} 1 & 0 & 1 \end{bmatrix}\begin{bmatrix} 1 \\ 0 \\ 1 \end{bmatrix} \end{bmatrix}$$

$$\cup \begin{bmatrix} \begin{pmatrix} 1 & 1 \\ 5 & 1 \\ 7 & 2 \end{pmatrix}\begin{pmatrix} 1 & 5 & 7 \\ 1 & 1 & 2 \end{pmatrix} & \begin{pmatrix} 1 & 1 \\ 5 & 1 \\ 7 & 2 \end{pmatrix}\begin{pmatrix} 8 & 2 \\ 1 & 0 \end{pmatrix} & \begin{pmatrix} 1 & 1 \\ 5 & 1 \\ 7 & 2 \end{pmatrix}\begin{pmatrix} 3 \\ 5 \end{pmatrix} \\ \hline \begin{pmatrix} 8 & 1 \\ 2 & 0 \end{pmatrix}\begin{pmatrix} 1 & 5 & 7 \\ 1 & 1 & 2 \end{pmatrix} & \begin{pmatrix} 8 & 1 \\ 2 & 0 \end{pmatrix}\begin{pmatrix} 8 & 2 \\ 1 & 0 \end{pmatrix} & \begin{pmatrix} 8 & 1 \\ 2 & 0 \end{pmatrix}\begin{pmatrix} 3 \\ 5 \end{pmatrix} \\ \hline \begin{pmatrix} 3 & 5 \end{pmatrix}\begin{pmatrix} 1 & 5 & 7 \\ 1 & 1 & 2 \end{pmatrix} & \begin{pmatrix} 3 & 5 \end{pmatrix}\begin{pmatrix} 8 & 2 \\ 1 & 0 \end{pmatrix} & \begin{pmatrix} 3 & 5 \end{pmatrix}\begin{pmatrix} 3 \\ 5 \end{pmatrix} \end{bmatrix} +$$

$$\begin{bmatrix} \begin{bmatrix} 3 & 1 & 0 \\ 0 & 2 & 0 \\ 1 & 0 & 1 \end{bmatrix}\begin{bmatrix} 3 & 0 & 1 \\ 1 & 2 & 0 \\ 0 & 0 & 1 \end{bmatrix} & \begin{bmatrix} 3 & 1 & 0 \\ 0 & 2 & 0 \\ 1 & 0 & 1 \end{bmatrix}\begin{bmatrix} 2 & 1 \\ 1 & 0 \\ 1 & 2 \end{bmatrix} & \begin{bmatrix} 3 & 1 & 0 \\ 0 & 2 & 0 \\ 1 & 0 & 1 \end{bmatrix}\begin{bmatrix} 7 \\ 0 \\ 1 \end{bmatrix} \\ \begin{bmatrix} 2 & 1 & 1 \\ 1 & 0 & 2 \end{bmatrix}\begin{bmatrix} 3 & 0 & 1 \\ 1 & 2 & 0 \\ 0 & 0 & 1 \end{bmatrix} & \begin{bmatrix} 2 & 1 & 1 \\ 1 & 0 & 2 \end{bmatrix}\begin{bmatrix} 2 & 1 \\ 1 & 0 \\ 1 & 2 \end{bmatrix} & \begin{bmatrix} 2 & 1 & 1 \\ 1 & 0 & 2 \end{bmatrix}\begin{bmatrix} 7 \\ 0 \\ 1 \end{bmatrix} \\ \begin{bmatrix} 7 & 0 & 1 \end{bmatrix}\begin{bmatrix} 3 & 0 & 1 \\ 1 & 2 & 0 \\ 0 & 0 & 1 \end{bmatrix} & \begin{bmatrix} 7 & 0 & 1 \end{bmatrix}\begin{bmatrix} 2 & 1 \\ 1 & 0 \\ 1 & 2 \end{bmatrix} & \begin{bmatrix} 7 & 0 & 1 \end{bmatrix}\begin{bmatrix} 7 \\ 0 \\ 1 \end{bmatrix} \end{bmatrix} +$$

$$\begin{bmatrix} \begin{pmatrix} 1 & 4 \\ 1 & 0 \\ 0 & 1 \end{pmatrix}\begin{pmatrix} 1 & 1 & 0 \\ 4 & 0 & 1 \end{pmatrix} & \begin{pmatrix} 1 & 4 \\ 1 & 0 \\ 0 & 1 \end{pmatrix}\begin{pmatrix} 0 & 2 \\ 1 & 1 \end{pmatrix} & \begin{pmatrix} 1 & 4 \\ 1 & 0 \\ 0 & 1 \end{pmatrix}\begin{pmatrix} 6 \\ 0 \end{pmatrix} \\ \begin{pmatrix} 0 & 1 \\ 2 & 1 \end{pmatrix}\begin{pmatrix} 1 & 1 & 0 \\ 4 & 0 & 1 \end{pmatrix} & \begin{pmatrix} 0 & 1 \\ 2 & 1 \end{pmatrix}\begin{pmatrix} 0 & 2 \\ 1 & 1 \end{pmatrix} & \begin{pmatrix} 0 & 1 \\ 2 & 1 \end{pmatrix}\begin{pmatrix} 6 \\ 0 \end{pmatrix} \\ \begin{pmatrix} 6 & 0 \end{pmatrix}\begin{pmatrix} 1 & 1 & 0 \\ 4 & 0 & 1 \end{pmatrix} & \begin{pmatrix} 6 & 0 \end{pmatrix}\begin{pmatrix} 0 & 2 \\ 1 & 1 \end{pmatrix} & \begin{pmatrix} 6 & 0 \end{pmatrix}\begin{pmatrix} 6 \\ 0 \end{pmatrix} \end{bmatrix}$$

$$= \left[\begin{array}{cc|cccc|c} 1 & 2 & 0 & 1 & 1 & 2 & 5 \\ 2 & 4 & 0 & 2 & 2 & 4 & 10 \\ \hline 0 & 0 & 0 & 0 & 0 & 0 & 0 \\ 1 & 2 & 0 & 1 & 1 & 2 & 5 \\ 1 & 2 & 0 & 1 & 1 & 2 & 5 \\ 2 & 4 & 0 & 2 & 2 & 4 & 10 \\ \hline 5 & 10 & 0 & 5 & 5 & 10 & 25 \end{array}\right] + \left[\begin{array}{cc|cccc|c} 1 & 0 & 1 & 2 & 0 & 1 & 2 \\ 0 & 1 & 2 & 1 & 1 & 0 & 1 \\ \hline 1 & 2 & 5 & 4 & 2 & 1 & 4 \\ 2 & 1 & 4 & 5 & 1 & 2 & 5 \\ 0 & 1 & 2 & 1 & 1 & 0 & 1 \\ 1 & 0 & 1 & 2 & 0 & 1 & 2 \\ \hline 2 & 1 & 4 & 5 & 1 & 2 & 5 \end{array}\right] +$$

$$\begin{bmatrix} 1 & 1 & 2 & 2 & 1 & 0 & 0 \\ 1 & 9 & 4 & 4 & 5 & 6 & 4 \\ \hline 2 & 4 & 5 & 4 & 2 & 1 & 1 \\ 2 & 4 & 4 & 5 & 4 & 2 & 1 \\ 1 & 5 & 2 & 4 & 5 & 4 & 2 \\ 0 & 6 & 1 & 2 & 4 & 5 & 3 \\ \hline 0 & 4 & 1 & 1 & 2 & 3 & 2 \end{bmatrix} \cup \begin{bmatrix} 2 & 6 & 9 & 9 & 2 & 8 \\ 6 & 26 & 37 & 41 & 10 & 20 \\ 9 & 37 & 53 & 58 & 14 & 31 \\ \hline 9 & 41 & 58 & 65 & 16 & 29 \\ 2 & 10 & 14 & 16 & 4 & 6 \\ \hline 8 & 20 & 31 & 29 & 6 & 34 \end{bmatrix} +$$

$$\begin{bmatrix} 10 & 2 & 3 & 7 & 3 & 21 \\ 2 & 4 & 0 & 2 & 0 & 0 \\ 3 & 0 & 2 & 3 & 3 & 8 \\ \hline 7 & 2 & 3 & 6 & 4 & 15 \\ 3 & 0 & 3 & 4 & 5 & 9 \\ \hline 21 & 0 & 8 & 15 & 9 & 50 \end{bmatrix} + \begin{bmatrix} 17 & 1 & 4 & 4 & 6 & 6 \\ 1 & 1 & 0 & 0 & 2 & 6 \\ 4 & 0 & 1 & 1 & 1 & 0 \\ \hline 4 & 0 & 1 & 1 & 1 & 0 \\ 6 & 2 & 1 & 1 & 5 & 12 \\ \hline 6 & 6 & 0 & 0 & 12 & 36 \end{bmatrix} =$$

$$\begin{bmatrix} 3 & 3 & 3 & 5 & 2 & 3 & 7 \\ 3 & 14 & 6 & 7 & 8 & 10 & 15 \\ \hline 3 & 6 & 10 & 8 & 4 & 2 & 5 \\ 5 & 7 & 8 & 11 & 6 & 6 & 11 \\ 2 & 8 & 4 & 6 & 7 & 6 & 8 \\ 3 & 10 & 2 & 6 & 6 & 10 & 15 \\ \hline 7 & 15 & 5 & 11 & 8 & 15 & 32 \end{bmatrix} \cup$$

$$\begin{bmatrix} 29 & 9 & 16 & 20 & 11 & 35 \\ 9 & 31 & 37 & 43 & 12 & 26 \\ 16 & 37 & 56 & 62 & 18 & 39 \\ \hline 20 & 43 & 62 & 72 & 21 & 44 \\ 11 & 12 & 18 & 21 & 14 & 27 \\ \hline 35 & 26 & 39 & 44 & 27 & 120 \end{bmatrix}.$$

We see the resultant is a symmetric superbimatrix. We see by the minor byproduct of superbimatrices with their respective transposes we can get more and more symmetric superbimatrices.

*Example 2.58:* Let $A = A_1 \cup A_2$ be a superbimatrix where

$$A_1 = \begin{bmatrix} 3 & 1 & 1 & 3 \\ 1 & 0 & 0 & 1 \\ \hline 2 & 1 & 0 & 2 \\ 5 & 2 & 1 & 1 \\ \hline 0 & 3 & 0 & 0 \end{bmatrix}$$

and

$$A_2 = \begin{bmatrix} 8 & 1 & 3 & 1 & 0 \\ \hline 0 & 0 & 1 & 2 & 5 \\ 1 & 0 & 0 & 1 & 0 \\ \hline 2 & 1 & 2 & 1 & 2 \\ \hline 5 & 1 & 0 & 0 & 3 \\ 6 & 0 & 1 & 1 & 0 \end{bmatrix}.$$

$A^T = (A_1 \cup A_2)^T = (A_1^T \cup A_2^T)$

$$= \begin{bmatrix} 3 & 1 & 2 & 5 & 0 \\ \hline 1 & 0 & 1 & 2 & 3 \\ \hline 1 & 0 & 0 & 1 & 0 \\ \hline 3 & 1 & 2 & 1 & 0 \end{bmatrix} \cup \begin{bmatrix} 8 & 0 & 1 & 2 & 5 & 6 \\ 1 & 0 & 0 & 1 & 1 & 0 \\ \hline 3 & 1 & 0 & 2 & 0 & 1 \\ 1 & 2 & 1 & 1 & 0 & 1 \\ \hline 0 & 5 & 0 & 2 & 3 & 0 \end{bmatrix}.$$

$$\begin{aligned} AA^T &= (A_1 \cup A_2)(A_1 \cup A_2)^T \\ &= (A_1 \cup A_2)(A_1^T \cup A_2^T) \\ &= A_1 A_1^T \cup A_2 A_2^T \end{aligned}$$

$$= \left[\begin{array}{ccc|c} 3 & 1 & 1 & 3 \\ 1 & 0 & 0 & 1 \\ \hline 2 & 1 & 0 & 2 \\ 5 & 2 & 1 & 1 \\ \hline 0 & 3 & 0 & 0 \end{array}\right] \left[\begin{array}{cc|ccc|c} 3 & 1 & 2 & 5 & 0 \\ \hline 1 & 0 & 1 & 2 & 3 \\ 1 & 0 & 0 & 1 & 0 \\ \hline 3 & 1 & 2 & 1 & 0 \end{array}\right] \cup$$

$$\left[\begin{array}{cc|cc|c} 8 & 1 & 3 & 1 & 0 \\ 0 & 0 & 1 & 2 & 5 \\ 1 & 0 & 0 & 1 & 0 \\ 2 & 1 & 2 & 1 & 2 \\ \hline 5 & 1 & 0 & 0 & 3 \\ 6 & 0 & 1 & 1 & 0 \end{array}\right] \left[\begin{array}{c|cc|cc|cc} 8 & 0 & 1 & 2 & 5 & 6 \\ \hline 1 & 0 & 0 & 1 & 1 & 0 \\ 3 & 1 & 0 & 2 & 0 & 1 \\ 1 & 2 & 1 & 1 & 0 & 1 \\ \hline 0 & 5 & 0 & 2 & 3 & 0 \end{array}\right]$$

$$= \left[\begin{array}{c|c|c} \begin{bmatrix}3\\1\end{bmatrix}\begin{bmatrix}3 & 1\end{bmatrix} & \begin{bmatrix}3\\1\end{bmatrix}\begin{bmatrix}2 & 5\end{bmatrix} & \begin{bmatrix}3\\1\end{bmatrix}[0] \\ \hline \begin{bmatrix}2\\5\end{bmatrix}\begin{bmatrix}3 & 1\end{bmatrix} & \begin{bmatrix}2\\5\end{bmatrix}\begin{bmatrix}2 & 5\end{bmatrix} & \begin{bmatrix}2\\5\end{bmatrix}[0] \\ \hline [0]\begin{bmatrix}3 & 1\end{bmatrix} & [0]\begin{bmatrix}2 & 5\end{bmatrix} & [0][0] \end{array}\right] +$$

$$\left[\begin{array}{c|c|c} \begin{pmatrix}1 & 1\\0 & 0\end{pmatrix}\begin{pmatrix}1 & 0\\1 & 0\end{pmatrix} & \begin{pmatrix}1 & 1\\0 & 0\end{pmatrix}\begin{pmatrix}1 & 2\\0 & 1\end{pmatrix} & \begin{pmatrix}1 & 1\\0 & 0\end{pmatrix}\begin{pmatrix}3\\0\end{pmatrix} \\ \hline \begin{pmatrix}1 & 0\\2 & 1\end{pmatrix}\begin{pmatrix}1 & 0\\1 & 0\end{pmatrix} & \begin{pmatrix}1 & 0\\2 & 1\end{pmatrix}\begin{pmatrix}1 & 2\\0 & 1\end{pmatrix} & \begin{pmatrix}1 & 0\\2 & 1\end{pmatrix}\begin{pmatrix}3\\0\end{pmatrix} \\ \hline \begin{pmatrix}3 & 0\end{pmatrix}\begin{pmatrix}1 & 0\\1 & 0\end{pmatrix} & \begin{pmatrix}3 & 0\end{pmatrix}\begin{pmatrix}1 & 2\\0 & 1\end{pmatrix} & \begin{pmatrix}3 & 0\end{pmatrix}\begin{pmatrix}3\\0\end{pmatrix} \end{array}\right] +$$

$$\left[\begin{array}{c|c|c} \begin{bmatrix}3\\1\end{bmatrix}\begin{bmatrix}3 & 1\end{bmatrix} & \begin{bmatrix}3\\1\end{bmatrix}\begin{bmatrix}2 & 1\end{bmatrix} & \begin{bmatrix}3\\1\end{bmatrix}[0] \\ \hline \begin{bmatrix}2\\1\end{bmatrix}\begin{bmatrix}3 & 1\end{bmatrix} & \begin{bmatrix}2\\1\end{bmatrix}\begin{bmatrix}2 & 1\end{bmatrix} & \begin{bmatrix}2\\1\end{bmatrix}[0] \\ \hline [0]\begin{bmatrix}3 & 1\end{bmatrix} & [0]\begin{bmatrix}2 & 1\end{bmatrix} & [0][0] \end{array}\right] \cup$$

$$\left[\begin{array}{c|c|c} (8\ 1)\begin{pmatrix}8\\1\end{pmatrix} & (8\ 1)\begin{pmatrix}0 & 1 & 2\\0 & 0 & 1\end{pmatrix} & (8\ 1)\begin{pmatrix}5 & 6\\1 & 0\end{pmatrix} \\ \hline \begin{pmatrix}0 & 0\\1 & 0\\2 & 1\end{pmatrix}\begin{pmatrix}8\\1\end{pmatrix} & \begin{pmatrix}0 & 0\\1 & 0\\2 & 1\end{pmatrix}\begin{pmatrix}0 & 1 & 2\\0 & 0 & 1\end{pmatrix} & \begin{pmatrix}0 & 0\\1 & 0\\2 & 1\end{pmatrix}\begin{pmatrix}5 & 6\\1 & 0\end{pmatrix} \\ \hline \begin{pmatrix}5 & 1\\6 & 0\end{pmatrix}\begin{pmatrix}8\\1\end{pmatrix} & \begin{pmatrix}5 & 1\\6 & 0\end{pmatrix}\begin{pmatrix}0 & 1 & 2\\0 & 0 & 1\end{pmatrix} & \begin{pmatrix}5 & 1\\6 & 0\end{pmatrix}\begin{pmatrix}5 & 6\\1 & 0\end{pmatrix} \end{array}\right] +$$

$$\left[\begin{array}{c|c|c} (3\ 1)\begin{pmatrix}3\\1\end{pmatrix} & (3\ 1)\begin{pmatrix}1 & 0 & 2\\2 & 1 & 1\end{pmatrix} & (3\ 1)\begin{pmatrix}0 & 1\\0 & 1\end{pmatrix} \\ \hline \begin{pmatrix}1 & 2\\0 & 1\\2 & 1\end{pmatrix}\begin{pmatrix}3\\1\end{pmatrix} & \begin{pmatrix}1 & 2\\0 & 1\\2 & 1\end{pmatrix}\begin{pmatrix}1 & 0 & 2\\2 & 1 & 1\end{pmatrix} & \begin{pmatrix}1 & 2\\0 & 1\\2 & 1\end{pmatrix}\begin{pmatrix}0 & 1\\0 & 1\end{pmatrix} \\ \hline \begin{pmatrix}0 & 0\\1 & 1\end{pmatrix}\begin{pmatrix}3\\1\end{pmatrix} & \begin{pmatrix}0 & 0\\1 & 1\end{pmatrix}\begin{pmatrix}1 & 0 & 2\\2 & 1 & 1\end{pmatrix} & \begin{pmatrix}0 & 0\\1 & 1\end{pmatrix}\begin{pmatrix}0 & 1\\0 & 1\end{pmatrix} \end{array}\right] +$$

$$\left[\begin{array}{c|c|c} [0][0] & 0(5\ 0\ 2) & [0](3\ 0) \\ \hline \begin{bmatrix}5\\0\\2\end{bmatrix}[0] & \begin{bmatrix}5\\0\\2\end{bmatrix}(5\ 0\ 2) & \begin{bmatrix}5\\0\\2\end{bmatrix}(3\ 0) \\ \hline \begin{pmatrix}3\\0\end{pmatrix}0 & \begin{pmatrix}3\\0\end{pmatrix}(5\ 0\ 2) & \begin{pmatrix}3\\0\end{pmatrix}(3\ 0) \end{array}\right] =$$

$$\begin{bmatrix} 9 & 3 & 6 & 15 & 0 \\ 3 & 1 & 2 & 5 & 0 \\ \hline 6 & 2 & 4 & 10 & 0 \\ 15 & 5 & 10 & 25 & 0 \\ \hline 0 & 0 & 0 & 0 & 0 \end{bmatrix} + \begin{bmatrix} 2 & 0 & 1 & 3 & 3 \\ 0 & 0 & 0 & 0 & 0 \\ \hline 1 & 0 & 1 & 2 & 3 \\ 3 & 0 & 2 & 5 & 6 \\ \hline 3 & 0 & 3 & 6 & 9 \end{bmatrix} + \begin{bmatrix} 9 & 3 & 6 & 3 & 0 \\ 3 & 1 & 2 & 1 & 0 \\ \hline 6 & 2 & 4 & 2 & 0 \\ 3 & 1 & 2 & 1 & 0 \\ \hline 0 & 0 & 0 & 0 & 0 \end{bmatrix}$$

$$\cup \begin{bmatrix} 65 & 0 & 8 & 17 & 41 & 48 \\ \hline 0 & 0 & 0 & 0 & 0 & 0 \\ 8 & 0 & 1 & 2 & 5 & 6 \\ 17 & 0 & 2 & 5 & 11 & 12 \\ \hline 41 & 0 & 5 & 11 & 26 & 30 \\ 48 & 0 & 6 & 12 & 30 & 36 \end{bmatrix} + \begin{bmatrix} 10 & 5 & 1 & 7 & 0 & 4 \\ \hline 5 & 5 & 2 & 4 & 0 & 3 \\ 1 & 2 & 1 & 1 & 0 & 1 \\ 7 & 4 & 1 & 5 & 0 & 3 \\ \hline 0 & 0 & 0 & 0 & 0 & 0 \\ 4 & 3 & 1 & 3 & 0 & 2 \end{bmatrix} +$$

$$\begin{bmatrix} 0 & 0 & 0 & 0 & 0 & 0 \\ \hline 0 & 25 & 0 & 10 & 15 & 0 \\ 0 & 0 & 0 & 0 & 0 & 0 \\ 0 & 10 & 0 & 4 & 6 & 0 \\ \hline 0 & 15 & 0 & 6 & 9 & 0 \\ 0 & 0 & 0 & 0 & 0 & 0 \end{bmatrix}$$

$$= \begin{bmatrix} 20 & 6 & 13 & 21 & 3 \\ 6 & 2 & 4 & 6 & 0 \\ \hline 13 & 4 & 9 & 14 & 3 \\ 21 & 6 & 14 & 31 & 6 \\ \hline 3 & 0 & 3 & 6 & 9 \end{bmatrix} \cup \begin{bmatrix} 75 & 5 & 9 & 24 & 41 & 52 \\ \hline 5 & 30 & 2 & 14 & 15 & 3 \\ 9 & 2 & 2 & 3 & 5 & 7 \\ 24 & 14 & 3 & 14 & 17 & 15 \\ \hline 41 & 15 & 5 & 17 & 35 & 30 \\ 52 & 3 & 7 & 15 & 30 & 38 \end{bmatrix}.$$

We see the resultant is symmetric superbimatrix.

**DEFINITION 2.22:** *Let $A = A_1 \cup A_2$ and $B = B_1 \cup B_2$ be any two semi superbimatrices. The minor byproduct $AB = (A_1 \cup A_2)(B_1 \cup B_2) = A_1 B_1 \cup A_2 B_2$ is defined if and only if*

1. *Number of rows of $B_1$ is equal to number of columns in $A_1$ ($A_1$ and $B_1$ are usual matrices).*

2. *Number of rows of the super matrix $A_2$ equal to the number of columns of $B_2$ and the vertical partition of $A_2$ and the horizontal partition of $B_2$ are identical; i.e., if in $A_2$ there is a partition between $r$ and $(r + 1)^{th}$ column then in $B_2$ we have a partition between $r$ and $(r + 1)^{th}$ row this is true of any $r$; $1 < r <$ number of columns in $A_2 =$ number of rows in $B_2$.*

Now we illustrate this by the following examples.

**Example 2.59:** Let $A = A_1 \cup A_2$ and $B = B_1 \cup B_2$ be any two semi superbimatrices. Here

$$A = A_1 \cup A_2$$

$$= \begin{bmatrix} 9 & 0 & 1 & 2 & 3 & 5 \\ 0 & 2 & 0 & 1 & 2 & 1 \\ 1 & 0 & 1 & 0 & 0 & 1 \\ 2 & 1 & 0 & 1 & 1 & 0 \\ 7 & 0 & 2 & 1 & 2 & 1 \end{bmatrix} \cup \left[\begin{array}{cc|cccc|c} 3 & 1 & 2 & 0 & 1 & 1 & 5 & 0 \\ 6 & 0 & 1 & 1 & 0 & 0 & 1 & 2 \\ \hline 1 & 1 & 0 & 1 & 1 & 0 & 0 & 1 \\ 2 & 0 & 0 & 1 & 0 & 1 & 0 & 1 \\ \hline 3 & 2 & 1 & 0 & 1 & 0 & 2 & 0 \\ 4 & 0 & 2 & 1 & 0 & 1 & 0 & 1 \end{array}\right].$$

$$B = B_1 \cup B_2 = \begin{bmatrix} 0 & 0 & 1 & 1 & 2 \\ 1 & 0 & 1 & 0 & 1 \\ 2 & 1 & 0 & 1 & 0 \\ 1 & 0 & 1 & 0 & 1 \\ 0 & 1 & 1 & 0 & 2 \\ 1 & 0 & 1 & 1 & 0 \end{bmatrix} \cup \left[\begin{array}{ccc|cc} 1 & 2 & 1 & 2 & 0 \\ 2 & 0 & 2 & 1 & 1 \\ \hline 0 & 1 & 1 & 0 & 1 \\ 1 & 0 & 0 & 1 & 0 \\ 0 & 1 & 1 & 1 & 1 \\ 1 & 0 & 1 & 0 & 0 \\ 1 & 0 & 1 & 1 & 0 \\ 2 & 1 & 1 & 0 & 1 \end{array}\right].$$

AB $= (A_1 \cup A_2)(B_1 \cup B_2)$

$$= \quad A_1 B_1 \cup A_2 B_2$$

$$= \begin{bmatrix} 9 & 0 & 1 & 2 & 3 & 5 \\ 0 & 2 & 0 & 1 & 2 & 1 \\ 1 & 0 & 1 & 0 & 0 & 1 \\ 2 & 1 & 0 & 1 & 1 & 0 \\ 7 & 0 & 2 & 1 & 2 & 1 \end{bmatrix} \begin{bmatrix} 0 & 0 & 1 & 1 & 2 \\ 1 & 0 & 1 & 0 & 1 \\ 2 & 1 & 0 & 1 & 0 \\ 1 & 0 & 1 & 0 & 1 \\ 0 & 1 & 1 & 0 & 2 \\ 1 & 0 & 1 & 1 & 0 \end{bmatrix} \cup$$

$$\begin{bmatrix} \begin{array}{cc|cccccc|c} 3 & 1 & 2 & 0 & 1 & 1 & 5 & 0 \\ 6 & 0 & 1 & 1 & 0 & 0 & 1 & 2 \\ \hline 1 & 1 & 0 & 1 & 1 & 0 & 0 & 1 \\ 2 & 0 & 0 & 1 & 0 & 1 & 0 & 1 \\ \hline 3 & 2 & 1 & 0 & 1 & 0 & 2 & 0 \\ 4 & 0 & 2 & 1 & 0 & 1 & 0 & 1 \end{array} \end{bmatrix} \begin{bmatrix} \begin{array}{ccc|cc} 1 & 2 & 1 & 2 & 0 \\ 2 & 0 & 2 & 1 & 1 \\ \hline 0 & 1 & 1 & 0 & 1 \\ 1 & 0 & 0 & 1 & 0 \\ 0 & 1 & 1 & 1 & 1 \\ 1 & 0 & 1 & 0 & 0 \\ 1 & 0 & 1 & 1 & 0 \\ \hline 2 & 1 & 1 & 0 & 1 \end{array} \end{bmatrix}$$

$$= \begin{bmatrix} 9 & 4 & 19 & 15 & 26 \\ 4 & 2 & 6 & 1 & 7 \\ 3 & 1 & 2 & 3 & 2 \\ 2 & 1 & 5 & 2 & 8 \\ 6 & 4 & 11 & 10 & 19 \end{bmatrix} \cup$$

$$\begin{bmatrix} \begin{pmatrix} 3 & 1 \\ 6 & 0 \end{pmatrix}\begin{pmatrix} 1 & 2 & 1 \\ 2 & 0 & 2 \end{pmatrix} & \begin{pmatrix} 3 & 1 \\ 6 & 0 \end{pmatrix}\begin{pmatrix} 2 & 0 \\ 1 & 1 \end{pmatrix} \\ \hline \begin{pmatrix} 1 & 1 \\ 2 & 0 \end{pmatrix}\begin{pmatrix} 1 & 2 & 1 \\ 2 & 0 & 2 \end{pmatrix} & \begin{pmatrix} 1 & 1 \\ 2 & 0 \end{pmatrix}\begin{pmatrix} 2 & 0 \\ 1 & 1 \end{pmatrix} \\ \hline \begin{pmatrix} 3 & 2 \\ 4 & 0 \end{pmatrix}\begin{pmatrix} 1 & 2 & 1 \\ 2 & 0 & 2 \end{pmatrix} & \begin{pmatrix} 3 & 2 \\ 4 & 0 \end{pmatrix}\begin{pmatrix} 2 & 0 \\ 1 & 1 \end{pmatrix} \end{bmatrix} +$$

$$\begin{bmatrix} \begin{pmatrix} 2 & 0 & 1 & 1 & 5 \\ 1 & 1 & 0 & 0 & 1 \end{pmatrix} \begin{pmatrix} 0 & 1 & 1 \\ 1 & 0 & 0 \\ 0 & 1 & 1 \\ 1 & 0 & 1 \\ 1 & 0 & 1 \end{pmatrix} & \begin{pmatrix} 2 & 0 & 1 & 1 & 5 \\ 1 & 1 & 0 & 0 & 1 \end{pmatrix} \begin{bmatrix} 0 & 1 \\ 1 & 0 \\ 1 & 1 \\ 0 & 0 \\ 1 & 0 \end{bmatrix} \\ \hline \begin{pmatrix} 0 & 1 & 1 & 0 & 0 \\ 0 & 1 & 0 & 1 & 0 \end{pmatrix} \begin{pmatrix} 0 & 1 & 1 \\ 1 & 0 & 0 \\ 0 & 1 & 1 \\ 1 & 0 & 1 \\ 1 & 0 & 1 \end{pmatrix} & \begin{pmatrix} 0 & 1 & 1 & 0 & 0 \\ 0 & 1 & 0 & 1 & 0 \end{pmatrix} \begin{bmatrix} 0 & 1 \\ 1 & 0 \\ 1 & 1 \\ 0 & 0 \\ 1 & 0 \end{bmatrix} \\ \hline \begin{pmatrix} 1 & 0 & 1 & 0 & 2 \\ 2 & 1 & 0 & 1 & 0 \end{pmatrix} \begin{pmatrix} 0 & 1 & 1 \\ 1 & 0 & 0 \\ 0 & 1 & 1 \\ 1 & 0 & 1 \\ 1 & 0 & 1 \end{pmatrix} & \begin{pmatrix} 1 & 0 & 1 & 0 & 2 \\ 2 & 1 & 0 & 1 & 0 \end{pmatrix} \begin{bmatrix} 0 & 1 \\ 1 & 0 \\ 1 & 1 \\ 0 & 0 \\ 1 & 0 \end{bmatrix} \end{bmatrix}$$

$$+ \begin{bmatrix} \begin{pmatrix} 0 \\ 2 \end{pmatrix} \begin{pmatrix} 2 & 1 & 1 \end{pmatrix} & \begin{pmatrix} 0 \\ 2 \end{pmatrix} \begin{pmatrix} 0 & 1 \end{pmatrix} \\ \hline \begin{pmatrix} 1 \\ 1 \end{pmatrix} \begin{pmatrix} 2 & 1 & 1 \end{pmatrix} & \begin{pmatrix} 1 \\ 1 \end{pmatrix} \begin{pmatrix} 0 & 1 \end{pmatrix} \\ \hline \begin{pmatrix} 0 \\ 1 \end{pmatrix} \begin{pmatrix} 2 & 1 & 1 \end{pmatrix} & \begin{pmatrix} 0 \\ 1 \end{pmatrix} \begin{pmatrix} 0 & 1 \end{pmatrix} \end{bmatrix} =$$

$$\begin{bmatrix} 9 & 4 & 19 & 15 & 26 \\ 4 & 2 & 6 & 1 & 7 \\ 3 & 1 & 2 & 3 & 2 \\ 2 & 1 & 5 & 2 & 8 \\ 6 & 4 & 11 & 10 & 19 \end{bmatrix} \cup \begin{bmatrix} 5 & 6 & 5 & 7 & 1 \\ 6 & 12 & 6 & 12 & 0 \\ \hline 3 & 2 & 3 & 3 & 1 \\ 2 & 4 & 2 & 4 & 0 \\ \hline 7 & 6 & 7 & 8 & 2 \\ 4 & 8 & 4 & 8 & 0 \end{bmatrix} +$$

$$\begin{bmatrix} 6 & 3 & 9 & 6 & 3 \\ 2 & 1 & 2 & 2 & 1 \\ \hline 1 & 1 & 1 & 2 & 1 \\ 2 & 0 & 1 & 1 & 0 \\ \hline 2 & 2 & 4 & 3 & 2 \\ 2 & 2 & 3 & 1 & 2 \end{bmatrix} + \begin{bmatrix} 0 & 0 & 0 & 0 & 0 \\ 4 & 2 & 2 & 0 & 2 \\ \hline 2 & 1 & 1 & 0 & 1 \\ 2 & 1 & 1 & 0 & 1 \\ \hline 0 & 0 & 0 & 0 & 0 \\ 2 & 1 & 1 & 0 & 1 \end{bmatrix} =$$

$$\begin{bmatrix} 9 & 4 & 19 & 15 & 26 \\ 4 & 2 & 6 & 1 & 7 \\ 3 & 1 & 2 & 3 & 2 \\ 2 & 1 & 5 & 2 & 8 \\ 6 & 4 & 11 & 10 & 19 \end{bmatrix} \cup \begin{bmatrix} 11 & 9 & 14 & 13 & 4 \\ 12 & 15 & 10 & 14 & 3 \\ \hline 6 & 4 & 5 & 5 & 3 \\ 6 & 5 & 4 & 5 & 1 \\ \hline 9 & 8 & 11 & 11 & 4 \\ 8 & 11 & 8 & 9 & 3 \end{bmatrix}$$

is a semi superbimatrix which is not symmetric.

**Example 2.60:** Let $A = A_1 \cup A_2$ and $B = B_1 \cup B_2$ be any two semi superbimatrices. We find their minor product.

Here

$$A = A_1 \cup A_2$$

$$= \begin{bmatrix} 0 & 3 & 1 & 2 \\ 1 & 0 & 2 & 0 \\ 0 & 1 & 1 & 1 \\ 1 & 0 & 1 & 0 \\ 2 & 1 & 0 & 1 \end{bmatrix} \cup \begin{bmatrix} 3 & 1 & 4 & 5 \\ 0 & 1 & 0 & 1 \\ 1 & 0 & 1 & 0 \\ \hline 1 & 0 & 0 & 1 \\ 0 & 1 & 1 & 0 \\ 1 & 1 & 1 & 0 \\ 1 & 0 & 1 & 1 \\ 0 & 1 & 1 & 1 \\ \hline 1 & 1 & 0 & 1 \\ 1 & 1 & 0 & 0 \end{bmatrix}$$

and
$$B = B_1 \cup B_2$$

$$= \begin{bmatrix} 1 & 4 & 3 & 2 \\ 2 & 1 & 0 & 1 \\ 3 & 0 & 1 & 0 \\ 4 & 2 & 1 & 2 \end{bmatrix} \cup \left[ \begin{array}{cccc|cc|ccc} 0 & 2 & 3 & 0 & 1 & 3 & 3 & 1 & 2 & 3 \\ 2 & 1 & 1 & 3 & 0 & 2 & 0 & 0 & 1 & 0 \\ \hline 3 & 0 & 2 & 1 & 3 & 1 & 1 & 1 & 0 & 1 \\ 1 & 3 & 0 & 2 & 2 & 0 & 2 & 1 & 2 & 0 \end{array} \right]$$

are the given semi superbimatrices.

$$AB = (A_1 \cup A_2)(B_1 \cup B_2)$$
$$= A_1 B_1 \cup A_2 B_2$$

$$= \begin{bmatrix} 0 & 3 & 1 & 2 \\ 1 & 0 & 2 & 0 \\ 0 & 1 & 1 & 1 \\ 1 & 0 & 1 & 0 \\ 2 & 1 & 0 & 1 \end{bmatrix} \begin{bmatrix} 1 & 4 & 3 & 2 \\ 2 & 1 & 0 & 1 \\ 3 & 0 & 1 & 0 \\ 4 & 2 & 1 & 2 \end{bmatrix} \cup$$

$$\left[ \begin{array}{cc|cc} 3 & 1 & 4 & 5 \\ 0 & 1 & 0 & 1 \\ 1 & 0 & 1 & 0 \\ \hline 1 & 0 & 0 & 1 \\ 0 & 1 & 1 & 0 \\ 1 & 1 & 1 & 0 \\ 1 & 0 & 1 & 1 \\ 0 & 1 & 1 & 1 \\ 1 & 1 & 0 & 1 \\ 1 & 1 & 0 & 0 \end{array} \right] \left[ \begin{array}{cccc|cc|ccc} 0 & 2 & 3 & 0 & 1 & 3 & 3 & 1 & 2 & 3 \\ 2 & 1 & 1 & 3 & 0 & 2 & 0 & 0 & 1 & 0 \\ \hline 3 & 0 & 2 & 1 & 3 & 1 & 1 & 1 & 0 & 1 \\ 1 & 3 & 0 & 2 & 2 & 0 & 2 & 1 & 2 & 0 \end{array} \right]$$

$$= \begin{bmatrix} 17 & 7 & 3 & 7 \\ 7 & 4 & 5 & 2 \\ 9 & 3 & 2 & 3 \\ 4 & 4 & 4 & 2 \\ 8 & 11 & 7 & 7 \end{bmatrix} \cup$$

$$\begin{bmatrix} \begin{pmatrix} 3 & 1 \\ 0 & 1 \\ 1 & 0 \end{pmatrix}\begin{pmatrix} 0 & 2 & 3 & 0 & 1 \\ 2 & 1 & 1 & 3 & 0 \end{pmatrix} & \begin{pmatrix} 3 & 1 \\ 0 & 1 \\ 1 & 0 \end{pmatrix}\begin{bmatrix} 3 & 3 \\ 2 & 0 \end{bmatrix} & \begin{pmatrix} 3 & 1 \\ 0 & 1 \\ 1 & 0 \end{pmatrix}\begin{pmatrix} 1 & 2 & 3 \\ 0 & 1 & 0 \end{pmatrix} \\ \begin{pmatrix} 1 & 0 \\ 0 & 1 \\ 1 & 1 \\ 1 & 0 \\ 0 & 1 \end{pmatrix}\begin{pmatrix} 0 & 2 & 3 & 0 & 1 \\ 2 & 1 & 1 & 3 & 0 \end{pmatrix} & \begin{pmatrix} 1 & 0 \\ 0 & 1 \\ 1 & 1 \\ 1 & 0 \\ 0 & 1 \end{pmatrix}\begin{bmatrix} 3 & 3 \\ 2 & 0 \end{bmatrix} & \begin{pmatrix} 1 & 0 \\ 0 & 1 \\ 1 & 1 \\ 1 & 0 \\ 0 & 1 \end{pmatrix}\begin{pmatrix} 1 & 2 & 3 \\ 0 & 1 & 0 \end{pmatrix} \\ \begin{pmatrix} 1 & 1 \\ 1 & 1 \end{pmatrix}\begin{pmatrix} 0 & 2 & 3 & 0 & 1 \\ 2 & 1 & 1 & 3 & 0 \end{pmatrix} & \begin{pmatrix} 1 & 1 \\ 1 & 1 \end{pmatrix}\begin{bmatrix} 3 & 3 \\ 2 & 0 \end{bmatrix} & \begin{pmatrix} 1 & 1 \\ 1 & 1 \end{pmatrix}\begin{pmatrix} 1 & 2 & 3 \\ 0 & 1 & 0 \end{pmatrix} \end{bmatrix}$$

$$+$$

$$\begin{bmatrix} \begin{bmatrix} 4 & 5 \\ 0 & 1 \\ 1 & 0 \end{bmatrix}\begin{bmatrix} 3 & 0 & 2 & 1 & 3 \\ 1 & 3 & 0 & 2 & 2 \end{bmatrix} & \begin{bmatrix} 4 & 5 \\ 0 & 1 \\ 1 & 0 \end{bmatrix}\begin{pmatrix} 1 & 1 \\ 0 & 2 \end{pmatrix} & \begin{bmatrix} 4 & 5 \\ 0 & 1 \\ 1 & 0 \end{bmatrix}\begin{pmatrix} 1 & 0 & 1 \\ 1 & 2 & 0 \end{pmatrix} \\ \begin{pmatrix} 0 & 1 \\ 1 & 0 \\ 1 & 0 \\ 1 & 1 \\ 1 & 1 \end{pmatrix}\begin{bmatrix} 3 & 0 & 2 & 1 & 3 \\ 1 & 3 & 0 & 2 & 2 \end{bmatrix} & \begin{pmatrix} 0 & 1 \\ 1 & 0 \\ 1 & 0 \\ 1 & 1 \\ 1 & 1 \end{pmatrix}\begin{pmatrix} 1 & 1 \\ 0 & 2 \end{pmatrix} & \begin{pmatrix} 0 & 1 \\ 1 & 0 \\ 1 & 0 \\ 1 & 1 \\ 1 & 1 \end{pmatrix}\begin{pmatrix} 1 & 0 & 1 \\ 1 & 2 & 0 \end{pmatrix} \\ \begin{pmatrix} 0 & 1 \\ 0 & 0 \end{pmatrix}\begin{bmatrix} 3 & 0 & 2 & 1 & 3 \\ 1 & 3 & 0 & 2 & 2 \end{bmatrix} & \begin{pmatrix} 0 & 1 \\ 0 & 0 \end{pmatrix}\begin{pmatrix} 1 & 1 \\ 0 & 2 \end{pmatrix} & \begin{pmatrix} 0 & 1 \\ 0 & 0 \end{pmatrix}\begin{pmatrix} 1 & 0 & 1 \\ 1 & 2 & 0 \end{pmatrix} \end{bmatrix}$$

$$= \begin{bmatrix} 17 & 7 & 3 & 7 \\ 7 & 4 & 5 & 2 \\ 9 & 3 & 2 & 3 \\ 4 & 4 & 4 & 2 \\ 8 & 11 & 7 & 7 \end{bmatrix} \cup \left[ \begin{array}{ccccc|ccc|ccc} 2 & 7 & 10 & 3 & 3 & 11 & 9 & 3 & 7 & 9 \\ 2 & 1 & 1 & 3 & 0 & 2 & 0 & 0 & 1 & 0 \\ 0 & 2 & 3 & 0 & 1 & 3 & 3 & 1 & 2 & 3 \\ \hline 0 & 2 & 3 & 0 & 1 & 3 & 3 & 1 & 2 & 3 \\ 2 & 1 & 1 & 3 & 0 & 2 & 0 & 0 & 1 & 0 \\ 2 & 3 & 4 & 3 & 1 & 5 & 3 & 1 & 3 & 3 \\ 0 & 2 & 3 & 0 & 1 & 3 & 3 & 1 & 2 & 3 \\ 2 & 1 & 1 & 3 & 0 & 2 & 0 & 0 & 1 & 0 \\ \hline 2 & 3 & 4 & 3 & 1 & 5 & 3 & 1 & 3 & 3 \\ 2 & 3 & 4 & 3 & 1 & 5 & 3 & 1 & 3 & 3 \end{array} \right] +$$

$$\left[ \begin{array}{ccccc|cc|ccc} 17 & 15 & 8 & 14 & 22 & 4 & 14 & 9 & 10 & 4 \\ 1 & 3 & 0 & 2 & 2 & 0 & 2 & 1 & 2 & 0 \\ 3 & 0 & 2 & 1 & 3 & 1 & 1 & 1 & 0 & 1 \\ \hline 1 & 3 & 0 & 2 & 2 & 0 & 2 & 1 & 2 & 0 \\ 3 & 0 & 2 & 1 & 3 & 1 & 1 & 1 & 0 & 1 \\ 3 & 0 & 2 & 1 & 3 & 1 & 1 & 1 & 0 & 1 \\ 4 & 3 & 2 & 3 & 5 & 1 & 3 & 2 & 2 & 1 \\ 4 & 3 & 2 & 3 & 5 & 1 & 3 & 2 & 2 & 1 \\ \hline 1 & 3 & 0 & 2 & 2 & 0 & 2 & 1 & 2 & 0 \\ 0 & 0 & 0 & 0 & 0 & 0 & 0 & 0 & 0 & 0 \end{array} \right]$$

$$= \begin{bmatrix} 17 & 7 & 3 & 7 \\ 7 & 4 & 5 & 2 \\ 9 & 3 & 2 & 3 \\ 4 & 4 & 4 & 2 \\ 8 & 11 & 7 & 7 \end{bmatrix} \cup$$

$$\begin{bmatrix} 19 & 22 & 18 & 17 & 25 & 15 & 23 & 12 & 17 & 13 \\ 3 & 4 & 1 & 5 & 2 & 2 & 2 & 1 & 3 & 0 \\ 3 & 2 & 5 & 1 & 4 & 4 & 4 & 2 & 2 & 4 \\ \hline 1 & 5 & 3 & 2 & 3 & 3 & 5 & 2 & 4 & 3 \\ 5 & 1 & 3 & 4 & 3 & 3 & 1 & 1 & 1 & 1 \\ 5 & 3 & 6 & 4 & 4 & 6 & 4 & 2 & 3 & 4 \\ 4 & 5 & 5 & 3 & 6 & 4 & 6 & 3 & 4 & 4 \\ 6 & 4 & 3 & 6 & 5 & 3 & 3 & 2 & 3 & 1 \\ \hline 3 & 6 & 4 & 5 & 3 & 5 & 5 & 2 & 5 & 3 \\ 2 & 3 & 4 & 3 & 1 & 5 & 3 & 1 & 3 & 3 \end{bmatrix}.$$

We see the resultant is also a semi superbimatrix. Thus the minor product of two (compatible under product) semi superbimatrices is a semi superbimatrix.

***Example 2.61:*** Let $A = A_1 \cup A_2$ and $B = B_1 \cup B_2$ be any two semi superbimatrices. We find the product AB. Here

$$A = A_1 \cup A_2$$

$$= \begin{bmatrix} 0 & 1 & 2 & 1 & 2 & 1 & 0 & 1 & 1 \\ 1 & 0 & 1 & 2 & 1 & 0 & 1 & 0 & 1 \\ 2 & 1 & 0 & 0 & 1 & 1 & 0 & 0 & 1 \\ \hline 1 & 5 & 6 & 5 & 6 & 1 & 1 & 2 & 1 \\ 1 & 0 & 0 & 0 & 1 & 0 & 1 & 1 & 0 \\ 2 & 3 & 1 & 1 & 0 & 1 & 0 & 0 & 1 \\ 1 & 0 & 2 & 1 & 1 & 0 & 1 & 0 & 1 \end{bmatrix} \cup \begin{bmatrix} 0 & 2 & 1 & 5 & 6 \\ 1 & 0 & 1 & 1 & 1 \\ 0 & 1 & 0 & 1 & 0 \\ 2 & 1 & 0 & 2 & 1 \\ 4 & 1 & 4 & 0 & 1 \\ 1 & 0 & 1 & 2 & 1 \end{bmatrix}$$

and $B = B_1 \cup B_2$ where

$$B_1 = \begin{bmatrix} 1 & 0 & | & 3 & 0 & 1 \\ 1 & 0 & | & 1 & 1 & 1 \\ 0 & 1 & | & 0 & 1 & 0 \\ \hline 1 & 1 & | & 1 & 2 & 1 \\ 0 & 1 & | & 3 & 1 & 2 \\ \hline 0 & 1 & | & 1 & 0 & 1 \\ 2 & 0 & | & 0 & 1 & 1 \\ 1 & 0 & | & 5 & 0 & 1 \\ 0 & 1 & | & 1 & 1 & 0 \end{bmatrix} \text{ and } B_2 = \begin{bmatrix} 3 & 0 & 1 & 2 \\ 1 & 1 & 2 & 0 \\ 2 & 0 & 1 & 1 \\ 5 & 1 & 0 & 0 \\ 1 & 0 & 1 & 0 \end{bmatrix}.$$

$$\begin{aligned} AB &= (A_1 \cup A_2)(B_1 \cup B_2) \\ &= A_1 B_1 \cup A_2 B_2 \end{aligned}$$

$$= \begin{bmatrix} 0 & 1 & 2 & | & 1 & 2 & | & 1 & 0 & 1 & 1 \\ 1 & 0 & 1 & | & 2 & 1 & | & 0 & 1 & 0 & 1 \\ 2 & 1 & 0 & | & 0 & 1 & | & 1 & 0 & 0 & 1 \\ \hline 1 & 5 & 6 & | & 5 & 6 & | & 1 & 1 & 2 & 1 \\ 1 & 0 & 0 & | & 0 & 1 & | & 0 & 1 & 1 & 0 \\ \hline 2 & 3 & 1 & | & 1 & 0 & | & 1 & 0 & 0 & 1 \\ 1 & 0 & 2 & | & 1 & 1 & | & 0 & 1 & 0 & 1 \end{bmatrix} \begin{bmatrix} 1 & 0 & | & 3 & 0 & 1 \\ 1 & 0 & | & 1 & 1 & 1 \\ 0 & 1 & | & 0 & 1 & 0 \\ \hline 1 & 1 & | & 1 & 2 & 1 \\ 0 & 1 & | & 3 & 1 & 2 \\ \hline 0 & 1 & | & 1 & 0 & 1 \\ 2 & 0 & | & 0 & 1 & 1 \\ 1 & 0 & | & 5 & 0 & 1 \\ 0 & 1 & | & 1 & 1 & 0 \end{bmatrix} \cup$$

$$\begin{bmatrix} 0 & 2 & 1 & 5 & 6 \\ 1 & 0 & 1 & 1 & 1 \\ 0 & 1 & 0 & 1 & 0 \\ 2 & 1 & 0 & 2 & 1 \\ 4 & 1 & 4 & 0 & 1 \\ 1 & 0 & 1 & 2 & 1 \end{bmatrix} \begin{bmatrix} 3 & 0 & 1 & 2 \\ 1 & 1 & 2 & 0 \\ 2 & 0 & 1 & 1 \\ 5 & 1 & 0 & 0 \\ 1 & 0 & 1 & 0 \end{bmatrix} =$$

$$\left[ \frac{\begin{pmatrix} 0 & 1 & 2 \\ 1 & 0 & 1 \\ 2 & 1 & 0 \end{pmatrix}\begin{bmatrix} 1 & 0 \\ 1 & 0 \\ 0 & 1 \end{bmatrix} \Bigg| \begin{pmatrix} 0 & 1 & 2 \\ 1 & 0 & 1 \\ 2 & 1 & 0 \end{pmatrix}\begin{pmatrix} 3 & 0 & 1 \\ 1 & 1 & 1 \\ 0 & 1 & 0 \end{pmatrix}}{\begin{bmatrix} 1 & 5 & 6 \\ 1 & 0 & 0 \\ 2 & 3 & 1 \\ 1 & 0 & 2 \end{bmatrix}\begin{bmatrix} 1 & 0 \\ 1 & 0 \\ 0 & 1 \end{bmatrix} \Bigg| \begin{bmatrix} 1 & 5 & 6 \\ 1 & 0 & 0 \\ 2 & 3 & 1 \\ 1 & 0 & 2 \end{bmatrix}\begin{pmatrix} 3 & 0 & 1 \\ 1 & 1 & 1 \\ 0 & 1 & 0 \end{pmatrix}} \right] +$$

$$\left[ \frac{\begin{pmatrix} 1 & 2 \\ 2 & 1 \\ 0 & 1 \end{pmatrix}\begin{bmatrix} 1 & 1 \\ 0 & 1 \end{bmatrix} \Bigg| \begin{bmatrix} 1 & 2 \\ 2 & 1 \\ 0 & 1 \end{bmatrix}\begin{bmatrix} 1 & 2 & 1 \\ 3 & 1 & 2 \end{bmatrix}}{\begin{bmatrix} 5 & 6 \\ 0 & 1 \\ 1 & 0 \\ 1 & 1 \end{bmatrix}\begin{bmatrix} 1 & 1 \\ 0 & 1 \end{bmatrix} \Bigg| \begin{bmatrix} 5 & 6 \\ 0 & 1 \\ 1 & 0 \\ 1 & 1 \end{bmatrix}\begin{bmatrix} 1 & 2 & 1 \\ 3 & 1 & 2 \end{bmatrix}} \right] +$$

$$\left[ \frac{\begin{pmatrix} 1 & 0 & 1 & 1 \\ 0 & 1 & 0 & 1 \\ 1 & 0 & 0 & 1 \end{pmatrix}\begin{bmatrix} 0 & 1 \\ 2 & 0 \\ 1 & 0 \\ 0 & 1 \end{bmatrix} \Bigg| \begin{pmatrix} 1 & 0 & 1 & 1 \\ 0 & 1 & 0 & 1 \\ 1 & 0 & 0 & 1 \end{pmatrix}\begin{pmatrix} 1 & 0 & 1 \\ 0 & 1 & 1 \\ 5 & 0 & 1 \\ 1 & 1 & 0 \end{pmatrix}}{\begin{bmatrix} 1 & 1 & 2 & 1 \\ 0 & 1 & 1 & 0 \\ 1 & 0 & 0 & 1 \\ 0 & 1 & 0 & 1 \end{bmatrix}\begin{bmatrix} 0 & 1 \\ 2 & 0 \\ 1 & 0 \\ 0 & 1 \end{bmatrix} \Bigg| \begin{bmatrix} 1 & 1 & 2 & 1 \\ 0 & 1 & 1 & 0 \\ 1 & 0 & 0 & 1 \\ 0 & 1 & 0 & 1 \end{bmatrix}\begin{pmatrix} 1 & 0 & 1 \\ 0 & 1 & 1 \\ 5 & 0 & 1 \\ 1 & 1 & 0 \end{pmatrix}} \right]$$

$\cup\ A_2 B_2.$

$$= \begin{bmatrix} 1 & 2 & 1 & 3 & 1 \\ 1 & 1 & 3 & 1 & 1 \\ 3 & 0 & 7 & 1 & 3 \\ \hline 6 & 6 & 8 & 11 & 6 \\ 1 & 0 & 3 & 0 & 1 \\ 5 & 1 & 9 & 4 & 5 \\ 1 & 2 & 3 & 2 & 1 \end{bmatrix} + \begin{bmatrix} 1 & 3 & 6 & 4 & 5 \\ 2 & 3 & 5 & 5 & 4 \\ 0 & 1 & 3 & 1 & 2 \\ \hline 5 & 11 & 23 & 16 & 17 \\ 0 & 1 & 3 & 1 & 2 \\ 1 & 1 & 1 & 2 & 1 \\ 1 & 2 & 4 & 3 & 3 \end{bmatrix} +$$

$$\begin{bmatrix} 1 & 2 & 7 & 1 & 2 \\ 2 & 1 & 1 & 2 & 1 \\ 0 & 2 & 2 & 1 & 1 \\ \hline 4 & 2 & 12 & 2 & 4 \\ 3 & 0 & 5 & 1 & 2 \\ 0 & 2 & 2 & 1 & 1 \\ 2 & 1 & 1 & 2 & 1 \end{bmatrix} \cup \begin{bmatrix} 35 & 7 & 11 & 1 \\ 11 & 1 & 3 & 3 \\ 6 & 2 & 2 & 0 \\ 17 & 3 & 5 & 4 \\ 22 & 1 & 11 & 12 \\ 16 & 2 & 3 & 3 \end{bmatrix}$$

$$= \begin{bmatrix} 3 & 7 & 14 & 8 & 8 \\ 5 & 5 & 9 & 8 & 6 \\ 3 & 3 & 12 & 3 & 6 \\ \hline 15 & 19 & 43 & 29 & 27 \\ 4 & 1 & 11 & 2 & 5 \\ 6 & 4 & 12 & 7 & 7 \\ 4 & 5 & 8 & 7 & 5 \end{bmatrix} \cup \begin{bmatrix} 35 & 7 & 11 & 1 \\ 11 & 1 & 3 & 3 \\ 6 & 2 & 2 & 0 \\ 17 & 3 & 5 & 4 \\ 22 & 1 & 11 & 12 \\ 16 & 2 & 3 & 3 \end{bmatrix}.$$

Clearly AB is again a semi superbimatrix. Thus the minor product of semi superbimatrices yields a semi superbimatrix.

Now we proceed on to find the minor product of a semi superbimatrix with its transpose.

*Example 2.62:* Let $A = A_1 \cup A_2$ be a semi superbimatrix where

$$A_1 = \begin{bmatrix} 0 & 1 & 2 & 3 \\ 1 & 0 & 1 & 1 \\ 0 & 1 & 0 & 1 \\ 1 & 2 & 0 & 0 \\ 5 & 1 & 2 & 1 \\ 0 & 1 & 1 & 0 \end{bmatrix}$$

and

$$A_2 = \left[ \begin{array}{cc|cccc} 1 & 0 & 1 & 5 & 1 & 1 \\ 2 & 1 & 0 & 6 & 1 & 0 \\ 0 & 1 & 2 & 1 & 0 & 1 \\ \hline 1 & 2 & 1 & 0 & 3 & 0 \\ 1 & 1 & 0 & 2 & 2 & 2 \\ 5 & 0 & 1 & 1 & 0 & 1 \\ 6 & 1 & 2 & 0 & 1 & 0 \end{array} \right].$$

$A^T = (A_1 \cup A_2)^T = (A_1^T \cup A_2^T)$

$$A_1^T = \begin{bmatrix} 0 & 1 & 0 & 1 & 5 & 0 \\ 1 & 0 & 1 & 2 & 1 & 1 \\ 2 & 1 & 0 & 0 & 2 & 1 \\ 3 & 1 & 1 & 0 & 1 & 0 \end{bmatrix}$$

and

$$A_2^T = \left[ \begin{array}{ccc|ccc|c} 1 & 2 & 0 & 1 & 1 & 5 & 6 \\ 0 & 1 & 1 & 2 & 1 & 0 & 1 \\ \hline 1 & 0 & 2 & 1 & 0 & 1 & 2 \\ 5 & 6 & 1 & 0 & 2 & 1 & 0 \\ 1 & 1 & 0 & 3 & 2 & 0 & 1 \\ 1 & 0 & 1 & 0 & 2 & 1 & 0 \end{array} \right].$$

$$\begin{aligned} (A_1 \cup A_2)(A_1 \cup A_2)^T &= (A_1 \cup A_2)(A_1^T \cup A_2^T) \\ &= A_1 A_1^T \cup A_2 A_2^T \end{aligned}$$

$$= \begin{bmatrix} 0 & 1 & 2 & 3 \\ 1 & 0 & 1 & 1 \\ 0 & 1 & 0 & 1 \\ 1 & 2 & 0 & 0 \\ 5 & 1 & 2 & 1 \\ 0 & 1 & 1 & 0 \end{bmatrix} \begin{bmatrix} 0 & 1 & 0 & 1 & 5 & 0 \\ 1 & 0 & 1 & 2 & 1 & 1 \\ 2 & 1 & 0 & 0 & 2 & 1 \\ 3 & 1 & 1 & 0 & 1 & 0 \end{bmatrix} \cup$$

$$\left[ \begin{array}{cc|ccc|c} 1 & 0 & 1 & 5 & 1 & 1 \\ 2 & 1 & 0 & 6 & 1 & 0 \\ 0 & 1 & 2 & 1 & 0 & 1 \\ \hline 1 & 2 & 1 & 0 & 3 & 0 \\ 1 & 1 & 0 & 2 & 2 & 2 \\ 5 & 0 & 1 & 1 & 0 & 1 \\ \hline 6 & 1 & 2 & 0 & 1 & 0 \end{array} \right] \left[ \begin{array}{ccc|ccc|c} 1 & 2 & 0 & 1 & 1 & 5 & 6 \\ 0 & 1 & 1 & 2 & 1 & 0 & 1 \\ \hline 1 & 0 & 2 & 1 & 0 & 1 & 2 \\ 5 & 6 & 1 & 0 & 2 & 1 & 0 \\ 1 & 1 & 0 & 3 & 2 & 0 & 1 \\ \hline 1 & 0 & 1 & 0 & 2 & 1 & 0 \end{array} \right]$$

$$= \begin{bmatrix} 14 & 5 & 4 & 2 & 8 & 3 \\ 5 & 3 & 1 & 1 & 8 & 1 \\ 4 & 1 & 2 & 2 & 2 & 1 \\ 2 & 1 & 2 & 5 & 7 & 2 \\ 8 & 8 & 2 & 7 & 31 & 3 \\ 3 & 1 & 1 & 2 & 3 & 2 \end{bmatrix} \cup$$

$$\left[ \begin{array}{c|c|c} \begin{pmatrix} 1 & 0 \\ 2 & 1 \\ 0 & 1 \end{pmatrix}\begin{pmatrix} 1 & 2 & 0 \\ 0 & 1 & 1 \end{pmatrix} & \begin{pmatrix} 1 & 0 \\ 2 & 1 \\ 0 & 1 \end{pmatrix}\begin{pmatrix} 1 & 1 & 5 \\ 2 & 1 & 0 \end{pmatrix} & \begin{pmatrix} 1 & 0 \\ 2 & 1 \\ 0 & 1 \end{pmatrix}\begin{pmatrix} 6 \\ 1 \end{pmatrix} \\ \hline \begin{pmatrix} 1 & 2 \\ 1 & 1 \\ 5 & 0 \end{pmatrix}\begin{pmatrix} 1 & 2 & 0 \\ 0 & 1 & 1 \end{pmatrix} & \begin{pmatrix} 1 & 2 \\ 1 & 1 \\ 5 & 0 \end{pmatrix}\begin{pmatrix} 1 & 1 & 5 \\ 2 & 1 & 0 \end{pmatrix} & \begin{pmatrix} 1 & 2 \\ 1 & 1 \\ 5 & 0 \end{pmatrix}\begin{pmatrix} 6 \\ 1 \end{pmatrix} \\ \hline \begin{pmatrix} 6 & 1 \end{pmatrix}\begin{pmatrix} 1 & 2 & 0 \\ 0 & 1 & 1 \end{pmatrix} & \begin{pmatrix} 6 & 1 \end{pmatrix}\begin{pmatrix} 1 & 1 & 5 \\ 2 & 1 & 0 \end{pmatrix} & \begin{pmatrix} 6 & 1 \end{pmatrix}\begin{pmatrix} 6 \\ 1 \end{pmatrix} \end{array} \right] +$$

$$\begin{bmatrix} \begin{pmatrix} 1 & 5 & 1 \\ 0 & 6 & 1 \\ 2 & 1 & 0 \end{pmatrix}\begin{pmatrix} 1 & 0 & 2 \\ 5 & 6 & 1 \\ 1 & 1 & 0 \end{pmatrix} & \begin{pmatrix} 1 & 5 & 1 \\ 0 & 6 & 1 \\ 2 & 1 & 0 \end{pmatrix}\begin{pmatrix} 1 & 0 & 1 \\ 0 & 2 & 1 \\ 3 & 2 & 0 \end{pmatrix} & \begin{pmatrix} 1 & 5 & 1 \\ 0 & 6 & 1 \\ 2 & 1 & 0 \end{pmatrix}\begin{pmatrix} 2 \\ 0 \\ 1 \end{pmatrix} \\ \hline \begin{pmatrix} 1 & 0 & 3 \\ 0 & 2 & 2 \\ 1 & 1 & 0 \end{pmatrix}\begin{pmatrix} 1 & 0 & 2 \\ 5 & 6 & 1 \\ 1 & 1 & 0 \end{pmatrix} & \begin{pmatrix} 1 & 0 & 3 \\ 0 & 2 & 2 \\ 1 & 1 & 0 \end{pmatrix}\begin{pmatrix} 1 & 0 & 1 \\ 0 & 2 & 1 \\ 3 & 2 & 0 \end{pmatrix} & \begin{pmatrix} 1 & 0 & 3 \\ 0 & 2 & 2 \\ 1 & 1 & 0 \end{pmatrix}\begin{pmatrix} 2 \\ 0 \\ 1 \end{pmatrix} \\ \hline \begin{pmatrix} 2 & 0 & 1 \end{pmatrix}\begin{pmatrix} 1 & 0 & 2 \\ 5 & 6 & 1 \\ 1 & 1 & 0 \end{pmatrix} & \begin{pmatrix} 2 & 0 & 1 \end{pmatrix}\begin{pmatrix} 1 & 0 & 1 \\ 0 & 2 & 1 \\ 3 & 2 & 0 \end{pmatrix} & \begin{pmatrix} 2 & 0 & 1 \end{pmatrix}\begin{pmatrix} 2 \\ 0 \\ 1 \end{pmatrix} \end{bmatrix}$$

$$+ \begin{bmatrix} \begin{pmatrix} 1 \\ 0 \\ 1 \end{pmatrix}\begin{pmatrix} 1 & 0 & 1 \end{pmatrix} & \begin{pmatrix} 1 \\ 0 \\ 1 \end{pmatrix}\begin{pmatrix} 0 & 2 & 1 \end{pmatrix} & \begin{pmatrix} 1 \\ 0 \\ 1 \end{pmatrix}(0) \\ \hline \begin{pmatrix} 0 \\ 2 \\ 1 \end{pmatrix}\begin{pmatrix} 1 & 0 & 1 \end{pmatrix} & \begin{pmatrix} 0 \\ 2 \\ 1 \end{pmatrix}\begin{pmatrix} 0 & 2 & 1 \end{pmatrix} & \begin{pmatrix} 0 \\ 2 \\ 1 \end{pmatrix}(0) \\ \hline (0)\begin{pmatrix} 1 & 0 & 1 \end{pmatrix} & (0)\begin{pmatrix} 0 & 2 & 1 \end{pmatrix} & (0)(0) \end{bmatrix}$$

$$= \begin{bmatrix} 14 & 5 & 4 & 2 & 8 & 3 \\ 5 & 3 & 1 & 1 & 8 & 1 \\ 4 & 1 & 2 & 2 & 2 & 1 \\ 2 & 1 & 2 & 5 & 7 & 2 \\ 8 & 8 & 2 & 7 & 31 & 3 \\ 3 & 1 & 1 & 2 & 3 & 2 \end{bmatrix} \cup \begin{bmatrix} 1 & 2 & 0 & 1 & 1 & 5 & 6 \\ 2 & 5 & 1 & 4 & 3 & 10 & 13 \\ 0 & 1 & 1 & 2 & 1 & 0 & 1 \\ \hline 1 & 4 & 2 & 5 & 3 & 5 & 8 \\ 1 & 3 & 1 & 3 & 2 & 5 & 7 \\ 5 & 10 & 0 & 5 & 5 & 25 & 30 \\ \hline 6 & 13 & 1 & 8 & 7 & 30 & 37 \end{bmatrix} +$$

$$\begin{bmatrix} 27 & 31 & 7 & 4 & 12 & 6 & 3 \\ 31 & 37 & 6 & 3 & 14 & 6 & 1 \\ 7 & 6 & 5 & 2 & 2 & 3 & 4 \\ \hline 4 & 3 & 2 & 10 & 6 & 1 & 5 \\ 12 & 14 & 2 & 6 & 8 & 2 & 2 \\ 6 & 6 & 3 & 1 & 2 & 2 & 2 \\ 3 & 1 & 4 & 5 & 2 & 2 & 5 \end{bmatrix} + \begin{bmatrix} 1 & 0 & 1 & 0 & 2 & 1 & 0 \\ 0 & 0 & 0 & 0 & 0 & 0 & 0 \\ 1 & 0 & 1 & 0 & 2 & 1 & 0 \\ \hline 0 & 0 & 0 & 0 & 0 & 0 & 0 \\ 2 & 0 & 2 & 0 & 4 & 2 & 0 \\ 1 & 0 & 1 & 0 & 2 & 1 & 0 \\ 0 & 0 & 0 & 0 & 0 & 0 & 0 \end{bmatrix}$$

$$= \begin{bmatrix} 14 & 5 & 4 & 2 & 8 & 3 \\ 5 & 3 & 1 & 1 & 8 & 1 \\ 4 & 1 & 2 & 2 & 2 & 1 \\ 2 & 1 & 2 & 5 & 7 & 2 \\ 8 & 8 & 2 & 7 & 31 & 3 \\ 3 & 1 & 1 & 2 & 3 & 2 \end{bmatrix} \cup \begin{bmatrix} 29 & 33 & 8 & 5 & 15 & 12 & 9 \\ 33 & 42 & 7 & 7 & 17 & 16 & 14 \\ 8 & 7 & 7 & 4 & 5 & 4 & 5 \\ \hline 5 & 7 & 4 & 15 & 9 & 6 & 13 \\ 15 & 17 & 5 & 9 & 14 & 9 & 9 \\ 12 & 16 & 4 & 6 & 9 & 28 & 32 \\ 9 & 14 & 5 & 13 & 9 & 32 & 42 \end{bmatrix}.$$

We see $AA^T$ is a symmetric semi superbimatrix. We give yet another example.

***Example 2.63:*** Let $A = A_1 \cup A_2$ be a semi superbimatrix where

$$A_2 = \begin{bmatrix} 0 & 1 & 2 & 3 & 4 & 1 \\ 2 & 3 & 4 & 1 & 0 & 0 \\ 3 & 4 & 1 & 0 & 1 & 0 \\ 4 & 1 & 0 & 1 & 0 & 3 \\ 1 & 0 & 1 & 0 & 3 & 4 \\ \hline 0 & 1 & 0 & 1 & 0 & 1 \\ 1 & 0 & 1 & 0 & 1 & 0 \\ 1 & 1 & 1 & 0 & 0 & 1 \end{bmatrix}$$

and

$$A_1 = \begin{bmatrix} 3 & 0 & 2 & 4 \\ 1 & 1 & 0 & 1 \\ 2 & 2 & 1 & 0 \\ 6 & 0 & 0 & 2 \\ 1 & 1 & 0 & 1 \end{bmatrix}.$$

$A^T = (A_1 \cup A_2)^T = (A_1^T \cup A_2^T)$

$$= \begin{bmatrix} 3 & 1 & 2 & 6 & 1 \\ 0 & 1 & 2 & 0 & 1 \\ 2 & 0 & 1 & 0 & 0 \\ 4 & 1 & 0 & 2 & 1 \end{bmatrix} \cup \left[\begin{array}{ccccc|ccc} 0 & 2 & 3 & 4 & 1 & 0 & 1 & 1 \\ 1 & 3 & 4 & 1 & 0 & 1 & 0 & 1 \\ 2 & 4 & 1 & 0 & 1 & 0 & 1 & 1 \\ 3 & 1 & 0 & 1 & 0 & 1 & 0 & 0 \\ 4 & 0 & 1 & 0 & 3 & 0 & 1 & 0 \\ 1 & 0 & 0 & 3 & 4 & 1 & 0 & 1 \end{array}\right].$$

$A^T = (A_1 \cup A_2)^T (A_1 \cup A_2)$
$= (A_1^T \cup A_2^T)(A_1 \cup A_2)$
$= A_1^T A_1 \cup A_2^T A_2$

$$= \begin{bmatrix} 3 & 1 & 2 & 6 & 1 \\ 0 & 1 & 2 & 0 & 1 \\ 2 & 0 & 1 & 0 & 0 \\ 4 & 1 & 0 & 2 & 1 \end{bmatrix} \begin{bmatrix} 3 & 0 & 2 & 4 \\ 1 & 1 & 0 & 1 \\ 2 & 2 & 1 & 0 \\ 6 & 0 & 0 & 2 \\ 1 & 1 & 0 & 1 \end{bmatrix} \cup$$

$$\left[\begin{array}{ccccc|ccc} 0 & 2 & 3 & 4 & 1 & 0 & 1 & 1 \\ 1 & 3 & 4 & 1 & 0 & 1 & 0 & 1 \\ \hline 2 & 4 & 1 & 0 & 1 & 0 & 1 & 1 \\ 3 & 1 & 0 & 1 & 0 & 1 & 0 & 0 \\ \hline 4 & 0 & 1 & 0 & 3 & 0 & 1 & 0 \\ 1 & 0 & 0 & 3 & 4 & 1 & 0 & 1 \end{array}\right] \left[\begin{array}{cc|cc|cc} 0 & 1 & 2 & 3 & 4 & 1 \\ 2 & 3 & 4 & 1 & 0 & 0 \\ 3 & 4 & 1 & 0 & 1 & 0 \\ 4 & 1 & 0 & 1 & 0 & 3 \\ 1 & 0 & 1 & 0 & 3 & 4 \\ \hline 0 & 1 & 0 & 1 & 0 & 1 \\ 1 & 0 & 1 & 0 & 1 & 0 \\ 1 & 1 & 1 & 0 & 0 & 1 \end{array}\right]$$

$$= \begin{bmatrix} 51 & 6 & 8 & 26 \\ 6 & 6 & 2 & 2 \\ 8 & 2 & 5 & 8 \\ 26 & 2 & 8 & 22 \end{bmatrix} \cup$$

$$\begin{bmatrix} \begin{pmatrix} 0 & 2 & 3 & 4 & 1 \\ 1 & 3 & 4 & 1 & 0 \end{pmatrix} \begin{bmatrix} 0 & 1 \\ 2 & 3 \\ 3 & 4 \\ 4 & 1 \\ 1 & 0 \end{bmatrix} & \begin{pmatrix} 0 & 2 & 3 & 4 & 1 \\ 1 & 3 & 4 & 1 & 0 \end{pmatrix} \begin{bmatrix} 2 & 3 \\ 4 & 1 \\ 1 & 0 \\ 0 & 1 \\ 1 & 0 \end{bmatrix} & \begin{pmatrix} 0 & 2 & 3 & 4 & 1 \\ 1 & 3 & 4 & 1 & 0 \end{pmatrix} \begin{bmatrix} 4 & 1 \\ 0 & 0 \\ 1 & 0 \\ 0 & 3 \\ 3 & 4 \end{bmatrix} \\ \begin{pmatrix} 2 & 4 & 1 & 0 & 1 \\ 3 & 1 & 0 & 1 & 0 \end{pmatrix} \begin{bmatrix} 0 & 1 \\ 2 & 3 \\ 3 & 4 \\ 4 & 1 \\ 1 & 0 \end{bmatrix} & \begin{pmatrix} 2 & 4 & 1 & 0 & 1 \\ 3 & 1 & 0 & 1 & 0 \end{pmatrix} \begin{bmatrix} 2 & 3 \\ 4 & 1 \\ 1 & 0 \\ 0 & 1 \\ 1 & 0 \end{bmatrix} & \begin{pmatrix} 2 & 4 & 1 & 0 & 1 \\ 3 & 1 & 0 & 1 & 0 \end{pmatrix} \begin{bmatrix} 4 & 1 \\ 0 & 0 \\ 1 & 0 \\ 0 & 3 \\ 3 & 4 \end{bmatrix} \\ \begin{pmatrix} 4 & 0 & 1 & 0 & 3 \\ 1 & 0 & 0 & 3 & 4 \end{pmatrix} \begin{bmatrix} 0 & 1 \\ 2 & 3 \\ 3 & 4 \\ 4 & 1 \\ 1 & 0 \end{bmatrix} & \begin{pmatrix} 4 & 0 & 1 & 0 & 3 \\ 1 & 0 & 0 & 3 & 4 \end{pmatrix} \begin{bmatrix} 2 & 3 \\ 4 & 1 \\ 1 & 0 \\ 0 & 1 \\ 1 & 0 \end{bmatrix} & \begin{pmatrix} 4 & 0 & 1 & 0 & 3 \\ 1 & 0 & 0 & 3 & 4 \end{pmatrix} \begin{bmatrix} 4 & 1 \\ 0 & 0 \\ 1 & 0 \\ 0 & 3 \\ 3 & 4 \end{bmatrix} \end{bmatrix}$$

$$+ \begin{bmatrix} \begin{pmatrix} 0 & 1 & 1 \\ 1 & 0 & 1 \end{pmatrix} \begin{bmatrix} 0 & 1 \\ 1 & 0 \\ 1 & 1 \end{bmatrix} & \begin{pmatrix} 0 & 1 & 1 \\ 1 & 0 & 1 \end{pmatrix} \begin{bmatrix} 0 & 1 \\ 1 & 0 \\ 1 & 0 \end{bmatrix} & \begin{pmatrix} 0 & 1 & 1 \\ 1 & 0 & 1 \end{pmatrix} \begin{bmatrix} 0 & 1 \\ 1 & 0 \\ 0 & 1 \end{bmatrix} \\ \begin{pmatrix} 0 & 1 & 1 \\ 1 & 0 & 0 \end{pmatrix} \begin{bmatrix} 0 & 1 \\ 1 & 0 \\ 1 & 1 \end{bmatrix} & \begin{pmatrix} 0 & 1 & 1 \\ 1 & 0 & 0 \end{pmatrix} \begin{bmatrix} 0 & 1 \\ 1 & 0 \\ 1 & 0 \end{bmatrix} & \begin{pmatrix} 0 & 1 & 1 \\ 1 & 0 & 0 \end{pmatrix} \begin{bmatrix} 0 & 1 \\ 1 & 0 \\ 0 & 1 \end{bmatrix} \\ \begin{pmatrix} 0 & 1 & 0 \\ 1 & 0 & 1 \end{pmatrix} \begin{bmatrix} 0 & 1 \\ 1 & 0 \\ 1 & 1 \end{bmatrix} & \begin{pmatrix} 0 & 1 & 0 \\ 1 & 0 & 1 \end{pmatrix} \begin{bmatrix} 0 & 1 \\ 1 & 0 \\ 1 & 0 \end{bmatrix} & \begin{pmatrix} 0 & 1 & 0 \\ 1 & 0 & 1 \end{pmatrix} \begin{bmatrix} 0 & 1 \\ 1 & 0 \\ 0 & 1 \end{bmatrix} \end{bmatrix}$$

$$= \begin{bmatrix} 51 & 6 & 8 & 26 \\ 6 & 6 & 2 & 2 \\ 8 & 2 & 5 & 8 \\ 26 & 2 & 8 & 22 \end{bmatrix} \cup$$

$$\begin{bmatrix} 30 & 22 & | & 12 & 6 & | & 6 & 16 \\ 22 & 27 & | & 18 & 7 & | & 8 & 4 \\ \hline 12 & 18 & | & 22 & 10 & | & 12 & 6 \\ 6 & 7 & | & 10 & 11 & | & 12 & 6 \\ \hline 6 & 8 & | & 12 & 12 & | & 26 & 16 \\ 16 & 4 & | & 6 & 6 & | & 16 & 26 \end{bmatrix} +$$

$$\begin{bmatrix} 2 & 1 & | & 2 & 0 & | & 1 & 1 \\ 1 & 2 & | & 1 & 1 & | & 0 & 2 \\ \hline 2 & 1 & | & 2 & 0 & | & 1 & 1 \\ 0 & 1 & | & 0 & 1 & | & 0 & 1 \\ \hline 1 & 0 & | & 1 & 0 & | & 1 & 0 \\ 1 & 2 & | & 1 & 1 & | & 0 & 2 \end{bmatrix} =$$

$$\begin{bmatrix} 51 & 6 & 8 & 26 \\ 6 & 6 & 2 & 2 \\ 8 & 2 & 5 & 8 \\ 26 & 2 & 8 & 22 \end{bmatrix} \cup \begin{bmatrix} 32 & 23 & | & 14 & 6 & | & 7 & 17 \\ 23 & 29 & | & 19 & 8 & | & 8 & 6 \\ \hline 14 & 19 & | & 24 & 10 & | & 13 & 7 \\ 6 & 8 & | & 10 & 12 & | & 12 & 7 \\ \hline 7 & 8 & | & 13 & 12 & | & 27 & 16 \\ 11 & 6 & | & 7 & 7 & | & 16 & 28 \end{bmatrix}.$$

Thus we see the resultant of the minor product of $A^T A$ is a symmetric semi superbimatrix. Now we find

$$\begin{aligned} AA^T &= (A_1 \cup A_2)(A_1 \cup A_2)^T \\ &= (A_1 \cup A_2)(A_1^T \cup A_2^T) \\ &= A_1 A_1^T \cup A_2 A_2^T \end{aligned}$$

$$= \begin{bmatrix} 3 & 0 & 2 & 4 \\ 1 & 1 & 0 & 1 \\ 2 & 2 & 1 & 0 \\ 6 & 0 & 0 & 2 \\ 1 & 1 & 0 & 1 \end{bmatrix} \begin{bmatrix} 3 & 1 & 2 & 6 & 1 \\ 0 & 1 & 2 & 0 & 1 \\ 2 & 0 & 1 & 0 & 0 \\ 4 & 1 & 0 & 2 & 1 \end{bmatrix}$$

$$\cup \left[ \begin{array}{cc|cc|cc} 0 & 1 & 2 & 3 & 4 & 1 \\ 2 & 3 & 4 & 1 & 0 & 0 \\ 3 & 4 & 1 & 0 & 1 & 0 \\ 4 & 1 & 0 & 1 & 0 & 3 \\ 1 & 0 & 1 & 0 & 3 & 4 \\ \hline 0 & 1 & 0 & 1 & 0 & 1 \\ 1 & 0 & 1 & 0 & 1 & 0 \\ 1 & 1 & 1 & 0 & 0 & 1 \end{array} \right] \left[ \begin{array}{ccccc|ccc} 0 & 2 & 3 & 4 & 1 & 0 & 1 & 1 \\ 1 & 3 & 4 & 1 & 0 & 1 & 0 & 1 \\ \hline 2 & 4 & 1 & 0 & 1 & 0 & 1 & 1 \\ 3 & 1 & 0 & 1 & 0 & 1 & 0 & 0 \\ \hline 4 & 0 & 1 & 0 & 3 & 0 & 1 & 0 \\ 1 & 0 & 0 & 3 & 4 & 1 & 0 & 1 \end{array} \right]$$

$$= \begin{bmatrix} 29 & 7 & 8 & 26 & 7 \\ 7 & 3 & 4 & 8 & 3 \\ 8 & 4 & 9 & 12 & 4 \\ 26 & 8 & 12 & 40 & 8 \\ 7 & 3 & 4 & 8 & 3 \end{bmatrix} \cup$$

$$\left[ \begin{array}{c|c} \begin{bmatrix} 0 & 1 \\ 2 & 3 \\ 3 & 4 \\ 4 & 1 \\ 1 & 0 \end{bmatrix} \begin{pmatrix} 0 & 2 & 3 & 4 & 1 \\ 1 & 3 & 4 & 1 & 0 \end{pmatrix} & \begin{bmatrix} 0 & 1 \\ 2 & 3 \\ 3 & 4 \\ 4 & 1 \\ 1 & 0 \end{bmatrix} \begin{bmatrix} 0 & 1 & 1 \\ 1 & 0 & 1 \end{bmatrix} \\ \hline \begin{bmatrix} 0 & 1 \\ 1 & 0 \\ 1 & 1 \end{bmatrix} \begin{pmatrix} 0 & 2 & 3 & 4 & 1 \\ 1 & 3 & 4 & 1 & 0 \end{pmatrix} & \begin{bmatrix} 0 & 1 \\ 1 & 0 \\ 1 & 1 \end{bmatrix} \begin{bmatrix} 0 & 1 & 1 \\ 1 & 0 & 1 \end{bmatrix} \end{array} \right] +$$

$$\left[\begin{array}{c|c} \begin{bmatrix} 2 & 3 \\ 4 & 1 \\ 1 & 0 \\ 0 & 1 \\ 1 & 0 \end{bmatrix} \begin{pmatrix} 2 & 4 & 1 & 0 & 1 \\ 3 & 1 & 0 & 1 & 0 \end{pmatrix} & \begin{bmatrix} 2 & 3 \\ 4 & 1 \\ 1 & 0 \\ 0 & 1 \\ 1 & 0 \end{bmatrix} \begin{bmatrix} 0 & 1 & 1 \\ 1 & 0 & 0 \end{bmatrix} \\ \hline \begin{bmatrix} 0 & 1 \\ 1 & 0 \\ 1 & 0 \end{bmatrix} \begin{pmatrix} 2 & 4 & 1 & 0 & 1 \\ 3 & 1 & 0 & 1 & 0 \end{pmatrix} & \begin{bmatrix} 0 & 1 \\ 1 & 0 \\ 1 & 0 \end{bmatrix} \begin{bmatrix} 0 & 1 & 1 \\ 1 & 0 & 0 \end{bmatrix} \end{array}\right] +$$

$$\left[\begin{array}{c|c} \begin{bmatrix} 4 & 1 \\ 0 & 0 \\ 1 & 0 \\ 0 & 3 \\ 3 & 4 \end{bmatrix} \begin{pmatrix} 4 & 0 & 1 & 0 & 3 \\ 1 & 0 & 0 & 3 & 4 \end{pmatrix} & \begin{bmatrix} 4 & 1 \\ 0 & 0 \\ 1 & 0 \\ 0 & 3 \\ 3 & 4 \end{bmatrix} \begin{bmatrix} 0 & 1 & 0 \\ 1 & 0 & 1 \end{bmatrix} \\ \hline \begin{bmatrix} 0 & 1 \\ 1 & 0 \\ 0 & 1 \end{bmatrix} \begin{pmatrix} 4 & 0 & 1 & 0 & 3 \\ 1 & 0 & 0 & 3 & 4 \end{pmatrix} & \begin{bmatrix} 0 & 1 \\ 1 & 0 \\ 0 & 1 \end{bmatrix} \begin{bmatrix} 0 & 1 & 0 \\ 1 & 0 & 1 \end{bmatrix} \end{array}\right]$$

$$= \begin{bmatrix} 29 & 7 & 8 & 26 & 7 \\ 7 & 3 & 4 & 8 & 3 \\ 8 & 4 & 9 & 12 & 4 \\ 26 & 8 & 12 & 40 & 8 \\ 7 & 3 & 4 & 8 & 3 \end{bmatrix} \cup \left[\begin{array}{ccccc|ccc} 1 & 3 & 4 & 1 & 0 & 1 & 0 & 1 \\ 3 & 13 & 18 & 11 & 2 & 3 & 2 & 5 \\ 4 & 18 & 25 & 16 & 3 & 4 & 3 & 7 \\ 1 & 11 & 16 & 17 & 4 & 1 & 4 & 5 \\ 0 & 2 & 3 & 4 & 1 & 0 & 1 & 1 \\ \hline 1 & 3 & 4 & 1 & 0 & 1 & 0 & 1 \\ 0 & 2 & 3 & 4 & 1 & 0 & 1 & 1 \\ 1 & 5 & 7 & 5 & 1 & 1 & 1 & 2 \end{array}\right] +$$

$$\begin{bmatrix} 13 & 11 & 2 & 3 & 2 & 3 & 2 & 2 \\ 11 & 17 & 4 & 1 & 4 & 1 & 4 & 4 \\ 2 & 4 & 1 & 0 & 1 & 0 & 1 & 1 \\ 3 & 1 & 0 & 1 & 0 & 1 & 0 & 0 \\ 2 & 4 & 1 & 0 & 1 & 0 & 1 & 1 \\ \hline 3 & 1 & 0 & 1 & 0 & 1 & 0 & 0 \\ 2 & 4 & 1 & 0 & 1 & 0 & 1 & 1 \\ 2 & 4 & 1 & 0 & 1 & 0 & 1 & 1 \end{bmatrix} +$$

$$\begin{bmatrix} 17 & 0 & 4 & 3 & 16 & 0 & 4 & 1 \\ 0 & 0 & 0 & 0 & 0 & 0 & 0 & 0 \\ 4 & 0 & 1 & 0 & 3 & 0 & 1 & 0 \\ 3 & 0 & 0 & 9 & 12 & 3 & 0 & 3 \\ 16 & 0 & 3 & 12 & 25 & 4 & 3 & 4 \\ \hline 0 & 0 & 0 & 3 & 4 & 1 & 0 & 1 \\ 4 & 0 & 1 & 0 & 3 & 0 & 1 & 0 \\ 1 & 0 & 0 & 3 & 4 & 1 & 0 & 1 \end{bmatrix}$$

$$= \begin{bmatrix} 29 & 7 & 8 & 26 & 7 \\ 7 & 3 & 4 & 8 & 3 \\ 8 & 4 & 9 & 12 & 4 \\ 26 & 8 & 12 & 40 & 8 \\ 7 & 3 & 4 & 8 & 3 \end{bmatrix} \cup \begin{bmatrix} 31 & 14 & 10 & 7 & 18 & 4 & 6 & 4 \\ 14 & 30 & 22 & 12 & 6 & 4 & 6 & 9 \\ 10 & 22 & 27 & 16 & 7 & 4 & 5 & 8 \\ 7 & 12 & 16 & 27 & 16 & 5 & 4 & 8 \\ 18 & 6 & 7 & 16 & 27 & 4 & 5 & 6 \\ \hline 4 & 4 & 4 & 5 & 4 & 3 & 0 & 2 \\ 6 & 6 & 5 & 4 & 5 & 0 & 3 & 2 \\ 4 & 9 & 8 & 8 & 6 & 2 & 2 & 4 \end{bmatrix}$$

is a symmetric semi superbimatrices. Clearly from this example 2.63 we see $AA^T \ne A^TA$, but both of them are symmetric super semi bimatrices.

**THEOREM 2.3:** *Let $A = A_1 \cup A_2$ be a semi superbimatrix. $A^T$ be the transpose of A. Then $AA^T$ and $A^TA$ are in general two distinct symmetric superbimatrices.*

The proof is left as an exercise to the reader; however we give a small hint for the interested reader to work out for the proof of the theorem.

*Hint:* Let $A = A_1 \cup A_2$ be a semi superbimatrix where $A_1$ is a m × n simple matrix and $A_2$ is a s × t super matrix.
Now
$$A^T = (A_1 \cup A_2)^T$$
$$= (A_1^T \cup A_2^T)$$

is such as $A_1^T$ is a n × m matrix and $A_2^T$ is a t × s supermatrix $AA^T = B_1 \cup B_2$ is a symmetric semi superbimatrix such that $B_1$ is a m × m symmetric simple matrix and $B_2$ is a s × s symmetric supermatrix. Thus $AA^T = B_1 \cup B_2$ is a symmetric semi superbimatrix.

On the other hand we see $A^TA = C_1 \cup C_2$ where $C_1$ is a n × n symmetric simple matrix and $C_2$ is a t × t symmetric supermatrix. Thus $A^TA$ is also a symmetric superbimatrix.

Clearly $AA^T \neq A^TA$ in general.

If both $A_1$ and $A_2$ in $A = A_1 \cup A_2$ are square matrices we want to find out what is $AA^T$ and $A^TA$ in case of semi superbimatrices. To this end we first given an example before we proceed to find the general rule.

*Example 2.64:* Let $A = A_1 \cup A_2$ be a mixed square semi superbimatrix where

$$A_1 = \begin{bmatrix} 3 & 0 & 1 & 1 & 0 \\ 1 & 1 & 0 & 1 & 1 \\ 0 & 0 & 1 & 0 & 1 \\ 2 & 1 & 0 & 1 & 0 \\ 1 & 0 & 1 & 0 & 1 \end{bmatrix}$$

and

$$A_2 = \begin{bmatrix} 0 & 1 & 2 & 3 & 4 & 5 \\ 1 & 2 & 3 & 4 & 5 & 0 \\ \hline 1 & 2 & 0 & 1 & 0 & 1 \\ 0 & 1 & 0 & 0 & 1 & 1 \\ 1 & 0 & 1 & 1 & 0 & 1 \\ 0 & 1 & 0 & 1 & 1 & 0 \end{bmatrix}.$$

We find $A^T = (A_1 \cup A_2)^T = (A_1^T \cup A_2^T)$

$$= \begin{bmatrix} 3 & 1 & 0 & 2 & 1 \\ 0 & 1 & 0 & 1 & 0 \\ 1 & 0 & 1 & 0 & 1 \\ 1 & 1 & 0 & 1 & 0 \\ 0 & 1 & 1 & 0 & 1 \end{bmatrix} \cup \begin{bmatrix} 0 & 1 & 1 & 0 & 1 & 0 \\ 1 & 2 & 2 & 1 & 0 & 1 \\ \hline 2 & 3 & 0 & 0 & 1 & 0 \\ 3 & 4 & 1 & 0 & 1 & 1 \\ 4 & 5 & 0 & 1 & 0 & 1 \\ 5 & 0 & 1 & 1 & 1 & 0 \end{bmatrix}$$

$$\begin{aligned} AA^T &= (A_1 \cup A_2)(A_1 \cup A_2)^T \\ &= (A_1 \cup A_2)(A_1^T \cup A_2^T) \\ &= A_1 A_1^T \cup A_2 A_2^T \end{aligned}$$

$$= \begin{bmatrix} 3 & 0 & 1 & 1 & 0 \\ 1 & 1 & 0 & 1 & 1 \\ 0 & 0 & 1 & 0 & 1 \\ 2 & 1 & 0 & 1 & 0 \\ 1 & 0 & 1 & 0 & 1 \end{bmatrix} \begin{bmatrix} 3 & 1 & 0 & 2 & 1 \\ 0 & 1 & 0 & 1 & 0 \\ 1 & 0 & 1 & 0 & 1 \\ 1 & 1 & 0 & 1 & 0 \\ 0 & 1 & 1 & 0 & 1 \end{bmatrix} \cup$$

$$\left[\begin{array}{cc|ccc|c} 0 & 1 & 2 & 3 & 4 & 5 \\ 1 & 2 & 3 & 4 & 5 & 0 \\ \hline 1 & 2 & 0 & 1 & 0 & 1 \\ 0 & 1 & 0 & 0 & 1 & 1 \\ 1 & 0 & 1 & 1 & 0 & 1 \\ 0 & 1 & 0 & 1 & 1 & 0 \end{array}\right] \left[\begin{array}{cc|cccc} 0 & 1 & 1 & 0 & 1 & 0 \\ 1 & 2 & 2 & 1 & 0 & 1 \\ 2 & 3 & 0 & 0 & 1 & 0 \\ 3 & 4 & 1 & 0 & 1 & 1 \\ 4 & 5 & 0 & 1 & 0 & 1 \\ \hline 5 & 0 & 1 & 1 & 1 & 0 \end{array}\right]$$

$$= \begin{bmatrix} 11 & 4 & 1 & 7 & 4 \\ 4 & 4 & 1 & 4 & 2 \\ 1 & 1 & 2 & 0 & 2 \\ 7 & 4 & 0 & 6 & 2 \\ 4 & 2 & 2 & 2 & 3 \end{bmatrix} \cup$$

$$\left[\begin{array}{c|c} \begin{pmatrix} 0 & 1 \\ 1 & 2 \end{pmatrix}\begin{pmatrix} 0 & 1 \\ 1 & 2 \end{pmatrix} & \begin{pmatrix} 0 & 1 \\ 1 & 2 \end{pmatrix}\begin{pmatrix} 1 & 0 & 1 & 0 \\ 2 & 1 & 0 & 1 \end{pmatrix} \\ \hline \begin{pmatrix} 1 & 2 \\ 0 & 1 \\ 1 & 0 \\ 0 & 1 \end{pmatrix}\begin{pmatrix} 0 & 1 \\ 1 & 2 \end{pmatrix} & \begin{pmatrix} 1 & 2 \\ 0 & 1 \\ 1 & 0 \\ 0 & 1 \end{pmatrix}\begin{pmatrix} 1 & 0 & 1 & 0 \\ 2 & 1 & 0 & 1 \end{pmatrix} \end{array}\right] +$$

$$\left[\begin{array}{c|c} \begin{pmatrix} 2 & 3 & 4 \\ 3 & 4 & 5 \end{pmatrix}\begin{pmatrix} 2 & 3 \\ 3 & 4 \\ 4 & 5 \end{pmatrix} & \begin{pmatrix} 2 & 3 & 4 \\ 3 & 4 & 5 \end{pmatrix}\begin{pmatrix} 0 & 0 & 1 & 0 \\ 1 & 0 & 1 & 1 \\ 0 & 1 & 0 & 1 \end{pmatrix} \\ \hline \begin{pmatrix} 0 & 1 & 0 \\ 0 & 0 & 1 \\ 1 & 1 & 0 \\ 0 & 1 & 1 \end{pmatrix}\begin{pmatrix} 2 & 3 \\ 3 & 4 \\ 4 & 5 \end{pmatrix} & \begin{pmatrix} 0 & 1 & 0 \\ 0 & 0 & 1 \\ 1 & 1 & 0 \\ 0 & 1 & 1 \end{pmatrix}\begin{pmatrix} 0 & 0 & 1 & 0 \\ 1 & 0 & 1 & 1 \\ 0 & 1 & 0 & 1 \end{pmatrix} \end{array}\right] +$$

$$\begin{bmatrix} \begin{pmatrix} 5 \\ 0 \end{pmatrix}(5 \quad 0) & \begin{pmatrix} 5 \\ 0 \end{pmatrix}(1 \quad 1 \quad 1 \quad 0) \\ \hline \begin{pmatrix} 1 \\ 1 \\ 1 \\ 0 \end{pmatrix}(5 \quad 0) & \begin{pmatrix} 1 \\ 1 \\ 1 \\ 0 \end{pmatrix}(1 \quad 1 \quad 1 \quad 0) \end{bmatrix}$$

$$= \begin{bmatrix} 11 & 4 & 1 & 7 & 4 \\ 4 & 4 & 1 & 4 & 2 \\ 1 & 1 & 2 & 0 & 2 \\ 7 & 4 & 0 & 6 & 2 \\ 4 & 2 & 2 & 2 & 3 \end{bmatrix} \cup \begin{bmatrix} 1 & 2 & 2 & 1 & 0 & 1 \\ 2 & 5 & 5 & 2 & 1 & 2 \\ \hline 2 & 5 & 5 & 2 & 1 & 2 \\ 1 & 2 & 2 & 1 & 0 & 1 \\ 0 & 1 & 1 & 0 & 1 & 0 \\ 1 & 2 & 2 & 1 & 0 & 1 \end{bmatrix} +$$

$$\begin{bmatrix} 29 & 38 & 3 & 4 & 5 & 7 \\ 38 & 50 & 4 & 5 & 7 & 9 \\ \hline 3 & 4 & 1 & 0 & 1 & 1 \\ 4 & 5 & 0 & 1 & 0 & 1 \\ 5 & 7 & 1 & 0 & 2 & 1 \\ 7 & 9 & 1 & 1 & 1 & 2 \end{bmatrix} + \begin{bmatrix} 25 & 0 & 5 & 5 & 5 & 0 \\ 0 & 0 & 0 & 0 & 0 & 0 \\ \hline 5 & 0 & 1 & 1 & 1 & 0 \\ 5 & 0 & 1 & 1 & 1 & 0 \\ 5 & 0 & 1 & 1 & 1 & 0 \\ 0 & 0 & 0 & 0 & 0 & 0 \end{bmatrix}$$

$$= \begin{bmatrix} 11 & 4 & 1 & 7 & 4 \\ 4 & 4 & 1 & 4 & 2 \\ 1 & 1 & 2 & 0 & 2 \\ 7 & 4 & 0 & 6 & 2 \\ 4 & 2 & 2 & 2 & 3 \end{bmatrix} \cup \begin{bmatrix} 55 & 40 & 10 & 10 & 10 & 8 \\ 40 & 55 & 9 & 7 & 8 & 11 \\ \hline 10 & 9 & 7 & 3 & 3 & 3 \\ 10 & 7 & 3 & 3 & 1 & 2 \\ 10 & 8 & 3 & 1 & 4 & 1 \\ 8 & 11 & 3 & 2 & 1 & 3 \end{bmatrix}.$$

We see the minor product yields a symmetric semi superbimatrix. Further we see $AA^T \neq A^TA$. The partition of the supermatrix component in $A^TA$ is different from $AA^T$. It is left as an exercise for the reader to find $A^TA$.

**Chapter Three**

# SUPER TRIMATRICES AND THEIR GENERALIZATIONS

In this chapter we introduce the notion of super trimatrices and give here some of their properties and operations on them.

**DEFINITION 3.1:** *$T = T_1 \cup T_2 \cup T_3$ is defined to be a super trimatrix if each of the $T_i$ is a supermatrix for $i = 1, 2, 3$. We demand either each $T_i$ must be a distinct matrix or each $T_i$ must have a distinct partition defined on it, $T_i \neq T_j$ if $i \neq j$, $1 \leq i, j \leq 3$.*

*Example 3.1:* Let $T = T_1 \cup T_2 \cup T_3$ where

$$T_1 = (1\ 0\ 3\ |\ 1\ 1\ 2\ 3\ 4\ |\ 0\ 1\ 3),$$
$$T_2 = (1\ 5\ |\ 5\ 1\ 3\ |\ 3\ 2\ 0\ 1\ 3)$$

and

$$T_3 = (1\ 1\ 1\ 1\ 1\ 1\ |\ 0\ 0\ 1\ 0\ |\ 1\ 1\ 1\ 1\ 1);$$

clearly T is a super trimatrix.

*Example 3.2:* Let $D = D_1 \cup D_2 \cup D_3$ where

$$D_1 = \begin{bmatrix} 1 \\ 2 \\ 3 \\ 4 \\ 7 \\ 8 \\ 9 \end{bmatrix}, D_2 = \begin{bmatrix} 3 \\ 1 \\ 0 \\ 1 \\ 1 \\ 0 \\ 1 \\ 3 \\ 2 \end{bmatrix} \text{ and } D_3 = \begin{bmatrix} 1 \\ 2 \\ 3 \\ 4 \\ 5 \\ 6 \\ 7 \\ 8 \\ 9 \end{bmatrix};$$

D is a super trimatrix.

*Example 3.3:* Let $T = T_1 \cup T_2 \cup T_3$ where

$$T_1 = \left[ \begin{array}{cc|c} 3 & 1 & 2 \\ 1 & 0 & 1 \\ 5 & 1 & 3 \end{array} \right],$$

$$T_2 = \left[ \begin{array}{cccc} 1 & 0 & 1 & 1 \\ 1 & 0 & 1 & 2 \\ \hline 0 & 1 & 0 & 1 \\ 2 & 1 & 0 & 2 \end{array} \right]$$

and

$$T_3 = \left[ \begin{array}{cc|ccc} 0 & 5 & 9 & 7 & 6 \\ 1 & 6 & 8 & 1 & 5 \\ 2 & 7 & 6 & 2 & 4 \\ \hline 3 & 8 & 5 & 3 & 2 \\ 4 & 9 & 4 & 0 & 1 \end{array} \right];$$

T is a super trimatrix.

*Example 3.4:* Let $U = U_1 \cup U_2 \cup U_3$ where

$$U_1 = \begin{bmatrix} 3 & 4 & 5 & 3 & 2 & 1 & 1 \\ 1 & 0 & 1 & 1 & 1 & 3 & 1 \\ 2 & 1 & 1 & 0 & 3 & 4 & 1 \end{bmatrix},$$

$$U_2 = \left[\begin{array}{cccccc|c} 4 & 9 & 4 & 2 & 1 & 5 & 1 \\ 7 & 6 & 2 & 1 & 3 & 6 & 0 \\ \hline 8 & 5 & 0 & 0 & 5 & 7 & 1 \\ 1 & 7 & 1 & 3 & 2 & 8 & 2 \end{array}\right]$$

and

$$U_3 = \begin{bmatrix} 5 & 1 & 3 & 1 \\ 3 & 2 & 1 & 1 \\ 2 & 0 & 0 & 1 \\ \hline 1 & 0 & 6 & 2 \\ 3 & 1 & 0 & 3 \\ 1 & 5 & 1 & 5 \end{bmatrix};$$

U is a super trimatrix.

*Example 3.5:* Let $V = V_1 \cup V_2 \cup V_3$ where

$$V_1 = \begin{bmatrix} 0 & 1 & 2 \\ \hline 1 & 2 & 0 \\ 2 & 1 & 0 \end{bmatrix},$$

$$V_2 = \begin{bmatrix} 3 & 1 \\ 1 & 0 \\ 2 & 1 \\ \hline 0 & 2 \\ 1 & 5 \\ 1 & 6 \end{bmatrix} \text{ and } V_3 = \begin{bmatrix} 1 \\ 2 \\ \hline 1 \\ 3 \\ 4 \\ 7 \end{bmatrix}.$$

V is a super trimatrix.

We have seen 5 examples of super trimatrices.

**DEFINITION 3.2**: *Let $T = T_1 \cup T_2 \cup T_3$ be a super trimatrix we call T to be a row super trimatrix if each of the $T_i$ is a row supermatrix; for i = 1, 2, 3.*

The super trimatrix given in example 3.1 is a row super trimatrix.

**DEFINITION 3.3:** *Let $D = D_1 \cup D_2 \cup D_3$ where D is a super trimatrix. If each of the $D_i$ is a column supermatrix, i = 1, 2, 3 then we call D to be column super trimatrix.*

The super trimatrix given in example 3.2 is a column super trimatrix.

**DEFINITION 3.4:** *Let $V = V_1 \cup V_2 \cup V_3$ be a super trimatrix. If each $V_i$ is a square supermatrix having a different order then V is a mixed square super trimatrix.*

The super trimatrix given in example 3.3 is a mixed square super trimatrix.

**Note:** If in a super trimatrix $T = T_1 \cup T_2 \cup T_3$ if each of the supermatrix is an m × m square supermatrix then we call T to be an m × m square super trimatrix.

***Example 3.6:*** Let $S = S_1 \cup S_2 \cup S_3$ be a super trimatrix where

$$S_1 = \begin{bmatrix} 5 & 0 & 1 & 2 & 3 \\ 0 & 1 & 2 & 3 & 5 \\ \hline 1 & 2 & 3 & 5 & 0 \\ 2 & 3 & 5 & 0 & 1 \\ 3 & 5 & 0 & 1 & 2 \end{bmatrix},$$

$$S_2 = \begin{bmatrix} 1 & 0 & 1 & 0 & 1 \\ 0 & 1 & 1 & 0 & 0 \\ \hline 1 & 0 & 1 & 1 & 1 \\ \hline 0 & 1 & 1 & 1 & 0 \\ 1 & 0 & 1 & 0 & 0 \end{bmatrix}$$

and

$$S_3 = \begin{bmatrix} 9 & 8 & 7 & 6 & 5 \\ \hline 4 & 3 & 2 & 1 & 0 \\ 0 & 1 & 2 & 3 & 4 \\ 5 & 6 & 7 & 8 & 9 \\ \hline 9 & 7 & 5 & 3 & 1 \end{bmatrix}.$$

Clearly S is a square super trimatrix. This is not a mixed square supermatrix.

It is in fact clear that S is a $5 \times 5$ square super trimatrix.

**DEFINITION 3.5:** *Let $M_1 \cup M_2 \cup M_3$ be a super trimatrix. If each $M_i$ is a rectangular supermatrix having a different order for $i = 1, 2, 3$ then we call M to be a mixed rectangular super trimatrix.*

The example given in 3.4 is a mixed rectangular super trimatrix.

*Example 3.7:* Let $T = T_1 \cup T_2 \cup T_3$ be a rectangular super trimatrix where

$$T_1 = \begin{bmatrix} 3 & 1 & 0 & 1 & 1 & 2 & 1 \\ 4 & 0 & 1 & 1 & 0 & 2 & 0 \\ \hline 6 & 2 & 0 & 1 & 1 & 3 & 1 \\ 7 & 1 & 3 & 0 & 0 & 0 & 1 \end{bmatrix}$$

$$T_2 = \begin{bmatrix} 1 & 3 & 6 & 9 & 6 & 3 & 5 \\ 0 & 4 & 7 & 8 & 5 & 2 & 7 \\ 2 & 5 & 8 & 7 & 4 & 1 & 8 \\ 1 & 2 & 0 & 0 & 1 & 3 & 9 \end{bmatrix}$$

and

$$T_3 = \begin{bmatrix} 0 & 4 & 0 & 1 & 1 & 0 & 0 \\ 1 & 5 & 1 & 0 & 0 & 0 & 0 \\ 2 & 6 & 0 & 1 & 0 & 1 & 1 \\ 3 & 1 & 1 & 0 & 1 & 1 & 0 \end{bmatrix}.$$

T is a 4 × 7 rectangular super trimatrix which is not a mixed rectangular super trimatrix.

*Example 3.8:* Let $N = N_1 \cup N_2 \cup N_3$ be a super trimatrix where

$$N_1 = \begin{bmatrix} 1 & 2 & 3 & 1 \\ 1 & 1 & 0 & 1 \\ 0 & 1 & 0 & 1 \\ 1 & 0 & 1 & 0 \\ 1 & 1 & 0 & 0 \\ 1 & 0 & 0 & 1 \\ 0 & 1 & 1 & 0 \\ 0 & 0 & 1 & 1 \\ 1 & 1 & 1 & 1 \end{bmatrix}, \quad N_2 = \begin{bmatrix} 1 & 2 \\ 2 & 1 \\ 1 & 3 \\ 3 & 1 \\ 1 & 1 \\ 2 & 2 \\ 3 & 3 \\ 4 & 4 \\ 5 & 5 \\ 6 & 6 \\ 1 & 0 \\ 0 & 1 \\ 1 & 1 \\ 0 & 0 \end{bmatrix}$$

and

$$N_3 = \begin{bmatrix} 0 & 1 & 2 & 3 & 4 \\ 1 & 2 & 3 & 4 & 0 \\ 2 & 3 & 4 & 0 & 1 \\ \hline 1 & 1 & 4 & 0 & 2 \\ 0 & 1 & 0 & 2 & 5 \\ 9 & 8 & 7 & 7 & 6 \\ \hline 6 & 6 & 5 & 5 & 4 \\ 4 & 4 & 3 & 3 & 2 \\ 2 & 2 & 1 & 1 & 5 \\ 5 & 5 & 4 & 4 & 3 \end{bmatrix}.$$

We see each $N_i$ is partitioned only horizontally and never vertically i.e., each $N_i$ is partitioned only along the rows for i = 1, 2, 3.

*Example 3.9:* Let $C = C_1 \cup C_2 \cup C_3$ be a super trimatrix; where

$$C_1 = \begin{bmatrix} 1 & 3 & 2 & 1 & 2 & 3 & 1 & 2 & 2 & 3 \\ 2 & 1 & 1 & 1 & 2 & 3 & 2 & 2 & 1 & 1 \\ 3 & 2 & 3 & 1 & 2 & 3 & 2 & 1 & 2 & 2 \end{bmatrix}$$

$$C_2 = \begin{bmatrix} 2 & 1 & 8 & 1 & 2 & 1 & 1 & 4 & 1 \\ 1 & 2 & 7 & 1 & 8 & 7 & 1 & -5 & 8 \\ 3 & 9 & 3 & 1 & 1 & 8 & 9 & 8 & 9 \\ 5 & 6 & 2 & 5 & 9 & 1 & 4 & 0 & 8 \\ 0 & 7 & 5 & 3 & 3 & 1 & 7 & -1 & 1 \end{bmatrix}$$

and

$$C_3 = \begin{bmatrix} 2 & 1 & 3 & 4 & 0 & 9 & 3 & 8 & 1 & 2 \\ 3 & 1 & 1 & 5 & 1 & 2 & 7 & 1 & 9 & 3 \end{bmatrix}.$$

Each $C_i$ is partitioned only vertically and never horizontally for i = 1, 2, 3.

**DEFINITION 3.6:** *Let $T = T_1 \cup T_2 \cup T_3$ be a super trimatrix. Of each $T_i$ is only partitioned horizontally i.e., only along the rows or in between the rows, then we call T to be a column super trivector.*

The super trimatrix given in example 3.8 is a column super trivector.

**DEFINITION 3.7:** *Let $T = T_1 \cup T_2 \cup T_3$ be a super trimatrix. If each of $T_i$ is only partitioned vertically i.e., only along the columns or in between the columns, then we call T to be a row super trivector.*

The super trimatrix C given in example 3.9 is a row super trivector.

Now having defined several types of super trimatrices we define operations on them. It is important to mention here that even it is very difficult to define addition of super trimatrices for if we need to define addition we need not only have the order to be the same but also the partition defined on them must be identical otherwise we cannot even define simple addition of super trimatrices.

*Example 3.10:* Let
$$T = T_1 \cup T_2 \cup T_3$$

$$= \begin{bmatrix} 3 \\ 1 \\ 2 \\ \hline 3 \\ 0 \\ 5 \\ 6 \end{bmatrix} \cup \begin{bmatrix} 3 & 1 & 2 & 4 \\ 0 & 1 & 1 & 1 \\ 1 & 1 & 1 & 0 \\ 1 & 0 & 1 & 1 \\ \hline 1 & 1 & 1 & 0 \\ 0 & 1 & 0 & 1 \end{bmatrix} \cup \left[ \begin{array}{ccc|cc} 3 & 0 & 1 & 3 & 5 \\ 1 & 2 & 0 & 1 & 1 \\ 0 & 1 & 0 & 1 & 0 \\ \hline 5 & 0 & 1 & 0 & 1 \\ 1 & 2 & 3 & 4 & 5 \end{array} \right]$$

be a super trimatrix and

$$S = S_1 \cup S_2 \cup S_3$$

$$= \begin{bmatrix} 4 \\ 0 \\ 1 \\ \overline{2} \\ 2 \\ 3 \\ 1 \end{bmatrix} \cup \begin{bmatrix} 0 & 1 & 1 & 2 \\ 1 & 1 & 1 & 3 \\ 2 & 0 & 3 & 1 \\ 1 & 1 & 5 & 0 \\ \overline{1} & 2 & 0 & 1 \\ 3 & 4 & 1 & 1 \end{bmatrix} \cup \left[ \begin{array}{cccccc} 0 & 1 & 2 & 3 & 4 \\ 0 & 1 & 2 & 0 & 1 \\ 2 & 1 & 0 & 2 & 1 \\ 1 & 2 & 1 & 2 & 2 \\ \hline 3 & 0 & 0 & 3 & 0 \end{array} \right]$$

another super trimatrix. We see in both T and S each $T_i$ and $S_i$ are supermatrices of same order with identical partitions defined on them $1 \leq i \leq 3$. Thus addition of S and T or T and S is defined

$$T + S = (T_1 \cup T_2 \cup T_3) + (S_1 \cup S_2 + S_3)$$
$$= (T_1 + S_1) \cup (T_2 + S_2) \cup (T_3 + S_3)$$

$$= \left\{ \begin{bmatrix} 3 \\ \overline{1} \\ 2 \\ \overline{3} \\ 0 \\ 5 \\ 6 \end{bmatrix} + \begin{bmatrix} 4 \\ 0 \\ 1 \\ \overline{2} \\ 2 \\ 3 \\ 1 \end{bmatrix} \right\} \cup \left\{ \begin{bmatrix} 3 & 1 & 2 & 4 \\ 0 & 1 & 1 & 1 \\ 1 & 1 & 1 & 0 \\ 1 & 0 & 1 & 1 \\ \hline 1 & 1 & 1 & 0 \\ 0 & 1 & 0 & 1 \end{bmatrix} + \begin{bmatrix} 0 & 1 & 1 & 2 \\ 1 & 1 & 1 & 3 \\ 2 & 0 & 3 & 1 \\ 1 & 1 & 5 & 0 \\ \hline 1 & 2 & 0 & 1 \\ 3 & 4 & 1 & 1 \end{bmatrix} \right\} \cup$$

$$\left\{ \begin{bmatrix} 3 & 0 & 1 & 3 & 5 \\ 1 & 2 & 0 & 1 & 1 \\ 0 & 1 & 0 & 1 & 0 \\ \hline 5 & 0 & 1 & 0 & 1 \\ 1 & 2 & 3 & 4 & 5 \end{bmatrix} + \begin{bmatrix} 0 & 1 & 2 & 3 & 4 \\ 0 & 1 & 2 & 0 & 1 \\ 2 & 1 & 0 & 2 & 1 \\ 1 & 2 & 1 & 2 & 2 \\ \hline 3 & 0 & 0 & 3 & 0 \end{bmatrix} \right\}$$

$$= \begin{bmatrix} 7 \\ 1 \\ 3 \\ \overline{5} \\ 2 \\ 8 \\ 7 \end{bmatrix} \cup \begin{bmatrix} 3 & 2 & 3 & 6 \\ 0 & 2 & 2 & 4 \\ 3 & 1 & 4 & 1 \\ 2 & 1 & 6 & 1 \\ \overline{2} & 3 & 1 & 1 \\ 3 & 5 & 1 & 2 \end{bmatrix} \cup \begin{bmatrix} 3 & 1 & 3 & 6 & 9 \\ 1 & 3 & 2 & 1 & 2 \\ 2 & 2 & 0 & 3 & 1 \\ \overline{6} & 2 & 2 & 2 & 3 \\ 4 & 2 & 3 & 7 & 5 \end{bmatrix}.$$

T + S also happens to be the same type of super trimatrix.

*Note:* If $T = T_1 \cup T_2 \cup T_3$ is a super trimatrix then

$$\begin{aligned} T + T &= (T_1 \cup T_2 \cup T_3) + (T_1 \cup T_2 \cup T_3) \\ &= (T_1 + T_1) \cup (T_2 + T_2) \cup (T_3 + T_3) \\ &= 2T_1 \cup 2T_2 \cup 2T_3 \\ &= 2T \end{aligned}$$

is a super trimatrix. Thus if we take $T + T + \ldots + T$(n times) $=$ $nT = nT_1 \cup nT_2 \cup nT_3$.

*Example 3.11:* Let

$$T = T_1 \cup T_2 \cup T_3$$

$$= \begin{bmatrix} 2 & 1 \\ \hline 0 & 1 \\ 1 & 2 \\ 2 & 1 \end{bmatrix} \cup \begin{bmatrix} 1 & 1 & 3 & 1 \\ 0 & 1 & 2 & 4 \\ 1 & 2 & 3 & 0 \\ 2 & 0 & 1 & 1 \end{bmatrix} \cup \begin{bmatrix} 1 & 0 & 1 & 1 & 0 & 1 \\ 2 & 1 & 1 & 2 & 0 & 0 \\ \hline 3 & 1 & 3 & 0 & 0 & 0 \\ 4 & 0 & 4 & 1 & 0 & 0 \\ 5 & 2 & 5 & 1 & 1 & 0 \end{bmatrix}$$

be a super trimatrix.

$$\begin{aligned} T + T &= (T_1 \cup T_2 \cup T_3) + (T_1 \cup T_2 \cup T_3) \\ &= (T_1 + T_1) \cup (T_2 + T_2) \cup (T_3 + T_3) \end{aligned}$$

$$= \left\{ \begin{bmatrix} 2 & 1 \\ \hline 0 & 1 \\ 1 & 2 \\ 2 & 1 \end{bmatrix} + \begin{bmatrix} 2 & 1 \\ \hline 0 & 1 \\ 1 & 2 \\ 2 & 1 \end{bmatrix} \right\} \cup \left\{ \begin{bmatrix} 1 & 1 & 3 & 1 \\ 0 & 1 & 2 & 4 \\ \hline 1 & 2 & 3 & 0 \\ 2 & 0 & 1 & 1 \end{bmatrix} + \begin{bmatrix} 1 & 1 & 3 & 1 \\ 0 & 1 & 2 & 4 \\ \hline 1 & 2 & 3 & 0 \\ 2 & 0 & 1 & 1 \end{bmatrix} \right\}$$

$$\cup \left\{ \begin{bmatrix} 1 & 0 & 1 & 1 & 0 & 1 \\ 2 & 1 & 1 & 2 & 0 & 0 \\ \hline 3 & 1 & 3 & 0 & 0 & 0 \\ 4 & 0 & 4 & 1 & 0 & 0 \\ 5 & 2 & 5 & 1 & 1 & 0 \end{bmatrix} + \begin{bmatrix} 1 & 0 & 1 & 1 & 0 & 1 \\ 2 & 1 & 1 & 2 & 0 & 0 \\ 3 & 1 & 3 & 0 & 0 & 0 \\ 4 & 0 & 4 & 1 & 0 & 0 \\ 5 & 2 & 5 & 1 & 1 & 0 \end{bmatrix} \right\}$$

$$= \begin{bmatrix} 4 & 2 \\ \hline 0 & 2 \\ 2 & 4 \\ 4 & 2 \end{bmatrix} \cup \begin{bmatrix} 2 & 2 & 6 & 2 \\ 0 & 2 & 4 & 8 \\ \hline 2 & 4 & 6 & 0 \\ 4 & 0 & 2 & 2 \end{bmatrix} \cup \begin{bmatrix} 2 & 0 & 2 & 2 & 0 & 2 \\ 4 & 2 & 2 & 4 & 0 & 0 \\ 6 & 2 & 6 & 0 & 0 & 0 \\ 8 & 0 & 8 & 2 & 0 & 0 \\ 10 & 4 & 10 & 2 & 2 & 0 \end{bmatrix}$$

$$= 2T_1 \cup 2T_2 \cup 2T_3.$$

It is easily verified $5T_1 \cup 5T_2 \cup 5T_3$

$$= \begin{bmatrix} 10 & 5 \\ \hline 0 & 5 \\ 5 & 10 \\ 10 & 5 \end{bmatrix} \cup \begin{bmatrix} 5 & 5 & 15 & 5 \\ 0 & 5 & 10 & 20 \\ \hline 5 & 10 & 15 & 0 \\ 10 & 0 & 15 & 0 \end{bmatrix} \cup \begin{bmatrix} 5 & 0 & 5 & 5 & 0 & 5 \\ 10 & 5 & 5 & 10 & 0 & 0 \\ \hline 15 & 5 & 15 & 0 & 0 & 0 \\ 20 & 0 & 20 & 5 & 0 & 0 \\ 25 & 10 & 25 & 5 & 5 & 0 \end{bmatrix}.$$

Now we can define product of two super trimatrices in many ways.

*Example 3.12:* Let
$$T = T_1 \cup T_2 \cup T_3$$

$$= \begin{bmatrix} 3 & 1 \\ 1 & 1 \\ 2 & 0 \\ 3 & 1 \\ 5 & 0 \\ 6 & 0 \end{bmatrix} \cup \left[ \begin{array}{cc|cc} 3 & 1 & 2 & 1 \\ 0 & 1 & 3 & 1 \\ \hline 1 & 1 & 0 & 2 \\ 0 & 0 & 1 & 0 \end{array} \right] \cup \left[ \begin{array}{cc|cccc} 3 & 1 & 2 & 1 & 0 & 1 & 1 \\ 1 & 0 & 1 & 1 & 1 & 1 & 2 \\ 0 & 1 & 0 & 1 & 1 & 0 & 1 \end{array} \right]$$

be a super trimatrix and

$$V = V_1 \cup V_2 \cup V_3$$

$$= \left[ \begin{array}{ccccc|ccc} 1 & 1 & 2 & 4 & 7 & 1 & 1 & 0 \\ 0 & 1 & 3 & 5 & 2 & 1 & 0 & 1 \end{array} \right] \cup \left[ \begin{array}{cc|cc} 0 & 1 & 2 & 3 \\ 1 & 1 & 1 & 1 \\ \hline 0 & 1 & 0 & 1 \\ 2 & 1 & 1 & 0 \end{array} \right] \cup \begin{bmatrix} 0 & 1 & 2 \\ 1 & 0 & 3 \\ 1 & 1 & 1 \\ 0 & 0 & 0 \\ 2 & 1 & 3 \\ 7 & 5 & 2 \\ 0 & 2 & 1 \end{bmatrix}$$

be another super trimatrix. Now the product

$$\begin{aligned} TV &= (T_1 \cup T_2 \cup T_3)(V_1 \cup V_2 \cup V_3) \\ &= T_1V_1 \cup T_2V_2 \cup T_3V_3 \end{aligned}$$

$$= \left\{ \begin{bmatrix} 3 & 1 \\ 1 & 1 \\ 2 & 0 \\ 3 & 1 \\ 5 & 0 \\ 6 & 0 \end{bmatrix} \left[ \begin{array}{ccccc|cccc} 1 & 1 & 2 & 4 & 7 & 1 & 1 & 0 & 1 & 0 \\ 0 & 1 & 3 & 5 & 2 & 1 & 0 & 1 & 1 & 0 \end{array} \right] \right\} \cup$$

$$\left\{ \left[ \begin{array}{cc|cc} 3 & 1 & 2 & 1 \\ 0 & 1 & 3 & 1 \\ \hline 1 & 1 & 0 & 2 \\ 0 & 0 & 1 & 0 \end{array} \right] \left[ \begin{array}{cc|cc} 0 & 1 & 2 & 3 \\ 1 & 1 & 1 & 1 \\ \hline 0 & 1 & 0 & 1 \\ 2 & 1 & 1 & 0 \end{array} \right] \right\}$$

$$\cup \left\{ \left[ \begin{array}{cc|cccccc} 3 & 1 & 2 & 1 & 0 & 1 & 1 \\ 1 & 0 & 1 & 1 & 1 & 1 & 2 \\ \hline 0 & 1 & 0 & 1 & 1 & 0 & 1 \end{array} \right] \left[ \begin{array}{ccc} 0 & 1 & 2 \\ 1 & 0 & 3 \\ \hline 1 & 1 & 1 \\ 0 & 0 & 0 \\ 2 & 1 & 3 \\ 7 & 5 & 2 \\ 0 & 2 & 1 \end{array} \right] \right\}$$

$$= \left[ \begin{array}{c} \left[ \begin{array}{cc} 3 & 1 \\ 1 & 1 \\ 2 & 0 \\ 3 & 1 \end{array} \right] \left[ \begin{array}{ccccc} 1 & 1 & 2 & 4 & 7 \\ 0 & 1 & 3 & 5 & 2 \end{array} \right] \\ \hline \left[ \begin{array}{cc} 5 & 0 \\ 6 & 0 \end{array} \right] \left[ \begin{array}{ccccc} 1 & 1 & 0 & 1 & 0 \\ 1 & 0 & 1 & 1 & 0 \end{array} \right] \end{array} \right]$$

$$\cup \left[ \begin{array}{cc|cc} \begin{pmatrix} 3 & 1 \\ 0 & 1 \end{pmatrix} \begin{pmatrix} 0 & 1 \\ 1 & 1 \end{pmatrix} & \begin{pmatrix} 2 & 1 \\ 3 & 1 \end{pmatrix} \begin{pmatrix} 0 & 1 \\ 2 & 1 \end{pmatrix} \\ \hline \begin{pmatrix} 1 & 1 \\ 0 & 0 \end{pmatrix} \begin{pmatrix} 0 & 1 \\ 2 & 1 \end{pmatrix} & \begin{pmatrix} 0 & 2 \\ 1 & 0 \end{pmatrix} \begin{pmatrix} 0 & 1 \\ 1 & 0 \end{pmatrix} \end{array} \right]$$

$$\cup \begin{bmatrix} \begin{bmatrix} 3 & 1 \\ 1 & 0 \\ 0 & 1 \end{bmatrix} \begin{bmatrix} 0 & 1 & 2 \\ 1 & 0 & 3 \end{bmatrix} \\ \hline \begin{bmatrix} 2 & 1 & 0 & 1 & 1 \\ 1 & 1 & 1 & 1 & 2 \\ 0 & 1 & 1 & 0 & 1 \end{bmatrix} \begin{bmatrix} 1 & 1 & 1 \\ 0 & 0 & 0 \\ 2 & 1 & 3 \\ 7 & 5 & 2 \\ 0 & 2 & 1 \end{bmatrix} \end{bmatrix}$$

$$= \begin{bmatrix} 3 & 4 & 9 & 17 & 23 \\ 1 & 2 & 5 & 9 & 9 \\ 2 & 4 & 4 & 8 & 14 \\ 9 & 4 & 9 & 17 & 23 \\ \hline 5 & 5 & 0 & 5 & 0 \\ 6 & 6 & 0 & 6 & 0 \end{bmatrix} \cup \begin{bmatrix} 1 & 4 & 2 & 3 \\ 1 & 1 & 2 & 4 \\ 2 & 2 & 2 & 0 \\ \hline 0 & 0 & 0 & 1 \end{bmatrix} \cup \begin{bmatrix} 1 & 3 & 9 \\ 0 & 1 & 2 \\ 1 & 0 & 3 \\ 9 & 9 & 5 \\ \hline 10 & 11 & 8 \\ 2 & 3 & 4 \end{bmatrix}.$$

Thus we are able to define a product but in all cases we may not be in a position to define a product of super trimatrices, which is a column super trimatrix.

*Example 3.13:* Let $A = A_1 \cup A_2 \cup A_3$ be a super trimatrix where

$$A_1 = \begin{bmatrix} 3 & 1 & 0 \\ 2 & 1 & 2 \\ 3 & 1 & 0 \\ \hline 1 & 1 & 1 \\ 0 & 0 & 2 \\ \hline 1 & 0 & 1 \\ 0 & 1 & 0 \end{bmatrix},$$

$$A_2 = \begin{bmatrix} 1 & 1 \\ 2 & 0 \\ 3 & 3 \\ \hline 1 & 4 \\ 2 & 1 \\ 3 & 1 \\ 0 & 1 \\ \hline 1 & 5 \end{bmatrix}$$

and

$$A_3 = \begin{bmatrix} 2 & 1 & 3 & 4 \\ 0 & 1 & 1 & 0 \\ \hline 1 & 0 & 1 & 0 \\ 0 & 1 & 0 & 1 \\ 1 & 1 & 1 & 1 \\ \hline 1 & 1 & 0 & 1 \\ 2 & 0 & 0 & 1 \\ 1 & 0 & 0 & 2 \\ 0 & 2 & 1 & 0 \end{bmatrix}$$

which is a column super trivector. Let $B = B_1 \cup B_2 \cup B_3$ be a row super trivector where

$$B_1 = \begin{bmatrix} 1 & 0 & 1 & 2 & 1 & 1 & 2 \\ 1 & 1 & 0 & 0 & 0 & 3 & 4 \\ 0 & 1 & 1 & 1 & 1 & 0 & 5 \end{bmatrix},$$

$$B_2 = \begin{bmatrix} 2 & 1 & 0 & 1 & 1 & 3 & 1 & 2 \\ 2 & 1 & 1 & 0 & 3 & 1 & 0 & 1 \end{bmatrix}$$

and

$$B_3 = \begin{bmatrix} 0 & 1 & | & 0 & 1 & 1 & | & 2 & 1 & 1 & 0 \\ 1 & 0 & | & 0 & 2 & 0 & | & 1 & 0 & 0 & 0 \\ 2 & 2 & | & 1 & 0 & 1 & | & 0 & 1 & 2 & 1 \\ 1 & 0 & | & 2 & 0 & 1 & | & 1 & 2 & 1 & 0 \end{bmatrix}.$$

Now

$$\begin{aligned}
AB &= (A_1 \cup A_2 \cup A_3)(B_1 \cup B_2 \cup B_3) \\
&= A_1B_1 \cup A_2B_2 \cup A_3B_3
\end{aligned}$$

$$= \begin{bmatrix} 3 & 1 & 0 \\ 2 & 1 & 2 \\ 3 & 1 & 0 \\ \hline 1 & 1 & 1 \\ 0 & 0 & 2 \\ \hline 1 & 0 & 1 \\ 0 & 1 & 0 \end{bmatrix} \begin{bmatrix} 1 & 0 & 1 & | & 2 & 1 & | & 1 & 2 \\ 1 & 1 & 0 & | & 0 & 0 & | & 3 & 4 \\ 0 & 1 & 1 & | & 1 & 1 & | & 0 & 5 \end{bmatrix} \cup$$

$$\left\{ \begin{bmatrix} 1 & 1 \\ 2 & 0 \\ 3 & 3 \\ 1 & 4 \\ \hline 2 & 1 \\ 3 & 1 \\ 0 & 1 \\ \hline 1 & 5 \end{bmatrix} \begin{bmatrix} 2 & 1 & 0 & 1 & | & 1 & 3 & 1 & | & 2 \\ 2 & 1 & 1 & 0 & | & 3 & 1 & 0 & | & 1 \end{bmatrix} \right\} \cup$$

$$\left\{ \begin{bmatrix} 2 & 1 & 3 & 4 \\ 0 & 1 & 1 & 0 \\ \hline 1 & 0 & 1 & 0 \\ 0 & 1 & 0 & 1 \\ 1 & 1 & 1 & 1 \\ \hline 1 & 1 & 0 & 1 \\ 2 & 0 & 0 & 1 \\ 1 & 0 & 0 & 2 \\ 0 & 2 & 1 & 0 \end{bmatrix} \begin{bmatrix} 0 & 1 & | & 0 & 1 & 1 & | & 2 & 1 & 1 & 0 \\ 1 & 0 & | & 0 & 2 & 0 & | & 1 & 0 & 0 & 0 \\ 2 & 2 & | & 1 & 0 & 1 & | & 0 & 1 & 2 & 1 \\ 1 & 0 & | & 2 & 0 & 1 & | & 1 & 2 & 1 & 0 \end{bmatrix} \right\} =$$

$$\begin{bmatrix} \begin{pmatrix} 3 & 1 & 0 \\ 2 & 1 & 2 \\ 3 & 1 & 0 \end{pmatrix}\begin{pmatrix} 1 & 0 & 1 \\ 1 & 1 & 0 \\ 0 & 1 & 1 \end{pmatrix} & \begin{pmatrix} 3 & 1 & 0 \\ 2 & 1 & 2 \\ 3 & 1 & 0 \end{pmatrix}\begin{pmatrix} 2 & 1 \\ 0 & 0 \\ 1 & 1 \end{pmatrix} & \begin{pmatrix} 3 & 1 & 0 \\ 2 & 1 & 2 \\ 3 & 1 & 0 \end{pmatrix}\begin{pmatrix} 1 & 2 \\ 3 & 4 \\ 0 & 5 \end{pmatrix} \\ \hline \begin{pmatrix} 1 & 1 & 1 \\ 0 & 0 & 2 \end{pmatrix}\begin{pmatrix} 1 & 0 & 1 \\ 1 & 1 & 0 \\ 0 & 1 & 1 \end{pmatrix} & \begin{pmatrix} 1 & 1 & 0 \\ 0 & 0 & 2 \end{pmatrix}\begin{pmatrix} 2 & 1 \\ 0 & 0 \\ 1 & 1 \end{pmatrix} & \begin{pmatrix} 1 & 1 & 1 \\ 0 & 0 & 2 \end{pmatrix}\begin{pmatrix} 1 & 2 \\ 3 & 4 \\ 0 & 5 \end{pmatrix} \\ \hline \begin{pmatrix} 1 & 0 & 1 \\ 0 & 1 & 0 \end{pmatrix}\begin{pmatrix} 1 & 0 & 1 \\ 1 & 1 & 0 \\ 0 & 1 & 1 \end{pmatrix} & \begin{pmatrix} 1 & 0 & 1 \\ 0 & 1 & 0 \end{pmatrix}\begin{pmatrix} 2 & 1 \\ 0 & 0 \\ 1 & 1 \end{pmatrix} & \begin{pmatrix} 1 & 0 & 1 \\ 0 & 1 & 0 \end{pmatrix}\begin{pmatrix} 1 & 2 \\ 3 & 4 \\ 0 & 5 \end{pmatrix} \end{bmatrix}$$

$$\cup \begin{bmatrix} \begin{pmatrix} 1 & 1 \\ 2 & 0 \\ 3 & 3 \\ 1 & 4 \end{pmatrix}\begin{pmatrix} 2 & 1 & 0 & 1 \\ 2 & 1 & 1 & 0 \end{pmatrix} & \begin{pmatrix} 1 & 1 \\ 2 & 0 \\ 3 & 3 \\ 1 & 4 \end{pmatrix}\begin{pmatrix} 1 & 3 & 1 \\ 3 & 1 & 0 \end{pmatrix} & \begin{pmatrix} 1 & 1 \\ 2 & 0 \\ 3 & 3 \\ 1 & 4 \end{pmatrix}\begin{pmatrix} 2 \\ 1 \end{pmatrix} \\ \hline \begin{pmatrix} 2 & 1 \\ 3 & 1 \\ 0 & 1 \end{pmatrix}\begin{pmatrix} 2 & 1 & 0 & 1 \\ 2 & 1 & 1 & 0 \end{pmatrix} & \begin{pmatrix} 2 & 1 \\ 3 & 1 \\ 0 & 1 \end{pmatrix}\begin{pmatrix} 1 & 3 & 1 \\ 3 & 1 & 0 \end{pmatrix} & \begin{pmatrix} 2 & 1 \\ 3 & 1 \\ 0 & 1 \end{pmatrix}\begin{pmatrix} 2 \\ 1 \end{pmatrix} \\ \hline \begin{pmatrix} 1 & 5 \end{pmatrix}\begin{pmatrix} 2 & 1 & 0 & 1 \\ 2 & 1 & 1 & 0 \end{pmatrix} & \begin{pmatrix} 1 & 5 \end{pmatrix}\begin{pmatrix} 1 & 3 & 1 \\ 3 & 1 & 0 \end{pmatrix} & \begin{pmatrix} 1 & 5 \end{pmatrix}\begin{pmatrix} 2 \\ 1 \end{pmatrix} \end{bmatrix} \cup$$

$$\begin{bmatrix} \begin{pmatrix} 2 & 1 & 3 & 4 \\ 0 & 1 & 1 & 0 \end{pmatrix} \begin{bmatrix} 0 & 1 \\ 1 & 0 \\ 2 & 2 \\ 1 & 0 \end{bmatrix} & \begin{pmatrix} 2 & 1 & 3 & 4 \\ 0 & 1 & 1 & 0 \end{pmatrix} \begin{pmatrix} 0 & 1 & 1 \\ 0 & 2 & 0 \\ 1 & 0 & 1 \\ 2 & 0 & 1 \end{pmatrix} & \begin{pmatrix} 2 & 1 & 3 & 4 \\ 0 & 1 & 1 & 0 \end{pmatrix} \begin{bmatrix} 2 & 1 & 1 & 0 \\ 1 & 0 & 0 & 0 \\ 0 & 1 & 2 & 1 \\ 1 & 2 & 1 & 0 \end{bmatrix} \\ \hline \begin{pmatrix} 1 & 0 & 1 & 0 \\ 0 & 1 & 0 & 1 \\ 1 & 1 & 1 & 1 \end{pmatrix} \begin{bmatrix} 0 & 1 \\ 1 & 0 \\ 2 & 2 \\ 1 & 0 \end{bmatrix} & \begin{pmatrix} 1 & 0 & 1 & 0 \\ 0 & 1 & 0 & 1 \\ 1 & 1 & 1 & 1 \end{pmatrix} \begin{pmatrix} 0 & 1 & 1 \\ 0 & 2 & 0 \\ 1 & 0 & 1 \\ 2 & 0 & 1 \end{pmatrix} & \begin{pmatrix} 1 & 0 & 1 & 0 \\ 0 & 1 & 0 & 1 \\ 1 & 1 & 1 & 1 \end{pmatrix} \begin{bmatrix} 2 & 1 & 1 & 0 \\ 1 & 0 & 0 & 0 \\ 0 & 1 & 2 & 1 \\ 1 & 2 & 1 & 0 \end{bmatrix} \\ \hline \begin{pmatrix} 1 & 0 & 1 & 0 \\ 0 & 1 & 0 & 1 \\ 1 & 1 & 1 & 1 \end{pmatrix} \begin{bmatrix} 0 & 1 \\ 1 & 0 \\ 2 & 2 \\ 1 & 0 \end{bmatrix} & \begin{pmatrix} 1 & 1 & 0 & 1 \\ 2 & 0 & 0 & 1 \\ 1 & 0 & 0 & 2 \\ 0 & 2 & 1 & 0 \end{pmatrix} \begin{pmatrix} 0 & 1 & 1 \\ 0 & 2 & 0 \\ 1 & 0 & 1 \\ 2 & 0 & 1 \end{pmatrix} & \begin{pmatrix} 1 & 1 & 0 & 1 \\ 2 & 0 & 0 & 1 \\ 1 & 0 & 0 & 2 \\ 0 & 2 & 1 & 0 \end{pmatrix} \begin{bmatrix} 2 & 1 & 1 & 0 \\ 1 & 0 & 0 & 0 \\ 0 & 1 & 2 & 1 \\ 1 & 2 & 1 & 0 \end{bmatrix} \end{bmatrix}$$

$$= \begin{bmatrix} \begin{array}{ccc|cc|cc} 4 & 1 & 3 & 6 & 3 & 6 & 10 \\ 3 & 3 & 4 & 6 & 4 & 5 & 18 \\ 4 & 1 & 3 & 6 & 3 & 6 & 10 \\ \hline 2 & 2 & 2 & 2 & 1 & 4 & 11 \\ 0 & 2 & 2 & 2 & 2 & 0 & 10 \\ \hline 1 & 1 & 2 & 3 & 2 & 1 & 7 \\ 1 & 1 & 0 & 0 & 0 & 3 & 4 \end{array} \end{bmatrix} \cup$$

$$\begin{bmatrix} \begin{array}{cccc|ccc|c} 4 & 2 & 1 & 1 & 4 & 4 & 1 & 3 \\ 4 & 2 & 0 & 2 & 2 & 6 & 2 & 4 \\ 12 & 6 & 3 & 3 & 12 & 12 & 3 & 9 \\ 6 & 5 & 4 & 1 & 13 & 7 & 1 & 6 \\ \hline 5 & 3 & 1 & 2 & 5 & 7 & 2 & 5 \\ 8 & 4 & 1 & 3 & 6 & 10 & 3 & 7 \\ 2 & 1 & 1 & 0 & 3 & 1 & 0 & 1 \\ \hline 12 & 6 & 5 & 1 & 16 & 8 & 1 & 7 \end{array} \end{bmatrix}$$

$$\cup \begin{bmatrix} 11 & 8 & 11 & 4 & 9 & 9 & 13 & 12 & 3 \\ 3 & 2 & 1 & 2 & 1 & 1 & 1 & 2 & 1 \\ \hline 2 & 3 & 1 & 1 & 2 & 2 & 2 & 3 & 1 \\ 2 & 0 & 2 & 2 & 1 & 2 & 2 & 1 & 0 \\ 4 & 3 & 3 & 3 & 3 & 4 & 4 & 4 & 1 \\ \hline 2 & 3 & 2 & 3 & 2 & 4 & 3 & 2 & 0 \\ 2 & 0 & 2 & 2 & 3 & 5 & 4 & 3 & 0 \\ 2 & 1 & 4 & 1 & 3 & 4 & 5 & 3 & 0 \\ 4 & 3 & 1 & 4 & 1 & 2 & 1 & 2 & 1 \end{bmatrix}.$$

We see the product of a column super trivector with a row super trivector when defined results in a super trimatrix which is not a super trivector.

We proceed on to define some more concepts before we define the product of super trimatrix with its transpose and so on.

**DEFINITION 3.8:** *Let $A = A_1 \cup A_2 \cup A_3$ be a super trimatrix. Then the transpose of A denoted by $A^T = (A_1 \cup A_2 \cup A_3)^T = A_1^T \cup A_2^T \cup A_3^T$ is again a super trimatrix.*

*Example 3.14:* Let $T = T_1 \cup T_2 \cup T_3$

$$= \begin{bmatrix} 0 & 1 & 2 & 3 & 4 \\ 1 & 1 & 0 & 1 & 2 \\ \hline 9 & 6 & 8 & 4 & 2 \\ 5 & 4 & 0 & 1 & 5 \end{bmatrix} \cup \begin{bmatrix} 3 & 1 & 2 & 0 & 5 & 7 & 8 \\ 1 & 6 & 2 & 3 & 4 & 5 & 9 \\ \hline 0 & 1 & 2 & 3 & 1 & 1 & 1 \\ 9 & 1 & 0 & 8 & 7 & 0 & 9 \\ 2 & 4 & 2 & 0 & 1 & 5 & 1 \end{bmatrix}$$

$$\cup \begin{bmatrix} 2 & 0 & 1 & 0 & 0 \\ 3 & 0 & 2 & 1 & 0 \\ 1 & 3 & 0 & 1 & 5 \\ 4 & 1 & 5 & 2 & 1 \end{bmatrix};$$

we see $T_1$ and $T_2$ are supermatrices where as $T_3$ is a simple(ordinary) matrix. So T is not a trimatrix or a super trimatrix.

*Example 3.15:* Let $T = T_1 \cup T_2 \cup T_3$

$$\begin{bmatrix} 4 & 1 & 1 \\ 0 & 2 & 3 \\ 3 & 5 & 4 \\ 1 & 1 & 0 \\ 0 & 9 & 5 \\ 9 & 2 & 6 \end{bmatrix} \cup \left[\begin{array}{ccc|cccc} 2 & 3 & 1 & 5 & 0 & 1 & 2 & 1 \\ 5 & 0 & 1 & 0 & 1 & 0 & 1 & 0 \\ \hline 7 & 2 & 3 & 3 & 2 & 1 & 2 & 0 \\ 2 & 1 & 3 & 2 & 1 & 1 & 1 & 2 \\ 0 & 5 & 1 & 1 & 0 & 2 & 3 & 4 \\ \hline 1 & 6 & 0 & 5 & 5 & 7 & 5 & 8 \end{array}\right].$$

$$\begin{bmatrix} 3 & 5 & 1 & 3 & 1 & 5 & 1 & 0 \\ 2 & 9 & 0 & 4 & 1 & 7 & 4 & 1 \\ 1 & 8 & 2 & 1 & 5 & 3 & 3 & 2 \\ 0 & 1 & 4 & 2 & 2 & 1 & 4 & 3 \end{bmatrix}.$$

We see $T_1$ and $T_3$ are simple matrices where as $T_2$ is a supermatrix. Thus T is not a super trimatrix or a trimatrix. So we define a new notion called semi super trimatrix which will accommodate both examples 3.14 and 3.15.

**DEFINITION 3.9:** *Let $T = T_1 \cup T_2 \cup T_3$, where some of the $T_i$'s are supermatrices and some of the $T_j$'s are ordinary matrices $1 \leq i, j \leq 3$. Then we call T to be a semi super trimatrix i.e., a semi super trimatrix T has at least one of the matrices $T_i$ to be a supermatrix, $1 \leq i \leq j \leq 3$.*

The matrices given in examples 3.14 and 3.15 are semi super trimatrices. Next consider the following examples.

*Example 3.16:* $T = T_1 \cup T_2 \cup T_3$

$$= \begin{bmatrix} 0 & 1 & | & 2 & 3 \\ 1 & 6 & | & 1 & 5 \\ \hline 2 & 1 & | & 2 & 3 \\ 3 & 5 & | & 3 & 4 \end{bmatrix} \cup \begin{bmatrix} 3 & 1 & | & 2 & 3 & | & 4 & 1 \\ 1 & 0 & | & 1 & 2 & | & 3 & 4 \\ \hline 2 & 1 & | & 5 & 6 & | & 7 & 8 \\ 3 & 2 & | & 6 & 9 & | & 8 & 7 \\ \hline 4 & 3 & | & 7 & 8 & | & 2 & 1 \\ 1 & 4 & | & 8 & 7 & | & 1 & 4 \end{bmatrix} \cup$$

$$\begin{bmatrix} 1 & 2 & | & 0 & 1 & | & 3 \\ 2 & 1 & | & 1 & 0 & | & 2 \\ \hline 0 & 1 & | & 5 & 7 & | & 9 \\ 1 & 0 & | & 7 & 1 & | & 2 \\ \hline 3 & 2 & | & 9 & 2 & | & 0 \end{bmatrix}$$

is a mixed square super trimatrix which is a symmetric super trimatrix.

Thus from this example we see if a super trimatrix is to be a symmetric super trimatrix then it should either be a square super trimatrix or a mixed square super trimatrix.

*Example 3.17:*

$$V = V_1 \cup V_2 \cup V_3$$

$$= \begin{bmatrix} 3 & | & 6 & 2 & 5 \\ \hline 6 & | & 1 & 0 & 1 \\ 2 & | & 0 & 9 & 2 \\ 5 & | & 1 & 2 & 7 \end{bmatrix} \cup \begin{bmatrix} 1 & 2 & | & 3 & 4 \\ 2 & 5 & | & 6 & 7 \\ \hline 3 & 6 & | & 8 & 9 \\ 4 & 7 & | & 9 & 0 \end{bmatrix} \cup \begin{bmatrix} 1 & 2 & 3 & | & 5 \\ 2 & 9 & 8 & | & 7 \\ 3 & 8 & 1 & | & 2 \\ \hline 5 & 7 & 2 & | & 9 \end{bmatrix};$$

V is a square super trimatrix. Further we see $V_1$, $V_2$ and $V_3$ are symmetrical supermatrices. Hence V is a square super trimatrix which is a symmetric super trimatrix.

**DEFINITION 3.10:** *Let $S = S_1 \cup S_2 \cup S_3$ be a square super trimatrix or a mixed square super trimatrix. If each of $S_i$ is a symmetric supermatrix then we call $S = S_1 \cup S_2 \cup S_3$ to be a symmetric trimatrix ($1 \leq i \leq 3$).*

The examples 3.16 and 3.17 are symmetric super trimatrices.

Now having defined the notion of symmetric super trimatrices we now proceed on to define quasi symmetric super trimatrices.

**DEFINITION 3.11:** *$V = V_1 \cup V_2 \cup V_3$ be a square super trimatrix or a mixed square super trimatrix or a mixed super trimatrix. We call V to be a quasi symmetric super trimatrix if at least one of the supermatrices $V_1$, $V_2$ or $V_3$ is a symmetric supermatrix.*

Now we illustrate this by some examples.

***Example 3.18:*** Let $T = T_1 \cup T_2 \cup T_3$

$$= \begin{bmatrix} 3 & 1 & 2 & 4 & 5 \\ 1 & 0 & 1 & 2 & 1 \\ \hline 0 & 1 & 3 & 0 & 9 \\ 3 & 1 & 5 & 1 & 8 \end{bmatrix} \cup \begin{bmatrix} 3 & 4 & 5 & 6 & 1 \\ 4 & 1 & 2 & 3 & 4 \\ \hline 5 & 2 & 8 & 9 & 1 \\ 6 & 3 & 9 & 7 & 2 \\ 1 & 4 & 1 & 2 & 6 \end{bmatrix} \cup$$

$$\begin{bmatrix} 6 & 0 & 9 & 2 \\ 6 & 2 & 1 & 6 \\ \hline 1 & 1 & 0 & 9 \\ 3 & 3 & 9 & 2 \\ 4 & 4 & 8 & 0 \\ \hline 6 & 3 & 2 & 1 \end{bmatrix};$$

T is super trimatrix which is a quasi symmetric super trimatrix as $T_2$ is a symmetric supermatrix but $T_1$ and $T_3$ are just supermatrices.

***Example 3.19:*** Let $V = V_1 \cup V_2 \cup V_3$ where

$$V_1 = \begin{bmatrix} 3 & 9 & 6 & | & 4 & 5 \\ 9 & 2 & 0 & | & 1 & 2 \\ 6 & 0 & 7 & | & 6 & 5 \\ \hline 4 & 1 & 6 & | & 9 & 2 \\ 5 & 2 & 5 & | & 2 & 8 \end{bmatrix},$$

$$V_2 = \begin{bmatrix} 0 & 1 & | & 2 & 3 & 4 & 5 & | & 6 \\ 1 & 4 & | & 2 & 3 & 4 & 5 & | & 0 \\ \hline 2 & 2 & | & 6 & 5 & 4 & 3 & | & 2 \\ 3 & 3 & | & 5 & 7 & 0 & 1 & | & 2 \\ 4 & 4 & | & 4 & 0 & 9 & 0 & | & 1 \\ 5 & 5 & | & 3 & 1 & 0 & 1 & | & 3 \\ \hline 6 & 0 & | & 2 & 2 & 1 & 3 & | & 7 \end{bmatrix}$$

and

$$V_3 = \begin{bmatrix} 3 & 1 & | & 0 & 9 & 6 & 4 & | & 3 & 2 & 9 \\ 1 & 2 & | & 1 & 8 & 7 & 6 & | & 4 & 5 & 3 \\ \hline 2 & 2 & | & 5 & 1 & 4 & 0 & | & 1 & 1 & 2 \\ 5 & 1 & | & 2 & 3 & 4 & 5 & | & 6 & 7 & 8 \\ 7 & 0 & | & 1 & 0 & 1 & 4 & | & 0 & 5 & 2 \end{bmatrix}.$$

Clearly V is a super trimatrix but V is only a quasi symmetric super trimatrix as $V_1$ and $V_2$ are symmetric supermatrices but $V_3$ is only a supermatrix.

Now having seen examples of quasi super trimatrices we proceed on to define the notion of quasi semi super trimatrices.

**DEFINITION 3.12:** *Let $V = V_1 \cup V_2 \cup V_3$ be a semi super trimatrix. If $V_1$, $V_2$ and $V_3$ are square matrices or super square matrices then we call V to be a mixed square semi super* trimatrix. If $V_1$, $V_2$ and $V_3$ are $n \times n$ square matrices or $n \times n$ supermatrices then we call V to be a square semi super trimatrix or an $n \times n$ semi super trimatrix.

**Example 3.20:** Let

$$T = T_1 \cup T_2 \cup T_3$$

$$= \begin{bmatrix} 3 & 1 \\ \hline 2 & 0 \end{bmatrix} \cup \begin{bmatrix} 3 & 0 \\ 1 & 2 \end{bmatrix} \cup \begin{bmatrix} 1 & 2 \\ 3 & 5 \end{bmatrix}$$

be a semi super trimatrix. Clearly T is a $2 \times 2$ or square semi super trimatrix.

**Example 3.21:** Let $V = V_1 \cup V_2 \cup V_3$ be a semi super trimatrix where

$$V_1 = \begin{bmatrix} 3 & 2 & 1 & 5 \\ 6 & 7 & 8 & 9 \\ 7 & 5 & 4 & 2 \\ 0 & 1 & 3 & 2 \end{bmatrix}.$$

$$V_2 = \begin{bmatrix} 3 & 1 & 0 & 1 & 2 & 5 \\ \hline 1 & 1 & 2 & 0 & 1 & 1 \\ 0 & 2 & 9 & 2 & 0 & 6 \\ 2 & 1 & 3 & 2 & 9 & 1 \\ 5 & 3 & 1 & 0 & 5 & 6 \\ \hline 1 & 2 & 3 & 4 & 2 & 1 \end{bmatrix}$$

and

$$V_3 = \begin{bmatrix} 3 & 1 & 2 \\ \hline 0 & 1 & 1 \\ \hline 1 & 0 & 1 \end{bmatrix}.$$

Clearly V is a mixed square semi super trimatrix.

Next we define mixed super trimatrix.

**DEFINITION 3.13:** *Let $T = T_1 \cup T_2 \cup T_3$ where some of $T_i$ is a square supermatrix or matrix and the rest are rectangular supermatrix or matrix then we call T to be a mixed super trimatrix.*

*Example 3.22:* Let $V = V_1 \cup V_2 \cup V_3$ where

$$V_1 = \begin{bmatrix} 3 & 1 & 4 & 5 \\ 7 & 0 & 2 & 3 \\ 8 & 4 & 8 & 2 \\ 9 & 7 & 1 & 0 \end{bmatrix},$$

$$V_2 = \left[\begin{array}{cccc|ccc|cc} 3 & 1 & 2 & 5 & 7 & 8 & 9 & 1 & 2 \\ 0 & 1 & 2 & 3 & 0 & 2 & 3 & 1 & 0 \\ 1 & 2 & 0 & 5 & 2 & 0 & 1 & 3 & 1 \\ \hline 1 & 1 & 1 & 0 & 3 & 1 & 0 & 1 & 5 \end{array}\right]$$

and

$$V_3 = \left[\begin{array}{cc|ccc|c} 3 & 4 & 5 & 6 & 7 & 8 \\ 4 & 2 & 1 & 0 & 5 & 6 \\ \hline 5 & 6 & 1 & 6 & 2 & 3 \\ 6 & 3 & 0 & 1 & 2 & 0 \\ 7 & 1 & 3 & 1 & 4 & 5 \\ \hline 8 & 2 & 0 & 5 & 2 & 1 \end{array}\right];$$

V is a semi super trimatrix which is a mixed semi super trimatrix.

*Example 3.23:* Let $S = S_1 \cup S_2 \cup S_3$ where

$$S_1 = \begin{bmatrix} 1 & 2 & 3 & 4 & 5 & 6 & 7 & 8 \\ 9 & 1 & 2 & 3 & 4 & 5 & 6 & 7 \\ 8 & 9 & 1 & 2 & 3 & 4 & 5 & 6 \end{bmatrix},$$

$$S_2 = \left[\begin{array}{cc|cc} 3 & 1 & 2 & 1 \\ 1 & 1 & 0 & 2 \\ 3 & 9 & 7 & 5 \\ \hline 1 & 8 & 6 & 4 \end{array}\right]$$

and

$$S_3 = \left[\begin{array}{cc|cc} 3 & 1 & 2 & 5 \\ \hline 1 & 1 & 0 & 0 \\ 9 & 9 & 2 & 2 \\ 2 & 1 & 3 & 5 \\ \hline 1 & 2 & 3 & 4 \\ 5 & 6 & 7 & 8 \\ 9 & 1 & 2 & 3 \end{array}\right].$$

S is a semi super trimatrix which is a mixed semi super trimatrix.

*Example 3.24:* Let $V = V_1 \cup V_2 \cup V_3$ be a semi super trimatrix where

$$V_1 = \begin{bmatrix} 1 & 3 & 5 & 7 & 9 & 2 & 4 & 6 & 8 \\ 2 & 4 & 6 & 8 & 1 & 3 & 5 & 7 & 9 \end{bmatrix},$$

$$V_2 = \left[\begin{array}{ccc|c} 3 & 1 & 2 & 4 \\ 0 & 1 & 0 & 2 \\ \hline 5 & 7 & 6 & 2 \\ 9 & 0 & 1 & 8 \\ 3 & 2 & 1 & 0 \\ 1 & 1 & 1 & 5 \end{array}\right]$$

and

$$V_3 = \begin{bmatrix} 1 & 2 & 3 & | & 4 & 5 \\ 6 & 3 & 1 & | & 2 & 1 \\ 7 & 0 & 1 & | & 1 & 1 \\ 8 & 5 & 6 & | & 0 & 1 \\ \hline 9 & 7 & 8 & | & 9 & 2 \\ 1 & 3 & 4 & | & 5 & 6 \end{bmatrix}.$$

We see the three matrices $V_1$, $V_2$ and $V_3$ are rectangular matrices of order $2 \times 9$, $6 \times 4$ and $6 \times 5$ respectively. V is called or defined as a mixed rectangular semi super trimatrix.

*Example 3.25:* Let $T = T_1 \cup T_2 \cup T_3$ be a semi super trimatrix where

$$T_1 = \begin{bmatrix} 3 & 1 & 2 & 6 & 7 & 8 & 1 \\ 0 & 1 & 2 & 3 & 4 & 5 & 6 \\ 7 & 8 & 9 & 1 & 2 & 3 & 0 \end{bmatrix},$$

$$T_2 = \begin{bmatrix} 1 & 2 & 3 & 4 & 5 & 6 & | & 7 \\ \hline 3 & 2 & 1 & 1 & 4 & 3 & | & 2 \\ 5 & 7 & 8 & 0 & 1 & 2 & | & 3 \end{bmatrix}$$

and

$$T_3 = \begin{bmatrix} 0 & | & 1 & 2 & | & 3 & 4 & 5 & 6 \\ 7 & | & 8 & 9 & | & 1 & 1 & 2 & 0 \\ 3 & | & 0 & 1 & | & 0 & 1 & 2 & 5 \end{bmatrix}.$$

We see each of the 3 matrices $T_1$, $T_2$ and $T_3$ are $3 \times 7$ matrices or $3 \times 7$ supermatrices. Thus T is a $3 \times 7$ rectangular semi super trimatrix.

*Example 3.26:* Let $T = T_1 \cup T_2 \cup T_3$ be a semi super trimatrix where

$$T_1 = \begin{bmatrix} 3 & 1 & 0 & 2 & 1 \\ 1 & 1 & 2 & 5 & 3 \end{bmatrix},$$

$$T_2 = \begin{bmatrix} 7 & 8 & 9 \\ 1 & 2 & 3 \\ 4 & 5 & 6 \\ 3 & 1 & 2 \\ 5 & 6 & 4 \\ 8 & 9 & 7 \end{bmatrix}$$

and

$$T_3 = \left[ \begin{array}{ccc|ccc|c} 3 & 1 & 2 & 3 & 4 & 5 & 1 \\ 2 & 0 & 1 & 0 & 1 & 1 & 2 \\ \hline 1 & 2 & 4 & 6 & 8 & 1 & 3 \\ 0 & 3 & 5 & 7 & 9 & 0 & 5 \end{array} \right].$$

T is a mixed rectangular semi super trimatrix.

Now we see just a column semi super trimatrix and row semi super trimatrix.

*Example 3.27:* Let $A = A_1 \cup A_2 \cup A_3$ where

$$A_1 = [1\ 2\ 3\ 4\ 5\ 6],$$
$$A_2 = [0\ 1\ 0\ 1\ |\ 0\ 1\ |\ 0\ 1\ 3]$$

and

$$A_3 = [9\ 8\ 7\ |\ 6\ 4\ 3\ 2\ |\ 1\ 0].$$

A is a semi super trimatrix which we call as simple row semi super trimatrix or simple row semi super trivector.

*Example 3.28:* Let $V = V_1 \cup V_2 \cup V_3$ where

$$V_1 = \begin{bmatrix} 1 \\ 2 \\ 3 \\ 4 \\ 5 \\ 6 \\ 7 \\ 8 \end{bmatrix}, \ V_2 = \begin{bmatrix} 3 \\ 1 \\ 2 \\ 5 \\ 6 \\ 1 \end{bmatrix} \text{ and } V_3 = \begin{bmatrix} 3 \\ 1 \\ 2 \\ 3 \\ 4 \\ 0 \\ 1 \\ 2 \end{bmatrix}.$$

V is a semi super trimatrix which is a simple column semi super trivector or a simple column semi super trimatrix. We define dual partition in case of super trimatrices. It is important to mention here the main difference between a simple row or column super trimatrix and a row or column super trivector is that when we say simple column or row super trimatrix each of its components in $T = T_1 \cup T_2 \cup T_3$ are just a simple $1 \times n_i$, $i = 1, 2, 3$ row matrices or simple $m_i \times 1$ column matrices.

We first see how the product is defined.

**DEFINITION 3.14:** *Let $T = T_1 \cup T_2 \cup T_3$ be a simple row semi super trimatrix, where $T_i = (a_1^i, \ldots, a_{n_i}^i)$; $1 \leq i \leq 3$. If $T_i$ is partitioned between rows r and r + 1, s and s + 1 and t and t + 1 and so on. We know*

$$T^T = (T_1 \cup T_2 \cup T_3)^T = T_1^T \cup T_2^T \cup T_3^T$$

*where*

$$T_i^T = \begin{bmatrix} a_1^i \\ a_2^i \\ \vdots \\ a_{n_i}^i \end{bmatrix};$$

*$i = 1, 2, 3$. Now each $T_i^T$ will be partitioned between the columns r and r + 1, s and s + 1 and t and t + 1 and so on. We call this type of partition carried on say from a $1 \times n_i$ row to any* $n_i \times 1$ *column vector to be a dual partition. Thus if a row supermatrix* $A = (a_1 \, a_2 \, \ldots \, a_n)$ *has some partition and if a column supermatrix*

$$B = \begin{pmatrix} b_1 \\ b_2 \\ \vdots \\ b_n \end{pmatrix}$$

*has the dual partition then we have a well defined product AB.*

This we will first illustrate by an example.

***Example 3.29:*** Let $T = T_1 \cup T_2 \cup T_3$ where $A_1 = [2\ 3\ 0\ |\ 1\ 5\ 4\ 7\ |\ 2\ 0]$, $A_2 = [1\ |\ 2\ 3\ 4\ 5\ |\ 0\ 1]$ and $A_3 = [1\ 1\ 1\ 0\ |\ 1\ 2\ 0\ 1\ 5\ |\ 0\ 1\ ]$ be a simple row super trimatrix. Let $B = B_1 \cup B_2 \cup B_3$ where

$$B_1 = \begin{bmatrix} 1 \\ 0 \\ 1 \\ \overline{2} \\ 2 \\ 0 \\ 1 \\ \overline{3} \\ 1 \end{bmatrix} \quad B_2 = \begin{bmatrix} 0 \\ \overline{2} \\ 0 \\ 1 \\ 5 \\ \overline{1} \\ 2 \end{bmatrix} \text{ and } B_3 = \begin{bmatrix} 3 \\ 1 \\ 2 \\ 0 \\ 1 \\ 1 \\ 0 \\ 1 \\ \overline{2} \\ 5 \\ 0 \end{bmatrix}$$

be a simple column super trimatrix. Then

AB     =    $(A_1 \cup A_2 \cup A_3) (B_1 \cup B_2 \cup B_3)$
         =    $A_1 B_1 \cup A_2 B_2 \cup A_3 B_3$

$$= \left\{ \begin{bmatrix} 2 & 3 & 0 \mid 1 & 5 & 4 & 7 \mid 2 & 0 \end{bmatrix} \begin{bmatrix} 1 \\ 0 \\ 1 \\ 2 \\ 2 \\ 0 \\ 1 \\ 3 \\ 1 \end{bmatrix} \right\} \cup$$

$$\left\{ \begin{bmatrix} 1 \mid 2 & 3 & 4 & 5 \mid 0 & 1 \end{bmatrix} \begin{bmatrix} 0 \\ 2 \\ 0 \\ 1 \\ 5 \\ 1 \\ 2 \end{bmatrix} \right\} \cup$$

$$\left\{ \begin{bmatrix} 1 & 1 & 1 & 0 \mid 1 & 2 & 0 & 1 & 5 \mid 0 & 1 \end{bmatrix} \begin{bmatrix} 3 \\ 1 \\ 2 \\ 0 \\ 1 \\ 1 \\ 0 \\ 1 \\ 2 \\ 5 \\ 0 \end{bmatrix} \right\}$$

$$= \left\{ [2 \ 3 \ 0] \begin{bmatrix} 1 \\ 0 \\ 1 \end{bmatrix} \middle| [1 \ 5 \ 4 \ 7] \begin{bmatrix} 2 \\ 2 \\ 0 \\ 1 \end{bmatrix} \middle| [2 \ 0] \begin{bmatrix} 3 \\ 1 \end{bmatrix} \right\} \cup$$

$$\left\{ [1][0] \middle| [2 \ 3 \ 4 \ 5] \begin{bmatrix} 2 \\ 0 \\ 1 \\ 5 \end{bmatrix} \middle| [0 \ 1] \begin{bmatrix} 1 \\ 2 \end{bmatrix} \right\} \cup$$

$$\left\{ [1 \ 1 \ 1 \ 0] \begin{bmatrix} 3 \\ 1 \\ 2 \\ 0 \end{bmatrix} \middle| [1 \ 2 \ 0 \ 1 \ 5] \begin{bmatrix} 1 \\ 1 \\ 0 \\ 1 \\ 2 \end{bmatrix} \middle| [0 \ 1] \begin{bmatrix} 5 \\ 0 \end{bmatrix} \right\}.$$

$$= [2 \mid 19 \mid 6] \cup [0 \mid 33 \mid 2] \cup [6 \mid 14 \mid 0].$$

Clearly this is a simple row super trimatrix.

Suppose we want to find

BA  =  $(B_1 \cup B_2 \cup B_3)(A_1 \cup A_2 \cup A_3)$
   =  $B_1A_1 \cup B_2A_2 \cup B_3A_3$.

$$= \left\{ \begin{bmatrix} 1 \\ 0 \\ \frac{1}{2} \\ 2 \\ 0 \\ \frac{1}{3} \\ 1 \end{bmatrix} \begin{bmatrix} 2 & 3 & 0 \mid 1 & 5 & 4 & 7 \mid 2 & 0 \end{bmatrix} \right\} \cup$$

$$\left\{ \begin{bmatrix} 0 \\ \frac{2}{0} \\ 1 \\ \frac{5}{1} \\ 2 \end{bmatrix} \begin{bmatrix} 1 \mid 2 & 3 & 4 & 5 \mid 0 & 1 \end{bmatrix} \right\} \cup$$

$$\left\{ \begin{bmatrix} 3 \\ 1 \\ 2 \\ \frac{0}{1} \\ 1 \\ 0 \\ 1 \\ \frac{2}{5} \\ 0 \end{bmatrix} \begin{bmatrix} 1 & 1 & 1 & 0 \mid 1 & 2 & 0 & 1 & 5 \mid 0 & 1 \end{bmatrix} \right\}$$

$$= \left\{ \begin{bmatrix} 1 \\ 0 \\ 1 \end{bmatrix} [2 \quad 3 \quad 0] \middle| \begin{bmatrix} 2 \\ 2 \\ 0 \\ 1 \end{bmatrix} [15 \quad 4 \quad 7] \middle| \begin{bmatrix} 3 \\ 1 \end{bmatrix} [2 \quad 0] \right\} \cup$$

$$\left\{ [0][1] \middle| \begin{bmatrix} 2 \\ 0 \\ 1 \\ 5 \end{bmatrix} [2 \quad 3 \quad 4 \quad 5] \middle| \begin{bmatrix} 1 \\ 2 \end{bmatrix} [0 \quad 1] \right\}$$

$$\left\{ \begin{bmatrix} 3 \\ 1 \\ 2 \\ 0 \end{bmatrix} [1 \quad 1 \quad 1 \quad 0] \middle| \begin{bmatrix} 1 \\ 1 \\ 0 \\ 1 \\ 2 \end{bmatrix} [1 \quad 2 \quad 0 \quad 1 \quad 5] \middle| \begin{bmatrix} 5 \\ 0 \end{bmatrix} [0 \quad 1] \right\}$$

$$= \begin{bmatrix} 2 & 3 & 0 & 2 & 10 & 8 & 14 & 6 & 0 \\ 0 & 0 & 0 & 2 & 10 & 8 & 14 & 2 & 0 \\ 2 & 3 & 0 & 0 & 0 & 0 & 0 & & \\ & & & 1 & 5 & 4 & 7 & & \end{bmatrix} \cup$$

$$\cdot \begin{bmatrix} 0 & 4 & 6 & 8 & 10 & 0 & 1 \\ & 0 & 0 & 0 & 0 & & \\ & 2 & 3 & 4 & 5 & 0 & 2 \\ & 10 & 15 & 20 & 25 & & \end{bmatrix} \cup$$

$$\begin{bmatrix} 3 & 3 & 3 & 0 & 1 & 2 & 0 & 1 & 5 & & \\ 1 & 1 & 1 & 0 & 1 & 2 & 0 & 1 & 5 & 0 & 5 \\ 2 & 2 & 2 & 0 & 0 & 0 & 0 & 0 & 0 & & \\ 0 & 0 & 0 & 0 & 1 & 2 & 0 & 1 & 5 & 0 & 5 \\ & & & & 2 & 4 & 0 & 2 & 10 & & \end{bmatrix}.$$

We see though we have the compatibility with respect to each cell, yet the resultant is some new structure which is never defined for they are not super trimatrices or vector or any known mathematical structure.

Thus as in case of matrices if AB is defined it may happen that BA is undefined likewise we see here AB is defined but BA is undefined.

Now we illustrate by a simple example the minor product of two super trimatrices.

***Example 3.30:*** Let $T = T_1 \cup T_2 \cup T_3$ and $S = S_1 \cup S_2 \cup S_3$ be any two super trivectors, where

$$T = T_1 \cup T_2 \cup T_3$$

$$= [3 \ 0 \mid 5 \ 1 \ 0 \ 2 \mid 3] \cup \begin{bmatrix} 0 & 1 & 0 & 1 & 2 & 3 \\ 1 & 2 & 6 & 0 & 4 & 5 \end{bmatrix} \cup$$

$$\begin{bmatrix} 1 & 2 & 3 & 2 & 0 & 1 & 3 & 2 & 1 \\ 1 & 1 & 0 & 1 & 5 & 0 & 1 & 0 & 1 \\ 1 & 1 & 5 & 0 & 4 & 0 & 7 & 3 & 0 \end{bmatrix}$$

be a row super trivector and
$$S = S_1 \cup S_2 \cup S_3$$

$$= \begin{bmatrix} 2 & 1 \\ 1 & 0 \\ \hline 1 & 2 \\ 2 & 1 \\ 3 & 4 \\ 4 & 3 \\ \hline 1 & 0 \end{bmatrix} \cup \begin{bmatrix} 1 & 2 & 3 & 4 \\ 0 & 1 & 2 & 5 \\ \hline 1 & 3 & 0 & 1 \\ 1 & 1 & 0 & 2 \\ 2 & 0 & 2 & 1 \\ 5 & 1 & 0 & 2 \end{bmatrix} \cup \begin{bmatrix} 1 & 0 & 1 & 1 & 1 \\ 0 & 1 & 0 & 1 & 0 \\ 2 & 0 & 1 & 2 & 0 \\ \hline 1 & 0 & 1 & 0 & 0 \\ 0 & 0 & 0 & 1 & 0 \\ 1 & 0 & 0 & 0 & 1 \\ 0 & 1 & 1 & 1 & 0 \\ \hline 3 & 1 & 0 & 0 & 1 \\ 1 & 0 & 1 & 0 & 1 \end{bmatrix}$$

be a column super trivector. Now

$$TS = (T_1 \cup T_2 \cup T_3)(S_1 \cup S_2 \cup S_3)$$
$$= T_1S_1 \cup T_2S_2 \cup T_3S_3$$

$$= \begin{bmatrix} 3 & 0 & | & 5 & 1 & 0 & 2 & | & 3 \end{bmatrix} \begin{bmatrix} 2 & 1 \\ 1 & 0 \\ \hline 1 & 2 \\ 2 & 1 \\ 3 & 4 \\ 4 & 3 \\ \hline 1 & 0 \end{bmatrix} \cup$$

$$\begin{bmatrix} 1 & | & 1 & 0 & | & 1 & 2 & 3 \\ 0 & | & 2 & 6 & | & 0 & 4 & 5 \end{bmatrix} \begin{bmatrix} 1 & 2 & 3 & 4 \\ \hline 0 & 1 & 2 & 5 \\ 1 & 3 & 0 & 1 \\ \hline 1 & 1 & 0 & 2 \\ 2 & 0 & 2 & 1 \\ 5 & 1 & 0 & 2 \end{bmatrix} \cup$$

$$\begin{bmatrix} 1 & 2 & 3 & | & 2 & 0 & 1 & 3 & | & 2 & 1 \\ 1 & 1 & 0 & | & 1 & 5 & 0 & 1 & | & 0 & 1 \\ 1 & 1 & 5 & | & 0 & 4 & 0 & 7 & | & 3 & 0 \end{bmatrix} \begin{bmatrix} 1 & 0 & 1 & 1 & 1 \\ 0 & 1 & 0 & 1 & 0 \\ 2 & 0 & 1 & 2 & 0 \\ \hline 1 & 0 & 1 & 0 & 0 \\ 0 & 0 & 0 & 1 & 0 \\ 1 & 0 & 0 & 0 & 1 \\ 0 & 1 & 1 & 1 & 0 \\ \hline 3 & 1 & 0 & 0 & 1 \\ 1 & 0 & 1 & 0 & 1 \end{bmatrix}.$$

$$= \left\{ \begin{bmatrix} 3 & 0 \end{bmatrix} \begin{bmatrix} 2 & 1 \\ 1 & 0 \end{bmatrix} + \begin{bmatrix} 5 & 1 & 0 & 2 \end{bmatrix} \begin{bmatrix} 1 & 2 \\ 2 & 1 \\ 3 & 4 \\ 4 & 3 \end{bmatrix} + 3\begin{bmatrix} 1 & 0 \end{bmatrix} \right\} \cup$$

$$\left\{ \begin{bmatrix} 1 \\ 0 \end{bmatrix} \begin{bmatrix} 1 & 2 & 3 & 4 \end{bmatrix} + \begin{bmatrix} 1 & 0 \\ 2 & 6 \end{bmatrix} \begin{bmatrix} 0 & 1 & 2 & 5 \\ 1 & 3 & 0 & 1 \end{bmatrix} + \right.$$

$$\left. \begin{bmatrix} 1 & 2 & 3 \\ 0 & 4 & 5 \end{bmatrix} \begin{bmatrix} 1 & 1 & 0 & 2 \\ 2 & 0 & 2 & 1 \\ 5 & 1 & 0 & 2 \end{bmatrix} \right\} \cup$$

$$\left\{ \begin{bmatrix} 1 & 2 & 3 \\ 1 & 1 & 0 \\ 1 & 1 & 5 \end{bmatrix} \begin{bmatrix} 1 & 0 & 1 & 1 & 1 \\ 0 & 1 & 0 & 1 & 0 \\ 2 & 0 & 1 & 2 & 0 \end{bmatrix} + \right.$$

$$\left. \begin{bmatrix} 2 & 0 & 1 & 3 \\ 1 & 5 & 0 & 1 \\ 0 & 4 & 0 & 7 \end{bmatrix} \begin{bmatrix} 1 & 0 & 1 & 0 & 0 \\ 0 & 0 & 0 & 1 & 0 \\ 1 & 0 & 0 & 0 & 1 \\ 0 & 1 & 1 & 1 & 0 \end{bmatrix} + \begin{bmatrix} 2 & 1 \\ 0 & 1 \\ 3 & 0 \end{bmatrix} \begin{bmatrix} 3 & 1 & 0 & 0 & 11 \\ 1 & 0 & 1 & 0 & 1 \end{bmatrix} \right\}$$

$$= \left\{ \begin{bmatrix} 6 & 3 \end{bmatrix} + \begin{bmatrix} 15 & 17 \end{bmatrix} + \begin{bmatrix} 3 & 0 \end{bmatrix} \right\} \cup$$

$$\left\{ \begin{bmatrix} 1 & 2 & 3 & 4 \\ 0 & 0 & 0 & 0 \end{bmatrix} + \begin{bmatrix} 0 & 1 & 2 & 5 \\ 6 & 20 & 4 & 16 \end{bmatrix} \right.$$

$$\left. + \begin{bmatrix} 20 & 4 & 4 & 10 \\ 33 & 5 & 8 & 14 \end{bmatrix} \right\} \cup$$

$$\left\{ \begin{bmatrix} 7 & 2 & 4 & 9 & 1 \\ 1 & 1 & 1 & 2 & 1 \\ 11 & 1 & 6 & 12 & 1 \end{bmatrix} + \begin{bmatrix} 3 & 3 & 5 & 3 & 1 \\ 1 & 1 & 2 & 6 & 0 \\ 0 & 7 & 7 & 11 & 0 \end{bmatrix} + \begin{bmatrix} 7 & 2 & 1 & 0 & 3 \\ 1 & 0 & 1 & 0 & 1 \\ 9 & 3 & 0 & 0 & 3 \end{bmatrix} \right\}$$

$$= \begin{bmatrix} 24 & 20 \end{bmatrix} \cup \begin{bmatrix} 21 & 7 & 9 & 19 \\ 39 & 25 & 12 & 30 \end{bmatrix} \cup \begin{bmatrix} 17 & 7 & 10 & 12 & 5 \\ 3 & 2 & 4 & 8 & 2 \\ 20 & 11 & 13 & 23 & 4 \end{bmatrix}.$$

We see the resultant is only a trimatrix. Thus the product TS of the row super trivector T with a compatible column super trivector S is only a trimatrix and not a super trimatrix.

*Example 3.31:* Let $T = T_1 \cup T_2 \cup T_3$ be a row super trivector and $S = S_1 \cup S_2 \cup S_3$ be a column super trivector.

$$\begin{aligned} TS &= [T_1 \cup T_2 \cup T_3][S_1 \cup S_2 \cup S_3] \\ &= T_1 S_1 \cup T_2 S_2 \cup T_3 S_3 \end{aligned}$$

where

$$T = \left[ \begin{array}{ccc|cccc|cc} 3 & 5 & 1 & 3 & 1 & 2 & 4 & 2 & 0 \\ 1 & 0 & 6 & 1 & 0 & 1 & 0 & 1 & 1 \end{array} \right] \cup$$

$$\left[ \begin{array}{cccc|cc|ccc} 3 & 0 & 1 & 1 & 1 & 1 & 0 & 1 & 0 \\ 1 & 6 & 0 & 0 & 0 & 0 & 1 & 0 & 1 \\ 0 & 1 & 0 & 1 & 2 & 1 & 1 & 1 & 0 \end{array} \right] \cup$$

$$\left[ \begin{array}{c|ccc|cccc} 1 & 3 & 3 & 1 & 0 & 3 & 1 & 2 & 1 \\ 2 & 4 & 1 & 0 & 1 & 0 & 1 & 0 & 0 \\ 1 & 0 & 1 & 1 & 0 & 1 & 0 & 1 & 0 \\ 0 & 1 & 2 & 2 & 1 & 0 & 1 & 0 & 1 \end{array} \right]$$

and

$$S = \begin{bmatrix} 3 & 1 \\ 1 & 0 \\ 1 & 3 \\ \hline 1 & 1 \\ 0 & 1 \\ 2 & 1 \\ 4 & 5 \\ \hline 1 & 0 \\ 1 & 4 \end{bmatrix} \cup \begin{bmatrix} 1 \\ 0 \\ 1 \\ \hline 1 \\ 5 \\ 0 \\ \hline 0 \\ 2 \\ 3 \end{bmatrix} \cup \begin{bmatrix} 3 & 0 & 1 \\ 0 & 1 & 0 \\ 1 & 1 & 2 \\ 0 & 1 & 1 \\ 1 & 1 & 0 \\ 1 & 0 & 1 \\ 0 & 1 & 0 \\ 0 & 0 & 0 \\ 1 & 1 & 1 \end{bmatrix}.$$

$$TS = \begin{bmatrix} 3 & 5 & 1 & | & 3 & 1 & 2 & 4 & | & 2 & 0 \\ 1 & 0 & 6 & | & 1 & 0 & 1 & 0 & | & 1 & 1 \end{bmatrix} \begin{bmatrix} 3 & 1 \\ 1 & 0 \\ 1 & 3 \\ \hline 1 & 1 \\ 0 & 1 \\ 2 & 1 \\ 4 & 5 \\ \hline 1 & 0 \\ 4 & 1 \end{bmatrix} \cup$$

$$\begin{bmatrix} 3 & 0 & 1 & 1 & | & 1 & 1 & | & 0 & 1 & 0 \\ 1 & 6 & 0 & 0 & | & 0 & 0 & | & 1 & 0 & 1 \\ 0 & 1 & 0 & 1 & | & 2 & 1 & | & 1 & 1 & 0 \end{bmatrix} \begin{bmatrix} 1 \\ 0 \\ 1 \\ 1 \\ \hline 5 \\ 0 \\ \hline 0 \\ 2 \\ 3 \end{bmatrix} \cup$$

$$\begin{bmatrix} 1 & 3 & 3 & 1 & 0 & 3 & 1 & 2 & 1 \\ 2 & 4 & 1 & 0 & 0 & 1 & 1 & 0 & 0 \\ 1 & 0 & 1 & 1 & 0 & 1 & 0 & 1 & 0 \\ 0 & 1 & 2 & 2 & 1 & 0 & 1 & 0 & 1 \end{bmatrix} \begin{bmatrix} 3 & 0 & 1 \\ 0 & 1 & 0 \\ 1 & 1 & 2 \\ 0 & 1 & 1 \\ \hline 1 & 1 & 0 \\ 1 & 0 & 1 \\ 0 & 1 & 0 \\ 0 & 0 & 6 \\ 1 & 1 & 1 \end{bmatrix}$$

$$= \left\{ \begin{bmatrix} 3 & 5 & 1 \\ 1 & 0 & 6 \end{bmatrix} \begin{bmatrix} 3 & 1 \\ 1 & 0 \\ 1 & 3 \end{bmatrix} + \begin{bmatrix} 3 & 1 & 2 & 4 \\ 1 & 0 & 1 & 0 \end{bmatrix} \begin{bmatrix} 1 & 1 \\ 0 & 1 \\ 2 & 1 \\ 4 & 5 \end{bmatrix} \right.$$

$$\left. + \begin{bmatrix} 2 & 0 \\ 1 & 1 \end{bmatrix} \begin{bmatrix} 1 & 0 \\ 1 & 4 \end{bmatrix} \right\} \cup$$

$$\left\{ \begin{bmatrix} 3 & 0 & 1 & 1 \\ 1 & 6 & 0 & 0 \\ 0 & 1 & 0 & 1 \end{bmatrix} \begin{bmatrix} 1 \\ 0 \\ 1 \\ 1 \end{bmatrix} + \begin{bmatrix} 1 & 1 \\ 0 & 0 \\ 2 & 1 \end{bmatrix} \begin{bmatrix} 5 \\ 0 \end{bmatrix} + \begin{bmatrix} 0 & 1 & 0 \\ 1 & 0 & 1 \\ 1 & 1 & 0 \end{bmatrix} \begin{bmatrix} 0 \\ 2 \\ 3 \end{bmatrix} \right\}$$

$$\cup \left\{ \begin{bmatrix} 1 \\ 2 \\ 1 \\ 0 \end{bmatrix} \begin{bmatrix} 3 & 0 & 1 \end{bmatrix} + \begin{bmatrix} 3 & 3 & 1 \\ 4 & 1 & 0 \\ 0 & 1 & 1 \\ 1 & 2 & 2 \end{bmatrix} \begin{bmatrix} 0 & 1 & 0 \\ 1 & 1 & 2 \\ 0 & 1 & 1 \end{bmatrix} + \right.$$

$$\left. \begin{bmatrix} 0 & 3 & 1 & 2 & 1 \\ 0 & 1 & 1 & 0 & 0 \\ 0 & 1 & 0 & 1 & 0 \\ 1 & 0 & 1 & 0 & 1 \end{bmatrix} \begin{bmatrix} 1 & 1 & 0 \\ 1 & 0 & 1 \\ 0 & 1 & 0 \\ 0 & 0 & 6 \\ 1 & 1 & 1 \end{bmatrix} \right\} =$$

$$\left\{ \begin{bmatrix} 15 & 6 \\ 9 & 19 \end{bmatrix} + \begin{bmatrix} 23 & 26 \\ 3 & 2 \end{bmatrix} + \begin{bmatrix} 2 & 0 \\ 2 & 4 \end{bmatrix} \right\} \cup \left\{ \begin{bmatrix} 5 \\ 1 \\ 1 \end{bmatrix} + \begin{bmatrix} 5 \\ 0 \\ 10 \end{bmatrix} + \begin{bmatrix} 2 \\ 3 \\ 2 \end{bmatrix} \right\} \cup$$

$$\left\{ \begin{bmatrix} 3 & 0 & 1 \\ 6 & 0 & 2 \\ 3 & 0 & 1 \\ 0 & 0 & 0 \end{bmatrix} + \begin{bmatrix} 3 & 7 & 7 \\ 1 & 5 & 2 \\ 1 & 2 & 3 \\ 2 & 5 & 6 \end{bmatrix} + \begin{bmatrix} 4 & 2 & 16 \\ 1 & 1 & 1 \\ 1 & 0 & 7 \\ 2 & 3 & 1 \end{bmatrix} \right\}$$

$$= \begin{bmatrix} 40 & 32 \\ 14 & 25 \end{bmatrix} \cup \begin{bmatrix} 12 \\ 4 \\ 13 \end{bmatrix} \cup \begin{bmatrix} 10 & 9 & 24 \\ 8 & 6 & 5 \\ 5 & 2 & 11 \\ 4 & 8 & 7 \end{bmatrix}.$$

We see TS is just a trimatrix which is not a super trimatrix.

Now we find the product of a row super trivector T with its transpose.

*Example 3.32:* Let $T = T_1 \cup T_2 \cup T_3$

$$= \begin{bmatrix} 3 & 0 & 4 & | & 7 & 2 & 3 & 1 & 5 & | & 1 & 2 \\ 1 & 2 & 5 & | & 0 & 1 & 0 & 1 & 0 & | & 2 & 0 \end{bmatrix} \cup$$

$$\begin{bmatrix} 3 & | & 1 & 0 & 1 & | & 1 & 2 & 3 & 4 & 5 \\ 1 & | & 0 & 0 & 1 & | & 0 & 1 & 1 & 0 & 1 \\ 2 & | & 0 & 1 & 5 & | & 0 & 0 & 1 & 0 & 1 \end{bmatrix}$$

$$\cup \begin{bmatrix} 1 & 2 & 3 & 4 & 1 & 0 & 2 & 2 & 0 & 1 \\ 2 & 1 & 4 & 3 & 1 & 1 & 1 & 0 & 1 & 0 \\ 5 & 0 & 1 & 2 & 0 & 1 & 1 & 3 & 0 & 1 \\ 0 & 1 & 0 & 0 & 0 & 0 & 1 & 1 & 2 & 2 \end{bmatrix}$$

be a row super trivector.

Now

$$\begin{aligned} T^T &= (T_1 \cup T_2 \cup T_3)^T \\ &= T_1^T \cup T_2^T \cup T_3^T \end{aligned}$$

$$= \begin{bmatrix} 3 & 1 \\ 0 & 2 \\ 4 & 5 \\ \hline 7 & 0 \\ 2 & 1 \\ 3 & 0 \\ 1 & 1 \\ 5 & 0 \\ \hline 1 & 2 \\ 2 & 0 \end{bmatrix} \cup \begin{bmatrix} 3 & 1 & 2 \\ 1 & 0 & 0 \\ 0 & 0 & 1 \\ 1 & 1 & 5 \\ 1 & 0 & 0 \\ 2 & 1 & 0 \\ 3 & 1 & 1 \\ 4 & 0 & 0 \\ 5 & 1 & 1 \end{bmatrix} \cup \begin{bmatrix} 1 & 2 & 5 & 0 \\ 2 & 1 & 0 & 1 \\ 3 & 4 & 1 & 0 \\ 4 & 3 & 2 & 0 \\ \hline 1 & 1 & 0 & 0 \\ 0 & 1 & 1 & 0 \\ 2 & 1 & 1 & 1 \\ \hline 2 & 0 & 3 & 1 \\ 0 & 1 & 0 & 2 \\ 1 & 0 & 1 & 2 \end{bmatrix}.$$

$$\begin{aligned} TT^T &= (T_1 \cup T_2 \cup T_3)(T_1 \cup T_2 \cup T_3)^T \\ &= (T_1 \cup T_2 \cup T_3)(T_1^T \cup T_2^T \cup T_3^T) \\ &= T_1 T_1^T \cup T_2 T_2^T \cup T_3 T_3^T \end{aligned}$$

$$= \begin{bmatrix} 3 & 0 & 4 & 7 & 2 & 3 & 1 & 5 & 1 & 2 \\ 1 & 2 & 5 & 0 & 1 & 0 & 1 & 0 & 2 & 0 \end{bmatrix} \begin{bmatrix} 3 & 1 \\ 0 & 2 \\ 4 & 5 \\ \hline 7 & 0 \\ 2 & 1 \\ 3 & 0 \\ 1 & 1 \\ 5 & 0 \\ \hline 1 & 2 \\ 2 & 0 \end{bmatrix} \cup$$

$$\begin{bmatrix} 3 & 1 & 0 & 1 & 1 & 2 & 3 & 4 & 5 \\ 1 & 0 & 0 & 1 & 0 & 1 & 1 & 0 & 1 \\ 2 & 0 & 1 & 5 & 0 & 0 & 1 & 0 & 1 \end{bmatrix} \begin{bmatrix} 3 & 1 & 2 \\ \hline 1 & 0 & 0 \\ 0 & 0 & 1 \\ 1 & 1 & 5 \\ \hline 1 & 0 & 0 \\ 2 & 1 & 0 \\ 3 & 1 & 1 \\ 4 & 0 & 0 \\ 5 & 1 & 1 \end{bmatrix} \cup$$

$$\begin{bmatrix} 1 & 2 & 3 & 4 & 1 & 0 & 2 & 2 & 0 & 1 \\ 2 & 1 & 4 & 3 & 1 & 1 & 1 & 0 & 1 & 0 \\ 5 & 0 & 1 & 2 & 0 & 1 & 1 & 3 & 0 & 1 \\ 0 & 1 & 0 & 0 & 0 & 0 & 1 & 1 & 2 & 2 \end{bmatrix} \begin{bmatrix} 1 & 2 & 5 & 0 \\ 2 & 1 & 0 & 1 \\ 3 & 4 & 1 & 0 \\ 4 & 3 & 2 & 0 \\ \hline 1 & 1 & 0 & 0 \\ 0 & 1 & 1 & 0 \\ 2 & 1 & 1 & 0 \\ \hline 2 & 0 & 3 & 1 \\ 0 & 1 & 0 & 2 \\ 1 & 0 & 1 & 2 \end{bmatrix}$$

$$= \left\{ \begin{bmatrix} 3 & 0 & 4 \\ 1 & 2 & 5 \end{bmatrix} \begin{bmatrix} 3 & 1 \\ 0 & 2 \\ 4 & 5 \end{bmatrix} + \begin{bmatrix} 7 & 2 & 3 & 1 & 5 \\ 0 & 1 & 0 & 1 & 0 \end{bmatrix} \begin{bmatrix} 7 & 0 \\ 2 & 1 \\ 3 & 0 \\ 1 & 1 \\ 5 & 0 \end{bmatrix} + \right.$$

$$\left. \begin{bmatrix} 1 & 2 \\ 2 & 0 \end{bmatrix} \begin{bmatrix} 1 & 2 \\ 2 & 0 \end{bmatrix} \right\} \cup$$

$$\left\{ \begin{bmatrix} 3 \\ 1 \\ 2 \end{bmatrix} \begin{bmatrix} 3 & 1 & 2 \end{bmatrix} + \begin{bmatrix} 1 & 0 & 1 \\ 0 & 0 & 1 \\ 0 & 1 & 5 \end{bmatrix} \begin{bmatrix} 1 & 0 & 0 \\ 0 & 0 & 1 \\ 1 & 1 & 5 \end{bmatrix} + \right.$$

$$\left. \begin{bmatrix} 1 & 2 & 3 & 4 & 5 \\ 0 & 1 & 1 & 0 & 1 \\ 0 & 0 & 1 & 0 & 1 \end{bmatrix} \begin{bmatrix} 1 & 0 & 0 \\ 2 & 1 & 0 \\ 3 & 1 & 1 \\ 4 & 0 & 0 \\ 5 & 1 & 1 \end{bmatrix} \right\} \cup$$

$$\left\{ \begin{bmatrix} 1 & 2 & 3 & 4 \\ 2 & 1 & 4 & 3 \\ 5 & 0 & 1 & 2 \\ 0 & 1 & 0 & 0 \end{bmatrix} \begin{bmatrix} 1 & 2 & 5 & 0 \\ 2 & 1 & 0 & 1 \\ 3 & 4 & 1 & 0 \\ 4 & 3 & 2 & 0 \end{bmatrix} + \begin{bmatrix} 1 & 0 & 2 \\ 1 & 1 & 1 \\ 0 & 1 & 1 \\ 0 & 0 & 1 \end{bmatrix} \begin{bmatrix} 1 & 1 & 0 & 0 \\ 0 & 1 & 1 & 0 \\ 2 & 1 & 1 & 1 \end{bmatrix} + \right.$$

$$\left. \begin{bmatrix} 2 & 0 & 1 \\ 0 & 1 & 0 \\ 3 & 0 & 1 \\ 1 & 2 & 2 \end{bmatrix} \begin{bmatrix} 2 & 0 & 3 & 1 \\ 0 & 1 & 0 & 2 \\ 1 & 0 & 1 & 2 \end{bmatrix} \right\}$$

$$\left\{ \begin{bmatrix} 25 & 23 \\ 23 & 30 \end{bmatrix} + \begin{bmatrix} 88 & 3 \\ 3 & 2 \end{bmatrix} + \begin{bmatrix} 5 & 2 \\ 2 & 4 \end{bmatrix} \right\} \cup$$

$$\left\{ \begin{bmatrix} 9 & 3 & 6 \\ 3 & 1 & 2 \\ 6 & 2 & 4 \end{bmatrix} + \begin{bmatrix} 2 & 1 & 5 \\ 1 & 1 & 5 \\ 5 & 5 & 26 \end{bmatrix} + \begin{bmatrix} 55 & 10 & 8 \\ 10 & 3 & 2 \\ 8 & 2 & 2 \end{bmatrix} \right\} \cup$$

$$\left\{ \begin{bmatrix} 30 & 28 & 16 & 2 \\ 28 & 29 & 20 & 1 \\ 16 & 20 & 30 & 0 \\ 2 & 1 & 0 & 1 \end{bmatrix} + \begin{bmatrix} 5 & 3 & 2 & 2 \\ 3 & 3 & 2 & 1 \\ 2 & 2 & 2 & 1 \\ 2 & 1 & 1 & 1 \end{bmatrix} + \begin{bmatrix} 5 & 0 & 7 & 4 \\ 0 & 1 & 0 & 2 \\ 7 & 0 & 10 & 5 \\ 4 & 2 & 5 & 9 \end{bmatrix} \right\}$$

$$= \begin{bmatrix} 118 & 28 \\ 28 & 36 \end{bmatrix} \cup \begin{bmatrix} 66 & 14 & 19 \\ 14 & 5 & 9 \\ 19 & 9 & 32 \end{bmatrix} \cup \begin{bmatrix} 40 & 31 & 25 & 8 \\ 31 & 33 & 22 & 4 \\ 25 & 22 & 42 & 6 \\ 8 & 4 & 6 & 11 \end{bmatrix}.$$

We see the resultant of $TT^T$ is a symmetric trimatrix which is not a super trimatrix.

***Example 3.33:*** Let $T = T_1 \cup T_2 \cup T_3$ be a row super trivector; to find the product of $T$ with $T^T$. Given

$$T = T_1 \cup T_2 \cup T_3$$

$$= \begin{bmatrix} 3 & 1 & 4 & 1 & | & 0 & 1 & | & 0 & 1 & 0 & 1 & 2 \\ 0 & 2 & 0 & 1 & | & 2 & 5 & | & 1 & 0 & 0 & 0 & 1 \end{bmatrix} \cup$$

$$\begin{bmatrix} 1 & 3 & 1 & | & 0 & 1 & 1 & 3 & | & 2 & 1 \\ 4 & 0 & 1 & | & 1 & 0 & 2 & 1 & | & 0 & 1 \\ 3 & 0 & 0 & | & 0 & 1 & 0 & 0 & | & 2 & 2 \end{bmatrix}$$

$$\cup \begin{bmatrix} 4 & 1 & 1 & 1 & 1 & 2 & 1 & 0 & 1 & 3 & 1 & 0 \\ 0 & 1 & 0 & 0 & 1 & 0 & 1 & 1 & 2 & 2 & 1 & 0 \\ 1 & 0 & 0 & 1 & 1 & 1 & 3 & 0 & 3 & 1 & 0 & 0 \\ 2 & 0 & 1 & 0 & 1 & 1 & 0 & 0 & 1 & 0 & 0 & 1 \end{bmatrix}.$$

Now

$$\begin{aligned} T^T &= (T_1 \cup T_2 \cup T_3)^T \\ &= T_1^T \cup T_2^T \cup T_3^T \end{aligned}$$

$$= \begin{bmatrix} 3 & 0 \\ 1 & 2 \\ 4 & 0 \\ 1 & 1 \\ \hline 0 & 2 \\ 1 & 5 \\ \hline 0 & 1 \\ 1 & 0 \\ 0 & 0 \\ 1 & 0 \\ 2 & 1 \end{bmatrix} \cup \begin{bmatrix} 1 & 4 & 3 \\ 3 & 0 & 0 \\ 1 & 1 & 0 \\ \hline 0 & 1 & 0 \\ 1 & 0 & 1 \\ 1 & 2 & 0 \\ \hline 3 & 1 & 0 \\ 2 & 0 & 2 \\ 1 & 1 & 2 \end{bmatrix} \cup \begin{bmatrix} 4 & 0 & 1 & 2 \\ 1 & 1 & 0 & 0 \\ 1 & 0 & 0 & 1 \\ 1 & 0 & 1 & 0 \\ \hline 1 & 1 & 1 & 1 \\ 2 & 0 & 1 & 1 \\ \hline 1 & 1 & 3 & 0 \\ 0 & 1 & 0 & 0 \\ 0 & 2 & 3 & 1 \\ \hline 3 & 2 & 1 & 0 \\ 1 & 1 & 0 & 0 \\ 0 & 0 & 0 & 1 \end{bmatrix}.$$

$$\begin{aligned} TT^T &= (T_1 \cup T_2 \cup T_3)(T_1 \cup T_2 \cup T_3)^T \\ &= (T_1 \cup T_2 \cup T_3)(T_1^T \cup T_2^T \cup T_3^T) \\ &= TT_1^T \cup T_2 T_2^T \cup T_3 T_3^T \end{aligned}$$

$$= \begin{bmatrix} 3 & 1 & 4 & 1 & | & 0 & 1 & | & 0 & 1 & 0 & 1 & 2 \\ 0 & 2 & 0 & 1 & | & 2 & 5 & | & 1 & 0 & 0 & 0 & 1 \end{bmatrix} \begin{bmatrix} 3 & 0 \\ 1 & 2 \\ 4 & 0 \\ 1 & 1 \\ \hline 0 & 2 \\ 1 & 5 \\ \hline 0 & 1 \\ 1 & 0 \\ 0 & 0 \\ 1 & 0 \\ 2 & 1 \end{bmatrix} \cup$$

$$\begin{bmatrix} 1 & 3 & 1 & | & 0 & 1 & 1 & 3 & | & 2 & 1 \\ 4 & 0 & 1 & | & 1 & 0 & 2 & 1 & | & 0 & 1 \\ 3 & 0 & 0 & | & 0 & 1 & 0 & 0 & | & 2 & 2 \end{bmatrix} \begin{bmatrix} 1 & 4 & 3 \\ 3 & 0 & 0 \\ 1 & 1 & 0 \\ \hline 0 & 1 & 0 \\ 1 & 0 & 1 \\ 1 & 2 & 0 \\ 3 & 1 & 0 \\ \hline 2 & 0 & 2 \\ 1 & 1 & 2 \end{bmatrix} \cup$$

$$\begin{bmatrix} 4 & 1 & 1 & 1 & 1 & 2 & 0 & 1 & 0 & 3 & 1 & 0 \\ 0 & 1 & 0 & 0 & 1 & 0 & 1 & 1 & 2 & 2 & 1 & 0 \\ 1 & 0 & 0 & 1 & 1 & 1 & 3 & 0 & 3 & 1 & 0 & 0 \\ 2 & 0 & 1 & 0 & 1 & 1 & 0 & 0 & 1 & 0 & 0 & 1 \end{bmatrix} \begin{bmatrix} 4 & 0 & 1 & 2 \\ 1 & 1 & 0 & 0 \\ 1 & 0 & 0 & 1 \\ 1 & 0 & 1 & 0 \\ 1 & 1 & 1 & 1 \\ \hline 2 & 0 & 1 & 1 \\ 0 & 1 & 3 & 0 \\ 1 & 1 & 0 & 0 \\ 0 & 2 & 3 & 1 \\ \hline 3 & 2 & 1 & 0 \\ 1 & 1 & 0 & 0 \\ 0 & 0 & 0 & 1 \end{bmatrix}$$

$$= \left\{ \begin{bmatrix} 3 & 1 & 4 & 1 \\ 0 & 2 & 0 & 1 \end{bmatrix} \begin{bmatrix} 3 & 0 \\ 1 & 2 \\ 4 & 0 \\ 1 & 1 \end{bmatrix} + \right.$$

$$\left. \begin{bmatrix} 0 & 1 \\ 2 & 5 \end{bmatrix} \begin{bmatrix} 0 & 2 \\ 1 & 5 \end{bmatrix} + \begin{bmatrix} 0 & 1 & 0 & 2 & 1 \\ 1 & 0 & 0 & 0 & 1 \end{bmatrix} \begin{bmatrix} 0 & 1 \\ 1 & 0 \\ 0 & 0 \\ 2 & 0 \\ 1 & 1 \end{bmatrix} \right\}$$

$$\cup \left\{ \begin{bmatrix} 1 & 3 & 1 \\ 4 & 0 & 1 \\ 3 & 0 & 0 \end{bmatrix} \begin{bmatrix} 1 & 4 & 3 \\ 3 & 0 & 0 \\ 1 & 1 & 0 \end{bmatrix} + \begin{bmatrix} 0 & 1 & 1 & 3 \\ 1 & 0 & 2 & 1 \\ 0 & 1 & 0 & 0 \end{bmatrix} \begin{bmatrix} 0 & 1 & 0 \\ 1 & 0 & 1 \\ 1 & 2 & 0 \\ 3 & 1 & 0 \end{bmatrix} + \right.$$

$$\left. \begin{bmatrix} 2 & 1 \\ 0 & 1 \\ 2 & 2 \end{bmatrix} \begin{bmatrix} 2 & 0 & 2 \\ 1 & 1 & 2 \end{bmatrix} \right\}$$

$$\cup \left\{ \begin{bmatrix} 4 & 1 & 1 & 1 & 1 \\ 0 & 1 & 0 & 0 & 1 \\ 1 & 0 & 0 & 1 & 1 \\ 2 & 0 & 1 & 0 & 1 \end{bmatrix} \begin{bmatrix} 4 & 0 & 1 & 2 \\ 1 & 1 & 0 & 0 \\ 1 & 0 & 0 & 1 \\ 1 & 0 & 1 & 0 \\ 1 & 1 & 1 & 1 \end{bmatrix} + \right.$$

$$\begin{bmatrix} 2 & 0 & 1 & 0 \\ 0 & 1 & 1 & 2 \\ 1 & 3 & 0 & 3 \\ 1 & 0 & 0 & 1 \end{bmatrix} \begin{bmatrix} 2 & 0 & 1 & 1 \\ 0 & 1 & 3 & 0 \\ 1 & 1 & 0 & 0 \\ 0 & 2 & 3 & 1 \end{bmatrix} + \begin{bmatrix} 3 & 1 & 0 \\ 2 & 1 & 0 \\ 1 & 0 & 0 \\ 0 & 0 & 1 \end{bmatrix} \begin{bmatrix} 3 & 2 & 1 & 0 \\ 1 & 1 & 0 & 0 \\ 0 & 0 & 0 & 1 \end{bmatrix} \right\} =$$

$$\left\{ \begin{bmatrix} 27 & 3 \\ 3 & 5 \end{bmatrix} + \begin{bmatrix} 1 & 5 \\ 5 & 29 \end{bmatrix} + \begin{bmatrix} 6 & 1 \\ 1 & 2 \end{bmatrix} \right\}$$

$$\cup \left\{ \begin{bmatrix} 11 & 5 & 3 \\ 5 & 17 & 12 \\ 3 & 12 & 9 \end{bmatrix} + \begin{bmatrix} 11 & 5 & 1 \\ 5 & 6 & 0 \\ 1 & 0 & 1 \end{bmatrix} + \begin{bmatrix} 5 & 1 & 6 \\ 1 & 1 & 2 \\ 6 & 2 & 8 \end{bmatrix} \right\} \cup$$

$$\left\{ \begin{bmatrix} 20 & 2 & 6 & 10 \\ 2 & 2 & 1 & 1 \\ 6 & 1 & 3 & 3 \\ 10 & 1 & 3 & 6 \end{bmatrix} + \begin{bmatrix} 5 & 1 & 2 & 2 \\ 1 & 6 & 9 & 2 \\ 2 & 9 & 19 & 4 \\ 2 & 2 & 4 & 2 \end{bmatrix} + \begin{bmatrix} 10 & 7 & 3 & 0 \\ 7 & 5 & 2 & 0 \\ 3 & 2 & 1 & 0 \\ 0 & 0 & 0 & 1 \end{bmatrix} \right\}$$

$$= \begin{bmatrix} 34 & 9 \\ 9 & 36 \end{bmatrix} + \begin{bmatrix} 27 & 11 & 10 \\ 11 & 24 & 14 \\ 10 & 14 & 18 \end{bmatrix} + \begin{bmatrix} 35 & 10 & 11 & 12 \\ 10 & 13 & 12 & 3 \\ 11 & 12 & 23 & 7 \\ 12 & 3 & 7 & 9 \end{bmatrix}.$$

is a trimatrix which is symmetric, clearly not a super trimatrix.

**DEFINITION 3.15:** *Let $T = T_1 \cup T_2 \cup T_3$ be a super trimatrix, if every $T_i$ is symmetric supermatrix $i = 1, 2, 3$ then we say $T$ is symmetric super trimatrix. Clearly if $T$ is a symmetric super trimatrix then $T$ should be either a mixed square symmetric super trimatrix or square symmetric super trimatrix.*

We illustrate this by the following examples.

*Example 3.34:* Let $T = T_1 \cup T_2 \cup T_3$ be a super trimatrix where

$$T_1 = \begin{bmatrix} 0 & 1 & 2 & 3 \\ \hline 1 & 4 & 5 & 6 \\ 2 & 5 & 0 & 2 \\ 3 & 6 & 2 & 7 \end{bmatrix},$$

$$T_2 = \begin{bmatrix} 0 & 1 & 0 & 1 \\ 1 & 2 & 0 & 2 \\ \hline 0 & 0 & 4 & 5 \\ 1 & 2 & 5 & 3 \end{bmatrix}$$

and

$$T_3 = \begin{bmatrix} 0 & 1 & 2 & 5 \\ 1 & 0 & 7 & 8 \\ 2 & 7 & 6 & 1 \\ \hline 5 & 8 & 1 & 5 \end{bmatrix}.$$

We see $T = T_1 \cup T_2 \cup T_3$ is a square symmetric super trimatrix. Further each $T_i$; $1 \le i \le 3$ are $4 \times 4$ symmetric super trimatrix.

*Example 3.35:* Let $V = V_1 \cup V_2 \cup V_3$ where

$$V_1 = \begin{bmatrix} 3 & 1 & 2 & 3 & 0 \\ 1 & 7 & 0 & 1 & 2 \\ \hline 2 & 0 & 5 & 0 & 1 \\ 3 & 1 & 0 & 2 & 3 \\ 0 & 2 & 1 & 3 & 1 \end{bmatrix},$$

$$V_2 = \begin{bmatrix} 1 & 0 & 1 & 2 & 3 & 0 \\ 0 & 5 & 0 & 1 & 1 & 0 \\ 1 & 0 & 7 & 2 & 0 & 2 \\ \hline 2 & 1 & 2 & 9 & 1 & 0 \\ 3 & 1 & 0 & 1 & 11 & 6 \\ 0 & 0 & 2 & 0 & 6 & 2 \end{bmatrix}$$

and

$$V_3 = \begin{bmatrix} 9 & 2 & 1 & 2 & 3 & 0 & 7 \\ 2 & 0 & 5 & 1 & 0 & 6 & 1 \\ \hline 1 & 5 & 3 & -1 & 2 & 0 & 6 \\ 2 & 1 & -1 & 7 & 1 & 1 & 5 \\ 3 & 0 & 2 & 1 & -3 & 5 & -1 \\ \hline 0 & 6 & 0 & 1 & 5 & 8 & 2 \\ 7 & 1 & 6 & 5 & -1 & 2 & 6 \end{bmatrix}$$

be a symmetric super trimatrix. Clearly V is a mixed square symmetric super trimatrix.

*Example 3.36:* Let $T = T_1 \cup T_2 \cup T_3$ where

$$T_1 = \begin{bmatrix} 3 & 1 & 3 & 0 & 5 \\ 1 & 7 & 2 & 1 & 0 \\ \hline 3 & 2 & 0 & 3 & 7 \\ 0 & 1 & 3 & 1 & 6 \\ 5 & 0 & 7 & 6 & 9 \end{bmatrix},$$

$$T_2 = \left[\begin{array}{cc|cccc} 0 & 1 & 2 & 3 & 0 & 6 \\ 1 & 2 & 1 & 0 & 1 & 2 \\ 2 & 1 & 9 & 6 & 0 & 3 \\ 3 & 0 & 6 & 1 & 2 & 1 \\ \hline 0 & 1 & 0 & 2 & 5 & 8 \\ 6 & 2 & 3 & 1 & 8 & 7 \end{array}\right]$$

and

$$T_3 = \left[\begin{array}{ccc|cc} 1 & 2 & 3 & 0 & 5 \\ 2 & 0 & 1 & 6 & 1 \\ 3 & 1 & 2 & 1 & 0 \\ \hline 0 & 6 & 1 & 8 & 9 \\ 5 & 1 & 0 & 9 & 6 \end{array}\right].$$

We see each of the matrices $T_1$, $T_2$ and $T_3$ are symmetric matrices.

Further they are also supermatrices but $T = T_1 \cup T_2 \cup T_3$ is not a symmetric super trimatrix as $T_2$ is only a symmetric matrix but $T_2$ is not symmetric supermatrix though $T_2$ is a supermatrix. Thus T is only a super trimatrix which is not symmetric.

From this we see each of the supermatrix which is symmetric must be partitioned, such that it is a symmetric supermatrix. Though $T_2$ is symmetric it is not symmetric supermatrix as the partition happens to yield a non symmetric supermatrix.

*Example 3.37:* Let $T = T_1 \cup T_2 \cup T_3$ where

$$T_1 = \left[\begin{array}{cccc|cc} 3 & 4 & 0 & 1 & 5 & 6 \\ 1 & 1 & 6 & 2 & 1 & 5 \\ 0 & 3 & 1 & 2 & 4 & 1 \end{array}\right]$$

$$T_2 = \begin{bmatrix} 0 & 1 & 2 & 3 & 4 & 5 \\ 1 & 6 & 0 & 1 & 2 & 1 \\ 2 & 0 & 7 & 6 & 0 & 2 \\ \hline 3 & 1 & 6 & 5 & 1 & 3 \\ 4 & 2 & 0 & 1 & 7 & 2 \\ 5 & 1 & 2 & 3 & 2 & 8 \end{bmatrix}$$

and

$$T_3 = \begin{bmatrix} 0 & 1 & 2 & 3 & 0 & 1 & 8 \\ 1 & 0 & 1 & 2 & 3 & 0 & 1 \\ 2 & 1 & 7 & 0 & 1 & 2 & 0 \\ 3 & 2 & 0 & 9 & 3 & 1 & 2 \\ 0 & 3 & 1 & 3 & 6 & 0 & 1 \\ 1 & 0 & 2 & 3 & 6 & 0 & 1 \\ 8 & 1 & 0 & 2 & 1 & 5 & 7 \end{bmatrix}.$$

We see T is a super trimatrix but it is not a symmetric trimatrix. Only one of the matrices $T_2$ alone is a super symmetric matrix. It is not even a square super trimatrix. Thus T is only a mixed super trimatrix.

*Example 3.38:* Let $T = T_1 \cup T_2 \cup T_3$ where

$$T_1 = \begin{bmatrix} 1 & 1 & 0 & 2 & 3 \\ 1 & 7 & 9 & 0 & 6 \\ \hline 0 & 9 & 1 & 2 & 1 \\ 2 & 0 & 2 & 0 & 2 \\ 3 & 6 & 1 & 1 & 1 \end{bmatrix},$$

$$T_2 = \left[\begin{array}{ccc|ccc} 4 & 1 & 0 & 2 & 3 & 1 \\ 1 & 0 & 8 & 9 & 6 & 3 \\ 0 & 8 & 7 & 1 & 2 & 3 \\ \hline 2 & 9 & 1 & 2 & 0 & 1 \\ 3 & 6 & 2 & 0 & 5 & 3 \\ 1 & 3 & 3 & 1 & 3 & 0 \end{array}\right]$$

and

$$T_3 = \left[\begin{array}{ccc|ccc} 3 & 1 & 0 & 1 & 3 & 1 \\ 1 & 2 & 1 & 2 & 3 & 4 \\ 0 & 1 & 5 & 1 & 2 & 3 \\ 1 & 2 & 1 & 0 & 1 & 2 \\ \hline 3 & 3 & 2 & 1 & 7 & 5 \\ 1 & 4 & 3 & 2 & 5 & 3 \end{array}\right].$$

T is a mixed square super trimatrix but T is not a symmetric mixed square super trimatrix. For $T_3$ is a supermatrix further $T_3$ is a symmetric matrix but after partition $T_3$ is not a symmetric super trimatrix. $T_1$ is not a symmetric matrix only a supermatrix. $T_2$ is a symmetric supermatrix. Thus T is only a mixed square super trimatrix.

*Example 3.39:* Let $V = V_1 \cup V_2 \cup V_3$ where

$$V_1 = \left[\begin{array}{ccc|c} 0 & 1 & 2 & 3 \\ 1 & 5 & 6 & 1 \\ 2 & 6 & 0 & 4 \\ \hline 3 & 1 & 4 & 7 \end{array}\right],$$

$$V_2 = \left[\begin{array}{cc|c} 3 & 0 & 1 \\ 0 & 1 & 2 \\ \hline 1 & 2 & 1 \end{array}\right]$$

and

$$V_3 = \begin{bmatrix} 7 & 1 & 2 & 3 & 4 & 5 & 6 \\ 1 & 9 & 1 & 2 & 3 & 4 & 5 \\ 2 & 1 & 8 & 1 & 2 & 3 & 4 \\ \hline 3 & 2 & 1 & 6 & 1 & 2 & 3 \\ 4 & 3 & 2 & 1 & 5 & 1 & 2 \\ \hline 5 & 4 & 3 & 2 & 1 & 4 & 1 \\ 6 & 5 & 4 & 3 & 2 & 1 & 3 \end{bmatrix}.$$

We see T is a mixed square super trimatrix which is also a mixed square symmetric super trimatrix.

Thus we are interested in studying those super trimatrices T = $T_1 \cup T_2 \cup T_3$ in which at least one of them is a symmetric supermatrix. To this end we give the following definition.

**DEFINITION 3.16:** *Let $T = T_1 \cup T_2 \cup T_3$ be a super trimatrix we say T is a quasi symmetric super trimatrix if at least one of the $T_i$ is a symmetric supermatrix $1 \leq i \leq 3$.*

1. It may so happen all the 3 matrices are symmetric matrices; yet all of them are not super symmetric, only one is a symmetric supermatrix.
2. It may so happen only one of the matrices $T_i$ alone is a super symmetric matrix where as others are rectangular supermatrices.
3. It may so happen T is a square super trimatrix or a mixed square super trimatrix where only one of the $T_i$ is a symmetric super trimatrix.

Thus in all these cases also we call T to be a quasi symmetric super trimatrix. Now we have seen if T is row super trivector then T. $T^T$ the product of T with its transpose yields a trimatrix which is not a super trimatrix but it is a symmetric trimatrix. Thus we have a method by which we can generate symmetric trimatrices, of course it may be square symmetric trimatrix or a mixed square symmetric trimatrix.

Now do we have any method of generating symmetric super trimatrices? The answer is yes and now we proceed on to generate them by a special product.

*Example 3.40:* Let $T = T_1 \cup T_2 \cup T_3$ be a super trimatrix where

$$T_1 = \begin{bmatrix} 0 & 1 & 2 & 1 & 1 & 1 \\ 1 & 0 & 1 & 1 & 0 & 1 \\ \hline 3 & 1 & 0 & 0 & 1 & 1 \\ \hline 4 & 0 & 2 & 1 & 0 & 0 \\ 5 & 2 & 0 & 0 & 1 & 0 \end{bmatrix},$$

$$T_2 = \begin{bmatrix} 1 & 2 & 3 & 0 & 1 & 2 & 5 \\ 0 & 1 & 1 & 1 & 0 & 0 & 1 \\ 1 & 0 & 0 & 0 & 1 & 1 & 0 \\ \hline 3 & 1 & 0 & 1 & 0 & 1 & 0 \\ 4 & 2 & 1 & 2 & 1 & 2 & 1 \\ 5 & 0 & 1 & 0 & 1 & 0 & 1 \end{bmatrix}$$

and

$$T_3 = \begin{bmatrix} 1 & 1 & 0 & 1 & 1 \\ 2 & 2 & 1 & 0 & 0 \\ 3 & 1 & 0 & 0 & 1 \\ \hline 0 & 6 & 1 & 0 & 2 \\ 1 & 1 & 1 & 0 & 1 \\ 2 & 0 & 1 & 0 & 2 \\ 0 & 1 & 0 & 1 & 0 \end{bmatrix}.$$

Now

$$T^T = (T_1 \cup T_2 \cup T_3)^T = T_1^T \cup T_2^T \cup T_3^T$$

where

$$T_1^T = \left[\begin{array}{cc|ccc}
0 & 1 & 3 & 4 & 5 \\
1 & 0 & 1 & 0 & 2 \\
2 & 1 & 0 & 2 & 0 \\
\hline
1 & 1 & 0 & 1 & 0 \\
1 & 0 & 1 & 0 & 1 \\
1 & 1 & 1 & 0 & 0
\end{array}\right],$$

$$T_2^T = \left[\begin{array}{c|cc|ccc}
1 & 0 & 1 & 3 & 4 & 5 \\
2 & 1 & 0 & 1 & 2 & 0 \\
3 & 1 & 0 & 0 & 1 & 1 \\
\hline
0 & 1 & 0 & 1 & 2 & 0 \\
1 & 0 & 1 & 0 & 1 & 1 \\
2 & 0 & 1 & 1 & 2 & 0 \\
5 & 1 & 0 & 0 & 1 & 1
\end{array}\right]$$

and

$$T_3^T = \left[\begin{array}{c|cc|ccc}
1 & 2 & 3 & 0 & 1 & 2 & 0 \\
1 & 2 & 1 & 6 & 1 & 0 & 1 \\
0 & 1 & 0 & 1 & 1 & 1 & 0 \\
1 & 0 & 0 & 0 & 0 & 0 & 1 \\
1 & 0 & 1 & 2 & 1 & 2 & 0
\end{array}\right].$$

$$\begin{aligned}
TT^T &= (T_1 \cup T_2 \cup T_3)(T_1 \cup T_2 \cup T_3)^T \\
&= (T_1 \cup T_2 \cup T_3)(T_1^T \cup T_2^T \cup T_3^T) \\
&= T_1 T_1^T \cup T_2 T_2^T \cup T_3 T_3^T
\end{aligned}$$

$$= \begin{bmatrix} 0 & 1 & 2 & 1 & 1 & 1 \\ 1 & 0 & 1 & 1 & 0 & 1 \\ \hline 3 & 1 & 0 & 0 & 1 & 1 \\ 4 & 0 & 2 & 1 & 0 & 0 \\ \hline 5 & 2 & 0 & 0 & 1 & 0 \end{bmatrix} \begin{bmatrix} 0 & 1 & 3 & 4 & 5 \\ 1 & 0 & 1 & 0 & 2 \\ 2 & 1 & 0 & 2 & 0 \\ \hline 1 & 1 & 0 & 1 & 0 \\ 1 & 0 & 1 & 0 & 1 \\ \hline 1 & 1 & 1 & 0 & 0 \end{bmatrix} \cup$$

$$\begin{bmatrix} 1 & 2 & 3 & 0 & 1 & 2 & 5 \\ \hline 0 & 1 & 1 & 1 & 0 & 0 & 1 \\ 1 & 0 & 0 & 0 & 1 & 1 & 0 \\ \hline 3 & 1 & 0 & 1 & 0 & 1 & 0 \\ 4 & 2 & 0 & 2 & 1 & 2 & 1 \\ 5 & 0 & 1 & 0 & 1 & 0 & 1 \end{bmatrix} \begin{bmatrix} 1 & 0 & 1 & 3 & 4 & 5 \\ 2 & 1 & 0 & 1 & 2 & 0 \\ 3 & 1 & 0 & 0 & 1 & 1 \\ \hline 0 & 1 & 0 & 1 & 2 & 0 \\ 1 & 0 & 1 & 0 & 1 & 1 \\ 2 & 0 & 1 & 1 & 2 & 0 \\ \hline 5 & 1 & 0 & 0 & 1 & 1 \end{bmatrix}$$

$$\cup \begin{bmatrix} 1 & 1 & 0 & 1 & 1 \\ \hline 2 & 2 & 1 & 0 & 0 \\ 3 & 1 & 0 & 0 & 1 \\ \hline 0 & 6 & 1 & 0 & 2 \\ 1 & 1 & 1 & 0 & 1 \\ 2 & 0 & 1 & 0 & 2 \\ 0 & 1 & 0 & 1 & 0 \end{bmatrix} \begin{bmatrix} 1 & 2 & 3 & 0 & 1 & 2 & 0 \\ \hline 1 & 2 & 1 & 6 & 1 & 0 & 1 \\ 0 & 1 & 0 & 1 & 1 & 1 & 0 \\ \hline 1 & 0 & 0 & 0 & 0 & 0 & 1 \\ 1 & 0 & 1 & 2 & 1 & 2 & 0 \end{bmatrix}$$

$$= \left\{ \begin{bmatrix} 0 \\ 1 \\ \hline 3 \\ 4 \\ \hline 5 \end{bmatrix} [0 \quad 1 \mid 3 \quad 4 \mid 5] + \right.$$

$$\begin{bmatrix} 1 & 2 \\ 0 & 1 \\ 1 & 0 \\ 0 & 2 \\ 2 & 0 \end{bmatrix} \begin{bmatrix} 1 & 0 & | & 1 & 0 & | & 2 \\ 2 & 1 & | & 0 & 2 & | & 0 \end{bmatrix} + \begin{bmatrix} 1 & 1 & 1 \\ 1 & 0 & 1 \\ 0 & 1 & 1 \\ 1 & 0 & 0 \\ 0 & 1 & 0 \end{bmatrix} \begin{bmatrix} 1 & 1 & | & 0 & 1 & | & 0 \\ 1 & 0 & | & 1 & 0 & | & 1 \\ 1 & 1 & | & 1 & 0 & | & 0 \end{bmatrix} \Bigg\}$$

$$\cup \Bigg\{ \begin{bmatrix} 1 & 2 & 3 \\ 0 & 1 & 1 \\ 1 & 0 & 0 \\ 3 & 1 & 0 \\ 4 & 2 & 0 \\ 5 & 0 & 1 \end{bmatrix} \begin{bmatrix} 1 & | & 0 & 1 & | & 3 & 4 & 5 \\ 2 & | & 1 & 0 & | & 1 & 2 & 0 \\ 3 & | & 1 & 0 & | & 0 & 1 & 1 \end{bmatrix} +$$

$$\begin{bmatrix} 0 & 1 & 2 \\ 1 & 0 & 0 \\ 0 & 1 & 1 \\ 1 & 0 & 1 \\ 2 & 1 & 2 \\ 0 & 1 & 0 \end{bmatrix} \begin{bmatrix} 0 & | & 1 & 0 & | & 1 & 2 & 0 \\ 1 & | & 0 & 1 & | & 0 & 1 & 1 \\ 2 & | & 0 & 1 & | & 1 & 2 & 0 \end{bmatrix} + \begin{bmatrix} 5 \\ 1 \\ 0 \\ 0 \\ 1 \\ 1 \end{bmatrix} \begin{bmatrix} 5 & | & 1 & 0 & | & 0 & 1 & 1 \end{bmatrix} \Bigg\}$$

$$\cup \Bigg\{ \begin{bmatrix} 1 \\ 2 \\ 3 \\ 0 \\ 1 \\ 2 \\ 0 \end{bmatrix} \begin{bmatrix} 1 & | & 2 & 3 & | & 0 & 1 & 2 & 0 \end{bmatrix} +$$

$$\begin{bmatrix} 1 & 0 \\ 2 & 1 \\ \hline 1 & 0 \\ 6 & 1 \\ 1 & 1 \\ 0 & 1 \\ 1 & 0 \end{bmatrix} \left[ \begin{array}{c|cc|cccc} 1 & 2 & 1 & 6 & 1 & 0 & 1 \\ 0 & 1 & 0 & 1 & 1 & 1 & 0 \end{array} \right] +$$

$$\begin{bmatrix} 1 & 1 \\ 0 & 0 \\ \hline 0 & 1 \\ 0 & 2 \\ 0 & 1 \\ 0 & 2 \\ 1 & 0 \end{bmatrix} \left[ \begin{array}{c|cc|cccc} 1 & 0 & 0 & 0 & 0 & 0 & 1 \\ 1 & 0 & 1 & 2 & 1 & 2 & 0 \end{array} \right] \Bigg\}$$

$$= \left\{ \left[ \begin{array}{c|c|c} \begin{bmatrix}0\\1\end{bmatrix}[0\ 1] & \begin{bmatrix}0\\1\end{bmatrix}[3\ 4] & \begin{bmatrix}0\\1\end{bmatrix}[5] \\ \hline \begin{bmatrix}3\\4\end{bmatrix}[0\ 1] & \begin{bmatrix}3\\4\end{bmatrix}[3\ 4] & \begin{bmatrix}3\\4\end{bmatrix}[5] \\ \hline [5][0\ 1] & [5][3\ 4] & [5][5] \end{array} \right] + \right.$$

$$\left[ \begin{array}{c|c|c} \begin{pmatrix}1 & 2\\0 & 1\end{pmatrix}\begin{pmatrix}1 & 0\\2 & 1\end{pmatrix} & \begin{pmatrix}1 & 2\\0 & 1\end{pmatrix}\begin{pmatrix}1 & 0\\0 & 2\end{pmatrix} & \begin{pmatrix}1 & 2\\0 & 1\end{pmatrix}\begin{pmatrix}2\\0\end{pmatrix} \\ \hline \begin{pmatrix}1 & 0\\0 & 2\end{pmatrix}\begin{pmatrix}1 & 0\\2 & 1\end{pmatrix} & \begin{pmatrix}1 & 0\\0 & 2\end{pmatrix}\begin{pmatrix}1 & 0\\0 & 2\end{pmatrix} & \begin{pmatrix}1 & 0\\0 & 2\end{pmatrix}\begin{pmatrix}2\\0\end{pmatrix} \\ \hline (2\ 0)\begin{pmatrix}1 & 0\\0 & 2\end{pmatrix} & (2\ 0)\begin{pmatrix}1 & 0\\0 & 2\end{pmatrix} & (2\ 0)\begin{pmatrix}2\\0\end{pmatrix} \end{array} \right] +$$

$$\left\{ \begin{array}{c|c|c} \begin{bmatrix} 1 & 1 & 1 \\ 1 & 0 & 1 \end{bmatrix} \begin{bmatrix} 1 & 1 \\ 1 & 0 \\ 1 & 1 \end{bmatrix} & \begin{bmatrix} 1 & 1 & 1 \\ 1 & 0 & 1 \end{bmatrix} \begin{bmatrix} 0 & 1 \\ 1 & 0 \\ 1 & 0 \end{bmatrix} & \begin{bmatrix} 1 & 1 & 1 \\ 1 & 0 & 1 \end{bmatrix} \begin{bmatrix} 0 \\ 1 \\ 0 \end{bmatrix} \\ \hline \begin{bmatrix} 0 & 1 & 1 \\ 1 & 0 & 0 \end{bmatrix} \begin{bmatrix} 1 & 1 \\ 1 & 0 \\ 1 & 1 \end{bmatrix} & \begin{bmatrix} 0 & 1 & 1 \\ 1 & 0 & 0 \end{bmatrix} \begin{bmatrix} 0 & 1 \\ 1 & 0 \\ 1 & 0 \end{bmatrix} & \begin{bmatrix} 0 & 1 & 1 \\ 1 & 0 & 0 \end{bmatrix} \begin{bmatrix} 0 \\ 1 \\ 0 \end{bmatrix} \\ \hline \begin{bmatrix} 0 & 1 & 0 \end{bmatrix} \begin{bmatrix} 1 & 1 \\ 1 & 0 \\ 1 & 1 \end{bmatrix} & \begin{bmatrix} 0 & 1 & 0 \end{bmatrix} \begin{bmatrix} 0 & 1 \\ 1 & 0 \\ 1 & 0 \end{bmatrix} & \begin{bmatrix} 0 & 1 & 0 \end{bmatrix} \begin{bmatrix} 0 \\ 1 \\ 0 \end{bmatrix} \end{array} \right\} \subset$$

$$\left\{ \begin{array}{c|c|c} \begin{bmatrix} 1 & 2 & 3 \end{bmatrix} \begin{bmatrix} 1 \\ 2 \\ 3 \end{bmatrix} & \begin{bmatrix} 1 & 2 & 3 \end{bmatrix} \begin{bmatrix} 0 & 1 \\ 1 & 0 \\ 1 & 0 \end{bmatrix} & \begin{bmatrix} 1 & 2 & 3 \end{bmatrix} \begin{bmatrix} 3 & 4 & 5 \\ 1 & 2 & 0 \\ 0 & 1 & 1 \end{bmatrix} \\ \hline \begin{bmatrix} 0 & 1 & 1 \\ 1 & 0 & 0 \end{bmatrix} \begin{bmatrix} 1 \\ 2 \\ 3 \end{bmatrix} & \begin{bmatrix} 0 & 1 & 1 \\ 1 & 0 & 0 \end{bmatrix} \begin{bmatrix} 0 & 1 \\ 1 & 0 \\ 1 & 0 \end{bmatrix} & \begin{bmatrix} 0 & 1 & 1 \\ 1 & 0 & 0 \end{bmatrix} \begin{bmatrix} 3 & 4 & 5 \\ 1 & 2 & 0 \\ 0 & 1 & 1 \end{bmatrix} \\ \hline \begin{bmatrix} 3 & 1 & 0 \\ 4 & 2 & 1 \\ 5 & 0 & 1 \end{bmatrix} \begin{bmatrix} 1 \\ 2 \\ 3 \end{bmatrix} & \begin{bmatrix} 3 & 1 & 0 \\ 4 & 2 & 0 \\ 5 & 0 & 1 \end{bmatrix} \begin{bmatrix} 0 & 1 \\ 1 & 0 \\ 1 & 0 \end{bmatrix} & \begin{bmatrix} 3 & 1 & 0 \\ 4 & 2 & 0 \\ 5 & 0 & 1 \end{bmatrix} \begin{bmatrix} 3 & 4 & 5 \\ 1 & 2 & 0 \\ 0 & 0 & 1 \end{bmatrix} \end{array} \right\} +$$

$$\left[ \begin{array}{c|c|c} (0 \ 1 \ 2) \begin{bmatrix} 0 \\ 1 \\ 2 \end{bmatrix} & (0 \ 1 \ 2) \begin{bmatrix} 1 & 0 \\ 0 & 1 \\ 0 & 1 \end{bmatrix} & (0 \ 1 \ 2) \begin{bmatrix} 1 & 2 & 0 \\ 0 & 1 & 1 \\ 1 & 2 & 0 \end{bmatrix} \\ \hline \begin{pmatrix} 1 & 0 & 0 \\ 0 & 1 & 1 \end{pmatrix} \begin{bmatrix} 0 \\ 1 \\ 2 \end{bmatrix} & \begin{pmatrix} 1 & 0 & 0 \\ 0 & 1 & 1 \end{pmatrix} \begin{bmatrix} 1 & 0 \\ 0 & 1 \\ 0 & 1 \end{bmatrix} & \begin{pmatrix} 1 & 0 & 0 \\ 0 & 1 & 1 \end{pmatrix} \begin{bmatrix} 1 & 2 & 0 \\ 0 & 1 & 1 \\ 1 & 2 & 0 \end{bmatrix} \\ \hline \begin{pmatrix} 1 & 0 & 1 \\ 2 & 1 & 2 \\ 0 & 1 & 0 \end{pmatrix} \begin{bmatrix} 0 \\ 1 \\ 2 \end{bmatrix} & \begin{pmatrix} 1 & 0 & 1 \\ 2 & 1 & 2 \\ 0 & 1 & 0 \end{pmatrix} \begin{bmatrix} 1 & 0 \\ 0 & 1 \\ 0 & 1 \end{bmatrix} & \begin{pmatrix} 1 & 0 & 1 \\ 2 & 1 & 2 \\ 0 & 1 & 0 \end{pmatrix} \begin{bmatrix} 1 & 2 & 0 \\ 0 & 1 & 1 \\ 1 & 2 & 0 \end{bmatrix} \end{array} \right] +$$

$$\left[\begin{array}{c|c|c} (5)(5) & (5)(1\ 0) & (5)(0\ 1\ 1) \\ \hline \begin{pmatrix}1\\0\end{pmatrix}(5) & \begin{pmatrix}1\\0\end{pmatrix}(1\ 0) & \begin{pmatrix}1\\0\end{pmatrix}(0\ 1\ 1) \\ \hline \begin{pmatrix}0\\1\\1\end{pmatrix}(5) & \begin{pmatrix}0\\1\\1\end{pmatrix}(1\ 0) & \begin{pmatrix}0\\1\\1\end{pmatrix}(0\ 1\ 1) \end{array}\right] \cup$$

$$\cup \left\{ \left[\begin{array}{c|c|c} (1)(1) & (1)(2\ 3) & (1)(0\ 1\ 2\ 0) \\ \hline \begin{pmatrix}2\\3\end{pmatrix}(1) & \begin{pmatrix}2\\3\end{pmatrix}(2\ 3) & \begin{pmatrix}2\\3\end{pmatrix}(0\ 1\ 2\ 0) \\ \hline \begin{pmatrix}0\\1\\2\\0\end{pmatrix}(1) & \begin{pmatrix}0\\1\\2\\0\end{pmatrix}(2\ 3) & \begin{pmatrix}0\\1\\2\\0\end{pmatrix}(0\ 1\ 2\ 0) \end{array}\right] \right\} +$$

$$\left[\begin{array}{c|c|c} (1\ 0)\begin{pmatrix}1\\0\end{pmatrix} & (1\ 0)\begin{pmatrix}2 & 1\\1 & 0\end{pmatrix} & (1\ 0)\begin{pmatrix}6 & 1 & 0 & 1\\1 & 1 & 1 & 0\end{pmatrix} \\ \hline \begin{pmatrix}2 & 1\\1 & 0\end{pmatrix}\begin{pmatrix}1\\0\end{pmatrix} & \begin{pmatrix}2 & 1\\1 & 0\end{pmatrix}\begin{pmatrix}2 & 1\\1 & 0\end{pmatrix} & \begin{pmatrix}2 & 1\\1 & 0\end{pmatrix}\begin{pmatrix}6 & 1 & 0 & 1\\1 & 1 & 1 & 0\end{pmatrix} \\ \hline \begin{pmatrix}6 & 1\\1 & 1\\0 & 1\\1 & 0\end{pmatrix}\begin{pmatrix}1\\0\end{pmatrix} & \begin{pmatrix}6 & 1\\1 & 1\\0 & 1\\1 & 0\end{pmatrix}\begin{pmatrix}2 & 1\\1 & 0\end{pmatrix} & \begin{pmatrix}6 & 1\\1 & 1\\0 & 1\\1 & 0\end{pmatrix}\begin{pmatrix}6 & 1 & 0 & 1\\1 & 1 & 1 & 0\end{pmatrix} \end{array}\right]$$

$$+\left[\begin{array}{c|c|c}
(1\ 1)\begin{pmatrix}1\\1\end{pmatrix} & (1\ 1)\begin{pmatrix}0 & 0\\0 & 1\end{pmatrix} & (1\ 1)\begin{pmatrix}0 & 0 & 0 & 1\\2 & 1 & 2 & 0\end{pmatrix} \\
\hline
\begin{pmatrix}0 & 0\\0 & 1\end{pmatrix}\begin{pmatrix}1\\1\end{pmatrix} & \begin{pmatrix}0 & 0\\0 & 1\end{pmatrix}\begin{pmatrix}0 & 0\\0 & 1\end{pmatrix} & \begin{pmatrix}0 & 0\\0 & 1\end{pmatrix}\begin{pmatrix}0 & 0 & 0 & 1\\2 & 1 & 2 & 0\end{pmatrix} \\
\hline
\begin{pmatrix}0 & 2\\0 & 1\\0 & 2\\0 & 1\end{pmatrix}\begin{pmatrix}1\\1\end{pmatrix} & \begin{pmatrix}0 & 2\\0 & 1\\0 & 2\\0 & 1\end{pmatrix}\begin{pmatrix}0 & 0\\0 & 1\end{pmatrix} & \begin{pmatrix}0 & 2\\0 & 1\\0 & 2\\0 & 1\end{pmatrix}\begin{pmatrix}0 & 0 & 0 & 1\\2 & 1 & 2 & 0\end{pmatrix}
\end{array}\right]\Bigg\}$$

$$\cup\left\{\left[\begin{array}{cc|ccc|c}
0 & 0 & 0 & 0 & 0 \\
0 & 1 & 3 & 4 & 5 \\
\hline
0 & 3 & 9 & 12 & 15 \\
0 & 4 & 12 & 16 & 20 \\
\hline
0 & 5 & 15 & 20 & 25
\end{array}\right]+\left[\begin{array}{cc|cc|c}
5 & 2 & 1 & 4 & 2 \\
2 & 1 & 0 & 2 & 0 \\
\hline
1 & 0 & 1 & 0 & 2 \\
4 & 2 & 0 & 4 & 0 \\
\hline
2 & 0 & 2 & 0 & 4
\end{array}\right]+\right.$$

$$\left.\left[\begin{array}{cc|cc|c}
3 & 2 & 2 & 1 & 1 \\
2 & 2 & 1 & 1 & 0 \\
\hline
2 & 1 & 2 & 0 & 1 \\
1 & 1 & 0 & 1 & 0 \\
\hline
1 & 0 & 1 & 0 & 1
\end{array}\right]\right\}\cup$$

$$\left\{\left[\begin{array}{c|cc|ccc}
14 & 5 & 1 & 5 & 11 & 8 \\
\hline
5 & 2 & 0 & 1 & 3 & 1 \\
1 & 0 & 1 & 3 & 4 & 5 \\
\hline
5 & 1 & 3 & 11 & 14 & 15 \\
11 & 3 & 4 & 14 & 20 & 20 \\
8 & 1 & 5 & 15 & 20 & 26
\end{array}\right]+\left[\begin{array}{c|cc|cc|c}
5 & 0 & 3 & 2 & 5 & 1 \\
\hline
0 & 1 & 0 & 1 & 2 & 0 \\
3 & 0 & 2 & 1 & 3 & 1 \\
\hline
2 & 1 & 1 & 2 & 4 & 0 \\
5 & 2 & 3 & 4 & 9 & 1 \\
1 & 0 & 1 & 0 & 1 & 1
\end{array}\right]+\right.$$

$$\left\{ \begin{bmatrix} 25 & 5 & 0 & 0 & 5 & 5 \\ 5 & 1 & 0 & 0 & 1 & 1 \\ 0 & 0 & 0 & 0 & 0 & 0 \\ 0 & 0 & 0 & 0 & 0 & 0 \\ 5 & 1 & 0 & 0 & 1 & 1 \\ 5 & 1 & 0 & 0 & 1 & 1 \end{bmatrix} \right\} \cup$$

$$\left\{ \begin{bmatrix} 1 & 2 & 3 & 0 & 1 & 2 & 0 \\ 2 & 4 & 6 & 0 & 2 & 4 & 0 \\ 3 & 6 & 9 & 0 & 3 & 6 & 0 \\ 0 & 0 & 0 & 0 & 0 & 0 & 0 \\ 1 & 2 & 3 & 0 & 1 & 2 & 0 \\ 2 & 4 & 6 & 0 & 2 & 4 & 0 \\ 0 & 0 & 0 & 0 & 0 & 0 & 0 \end{bmatrix} + \begin{bmatrix} 1 & 2 & 1 & 6 & 1 & 0 & 1 \\ 2 & 5 & 2 & 13 & 3 & 1 & 2 \\ 1 & 2 & 1 & 6 & 1 & 0 & 1 \\ 6 & 13 & 6 & 37 & 7 & 1 & 6 \\ 1 & 3 & 1 & 7 & 2 & 1 & 1 \\ 0 & 1 & 0 & 1 & 1 & 1 & 0 \\ 1 & 2 & 1 & 6 & 1 & 0 & 1 \end{bmatrix} + \right.$$

$$\left. + \begin{bmatrix} 2 & 0 & 1 & 2 & 1 & 2 & 1 \\ 0 & 0 & 0 & 0 & 0 & 0 & 0 \\ 1 & 0 & 1 & 2 & 1 & 2 & 0 \\ 2 & 0 & 2 & 4 & 2 & 4 & 0 \\ 1 & 0 & 1 & 2 & 1 & 2 & 0 \\ 2 & 0 & 2 & 4 & 2 & 4 & 0 \\ 1 & 0 & 0 & 0 & 0 & 0 & 1 \end{bmatrix} \right\} =$$

$$\begin{bmatrix} 8 & 4 & 3 & 5 & 3 \\ 4 & 4 & 4 & 7 & 5 \\ 3 & 4 & 12 & 12 & 18 \\ 5 & 7 & 12 & 21 & 20 \\ 3 & 5 & 18 & 20 & 30 \end{bmatrix} \cup \begin{bmatrix} 44 & 10 & 4 & 7 & 21 & 14 \\ 10 & 4 & 0 & 2 & 6 & 2 \\ 4 & 0 & 3 & 4 & 7 & 6 \\ 7 & 2 & 4 & 13 & 18 & 15 \\ 21 & 6 & 7 & 18 & 30 & 22 \\ 14 & 2 & 6 & 15 & 22 & 28 \end{bmatrix} \cup$$

$$\begin{bmatrix} 4 & 4 & 5 & 8 & 3 & 4 & 2 \\ 4 & 9 & 8 & 13 & 5 & 5 & 2 \\ 5 & 8 & 11 & 8 & 5 & 8 & 1 \\ \hline 8 & 13 & 8 & 41 & 9 & 5 & 6 \\ 3 & 5 & 5 & 9 & 4 & 5 & 1 \\ 4 & 5 & 8 & 5 & 5 & 9 & 0 \\ 2 & 2 & 1 & 6 & 1 & 0 & 2 \end{bmatrix} \Bigg\}$$

$$= S_1 \cup S_2 \cup S_3 = S.$$

We see S is a symmetric super trimatrix. Thus using the product of T with $T^T$ we get a symmetric super trimatrix.

Now for the same trimatrix $T = T_1 \cup T_2 \cup T_3$ we find the product

$$\begin{aligned} TT^T &= (T_1 \cup T_2 \cup T_3)^T (T_1 \cup T_2 \cup T_3) \\ &= (T_1^T \cup T_2^T \cup T_3^T)(T_1 \cup T_2 \cup T_3) \\ &= (T_1^T T_1 \cup T_2^T T_2 \cup T_3^T T_3) \end{aligned}$$

$$= \begin{bmatrix} 0 & 1 & 3 & 4 & 5 \\ 1 & 0 & 1 & 0 & 2 \\ 2 & 1 & 0 & 2 & 0 \\ \hline 1 & 1 & 0 & 1 & 0 \\ 1 & 0 & 1 & 0 & 1 \\ 1 & 1 & 1 & 0 & 0 \end{bmatrix} \begin{bmatrix} 0 & 1 & 2 & 1 & 1 & 1 \\ \hline 1 & 0 & 1 & 1 & 0 & 1 \\ \hline 3 & 1 & 0 & 0 & 1 & 1 \\ 4 & 0 & 2 & 1 & 0 & 0 \\ 5 & 2 & 0 & 0 & 1 & 0 \end{bmatrix} \cup$$

$$\begin{bmatrix} 1 & 0 & 1 & 3 & 4 & 5 \\ 2 & 1 & 0 & 1 & 2 & 0 \\ 3 & 1 & 0 & 0 & 1 & 1 \\ \hline 0 & 1 & 0 & 1 & 2 & 0 \\ 1 & 0 & 1 & 0 & 1 & 1 \\ 2 & 0 & 1 & 1 & 2 & 0 \\ 5 & 1 & 0 & 0 & 1 & 1 \end{bmatrix} \begin{bmatrix} 1 & 2 & 3 & 0 & 1 & 2 & 5 \\ \hline 0 & 1 & 1 & 1 & 0 & 0 & 1 \\ 1 & 0 & 0 & 0 & 1 & 1 & 0 \\ \hline 3 & 1 & 0 & 1 & 0 & 1 & 0 \\ 4 & 2 & 1 & 2 & 1 & 2 & 1 \\ 5 & 0 & 1 & 0 & 1 & 0 & 1 \end{bmatrix}$$

$$\cup \begin{bmatrix} 1 & 2 & 3 & 0 & 1 & 2 & 0 \\ \hline 1 & 2 & 1 & 6 & 1 & 0 & 1 \\ 0 & 1 & 0 & 1 & 1 & 1 & 0 \\ 1 & 0 & 0 & 0 & 0 & 0 & 1 \\ 1 & 0 & 1 & 2 & 1 & 2 & 0 \end{bmatrix} \begin{bmatrix} 1 & 1 & 0 & 1 & 1 \\ 2 & 2 & 1 & 0 & 0 \\ 3 & 1 & 0 & 0 & 1 \\ \hline 0 & 6 & 1 & 0 & 2 \\ 1 & 1 & 1 & 0 & 1 \\ 2 & 0 & 1 & 0 & 2 \\ 0 & 1 & 0 & 1 & 0 \end{bmatrix} =$$

$$\left\{ \begin{bmatrix} 0 & 1 \\ 1 & 0 \\ \hline 2 & 1 \\ 1 & 1 \\ 1 & 0 \\ 1 & 1 \end{bmatrix} \begin{bmatrix} 0 & 1 & 2 & 1 & 1 & 1 \\ 1 & 0 & 1 & 1 & 0 & 1 \end{bmatrix} + \begin{bmatrix} 3 & 4 \\ 1 & 0 \\ \hline 0 & 2 \\ 0 & 1 \\ 1 & 0 \\ 1 & 0 \end{bmatrix} \begin{bmatrix} 3 & 1 & 0 & 0 & 1 & 1 \\ 4 & 0 & 2 & 1 & 0 & 0 \end{bmatrix} \right.$$

$$\left. + \begin{bmatrix} \frac{5}{2} \\ 0 \\ \hline 0 \\ 1 \\ 0 \end{bmatrix} \begin{bmatrix} 5 & 2 & 0 & 0 & 1 & 0 \end{bmatrix} \right\} \cup$$

$$\left\{ \begin{bmatrix} 1 \\ 2 \\ 3 \\ \hline 0 \\ 1 \\ 2 \\ \hline 5 \end{bmatrix} \begin{bmatrix} 1 & 2 & 3 \mid 0 & 1 & 2 \mid 5 \end{bmatrix} + \right.$$

$$\begin{bmatrix} 0 & 1 \\ 1 & 0 \\ 1 & 0 \\ \hline 1 & 0 \\ 0 & 1 \\ \hline 0 & 1 \\ 1 & 0 \end{bmatrix} \begin{bmatrix} 0 & 1 & 1 \mid 1 & 0 & 0 \mid 1 \\ 1 & 0 & 0 \mid 0 & 1 & 1 \mid 0 \end{bmatrix} +$$

$$\left. \begin{bmatrix} 3 & 4 & 5 \\ 1 & 2 & 0 \\ 0 & 1 & 1 \\ \hline 1 & 2 & 0 \\ 0 & 1 & 1 \\ \hline 1 & 2 & 0 \\ 0 & 1 & 1 \end{bmatrix} \begin{bmatrix} 3 & 1 & 0 \mid 1 & 0 & 1 \mid 0 \\ 4 & 2 & 1 \mid 2 & 1 & 2 \mid 1 \\ 5 & 0 & 1 \mid 0 & 1 & 0 \mid 1 \end{bmatrix} \right\}$$

$$\cup \left\{ \begin{bmatrix} 1 \\ 1 \\ \hline 0 \\ \hline 1 \\ 1 \end{bmatrix} \begin{bmatrix} 1 \mid 1 & 0 \mid 1 & 1 \end{bmatrix} + \begin{bmatrix} 2 & 3 \\ 2 & 1 \\ 1 & 0 \\ \hline 0 & 0 \\ 0 & 1 \end{bmatrix} \begin{bmatrix} 2 \mid 2 & 1 \mid 0 & 0 \\ 3 \mid 1 & 0 \mid 0 & 1 \end{bmatrix} + \right.$$

$$\left. \begin{bmatrix} \begin{array}{ccc|c} 0 & 1 & 2 & 0 \\ \hline 6 & 1 & 0 & 1 \\ 1 & 1 & 1 & 0 \\ \hline 0 & 0 & 0 & 1 \\ 2 & 1 & 2 & 0 \end{array} \end{bmatrix} \begin{bmatrix} \begin{array}{c|cc|cc} 0 & 6 & 1 & 0 & 2 \\ 1 & 1 & 1 & 0 & 1 \\ 2 & 0 & 1 & 0 & 2 \\ 0 & 1 & 0 & 1 & 0 \end{array} \end{bmatrix} \right\}$$

$$= \left\{ \begin{bmatrix} \begin{array}{c|c|c} [0\ 1]\begin{bmatrix}0\\1\end{bmatrix} & [0\ 1]\begin{bmatrix}1 & 2\\0 & 1\end{bmatrix} & [0\ 1]\begin{bmatrix}1 & 1 & 1\\1 & 0 & 1\end{bmatrix} \\ \hline \begin{bmatrix}1 & 0\\2 & 1\end{bmatrix}\begin{bmatrix}0\\1\end{bmatrix} & \begin{bmatrix}1 & 0\\2 & 1\end{bmatrix}\begin{bmatrix}1 & 2\\0 & 1\end{bmatrix} & \begin{bmatrix}1 & 0\\2 & 1\end{bmatrix}\begin{bmatrix}1 & 1 & 1\\1 & 0 & 1\end{bmatrix} \\ \hline \begin{bmatrix}1 & 1\\1 & 0\\1 & 1\end{bmatrix}\begin{bmatrix}0\\1\end{bmatrix} & \begin{bmatrix}1 & 1\\1 & 0\\1 & 1\end{bmatrix}\begin{bmatrix}1 & 2\\0 & 1\end{bmatrix} & \begin{bmatrix}1 & 1\\1 & 0\\1 & 1\end{bmatrix}\begin{bmatrix}1 & 1 & 1\\1 & 0 & 1\end{bmatrix} \end{array} \end{bmatrix} \right.$$

$$+ \begin{bmatrix} \begin{array}{c|c|c} (3\ 4)\begin{pmatrix}3\\4\end{pmatrix} & (3\ 4)\begin{pmatrix}1 & 0\\0 & 2\end{pmatrix} & (3\ 4)\begin{pmatrix}0 & 1 & 1\\1 & 0 & 0\end{pmatrix} \\ \hline \begin{pmatrix}1 & 0\\0 & 2\end{pmatrix}\begin{pmatrix}3\\4\end{pmatrix} & \begin{pmatrix}1 & 0\\0 & 2\end{pmatrix}\begin{pmatrix}1 & 0\\0 & 2\end{pmatrix} & \begin{pmatrix}1 & 0\\0 & 2\end{pmatrix}\begin{pmatrix}0 & 1 & 1\\1 & 0 & 0\end{pmatrix} \\ \hline \begin{pmatrix}0 & 1\\1 & 0\\1 & 0\end{pmatrix}\begin{pmatrix}3\\4\end{pmatrix} & \begin{pmatrix}0 & 1\\1 & 0\\1 & 0\end{pmatrix}\begin{pmatrix}1 & 0\\0 & 2\end{pmatrix} & \begin{pmatrix}0 & 1\\1 & 0\\1 & 0\end{pmatrix}\begin{pmatrix}0 & 1 & 1\\1 & 0 & 0\end{pmatrix} \end{array} \end{bmatrix}$$

$$\left. + \begin{bmatrix} \begin{array}{c|c|c} (5)(5) & (5)(2\ 0) & (5)(0\ 1\ 0) \\ \hline \begin{pmatrix}2\\0\end{pmatrix}(5) & \begin{pmatrix}2\\0\end{pmatrix}(2\ 0) & \begin{pmatrix}2\\0\end{pmatrix}(0\ 1\ 0) \\ \hline \begin{pmatrix}0\\1\\0\end{pmatrix} & \begin{pmatrix}0\\1\\0\end{pmatrix}(2\ 0) & \begin{pmatrix}0\\1\\0\end{pmatrix}(0\ 1\ 0) \end{array} \end{bmatrix} \right\} \cup$$

$$\left\{\left[\begin{array}{c|c|c}\begin{pmatrix}1\\2\\3\end{pmatrix}(1\ 2\ 3) & \begin{pmatrix}1\\2\\3\end{pmatrix}(0\ 1\ 2) & \begin{pmatrix}1\\2\\3\end{pmatrix}(5)\\\hline\begin{pmatrix}0\\1\\2\end{pmatrix}(1\ 2\ 3) & \begin{pmatrix}0\\1\\2\end{pmatrix}(0\ 1\ 2) & \begin{pmatrix}0\\1\\2\end{pmatrix}(5)\\\hline(5)(1\ 2\ 3) & (5)(0\ 1\ 2) & (5)(5)\end{array}\right]+\right.$$

$$\left[\begin{array}{c|c|c}\begin{pmatrix}0&1\\1&0\\1&0\end{pmatrix}\begin{bmatrix}0&1&1\\1&0&0\end{bmatrix} & \begin{pmatrix}0&1\\1&0\\1&0\end{pmatrix}\begin{pmatrix}1&0&0\\0&1&1\end{pmatrix} & \begin{pmatrix}0&1\\1&0\\1&0\end{pmatrix}\begin{pmatrix}1\\0\end{pmatrix}\\\hline\begin{bmatrix}1&0\\0&1\\0&1\end{bmatrix}\begin{bmatrix}0&1&1\\1&0&0\end{bmatrix} & \begin{bmatrix}1&0\\0&1\\0&1\end{bmatrix}\begin{pmatrix}1&0&0\\0&1&1\end{pmatrix} & \begin{pmatrix}1&0\\0&1\\0&1\end{pmatrix}\begin{pmatrix}1\\0\end{pmatrix}\\\hline(1\ 0)\begin{bmatrix}0&1&1\\1&0&0\end{bmatrix} & (1\ 0)\begin{pmatrix}1&0&0\\0&1&1\end{pmatrix} & (1\ 0)\begin{pmatrix}1\\0\end{pmatrix}\end{array}\right]+$$

$$\left.\left[\begin{array}{c|c|c}\begin{bmatrix}3&4&5\\1&2&0\\0&1&1\end{bmatrix}\begin{bmatrix}3&1&0\\4&2&1\\5&0&1\end{bmatrix} & \begin{bmatrix}3&4&5\\1&2&0\\0&1&1\end{bmatrix}\begin{bmatrix}1&0&1\\2&1&2\\0&1&0\end{bmatrix} & \begin{bmatrix}3&4&5\\1&2&0\\0&1&1\end{bmatrix}\begin{bmatrix}0\\1\\1\end{bmatrix}\\\hline\begin{bmatrix}1&2&0\\0&1&1\\1&2&0\end{bmatrix}\begin{bmatrix}3&1&0\\4&2&1\\5&0&1\end{bmatrix} & \begin{bmatrix}1&2&0\\0&1&1\\1&2&0\end{bmatrix}\begin{bmatrix}1&0&1\\2&1&2\\0&1&0\end{bmatrix} & \begin{bmatrix}1&2&0\\0&1&1\\1&2&0\end{bmatrix}\begin{bmatrix}0\\1\\1\end{bmatrix}\\\hline[0\ 1\ 1]\begin{bmatrix}3&1&0\\4&2&1\\5&0&1\end{bmatrix} & [0\ 1\ 1]\begin{bmatrix}1&0&1\\2&1&2\\0&1&0\end{bmatrix} & [0\ 1\ 1]\begin{bmatrix}0\\1\\1\end{bmatrix}\end{array}\right]\right\}$$

$$\cup \left\{ \left[ \begin{array}{c|c|c} (1)(1) & (1)(1\ 0) & (1)(1\ 1) \\ \hline \begin{bmatrix} 1 \\ 0 \end{bmatrix}[1] & \begin{pmatrix} 1 \\ 0 \end{pmatrix}(1\ 0) & \begin{pmatrix} 1 \\ 0 \end{pmatrix}(1\ 1) \\ \hline \begin{pmatrix} 1 \\ 1 \end{pmatrix}(1) & \begin{pmatrix} 1 \\ 1 \end{pmatrix}(1\ 0) & \begin{pmatrix} 1 \\ 1 \end{pmatrix}(1\ 1) \end{array} \right] \right. +$$

$$\left[ \begin{array}{c|c|c} (2\ 3)\begin{pmatrix} 2 \\ 3 \end{pmatrix} & (2\ 3)\begin{pmatrix} 2 & 1 \\ 1 & 0 \end{pmatrix} & (2\ 3)\begin{pmatrix} 0 & 0 \\ 0 & 1 \end{pmatrix} \\ \hline \begin{pmatrix} 2 & 1 \\ 1 & 0 \end{pmatrix}\begin{pmatrix} 2 \\ 3 \end{pmatrix} & \begin{pmatrix} 2 & 1 \\ 1 & 0 \end{pmatrix}\begin{pmatrix} 2 & 1 \\ 1 & 0 \end{pmatrix} & \begin{pmatrix} 2 & 1 \\ 1 & 0 \end{pmatrix}\begin{pmatrix} 0 & 0 \\ 0 & 1 \end{pmatrix} \\ \hline \begin{pmatrix} 0 & 0 \\ 0 & 1 \end{pmatrix}\begin{pmatrix} 2 \\ 3 \end{pmatrix} & \begin{pmatrix} 0 & 0 \\ 0 & 1 \end{pmatrix}\begin{pmatrix} 2 & 1 \\ 1 & 0 \end{pmatrix} & \begin{pmatrix} 0 & 0 \\ 0 & 1 \end{pmatrix}\begin{pmatrix} 0 & 0 \\ 0 & 1 \end{pmatrix} \end{array} \right]$$

$$+ \left[ \begin{array}{c|c|c} (0\ 1\ 2\ 0)\begin{bmatrix} 0 \\ 1 \\ 2 \\ 0 \end{bmatrix} & (0\ 1\ 2\ 0)\begin{bmatrix} 6 & 1 \\ 1 & 1 \\ 0 & 1 \\ 1 & 0 \end{bmatrix} & (0\ 1\ 2\ 0)\begin{bmatrix} 0 & 2 \\ 0 & 1 \\ 0 & 2 \\ 1 & 0 \end{bmatrix} \\ \hline \begin{pmatrix} 6 & 1 & 0 & 1 \\ 1 & 1 & 1 & 0 \end{pmatrix}\begin{bmatrix} 0 \\ 1 \\ 2 \\ 0 \end{bmatrix} & \begin{pmatrix} 6 & 1 & 0 & 1 \\ 1 & 1 & 1 & 0 \end{pmatrix}\begin{bmatrix} 6 & 1 \\ 1 & 1 \\ 0 & 1 \\ 1 & 0 \end{bmatrix} & \begin{pmatrix} 6 & 1 & 0 & 1 \\ 1 & 1 & 1 & 0 \end{pmatrix}\begin{bmatrix} 0 & 2 \\ 0 & 1 \\ 0 & 2 \\ 1 & 0 \end{bmatrix} \\ \hline \begin{pmatrix} 0 & 0 & 0 & 1 \\ 2 & 1 & 2 & 0 \end{pmatrix}\begin{bmatrix} 0 \\ 1 \\ 2 \\ 0 \end{bmatrix} & \begin{pmatrix} 0 & 0 & 0 & 1 \\ 2 & 1 & 2 & 0 \end{pmatrix}\begin{bmatrix} 6 & 1 \\ 1 & 1 \\ 0 & 1 \\ 1 & 0 \end{bmatrix} & \begin{pmatrix} 0 & 0 & 0 & 1 \\ 2 & 1 & 2 & 0 \end{pmatrix}\begin{bmatrix} 0 & 2 \\ 0 & 1 \\ 0 & 2 \\ 1 & 0 \end{bmatrix} \end{array} \right] \left. \right\}$$

$$= \left\{ \left[ \begin{array}{c|cc|ccc} 1 & 0 & 1 & 1 & 0 & 1 \\ \hline 0 & 1 & 2 & 1 & 1 & 1 \\ 1 & 2 & 5 & 3 & 2 & 3 \\ \hline 1 & 1 & 3 & 2 & 1 & 2 \\ 0 & 1 & 2 & 1 & 1 & 1 \\ 1 & 1 & 3 & 2 & 1 & 2 \end{array} \right] + \left[ \begin{array}{c|cc|ccc} 25 & 3 & 8 & 4 & 3 & 3 \\ \hline 3 & 1 & 0 & 0 & 1 & 1 \\ 8 & 0 & 4 & 2 & 0 & 0 \\ \hline 4 & 0 & 2 & 1 & 0 & 0 \\ 3 & 1 & 0 & 0 & 1 & 1 \\ 3 & 1 & 0 & 0 & 1 & 1 \end{array} \right] + \right.$$

$$\left. \left[ \begin{array}{c|cc|ccc} 25 & 10 & 0 & 0 & 5 & 0 \\ \hline 10 & 4 & 0 & 0 & 2 & 0 \\ 0 & 0 & 0 & 0 & 0 & 0 \\ \hline 0 & 0 & 0 & 0 & 0 & 0 \\ 5 & 2 & 0 & 0 & 1 & 0 \\ 0 & 0 & 0 & 0 & 0 & 0 \end{array} \right] \right\}$$

$$\cup \left\{ \left[ \begin{array}{ccc|ccc|c} 1 & 2 & 3 & 0 & 1 & 2 & 5 \\ 2 & 4 & 6 & 0 & 2 & 4 & 10 \\ 3 & 6 & 9 & 0 & 3 & 6 & 15 \\ \hline 0 & 0 & 0 & 0 & 0 & 0 & 0 \\ 1 & 2 & 3 & 0 & 1 & 2 & 5 \\ 2 & 4 & 6 & 0 & 2 & 4 & 10 \\ \hline 5 & 10 & 15 & 0 & 5 & 10 & 25 \end{array} \right] + \left[ \begin{array}{ccc|ccc|c} 1 & 0 & 0 & 0 & 1 & 1 & 0 \\ 0 & 1 & 1 & 1 & 0 & 0 & 1 \\ 0 & 1 & 1 & 1 & 0 & 0 & 1 \\ \hline 0 & 1 & 1 & 1 & 0 & 0 & 1 \\ 1 & 0 & 0 & 0 & 1 & 1 & 0 \\ 1 & 0 & 0 & 0 & 1 & 1 & 0 \\ \hline 0 & 1 & 1 & 1 & 0 & 0 & 1 \end{array} \right] + \right.$$

$$\left. \left[ \begin{array}{ccc|ccc|c} 50 & 11 & 9 & 11 & 9 & 11 & 9 \\ 11 & 5 & 2 & 5 & 2 & 5 & 2 \\ 9 & 2 & 2 & 2 & 2 & 2 & 2 \\ \hline 11 & 5 & 2 & 5 & 2 & 5 & 2 \\ 9 & 2 & 2 & 2 & 2 & 2 & 2 \\ 11 & 5 & 2 & 5 & 2 & 5 & 2 \\ \hline 9 & 2 & 2 & 2 & 2 & 2 & 2 \end{array} \right] \right\} \cup$$

$$\left\{\begin{bmatrix} 1 & 1 & 0 & 1 & 1 \\ 1 & 1 & 0 & 1 & 1 \\ 0 & 0 & 0 & 0 & 0 \\ 1 & 1 & 0 & 1 & 1 \\ 1 & 1 & 0 & 1 & 1 \end{bmatrix} + \begin{bmatrix} 13 & 7 & 2 & 0 & 3 \\ 7 & 5 & 2 & 0 & 1 \\ 2 & 2 & 1 & 0 & 0 \\ 0 & 0 & 0 & 0 & 0 \\ 3 & 1 & 0 & 0 & 1 \end{bmatrix} + \begin{bmatrix} 5 & 1 & 3 & 0 & 5 \\ 1 & 38 & 7 & 1 & 13 \\ 3 & 7 & 3 & 0 & 5 \\ 0 & 1 & 0 & 1 & 0 \\ 5 & 13 & 5 & 0 & 9 \end{bmatrix}\right\}$$

$$= \begin{bmatrix} 51 & 13 & 9 & 5 & 8 & 4 \\ 13 & 6 & 2 & 1 & 4 & 1 \\ 9 & 2 & 9 & 5 & 2 & 3 \\ 5 & 1 & 5 & 3 & 1 & 2 \\ 8 & 4 & 2 & 1 & 1 & 2 \\ 4 & 1 & 3 & 2 & 2 & 3 \end{bmatrix} \cup \begin{bmatrix} 52 & 13 & 12 & 11 & 11 & 14 & 14 \\ 13 & 11 & 9 & 6 & 4 & 9 & 13 \\ 12 & 9 & 12 & 3 & 5 & 8 & 18 \\ 11 & 6 & 3 & 6 & 2 & 5 & 3 \\ 11 & 4 & 5 & 2 & 4 & 5 & 7 \\ 14 & 9 & 8 & 5 & 5 & 10 & 12 \\ 14 & 13 & 18 & 3 & 7 & 12 & 28 \end{bmatrix}$$

$$\cup \begin{bmatrix} 19 & 9 & 5 & 1 & 9 \\ 9 & 44 & 9 & 2 & 15 \\ 5 & 9 & 4 & 0 & 5 \\ 1 & 2 & 0 & 2 & 1 \\ 9 & 15 & 5 & 1 & 11 \end{bmatrix} = P_1 \cup P_2 \cup P_3 = P.$$

We see P is also a symmetric super trimatrix. However we see P and S i.e., $TT^T$ and $T^TT$ have no relation.

Thus for a given super trimatrix T we can obtain two symmetric super trimatrices. Now we proceed on to define the notion of semi super trimatrix and the types of semi super trimatrices.

**DEFINITION 3.17:** *Let $T = T_1 \cup T_2 \cup T_3$, where at least one of the $T_i$ is a supermatrix and at least one of the $T_j$ $(i \neq j)$ is just a matrix and not a supermatrix $(1 \leq i, j \leq 3)$ then we call T to be a semi super trimatrix.*

***Example 3.41:*** Let $T = T_1 \cup T_2 \cup T_3$ where $T_1 = [\,0\ 1\ 2\ 3\ 4\,]$; $T_2 = [2\ 1\ |\ 0\ 5\ 7\ |\ 1\ 1\ 1\ 3\,]$ and $T_3 = [\,1\ 1\ 0\ 3\ |\ 8\ 9\ 3\ |\ 1\ 2\ 5\ 7\ 1]$. T is a semi super trimatrix for $T_1$ is just a row vector where as $T_2$ and $T_3$ are super row vectors.

***Example 3.42:*** Let $S = S_1 \cup S_2 \cup S_3$ where

$$S_1 = \begin{bmatrix} 3 \\ 1 \\ \hline 0 \\ 1 \\ 1 \\ \hline 0 \end{bmatrix},\ S_2 = \begin{bmatrix} 2 \\ 3 \\ 4 \\ 5 \\ 6 \end{bmatrix}$$

and

$$S_3 = \begin{bmatrix} 1 \\ 0 \\ 1 \\ 1 \\ \hline 1 \\ 1 \\ 1 \end{bmatrix}.$$

S is semi super trimatrix as $S_1$ and $S_3$ are column super vectors where as $S_2$ is just a column vector.

***Example 3.43:*** Let $Q = Q_1 \cup Q_2 \cup Q_3$ where

$$Q_1 = [1\ 1\ 1\ 1\ 1\ 0],$$
$$Q_2 = [2\ 3\ 0\ 4\ 1\ 5\ 7]$$

and

$$Q_3 = [3\ 1\ |\ 7\ 0\ 5\ |\ 1\ 1\ 4\ 3\ |\ 0];$$

Q is a semi super trimatrix for $Q_1$ and $Q_2$ are just row vectors where as $Q_3$ is a row super vector.

*Example 3.44:* Let $V = V_1 \cup V_2 \cup V_3$ where

$$V_1 = \begin{bmatrix} 1 \\ 0 \\ 1 \\ 2 \end{bmatrix}, V_2 = \begin{bmatrix} 1 \\ 1 \\ 0 \\ \hline 1 \\ 2 \\ 5 \end{bmatrix} \text{ and } V_3 = \begin{bmatrix} 1 \\ 0 \\ 2 \\ 1 \end{bmatrix};$$

V is a semi super trimatrix as $V_1$ and $V_3$ are just column vectors where as only $V_2$ is a super column vector.

*Example 3.45:* Let $V = V_1 \cup V_2 \cup V_3$ where

$$V_1 = \begin{bmatrix} 2 & 3 \\ 1 & 5 \end{bmatrix}, V_2 = \begin{bmatrix} 3 & 1 & 0 & 2 \\ 1 & 1 & 5 & 6 \\ 3 & 4 & 8 & 1 \\ 1 & 2 & 3 & 4 \end{bmatrix} \text{ and } V_3 = \begin{bmatrix} 1 & 0 & 2 \\ 2 & 0 & 1 \\ 1 & 1 & 0 \end{bmatrix}.$$

V is a semi super trimatrix as $V_1$ and $V_3$ are just square matrices where as $V_2$ is a square supermatrix.

*Example 3.46:* Let $W = W_1 \cup W_2 \cup W_3$ where

$$W_1 = \begin{bmatrix} 1 & 2 & 3 \\ 0 & 1 & 2 \\ 5 & 6 & 7 \end{bmatrix}, W_2 = \begin{bmatrix} 3 & 1 & 0 & 2 \\ 1 & 1 & 5 & 6 \\ 3 & 4 & 8 & 1 \\ 1 & 2 & 3 & 4 \end{bmatrix}$$

and

$$W_3 = \begin{bmatrix} 1 & 0 & 2 \\ 2 & 0 & 1 \\ 1 & 1 & 0 \end{bmatrix}.$$

W is a semi super trimatrix as $W_1$ and $W_3$ are just square matrices where as $W_2$ is a square supermatrix.

*Example 3.47:* Let $W = W_1 \cup W_2 \cup W_3$ where

$$W_1 = \begin{bmatrix} 1 & 2 & 3 \\ 0 & 1 & 2 \\ \hline 5 & 6 & 7 \end{bmatrix} \quad W_2 = \begin{bmatrix} 2 & 1 \\ 0 & 5 \end{bmatrix}$$

and

$$W_3 = \begin{bmatrix} 1 & 2 & 3 & 4 & 5 \\ 6 & 7 & 8 & 9 & 0 \\ \hline 0 & 9 & 8 & 7 & 6 \\ 5 & 4 & 3 & 2 & 1 \\ \hline 1 & 1 & 0 & 1 & 1 \end{bmatrix}.$$

W is a semi super trimatrix as $W_2$ is just a matrix where as $W_1$ and $W_3$ are supermatrices.

*Example 3.48:* Let $S = S_1 \cup S_2 \cup S_3$ where

$$S_1 = \begin{bmatrix} 3 & 1 & 0 & 1 & 1 \\ 1 & 2 & 5 & 6 & 0 \\ 3 & 6 & 7 & 1 & 2 \end{bmatrix}, S_2 = \begin{bmatrix} 3 & 1 & 2 \\ 0 & 1 & 5 \\ \hline 3 & 1 & 1 \\ 1 & 1 & 1 \\ 5 & 6 & 7 \end{bmatrix}$$

and

$$S_3 = \begin{bmatrix} 2 & 1 & 3 & 4 & 5 & 6 & 7 \\ 1 & 2 & 3 & 4 & 8 & 9 & 0 \end{bmatrix}.$$

S is semi super trimatrix, as $S_1$ is just a rectangular matrix but $S_2$ and $S_3$ are rectangular supermatrices.

*Example 3.49:* Let $T = T_1 \cup T_2 \cup T_3$ where

$$T_1 = \begin{bmatrix} 3 & 1 & 2 \\ 5 & 6 & 7 \\ 8 & 9 & 0 \end{bmatrix}, T_2 = \left[ \begin{array}{ccc|cc} 1 & 1 & 7 & 3 & 8 \\ 4 & 0 & 6 & 2 & 0 \\ \hline 6 & 9 & 5 & 1 & 5 \\ 7 & 8 & 4 & 1 & 6 \end{array} \right]$$

and

$$T_3 = \left[ \begin{array}{cc|cc} 1 & 2 & 3 & 4 \\ 5 & 6 & 7 & 8 \\ \hline 9 & 0 & 1 & 2 \\ 3 & 4 & 5 & 6 \end{array} \right].$$

T is a semi super trimatrix as $T_1$ is a square matrix, $T_2$ and $T_3$ are supermatrices.

Just as in case of semi superbimatrix we can in case of semi super trimatrices also define 5 types of semi super trimatrices. If in the semi super trimatrix $T_1 \cup T_2 \cup T_3$ all the 3 matrices are just column vectors we call T to be a semi super column vector.

Examples 3.42 and 3.44 are semi super column trimatrices. If in the row semi super trimatrix $V = V_1 \cup V_2 \cup V_3$, $V_i$'s are matrices $1 \le i \le 3$; some of them super row vector; then we call V to be a semi super row trivector. The examples 3.41 and 3.43 are semi super row trivector. Let $V = V_1 \cup V_2 \cup V_3$ be a semi super trimatrix where each of the $V_i$ is an n × n square matrices $1 \le i \le 3$, some of which are super square matrices and others just square matrices. We call $V = V_1 \cup V_2 \cup V_3$ to be a semi square super trimatrix.

We give the following example.

***Example 3.50:*** Let $U = U_1 \cup U_2 \cup U_3$ where

$$U_1 = \left[ \begin{array}{ccc|c} 1 & 2 & 3 & 4 \\ 5 & 6 & 7 & 8 \\ 9 & 0 & 1 & 2 \\ \hline 3 & 4 & 5 & 6 \end{array} \right],$$

$$U_2 = \begin{bmatrix} 0 & 1 & 0 & 1 \\ 1 & 0 & 1 & 0 \\ 1 & 1 & 1 & 1 \\ 0 & 1 & 0 & 1 \end{bmatrix}$$

and

$$U_3 = \left[ \begin{array}{cc|cc} 1 & 2 & 3 & 0 \\ 0 & 1 & 2 & 0 \\ \hline 3 & 1 & 0 & 3 \\ 1 & 1 & 0 & 5 \end{array} \right]$$

be a semi super trimatrix. U is a 4 × 4 square semi super trimatrix.

Next we consider a semi super trimatrix $T = T_1 \cup T_2 \cup T_3$, where $T_i$'s are square matrices of different order; some of the square matrices are super square matrices. We call T to be a mixed semi super square matrix. The examples 3.45 and 3.46 are mixed semi super square trimatrices.

Now we proceed on to define the notion of a mixed semi super rectangular trimatrix and a semi super rectangular trimatrix.

Let $V = V_1 \cup V_2 \cup V_3$ where $V_i$'s are m × n (m ≠ n) rectangular matrices some just ordinary and other supermatrices. Thus we call V to be an m × n rectangular semi super trimatrix.

We illustrate this by a simple example.

*Example 3.51:* Let $M = M_1 \cup M_2 \cup M_3$ where

$$M_1 = \left[ \begin{array}{cc|ccccc} 3 & 1 & 5 & 2 & 0 & 1 & 3 \\ 0 & 0 & 1 & 0 & 2 & 0 & 1 \\ \hline 1 & 5 & 6 & 2 & 0 & 5 & 0 \end{array} \right],$$

$$M_2 = \begin{bmatrix} 1 & 2 & 1 & 0 & 2 & 0 & 1 \\ 0 & 1 & 0 & 1 & 0 & 2 & 1 \\ 2 & 0 & 2 & 2 & 1 & 1 & 0 \end{bmatrix}$$

and

$$M_3 = \left[ \begin{array}{cccc|ccc} 1 & 0 & 1 & 1 & 1 & 1 & 0 \\ 1 & 1 & 1 & 1 & 1 & 0 & 0 \\ 0 & 1 & 1 & 1 & 0 & 1 & 1 \end{array} \right]$$

be a semi super trimatrix. We see M is a $3 \times 7$ rectangular semi super trimatrix.

A semi super trimatrix $V = V_1 \cup V_2 \cup V_3$ is said to be a mixed rectangular semi super trimatrix if $V_i$'s are rectangular matrices or rectangular supermatrices of different orders. The semi supermatrices given in the example 3.48 is a mixed rectangular semi super trimatrix.

A semi super trimatrix W is said to be a mixed semi super trimatrix if in $W = W_1 \cup W_2 \cup W_3$ some of the matrices $W_i$'s are square matrices or square supermatrices and some of the $W_j$'s are rectangular matrices or rectangular supermatrices.

*Example 3.52:* Let $T = T_1 \cup T_2 \cup T_3$ where

$$T_1 = \left[ \begin{array}{c|ccc} 2 & 1 & 0 & 3 \\ \hline 5 & 0 & 2 & 1 \\ 0 & 1 & 0 & 2 \\ 1 & 0 & 2 & 0 \end{array} \right],$$

$$T_2 = \begin{bmatrix} 3 & 0 & 1 & 4 & 7 & 1 \\ 1 & 0 & 2 & 5 & 8 & 4 \\ 2 & 7 & 3 & 6 & 9 & 2 \end{bmatrix}$$

and

$$T_3 = \left[\begin{array}{cccc} 3 & 1 & 4 & 5 \\ 0 & 1 & 2 & 3 \\ 4 & 5 & 6 & 7 \\ \hline 8 & 9 & 1 & 0 \\ 1 & 1 & 1 & 2 \\ \hline 1 & 3 & 1 & 4 \\ 1 & 5 & 1 & 6 \end{array}\right];$$

T is a mixed semi super trimatrix.

*Example 3.53:* Let $S = S_1 \cup S_2 \cup S_3$ where

$$S_1 = \left[\begin{array}{cc|cccccc} 3 & 1 & 5 & 9 & 2 & 1 & 3 \\ 0 & 2 & 6 & 1 & 0 & 1 & 4 \\ \hline 1 & 3 & 7 & 1 & 1 & 1 & 5 \\ 2 & 4 & 8 & 1 & 2 & 1 & 6 \end{array}\right],$$

$$S_2 = \left[\begin{array}{ccccc} 3 & 1 & 3 & 5 & 2 \\ 1 & 3 & 0 & 1 & 1 \\ 2 & 0 & 1 & 0 & 1 \\ 3 & 1 & 1 & 1 & 0 \\ 0 & 2 & 5 & 2 & 6 \end{array}\right]$$

and

$$S_3 = \left[\begin{array}{c|cc} 1 & 2 & 3 \\ 1 & 0 & 2 \\ \hline 2 & 0 & 1 \\ 0 & 1 & 0 \\ 5 & 6 & 7 \\ 8 & 9 & 2 \\ 1 & 0 & 3 \end{array}\right].$$

S is a mixed semi super trimatrix.

Now having seen the 5 types of semi super trimatrices now we define the notion of semi super trimatrix which is symmetric and a quasi symmetric semi super trimatrix.

**Example 3.54:** Let $T = T_1 \cup T_2 \cup T_3$ where T is a semi super trimatrix. We say T is a symmetric semi super trimatrix if each of the $T_i$ is a symmetric matrix or a symmetric supermatrix.

**Example 3.55:** Let $V = V_1 \cup V_2 \cup V_3$ where

$$V_1 = \begin{bmatrix} 2 & 0 & | & 1 & 2 & 3 \\ 0 & 5 & | & 3 & 2 & 1 \\ \hline 1 & 3 & | & 7 & 5 & 6 \\ 2 & 2 & | & 5 & 8 & 9 \\ 3 & 1 & | & 6 & 9 & 0 \end{bmatrix},$$

$$V_2 = \begin{bmatrix} 1 & 2 & | & 3 & 4 & 5 & 6 \\ 2 & 0 & | & 9 & 8 & 7 & 5 \\ \hline 3 & 9 & | & 2 & 7 & 6 & 4 \\ 4 & 8 & | & 7 & 3 & 1 & 2 \\ 5 & 7 & | & 6 & 1 & 7 & 8 \\ \hline 6 & 5 & | & 4 & 2 & 8 & 1 \end{bmatrix}$$

and

$$V_3 = \begin{bmatrix} 1 & 2 & 3 & 0 \\ 2 & 5 & 0 & 1 \\ 3 & 0 & 7 & 2 \\ 0 & 1 & 2 & 6 \end{bmatrix}$$

is a symmetric semi super trimatrix as $V_1$ and $V_2$ are symmetric supermatrices and $V_3$ is just a symmetric matrix.

**Example 3.56:** Let $T = T_1 \cup T_2 \cup T_3$ where

$$T_1 = \begin{bmatrix} 0 & -1 & 1 & -1 & 1 \\ -1 & 2 & 0 & 1 & 0 \\ 1 & 0 & 5 & 1 & -1 \\ -1 & 1 & 1 & 7 & 0 \\ 1 & 0 & -1 & 0 & 8 \end{bmatrix},$$

$$T_2 = \begin{bmatrix} 3 & 0 & 8 \\ 0 & 5 & 6 \\ 8 & 6 & 3 \end{bmatrix}$$

and

$$T_3 = \left[\begin{array}{cccc|cc} 0 & 1 & 2 & 3 & 4 & 5 \\ 1 & 2 & 0 & 1 & 0 & 1 \\ 2 & 0 & 7 & 5 & 6 & 2 \\ 3 & 1 & 5 & 0 & 7 & 1 \\ \hline 4 & 0 & 6 & 7 & 9 & 0 \\ 5 & 1 & 2 & 1 & 0 & 8 \end{array}\right].$$

Here $T_1$ and $T_2$ are symmetric matrices where as $T_3$ is a symmetric supermatrix. Thus T is a symmetric semi supermatrix.

Now we proceed on to define the notion of quasi symmetric semi supermatrix and illustrate it by some simple examples.

**DEFINITION 3.18:** *Let $T = T_1 \cup T_2 \cup T_3$ be a semi super trimatrix. If one of the $T_i$'s is symmetric supermatrix and one of the $T_j$'s is a symmetric matrix then we call T to be a quasi symmetric semi super trimatrix; i.e., the other $T_k$ can be a square matrix which is not symmetric or a square supermatrix which is not symmetric or a rectangular supermatrix or a rectangular ordinary matrix ($1 \leq i, j, k \leq 3$).*

*Example 3.57:* Let $T = T_1 \cup T_2 \cup T_3$ be a semi super trimatrix where

$$T_1 = \begin{bmatrix} 3 & 0 & 8 & 4 & 5 \\ 0 & 1 & 2 & 3 & 4 \\ \hline 8 & 2 & -1 & 5 & 1 \\ 4 & 3 & 5 & 0 & 2 \\ 5 & 4 & 1 & 2 & 8 \end{bmatrix},$$

$$T_2 = \begin{bmatrix} 3 & 1 & 2 \\ \hline 0 & 4 & 3 \\ 5 & 0 & 1 \\ 4 & 8 & 0 \\ 1 & 1 & 2 \\ \hline 3 & 0 & 1 \\ 6 & 8 & 5 \end{bmatrix}$$

and

$$T_3 = \begin{bmatrix} 3 & 4 & 0 & 1 \\ 4 & 8 & 5 & 6 \\ \hline 0 & 5 & 1 & 2 \\ 1 & 6 & 2 & 7 \end{bmatrix}.$$

T is a quasi symmetric super trimatrix.

*Example 3.58:* Let $V = V_1 \cup V_2 \cup V_3$ where

$$V_1 = \begin{bmatrix} 3 & 0 & 1 & 2 \\ 0 & 5 & 6 & 3 \\ 1 & 6 & 7 & 1 \\ 2 & 3 & 1 & 8 \end{bmatrix},$$

$$V_2 = \begin{bmatrix} 3 & 1 & 0 & 1 \\ 1 & 0 & 2 & 3 \\ 1 & 6 & 2 & 1 \\ 7 & 5 & 4 & 3 \\ 1 & 1 & 2 & 1 \\ 0 & 3 & 5 & 1 \\ 6 & 8 & 3 & 4 \\ 1 & 1 & 3 & 8 \\ 0 & 1 & 0 & 1 \end{bmatrix}$$

and

$$V_3 = \left[ \begin{array}{ccc|ccc} 1 & 2 & 3 & 4 & 5 & 6 \\ 2 & 9 & 8 & 7 & 6 & 5 \\ 3 & 8 & 1 & 0 & 2 & 4 \\ \hline 4 & 7 & 0 & 9 & 0 & 1 \\ 5 & 6 & 0 & 0 & 3 & 4 \\ 6 & 5 & 4 & 1 & 4 & 7 \end{array} \right],$$

we see $V_1$ is a symmetric square matrix. $V_3$ is a symmetric square supermatrix where as $V_2$ is just a rectangular matrix. Thus V is a quasi symmetric semi super trimatrix.

*Example 3.59:* Let $S = S_1 \cup S_2 \cup S_3$ where $S_1$ is a square matrix given by

$$S_1 = \begin{bmatrix} 9 & 8 & 7 & 6 & 4 & 3 \\ 2 & 1 & 0 & 1 & 2 & 3 \\ 4 & 5 & 6 & 7 & 8 & 9 \\ 1 & 0 & 1 & 5 & 2 & 3 \\ 4 & 5 & 3 & 1 & 0 & 1 \\ 0 & 7 & 0 & 3 & 6 & 2 \end{bmatrix},$$

$$S_2 = \begin{bmatrix} 0 & 3 & 0 & 1 & 5 \\ 3 & 2 & 1 & 0 & 2 \\ 0 & 1 & 7 & 9 & 3 \\ \hline 1 & 0 & 9 & 1 & 0 \\ 5 & 2 & 3 & 0 & 7 \end{bmatrix}$$

and

$$S_3 = \begin{bmatrix} 3 & 4 & 5 & 1 \\ 4 & 2 & 0 & 1 \\ 5 & 0 & 3 & 11 \\ 1 & 1 & 1 & 0 \end{bmatrix}$$

is not a quasi symmetric semi trimatrix as $S_2$ is a symmetric as a matrix but not a symmetric supermatrix because of the partition. $S_3$ is a symmetric matrix and $S_1$ just a non symmetric square matrix.

*Example 3.60:* Let $T = T_1 \cup T_2 \cup T_3$ where

$$T_1 = \begin{bmatrix} 0 & 8 & 9 & 7 & 6 \\ \hline 8 & 1 & 5 & 6 & 7 \\ 9 & 5 & 2 & 3 & 4 \\ 7 & 6 & 3 & 5 & 7 \\ \hline 6 & 7 & 4 & 7 & 8 \end{bmatrix},$$

$$T_2 = \begin{bmatrix} 0 & 1 & 2 & 3 & 0 & 5 \\ 1 & 6 & 5 & 4 & 3 & 2 \\ 2 & 5 & 7 & 1 & 0 & 9 \\ 3 & 4 & 1 & 8 & 6 & 1 \\ 0 & 3 & 0 & 6 & 5 & 8 \\ 5 & 2 & 9 & 1 & 8 & 3 \end{bmatrix}$$

and

$$T_3 = \begin{bmatrix} 2 & 1 & 0 \\ 1 & 8 & 2 \\ 5 & 0 & 9 \end{bmatrix}.$$

T is a quasi symmetric semi super trimatrix.

Now having seen examples of quasi symmetric semi super trimatrix we now proceed on to define some of its properties.

A matrix T of the form

$$T = \begin{bmatrix} 3 & 1 & 5 \\ \hline 1 & 0 & 1 \\ 1 & 1 & 1 \\ 2 & 6 & 2 \\ 0 & 1 & 2 \\ \hline 5 & 6 & 7 \\ 3 & 1 & 2 \end{bmatrix}$$

will be known as a special super column vector. Likewise

$$S = \begin{bmatrix} 3 & 1 & 4 & 7 & 1 & 0 & 7 & 1 & 1 \\ 1 & 0 & 5 & 8 & 0 & 4 & 8 & 2 & 0 \\ 5 & 2 & 0 & 9 & 2 & 5 & 0 & 3 & 0 \\ 6 & 3 & 6 & 0 & 3 & 6 & 9 & 4 & 3 \end{bmatrix}$$

is a supermatrix which will be known as the special super row vector. We see these matrices are partitioned either horizontally or vertically never both vertically and horizontally.

We see the matrix

$$S = \begin{bmatrix} 2 & 3 & 4 & | & 0 \\ 1 & 1 & 5 & | & 3 \\ \hline 0 & 1 & 1 & | & 8 \\ 1 & 1 & 1 & | & 0 \\ 2 & 3 & 4 & | & 5 \\ \hline 6 & 7 & 8 & | & 9 \\ 1 & 0 & 8 & | & 0 \end{bmatrix}$$

is not a column super vector for it is divided both vertically and horizontally.

Thus a special column super vector is always divided or partitioned only horizontally and a special super row vector is always partitioned only vertically. Thus

$$P = \begin{bmatrix} 0 & 2 & | & 1 & 0 & 9 & 8 & | & 6 & 4 & 3 \\ 1 & 1 & | & 1 & 3 & 2 & 1 & | & 3 & 5 & 7 \\ 2 & 0 & | & 2 & 5 & 9 & 0 & | & 4 & 3 & 2 \end{bmatrix}$$

is not a special super row vector.

Now we give conditions under which the product of a special column super vector is compatible with a special row super vector and so on.

We define both the minor product as well as the major product of these supermatrices when specially these supermatrices are super trimatrices.

**DEFINITION 3.19:** *Let $T = T_1 \cup T_2 \cup T_3$ be a semi super trimatrix if at least one of the $T_i$ is a column super vector and one of the $T_j$'s is a simple column vector i.e., $T_j$ is a $m \times n$ matrix with $m > n$ then we call T to be a special column semi super trimatrix or vector ($1 \leq i, j \leq 3$).*

We first illustrate this by a simple example before we define more concepts.

***Example 3.61:*** Let $T = T_1 \cup T_2 \cup T_3$ where

$$T_1 = \left[\begin{array}{cc|cccc|cc} 3 & 1 & 5 & 7 & 9 & 1 & 1 & 0 \\ 2 & 0 & 6 & 8 & 10 & 1 & 2 & 9 \end{array}\right],$$

$$T_2 = \left[\begin{array}{ccccc} 3 & 1 & 2 & 5 & 1 \\ 0 & 1 & 2 & 3 & 4 \\ \hline 5 & 6 & 7 & 8 & 9 \\ 9 & 8 & 7 & 6 & 5 \\ 4 & 3 & 2 & 1 & 0 \\ 1 & 1 & 0 & 1 & 1 \\ \hline 0 & 1 & 1 & 1 & 0 \\ 3 & 0 & 0 & 2 & 1 \end{array}\right]$$

and

$$T_3 = \left[\begin{array}{ccc} 2 & 1 & 5 \\ 6 & 7 & 8 \\ 9 & 0 & 1 \\ 2 & 3 & 4 \\ 6 & 7 & 8 \\ 9 & 1 & 3 \end{array}\right].$$

T is a special column semi super trivector. Clearly we see all the 3 matrices need not be an m × n matrix with n > m.

***Example 3.62:*** Let $S = S_1 \cup S_2 \cup S_3$ where

$$S_1 = \left[\begin{array}{c|ccc} 3 & 1 & 2 & 5 \\ 6 & 0 & 1 & 1 \\ 0 & 7 & 2 & 3 \\ 0 & 1 & 1 & 5 \end{array}\right],$$

$$S_2 = \begin{bmatrix} 3 & 1 & 2 \\ 0 & 1 & 1 \\ \hline 6 & 1 & 1 \\ 1 & 0 & 1 \\ 1 & 1 & 0 \\ 3 & 2 & 5 \\ \hline 6 & 7 & 8 \\ 1 & 3 & 4 \end{bmatrix}$$

and

$$S_3 = \begin{bmatrix} 3 & 1 & 2 & 3 & 5 \\ 1 & 0 & 2 & 6 & 1 \\ 0 & 1 & 0 & 7 & 8 \\ 8 & 7 & 0 & 0 & 7 \\ 1 & 1 & 0 & 8 & 4 \\ 1 & 2 & 3 & 4 & 5 \\ 6 & 7 & 8 & 9 & 0 \\ 1 & 3 & 0 & 1 & 1 \\ 8 & 1 & 9 & 1 & 0 \end{bmatrix};$$

S is a special semi super column trivector or matrix we see $S_1$ is a square supermatrix but $S_2$ is a special column super trivector and $S_3$ is an m × n matrix with m > n.

**DEFINITION 3.20:** *Let $P = P_1 \cup P_2 \cup P_3$ where at least one of the $P_i$'s is a special row super vector; at least one of the $P_j$'s is a special row vector i.e., $P_j$ is an m × n matrix in which n > m. 1 ≤ i, j ≤ 3.*

*Then we call P to be a special semi super row trivector.*

We illustrate this by the following examples.

**Example 3.63:** Let $P = P_1 \cup P_2 \cup P_3$ where

$$P_1 = \begin{bmatrix} 3 & 1 & 2 & 0 \\ 1 & 9 & 8 & 1 \\ 2 & 7 & 6 & 3 \\ 5 & 3 & 2 & 8 \end{bmatrix},$$

$$P_2 = \begin{bmatrix} 0 & 1 & 2 & 3 & 4 & 5 & 6 & 7 & 8 & 9 \\ 1 & 0 & 1 & 2 & 0 & 7 & 0 & 6 & 7 & 6 \\ 3 & 1 & 1 & 0 & 6 & 1 & 1 & 0 & 0 & 1 \end{bmatrix}$$

and

$$P_3 = \begin{bmatrix} 3 & 1 & 7 & 9 & 3 & 4 & 3 & 1 & 0 \\ 5 & 3 & 8 & 6 & 8 & 1 & 5 & 1 & 7 \end{bmatrix}.$$

P is a special semi super row trivector.

**Example 3.64:** Let $S = S_1 \cup S_2 \cup S_3$ where

$$S_2 = \begin{bmatrix} 3 & 4 & 6 & 8 & 0 \\ 1 & 5 & 7 & 9 & 1 \end{bmatrix},$$

$$S_1 = \begin{bmatrix} 3 & 1 & 2 & 5 \\ 0 & 1 & 0 & 1 \\ 7 & 6 & 2 & 5 \\ 3 & 0 & 1 & 4 \\ 1 & 1 & 0 & 3 \\ 1 & 0 & 1 & 1 \\ 1 & 1 & 1 & 1 \end{bmatrix}$$

and

$$S_3 = \begin{bmatrix} 3 & 9 & 8 & 7 & 6 & 3 & 4 & 2 & 3 & 4 & 0 & 1 \\ 1 & 1 & 2 & 3 & 4 & 5 & 6 & 7 & 8 & 9 & 8 & 7 \\ 2 & 6 & 5 & 4 & 3 & 2 & 1 & 0 & 1 & 2 & 3 & 1 \\ 4 & 7 & 8 & 4 & 0 & 1 & 2 & 3 & 5 & 1 & 0 & 1 \end{bmatrix}.$$

S is a special semi super row trimatrix or trivector, as $S_2$ and $S_3$ are rectangular m × n matrices with m > n.

Now we proceed on to illustrate by examples the products of semi super trimatrices.

***Example 3.65:*** Let $T = T_1 \cup T_2 \cup T_3$ be a special semi super row trivector and $P = P_1 \cup P_2 \cup P_3$ a special semi super column trivector. Here $T = T_1 \cup T_2 \cup T_3$

$$= \begin{bmatrix} 3 & 1 & 0 \\ 5 & 4 & 1 \\ 0 & 1 & 0 \end{bmatrix} \cup \begin{bmatrix} 2 & 0 & 5 & 7 & 3 & 1 & 4 \\ 1 & 2 & 1 & 1 & 0 & 0 & 2 \\ 3 & 1 & 2 & 0 & 1 & 1 & 0 \end{bmatrix} \cup$$

$$\begin{bmatrix} 1 & 3 & 5 & 1 & 0 & 0 & 1 & 1 \\ 2 & 4 & 0 & 1 & 1 & 0 & 0 & 0 \end{bmatrix}$$

is special row trimatrix.

$$P = P_1 \cup P_2 \cup P_3$$

$$= \begin{bmatrix} 0 & 1 & 1 \\ 1 & 0 & 1 \\ 1 & 1 & 0 \end{bmatrix} \cup \begin{bmatrix} 1 & 0 \\ 0 & 1 \\ 2 & 1 \\ 1 & 0 \\ 1 & 1 \\ 0 & 1 \\ 2 & 0 \end{bmatrix} \cup \begin{bmatrix} 1 & 0 & 1 \\ 0 & 1 & 0 \\ 1 & 1 & 1 \\ 0 & 1 & 2 \\ 1 & 0 & 0 \\ 0 & 1 & 1 \\ 1 & 0 & 1 \\ 0 & 0 & 0 \end{bmatrix}$$

is a special semi super column trimatrix.

$$
\begin{aligned}
TP &= (T_1 \cup T_2 \cup T_3)(P_1 \cup P_2 \cup P_3) \\
&= T_1P_1 \cup T_2P_2 \cup T_3P_3
\end{aligned}
$$

$$
= \begin{bmatrix} 3 & 1 & 0 \\ 5 & 4 & 1 \\ 0 & 1 & 0 \end{bmatrix} \begin{bmatrix} 0 & 1 & 1 \\ 1 & 0 & 1 \\ 1 & 1 & 0 \end{bmatrix} \cup
$$

$$
\begin{bmatrix} 2 & 0 & | & 5 & | & 7 & 3 & 1 & 4 \\ 1 & 2 & | & 1 & | & 1 & 0 & 0 & 2 \\ 3 & 1 & | & 2 & | & 0 & 1 & 1 & 0 \end{bmatrix} \begin{bmatrix} 1 & 0 \\ 0 & 1 \\ \hline 2 & 1 \\ \hline 1 & 0 \\ 1 & 1 \\ 0 & 1 \\ 2 & 0 \end{bmatrix} \cup
$$

$$
\begin{bmatrix} 1 & 3 & 5 & | & 1 & | & 0 & 0 & 1 & 1 \\ 2 & 4 & 0 & | & 1 & | & 0 & 0 & 0 & 0 \end{bmatrix} \begin{bmatrix} 1 & 0 & 1 \\ 0 & 1 & 0 \\ 1 & 1 & 1 \\ \hline 0 & 1 & 2 \\ \hline 1 & 0 & 0 \\ 0 & 1 & 1 \\ 1 & 0 & 1 \\ 0 & 0 & 0 \end{bmatrix}
$$

$$
= \begin{bmatrix} 1 & 3 & 4 \\ 5 & 6 & 9 \\ 1 & 0 & 1 \end{bmatrix} \cup \left\{ \begin{bmatrix} 2 & 0 \\ 1 & 2 \\ 3 & 1 \end{bmatrix} \begin{bmatrix} 1 & 0 \\ 0 & 1 \end{bmatrix} + \right.
$$

$$+ \begin{bmatrix} 5 \\ 1 \\ 2 \end{bmatrix} \begin{bmatrix} 2 & 1 \end{bmatrix} + \begin{bmatrix} 7 & 3 & 1 & 4 \\ 1 & 0 & 0 & 2 \\ 0 & 1 & 1 & 0 \end{bmatrix} \begin{bmatrix} 1 & 0 \\ 1 & 1 \\ 0 & 1 \\ 2 & 0 \end{bmatrix} \Bigg\}$$

$$\cup \left\{ \begin{bmatrix} 1 & 3 & 5 \\ 2 & 4 & 0 \end{bmatrix} \begin{bmatrix} 1 & 0 & 1 \\ 0 & 1 & 0 \\ 1 & 1 & 1 \end{bmatrix} + \begin{bmatrix} 1 \\ 1 \end{bmatrix} \begin{bmatrix} 0 & 1 & 2 \end{bmatrix} + \right.$$

$$\begin{bmatrix} 0 & 0 & 1 & 1 \\ 1 & 0 & 0 & 0 \end{bmatrix} \begin{bmatrix} 1 & 0 & 0 \\ 0 & 1 & 1 \\ 1 & 0 & 1 \\ 0 & 0 & 0 \end{bmatrix} \Bigg\}$$

$$= \begin{bmatrix} 1 & 3 & 4 \\ 5 & 6 & 9 \\ 1 & 0 & 1 \end{bmatrix} \cup \left\{ \begin{bmatrix} 2 & 0 \\ 1 & 2 \\ 3 & 1 \end{bmatrix} + \begin{bmatrix} 10 & 5 \\ 2 & 1 \\ 4 & 2 \end{bmatrix} + \begin{bmatrix} 18 & 4 \\ 5 & 0 \\ 1 & 2 \end{bmatrix} \right\}$$

$$\cup \left\{ \begin{bmatrix} 6 & 8 & 6 \\ 2 & 4 & 2 \end{bmatrix} + \begin{bmatrix} 0 & 1 & 2 \\ 0 & 1 & 2 \end{bmatrix} + \begin{bmatrix} 1 & 0 & 1 \\ 1 & 0 & 0 \end{bmatrix} \right\}$$

$$= \begin{bmatrix} 1 & 3 & 4 \\ 5 & 6 & 9 \\ 1 & 0 & 1 \end{bmatrix} \cup \begin{bmatrix} 30 & 9 \\ 8 & 3 \\ 7 & 5 \end{bmatrix} \cup \begin{bmatrix} 7 & 9 & 9 \\ 3 & 5 & 4 \end{bmatrix}.$$

We see TP is just a trimatrix which is not a semi super trimatrix.

***Example 3.66:*** Let $T = T_1 \cup T_2 \cup T_3$ be a special semi super row trimatrix and $V = V_1 \cup V_2 \cup V_3$ be a special semi super column trimatrix.

$$T = T_1 \cup T_2 \cup T_3$$

$$= \begin{bmatrix} 2 & 3 & 4 & 5 \\ 1 & 2 & 3 & 4 \\ 0 & 1 & 1 & 2 \end{bmatrix} \cup$$

$$\begin{bmatrix} 1 & 5 & | & 1 & 0 & 1 & 0 & | & 1 & 0 & 1 \\ 2 & 0 & | & 0 & 1 & 1 & 1 & | & 0 & 0 & 1 \\ 3 & 1 & | & 0 & 1 & 0 & 0 & | & 0 & 1 & 1 \\ 4 & 1 & | & 1 & 0 & 1 & 1 & | & 0 & 0 & 0 \end{bmatrix}$$

$$\cup \begin{bmatrix} 1 & 1 & 0 & | & 0 & 1 & | & 1 & 1 & 0 & 1 \\ 3 & 1 & 1 & | & 1 & 0 & | & 1 & 0 & 1 & 0 \\ 0 & 1 & 0 & | & 1 & 1 & | & 0 & 1 & 0 & 1 \end{bmatrix}$$

be a special semi super row trimatrix.

$$V = V_1 \cup V_2 \cup V_3$$

$$= \begin{bmatrix} 0 & 1 & 1 \\ 1 & 0 & 0 \\ 2 & 1 & 1 \\ 3 & 1 & 0 \end{bmatrix} \cup \begin{bmatrix} 1 & 0 \\ 0 & 1 \\ \hline 1 & 1 \\ 1 & 0 \\ 0 & 1 \\ 0 & 0 \\ \hline 2 & 1 \\ 0 & 1 \\ 1 & 0 \end{bmatrix} \cup \begin{bmatrix} 1 & 1 & 0 & 1 \\ 0 & 1 & 0 & 0 \\ 1 & 0 & 0 & 1 \\ \hline 1 & 1 & 0 & 0 \\ 0 & 0 & 1 & 1 \\ 0 & 1 & 1 & 0 \\ \hline 1 & 0 & 0 & 0 \\ 0 & 1 & 0 & 0 \\ 0 & 0 & 0 & 1 \end{bmatrix}$$

be a special semi super column trimatrix. To find

$$\begin{aligned} TV &= (T = T_1 \cup T_2 \cup T_3)(V = V_1 \cup V_2 \cup V_3) \\ &= T_1 V_1 \cup T_2 V_2 \cup T_3 V_3 \end{aligned}$$

$$= \begin{bmatrix} 2 & 3 & 4 & 5 \\ 1 & 2 & 3 & 4 \\ 0 & 1 & 1 & 2 \end{bmatrix} \begin{bmatrix} 0 & 1 & 1 \\ 1 & 0 & 0 \\ 2 & 1 & 1 \\ 3 & 1 & 0 \end{bmatrix} \cup$$

$$\begin{bmatrix} 1 & 5 & | & 1 & 0 & 1 & 0 & | & 1 & 0 & 1 \\ 2 & 0 & | & 0 & 1 & 1 & 1 & | & 0 & 0 & 1 \\ 3 & 1 & | & 0 & 1 & 0 & 0 & | & 0 & 1 & 1 \\ 4 & 1 & | & 1 & 0 & 1 & 1 & | & 0 & 0 & 0 \end{bmatrix} \begin{bmatrix} 1 & 0 \\ 0 & 1 \\ \hline 1 & 1 \\ 1 & 0 \\ 0 & 1 \\ 0 & 0 \\ \hline 2 & 1 \\ 0 & 1 \\ 1 & 0 \end{bmatrix}$$

$$\cup \begin{bmatrix} 1 & 1 & 0 & | & 0 & 1 & | & 1 & 1 & 0 & 1 \\ 3 & 1 & 1 & | & 1 & 0 & | & 1 & 0 & 1 & 0 \\ 0 & 1 & 0 & | & 1 & 1 & | & 0 & 1 & 0 & 1 \end{bmatrix} \begin{bmatrix} 1 & 1 & 0 & 1 \\ 0 & 1 & 0 & 0 \\ 1 & 0 & 0 & 1 \\ \hline 1 & 1 & 0 & 0 \\ 0 & 0 & 1 & 1 \\ \hline 0 & 1 & 1 & 0 \\ 1 & 0 & 0 & 0 \\ 0 & 1 & 0 & 0 \\ 0 & 0 & 0 & 1 \end{bmatrix}$$

$$= \begin{bmatrix} 26 & 11 & 6 \\ 20 & 8 & 4 \\ 9 & 3 & 1 \end{bmatrix} \cup \left\{ \begin{bmatrix} 1 & 5 \\ 2 & 0 \\ 3 & 1 \\ 4 & 1 \end{bmatrix} \begin{bmatrix} 1 & 0 \\ 0 & 1 \end{bmatrix} + \right.$$

$$\left. \begin{bmatrix} 1 & 0 & 1 & 0 \\ 0 & 1 & 1 & 1 \\ 0 & 1 & 0 & 0 \\ 1 & 0 & 1 & 1 \end{bmatrix} \begin{bmatrix} 1 & 1 \\ 1 & 0 \\ 0 & 1 \\ 0 & 0 \end{bmatrix} + \begin{bmatrix} 1 & 0 & 1 \\ 0 & 0 & 1 \\ 0 & 1 & 1 \\ 0 & 0 & 0 \end{bmatrix} \begin{bmatrix} 2 & 1 \\ 0 & 1 \\ 1 & 0 \end{bmatrix} \right\}$$

$$\cup \left\{ \begin{bmatrix} 1 & 1 & 0 \\ 3 & 1 & 1 \\ 0 & 1 & 0 \end{bmatrix} \begin{bmatrix} 1 & 1 & 0 & 1 \\ 0 & 1 & 0 & 0 \\ 1 & 0 & 0 & 1 \end{bmatrix} + \begin{bmatrix} 0 & 1 \\ 1 & 0 \\ 1 & 1 \end{bmatrix} \begin{bmatrix} 1 & 1 & 0 & 0 \\ 0 & 0 & 1 & 1 \end{bmatrix} + \right.$$

$$\left. \begin{bmatrix} 1 & 1 & 0 & 1 \\ 1 & 0 & 1 & 0 \\ 0 & 1 & 0 & 1 \end{bmatrix} \begin{bmatrix} 0 & 1 & 1 & 0 \\ 1 & 0 & 0 & 0 \\ 0 & 1 & 0 & 0 \\ 0 & 0 & 0 & 1 \end{bmatrix} \right\}$$

$$= \begin{bmatrix} 26 & 11 & 6 \\ 20 & 8 & 4 \\ 9 & 3 & 1 \end{bmatrix} \cup \left\{ \begin{bmatrix} 1 & 5 \\ 2 & 0 \\ 3 & 1 \\ 4 & 1 \end{bmatrix} + \begin{bmatrix} 1 & 2 \\ 1 & 1 \\ 1 & 0 \\ 1 & 2 \end{bmatrix} + \begin{bmatrix} 3 & 1 \\ 1 & 0 \\ 1 & 1 \\ 0 & 0 \end{bmatrix} \right\} \cup$$

$$\left\{ \begin{bmatrix} 1 & 2 & 0 & 1 \\ 4 & 4 & 0 & 4 \\ 0 & 1 & 0 & 0 \end{bmatrix} + \begin{bmatrix} 0 & 0 & 1 & 1 \\ 1 & 1 & 0 & 0 \\ 1 & 1 & 1 & 1 \end{bmatrix} + \begin{bmatrix} 1 & 1 & 1 & 1 \\ 0 & 2 & 1 & 0 \\ 1 & 0 & 1 & 1 \end{bmatrix} \right\}$$

$$= \begin{bmatrix} 26 & 11 & 6 \\ 20 & 8 & 4 \\ 9 & 3 & 1 \end{bmatrix} \cup \begin{bmatrix} 5 & 8 \\ 4 & 1 \\ 5 & 2 \\ 5 & 3 \end{bmatrix} \cup \begin{bmatrix} 2 & 3 & 2 & 3 \\ 5 & 7 & 1 & 4 \\ 2 & 2 & 1 & 2 \end{bmatrix}.$$

Clearly this is only a mixed trimatrix which is not a semi super trimatrix.

Now we proceed on to define the notion of minor product of special semi super row trivector and a special semi super column trivector.

**DEFINITION 3.21:** *Let $T = T_1 \cup T_2 \cup T_3$ be any special semi super row trimatrix and $V = V_1 \cup V_2 \cup V_3$ be another special semi super column trimatrix. The minor product $TV = (T = T_1 \cup T_2 \cup T_3)(V = V_1 \cup V_2 \cup V_3) = T_1V_1 \cup T_2V_2 \cup T_3V_3$ is defined if each of the product $T_iV_i$, $1 \leq i \leq 3$ is defined.*

**Note:** In the case of minor product of two special semi super trimatrices the resultant is a only trimatrix and not even a super trimatrix or a semi super trimatrix.

***Example 3.67:*** Let $T = T_1 \cup T_2 \cup T_3$ be a special semi super row trivector to find the minor product of $TT^T$.
Given $T = T_1 \cup T_2 \cup T_3$

$$= \begin{bmatrix} 3 & 1 & 0 & 2 & 5 \\ 1 & 0 & 1 & 0 & 1 \\ 0 & 1 & 0 & 1 & 0 \\ 3 & 0 & 2 & 0 & 1 \end{bmatrix} \cup$$

$$\begin{bmatrix} 1 & 2 & 0 & 1 & 3 & 1 & 2 & 5 & 0 \\ 0 & 1 & 2 & 0 & 1 & 0 & 1 & 0 & 1 \\ 3 & 0 & 1 & 1 & 0 & 1 & 0 & 1 & 0 \end{bmatrix} \cup$$

$$\begin{bmatrix} 1 & 1 & 3 & 0 & 1 & 0 & 1 & 1 & 0 & 1 & 0 \\ 2 & 0 & 1 & 0 & 0 & 1 & 1 & 0 & 1 & 1 & 1 \\ 3 & 1 & 0 & 1 & 0 & 1 & 0 & 1 & 1 & 1 & 0 \\ 0 & 0 & 0 & 1 & 1 & 0 & 1 & 1 & 1 & 0 & 1 \end{bmatrix}$$

is a special semi super row trivector.

$$\begin{aligned} T^T &= (T = T_1 \cup T_2 \cup T_3)^T \\ &= T_1^T \cup T_2^T \cup T_3^T \end{aligned}$$

$$= \begin{bmatrix} 3 & 1 & 0 & 3 \\ 1 & 0 & 1 & 0 \\ 0 & 1 & 0 & 2 \\ 2 & 0 & 1 & 0 \\ 5 & 1 & 0 & 1 \end{bmatrix} \cup \begin{bmatrix} 1 & 0 & 3 \\ 2 & 1 & 0 \\ \hline 0 & 2 & 1 \\ 1 & 0 & 1 \\ 3 & 1 & 0 \\ \hline 1 & 0 & 1 \\ 2 & 1 & 0 \\ 5 & 0 & 1 \\ 0 & 1 & 0 \end{bmatrix} \cup$$

$$\begin{bmatrix} 1 & 2 & 3 & 0 \\ 1 & 0 & 1 & 0 \\ 3 & 1 & 0 & 0 \\ 0 & 0 & 1 & 1 \\ 1 & 0 & 0 & 1 \\ \hline 0 & 1 & 1 & 0 \\ 1 & 1 & 0 & 1 \\ 1 & 0 & 1 & 1 \\ \hline 0 & 1 & 1 & 1 \\ 1 & 1 & 1 & 0 \\ 0 & 1 & 0 & 1 \end{bmatrix}.$$

$$\begin{aligned} \text{Now } TT^T &= (T_1 \cup T_2 \cup T_3)(T_1^T \cup T_2^T \cup T_3^T) \\ &= T_1 T_1^T \cup T_2 T_2^T \cup T_3 T_3^T \end{aligned}$$

$$= \begin{bmatrix} 3 & 1 & 0 & 2 & 5 \\ 1 & 0 & 1 & 0 & 1 \\ 0 & 1 & 0 & 1 & 0 \\ 3 & 0 & 2 & 0 & 1 \end{bmatrix} \begin{bmatrix} 3 & 1 & 0 & 3 \\ 1 & 0 & 1 & 0 \\ 0 & 1 & 0 & 2 \\ 2 & 0 & 1 & 0 \\ 5 & 1 & 0 & 1 \end{bmatrix} \cup$$

$$\begin{bmatrix} 1 & 2 & | & 0 & 1 & 3 & | & 1 & 2 & 5 & 0 \\ 0 & 1 & | & 2 & 0 & 1 & | & 0 & 1 & 0 & 1 \\ 3 & 0 & | & 1 & 1 & 0 & | & 1 & 0 & 1 & 0 \end{bmatrix} \begin{bmatrix} 1 & 0 & 3 \\ 2 & 1 & 0 \\ \hline 0 & 2 & 1 \\ 1 & 0 & 1 \\ 3 & 1 & 0 \\ \hline 1 & 0 & 1 \\ 2 & 1 & 0 \\ 5 & 0 & 1 \\ 0 & 1 & 0 \end{bmatrix}$$

$$\begin{bmatrix} 1 & 1 & 3 & 0 & 1 & | & 0 & 1 & 1 & | & 0 & 1 & | & 0 \\ 2 & 0 & 1 & 0 & 0 & | & 1 & 1 & 0 & | & 1 & 1 & | & 1 \\ 3 & 1 & 0 & 1 & 0 & | & 1 & 0 & 1 & | & 1 & 1 & | & 0 \\ 0 & 0 & 0 & 1 & 1 & | & 0 & 1 & 1 & | & 1 & 0 & | & 1 \end{bmatrix} \begin{bmatrix} 1 & 2 & 3 & 0 \\ 1 & 0 & 1 & 0 \\ 3 & 1 & 0 & 0 \\ 0 & 0 & 1 & 1 \\ 1 & 0 & 0 & 1 \\ \hline 0 & 1 & 1 & 0 \\ 1 & 1 & 0 & 1 \\ 1 & 0 & 1 & 1 \\ \hline 0 & 1 & 1 & 1 \\ 1 & 1 & 1 & 0 \\ \hline 0 & 1 & 0 & 1 \end{bmatrix}$$

$$= \begin{bmatrix} 39 & 8 & 3 & 14 \\ 8 & 3 & 0 & 6 \\ 3 & 0 & 2 & 0 \\ 14 & 6 & 0 & 14 \end{bmatrix} \cup \left\{ \begin{bmatrix} 1 & 2 \\ 0 & 1 \\ 3 & 0 \end{bmatrix} \begin{bmatrix} 1 & 0 & 3 \\ 2 & 1 & 0 \end{bmatrix} + \right.$$

$$\left. \begin{bmatrix} 0 & 1 & 3 \\ 2 & 0 & 1 \\ 1 & 1 & 0 \end{bmatrix} \begin{bmatrix} 0 & 2 & 1 \\ 1 & 0 & 1 \\ 3 & 1 & 0 \end{bmatrix} + \begin{bmatrix} 1 & 2 & 5 & 0 \\ 0 & 1 & 0 & 1 \\ 1 & 0 & 1 & 0 \end{bmatrix} \begin{bmatrix} 1 & 0 & 1 \\ 2 & 1 & 0 \\ 5 & 0 & 1 \\ 0 & 1 & 0 \end{bmatrix} \right\} \cup$$

$$\left\{ \begin{bmatrix} 1 & 1 & 3 & 0 & 1 \\ 2 & 0 & 1 & 0 & 0 \\ 3 & 1 & 0 & 1 & 0 \\ 0 & 0 & 0 & 1 & 1 \end{bmatrix} \begin{bmatrix} 1 & 2 & 3 & 0 \\ 1 & 0 & 1 & 0 \\ 3 & 1 & 0 & 0 \\ 0 & 0 & 1 & 1 \\ 1 & 0 & 0 & 1 \end{bmatrix} \right.$$

$$+ \begin{bmatrix} 0 & 1 & 1 \\ 1 & 1 & 0 \\ 1 & 0 & 1 \\ 0 & 1 & 1 \end{bmatrix} \begin{bmatrix} 0 & 1 & 1 & 0 \\ 1 & 1 & 0 & 1 \\ 1 & 0 & 1 & 1 \end{bmatrix} + \begin{bmatrix} 0 & 1 \\ 1 & 1 \\ 1 & 1 \\ 1 & 0 \end{bmatrix} \begin{bmatrix} 0 & 1 & 1 & 1 \\ 1 & 1 & 1 & 0 \end{bmatrix} +$$

$$\left. \begin{bmatrix} 0 \\ 1 \\ 0 \\ 1 \end{bmatrix} \begin{bmatrix} 0 & 1 & 0 & 1 \end{bmatrix} \right\} = \begin{bmatrix} 39 & 8 & 3 & 14 \\ 8 & 3 & 0 & 6 \\ 3 & 0 & 2 & 0 \\ 14 & 6 & 0 & 14 \end{bmatrix} \cup \left\{ \begin{bmatrix} 5 & 2 & 3 \\ 2 & 1 & 0 \\ 3 & 0 & 9 \end{bmatrix} + \right.$$

$$\left. \begin{bmatrix} 10 & 3 & 1 \\ 3 & 5 & 2 \\ 1 & 2 & 2 \end{bmatrix} + \begin{bmatrix} 30 & 2 & 6 \\ 2 & 2 & 0 \\ 6 & 0 & 2 \end{bmatrix} \right\} \cup$$

$$\left\{ \begin{bmatrix} 12 & 5 & 4 & 1 \\ 5 & 5 & 6 & 0 \\ 4 & 6 & 11 & 1 \\ 1 & 0 & 1 & 2 \end{bmatrix} + \begin{bmatrix} 2 & 1 & 1 & 2 \\ 1 & 2 & 1 & 1 \\ 1 & 1 & 2 & 1 \\ 2 & 1 & 1 & 2 \end{bmatrix} + \right.$$

$$\left. \begin{bmatrix} 1 & 1 & 1 & 0 \\ 1 & 2 & 2 & 1 \\ 1 & 2 & 2 & 1 \\ 0 & 1 & 1 & 1 \end{bmatrix} + \begin{bmatrix} 0 & 0 & 0 & 0 \\ 0 & 1 & 0 & 1 \\ 0 & 0 & 0 & 0 \\ 0 & 1 & 0 & 1 \end{bmatrix} \right\}$$

$$= \begin{bmatrix} 39 & 8 & 3 & 14 \\ 8 & 3 & 0 & 6 \\ 3 & 0 & 2 & 0 \\ 14 & 6 & 0 & 14 \end{bmatrix} \cup \begin{bmatrix} 45 & 7 & 10 \\ 7 & 8 & 2 \\ 10 & 2 & 13 \end{bmatrix} \cup \begin{bmatrix} 15 & 7 & 6 & 3 \\ 7 & 10 & 9 & 3 \\ 6 & 9 & 15 & 3 \\ 3 & 3 & 3 & 6 \end{bmatrix}.$$

We see the minor product of $TT^T$ is a symmetric trimatrix and is not a super trimatrix. In fact $TT^T$ is a mixed symmetric trimatrix. Thus the minor product of special semi super row trivector with its transpose yields a symmetric trimatrix.

***Example 3.68:*** Let $T = T_1 \cup T_2 \cup T_3$ be a special semi super row trivector. To find the product of $TT^T$. Given

$$T = \begin{bmatrix} 0 & 1 & 2 & 3 \\ 1 & 2 & 3 & 0 \\ 2 & 3 & 1 & 0 \end{bmatrix} \cup \begin{bmatrix} 1 & 2 & 3 & 0 & 1 & 1 & 1 & 1 \\ 0 & 0 & 1 & 1 & 0 & 1 & 0 & 1 \\ 2 & 0 & 1 & 0 & 1 & 0 & 1 & 0 \\ 1 & 1 & 0 & 1 & 0 & 0 & 1 & 0 \end{bmatrix} \cup$$

$$\begin{bmatrix} 1 & 0 & 2 & 1 & 1 & 1 & 1 & 0 & 0 & 1 & 0 \\ 0 & 1 & 0 & 2 & 1 & 0 & 0 & 0 & 1 & 1 & 1 \\ 2 & 2 & 1 & 0 & 0 & 1 & 0 & 1 & 0 & 1 & 1 \end{bmatrix},$$

$$\begin{aligned} T^T &= (T_1 \cup T_2 \cup T_3)^T \\ &= T_1^T \cup T_2^T \cup T_3^T \end{aligned}$$

$$= \begin{bmatrix} 0 & 1 & 2 \\ 1 & 2 & 3 \\ 2 & 3 & 1 \\ 3 & 0 & 0 \end{bmatrix} \cup \begin{bmatrix} 1 & 0 & 2 & 1 \\ 2 & 0 & 0 & 1 \\ \hline 3 & 1 & 1 & 0 \\ \hline 0 & 1 & 0 & 1 \\ 1 & 0 & 1 & 0 \\ 1 & 1 & 0 & 0 \\ 1 & 0 & 1 & 1 \\ 1 & 1 & 0 & 0 \end{bmatrix} \cup \begin{bmatrix} 1 & 0 & 2 \\ 0 & 1 & 2 \\ 2 & 0 & 1 \\ 1 & 2 & 0 \\ 1 & 1 & 0 \\ \hline 1 & 0 & 1 \\ 1 & 0 & 0 \\ 0 & 0 & 1 \\ \hline 0 & 1 & 0 \\ 1 & 1 & 1 \\ 0 & 1 & 1 \end{bmatrix}.$$

$$\begin{aligned}
TT^T &= (T_1 \cup T_2 \cup T_3)(T_1^T \cup T_2^T \cup T_3^T) \\
&= T_1 T_1^T \cup T_2 T_2^T \cup T_3 T_3^T
\end{aligned}$$

$$= \begin{bmatrix} 0 & 1 & 2 & 3 \\ 1 & 2 & 3 & 0 \\ 2 & 3 & 1 & 0 \end{bmatrix} \begin{bmatrix} 0 & 1 & 2 \\ 1 & 2 & 3 \\ 2 & 3 & 1 \\ 3 & 0 & 0 \end{bmatrix} \cup$$

$$\begin{bmatrix} 1 & 2 & 3 & 0 & 1 & 1 & 1 & 1 \\ 0 & 6 & 1 & 1 & 0 & 1 & 0 & 1 \\ 2 & 0 & 1 & 0 & 1 & 0 & 1 & 0 \\ 1 & 1 & 0 & 1 & 0 & 0 & 1 & 0 \end{bmatrix} \begin{bmatrix} 1 & 0 & 2 & 1 \\ 2 & 6 & 0 & 1 \\ \hline 3 & 1 & 1 & 0 \\ \hline 0 & 1 & 0 & 1 \\ 1 & 0 & 1 & 0 \\ 1 & 1 & 0 & 0 \\ 1 & 0 & 1 & 1 \\ 1 & 1 & 0 & 0 \end{bmatrix}$$

$$\begin{bmatrix} 1 & 0 & 2 & 1 \\ 0 & 1 & 0 & 2 \\ 2 & 2 & 1 & 0 \end{bmatrix} \begin{bmatrix} 1 & 1 & 1 \\ 1 & 0 & 0 \\ 0 & 1 & 0 \end{bmatrix} \begin{bmatrix} 0 & 0 & 1 \\ 0 & 1 & 1 \\ 1 & 0 & 1 \end{bmatrix} \begin{bmatrix} 0 \\ 1 \\ 1 \end{bmatrix} \begin{bmatrix} 1 & 0 & 2 \\ 0 & 1 & 2 \\ 2 & 0 & 1 \\ 1 & 2 & 0 \\ 1 & 1 & 0 \\ \hline 1 & 0 & 1 \\ 1 & 0 & 0 \\ 0 & 0 & 1 \\ \hline 0 & 1 & 0 \\ 1 & 1 & 1 \\ \hline 0 & 1 & 1 \end{bmatrix}$$

$$= \begin{bmatrix} 14 & 8 & 5 \\ 8 & 14 & 11 \\ 5 & 11 & 14 \end{bmatrix} \cup \left\{ \begin{bmatrix} 1 & 2 \\ 0 & 6 \\ 2 & 0 \\ 1 & 1 \end{bmatrix} \begin{bmatrix} 1 & 0 & 2 & 1 \\ 2 & 6 & 0 & 1 \end{bmatrix} + \right.$$

$$+ \begin{pmatrix} 3 \\ 1 \\ 1 \\ 0 \end{pmatrix} \begin{bmatrix} 3 & 1 & 1 & 0 \end{bmatrix} + \begin{bmatrix} 0 & 1 & 1 & 1 & 1 \\ 1 & 0 & 1 & 0 & 1 \\ 0 & 1 & 0 & 1 & 0 \\ 1 & 0 & 0 & 1 & 0 \end{bmatrix} \begin{bmatrix} 0 & 1 & 0 & 1 \\ 1 & 0 & 1 & 0 \\ 1 & 1 & 0 & 0 \\ 1 & 0 & 1 & 1 \\ 1 & 1 & 0 & 0 \end{bmatrix} \right\} \cup$$

$$\left\{ \begin{bmatrix} 1 & 0 & 2 & 1 & 1 \\ 0 & 1 & 0 & 2 & 1 \\ 2 & 2 & 1 & 0 & 0 \end{bmatrix} \begin{bmatrix} 1 & 0 & 2 \\ 0 & 1 & 2 \\ 2 & 0 & 1 \\ 1 & 2 & 0 \\ 1 & 1 & 0 \end{bmatrix} + \right.$$

$$+ \begin{bmatrix} 1 & 1 & 0 \\ 0 & 0 & 0 \\ 1 & 0 & 1 \end{bmatrix} \begin{bmatrix} 1 & 0 & 1 \\ 1 & 0 & 0 \\ 0 & 0 & 1 \end{bmatrix} + \begin{bmatrix} 0 & 1 \\ 1 & 1 \\ 0 & 1 \end{bmatrix} \begin{bmatrix} 0 & 1 & 0 \\ 1 & 1 & 1 \end{bmatrix} + \begin{bmatrix} 0 \\ 1 \\ 1 \end{bmatrix} [0 \; 1 \; 1] \Bigg\}$$

$$= \begin{bmatrix} 14 & 8 & 5 \\ 8 & 14 & 11 \\ 5 & 11 & 14 \end{bmatrix} \cup \Bigg\{ \begin{bmatrix} 5 & 12 & 2 & 3 \\ 12 & 36 & 0 & 6 \\ 2 & 0 & 4 & 2 \\ 3 & 6 & 2 & 2 \end{bmatrix} +$$

$$\begin{bmatrix} 9 & 3 & 3 & 0 \\ 3 & 1 & 1 & 0 \\ 3 & 1 & 1 & 0 \\ 0 & 0 & 0 & 0 \end{bmatrix} + \begin{bmatrix} 4 & 2 & 2 & 1 \\ 2 & 3 & 0 & 1 \\ 2 & 0 & 2 & 1 \\ 1 & 1 & 1 & 2 \end{bmatrix} \Bigg\}$$

$$\cup \Bigg\{ \begin{bmatrix} 7 & 3 & 4 \\ 3 & 6 & 2 \\ 4 & 2 & 9 \end{bmatrix} + \begin{bmatrix} 2 & 0 & 1 \\ 0 & 0 & 0 \\ 1 & 0 & 2 \end{bmatrix} + \begin{bmatrix} 1 & 1 & 1 \\ 1 & 2 & 1 \\ 1 & 1 & 1 \end{bmatrix} + \begin{bmatrix} 0 & 0 & 0 \\ 0 & 1 & 1 \\ 0 & 1 & 1 \end{bmatrix} \Bigg\}$$

$$= \begin{bmatrix} 14 & 8 & 5 \\ 8 & 14 & 11 \\ 5 & 11 & 14 \end{bmatrix} \cup \begin{bmatrix} 18 & 17 & 7 & 4 \\ 17 & 40 & 1 & 7 \\ 7 & 1 & 7 & 3 \\ 4 & 7 & 3 & 4 \end{bmatrix} \cup \begin{bmatrix} 10 & 4 & 6 \\ 4 & 9 & 4 \\ 6 & 4 & 13 \end{bmatrix}.$$

This is clearly a symmetric trimatrix which is neither semi super or super.

Now we proceed on to illustrate major product of semi super trimatrices.

*Example 3.69:* Let $T = T_1 \cup T_2 \cup T_3$ be a special semi super row trivector. $V = V_1 \cup V_2 \cup V_3$ be a special semi super column trivector.

$$\begin{aligned} TV &= (T = T_1 \cup T_2 \cup T_3)(V = V_1 \cup V_2 \cup V_3) \\ &= T_1V_1 \cup T_2V_2 \cup T_3V_3 \end{aligned}$$

here

$$T_1 = \begin{bmatrix} 7 & 0 & 2 & 1 & 5 \\ 1 & 1 & 0 & 0 & 1 \\ 2 & 0 & 1 & 0 & 1 \end{bmatrix},$$

$$T_2 = \left[\begin{array}{c|cc|ccc|cccc} 1 & 0 & 3 & 1 & 0 & 1 & 1 & 0 & 1 & 2 \\ 0 & 1 & 0 & 0 & 1 & 1 & 0 & 1 & 1 & 2 \\ 7 & 0 & 1 & 1 & 1 & 0 & 0 & 2 & 1 & 1 \end{array}\right]$$

and

$$T_3 = \left[\begin{array}{ccc|cc|c|ccccc} 0 & 1 & 2 & 3 & 1 & 0 & 3 & 4 & 5 & 0 & 1 \\ 1 & 0 & 1 & 1 & 0 & 1 & 1 & 0 & 1 & 1 & 1 \\ 2 & 1 & 0 & 1 & 1 & 0 & 0 & 1 & 0 & 0 & 0 \\ 3 & 0 & 1 & 0 & 1 & 1 & 0 & 0 & 1 & 0 & 1 \\ 0 & 1 & 0 & 1 & 0 & 0 & 1 & 1 & 0 & 1 & 0 \end{array}\right]$$

be the special semi super row trivector. $V = V_1 \cup V_2 \cup V_3$;

$$V_1 = \begin{bmatrix} 1 & 0 & 2 \\ 0 & 1 & 1 \\ 1 & 0 & 0 \\ 0 & 1 & 0 \\ 1 & 0 & 1 \end{bmatrix}, V_2 = \left[\begin{array}{ccc} 0 & 1 & 1 & 0 \\ 1 & 0 & 1 & 0 \\ 0 & 1 & 1 & 1 \\ 0 & 0 & 1 & 0 \\ \hline 1 & 0 & 0 & 1 \\ 0 & 1 & 0 & 1 \\ \hline 0 & 1 & 1 & 0 \\ 1 & 0 & 0 & 1 \\ 0 & 1 & 0 & 1 \\ 0 & 0 & 1 & 0 \end{array}\right], V_3 = \left[\begin{array}{ccc} 0 & 1 & 0 \\ 1 & 0 & 1 \\ 0 & 1 & 1 \\ \hline 1 & 1 & 1 \\ 0 & 1 & 1 \\ 0 & 0 & 1 \\ \hline 1 & 0 & 1 \\ 1 & 1 & 0 \\ 0 & 1 & 0 \\ 1 & 0 & 1 \\ 0 & 1 & 0 \end{array}\right]$$

be the semi super column trivector.

$$TV = \begin{bmatrix} 7 & 0 & 2 & 1 & 5 \\ 1 & 1 & 0 & 0 & 1 \\ 2 & 0 & 1 & 0 & 1 \end{bmatrix} \begin{bmatrix} 1 & 0 & 2 \\ 0 & 1 & 1 \\ 1 & 0 & 0 \\ 0 & 1 & 0 \\ 1 & 0 & 1 \end{bmatrix} \cup$$

$$\begin{bmatrix} 1 & 0 & 3 & 1 & 0 & 1 & 1 & 0 & 1 & 2 \\ 0 & 1 & 0 & 0 & 1 & 1 & 0 & 1 & 1 & 2 \\ 7 & 0 & 1 & 1 & 1 & 0 & 0 & 2 & 1 & 1 \end{bmatrix} \begin{bmatrix} 0 & 1 & 1 & 1 \\ 1 & 0 & 1 & 0 \\ 0 & 1 & 1 & 1 \\ 0 & 0 & 1 & 0 \\ 1 & 0 & 0 & 1 \\ 0 & 1 & 0 & 1 \\ 0 & 1 & 1 & 0 \\ 1 & 0 & 0 & 1 \\ 0 & 1 & 0 & 1 \\ 0 & 0 & 1 & 0 \end{bmatrix}$$

$$\cup \begin{bmatrix} 0 & 1 & 2 & 3 & 1 & 0 & 3 & 4 & 5 & 0 & 1 \\ 1 & 0 & 1 & 1 & 0 & 1 & 0 & 1 & 1 & 1 & 1 \\ 2 & 1 & 0 & 1 & 1 & 0 & 0 & 1 & 0 & 0 & 0 \\ 3 & 0 & 1 & 0 & 1 & 1 & 0 & 0 & 1 & 0 & 1 \\ 0 & 1 & 0 & 1 & 0 & 0 & 1 & 1 & 0 & 1 & 0 \end{bmatrix} \begin{bmatrix} 0 & 1 & 0 \\ 1 & 0 & 1 \\ 0 & 1 & 1 \\ 1 & 1 & 1 \\ 0 & 1 & 1 \\ 0 & 0 & 1 \\ 1 & 0 & 1 \\ 1 & 1 & 0 \\ 0 & 1 & 0 \\ 1 & 0 & 1 \\ 0 & 1 & 0 \end{bmatrix}$$

$$= \begin{bmatrix} 14 & 1 & 19 \\ 2 & 1 & 4 \\ 4 & 0 & 5 \end{bmatrix} \cup \left\{ \begin{bmatrix} 1 \\ 0 \\ 7 \end{bmatrix} \begin{bmatrix} 0 & 1 & 1 & 1 \end{bmatrix} + \right.$$

$$\begin{bmatrix} 0 & 3 \\ 1 & 0 \\ 0 & 1 \end{bmatrix} \begin{bmatrix} 1 & 0 & 1 & 0 \\ 0 & 1 & 1 & 1 \end{bmatrix} + \begin{bmatrix} 1 & 0 & 1 \\ 0 & 1 & 1 \\ 1 & 1 & 0 \end{bmatrix} \begin{bmatrix} 0 & 0 & 1 & 0 \\ 1 & 0 & 0 & 1 \\ 0 & 1 & 0 & 1 \end{bmatrix}$$

$$+ \begin{bmatrix} 1 & 0 & 1 & 2 \\ 0 & 1 & 1 & 2 \\ 0 & 2 & 1 & 1 \end{bmatrix} \begin{bmatrix} 0 & 1 & 1 & 0 \\ 1 & 0 & 0 & 1 \\ 0 & 1 & 0 & 1 \\ 0 & 0 & 1 & 0 \end{bmatrix} \right\} \cup$$

$$\left\{ \begin{bmatrix} 0 & 1 & 2 \\ 1 & 0 & 1 \\ 2 & 1 & 0 \\ 3 & 0 & 1 \\ 0 & 1 & 0 \end{bmatrix} \begin{bmatrix} 0 & 1 & 0 \\ 1 & 0 & 1 \\ 0 & 1 & 1 \end{bmatrix} + \begin{bmatrix} 3 & 1 \\ 1 & 0 \\ 1 & 1 \\ 0 & 1 \\ 1 & 0 \end{bmatrix} \begin{bmatrix} 1 & 1 & 1 \\ 0 & 1 & 1 \end{bmatrix} + \right.$$

$$\left. \begin{bmatrix} 0 \\ 1 \\ 0 \\ 1 \\ 0 \end{bmatrix} \begin{bmatrix} 0 & 0 & 1 \end{bmatrix} + \begin{bmatrix} 3 & 4 & 5 & 0 & 1 \\ 1 & 0 & 1 & 1 & 1 \\ 0 & 1 & 0 & 0 & 0 \\ 0 & 0 & 1 & 0 & 1 \\ 1 & 1 & 0 & 1 & 0 \end{bmatrix} \begin{bmatrix} 1 & 0 & 1 \\ 1 & 1 & 0 \\ 0 & 1 & 0 \\ 1 & 0 & 1 \\ 0 & 1 & 0 \end{bmatrix} \right\}$$

$$= \begin{bmatrix} 14 & 1 & 19 \\ 2 & 1 & 4 \\ 4 & 0 & 5 \end{bmatrix} \cup \left\{ \begin{bmatrix} 0 & 1 & 1 & 1 \\ 0 & 0 & 0 & 0 \\ 0 & 7 & 7 & 7 \end{bmatrix} + \begin{bmatrix} 0 & 3 & 3 & 3 \\ 1 & 0 & 1 & 0 \\ 0 & 1 & 1 & 1 \end{bmatrix} + \right.$$

$$\left\{ \begin{bmatrix} 0 & 1 & 1 & 1 \\ 1 & 1 & 0 & 2 \\ 1 & 0 & 1 & 1 \end{bmatrix} + \begin{bmatrix} 0 & 2 & 3 & 1 \\ 1 & 1 & 2 & 2 \\ 2 & 1 & 1 & 3 \end{bmatrix} \right\} \cup$$

$$\left\{ \begin{bmatrix} 1 & 2 & 3 \\ 0 & 2 & 1 \\ 1 & 2 & 1 \\ 0 & 4 & 1 \\ 1 & 0 & 1 \end{bmatrix} + \begin{bmatrix} 3 & 4 & 4 \\ 1 & 1 & 1 \\ 1 & 2 & 2 \\ 0 & 1 & 1 \\ 1 & 1 & 1 \end{bmatrix} + \begin{bmatrix} 0 & 0 & 0 \\ 0 & 0 & 1 \\ 0 & 0 & 0 \\ 0 & 0 & 1 \\ 0 & 0 & 0 \end{bmatrix} + \begin{bmatrix} 7 & 10 & 3 \\ 2 & 2 & 2 \\ 1 & 1 & 0 \\ 0 & 2 & 0 \\ 3 & 1 & 2 \end{bmatrix} \right\}$$

$$= \begin{bmatrix} 14 & 1 & 19 \\ 2 & 1 & 4 \\ 4 & 0 & 5 \end{bmatrix} \cup \begin{bmatrix} 0 & 7 & 8 & 6 \\ 3 & 2 & 3 & 4 \\ 3 & 9 & 9 & 12 \end{bmatrix} \cup \begin{bmatrix} 11 & 16 & 10 \\ 3 & 5 & 5 \\ 3 & 5 & 3 \\ 0 & 7 & 3 \\ 5 & 2 & 4 \end{bmatrix}.$$

Thus the product results only in a usual trimatrix.

Now we give an illustration of the major product.

***Example 3.70:*** Let $S = S_1 \cup S_2 \cup S_3$ be a special semi super column trivector. $V = V_1 \cup V_2 \cup V_3$ be another special semi super row trivector. To find the product SV. Given

$$S = \begin{bmatrix} 3 & 1 & 0 & 1 \\ 4 & 0 & 1 & 0 \\ 5 & 1 & 0 & 0 \\ 0 & 1 & 1 & 0 \\ 1 & 0 & 0 & 1 \end{bmatrix} \cup \left[ \begin{array}{cccc} 1 & 0 & 1 & 2 \\ 0 & 1 & 0 & 1 \\ 1 & 1 & 0 & 0 \\ \hline 2 & 1 & 0 & 1 \\ 1 & 0 & 2 & 1 \\ 1 & 1 & 1 & 1 \\ 0 & 1 & 1 & 0 \\ \hline 2 & 0 & 0 & 1 \end{array} \right] \cup \left[ \begin{array}{ccc|cc} 1 & 0 & 1 & 1 & 1 \\ 0 & 1 & 0 & 1 & 1 \\ 1 & 1 & 0 & 1 & 0 \\ \hline 1 & 2 & 1 & 1 & 0 \\ 1 & 0 & 1 & 0 & 1 \\ 0 & 1 & 0 & 1 & 0 \\ 1 & 0 & 1 & 1 & 2 \\ 0 & 1 & 0 & 1 & 0 \\ 3 & 1 & 1 & 0 & 1 \\ 0 & 2 & 0 & 1 & 0 \end{array} \right]$$

is a special semi super column trivector. Here

$$V = \begin{bmatrix} 0 & 2 & 3 & 1 & 1 \\ 0 & 0 & 0 & 1 & 0 \\ 1 & 1 & 1 & 1 & 0 \\ 1 & 0 & 0 & 0 & 1 \end{bmatrix} \cup \left[\begin{array}{cccc|cc|ccc} 0 & 1 & 1 & 0 & 0 & 1 & 1 & 0 & 1 \\ 1 & 0 & 1 & 1 & 1 & 0 & 0 & 1 & 0 \\ 2 & 0 & 0 & 1 & 1 & 1 & 1 & 1 & 0 \\ 0 & 1 & 0 & 0 & 0 & 0 & 1 & 0 & 1 \end{array}\right]$$

$$\cup \left[\begin{array}{ccc|cccc|cc} 0 & 2 & 1 & 3 & 1 & 0 & 1 & 3 & 1 \\ 1 & 0 & 0 & 1 & 0 & 1 & 0 & 0 & 0 & 1 \\ 1 & 0 & 1 & 1 & 0 & 0 & 1 & 0 & 1 & 0 \\ 0 & 0 & 0 & 2 & 0 & 0 & 1 & 0 & 2 & 0 \\ 1 & 1 & 1 & 1 & 1 & 1 & 0 & 1 & 1 & 1 \end{array}\right]$$

be the special semi super row trivector.

$$\begin{aligned} SV &= (S_1 \cup S_2 \cup S_3)(V_1 \cup V_2 \cup V_3) \\ &= S_1V_1 \cup S_2V_2 \cup S_3V_3 \end{aligned}$$

$$= \begin{bmatrix} 3 & 1 & 0 & 1 \\ 4 & 0 & 1 & 0 \\ 5 & 1 & 0 & 0 \\ 0 & 1 & 1 & 0 \\ 1 & 0 & 0 & 1 \end{bmatrix} \begin{bmatrix} 0 & 2 & 3 & 1 & 1 \\ 0 & 0 & 0 & 1 & 0 \\ 1 & 1 & 1 & 1 & 0 \\ 1 & 0 & 0 & 0 & 1 \end{bmatrix} \cup$$

$$\begin{bmatrix} 1 & 0 & 1 & 2 \\ 0 & 1 & 0 & 1 \\ 1 & 1 & 0 & 0 \\ \hline 2 & 1 & 0 & 1 \\ 1 & 0 & 2 & 1 \\ 1 & 1 & 1 & 1 \\ \hline 0 & 1 & 1 & 0 \\ 2 & 0 & 0 & 1 \end{bmatrix} \left[\begin{array}{cccc|cc|ccc} 0 & 1 & 1 & 0 & 0 & 1 & 1 & 0 & 1 \\ 1 & 0 & 1 & 1 & 1 & 0 & 0 & 1 & 0 \\ 2 & 0 & 0 & 1 & 1 & 1 & 1 & 1 & 0 \\ 0 & 1 & 0 & 0 & 0 & 0 & 1 & 0 & 1 \end{array}\right] \cup$$

$$\begin{bmatrix} 1 & 0 & 1 & 1 & 1 \\ \hline 0 & 1 & 0 & 1 & 1 \\ 1 & 1 & 0 & 1 & 0 \\ \hline 1 & 2 & 1 & 1 & 0 \\ 1 & 0 & 1 & 0 & 1 \\ 0 & 1 & 0 & 1 & 0 \\ \hline 1 & 0 & 1 & 1 & 2 \\ 0 & 1 & 0 & 1 & 0 \\ 3 & 1 & 1 & 0 & 1 \\ 0 & 2 & 0 & 1 & 0 \end{bmatrix} \begin{bmatrix} 0 & 2 & 1 & 3 & 1 & 0 & 0 & 1 & 3 & 1 \\ 1 & 0 & 0 & 1 & 0 & 1 & 0 & 0 & 0 & 1 \\ 1 & 0 & 1 & 1 & 0 & 0 & 1 & 0 & 1 & 0 \\ 0 & 0 & 0 & 2 & 0 & 0 & 1 & 0 & 2 & 0 \\ 1 & 1 & 1 & 1 & 1 & 1 & 0 & 1 & 1 & 1 \end{bmatrix}$$

$$= \begin{bmatrix} 1 & 6 & 9 & 4 & 4 \\ 1 & 9 & 13 & 5 & 4 \\ 0 & 10 & 15 & 6 & 5 \\ 1 & 1 & 1 & 2 & 0 \\ 1 & 2 & 3 & 1 & 2 \end{bmatrix} \cup$$

$$\left[ \begin{array}{c|c} \begin{pmatrix} 1 & 0 & 1 & 2 \\ 0 & 1 & 0 & 1 \\ 1 & 1 & 0 & 0 \end{pmatrix} \begin{bmatrix} 0 & 1 & 1 & 0 \\ 1 & 0 & 1 & 1 \\ 2 & 0 & 0 & 1 \\ 0 & 1 & 0 & 0 \end{bmatrix} & \begin{pmatrix} 1 & 0 & 1 & 2 \\ 0 & 1 & 0 & 1 \\ 1 & 1 & 0 & 0 \end{pmatrix} \begin{pmatrix} 0 & 1 \\ 1 & 0 \\ 1 & 1 \\ 0 & 0 \end{pmatrix} \\ \hline \begin{pmatrix} 2 & 1 & 0 & 1 \\ 1 & 0 & 2 & 1 \\ 1 & 1 & 1 & 1 \\ 0 & 1 & 1 & 0 \end{pmatrix} \begin{bmatrix} 0 & 1 & 1 & 0 \\ 1 & 0 & 1 & 1 \\ 2 & 0 & 0 & 1 \\ 0 & 1 & 0 & 0 \end{bmatrix} & \begin{pmatrix} 2 & 1 & 0 & 1 \\ 1 & 0 & 2 & 1 \\ 1 & 1 & 1 & 1 \\ 0 & 1 & 1 & 0 \end{pmatrix} \begin{pmatrix} 0 & 1 \\ 1 & 0 \\ 1 & 1 \\ 0 & 0 \end{pmatrix} \\ \hline \begin{pmatrix} 2 & 0 & 0 & 1 \end{pmatrix} \begin{bmatrix} 0 & 1 & 1 & 0 \\ 1 & 0 & 1 & 1 \\ 2 & 0 & 0 & 1 \\ 0 & 1 & 0 & 0 \end{bmatrix} & \begin{pmatrix} 2 & 0 & 0 & 1 \end{pmatrix} \begin{pmatrix} 0 & 1 \\ 1 & 0 \\ 1 & 1 \\ 0 & 0 \end{pmatrix} \end{array} \right]$$

$$\subset \left[ \begin{array}{c} \begin{pmatrix} 1 & 0 & 1 & 2 \\ 0 & 1 & 0 & 1 \\ 1 & 1 & 0 & 0 \end{pmatrix} \begin{bmatrix} 1 & 0 & 1 \\ 0 & 1 & 0 \\ 1 & 1 & 0 \\ 1 & 0 & 1 \end{bmatrix} \\ \hline \begin{pmatrix} 2 & 1 & 0 & 1 \\ 1 & 0 & 2 & 1 \\ 1 & 1 & 1 & 1 \\ 0 & 1 & 1 & 0 \end{pmatrix} \begin{bmatrix} 1 & 0 & 1 \\ 0 & 1 & 0 \\ 1 & 1 & 0 \\ 1 & 0 & 1 \end{bmatrix} \\ \hline \begin{pmatrix} 2 & 0 & 0 & 1 \end{pmatrix} \begin{bmatrix} 1 & 0 & 1 \\ 0 & 1 & 0 \\ 1 & 1 & 0 \\ 1 & 0 & 1 \end{bmatrix} \end{array} \right]$$

$$
\begin{bmatrix}
\begin{pmatrix} 1 & 0 & 1 & 1 & 1 \end{pmatrix} \begin{bmatrix} 0 & 2 & 1 \\ 1 & 0 & 0 \\ 1 & 0 & 1 \\ 0 & 0 & 0 \\ 1 & 1 & 1 \end{bmatrix} & \begin{bmatrix} 1 & 0 & 1 & 1 & 1 \end{bmatrix} \begin{bmatrix} 3 & 1 & 0 & 0 & 1 \\ 1 & 0 & 1 & 0 & 0 \\ 1 & 0 & 0 & 1 & 0 \\ 2 & 0 & 0 & 1 & 0 \\ 1 & 1 & 1 & 0 & 1 \end{bmatrix} \\
\begin{pmatrix} 0 & 1 & 0 & 1 & 1 \\ 1 & 1 & 0 & 1 & 0 \end{pmatrix} \begin{bmatrix} 0 & 2 & 1 \\ 1 & 0 & 0 \\ 1 & 0 & 1 \\ 0 & 0 & 0 \\ 1 & 1 & 1 \end{bmatrix} & \begin{bmatrix} 0 & 1 & 0 & 1 & 1 \\ 1 & 1 & 0 & 1 & 0 \end{bmatrix} \begin{bmatrix} 3 & 1 & 0 & 0 & 1 \\ 1 & 0 & 1 & 0 & 0 \\ 1 & 0 & 0 & 1 & 0 \\ 2 & 0 & 0 & 1 & 0 \\ 1 & 1 & 1 & 0 & 1 \end{bmatrix} \\
\begin{pmatrix} 1 & 2 & 1 & 1 & 0 \\ 1 & 0 & 1 & 0 & 1 \\ 0 & 1 & 0 & 1 & 0 \end{pmatrix} \begin{bmatrix} 0 & 2 & 1 \\ 1 & 0 & 0 \\ 1 & 0 & 1 \\ 0 & 0 & 0 \\ 1 & 1 & 1 \end{bmatrix} & \begin{bmatrix} 1 & 2 & 1 & 1 & 0 \\ 1 & 0 & 1 & 0 & 1 \\ 0 & 1 & 0 & 1 & 0 \end{bmatrix} \begin{bmatrix} 3 & 1 & 0 & 0 & 1 \\ 1 & 0 & 1 & 0 & 0 \\ 1 & 0 & 0 & 1 & 0 \\ 2 & 0 & 0 & 1 & 0 \\ 1 & 1 & 1 & 0 & 1 \end{bmatrix} \\
\begin{pmatrix} 1 & 0 & 1 & 1 & 2 \\ 0 & 1 & 0 & 1 & 0 \\ 3 & 1 & 1 & 0 & 1 \\ 0 & 2 & 0 & 1 & 0 \end{pmatrix} \begin{bmatrix} 0 & 2 & 1 \\ 1 & 0 & 0 \\ 1 & 0 & 1 \\ 0 & 0 & 0 \\ 1 & 1 & 1 \end{bmatrix} & \begin{bmatrix} 1 & 0 & 1 & 1 & 2 \\ 0 & 1 & 0 & 1 & 0 \\ 3 & 1 & 1 & 0 & 1 \\ 0 & 2 & 0 & 1 & 0 \end{bmatrix} \begin{bmatrix} 3 & 1 & 0 & 0 & 1 \\ 1 & 0 & 1 & 0 & 0 \\ 1 & 0 & 0 & 1 & 0 \\ 2 & 0 & 0 & 1 & 0 \\ 1 & 1 & 1 & 0 & 1 \end{bmatrix}
\end{bmatrix}
$$

$$
\left[\begin{array}{c}
(1\ 0\ 1\ 1\ 1)\begin{bmatrix}3 & 1\\ 0 & 1\\ 1 & 0\\ 2 & 0\\ 1 & 1\end{bmatrix}\\
\begin{pmatrix}0 & 1 & 0 & 1 & 1\\ 1 & 1 & 0 & 1 & 0\end{pmatrix}\begin{bmatrix}3 & 1\\ 0 & 1\\ 1 & 0\\ 2 & 0\\ 1 & 1\end{bmatrix}\\
\begin{pmatrix}1 & 2 & 1 & 1 & 0\\ 1 & 0 & 1 & 0 & 1\\ 0 & 1 & 0 & 1 & 0\end{pmatrix}\begin{bmatrix}3 & 1\\ 0 & 1\\ 1 & 0\\ 2 & 0\\ 1 & 1\end{bmatrix}\\
\begin{pmatrix}1 & 0 & 1 & 1 & 2\\ 0 & 1 & 0 & 1 & 0\\ 3 & 1 & 1 & 0 & 1\\ 0 & 2 & 0 & 1 & 0\end{pmatrix}\begin{bmatrix}3 & 1\\ 0 & 1\\ 1 & 0\\ 2 & 0\\ 1 & 1\end{bmatrix}
\end{array}\right]
=\begin{bmatrix}1 & 6 & 9 & 4 & 4\\ 1 & 9 & 13 & 5 & 4\\ 0 & 10 & 15 & 6 & 5\\ 1 & 1 & 1 & 2 & 0\\ 1 & 2 & 3 & 1 & 2\end{bmatrix}\subset
$$

$$
\left[\begin{array}{cccc|cc|ccc}
2 & 3 & 1 & 1 & 1 & 2 & 4 & 1 & 3\\
1 & 1 & 1 & 1 & 1 & 0 & 1 & 1 & 1\\
1 & 1 & 2 & 1 & 1 & 1 & 1 & 1 & 1\\
\hline
1 & 3 & 3 & 1 & 1 & 2 & 3 & 1 & 3\\
4 & 2 & 1 & 2 & 2 & 3 & 4 & 2 & 2\\
3 & 2 & 2 & 2 & 2 & 2 & 3 & 2 & 2\\
3 & 0 & 1 & 2 & 2 & 1 & 1 & 2 & 0\\
\hline
0 & 3 & 2 & 0 & 0 & 2 & 3 & 0 & 3
\end{array}\right]\subset
$$

$$\begin{bmatrix} 2 & 3 & 3 & | & 7 & 2 & 1 & 2 & 2 & | & 7 & 2 \\ 2 & 1 & 1 & | & 4 & 1 & 2 & 1 & 1 & | & 3 & 2 \\ 1 & 2 & 1 & | & 6 & 1 & 1 & 1 & 1 & | & 5 & 2 \\ \hline 3 & 2 & 2 & | & 8 & 1 & 2 & 2 & 1 & | & 6 & 3 \\ 2 & 3 & 3 & | & 5 & 2 & 1 & 1 & 2 & | & 5 & 2 \\ 1 & 0 & 0 & | & 3 & 0 & 1 & 1 & 0 & | & 2 & 1 \\ \hline 3 & 4 & 4 & | & 8 & 3 & 2 & 2 & 3 & | & 8 & 3 \\ 1 & 0 & 0 & | & 3 & 0 & 1 & 1 & 0 & | & 2 & 1 \\ 3 & 7 & 5 & | & 12 & 4 & 2 & 1 & 4 & | & 11 & 5 \\ 2 & 0 & 0 & | & 4 & 0 & 2 & 1 & 0 & | & 2 & 2 \end{bmatrix}.$$

We see the major product yields a semi super trimatrix. Thus using major product of compatible special semi super trivector we can get more and more semi super trimatrices which are not trivectors. We give yet another example of how the major product is determined.

*Example 3.71:* Let

$$T = T_1 \cup T_2 \cup T_3$$

$$= \begin{bmatrix} 3 & 1 & 0 & 1 & 2 \\ 1 & 1 & 0 & 1 & 0 \\ 2 & 3 & 0 & 1 & 4 \end{bmatrix} \cup \begin{bmatrix} 2 & 1 & 5 \\ 1 & 0 & 1 \\ 5 & 1 & 1 \\ \hline 0 & 1 & 1 \\ 1 & 0 & 0 \\ 0 & 1 & 0 \\ 1 & 1 & 0 \\ \hline 0 & 1 & 2 \\ 2 & 1 & 0 \\ 3 & 1 & 4 \\ 0 & 0 & 1 \end{bmatrix} \cup \begin{bmatrix} 3 & 4 \\ 1 & 0 \\ 5 & 2 \\ \hline 3 & 0 \\ 1 & 1 \\ \hline 1 & 1 \\ 2 & 1 \\ \hline 3 & 1 \\ 4 & 5 \\ 1 & 6 \\ \hline 1 & 2 \end{bmatrix}$$

be a special semi super column trivector and

$$S = S_1 \cup S_2 \cup S_3$$

$$= \begin{bmatrix} 1 & 0 & 1 & 0 \\ 2 & 1 & 0 & 1 \\ 3 & 0 & 1 & 2 \\ 0 & 0 & 1 & 0 \\ 7 & 1 & 0 & 1 \end{bmatrix} \cup$$

$$\begin{bmatrix} 0 & 1 & 0 & 1 & | & 3 & 1 & 1 & | & 3 & 0 & 1 & 4 & 5 & 0 \\ 1 & 1 & 0 & 2 & | & 1 & 1 & 1 & | & 1 & 1 & 0 & 1 & 1 & 0 \\ 2 & 1 & 0 & 0 & | & 0 & 2 & 1 & | & 2 & 1 & 0 & 1 & 0 & 0 \end{bmatrix}$$

$$\cup \begin{bmatrix} 1 & 3 & | & 0 & 1 & | & 1 & 1 & 0 & | & 1 & 0 & 1 \\ 2 & 4 & | & 1 & 0 & | & 0 & 0 & 1 & | & 1 & 1 & 0 \end{bmatrix}$$

be a special super row trivector.

TS = $(T_1 \cup T_2 \cup T_3)(S_1 \cup S_2 \cup S_3)$
    = $T_1 S_1 \cup T_2 S_2 \cup T_3 S_3$

$$= \begin{bmatrix} 3 & 1 & 0 & 1 & 2 \\ 1 & 1 & 0 & 1 & 0 \\ 2 & 3 & 0 & 1 & 4 \end{bmatrix} \begin{bmatrix} 1 & 0 & 1 & 0 \\ 2 & 1 & 0 & 1 \\ 3 & 0 & 1 & 2 \\ 0 & 0 & 1 & 0 \\ 7 & 1 & 0 & 1 \end{bmatrix}$$

$$\cup \begin{bmatrix} 2 & 1 & 5 \\ 1 & 0 & 1 \\ 5 & 1 & 1 \\ \hline 0 & 1 & 1 \\ 1 & 0 & 0 \\ 0 & 1 & 0 \\ 1 & 1 & 0 \\ \hline 0 & 1 & 2 \\ 2 & 1 & 0 \\ 3 & 1 & 4 \\ 0 & 0 & 1 \end{bmatrix} \begin{bmatrix} 0 & 1 & 0 & 1 & 3 & 1 & 1 & 3 & 0 & 1 & 4 & 5 & 0 \\ 1 & 1 & 0 & 2 & 1 & 1 & 1 & 1 & 1 & 0 & 1 & 1 & 0 \\ 2 & 1 & 0 & 0 & 0 & 2 & 1 & 2 & 1 & 0 & 1 & 0 & 0 \end{bmatrix} \cup$$

$$\begin{bmatrix} 3 & 4 \\ 1 & 0 \\ 5 & 2 \\ \hline 3 & 0 \\ 1 & 1 \\ 1 & 1 \\ 2 & 1 \\ \hline 3 & 1 \\ 4 & 5 \\ 1 & 6 \\ 1 & 2 \end{bmatrix} \begin{bmatrix} 1 & 3 & 0 & 1 & 1 & 1 & 0 & 1 & 0 & 1 \\ 2 & 4 & 1 & 0 & 0 & 0 & 1 & 1 & 1 & 0 \end{bmatrix}$$

$$= \begin{bmatrix} 19 & 3 & 4 & 3 \\ 3 & 1 & 2 & 1 \\ 34 & 7 & 3 & 7 \end{bmatrix} \cup$$

$$\left[ \begin{array}{c|c} \begin{pmatrix} 2 & 1 & 5 \\ 1 & 0 & 1 \\ 5 & 1 & 1 \end{pmatrix} \begin{bmatrix} 0 & 1 & 0 & 1 \\ 1 & 1 & 0 & 2 \\ 2 & 1 & 0 & 0 \end{bmatrix} & \begin{pmatrix} 2 & 1 & 5 \\ 1 & 0 & 1 \\ 5 & 1 & 1 \end{pmatrix} \begin{pmatrix} 3 & 1 & 1 \\ 1 & 1 & 1 \\ 0 & 2 & 1 \end{pmatrix} \\ \hline \begin{pmatrix} 0 & 1 & 1 \\ 1 & 0 & 0 \\ 0 & 1 & 0 \\ 1 & 1 & 0 \end{pmatrix} \begin{bmatrix} 0 & 1 & 0 & 1 \\ 1 & 1 & 0 & 2 \\ 2 & 1 & 0 & 0 \end{bmatrix} & \begin{pmatrix} 0 & 1 & 1 \\ 1 & 0 & 0 \\ 0 & 1 & 0 \\ 1 & 1 & 0 \end{pmatrix} \begin{pmatrix} 3 & 1 & 1 \\ 1 & 1 & 1 \\ 0 & 2 & 1 \end{pmatrix} \\ \hline \begin{pmatrix} 0 & 1 & 2 \\ 2 & 1 & 0 \\ 3 & 1 & 4 \\ 0 & 0 & 1 \end{pmatrix} \begin{bmatrix} 0 & 1 & 0 & 1 \\ 1 & 1 & 0 & 2 \\ 2 & 1 & 0 & 0 \end{bmatrix} & \begin{pmatrix} 0 & 1 & 2 \\ 2 & 1 & 0 \\ 3 & 1 & 4 \\ 0 & 0 & 1 \end{pmatrix} \begin{pmatrix} 3 & 1 & 1 \\ 1 & 1 & 1 \\ 0 & 2 & 1 \end{pmatrix} \end{array} \right]$$

$$\subset \left[ \begin{array}{c} \begin{pmatrix} 2 & 1 & 5 \\ 1 & 0 & 1 \\ 5 & 1 & 1 \end{pmatrix} \begin{bmatrix} 3 & 0 & 1 & 4 & 5 & 0 \\ 1 & 1 & 0 & 1 & 1 & 0 \\ 2 & 1 & 0 & 1 & 0 & 0 \end{bmatrix} \\ \hline \begin{pmatrix} 0 & 1 & 1 \\ 1 & 0 & 0 \\ 0 & 1 & 0 \\ 1 & 1 & 0 \end{pmatrix} \begin{bmatrix} 3 & 0 & 1 & 4 & 5 & 0 \\ 1 & 1 & 0 & 1 & 1 & 0 \\ 2 & 1 & 0 & 1 & 0 & 0 \end{bmatrix} \\ \hline \begin{pmatrix} 0 & 1 & 2 \\ 2 & 1 & 0 \\ 3 & 1 & 4 \\ 0 & 0 & 1 \end{pmatrix} \begin{bmatrix} 3 & 0 & 1 & 4 & 5 & 0 \\ 1 & 1 & 0 & 1 & 1 & 0 \\ 2 & 1 & 0 & 1 & 0 & 0 \end{bmatrix} \end{array} \right]$$

$$\begin{bmatrix} \begin{bmatrix} 3 & 4 \\ 1 & 0 \\ 5 & 2 \end{bmatrix}\begin{bmatrix} 1 & 3 \\ 2 & 4 \end{bmatrix} & \begin{bmatrix} 3 & 4 \\ 1 & 0 \\ 5 & 2 \end{bmatrix}\begin{bmatrix} 0 & 1 \\ 1 & 0 \end{bmatrix} & \begin{bmatrix} 3 & 4 \\ 1 & 0 \\ 5 & 2 \end{bmatrix}\begin{bmatrix} 1 & 1 & 0 \\ 0 & 0 & 1 \end{bmatrix} & \begin{bmatrix} 3 & 4 \\ 1 & 0 \\ 5 & 2 \end{bmatrix}\begin{bmatrix} 1 & 0 & 1 \\ 1 & 1 & 0 \end{bmatrix} \\ \begin{bmatrix} 3 & 0 \\ 1 & 1 \end{bmatrix}\begin{bmatrix} 1 & 3 \\ 2 & 4 \end{bmatrix} & \begin{bmatrix} 3 & 0 \\ 1 & 1 \end{bmatrix}\begin{bmatrix} 0 & 1 \\ 1 & 0 \end{bmatrix} & \begin{bmatrix} 3 & 0 \\ 1 & 1 \end{bmatrix}\begin{bmatrix} 1 & 1 & 0 \\ 0 & 0 & 1 \end{bmatrix} & \begin{bmatrix} 3 & 0 \\ 1 & 1 \end{bmatrix}\begin{bmatrix} 1 & 0 & 1 \\ 1 & 1 & 0 \end{bmatrix} \\ \begin{bmatrix} 1 & 1 \\ 2 & 1 \\ 3 & 1 \\ 4 & 5 \\ 1 & 6 \end{bmatrix}\begin{bmatrix} 1 & 3 \\ 2 & 4 \end{bmatrix} & \begin{bmatrix} 1 & 1 \\ 2 & 1 \\ 3 & 1 \\ 4 & 5 \\ 1 & 6 \end{bmatrix}\begin{bmatrix} 0 & 1 \\ 1 & 0 \end{bmatrix} & \begin{bmatrix} 1 & 1 \\ 2 & 1 \\ 3 & 1 \\ 4 & 5 \\ 1 & 6 \end{bmatrix}\begin{bmatrix} 1 & 1 & 0 \\ 0 & 0 & 1 \end{bmatrix} & \begin{bmatrix} 1 & 1 \\ 2 & 1 \\ 3 & 1 \\ 4 & 5 \\ 1 & 6 \end{bmatrix}\begin{bmatrix} 1 & 0 & 1 \\ 1 & 1 & 0 \end{bmatrix} \\ \begin{bmatrix} 1 & 2 \end{bmatrix}\begin{bmatrix} 1 & 3 \\ 2 & 4 \end{bmatrix} & \begin{bmatrix} 1 & 2 \end{bmatrix}\begin{bmatrix} 0 & 1 \\ 1 & 0 \end{bmatrix} & \begin{bmatrix} 1 & 2 \end{bmatrix}\begin{bmatrix} 1 & 1 & 0 \\ 0 & 0 & 1 \end{bmatrix} & \begin{bmatrix} 1 & 2 \end{bmatrix}\begin{bmatrix} 1 & 0 & 1 \\ 1 & 1 & 0 \end{bmatrix} \end{bmatrix}$$

$$= \begin{bmatrix} 19 & 3 & 4 & 3 \\ 3 & 1 & 2 & 1 \\ 33 & 7 & 3 & 7 \end{bmatrix} \cup$$

$$\begin{bmatrix} 11 & 8 & 0 & 4 & 7 & 13 & 8 & 17 & 6 & 2 & 14 & 11 & 0 \\ 2 & 2 & 0 & 1 & 3 & 3 & 2 & 5 & 1 & 1 & 5 & 5 & 0 \\ 3 & 2 & 0 & 7 & 16 & 8 & 7 & 18 & 2 & 5 & 22 & 26 & 0 \\ \hline 3 & 2 & 0 & 2 & 1 & 3 & 2 & 3 & 2 & 0 & 2 & 1 & 0 \\ 0 & 1 & 0 & 1 & 3 & 1 & 1 & 3 & 0 & 1 & 4 & 5 & 0 \\ 1 & 1 & 0 & 2 & 1 & 1 & 1 & 1 & 1 & 0 & 1 & 1 & 0 \\ 1 & 2 & 0 & 3 & 4 & 2 & 2 & 4 & 1 & 1 & 5 & 6 & 0 \\ \hline 5 & 3 & 0 & 2 & 1 & 5 & 3 & 5 & 3 & 0 & 3 & 1 & 0 \\ 1 & 3 & 0 & 4 & 7 & 3 & 3 & 7 & 1 & 2 & 9 & 11 & 0 \\ 9 & 8 & 0 & 5 & 10 & 12 & 8 & 18 & 5 & 3 & 17 & 16 & 0 \\ 2 & 1 & 0 & 0 & 0 & 2 & 1 & 2 & 1 & 0 & 1 & 0 & 0 \end{bmatrix}$$

$$\cup \begin{bmatrix} 11 & 25 & 4 & 3 & 3 & 3 & 4 & 7 & 4 & 3 \\ 1 & 3 & 0 & 1 & 1 & 1 & 0 & 1 & 0 & 1 \\ 9 & 23 & 2 & 5 & 5 & 5 & 2 & 7 & 2 & 5 \\ \hline 3 & 9 & 0 & 3 & 3 & 3 & 0 & 3 & 0 & 3 \\ 3 & 7 & 1 & 1 & 1 & 1 & 1 & 2 & 1 & 1 \\ \hline 3 & 7 & 1 & 1 & 1 & 1 & 1 & 2 & 1 & 1 \\ 4 & 10 & 1 & 2 & 2 & 2 & 1 & 3 & 1 & 2 \\ 5 & 13 & 1 & 3 & 3 & 3 & 1 & 4 & 1 & 3 \\ 14 & 32 & 5 & 4 & 4 & 4 & 5 & 9 & 5 & 4 \\ 13 & 27 & 6 & 1 & 1 & 1 & 6 & 7 & 6 & 1 \\ \hline 5 & 11 & 2 & 1 & 1 & 1 & 2 & 3 & 2 & 1 \end{bmatrix}.$$

We see the resultant is a semi super trimatrix.

*Example 3.72:* Let

$$T = T_1 \cup T_2 \cup T_3$$

$$= \begin{bmatrix} 0 & 1 & 0 & 1 & 0 \\ 1 & 0 & 1 & 0 & 1 \\ 1 & 1 & 1 & 0 & 0 \\ 0 & 0 & 0 & 1 & 1 \end{bmatrix} \cup \begin{bmatrix} 3 & 1 & 0 & 1 \\ 1 & 0 & 1 & 2 \\ 1 & 1 & 1 & 0 \\ \hline 1 & 0 & 1 & 0 \\ 0 & 1 & 0 & 0 \\ \hline 1 & 0 & 1 & 1 \\ 0 & 1 & 1 & 0 \\ 1 & 0 & 0 & 1 \\ 1 & 0 & 1 & 0 \\ 1 & 0 & 0 & 0 \\ 0 & 1 & 0 & 0 \end{bmatrix}$$

$$\cup \begin{bmatrix} 0 & 1 & 3 & 1 & 2 \\ 1 & 0 & 1 & 0 & 0 \\ 0 & 1 & 0 & 0 & 1 \\ \hline 1 & 1 & 0 & 0 & 0 \\ 2 & 1 & 1 & 1 & 0 \\ 0 & 0 & 1 & 1 & 0 \\ \hline 1 & 0 & 1 & 0 & 0 \\ 1 & 1 & 0 & 1 & 1 \\ 0 & 1 & 0 & 1 & 0 \end{bmatrix}$$

be a special semi super column trivector.

$$\begin{aligned} T^T &= (T = T_1 \cup T_2 \cup T_3)^T \\ &= T_1^T \cup T_2^T \cup T_3^T \end{aligned}$$

$$= \begin{bmatrix} 0 & 1 & 1 & 0 \\ 1 & 0 & 1 & 0 \\ 0 & 1 & 1 & 0 \\ 1 & 0 & 0 & 1 \\ 0 & 1 & 0 & 1 \end{bmatrix} \cup$$

$$\begin{bmatrix} 3 & 1 & 1 & 1 & 0 & 1 & 0 & 1 & 1 & 1 & 0 \\ 1 & 0 & 1 & 0 & 1 & 0 & 1 & 0 & 0 & 0 & 1 \\ 0 & 1 & 1 & 1 & 0 & 1 & 1 & 0 & 1 & 0 & 0 \\ 1 & 2 & 0 & 0 & 0 & 1 & 0 & 1 & 0 & 0 & 0 \end{bmatrix} \cup$$

$$\begin{bmatrix} 0 & 1 & 0 & 1 & 2 & 0 & 1 & 1 & 0 \\ 1 & 0 & 1 & 1 & 1 & 0 & 0 & 1 & 1 \\ 3 & 1 & 0 & 0 & 1 & 1 & 1 & 0 & 0 \\ 1 & 0 & 0 & 0 & 1 & 1 & 0 & 1 & 1 \\ 2 & 0 & 1 & 0 & 0 & 0 & 0 & 1 & 0 \end{bmatrix}.$$

$$TT^T = (T_1 \cup T_2 \cup T_3)(T_1 \cup T_2 \cup T_3)^T$$

$$
\begin{aligned}
&= (T_1 \cup T_2 \cup T_3)(T_1^T \cup T_2^T \cup T_3^T) \\
&= T_1 T_1^T \cup T_2 T_2^T \cup T_3 T_3^T
\end{aligned}
$$

$$
= \begin{bmatrix} 0 & 1 & 0 & 1 & 0 \\ 1 & 0 & 1 & 0 & 1 \\ 1 & 1 & 1 & 0 & 0 \\ 0 & 0 & 0 & 1 & 1 \end{bmatrix} \begin{bmatrix} 0 & 1 & 1 & 0 \\ 1 & 0 & 1 & 0 \\ 0 & 1 & 1 & 0 \\ 1 & 0 & 0 & 1 \\ 0 & 1 & 0 & 1 \end{bmatrix}
$$

$$
\cup \left\{ \begin{bmatrix} 3 & 1 & 0 & 1 \\ 1 & 0 & 1 & 2 \\ 1 & 1 & 1 & 0 \\ \hline 1 & 0 & 1 & 0 \\ 0 & 1 & 0 & 0 \\ \hline 1 & 0 & 1 & 1 \\ 0 & 1 & 1 & 0 \\ 1 & 0 & 0 & 1 \\ 1 & 0 & 1 & 0 \\ 1 & 0 & 0 & 0 \\ 0 & 1 & 0 & 0 \end{bmatrix} \begin{bmatrix} 3 & 1 & 1 & 1 & 0 & 1 & 0 & 1 & 1 & 1 & 0 \\ 1 & 0 & 1 & 0 & 1 & 0 & 1 & 0 & 0 & 0 & 1 \\ 0 & 1 & 1 & 1 & 0 & 1 & 1 & 0 & 1 & 0 & 0 \\ 1 & 2 & 0 & 0 & 0 & 1 & 0 & 1 & 0 & 0 & 0 \end{bmatrix} \right\}
$$

$$
\cup \left\{ \begin{bmatrix} 0 & 1 & 3 & 1 & 2 \\ 1 & 0 & 1 & 0 & 0 \\ 0 & 1 & 0 & 0 & 1 \\ 1 & 1 & 0 & 0 & 0 \\ \hline 2 & 1 & 1 & 1 & 0 \\ 0 & 0 & 1 & 1 & 0 \\ \hline 1 & 0 & 1 & 0 & 0 \\ 1 & 1 & 0 & 1 & 1 \\ 0 & 1 & 0 & 1 & 0 \end{bmatrix} \begin{bmatrix} 0 & 1 & 0 & 1 & 2 & 0 & 1 & 1 & 0 \\ 1 & 0 & 1 & 1 & 1 & 0 & 0 & 1 & 1 \\ 3 & 1 & 0 & 0 & 1 & 1 & 1 & 0 & 0 \\ 1 & 0 & 0 & 0 & 1 & 1 & 0 & 1 & 1 \\ 2 & 0 & 1 & 0 & 0 & 0 & 0 & 1 & 0 \end{bmatrix} \right\}
$$

$$= \begin{bmatrix} 2 & 0 & 1 & 1 \\ 0 & 3 & 2 & 1 \\ 1 & 2 & 3 & 0 \\ 1 & 1 & 0 & 2 \end{bmatrix} \cup$$

$$\left[ \begin{array}{c|c} \begin{pmatrix} 3 & 1 & 0 & 1 \\ 1 & 0 & 1 & 2 \\ 1 & 1 & 1 & 0 \end{pmatrix} \begin{pmatrix} 3 & 1 & 1 \\ 1 & 0 & 1 \\ 0 & 1 & 1 \\ 1 & 2 & 0 \end{pmatrix} & \begin{pmatrix} 3 & 1 & 0 & 1 \\ 1 & 0 & 1 & 2 \\ 1 & 1 & 1 & 0 \end{pmatrix} \begin{pmatrix} 1 & 0 \\ 0 & 1 \\ 1 & 0 \\ 0 & 0 \end{pmatrix} \\ \hline \begin{pmatrix} 1 & 0 & 1 & 0 \\ 0 & 1 & 0 & 0 \end{pmatrix} \begin{pmatrix} 3 & 1 & 1 \\ 1 & 0 & 1 \\ 0 & 1 & 1 \\ 1 & 2 & 0 \end{pmatrix} & \begin{pmatrix} 1 & 0 & 1 & 0 \\ 0 & 1 & 0 & 0 \end{pmatrix} \begin{pmatrix} 1 & 0 \\ 0 & 1 \\ 1 & 0 \\ 0 & 0 \end{pmatrix} \\ \hline \begin{bmatrix} 1 & 0 & 1 & 1 \\ 0 & 1 & 1 & 0 \\ 1 & 0 & 0 & 1 \\ 1 & 0 & 1 & 0 \\ 1 & 0 & 0 & 0 \\ 0 & 1 & 0 & 0 \end{bmatrix} \begin{pmatrix} 3 & 1 & 1 \\ 1 & 0 & 1 \\ 0 & 1 & 1 \\ 1 & 2 & 0 \end{pmatrix} & \begin{bmatrix} 1 & 0 & 1 & 1 \\ 0 & 1 & 1 & 0 \\ 1 & 0 & 0 & 1 \\ 1 & 0 & 1 & 0 \\ 1 & 0 & 0 & 0 \\ 0 & 1 & 0 & 0 \end{bmatrix} \begin{pmatrix} 1 & 0 \\ 0 & 1 \\ 1 & 0 \\ 0 & 0 \end{pmatrix} \end{array} \right.$$

$$\left| \begin{pmatrix} 3 & 1 & 0 & 1 \\ 1 & 0 & 1 & 2 \\ 1 & 1 & 1 & 0 \end{pmatrix} \begin{bmatrix} 1 & 0 & 1 & 1 & 1 & 0 \\ 0 & 1 & 0 & 0 & 0 & 1 \\ 1 & 1 & 0 & 1 & 0 & 0 \\ 1 & 0 & 1 & 0 & 0 & 0 \end{bmatrix} \right|$$

$$\left| \begin{pmatrix} 1 & 0 & 1 & 0 \\ 0 & 1 & 0 & 0 \end{pmatrix} \begin{bmatrix} 1 & 0 & 1 & 1 & 1 & 0 \\ 0 & 1 & 0 & 0 & 0 & 1 \\ 1 & 1 & 0 & 1 & 0 & 0 \\ 1 & 0 & 1 & 0 & 0 & 0 \end{bmatrix} \right| \cup$$

$$\left[ \begin{bmatrix} 1 & 0 & 1 & 1 \\ 0 & 1 & 1 & 0 \\ 1 & 0 & 0 & 1 \\ 1 & 0 & 1 & 0 \\ 1 & 0 & 0 & 0 \\ 0 & 1 & 0 & 0 \end{bmatrix} \begin{bmatrix} 1 & 0 & 1 & 1 & 1 & 0 \\ 0 & 1 & 0 & 0 & 0 & 1 \\ 1 & 1 & 0 & 1 & 0 & 0 \\ 1 & 0 & 1 & 0 & 0 & 0 \end{bmatrix} \right]$$

$$\left[ \begin{pmatrix} 0 & 1 & 3 & 1 & 2 \\ 1 & 0 & 1 & 0 & 0 \\ 0 & 1 & 0 & 0 & 1 \\ 1 & 1 & 0 & 0 & 0 \end{pmatrix} \begin{pmatrix} 0 & 1 & 0 & 1 \\ 1 & 0 & 1 & 1 \\ 3 & 1 & 0 & 0 \\ 1 & 0 & 0 & 0 \\ 2 & 0 & 1 & 0 \end{pmatrix} \middle| \begin{pmatrix} 0 & 1 & 3 & 1 & 2 \\ 1 & 0 & 1 & 0 & 0 \\ 0 & 1 & 0 & 0 & 1 \\ 1 & 1 & 0 & 0 & 0 \end{pmatrix} \begin{bmatrix} 2 & 0 \\ 1 & 0 \\ 1 & 1 \\ 1 & 1 \\ 0 & 0 \end{bmatrix} \right.$$

$$\left| \begin{pmatrix} 2 & 1 & 1 & 1 & 0 \\ 0 & 0 & 1 & 1 & 0 \end{pmatrix} \begin{pmatrix} 0 & 1 & 0 & 1 \\ 1 & 0 & 1 & 1 \\ 3 & 1 & 0 & 0 \\ 1 & 0 & 0 & 0 \\ 2 & 0 & 1 & 0 \end{pmatrix} \right| \begin{pmatrix} 2 & 1 & 1 & 1 & 0 \\ 0 & 0 & 1 & 1 & 0 \end{pmatrix} \begin{bmatrix} 2 & 0 \\ 1 & 0 \\ 1 & 1 \\ 1 & 1 \\ 0 & 0 \end{bmatrix}$$

$$\left. \begin{pmatrix} 1 & 0 & 1 & 0 & 0 \\ 1 & 1 & 0 & 1 & 1 \\ 0 & 1 & 0 & 1 & 0 \end{pmatrix} \begin{pmatrix} 0 & 1 & 0 & 1 \\ 1 & 0 & 1 & 1 \\ 3 & 1 & 0 & 0 \\ 1 & 0 & 0 & 0 \\ 2 & 0 & 1 & 0 \end{pmatrix} \middle| \begin{pmatrix} 1 & 0 & 1 & 0 & 0 \\ 1 & 1 & 0 & 1 & 1 \\ 0 & 1 & 0 & 1 & 0 \end{pmatrix} \begin{bmatrix} 2 & 0 \\ 1 & 0 \\ 1 & 1 \\ 1 & 1 \\ 0 & 0 \end{bmatrix} \right]$$

$$\begin{bmatrix} \begin{pmatrix} 0 & 1 & 3 & 1 & 2 \\ 1 & 0 & 1 & 0 & 0 \\ 0 & 1 & 0 & 0 & 1 \\ 1 & 1 & 0 & 0 & 0 \end{pmatrix} \begin{pmatrix} 1 & 1 & 0 \\ 0 & 1 & 1 \\ 1 & 0 & 0 \\ 0 & 1 & 1 \\ 0 & 1 & 0 \end{pmatrix} \\ \hline \begin{pmatrix} 2 & 1 & 1 & 1 & 0 \\ 0 & 0 & 1 & 1 & 0 \end{pmatrix} \begin{pmatrix} 1 & 1 & 0 \\ 0 & 1 & 1 \\ 1 & 0 & 0 \\ 0 & 1 & 1 \\ 0 & 1 & 0 \end{pmatrix} \\ \hline \begin{pmatrix} 1 & 0 & 1 & 0 & 0 \\ 1 & 1 & 0 & 1 & 1 \\ 0 & 1 & 0 & 1 & 0 \end{pmatrix} \begin{pmatrix} 1 & 1 & 0 \\ 0 & 1 & 1 \\ 1 & 0 & 0 \\ 0 & 1 & 1 \\ 0 & 1 & 0 \end{pmatrix} \end{bmatrix} = \begin{bmatrix} 2 & 0 & 1 & 1 \\ 0 & 3 & 2 & 1 \\ 1 & 2 & 3 & 0 \\ 1 & 1 & 0 & 2 \end{bmatrix} \cup$$

$$\begin{bmatrix} 11 & 5 & 4 & 3 & 1 & 4 & 1 & 4 & 3 & 3 & 1 \\ 5 & 6 & 2 & 2 & 0 & 4 & 1 & 3 & 2 & 1 & 0 \\ 4 & 2 & 3 & 2 & 1 & 2 & 2 & 1 & 2 & 1 & 1 \\ \hline 3 & 2 & 2 & 2 & 0 & 2 & 1 & 1 & 2 & 1 & 0 \\ 1 & 0 & 1 & 0 & 1 & 0 & 1 & 0 & 0 & 0 & 1 \\ \hline 4 & 4 & 2 & 2 & 0 & 3 & 1 & 2 & 2 & 1 & 0 \\ 1 & 1 & 2 & 1 & 1 & 1 & 2 & 0 & 1 & 0 & 1 \\ 4 & 3 & 1 & 1 & 0 & 2 & 0 & 2 & 1 & 1 & 0 \\ 3 & 2 & 2 & 2 & 0 & 2 & 1 & 1 & 2 & 1 & 1 \\ 3 & 1 & 1 & 1 & 1 & 1 & 0 & 1 & 1 & 1 & 0 \\ 1 & 0 & 1 & 0 & 1 & 0 & 1 & 0 & 1 & 0 & 1 \end{bmatrix} \cup$$

$$\begin{bmatrix} 15 & 3 & 3 & 1 & 5 & 4 & 3 & 4 & 2 \\ 3 & 2 & 0 & 1 & 3 & 1 & 2 & 1 & 0 \\ 3 & 0 & 2 & 1 & 1 & 0 & 0 & 2 & 1 \\ 1 & 1 & 1 & 2 & 3 & 0 & 1 & 2 & 1 \\ \hline 5 & 3 & 1 & 3 & 7 & 2 & 3 & 4 & 2 \\ 4 & 1 & 0 & 0 & 2 & 2 & 1 & 1 & 1 \\ \hline 3 & 2 & 0 & 1 & 3 & 1 & 2 & 1 & 0 \\ 4 & 1 & 2 & 2 & 4 & 1 & 1 & 4 & 2 \\ 2 & 0 & 1 & 1 & 2 & 1 & 0 & 2 & 2 \end{bmatrix}.$$

We see the resultant is a symmetric semi super trimatrix. Thus the major product of a special semi super column trimatrix T with its transpose $T^T$ yields a symmetric semi super trimatrix.

***Example 3.73:*** Let $P = P_1 \cup P_2 \cup P_3$ be a special semi super row trimatrix. To find the value of $P^T P$. Given

$$P = \begin{bmatrix} 1 & 0 & 1 & 1 \\ 0 & 1 & 0 & 1 \\ 1 & 0 & 1 & 0 \\ 0 & 1 & 1 & 0 \\ 1 & 1 & 0 & 0 \end{bmatrix} \cup$$

$$\begin{bmatrix} 3 & 1 & 1 & 1 & 0 & 1 & 0 & 1 & 1 & 0 & 1 \\ 1 & 0 & 1 & 0 & 0 & 0 & 1 & 0 & 1 & 0 & 0 \\ 0 & 1 & 0 & 0 & 1 & 1 & 1 & 1 & 0 & 0 & 0 \\ 1 & 3 & 0 & 1 & 1 & 0 & 1 & 0 & 1 & 0 & 1 \end{bmatrix} \cup$$

$$\begin{bmatrix} 1 & 1 & 1 & 1 & 1 & 1 & 0 & 1 & 1 & 1 \\ 0 & 1 & 0 & 0 & 1 & 0 & 1 & 0 & 0 & 1 \\ 1 & 0 & 1 & 1 & 0 & 1 & 0 & 1 & 0 & 0 \end{bmatrix}$$

$$P^T = \begin{bmatrix} 1 & 0 & 1 & 0 & 1 \\ 0 & 1 & 0 & 1 & 1 \\ 1 & 0 & 1 & 1 & 0 \\ 1 & 1 & 0 & 0 & 0 \end{bmatrix} \cup$$

$$\begin{bmatrix} 3 & 1 & 0 & 1 \\ 1 & 0 & 1 & 3 \\ \hline 1 & 1 & 0 & 0 \\ 1 & 0 & 0 & 1 \\ 0 & 0 & 1 & 1 \\ \hline 1 & 0 & 1 & 0 \\ 0 & 1 & 1 & 1 \\ 1 & 0 & 1 & 0 \\ 1 & 1 & 0 & 1 \\ 0 & 0 & 0 & 0 \\ 1 & 0 & 0 & 1 \end{bmatrix} \cup \begin{bmatrix} 1 & 0 & 1 \\ 1 & 1 & 0 \\ 1 & 0 & 1 \\ \hline 1 & 0 & 1 \\ 1 & 1 & 0 \\ 1 & 0 & 1 \\ \hline 0 & 1 & 0 \\ 1 & 0 & 1 \\ 1 & 0 & 0 \\ 1 & 1 & 0 \end{bmatrix}.$$

$$\begin{aligned}
P^T P &= (P_1 \cup P_2 \cup P_3)^T (P_1 \cup P_2 \cup P_3) \\
&= (P_1^T \cup P_2^T \cup P_3^T)(P_1 \cup P_2 \cup P_3) \\
&= P_1^T P_1 \cup P_2^T P_2 \cup P_3 P_3^T
\end{aligned}$$

$$= \begin{bmatrix} 1 & 0 & 1 & 0 & 1 \\ 0 & 1 & 0 & 1 & 1 \\ 1 & 0 & 1 & 1 & 0 \\ 1 & 1 & 0 & 0 & 0 \end{bmatrix} \begin{bmatrix} 1 & 0 & 1 & 1 \\ 0 & 1 & 0 & 1 \\ 1 & 0 & 1 & 0 \\ 0 & 1 & 1 & 0 \\ 1 & 1 & 0 & 0 \end{bmatrix} \cup$$

$$\begin{bmatrix} 3 & 1 & 0 & 1 \\ 1 & 0 & 1 & 3 \\ \hline 1 & 1 & 0 & 0 \\ 1 & 0 & 0 & 1 \\ 0 & 0 & 1 & 1 \\ \hline 1 & 0 & 1 & 0 \\ 0 & 1 & 1 & 1 \\ 1 & 0 & 1 & 0 \\ 1 & 1 & 0 & 1 \\ 0 & 0 & 0 & 0 \\ 1 & 0 & 0 & 1 \end{bmatrix} \begin{bmatrix} 3 & 1 & 1 & 1 & 0 & 1 & 0 & 1 & 1 & 0 & 1 \\ 1 & 0 & 1 & 0 & 0 & 0 & 1 & 0 & 1 & 0 & 0 \\ 0 & 1 & 0 & 0 & 1 & 1 & 1 & 1 & 0 & 0 & 0 \\ 1 & 3 & 0 & 1 & 1 & 0 & 1 & 0 & 1 & 0 & 1 \end{bmatrix} \cup$$

$$\begin{bmatrix} 1 & 0 & 1 \\ 1 & 1 & 0 \\ 1 & 0 & 1 \\ \hline 1 & 0 & 1 \\ 1 & 1 & 0 \\ \hline 1 & 0 & 1 \\ 0 & 1 & 0 \\ 1 & 0 & 1 \\ 1 & 0 & 0 \\ 1 & 1 & 0 \end{bmatrix} \begin{bmatrix} 1 & 1 & 1 & 1 & 1 & 1 & 0 & 1 & 1 & 1 \\ 0 & 1 & 0 & 0 & 1 & 0 & 1 & 0 & 0 & 1 \\ 1 & 0 & 1 & 1 & 0 & 1 & 0 & 1 & 0 & 0 \end{bmatrix}$$

$$= \begin{bmatrix} 3 & 1 & 2 & 1 \\ 1 & 3 & 1 & 1 \\ 2 & 1 & 3 & 1 \\ 1 & 1 & 1 & 2 \end{bmatrix} \cup$$

$$\left[ \begin{array}{c|c} \begin{pmatrix} 3 & 1 & 0 & 1 \\ 1 & 0 & 1 & 3 \end{pmatrix} \begin{bmatrix} 3 & 1 \\ 1 & 0 \\ 0 & 1 \\ 1 & 3 \end{bmatrix} & \begin{pmatrix} 3 & 1 & 0 & 1 \\ 1 & 0 & 1 & 3 \end{pmatrix} \begin{bmatrix} 1 & 1 & 0 \\ 1 & 0 & 0 \\ 0 & 0 & 1 \\ 0 & 1 & 1 \end{bmatrix} \\ \hline \begin{pmatrix} 1 & 1 & 0 & 0 \\ 1 & 0 & 0 & 1 \\ 0 & 0 & 1 & 1 \end{pmatrix} \begin{bmatrix} 3 & 1 \\ 1 & 0 \\ 0 & 1 \\ 1 & 3 \end{bmatrix} & \begin{pmatrix} 1 & 1 & 0 & 0 \\ 1 & 0 & 0 & 1 \\ 0 & 0 & 1 & 1 \end{pmatrix} \begin{bmatrix} 1 & 1 & 0 \\ 1 & 0 & 0 \\ 0 & 0 & 1 \\ 0 & 1 & 1 \end{bmatrix} \\ \hline \begin{pmatrix} 1 & 0 & 1 & 0 \\ 0 & 1 & 1 & 1 \\ 1 & 0 & 1 & 0 \\ 1 & 1 & 0 & 1 \\ 0 & 0 & 0 & 0 \\ 1 & 0 & 0 & 1 \end{pmatrix} \begin{bmatrix} 3 & 1 \\ 1 & 0 \\ 0 & 1 \\ 1 & 3 \end{bmatrix} & \begin{pmatrix} 1 & 0 & 1 & 0 \\ 0 & 1 & 1 & 1 \\ 1 & 0 & 1 & 0 \\ 1 & 1 & 0 & 1 \\ 0 & 0 & 0 & 0 \\ 1 & 0 & 0 & 1 \end{pmatrix} \begin{bmatrix} 1 & 1 & 0 \\ 1 & 0 & 0 \\ 0 & 0 & 1 \\ 0 & 1 & 1 \end{bmatrix} \end{array} \right]$$

$$\subset \left[ \begin{array}{c} \begin{pmatrix} 3 & 1 & 0 & 1 \\ 1 & 0 & 1 & 3 \end{pmatrix} \begin{bmatrix} 1 & 0 & 1 & 1 & 0 & 1 \\ 0 & 1 & 0 & 1 & 0 & 0 \\ 1 & 1 & 1 & 0 & 0 & 0 \\ 0 & 1 & 0 & 1 & 0 & 1 \end{bmatrix} \\ \hline \begin{pmatrix} 1 & 1 & 0 & 0 \\ 1 & 0 & 0 & 1 \\ 0 & 0 & 1 & 1 \end{pmatrix} \begin{bmatrix} 1 & 0 & 1 & 1 & 0 & 1 \\ 0 & 1 & 0 & 1 & 0 & 0 \\ 1 & 1 & 1 & 0 & 0 & 0 \\ 0 & 1 & 0 & 1 & 0 & 1 \end{bmatrix} \\ \hline \begin{pmatrix} 1 & 0 & 1 & 0 \\ 0 & 1 & 1 & 1 \\ 1 & 0 & 1 & 0 \\ 1 & 1 & 0 & 1 \\ 0 & 0 & 0 & 0 \\ 1 & 0 & 0 & 1 \end{pmatrix} \begin{bmatrix} 1 & 0 & 1 & 1 & 0 & 1 \\ 0 & 1 & 0 & 1 & 0 & 0 \\ 1 & 1 & 1 & 0 & 0 & 0 \\ 0 & 1 & 0 & 1 & 0 & 1 \end{bmatrix} \end{array} \right]$$

$$
\begin{bmatrix}
\begin{bmatrix} 1 & 0 & 1 \\ 1 & 1 & 0 \\ 1 & 0 & 1 \end{bmatrix}\begin{bmatrix} 1 & 1 & 1 \\ 0 & 1 & 0 \\ 1 & 0 & 1 \end{bmatrix} & \begin{bmatrix} 1 & 0 & 1 \\ 1 & 1 & 0 \\ 1 & 0 & 1 \end{bmatrix}\begin{bmatrix} 1 & 1 \\ 0 & 1 \\ 1 & 0 \end{bmatrix} & \begin{bmatrix} 1 & 0 & 1 \\ 1 & 1 & 0 \\ 1 & 0 & 1 \end{bmatrix}\begin{bmatrix} 1 & 0 & 1 & 1 & 1 \\ 0 & 1 & 0 & 0 & 1 \\ 1 & 0 & 1 & 0 & 0 \end{bmatrix} \\
\begin{bmatrix} 1 & 0 & 1 \\ 1 & 1 & 0 \end{bmatrix}\begin{bmatrix} 1 & 1 & 1 \\ 0 & 1 & 0 \\ 1 & 0 & 1 \end{bmatrix} & \begin{bmatrix} 1 & 0 & 1 \\ 1 & 1 & 0 \end{bmatrix}\begin{bmatrix} 1 & 1 \\ 0 & 1 \\ 1 & 0 \end{bmatrix} & \begin{bmatrix} 1 & 0 & 1 \\ 1 & 1 & 0 \end{bmatrix}\begin{bmatrix} 1 & 0 & 1 & 1 & 1 \\ 0 & 1 & 0 & 0 & 1 \\ 1 & 0 & 1 & 0 & 0 \end{bmatrix} \\
\begin{bmatrix} 1 & 0 & 1 \\ 0 & 1 & 0 \\ 1 & 0 & 1 \\ 1 & 0 & 0 \\ 1 & 1 & 0 \end{bmatrix}\begin{bmatrix} 1 & 1 & 1 \\ 0 & 1 & 0 \\ 1 & 0 & 1 \end{bmatrix} & \begin{bmatrix} 1 & 0 & 1 \\ 0 & 1 & 0 \\ 1 & 0 & 1 \\ 1 & 0 & 0 \\ 1 & 1 & 0 \end{bmatrix}\begin{bmatrix} 1 & 1 \\ 0 & 1 \\ 1 & 0 \end{bmatrix} & \begin{bmatrix} 1 & 0 & 1 \\ 0 & 1 & 0 \\ 1 & 0 & 1 \\ 1 & 0 & 0 \\ 1 & 1 & 0 \end{bmatrix}\begin{bmatrix} 1 & 0 & 1 & 1 & 1 \\ 0 & 1 & 0 & 0 & 1 \\ 1 & 0 & 1 & 0 & 0 \end{bmatrix}
\end{bmatrix}
$$

$$
= \begin{bmatrix} 3 & 1 & 2 & 1 \\ 1 & 3 & 1 & 1 \\ 2 & 1 & 3 & 1 \\ 1 & 1 & 1 & 2 \end{bmatrix} \cup
$$

$$
\begin{bmatrix}
11 & 6 & 4 & 4 & 1 & 3 & 2 & 3 & 5 & 0 & 4 \\
6 & 11 & 1 & 4 & 4 & 2 & 4 & 2 & 4 & 0 & 4 \\
4 & 1 & 2 & 1 & 0 & 1 & 1 & 1 & 2 & 0 & 1 \\
4 & 4 & 1 & 2 & 1 & 1 & 1 & 1 & 2 & 0 & 2 \\
1 & 4 & 0 & 1 & 2 & 1 & 2 & 1 & 1 & 0 & 1 \\
3 & 2 & 1 & 1 & 1 & 2 & 1 & 2 & 1 & 0 & 1 \\
2 & 4 & 1 & 1 & 2 & 1 & 3 & 1 & 2 & 0 & 1 \\
3 & 2 & 1 & 1 & 1 & 2 & 1 & 2 & 1 & 0 & 1 \\
5 & 4 & 2 & 2 & 1 & 1 & 2 & 1 & 3 & 0 & 2 \\
0 & 0 & 0 & 0 & 0 & 0 & 0 & 0 & 0 & 0 & 0 \\
4 & 4 & 1 & 2 & 1 & 1 & 1 & 1 & 2 & 0 & 2
\end{bmatrix} \cup
$$

$$\begin{bmatrix} 2 & 1 & 2 & 2 & 1 & 2 & 0 & 2 & 1 & 1 \\ 1 & 2 & 1 & 1 & 2 & 1 & 1 & 1 & 1 & 2 \\ 2 & 1 & 2 & 2 & 1 & 2 & 0 & 2 & 1 & 1 \\ \hline 2 & 1 & 2 & 2 & 1 & 2 & 0 & 2 & 1 & 1 \\ 1 & 2 & 1 & 1 & 2 & 1 & 1 & 1 & 1 & 2 \\ \hline 2 & 1 & 2 & 2 & 1 & 2 & 0 & 2 & 1 & 1 \\ 0 & 1 & 0 & 0 & 1 & 0 & 1 & 0 & 0 & 1 \\ 2 & 1 & 2 & 2 & 1 & 2 & 0 & 2 & 1 & 1 \\ 1 & 1 & 1 & 1 & 1 & 1 & 0 & 1 & 1 & 1 \\ 1 & 2 & 1 & 1 & 2 & 1 & 1 & 1 & 1 & 2 \end{bmatrix}.$$

Clearly $P^TP$ is again a symmetric semi super trimatrix. Thus the major product of a transpose of a special semi supermatrix P with itself yields a symmetric semi supermatrix. Now we proceed on to define the notion of major product of semi supermatrices.

***Example 3.74:*** Let $P = P_1 \cup P_2 \cup P_3$ and $T = T_1 \cup T_2 \cup T_3$ be two semi super trimatrices for which major product PT is defined. We find this product for P and T.

Given

$$P = \begin{bmatrix} 3 & 1 & 0 & 1 \\ 1 & 2 & 0 & 2 \\ 0 & 1 & 6 & 3 \\ 1 & 0 & 0 & 0 \end{bmatrix} \cup \begin{bmatrix} 1 & 1 & 1 & 0 & 3 & 0 & 1 \\ \hline 0 & 0 & 1 & 0 & 1 & 1 & 1 \\ 3 & 1 & 0 & 0 & 1 & 0 & 1 \\ 1 & 0 & 0 & 1 & 0 & 0 & 1 \\ 1 & 1 & 0 & 1 & 0 & 0 & 0 \\ \hline 0 & 1 & 0 & 1 & 0 & 1 & 1 \\ 1 & 1 & 1 & 0 & 1 & 0 & 2 \end{bmatrix} \cup$$

$$\begin{bmatrix} 1 & 0 & 1 & 1 & 0 & 1 \\ 0 & 1 & 0 & 1 & 0 & 1 \\ 1 & 1 & 1 & 0 & 0 & 0 \\ \hline 0 & 1 & 0 & 1 & 0 & 1 \\ 1 & 0 & 1 & 1 & 1 & 0 \\ \hline 0 & 1 & 0 & 0 & 0 & 1 \\ 0 & 0 & 1 & 0 & 1 & 0 \end{bmatrix}$$

is a semi super trimatrix.
Now

$$T = \begin{bmatrix} 1 & 2 & 1 \\ 0 & 1 & 0 \\ 1 & 0 & 0 \\ 1 & 0 & 1 \end{bmatrix} \cup \begin{bmatrix} 0 & 1 & 5 & 3 & 1 & 2 & 0 \\ \hline 0 & 1 & 0 & 0 & 1 & 1 & 2 \\ 1 & 0 & 1 & 1 & 0 & 1 & 1 \\ 1 & 1 & 0 & 1 & 0 & 0 & 1 \\ \hline 0 & 0 & 1 & 0 & 1 & 0 & 1 \\ 1 & 1 & 1 & 1 & 1 & 1 & 2 \\ 0 & 0 & 0 & 0 & 1 & 2 & 0 \end{bmatrix}$$

$$\cup \begin{bmatrix} 1 & 0 & 1 & 1 & 1 & 1 & 1 \\ 1 & 1 & 1 & 0 & 1 & 0 & 1 \\ 1 & 0 & 1 & 1 & 0 & 1 & 2 \\ 1 & 1 & 0 & 1 & 1 & 0 & 2 \\ \hline 0 & 1 & 1 & 2 & 1 & 1 & 1 \\ 2 & 0 & 1 & 1 & 0 & 0 & 1 \end{bmatrix}.$$

$$PT = (P_1 \cup P_2 \cup P_3)(T_1 \cup T_2 \cup T_3)$$

$$= P_1T_1 \cup P_2T_2 \cup P_3T_3$$

$$= \begin{bmatrix} 3 & 1 & 0 & 1 \\ 1 & 2 & 0 & 2 \\ 0 & 1 & 6 & 3 \\ 1 & 0 & 0 & 0 \end{bmatrix} \begin{bmatrix} 1 & 2 & 1 \\ 0 & 1 & 0 \\ 1 & 0 & 0 \\ 1 & 0 & 1 \end{bmatrix} \cup (Z_{ij}^1) \cup (Z_{ij}^2)$$

$$= \begin{bmatrix} 4 & 7 & 4 \\ 3 & 4 & 3 \\ 9 & 1 & 3 \\ 1 & 2 & 1 \end{bmatrix} \cup (Z_{ij}^1) \cup (Z_{ij}^2).$$

We calculate

$$(Z_{ij}^1) = P_2 T_2 = \left[\begin{array}{c|ccc|ccc} 1 & 1 & 1 & 0 & 3 & 0 & 1 \\ \hline 0 & 0 & 1 & 0 & 1 & 1 & 1 \\ 3 & 1 & 0 & 0 & 1 & 0 & 1 \\ 1 & 0 & 0 & 1 & 0 & 0 & 1 \\ 1 & 1 & 0 & 1 & 0 & 0 & 0 \\ \hline 0 & 1 & 0 & 1 & 0 & 1 & 1 \\ 1 & 1 & 1 & 0 & 1 & 0 & 2 \end{array}\right]$$

$$\left[\begin{array}{c|cccc|cc} 0 & 1 & 5 & 3 & 1 & 2 & 0 \\ \hline 0 & 1 & 0 & 0 & 1 & 1 & 2 \\ 1 & 0 & 1 & 1 & 0 & 1 & 1 \\ 1 & 1 & 0 & 1 & 0 & 0 & 1 \\ \hline 0 & 0 & 1 & 0 & 1 & 0 & 1 \\ 1 & 1 & 1 & 1 & 1 & 1 & 2 \\ 0 & 0 & 0 & 0 & 1 & 2 & 0 \end{array}\right]$$

$$Z_{11}^1 = [1 \mid 1\ 1\ 0 \mid 3\ 0\ 1] \begin{bmatrix} 0 \\ \overline{0} \\ 1 \\ 1 \\ \overline{0} \\ 1 \\ 0 \end{bmatrix}$$

$$= (1)\,(0) + (1\ 1\ 0) \begin{bmatrix} 0 \\ 1 \\ 1 \end{bmatrix} + (3\ 0\ 1) \begin{bmatrix} 0 \\ 1 \\ 0 \end{bmatrix}$$

$$= (0) + (1) + (0) = 1.$$

$$Z_{12}^1 = [1 \mid 1\ 1\ 0 \mid 3\ 0\ 1] \begin{bmatrix} 1 & 5 & 3 & 1 \\ \hline 1 & 0 & 0 & 1 \\ 0 & 1 & 1 & 0 \\ 1 & 0 & 1 & 0 \\ \hline 0 & 1 & 0 & 1 \\ 1 & 1 & 1 & 1 \\ 0 & 0 & 0 & 1 \end{bmatrix}$$

$$= [1]\,[1\ 5\ 3\ 1] + [1\ 1\ 0] \begin{bmatrix} 1 & 0 & 0 & 1 \\ 0 & 1 & 1 & 0 \\ 1 & 0 & 1 & 0 \end{bmatrix} + [3\ 0\ 1] \begin{bmatrix} 0 & 1 & 0 & 1 \\ 1 & 1 & 1 & 1 \\ 0 & 0 & 0 & 1 \end{bmatrix}$$

$$= [1\ 5\ 3\ 1] + [1\ 1\ 1\ 1] + [0\ 3\ 0\ 4]$$

$$= [2\ 9\ 4\ 6].$$

$$Z^1_{13} = [1 \mid 1\ 1\ 0 \mid 3\ 0\ 1] \begin{bmatrix} 2 & 0 \\ \hline 1 & 2 \\ 1 & 1 \\ 0 & 1 \\ \hline 0 & 1 \\ 1 & 2 \\ 2 & 0 \end{bmatrix}$$

$$= [1][2\ 0] + [1\ 1\ 0] \begin{bmatrix} 1 & 2 \\ 1 & 1 \\ 0 & 1 \end{bmatrix} + [3\ 0\ 1] \begin{bmatrix} 0 & 1 \\ 1 & 2 \\ 2 & 0 \end{bmatrix}$$

$$= [2\ 0] + [2\ 3] + [2\ 3]$$

$$= [6\ 6].$$

$$Z^1_{21} = \begin{bmatrix} 0 & 0 & 1 & 0 & 1 & 1 & 1 \\ 3 & 1 & 0 & 0 & 1 & 0 & 1 \\ 1 & 0 & 0 & 1 & 0 & 0 & 1 \\ 1 & 1 & 0 & 1 & 0 & 0 & 0 \end{bmatrix} \begin{bmatrix} 0 \\ \hline 0 \\ 1 \\ 1 \\ \hline 0 \\ 1 \\ 0 \end{bmatrix}$$

$$= \begin{bmatrix} 0 \\ 3 \\ 1 \\ 1 \end{bmatrix}[0] + \begin{bmatrix} 0 & 1 & 0 \\ 1 & 0 & 0 \\ 0 & 0 & 1 \\ 1 & 0 & 1 \end{bmatrix} \begin{bmatrix} 0 \\ 1 \\ 1 \end{bmatrix} + \begin{bmatrix} 1 & 1 & 1 \\ 1 & 0 & 1 \\ 0 & 0 & 1 \\ 0 & 0 & 0 \end{bmatrix} \begin{bmatrix} 0 \\ 1 \\ 0 \end{bmatrix}$$

$$= \begin{bmatrix} 0 \\ 0 \\ 0 \\ 0 \end{bmatrix} + \begin{bmatrix} 1 \\ 0 \\ 1 \\ 1 \end{bmatrix} + \begin{bmatrix} 1 \\ 0 \\ 0 \\ 0 \end{bmatrix} = \begin{bmatrix} 2 \\ 0 \\ 1 \\ 1 \end{bmatrix}.$$

$$Z_{22}^1 = \begin{bmatrix} 0 & 0 & 1 & 0 & 1 & 1 & 1 \\ 3 & 1 & 0 & 0 & 1 & 0 & 1 \\ 1 & 0 & 0 & 1 & 0 & 0 & 1 \\ 1 & 1 & 0 & 1 & 0 & 0 & 0 \end{bmatrix} \begin{bmatrix} 1 & 5 & 3 & 1 \\ 1 & 0 & 0 & 1 \\ 0 & 1 & 1 & 0 \\ 1 & 0 & 1 & 0 \\ 0 & 1 & 0 & 1 \\ 1 & 1 & 1 & 1 \\ 0 & 0 & 0 & 1 \end{bmatrix}$$

$$= \begin{bmatrix} 0 \\ 3 \\ 1 \\ 1 \end{bmatrix} \begin{bmatrix} 1 & 5 & 3 & 1 \end{bmatrix} + \begin{bmatrix} 0 & 1 & 0 \\ 1 & 0 & 0 \\ 0 & 0 & 1 \\ 1 & 0 & 1 \end{bmatrix} \begin{bmatrix} 1 & 0 & 0 & 1 \\ 0 & 1 & 1 & 0 \\ 1 & 0 & 1 & 0 \end{bmatrix} +$$

$$\begin{bmatrix} 1 & 1 & 1 \\ 1 & 0 & 1 \\ 0 & 0 & 1 \\ 0 & 0 & 0 \end{bmatrix} \begin{bmatrix} 0 & 1 & 0 & 1 \\ 1 & 1 & 1 & 1 \\ 0 & 0 & 0 & 1 \end{bmatrix}.$$

$$= \begin{bmatrix} 0 & 0 & 0 & 0 \\ 3 & 15 & 9 & 3 \\ 1 & 5 & 3 & 1 \\ 1 & 5 & 3 & 1 \end{bmatrix} + \begin{bmatrix} 0 & 1 & 1 & 0 \\ 1 & 0 & 0 & 1 \\ 1 & 0 & 1 & 0 \\ 2 & 0 & 1 & 1 \end{bmatrix} + \begin{bmatrix} 1 & 2 & 1 & 3 \\ 0 & 1 & 0 & 2 \\ 0 & 0 & 0 & 1 \\ 0 & 0 & 0 & 0 \end{bmatrix}$$

$$= \begin{bmatrix} 1 & 3 & 2 & 3 \\ 4 & 16 & 9 & 6 \\ 2 & 5 & 4 & 2 \\ 3 & 5 & 4 & 2 \end{bmatrix}.$$

$$Z_{23}^1 = \begin{bmatrix} 0 & 0 & 1 & 0 & 1 & 1 & 1 \\ 3 & 1 & 0 & 0 & 1 & 0 & 1 \\ 1 & 0 & 0 & 1 & 0 & 0 & 1 \\ 1 & 1 & 0 & 1 & 0 & 0 & 0 \end{bmatrix} \begin{bmatrix} 2 & 0 \\ \hline 1 & 2 \\ 1 & 1 \\ 0 & 1 \\ \hline 0 & 1 \\ 1 & 2 \\ 0 & 1 \end{bmatrix}$$

$$= \begin{bmatrix} 0 \\ 3 \\ 1 \\ 1 \end{bmatrix} [2 \; 0] + \begin{bmatrix} 0 & 1 & 0 \\ 1 & 0 & 0 \\ 0 & 0 & 1 \\ 1 & 0 & 1 \end{bmatrix} \begin{bmatrix} 1 & 2 \\ 1 & 1 \\ 0 & 1 \end{bmatrix} + \begin{bmatrix} 1 & 1 & 1 \\ 1 & 0 & 1 \\ 0 & 0 & 1 \\ 0 & 0 & 0 \end{bmatrix} \begin{bmatrix} 0 & 1 \\ 1 & 2 \\ 2 & 0 \end{bmatrix}$$

$$= \begin{bmatrix} 0 & 0 \\ 6 & 0 \\ 2 & 0 \\ 2 & 0 \end{bmatrix} + \begin{bmatrix} 1 & 1 \\ 1 & 2 \\ 0 & 1 \\ 1 & 3 \end{bmatrix} + \begin{bmatrix} 3 & 3 \\ 2 & 1 \\ 2 & 0 \\ 0 & 0 \end{bmatrix} = \begin{bmatrix} 4 & 4 \\ 9 & 3 \\ 4 & 1 \\ 3 & 3 \end{bmatrix}.$$

$$Z_{31}^1 = \begin{bmatrix} 0 & 1 & 0 & 1 & 0 & 1 & 1 \\ 1 & 1 & 1 & 0 & 1 & 0 & 2 \end{bmatrix} \begin{bmatrix} 0 \\ \hline 0 \\ 1 \\ 1 \\ \hline 0 \\ 1 \\ 0 \end{bmatrix}$$

$$= \begin{bmatrix} 0 \\ 1 \end{bmatrix}[0] + \begin{bmatrix} 1 & 0 & 1 \\ 1 & 1 & 0 \end{bmatrix}\begin{bmatrix} 0 \\ 1 \\ 1 \end{bmatrix} + \begin{bmatrix} 0 & 1 & 1 \\ 1 & 0 & 2 \end{bmatrix}\begin{bmatrix} 0 \\ 1 \\ 0 \end{bmatrix}$$

$$= \begin{bmatrix} 0 \\ 0 \end{bmatrix} + \begin{bmatrix} 1 \\ 1 \end{bmatrix} + \begin{bmatrix} 1 \\ 0 \end{bmatrix} = \begin{bmatrix} 2 \\ 1 \end{bmatrix}.$$

$$Z^1_{32} = \begin{bmatrix} 0 & | & 1 & 0 & 1 & | & 0 & 1 & 1 \\ 1 & | & 1 & 1 & 0 & | & 1 & 0 & 2 \end{bmatrix} \begin{bmatrix} 1 & 5 & 3 & 1 \\ \hline 1 & 0 & 0 & 1 \\ 0 & 1 & 1 & 0 \\ 1 & 0 & 1 & 0 \\ \hline 0 & 1 & 0 & 1 \\ 1 & 1 & 1 & 1 \\ 1 & 0 & 1 & 0 \end{bmatrix}$$

$$= \begin{bmatrix} 0 \\ 1 \end{bmatrix}\begin{bmatrix} 1 & 5 & 3 & 1 \end{bmatrix} +$$

$$\begin{bmatrix} 1 & 0 & 1 \\ 1 & 1 & 0 \end{bmatrix}\begin{bmatrix} 1 & 0 & 0 & 1 \\ 0 & 1 & 1 & 0 \\ 1 & 0 & 1 & 0 \end{bmatrix} + \begin{bmatrix} 0 & 1 & 1 \\ 1 & 0 & 2 \end{bmatrix}\begin{bmatrix} 0 & 1 & 0 & 1 \\ 1 & 1 & 1 & 1 \\ 0 & 0 & 0 & 1 \end{bmatrix}$$

$$= \begin{bmatrix} 0 & 0 & 0 & 0 \\ 1 & 5 & 3 & 1 \end{bmatrix} + \begin{bmatrix} 2 & 0 & 1 & 1 \\ 1 & 1 & 1 & 1 \end{bmatrix} + \begin{bmatrix} 1 & 1 & 1 & 2 \\ 0 & 1 & 0 & 3 \end{bmatrix}$$

$$= \begin{bmatrix} 3 & 1 & 2 & 3 \\ 2 & 7 & 4 & 5 \end{bmatrix}.$$

$$Z_{33}^1 = \begin{bmatrix} 0 & | & 1 & 0 & 1 & | & 0 & 1 & 1 \\ 1 & | & 1 & 1 & 0 & | & 1 & 0 & 2 \end{bmatrix} \begin{bmatrix} 2 & 0 \\ \hline 1 & 2 \\ 1 & 1 \\ 0 & 1 \\ \hline 0 & 1 \\ 1 & 2 \\ 2 & 0 \end{bmatrix}$$

$$= \begin{bmatrix} 0 \\ 1 \end{bmatrix} [2 \ 0] + \begin{bmatrix} 1 & 0 & 1 \\ 1 & 1 & 0 \end{bmatrix} \begin{bmatrix} 1 & 2 \\ 1 & 1 \\ 0 & 1 \end{bmatrix} + \begin{bmatrix} 0 & 1 & 1 \\ 1 & 0 & 2 \end{bmatrix} \begin{bmatrix} 0 & 1 \\ 1 & 2 \\ 2 & 0 \end{bmatrix}$$

$$= \begin{bmatrix} 0 & 0 \\ 2 & 0 \end{bmatrix} + \begin{bmatrix} 1 & 3 \\ 2 & 3 \end{bmatrix} + \begin{bmatrix} 3 & 2 \\ 4 & 1 \end{bmatrix} = \begin{bmatrix} 4 & 5 \\ 8 & 4 \end{bmatrix}.$$

$$Z_{ij}^1 = \begin{bmatrix} Z_{11} & Z_{12} & Z_{13} \\ Z_{21} & Z_{22} & Z_{23} \\ Z_{31} & Z_{32} & Z_{33} \end{bmatrix}$$

$$= \begin{bmatrix} 1 & | & 2 & 9 & 4 & 6 & | & 6 & 6 \\ 2 & | & 1 & 3 & 2 & 3 & | & 4 & 4 \\ 0 & | & 4 & 16 & 9 & 6 & | & 9 & 3 \\ \hline 1 & | & 2 & 5 & 4 & 2 & | & 4 & 1 \\ 1 & | & 3 & 5 & 4 & 2 & | & 3 & 3 \\ \hline 2 & | & 3 & 1 & 2 & 3 & | & 4 & 5 \\ 1 & | & 2 & 7 & 4 & 5 & | & 8 & 4 \end{bmatrix}.$$

In the similar manner $Z_{ij}^2$ is also calculated. We see the major product of two semi super trimatrix yields a semi super trimatrix provided the product is defined.

Now having defined this we can also as in case of special semi super trivectors define the product of a semi super trimatrix with its transpose.

***Example 3.75:*** Let $T = T_1 \cup T_2 \cup T_3$ be a semi super trimatrix where

$$T = \begin{bmatrix} 2 & 1 & 0 & 0 & 0 \\ 0 & 1 & 0 & 1 & 0 \\ 1 & 0 & 0 & 0 & 1 \\ 0 & 1 & 1 & 0 & 0 \\ 3 & 0 & 1 & 0 & 2 \end{bmatrix} \cup \left[ \begin{array}{c|cccc|cc} 1 & 0 & 1 & 1 & 2 & 3 & 4 \\ 0 & 1 & 0 & 0 & 1 & 0 & 1 \\ 1 & 0 & 1 & 0 & 1 & 1 & 0 \\ 3 & 1 & 0 & 1 & 0 & 1 & 1 \\ 0 & 1 & 1 & 1 & 0 & 0 & 0 \\ 1 & 0 & 1 & 0 & 1 & 1 & 1 \\ \hline 1 & 0 & 0 & 1 & 1 & 1 & 0 \\ 0 & 1 & 0 & 0 & 1 & 0 & 1 \end{array} \right]$$

$$\cup \left[ \begin{array}{c|cccc|cc} 1 & 1 & 3 & 1 & 0 & 1 & 0 \\ 2 & 0 & 0 & 1 & 1 & 0 & 1 \\ \hline 0 & 1 & 0 & 1 & 1 & 1 & 1 \\ 1 & 1 & 0 & 1 & 0 & 1 & 1 \\ 1 & 0 & 1 & 1 & 0 & 0 & 1 \\ 0 & 0 & 1 & 0 & 0 & 1 & 0 \\ \hline 3 & 1 & 0 & 3 & 1 & 1 & 0 \end{array} \right].$$

$$T^T = \begin{bmatrix} 2 & 0 & 1 & 0 & 3 \\ 1 & 1 & 0 & 1 & 0 \\ 0 & 0 & 0 & 1 & 1 \\ 0 & 1 & 0 & 0 & 0 \\ 0 & 0 & 1 & 0 & 2 \end{bmatrix}$$

$$\cup \begin{bmatrix} 1 & 0 & 1 & 3 & 0 & 1 & 1 & 0 \\ \hline 0 & 1 & 0 & 1 & 1 & 0 & 0 & 1 \\ 1 & 0 & 1 & 0 & 1 & 1 & 0 & 0 \\ 1 & 0 & 0 & 1 & 1 & 0 & 1 & 0 \\ 2 & 1 & 1 & 0 & 0 & 1 & 1 & 1 \\ \hline 3 & 0 & 1 & 1 & 0 & 1 & 1 & 0 \\ 4 & 1 & 0 & 1 & 0 & 1 & 0 & 1 \end{bmatrix} \cup \begin{bmatrix} 1 & 2 & 0 & 1 & 1 & 0 & 3 \\ \hline 1 & 0 & 1 & 1 & 0 & 0 & 1 \\ 3 & 0 & 0 & 0 & 1 & 1 & 0 \\ 1 & 1 & 1 & 1 & 1 & 0 & 3 \\ 0 & 1 & 1 & 0 & 0 & 0 & 1 \\ \hline 1 & 0 & 1 & 1 & 0 & 1 & 1 \\ 0 & 1 & 1 & 1 & 1 & 0 & 0 \end{bmatrix}$$

$$TT^T = \begin{bmatrix} 2 & 1 & 0 & 0 & 0 \\ 0 & 1 & 0 & 1 & 0 \\ 1 & 0 & 0 & 0 & 1 \\ 0 & 1 & 1 & 0 & 0 \\ 3 & 0 & 1 & 0 & 2 \end{bmatrix} \begin{bmatrix} 2 & 0 & 1 & 0 & 3 \\ 1 & 1 & 0 & 1 & 0 \\ 0 & 0 & 0 & 1 & 1 \\ 0 & 1 & 0 & 0 & 0 \\ 0 & 0 & 1 & 0 & 2 \end{bmatrix} \cup$$

$$\begin{bmatrix} 1 & 0 & 1 & 1 & 2 & 3 & 4 \\ \hline 0 & 1 & 0 & 0 & 1 & 0 & 1 \\ 1 & 0 & 1 & 0 & 1 & 1 & 0 \\ 3 & 1 & 0 & 1 & 0 & 1 & 1 \\ 0 & 1 & 1 & 1 & 0 & 0 & 0 \\ 1 & 0 & 1 & 0 & 1 & 1 & 1 \\ \hline 1 & 0 & 0 & 1 & 1 & 1 & 0 \\ 0 & 1 & 0 & 0 & 1 & 0 & 1 \end{bmatrix} \begin{bmatrix} 1 & 0 & 1 & 3 & 0 & 1 & 1 & 0 \\ \hline 0 & 1 & 0 & 1 & 1 & 0 & 0 & 1 \\ 1 & 0 & 1 & 0 & 1 & 1 & 0 & 0 \\ 1 & 0 & 0 & 1 & 1 & 0 & 1 & 0 \\ 2 & 1 & 1 & 0 & 0 & 1 & 1 & 1 \\ \hline 3 & 0 & 1 & 1 & 0 & 1 & 1 & 0 \\ 4 & 1 & 0 & 1 & 0 & 1 & 0 & 1 \end{bmatrix} \cup$$

$$\begin{bmatrix} 1 & 1 & 3 & 1 & 0 & 1 & 0 \\ 2 & 0 & 0 & 1 & 1 & 0 & 1 \\ \hline 0 & 1 & 0 & 1 & 1 & 1 & 1 \\ 1 & 1 & 0 & 1 & 0 & 1 & 1 \\ 1 & 0 & 1 & 1 & 0 & 0 & 1 \\ 0 & 0 & 1 & 0 & 0 & 1 & 0 \\ \hline 3 & 1 & 0 & 3 & 1 & 1 & 0 \end{bmatrix} \begin{bmatrix} 1 & 2 & 0 & 1 & 1 & 0 & 3 \\ \hline 1 & 0 & 1 & 1 & 0 & 0 & 1 \\ 3 & 0 & 0 & 0 & 1 & 1 & 0 \\ 1 & 1 & 1 & 1 & 1 & 0 & 3 \\ 0 & 1 & 1 & 0 & 0 & 0 & 1 \\ \hline 1 & 0 & 1 & 1 & 0 & 1 & 1 \\ 0 & 1 & 1 & 1 & 1 & 0 & 0 \end{bmatrix} \cup$$

$$= \begin{bmatrix} 5 & 1 & 2 & 1 & 6 \\ 1 & 2 & 0 & 1 & 0 \\ 2 & 0 & 2 & 0 & 5 \\ 1 & 1 & 0 & 2 & 1 \\ 6 & 0 & 5 & 1 & 14 \end{bmatrix} \cup \begin{bmatrix} Z^1_{11} & Z^1_{12} & Z^1_{13} \\ Z^1_{21} & Z^1_{22} & Z^1_{23} \\ Z^1_{31} & Z^1_{32} & Z^1_{33} \end{bmatrix} \cup$$

$$\begin{bmatrix} Z^2_{11} & Z^2_{12} & Z^2_{13} \\ Z^2_{21} & Z^2_{22} & Z^2_{23} \\ Z^2_{31} & Z^2_{32} & Z^2_{33} \end{bmatrix}$$

where the calculations of $Z^1_{ij}$ and $Z^2_{ij}$ will be shown explicitly and then filled in to form the supermatrix.

$$Z^1_{11} = [1 \mid 0\ 1\ 1\ 2 \mid 3\ 4] \begin{bmatrix} 1 \\ 0 \\ 1 \\ 1 \\ 2 \\ 3 \\ 4 \end{bmatrix}$$

$$= [1][1] + [0\ 1\ 1\ 2] \begin{bmatrix} 0 \\ 1 \\ 1 \\ 2 \end{bmatrix} + [3\ 4] \begin{bmatrix} 3 \\ 4 \end{bmatrix}$$

$$= [1 + 6 + 25] = [32]$$

$$Z_{12}^1 = [1 \mid 0\ 1\ 1\ 2 \mid 3\ 4] \begin{bmatrix} 0 & 1 & 3 & 0 & 1 \\ \hline 1 & 0 & 1 & 1 & 0 \\ 0 & 1 & 0 & 1 & 1 \\ 0 & 0 & 1 & 1 & 0 \\ 1 & 1 & 0 & 0 & 1 \\ \hline 0 & 1 & 1 & 0 & 1 \\ 1 & 0 & 1 & 0 & 1 \end{bmatrix}$$

$$= [1][0\ 1\ 3\ 0\ 1] + [0\ 1\ 1\ 2] \begin{bmatrix} 1 & 0 & 1 & 1 & 0 \\ 0 & 1 & 0 & 1 & 1 \\ 0 & 0 & 1 & 1 & 0 \\ 1 & 1 & 0 & 0 & 1 \end{bmatrix} +$$

$$[3\ 4] \begin{bmatrix} 0 & 1 & 1 & 0 & 1 \\ 1 & 0 & 1 & 0 & 1 \end{bmatrix}$$

$$= [0\ 1\ 3\ 0\ 1] + [2\ 3\ 1\ 2\ 3] + [4\ 3\ 7\ 0\ 7]$$

$$= [6\ 7\ 11\ 2\ 11].$$

$$Z_{13}^1 = [1 \mid 0\ 1\ 1\ 2 \mid 3\ 4] \begin{bmatrix} 1 & 0 \\ \hline 0 & 1 \\ 0 & 0 \\ 1 & 0 \\ 1 & 1 \\ \hline 1 & 0 \\ 0 & 1 \end{bmatrix}$$

$$= [1][1\ 0] + [0\ 1\ 1\ 2] \begin{bmatrix} 0 & 1 \\ 0 & 0 \\ 1 & 0 \\ 1 & 1 \end{bmatrix} + [3\ 4] \begin{bmatrix} 1 & 0 \\ 0 & 1 \end{bmatrix}$$

$$= [1\ 0] + [3\ 2] + [3\ 4] = [7\ 6].$$

$$Z_{21}^1 = \begin{bmatrix} 0 & 1 & 0 & 0 & 1 & 0 & 1 \\ 1 & 0 & 1 & 0 & 1 & 1 & 0 \\ 3 & 1 & 0 & 1 & 0 & 1 & 1 \\ 0 & 1 & 1 & 1 & 0 & 0 & 0 \\ 1 & 0 & 1 & 0 & 1 & 1 & 1 \end{bmatrix} \begin{bmatrix} 1 \\ 0 \\ 1 \\ 1 \\ 2 \\ 3 \\ 4 \end{bmatrix}$$

$$= \begin{bmatrix} 0 \\ 1 \\ 3 \\ 0 \\ 1 \end{bmatrix}[1] + \begin{bmatrix} 1 & 0 & 0 & 1 \\ 0 & 1 & 0 & 1 \\ 1 & 0 & 1 & 0 \\ 1 & 1 & 1 & 0 \\ 0 & 1 & 0 & 1 \end{bmatrix} \begin{bmatrix} 0 \\ 1 \\ 1 \\ 2 \end{bmatrix} + \begin{bmatrix} 0 & 1 \\ 1 & 0 \\ 1 & 1 \\ 0 & 0 \\ 1 & 1 \end{bmatrix} \begin{bmatrix} 3 \\ 4 \end{bmatrix}$$

$$= \begin{bmatrix} 0 \\ 1 \\ 3 \\ 0 \\ 1 \end{bmatrix} + \begin{bmatrix} 2 \\ 3 \\ 1 \\ 2 \\ 3 \end{bmatrix} + \begin{bmatrix} 4 \\ 3 \\ 7 \\ 0 \\ 7 \end{bmatrix} = \begin{bmatrix} 6 \\ 7 \\ 11 \\ 2 \\ 11 \end{bmatrix}.$$

$$Z_{22}^1 = \begin{bmatrix} 0 & 1 & 0 & 0 & 1 & 0 & 1 \\ 1 & 0 & 1 & 0 & 1 & 1 & 0 \\ 3 & 1 & 0 & 1 & 0 & 1 & 1 \\ 0 & 1 & 1 & 1 & 0 & 0 & 0 \\ 1 & 0 & 1 & 0 & 1 & 1 & 1 \end{bmatrix} \begin{bmatrix} 0 & 1 & 3 & 0 & 1 \\ 1 & 0 & 1 & 1 & 0 \\ 0 & 1 & 0 & 1 & 1 \\ 0 & 0 & 1 & 1 & 0 \\ 1 & 1 & 0 & 0 & 1 \\ 0 & 1 & 1 & 0 & 1 \\ 1 & 0 & 1 & 0 & 1 \end{bmatrix} =$$

$$= \begin{bmatrix} 0 \\ 1 \\ 3 \\ 0 \\ 1 \end{bmatrix} [0\ 1\ 3\ 0\ 1] + \begin{bmatrix} 1 & 0 & 0 & 1 \\ 0 & 1 & 0 & 1 \\ 1 & 0 & 1 & 0 \\ 1 & 1 & 1 & 0 \\ 0 & 1 & 0 & 1 \end{bmatrix} \begin{bmatrix} 1 & 0 & 1 & 1 & 0 \\ 0 & 1 & 0 & 1 & 1 \\ 0 & 0 & 1 & 1 & 0 \\ 1 & 1 & 0 & 0 & 1 \end{bmatrix}$$

$$+ \begin{bmatrix} 0 & 1 \\ 1 & 0 \\ 1 & 1 \\ 0 & 0 \\ 1 & 1 \end{bmatrix} \begin{bmatrix} 0 & 1 & 1 & 0 & 1 \\ 1 & 0 & 1 & 0 & 1 \end{bmatrix}$$

$$= \begin{bmatrix} 0 & 0 & 0 & 0 & 0 \\ 0 & 1 & 3 & 0 & 1 \\ 0 & 3 & 9 & 0 & 3 \\ 0 & 0 & 0 & 0 & 0 \\ 0 & 1 & 3 & 0 & 1 \end{bmatrix} + \begin{bmatrix} 2 & 1 & 1 & 1 & 1 \\ 1 & 2 & 0 & 1 & 2 \\ 1 & 0 & 2 & 2 & 0 \\ 1 & 1 & 2 & 3 & 1 \\ 1 & 2 & 0 & 1 & 2 \end{bmatrix} + \begin{bmatrix} 1 & 0 & 1 & 0 & 1 \\ 0 & 1 & 1 & 0 & 1 \\ 1 & 1 & 2 & 0 & 2 \\ 0 & 0 & 0 & 0 & 0 \\ 1 & 1 & 2 & 0 & 2 \end{bmatrix}$$

$$= \begin{bmatrix} 3 & 1 & 2 & 1 & 2 \\ 1 & 4 & 4 & 1 & 4 \\ 2 & 4 & 13 & 2 & 5 \\ 1 & 1 & 2 & 3 & 1 \\ 2 & 4 & 5 & 1 & 5 \end{bmatrix}.$$

$$Z_{23}^1 = \begin{bmatrix} 0 & 1 & 0 & 0 & 1 & 0 & 1 \\ 1 & 0 & 1 & 0 & 1 & 1 & 0 \\ 3 & 1 & 0 & 1 & 0 & 1 & 1 \\ 0 & 1 & 1 & 1 & 0 & 0 & 0 \\ 1 & 0 & 1 & 0 & 1 & 1 & 1 \end{bmatrix} \begin{bmatrix} 1 & 0 \\ 0 & 1 \\ 0 & 0 \\ 1 & 0 \\ 1 & 1 \\ 1 & 0 \\ 0 & 1 \end{bmatrix}$$

$$= \begin{bmatrix} 0 \\ 1 \\ 3 \\ 0 \\ 1 \end{bmatrix} [1\ 0] + \begin{bmatrix} 1 & 0 & 0 & 1 \\ 0 & 1 & 0 & 1 \\ 1 & 0 & 1 & 0 \\ 1 & 1 & 1 & 0 \\ 0 & 1 & 0 & 1 \end{bmatrix} \begin{bmatrix} 0 & 1 \\ 0 & 0 \\ 1 & 0 \\ 1 & 1 \end{bmatrix} + \begin{bmatrix} 0 & 1 \\ 1 & 0 \\ 1 & 1 \\ 0 & 0 \\ 1 & 1 \end{bmatrix} \begin{bmatrix} 1 & 0 \\ 0 & 1 \end{bmatrix}$$

$$= \begin{bmatrix} 0 & 0 \\ 1 & 0 \\ 3 & 0 \\ 0 & 0 \\ 1 & 0 \end{bmatrix} + \begin{bmatrix} 1 & 2 \\ 1 & 1 \\ 1 & 1 \\ 1 & 1 \\ 1 & 1 \end{bmatrix} + \begin{bmatrix} 0 & 1 \\ 1 & 0 \\ 1 & 1 \\ 0 & 0 \\ 1 & 1 \end{bmatrix} = \begin{bmatrix} 1 & 3 \\ 3 & 1 \\ 5 & 2 \\ 1 & 1 \\ 3 & 2 \end{bmatrix}.$$

$$Z_{31}^{1} = \begin{bmatrix} 1 & 0 & 0 & 1 & 1 & 1 & 0 \\ 0 & 1 & 0 & 0 & 1 & 0 & 1 \end{bmatrix} \begin{bmatrix} 1 \\ 0 \\ 1 \\ 1 \\ 2 \\ 3 \\ 4 \end{bmatrix}$$

$$= \begin{bmatrix} 1 \\ 0 \end{bmatrix} [1] + \begin{bmatrix} 0 & 0 & 1 & 1 \\ 1 & 0 & 0 & 1 \end{bmatrix} \begin{bmatrix} 0 \\ 1 \\ 1 \\ 2 \end{bmatrix} + \begin{bmatrix} 1 & 0 \\ 0 & 1 \end{bmatrix} \begin{bmatrix} 3 \\ 4 \end{bmatrix}$$

$$= \begin{bmatrix} 1 \\ 0 \end{bmatrix} + \begin{bmatrix} 3 \\ 2 \end{bmatrix} + \begin{bmatrix} 3 \\ 4 \end{bmatrix} = \begin{bmatrix} 7 \\ 6 \end{bmatrix}.$$

$$Z_{32}^1 = \begin{bmatrix} 1 & 0 & 0 & 1 & 1 & 1 & 0 \\ 0 & 1 & 0 & 0 & 1 & 0 & 1 \end{bmatrix} \begin{bmatrix} 0 & 1 & 3 & 0 & 1 \\ 1 & 0 & 1 & 1 & 0 \\ 0 & 1 & 0 & 1 & 1 \\ 0 & 0 & 1 & 1 & 0 \\ 1 & 1 & 0 & 0 & 1 \\ \hline 0 & 1 & 1 & 0 & 1 \\ 1 & 0 & 1 & 0 & 1 \end{bmatrix}$$

$$= \begin{bmatrix} 1 \\ 0 \end{bmatrix} [0\ 1\ 3\ 0\ 1] + \begin{bmatrix} 0 & 0 & 1 & 1 \\ 1 & 0 & 0 & 1 \end{bmatrix} \begin{bmatrix} 1 & 0 & 1 & 1 & 0 \\ 0 & 1 & 0 & 1 & 1 \\ 0 & 0 & 1 & 1 & 0 \\ 1 & 1 & 0 & 0 & 1 \end{bmatrix}$$

$$+ \begin{bmatrix} 1 & 0 \\ 0 & 1 \end{bmatrix} \begin{bmatrix} 0 & 1 & 1 & 0 & 1 \\ 1 & 0 & 1 & 0 & 1 \end{bmatrix}$$

$$= \begin{bmatrix} 0 & 1 & 3 & 0 & 1 \\ 0 & 0 & 0 & 0 & 0 \end{bmatrix} + \begin{bmatrix} 1 & 1 & 1 & 1 & 1 \\ 2 & 1 & 1 & 1 & 1 \end{bmatrix} + \begin{bmatrix} 0 & 1 & 1 & 0 & 1 \\ 1 & 0 & 1 & 0 & 1 \end{bmatrix}$$

$$= \begin{bmatrix} 1 & 3 & 5 & 1 & 3 \\ 3 & 1 & 2 & 1 & 2 \end{bmatrix}.$$

$$Z_{33}^1 = \begin{bmatrix} 1 & 0 & 0 & 1 & 1 & 1 & 0 \\ 0 & 1 & 0 & 0 & 1 & 0 & 1 \end{bmatrix} \begin{bmatrix} 1 & 0 \\ 0 & 1 \\ 0 & 0 \\ 1 & 0 \\ 1 & 1 \\ \hline 1 & 0 \\ 0 & 1 \end{bmatrix}$$

$$= \begin{bmatrix} 1 \\ 0 \end{bmatrix} \begin{bmatrix} 1 & 0 \end{bmatrix} + \begin{bmatrix} 0 & 0 & 1 & 1 \\ 1 & 0 & 0 & 1 \end{bmatrix} \begin{bmatrix} 0 & 1 \\ 0 & 0 \\ 1 & 0 \\ 1 & 1 \end{bmatrix} + \begin{bmatrix} 1 & 0 \\ 0 & 1 \end{bmatrix} \begin{bmatrix} 1 & 0 \\ 0 & 1 \end{bmatrix}$$

$$= \begin{bmatrix} 1 & 0 \\ 0 & 0 \end{bmatrix} + \begin{bmatrix} 2 & 1 \\ 1 & 2 \end{bmatrix} + \begin{bmatrix} 1 & 0 \\ 0 & 1 \end{bmatrix}$$

$$= \begin{bmatrix} 4 & 1 \\ 1 & 3 \end{bmatrix}.$$

Thus $\begin{bmatrix} Z^1_{11} & Z^1_{12} & Z^1_{13} \\ Z^1_{21} & Z^1_{22} & Z^1_{23} \\ Z^1_{31} & Z^1_{32} & Z^1_{33} \end{bmatrix} = \left[\begin{array}{c|cccc|cc} 32 & 6 & 7 & 11 & 2 & 11 & 7 & 6 \\ \hline 6 & 3 & 1 & 2 & 1 & 2 & 1 & 3 \\ 7 & 1 & 4 & 4 & 1 & 4 & 3 & 1 \\ 11 & 2 & 4 & 13 & 2 & 5 & 5 & 2 \\ 2 & 1 & 1 & 2 & 3 & 1 & 1 & 1 \\ 11 & 2 & 4 & 5 & 1 & 5 & 3 & 2 \\ \hline 7 & 1 & 3 & 5 & 1 & 3 & 4 & 1 \\ 6 & 3 & 1 & 2 & 1 & 2 & 1 & 3 \end{array}\right].$

Now we find the values of $Z^2_{ij}$ $(1 < i, j < 3)$.

$$Z^2_{11} = \begin{bmatrix} 1 & | & 1 & 3 & 1 & 0 & | & 1 & 0 \\ 2 & | & 0 & 0 & 1 & 1 & | & 0 & 1 \end{bmatrix} \begin{bmatrix} 1 & 2 \\ 1 & 0 \\ 3 & 0 \\ 1 & 1 \\ 0 & 1 \\ \hline 1 & 0 \\ 0 & 1 \end{bmatrix}$$

$$= \begin{bmatrix} 1 \\ 2 \end{bmatrix} \begin{bmatrix} 1 & 2 \end{bmatrix} + \begin{bmatrix} 1 & 3 & 1 & 0 \\ 0 & 0 & 1 & 1 \end{bmatrix} \begin{bmatrix} 1 & 0 \\ 3 & 0 \\ 1 & 1 \\ 0 & 1 \end{bmatrix} + \begin{bmatrix} 1 & 0 \\ 0 & 1 \end{bmatrix} \begin{bmatrix} 1 & 0 \\ 0 & 1 \end{bmatrix}$$

$$= \begin{bmatrix} 1 & 2 \\ 2 & 4 \end{bmatrix} + \begin{bmatrix} 11 & 1 \\ 1 & 2 \end{bmatrix} + \begin{bmatrix} 1 & 0 \\ 0 & 1 \end{bmatrix}$$

$$= \begin{bmatrix} 13 & 3 \\ 3 & 7 \end{bmatrix}.$$

$$Z_{12}^2 = \begin{bmatrix} 1 & 1 & 3 & 1 & 0 & 1 & 0 \\ 2 & 0 & 0 & 1 & 1 & 0 & 1 \end{bmatrix} \begin{bmatrix} 0 & 1 & 1 & 0 \\ 1 & 1 & 0 & 0 \\ 0 & 0 & 1 & 1 \\ 1 & 1 & 1 & 0 \\ 1 & 0 & 0 & 0 \\ 1 & 1 & 0 & 1 \\ 1 & 1 & 1 & 0 \end{bmatrix}$$

$$= \begin{bmatrix} 1 \\ 2 \end{bmatrix} \begin{bmatrix} 0 & 1 & 1 & 0 \end{bmatrix} + \begin{bmatrix} 1 & 3 & 1 & 0 \\ 0 & 0 & 1 & 1 \end{bmatrix} \begin{bmatrix} 1 & 1 & 0 & 0 \\ 0 & 0 & 1 & 1 \\ 1 & 1 & 1 & 0 \\ 1 & 0 & 0 & 0 \end{bmatrix} +$$

$$\begin{bmatrix} 1 & 0 \\ 0 & 1 \end{bmatrix} \begin{bmatrix} 1 & 1 & 0 & 1 \\ 1 & 1 & 1 & 0 \end{bmatrix}$$

$$= \begin{bmatrix} 0 & 1 & 1 & 0 \\ 0 & 2 & 2 & 0 \end{bmatrix} + \begin{bmatrix} 2 & 2 & 4 & 3 \\ 2 & 1 & 1 & 0 \end{bmatrix} + \begin{bmatrix} 1 & 1 & 0 & 1 \\ 1 & 1 & 1 & 0 \end{bmatrix}$$

$$= \begin{bmatrix} 3 & 4 & 5 & 4 \\ 3 & 4 & 4 & 0 \end{bmatrix}.$$

$$Z_{13}^2 = \begin{bmatrix} 1 & 1 & 3 & 1 & 0 & 1 & 0 \\ 2 & 0 & 0 & 1 & 1 & 0 & 1 \end{bmatrix} \begin{bmatrix} 3 \\ 1 \\ 0 \\ 3 \\ 1 \\ 1 \\ 0 \end{bmatrix}$$

$$= \begin{bmatrix} 1 \\ 2 \end{bmatrix}[3] + \begin{bmatrix} 1 & 3 & 1 & 0 \\ 0 & 0 & 1 & 1 \end{bmatrix} \begin{bmatrix} 1 \\ 0 \\ 3 \\ 1 \end{bmatrix} + \begin{bmatrix} 1 & 0 \\ 0 & 1 \end{bmatrix} \begin{bmatrix} 1 \\ 0 \end{bmatrix}$$

$$= \begin{bmatrix} 3 \\ 6 \end{bmatrix} + \begin{bmatrix} 4 \\ 4 \end{bmatrix} + \begin{bmatrix} 1 \\ 0 \end{bmatrix} = \begin{bmatrix} 8 \\ 10 \end{bmatrix}.$$

$$Z_{21}^2 = \begin{bmatrix} 0 & 1 & 0 & 1 & 1 & 1 & 1 \\ 1 & 1 & 0 & 1 & 0 & 1 & 1 \\ 1 & 0 & 1 & 1 & 0 & 0 & 1 \\ 0 & 0 & 1 & 0 & 0 & 1 & 0 \end{bmatrix} \begin{bmatrix} 1 & 2 \\ 1 & 0 \\ 3 & 0 \\ 1 & 1 \\ 0 & 1 \\ 1 & 0 \\ 0 & 1 \end{bmatrix}$$

$$= \begin{bmatrix} 0 \\ 1 \\ 1 \\ 0 \end{bmatrix} [1 \quad 2] + \begin{bmatrix} 1 & 0 & 1 & 1 \\ 1 & 0 & 1 & 0 \\ 0 & 1 & 1 & 0 \\ 0 & 1 & 0 & 0 \end{bmatrix} \begin{bmatrix} 1 & 0 \\ 3 & 0 \\ 1 & 1 \\ 0 & 1 \end{bmatrix}$$

$$+ \begin{bmatrix} 1 & 1 \\ 1 & 1 \\ 0 & 1 \\ 1 & 0 \end{bmatrix} \begin{bmatrix} 1 & 0 \\ 0 & 1 \end{bmatrix} =$$

$$= \begin{bmatrix} 0 & 0 \\ 1 & 2 \\ 1 & 2 \\ 0 & 0 \end{bmatrix} + \begin{bmatrix} 2 & 2 \\ 2 & 1 \\ 4 & 1 \\ 3 & 0 \end{bmatrix} + \begin{bmatrix} 1 & 1 \\ 1 & 1 \\ 0 & 1 \\ 1 & 0 \end{bmatrix}$$

$$= \begin{bmatrix} 3 & 3 \\ 4 & 4 \\ 5 & 4 \\ 4 & 0 \end{bmatrix}.$$

$$Z_{22}^2 = \begin{bmatrix} 0 & 1 & 0 & 1 & 1 & 1 & 1 \\ 1 & 1 & 0 & 1 & 0 & 1 & 1 \\ 1 & 0 & 1 & 1 & 0 & 0 & 1 \\ 0 & 0 & 1 & 0 & 0 & 1 & 0 \end{bmatrix} \begin{bmatrix} 0 & 1 & 1 & 0 \\ 1 & 1 & 0 & 0 \\ 0 & 0 & 1 & 1 \\ 1 & 1 & 1 & 0 \\ 1 & 0 & 0 & 0 \\ 1 & 1 & 0 & 1 \\ 1 & 1 & 1 & 0 \end{bmatrix}$$

$$= \begin{bmatrix} 0 \\ 1 \\ 1 \\ 0 \end{bmatrix} \begin{bmatrix} 0 & 1 & 1 & 0 \end{bmatrix} + \begin{bmatrix} 1 & 0 & 1 & 1 \\ 1 & 0 & 1 & 0 \\ 0 & 1 & 1 & 0 \\ 0 & 1 & 0 & 0 \end{bmatrix} \begin{bmatrix} 1 & 1 & 0 & 0 \\ 0 & 0 & 1 & 1 \\ 1 & 1 & 1 & 0 \\ 1 & 0 & 0 & 0 \end{bmatrix}$$

$$+ \begin{bmatrix} 1 & 1 \\ 1 & 1 \\ 0 & 1 \\ 1 & 0 \end{bmatrix} \begin{bmatrix} 1 & 1 & 0 & 1 \\ 1 & 1 & 1 & 0 \end{bmatrix}$$

$$= \begin{bmatrix} 0 & 0 & 0 & 0 \\ 0 & 1 & 1 & 0 \\ 0 & 1 & 1 & 0 \\ 0 & 0 & 0 & 0 \end{bmatrix} + \begin{bmatrix} 3 & 2 & 1 & 0 \\ 2 & 2 & 1 & 0 \\ 1 & 1 & 2 & 1 \\ 0 & 0 & 1 & 1 \end{bmatrix} + \begin{bmatrix} 2 & 2 & 1 & 1 \\ 2 & 2 & 1 & 1 \\ 1 & 1 & 1 & 0 \\ 1 & 1 & 0 & 1 \end{bmatrix}$$

$$= \begin{bmatrix} 5 & 4 & 2 & 1 \\ 4 & 5 & 3 & 1 \\ 2 & 3 & 4 & 1 \\ 1 & 1 & 1 & 2 \end{bmatrix}.$$

$$Z_{23}^2 = \begin{bmatrix} 0 & 1 & 0 & 1 & 1 & 1 & 1 \\ 1 & 1 & 0 & 1 & 0 & 1 & 1 \\ 1 & 0 & 1 & 1 & 0 & 0 & 1 \\ 0 & 0 & 1 & 0 & 0 & 1 & 0 \end{bmatrix} \begin{bmatrix} 3 \\ 1 \\ 0 \\ 3 \\ 1 \\ 1 \\ 0 \end{bmatrix}$$

$$= \begin{bmatrix} 0 \\ 1 \\ 1 \\ 0 \end{bmatrix} [3] + \begin{bmatrix} 1 & 0 & 1 & 1 \\ 1 & 0 & 1 & 0 \\ 0 & 1 & 1 & 0 \\ 0 & 1 & 0 & 0 \end{bmatrix} \begin{bmatrix} 1 \\ 0 \\ 3 \\ 1 \end{bmatrix} + \begin{bmatrix} 1 & 1 \\ 1 & 1 \\ 0 & 1 \\ 1 & 0 \end{bmatrix} \begin{bmatrix} 1 \\ 0 \end{bmatrix}$$

$$= \begin{bmatrix} 0 \\ 3 \\ 3 \\ 0 \end{bmatrix} + \begin{bmatrix} 5 \\ 4 \\ 3 \\ 0 \end{bmatrix} + \begin{bmatrix} 1 \\ 1 \\ 0 \\ 1 \end{bmatrix} = \begin{bmatrix} 6 \\ 8 \\ 6 \\ 1 \end{bmatrix}.$$

$$Z_{31}^2 = \begin{bmatrix} 3 \mid 1 & 0 & 3 & 1 \mid 1 & 0 \end{bmatrix} \begin{bmatrix} 1 & 2 \\ \hline 1 & 0 \\ 3 & 0 \\ 1 & 1 \\ 0 & 1 \\ \hline 1 & 0 \\ 0 & 1 \end{bmatrix}$$

$$[3] + \begin{bmatrix} 1 & 2 \end{bmatrix} + \begin{bmatrix} 10 & 31 \end{bmatrix} \begin{bmatrix} 1 & 0 \\ 3 & 0 \\ 1 & 1 \\ 0 & 1 \end{bmatrix} + \begin{bmatrix} 1 & 0 \end{bmatrix} \begin{bmatrix} 1 & 0 \\ 0 & 1 \end{bmatrix}$$

$$= \begin{bmatrix} 3 & 6 \end{bmatrix} + \begin{bmatrix} 4 & 4 \end{bmatrix} + \begin{bmatrix} 1 & 0 \end{bmatrix} = \begin{bmatrix} 8 & 10 \end{bmatrix}.$$

$$Z_{32}^2 = \begin{bmatrix} 3 \mid 1 & 0 & 3 & 1 \mid 1 & 0 \end{bmatrix} \begin{bmatrix} 0 & 1 & 1 & 0 \\ \hline 1 & 1 & 0 & 0 \\ 0 & 0 & 1 & 1 \\ 1 & 1 & 1 & 0 \\ 1 & 0 & 0 & 0 \\ \hline 1 & 1 & 0 & 1 \\ 1 & 1 & 1 & 0 \end{bmatrix}$$

$$= [3]\begin{bmatrix} 0 & 1 & 1 & 0 \end{bmatrix} + \begin{bmatrix} 1 & 0 & 3 & 1 \end{bmatrix} \begin{bmatrix} 1 & 1 & 0 & 0 \\ 0 & 0 & 1 & 1 \\ 1 & 1 & 1 & 0 \\ 1 & 0 & 0 & 0 \end{bmatrix} +$$

$$\begin{bmatrix} 1 & 0 \end{bmatrix} \begin{bmatrix} 1 & 1 & 0 & 1 \\ 1 & 1 & 1 & 0 \end{bmatrix}$$

$$= \begin{bmatrix} 0 & 3 & 3 & 0 \end{bmatrix} + \begin{bmatrix} 5 & 4 & 3 & 0 \end{bmatrix} + \begin{bmatrix} 1 & 1 & 0 & 1 \end{bmatrix}$$

$$= \begin{bmatrix} 6 & 8 & 6 & 1 \end{bmatrix}.$$

$$Z_{33}^2 = \begin{bmatrix} 3 \mid 1 & 0 & 3 & 1 \mid 1 & 0 \end{bmatrix} \begin{bmatrix} 3 \\ \hline 1 \\ 0 \\ 3 \\ 1 \\ \hline 1 \\ 0 \end{bmatrix}$$

$$= [3][3] + \begin{bmatrix} 1 & 0 & 3 & 1 \end{bmatrix} \begin{bmatrix} 1 \\ 0 \\ 3 \\ 1 \end{bmatrix} + \begin{bmatrix} 1 & 0 \end{bmatrix} \begin{bmatrix} 1 \\ 0 \end{bmatrix}$$

$$= [9] + [11] + [1] = [21].$$

Now

$$\begin{bmatrix} Z_{11}^2 & Z_{12}^2 & Z_{13}^2 \\ Z_{21}^2 & Z_{22}^2 & Z_{23}^2 \\ Z_{31}^2 & Z_{32}^2 & Z_{33}^2 \end{bmatrix}$$

$$= \begin{bmatrix} 13 & 3 & 3 & 4 & 5 & 4 & 8 \\ 3 & 7 & 3 & 4 & 4 & 0 & 10 \\ \hline 3 & 3 & 5 & 4 & 2 & 1 & 6 \\ 4 & 4 & 4 & 5 & 3 & 1 & 8 \\ 5 & 4 & 2 & 3 & 4 & 1 & 6 \\ 4 & 0 & 1 & 1 & 1 & 2 & 1 \\ \hline 8 & 10 & 6 & 8 & 6 & 1 & 21 \end{bmatrix}.$$

Thus

$$TT^T = \begin{bmatrix} 5 & 1 & 2 & 1 & 6 \\ 1 & 2 & 0 & 1 & 0 \\ 2 & 0 & 2 & 0 & 5 \\ 1 & 1 & 0 & 2 & 1 \\ 6 & 0 & 5 & 1 & 14 \end{bmatrix} \cup$$

$$\begin{bmatrix} 32 & 6 & 7 & 11 & 2 & 11 & 7 & 6 \\ \hline 6 & 3 & 1 & 2 & 1 & 2 & 1 & 3 \\ 7 & 1 & 4 & 4 & 1 & 4 & 3 & 1 \\ 11 & 2 & 4 & 13 & 2 & 5 & 5 & 2 \\ 2 & 1 & 1 & 2 & 3 & 1 & 1 & 1 \\ 11 & 2 & 4 & 5 & 1 & 5 & 3 & 2 \\ \hline 7 & 1 & 3 & 5 & 1 & 3 & 4 & 1 \\ 6 & 3 & 1 & 2 & 1 & 2 & 1 & 3 \end{bmatrix}$$

$$\cup \begin{bmatrix} 13 & 3 & 3 & 4 & 5 & 4 & 8 \\ 3 & 7 & 3 & 4 & 4 & 0 & 10 \\ \hline 3 & 3 & 5 & 4 & 2 & 1 & 6 \\ 4 & 4 & 4 & 5 & 3 & 1 & 8 \\ 5 & 4 & 2 & 3 & 4 & 1 & 6 \\ 4 & 0 & 1 & 1 & 1 & 2 & 1 \\ \hline 8 & 10 & 6 & 8 & 6 & 1 & 21 \end{bmatrix}.$$

We see the product gives us a symmetric semi super trimatrix. The interested reader can find $T^TT$ which will give yet another symmetric semi super trimatrix.

We have given the explicit working of the product, the main motivation for it is that while calculating for a general case we have to give lots of notations and we felt it would only confuse the reader. It is certainly easy to find the major products using these illustrations for after all what we are interested is that the reader should be in a positions to work it out. The theory behind it is not very difficult but the notational representations is little cumbersome.

Chapter Four

# SUPER n-MATRICES AND THEIR PROPERTIES

In this chapter we for the first time introduce the notion of super n-matrices (n-an integer and n > 3) and give a few of its properties. When n = 1 we get the supermatrix, n = 2, it is the superbimatrix studied in chapter two. The study of super tri matrix has been studied in chapter three, the case when n = 3; when we say super n-matrix we mean n a positive integer and n > 3. Here we also define the notion of semi super n-matrices and show how the product defined using them at times yields only an ordinary or elementary n-matrix. Further some of the products induce a symmetric super n-matrix or a quasi symmetric super n-matrix.

**DEFINITION 4.1:** *Let $V = V_1 \cup V_2 \cup \ldots \cup V_n$ (n > 3) denote n distinct super matrices, i.e., each $V_i$ is a super matrix $1 \leq i \leq n$. '$\cup$' is just a symbol. We define V to be a super n-matrix.*

*Example 4.1:* Let $V = V_1 \cup V_2 \cup V_3 \cup V_4$ where

$$V_1 = [1\ 0\ |\ 2\ 3\ 4\ |\ 5],$$

$$V_2 = \begin{bmatrix} 2 \\ 0 \\ 1 \\ \hline 1 \\ 1 \\ \hline 3 \\ 2 \\ 5 \end{bmatrix},$$

$$V_3 = \begin{bmatrix} 3 & 0 \\ 1 & 2 \end{bmatrix}$$

and

$$V_4 = \begin{bmatrix} 3 & 0 & 1 \\ 1 & 1 & 1 \\ \hline 2 & 0 & 2 \\ 5 & 3 & 5 \end{bmatrix}.$$

V is a super 4-matrix. Here n = 4.

***Example 4.2:*** Let $T = T_1 \cup T_2 \cup T_3 \cup T_4 \cup T_5$ where

$$T_1 = \begin{bmatrix} 3 & 0 & 2 & 1 & 5 & 3 & 1 \\ 1 & 1 & 2 & 0 & 7 & 8 & 0 \end{bmatrix},$$

$$T_2 = \begin{bmatrix} 3 & 1 & 3 & 4 \\ 0 & 1 & 1 & 0 \\ 1 & 0 & 0 & 1 \\ 3 & 1 & 5 & 2 \\ \hline 1 & 2 & 3 & 4 \\ 9 & 8 & 3 & 1 \end{bmatrix},$$

$$T_3 = \begin{bmatrix} 3 & 1 & | & 2 & 4 & 5 \\ 0 & 1 & | & 0 & 5 & 1 \\ \hline 6 & 3 & | & 7 & 8 & 9 \\ 1 & 1 & | & 0 & 1 & 2 \\ 2 & 0 & | & 2 & 0 & 3 \end{bmatrix},$$

$$T_4 = \begin{bmatrix} 8 & 5 & | & 3 & 1 & 2 \\ \hline 1 & 0 & | & 1 & 1 & 2 \\ 8 & 4 & | & 0 & 3 & 1 \\ 7 & 1 & | & 0 & 1 & 1 \\ 8 & 3 & | & 8 & 1 & 5 \end{bmatrix}$$

and

$$T_5 = \begin{bmatrix} 3 & 4 & | & 1 & 0 \\ 1 & 3 & | & 0 & 1 \\ \hline 5 & 1 & | & 1 & 0 \\ 0 & 2 & | & 3 & 1 \end{bmatrix}.$$

T is a super 5-matrix. In this case we have n = 5.

***Example 4.3:*** Let $T = T_1 \cup T_2 \cup T_3 \cup T_4 \cup T_5 \cup T_6$ where

$$T_1 = [0 \mid 1\ 0\ 1 \mid 2\ 3\ 4],$$

$$T_2 = [1 \mid 2\ 3\ 4\ 5 \mid 5\ 7\ 8 \mid 9\ 0],$$

$$T_3 = [0\ 1\ 3 \mid 4\ 5 \mid 7\ 8\ 9\ 10],$$

$$T_4 = [6 \mid 1\ 2 \mid 3\ 0\ 1\ 4\ 6\ 1],$$

$$T_5 = [3\ 1\ 0 \mid 2\ 2\ 5\ 0\ 1]$$

and

$$T_6 = [1 \mid 2\ 3 \mid 4\ 5\ 6 \mid 7\ 1\ 8\ 1].$$

T is a super 6-matrix. We see each of the super matrices $T_i$'s are row super vectors.

**DEFINITION 4.2:** *Let $T = T_1 \cup T_2 \cup ... \cup T_n$ ($n > 3$) be a super n-matrix. If each of the $T_i$ is a super row vector $i = 1, 2, ..., n$ then we call T to be a row super n-vector or super row n-vector.*

The example 4.3 is a super row 6-vector.

*Example 4.4:* Let $S = S_1 \cup S_2 \cup S_3 \cup S_4$ where

$$S_1 = \begin{bmatrix} \frac{1}{2} \\ \frac{3}{4} \\ 5 \\ 6 \end{bmatrix}, S_2 = \begin{bmatrix} 0 \\ \frac{1}{2} \\ 4 \\ 7 \\ \frac{0}{2} \\ \frac{5}{1} \end{bmatrix}, S_3 = \begin{bmatrix} 1 \\ 2 \\ 5 \\ 3 \\ \frac{1}{7} \\ 6 \\ 4 \end{bmatrix} \text{ and } S_4 = \begin{bmatrix} \frac{1}{2} \\ \frac{-1}{3} \\ 2 \\ \frac{-3}{1} \\ 0 \\ 1 \end{bmatrix}.$$

S is a super 4-matrix or a super 4-column vector.

**DEFINITION 4.3:** *Let $T = T_1 \cup T_2 \cup T_3 \cup ... \cup T_n$ be a super n-matrix ($n > 3$). If each of the $T_i$ is a super column vector ($1 \leq i \leq n$) then we call T to be a super n-column vector or column super n-vector or column super n-matrix.*

Example 4.4 is a super 4-column vector.

*Example 4.5:* Let $V = V_1 \cup V_2 \cup V_3 \cup V_4 \cup V_5$ where

$$V_1 = \begin{bmatrix} 2 & | & 1 \\ 0 & | & 1 \end{bmatrix}, V_2 = \begin{bmatrix} 1 & 1 \\ \hline 0 & 2 \end{bmatrix},$$

$$V_3 = \begin{bmatrix} 1 & 2 & 3 \\ 4 & 5 & 6 \\ 7 & 8 & 9 \end{bmatrix}$$

$$V_4 = \begin{bmatrix} 1 & 2 & 3 & 4 \\ 0 & 1 & 2 & 5 \\ 7 & 8 & 1 & 0 \\ 9 & 6 & 4 & 2 \end{bmatrix}$$

and

$$V_5 = \begin{bmatrix} 1 & 2 & 0 & 1 & 1 \\ 0 & 1 & 2 & 0 & 1 \\ 1 & 4 & 0 & 1 & 3 \\ 2 & 5 & 1 & 2 & 1 \\ 3 & 6 & 1 & 0 & 0 \end{bmatrix}.$$

V is a super 5-matrix. We see each of the super matrices $V_i$; $1 \leq i \leq 5$ are square super matrices.

**DEFINITION 4.4:** *Let $K = K_1 \cup K_2 \cup K_3 \cup ... \cup K_n$ ($n > 3$) be a super n-matrix if each of the $K_i$ are $t \times t$ square matrices for $i = 1, 2, ..., n$ then we call K to be a $t \times t$ square super n-matrix. If on the other hand each of $K_i$ is a $m_i \times m_i$ square matrix $i = 1, 2, ..., n$ then we call K to be a mixed square super n-matrix.*

The example 4.5 is a mixed super square n-matrix ($n = 5$).

*Example 4.6:* Let $P = P_1 \cup P_2 \cup P_3 \cup P_4$ where

$$P_1 = \begin{bmatrix} 1 & 2 & 3 & 4 \\ 5 & 6 & 7 & 8 \\ 9 & 0 & 1 & 2 \\ 3 & 4 & 5 & 6 \end{bmatrix},$$

$$P_2 = \left[\begin{array}{cc|cc} 0 & 1 & 2 & 3 \\ 1 & 0 & 1 & 1 \\ 0 & 1 & 0 & 1 \\ \hline 1 & 0 & 0 & 1 \end{array}\right],$$

$$P_3 = \left[\begin{array}{c|ccc} 1 & 2 & 3 & 4 \\ \hline 0 & 1 & 0 & 1 \\ 1 & 1 & 1 & 1 \\ 1 & 1 & 0 & 0 \end{array}\right]$$

and

$$P_4 = \left[\begin{array}{ccc|c} 1 & 0 & 2 & 4 \\ 1 & 1 & 0 & 1 \\ \hline 1 & 0 & 1 & 1 \\ 0 & 1 & 1 & 1 \end{array}\right].$$

P is a super 4-matrix. Infact P is a $4 \times 4$ square super 4-matrix.

*Example 4.7:* Let $S = S_1 \cup S_2 \cup S_3 \cup S_4 \cup S_5$ where

$$S_1 = \left[\begin{array}{c|ccc|cc} 3 & 1 & 0 & 1 & 0 & 9 \\ 0 & 1 & 1 & 1 & 7 & 1 \\ \hline 1 & 1 & 1 & 0 & 8 & 8 \end{array}\right],$$

$$S_2 = \left[\begin{array}{cc|cc} 3 & 4 & 1 & 2 \\ \hline 1 & 1 & 1 & 0 \\ 0 & 0 & 0 & 1 \\ 1 & 1 & 0 & 1 \\ \hline 0 & 1 & 0 & 9 \\ 1 & 1 & 2 & 5 \end{array}\right],$$

$$S_3 = \begin{bmatrix} 2 & 1 & 0 & 4 & 7 & 5 \\ 3 & 1 & 1 & 9 & 6 & 6 \\ 1 & 1 & 1 & 8 & 4 & 7 \end{bmatrix},$$

$$S_4 = \begin{bmatrix} 1 & 7 & 9 & 3 \\ 2 & 8 & 8 & 1 \\ 3 & 9 & 7 & 2 \\ \hline 4 & 0 & 6 & 0 \\ 5 & 1 & 5 & 2 \\ 6 & 2 & 4 & 5 \end{bmatrix}$$

and

$$S_5 = \begin{bmatrix} 1 & 1 & 7 & 4 & 5 & 1 & 1 & 4 \\ 2 & 1 & 8 & 3 & 6 & 2 & 1 & 0 \\ 3 & 1 & 9 & 2 & 7 & 3 & 1 & 8 \\ 4 & 1 & 6 & 1 & 9 & 0 & 1 & 7 \end{bmatrix}.$$

S is a super 5-matrix. We see each of the $S_i$; $1 \le i \le 5$ are rectangular super matrices of different order.

**DEFINITION 4.5:** *Let $S = S_1 \cup S_2 \cup S_3 \cup … \cup S_n$ (n >3) be a super n-matrix. If each of the $S_i$'s are rectangular super matrices of different order then we call S to be a mixed super rectangular n-matrix or a mixed rectangular super n-matrix.*

The super n-matrix given in example 4.7 is a mixed rectangular n-matrix. If in the super n-matrix $S = S_1 \cup S_2 \cup … \cup S_n$ each of the $S_i$'s are m × t (m ≠ t) rectangular super matrices then we call S to be an m × t rectangular n-matrix.

The super n-matrix given in example 4.7 is a mixed rectangular n-matrix.

**DEFINITION 4.6:** *Let $T = T_1 \cup T_2 \cup T_3 \cup … \cup T_n$ (n > 3) be a super n-matrix. If some of the $T_i$'s are square supermatrices and* *some of the $T_j$'s are rectangular super matrices ($i = j$); $i \leq j$, $i \leq n$, then we call T to be a mixed super n-matrix.*

***Example 4.8:*** Let $K = K_1 \cup K_2 \cup K_3 \cup K_4 \cup K_5$ be a super 5-matrix where

$$K_1 = \left[\begin{array}{c|cc} 3 & 1 & 2 \\ 1 & 1 & 0 \\ \hline 0 & 1 & 1 \end{array}\right],$$

$$K_2 = \left[\begin{array}{c|cccc} 3 & 0 & 1 & 1 & 1 \\ \hline 0 & 1 & 2 & 0 & 1 \\ 1 & 1 & 0 & 1 & 2 \end{array}\right],$$

$$K_3 = \left[\begin{array}{ccc|c} 3 & 1 & 0 & 1 \\ \hline 0 & 2 & 1 & 0 \\ 1 & 1 & 4 & 3 \\ 5 & 6 & 1 & 2 \end{array}\right],$$

$$K_4 = \left[\begin{array}{ccc|ccc|c} 1 & 6 & 5 & 0 & 1 & 0 & 3 \\ 2 & 7 & 4 & 9 & 1 & 1 & 4 \\ \hline 3 & 8 & 3 & 8 & 0 & 6 & 2 \\ 4 & 9 & 2 & 7 & 2 & 0 & 1 \\ 5 & 0 & 1 & 6 & 3 & 1 & 5 \end{array}\right]$$

and

$$K_5 = \left[\begin{array}{cc|ccc|cc|c} 1 & 1 & 3 & 5 & 7 & 9 & 1 & 7 \\ 0 & 2 & 4 & 6 & 8 & 0 & 4 & 8 \end{array}\right].$$

K is a mixed super n-matrix here n= 5.

**DEFINITION 4.7:** *Let $T = T_1 \cup T_2 \cup ... \cup T_n$ ($n > 3$) be a super n-matrix we say T is a special super row n-vector if each of the $T_i$'s is a $m_i \times n_i$ matrix ($n_i > m_i$); $1 < i < n$; with partitions done* *only vertically i.e., only between the columns. No partition is made between the rows. If $P = P_1 \cup P_2 \cup ... \cup P_n$ ($n > 3$) be a super n-matrix we say P is a special super column n vector if each $P_i$ is a $t_i \times s_i$ matrix with $t_i > s_i$ ($1 \le i \le n$) and each column matrix is partitioned only in between the rows i.e., horizontally and never partitioned in between the columns.*

**Example 4.9:** Let $S = S_1 \cup S_2 \cup S_3 \cup S_4 \cup S_5$ where

$$S_1 = \begin{bmatrix} 1 & 3 & 1 & 0 & 3 \\ 2 & 4 & 1 & 1 & 3 \end{bmatrix},$$

$$S_2 = \begin{bmatrix} 3 & 1 & 4 & 7 & 1 & 0 & 1 \\ 1 & 0 & 5 & 8 & 1 & 3 & 1 \\ 2 & 1 & 6 & 9 & 2 & 7 & 4 \end{bmatrix},$$

$$S_3 = \begin{bmatrix} 1 & 1 & 2 & 3 & 5 & 2 & 0 & 1 & 4 \\ 0 & 1 & 1 & 0 & 2 & 1 & 1 & 7 & 2 \\ 2 & 4 & 3 & 0 & 1 & 0 & 1 & 1 & 0 \\ 3 & 0 & 1 & 1 & 1 & 1 & 1 & 0 & 1 \end{bmatrix},$$

$$S_4 = \begin{bmatrix} 1 & 0 & 1 & 2 & 3 & 1 & 1 & 2 & 3 & 1 \\ 1 & 1 & 1 & 0 & 1 & 1 & 0 & 1 & 2 & 3 \\ 2 & 1 & 1 & 5 & 1 & 2 & 0 & 1 & 1 & 0 \end{bmatrix},$$

and

$$S_5 = \begin{bmatrix} 1 & 1 & 7 & 3 & 2 & 7 & 1 & 5 & 1 & 1 & 6 \\ 2 & 1 & 8 & 2 & 3 & 8 & 2 & 6 & 2 & 0 & 2 \\ 1 & 1 & 6 & 5 & 4 & 9 & 3 & 7 & 3 & 7 & 3 \\ 0 & 0 & 4 & 1 & 6 & 0 & 4 & 8 & 6 & 2 & 1 \end{bmatrix}.$$

S is a super 5-matrix which is a special super row 5-matrix.

**Example 4.10:** Let $T = T_1 \cup T_2 \cup T_3 \cup T_4$ where

$$T_1 = \begin{bmatrix} 3 & 1 & 0 & 1 \\ 1 & 1 & 2 & 3 \\ 5 & 6 & 7 & 8 \\ 9 & 10 & 1 & 1 \\ 1 & 2 & 1 & 1 \\ 0 & 1 & 0 & 1 \\ 1 & 1 & 1 & 0 \\ 0 & 1 & 1 & 0 \end{bmatrix},$$

$$T_2 = \begin{bmatrix} 2 & 1 & 0 \\ 3 & 1 & 1 \\ 9 & 8 & 7 \\ 6 & 5 & 4 \\ 3 & 2 & 1 \\ 0 & 1 & 1 \\ 1 & 0 & 0 \\ 0 & 1 & 1 \end{bmatrix}, \quad T_3 = \begin{bmatrix} 1 & 2 \\ 3 & 1 \\ 1 & 5 \\ 5 & 8 \\ 7 & 6 \\ 6 & 8 \\ 9 & 1 \\ 1 & 4 \\ 4 & 3 \\ 3 & 6 \end{bmatrix}$$

and

$$T_4 = \begin{bmatrix} 1 & 3 & 1 \\ 1 & 1 & 1 \\ 0 & 1 & 2 \\ 3 & 1 & 1 \\ 1 & 1 & 1 \\ 1 & 1 & 4 \\ 4 & 2 & 3 \\ 3 & 0 & 1 \\ 1 & 1 & 5 \\ 0 & 1 & 2 \end{bmatrix}.$$

T is a super 4-matrix which is a special super column 4-matrix.

*Example 4.11:* Let $S = S_1 \cup S_2 \cup S_3 \cup S_4$ where

$$S_1 = \begin{bmatrix} 2 & 1 & 4 & 1 \\ \hline 0 & 3 & 1 & 5 \\ 1 & 2 & 3 & 4 \\ 5 & 6 & 7 & 8 \\ \hline 9 & 0 & 1 & 8 \\ \hline 1 & 1 & 1 & 4 \\ 4 & 4 & 1 & 2 \end{bmatrix}, \quad S_2 = \begin{bmatrix} 3 & 1 & 1 \\ 0 & 2 & 4 \\ \hline 1 & 3 & 1 \\ 1 & 4 & 6 \\ 7 & 8 & 9 \\ \hline 1 & 0 & 4 \\ 4 & 1 & 2 \\ \hline 1 & 1 & 6 \end{bmatrix},$$

$$S_3 = \begin{bmatrix} 1 & 2 & 3 & 4 & 5 \\ 6 & 7 & 8 & 9 & 0 \\ \hline 3 & 2 & 1 & 4 & 8 \\ 1 & 1 & 1 & 4 & 1 \\ \hline 7 & 0 & 8 & 1 & 3 \\ 3 & 1 & 2 & 5 & 6 \\ 1 & 1 & 0 & 1 & 1 \end{bmatrix}$$

and

$$S_4 = \begin{bmatrix} 1 & 2 & 3 & 1 \\ \hline 1 & 1 & 0 & 1 \\ 1 & 0 & 1 & 1 \\ 2 & 1 & 1 & 1 \\ \hline 4 & 2 & 3 & 1 \\ 1 & 1 & 0 & 1 \\ 1 & 5 & 0 & 3 \\ 1 & 7 & 2 & 3 \end{bmatrix}.$$

Clearly S is a super 4 matrix but S is not a special super column 4 matrix for $S_1$ and $S_3$ are partitioned vertically also. $S_1$ and $S_3$ are only super matrices and are not special super column vectors.

***Example 4.12:*** Let $V = V_1 \cup V_2 \cup V_3 \cup V_4$ be a super 4-matrix; where

$$V_1 = \left[\begin{array}{cc|cccc|c} 2 & 1 & 0 & 1 & 7 & 1 & 8 \\ 0 & 1 & 1 & 4 & 6 & 1 & 7 \\ \hline 1 & 1 & 3 & 5 & 8 & 0 & 6 \\ 0 & 2 & 4 & 6 & 1 & 9 & 4 \end{array}\right],$$

$$V_2 = \left[\begin{array}{cc|ccc} 1 & 3 & 5 & 7 & 9 \\ 2 & 4 & 6 & 8 & 0 \end{array}\right]$$

$$V_3 = \left[\begin{array}{cc|ccc|cc} 1 & 9 & 3 & 0 & 6 & 1 & 8 \\ 2 & 8 & 2 & 1 & 5 & 2 & 9 \\ 3 & 7 & 1 & 2 & 4 & 5 & 1 \end{array}\right]$$

and

$$V_4 = \left[\begin{array}{c|ccc|ccc} 1 & 3 & 1 & 0 & 1 & 1 & 1 \\ \hline 2 & 1 & 5 & 1 & 0 & 1 & 1 \end{array}\right].$$

V is not a special super row vector for the super matrix $V_1$ and $V_4$ are partitioned horizontally also.

Now we can define minor product of special super n-matrices. We first illustrate it and then define the concept.

***Example 4.13:*** Let $T = T_1 \cup T_2 \cup T_3 \cup T_4$ be a special super column 4-matrix and $V = V_1 \cup V_2 \cup V_3 \cup V_4$ be a special super row 4 matrix the minor product TV is defined as follows.

$$T = T_1 \cup T_2 \cup T_3 \cup T_4$$

$$= \begin{bmatrix} 1 & 2 \\ 0 & 1 \\ \hline 1 & 1 \\ 2 & 3 \\ 3 & 4 \\ 1 & 0 \end{bmatrix} \cup \begin{bmatrix} 1 & 2 & 3 \\ 0 & 1 & 1 \\ \hline 1 & 1 & 1 \\ 1 & 1 & 2 \\ \hline 0 & 1 & 1 \\ \hline 1 & 1 & 1 \end{bmatrix} \cup \begin{bmatrix} 0 & 1 & 2 & 3 \\ 1 & 1 & 1 & 1 \\ 1 & 2 & 3 & 4 \\ 1 & 1 & 1 & 0 \\ 0 & 1 & 1 & 1 \\ 1 & 1 & 1 & 4 \\ 0 & 0 & 1 & 0 \end{bmatrix}$$

$$\cup \begin{bmatrix} 1 & 2 & 0 & 1 & 1 \\ 0 & 1 & 1 & 0 & 1 \\ \hline 1 & 2 & 3 & 0 & 1 \\ 1 & 0 & 1 & 1 & 0 \\ 0 & 1 & 0 & 0 & 1 \\ 1 & 1 & 0 & 1 & 0 \\ 0 & 1 & 0 & 0 & 1 \\ 1 & 1 & 1 & 0 & 1 \end{bmatrix}$$

be the given special super column 4 vector. $S = S_1 \cup S_2 \cup S_3 \cup S_4$

$$= \begin{bmatrix} 0 & 2 & | & 1 & 0 & 1 & 1 \\ 1 & 3 & | & 1 & 1 & 1 & 3 \end{bmatrix} \cup$$

$$\begin{bmatrix} 0 & 1 & 1 & 1 & | & 1 & 0 & 1 & | & 1 \\ 1 & 0 & 1 & 0 & | & 1 & 1 & 0 & | & 1 \\ 1 & 0 & 1 & 1 & | & 0 & 1 & 1 & | & 0 \end{bmatrix} \cup$$

$$\begin{bmatrix} 1 & 2 & 0 & | & 1 & 1 & 1 & 0 & 1 & | & 1 & 0 \\ 0 & 1 & 1 & | & 1 & 1 & 0 & 1 & 0 & | & 0 & 1 \\ 0 & 1 & 0 & | & 1 & 1 & 0 & 1 & 1 & | & 0 & 0 \\ 1 & 1 & 1 & | & 0 & 0 & 1 & 0 & 0 & | & 0 & 0 \end{bmatrix} \cup$$

$$\begin{bmatrix} 1 & 1 & 0 & 0 & 0 & 0 & 0 & 1 & 1 \\ 1 & 0 & 1 & 0 & 0 & 0 & 1 & 0 & 1 \\ 0 & 0 & 0 & 1 & 0 & 0 & 1 & 1 & 0 \\ 1 & 0 & 0 & 0 & 1 & 0 & 0 & 0 & 1 \\ 1 & 0 & 0 & 0 & 0 & 1 & 1 & 1 & 1 \end{bmatrix}$$

be the given special super row 4-vector.

$$\begin{aligned}
TS &= (T_1 \cup T_2 \cup T_3 \cup T_4)(S_1 \cup S_2 \cup S_3 \cup S_4) \\
&= T_1 S_1 \cup T_2 S_2 \cup T_3 S_3 \cup T_4 S_4
\end{aligned}$$

$$= \begin{bmatrix} 1 & 2 \\ 0 & 1 \\ \hline 1 & 1 \\ 2 & 3 \\ 3 & 4 \\ 1 & 0 \end{bmatrix} \begin{bmatrix} 0 & 2 & 1 & 0 & 1 & 1 \\ 1 & 3 & 1 & 1 & 1 & 3 \end{bmatrix} \cup$$

$$\begin{bmatrix} 1 & 2 & 3 \\ 0 & 1 & 1 \\ \hline 1 & 1 & 1 \\ \hline 1 & 1 & 2 \\ 0 & 1 & 1 \\ \hline 1 & 1 & 1 \end{bmatrix} \begin{bmatrix} 0 & 1 & 1 & 1 & 1 & 0 & 1 & 1 \\ 1 & 0 & 1 & 0 & 1 & 1 & 0 & 1 \\ 1 & 0 & 1 & 1 & 0 & 1 & 1 & 0 \end{bmatrix} \cup$$

$$\begin{bmatrix} 0 & 1 & 2 & 3 \\ \hline 1 & 1 & 1 & 1 \\ 1 & 2 & 3 & 4 \\ \hline 1 & 1 & 1 & 0 \\ 0 & 1 & 1 & 1 \\ 1 & 1 & 1 & 4 \\ 0 & 0 & 1 & 0 \end{bmatrix} \begin{bmatrix} 1 & 2 & 0 & 1 & 1 & 1 & 0 & 1 & 1 & 0 \\ 0 & 1 & 1 & 1 & 1 & 0 & 1 & 0 & 0 & 1 \\ 0 & 1 & 0 & 1 & 1 & 0 & 1 & 1 & 0 & 0 \\ 1 & 1 & 1 & 0 & 0 & 1 & 0 & 0 & 0 & 0 \end{bmatrix} \cup$$

$$\begin{bmatrix} 1 & 2 & 0 & 1 & 1 \\ 0 & 1 & 1 & 0 & 1 \\ \hline 1 & 2 & 3 & 0 & 1 \\ \hline 1 & 0 & 1 & 1 & 0 \\ 0 & 1 & 0 & 0 & 1 \\ 1 & 1 & 0 & 1 & 0 \\ 0 & 1 & 0 & 0 & 1 \\ 1 & 1 & 1 & 0 & 1 \end{bmatrix} \begin{bmatrix} 1 & 1 & 0 & 0 & 0 & 0 & 0 & 1 & 1 \\ 1 & 0 & 1 & 0 & 0 & 0 & 1 & 0 & 1 \\ 0 & 0 & 0 & 1 & 0 & 0 & 1 & 1 & 0 \\ 1 & 0 & 0 & 0 & 1 & 0 & 0 & 0 & 1 \\ 1 & 0 & 0 & 0 & 0 & 1 & 1 & 1 & 1 \end{bmatrix}$$

$$= \begin{bmatrix} \begin{bmatrix} 1 & 2 \\ 0 & 1 \end{bmatrix} \begin{bmatrix} 0 & 2 \\ 1 & 3 \end{bmatrix} & \begin{bmatrix} 1 & 2 \\ 0 & 1 \end{bmatrix} \begin{bmatrix} 1 & 0 & 1 & 1 \\ 1 & 1 & 1 & 3 \end{bmatrix} \\ \hline \begin{bmatrix} 1 & 1 \\ 2 & 3 \\ 3 & 4 \\ 1 & 0 \end{bmatrix} \begin{bmatrix} 0 & 2 \\ 1 & 3 \end{bmatrix} & \begin{bmatrix} 1 & 1 \\ 2 & 3 \\ 3 & 4 \\ 1 & 0 \end{bmatrix} \begin{bmatrix} 1 & 0 & 1 & 1 \\ 1 & 1 & 1 & 3 \end{bmatrix} \end{bmatrix} \cup$$

$$\begin{bmatrix} \begin{bmatrix} 1 & 2 & 3 \\ 0 & 1 & 1 \\ 1 & 1 & 1 \end{bmatrix} \begin{bmatrix} 0 & 1 & 1 & 1 \\ 1 & 0 & 1 & 0 \\ 1 & 0 & 1 & 1 \end{bmatrix} & \begin{bmatrix} 1 & 2 & 3 \\ 0 & 1 & 1 \\ 1 & 1 & 1 \end{bmatrix} \begin{bmatrix} 1 & 0 & 1 \\ 1 & 1 & 0 \\ 0 & 1 & 1 \end{bmatrix} & \begin{bmatrix} 1 & 2 & 3 \\ 0 & 1 & 1 \\ 1 & 1 & 1 \end{bmatrix} \begin{bmatrix} 1 \\ 1 \\ 0 \end{bmatrix} \\ \hline \begin{bmatrix} 1 & 1 & 2 \\ 0 & 1 & 1 \end{bmatrix} \begin{bmatrix} 0 & 1 & 1 & 1 \\ 1 & 0 & 1 & 0 \\ 1 & 0 & 1 & 1 \end{bmatrix} & \begin{bmatrix} 1 & 1 & 2 \\ 0 & 1 & 1 \end{bmatrix} \begin{bmatrix} 1 & 0 & 1 \\ 1 & 1 & 0 \\ 0 & 1 & 1 \end{bmatrix} & \begin{bmatrix} 1 & 1 & 2 \\ 0 & 1 & 1 \end{bmatrix} \begin{bmatrix} 1 \\ 1 \\ 0 \end{bmatrix} \\ \hline \begin{bmatrix} 1 & 1 & 1 \end{bmatrix} \begin{bmatrix} 0 & 1 & 1 & 1 \\ 1 & 0 & 1 & 0 \\ 1 & 0 & 1 & 1 \end{bmatrix} & \begin{bmatrix} 1 & 1 & 1 \end{bmatrix} \begin{bmatrix} 1 & 0 & 1 \\ 1 & 1 & 0 \\ 0 & 1 & 1 \end{bmatrix} & \begin{bmatrix} 1 & 1 & 1 \end{bmatrix} \begin{bmatrix} 1 \\ 1 \\ 0 \end{bmatrix} \end{bmatrix}$$

$$\subset \left[ \begin{array}{c|c} (0\ 1\ 2\ 3)\begin{bmatrix} 1 & 2 & 0 \\ 0 & 1 & 1 \\ 0 & 1 & 0 \\ 1 & 1 & 1 \end{bmatrix} & (0\ 1\ 2\ 3)\begin{bmatrix} 1 & 1 & 1 & 0 & 1 \\ 1 & 1 & 0 & 1 & 0 \\ 1 & 1 & 0 & 1 & 1 \\ 0 & 0 & 1 & 0 & 0 \end{bmatrix} \\ \hline \begin{pmatrix} 1 & 1 & 1 & 1 \\ 1 & 2 & 3 & 4 \end{pmatrix}\begin{pmatrix} 1 & 2 & 0 \\ 0 & 1 & 1 \\ 0 & 1 & 0 \\ 1 & 1 & 1 \end{pmatrix} & \begin{pmatrix} 1 & 1 & 1 & 1 \\ 1 & 2 & 3 & 4 \end{pmatrix}\begin{bmatrix} 1 & 1 & 1 & 0 & 1 \\ 1 & 1 & 0 & 1 & 0 \\ 1 & 1 & 0 & 1 & 1 \\ 0 & 0 & 1 & 0 & 0 \end{bmatrix} \\ \hline \begin{pmatrix} 1 & 1 & 1 & 0 \\ 0 & 1 & 1 & 1 \\ 1 & 1 & 1 & 4 \\ 0 & 0 & 1 & 0 \end{pmatrix}\begin{pmatrix} 1 & 2 & 0 \\ 0 & 1 & 1 \\ 0 & 1 & 0 \\ 1 & 1 & 1 \end{pmatrix} & \begin{pmatrix} 1 & 1 & 1 & 0 \\ 0 & 1 & 1 & 1 \\ 1 & 1 & 1 & 4 \\ 0 & 0 & 1 & 0 \end{pmatrix}\begin{bmatrix} 1 & 1 & 1 & 0 & 1 \\ 1 & 1 & 0 & 1 & 0 \\ 1 & 1 & 0 & 1 & 1 \\ 0 & 0 & 1 & 0 & 0 \end{bmatrix} \end{array} \right]$$

$$\left[ \begin{array}{c} (0\ 1\ 2\ 3)\begin{bmatrix} 1 & 0 \\ 0 & 1 \\ 0 & 0 \\ 0 & 0 \end{bmatrix} \\ \hline \begin{pmatrix} 1 & 1 & 1 & 1 \\ 1 & 2 & 3 & 4 \end{pmatrix}\begin{bmatrix} 1 & 0 \\ 0 & 1 \\ 0 & 0 \\ 0 & 0 \end{bmatrix} \\ \hline \begin{pmatrix} 1 & 1 & 1 & 0 \\ 0 & 1 & 1 & 1 \\ 1 & 1 & 1 & 4 \\ 0 & 0 & 1 & 0 \end{pmatrix}\begin{bmatrix} 1 & 0 \\ 0 & 1 \\ 0 & 0 \\ 0 & 0 \end{bmatrix} \end{array} \right] \subset$$

$$\begin{bmatrix} \begin{pmatrix} 1 & 2 & 0 & 1 & 1 \\ 0 & 1 & 1 & 0 & 1 \end{pmatrix} \begin{bmatrix} 1 & 1 \\ 1 & 0 \\ 0 & 0 \\ 1 & 0 \\ 1 & 0 \end{bmatrix} & \begin{pmatrix} 1 & 2 & 0 & 1 & 1 \\ 0 & 1 & 1 & 0 & 1 \end{pmatrix} \begin{bmatrix} 0 & 0 & 0 & 0 \\ 1 & 0 & 0 & 0 \\ 0 & 1 & 0 & 0 \\ 0 & 0 & 1 & 0 \\ 0 & 0 & 0 & 1 \end{bmatrix} \\ \hline \begin{bmatrix} 1 & 2 & 3 & 0 & 1 \end{bmatrix} \begin{bmatrix} 1 & 1 \\ 1 & 0 \\ 0 & 0 \\ 1 & 0 \\ 1 & 0 \end{bmatrix} & \begin{bmatrix} 1 & 2 & 3 & 0 & 1 \end{bmatrix} \begin{bmatrix} 0 & 0 & 0 & 0 \\ 1 & 0 & 0 & 0 \\ 0 & 1 & 0 & 0 \\ 0 & 0 & 1 & 0 \\ 0 & 0 & 0 & 1 \end{bmatrix} \\ \hline \begin{pmatrix} 1 & 0 & 1 & 1 & 0 \\ 0 & 1 & 0 & 0 & 1 \\ 1 & 1 & 0 & 1 & 0 \\ 0 & 1 & 0 & 0 & 1 \\ 1 & 1 & 1 & 0 & 1 \end{pmatrix} \begin{bmatrix} 1 & 1 \\ 1 & 0 \\ 0 & 0 \\ 1 & 0 \\ 1 & 0 \end{bmatrix} & \begin{pmatrix} 1 & 0 & 1 & 1 & 0 \\ 0 & 1 & 0 & 0 & 1 \\ 1 & 1 & 0 & 1 & 0 \\ 0 & 1 & 0 & 0 & 1 \\ 1 & 1 & 1 & 0 & 1 \end{pmatrix} \begin{bmatrix} 0 & 0 & 0 & 0 \\ 1 & 0 & 0 & 0 \\ 0 & 1 & 0 & 0 \\ 0 & 0 & 1 & 0 \\ 0 & 0 & 0 & 1 \end{bmatrix} \end{bmatrix}$$

$$\begin{bmatrix} \begin{pmatrix} 1 & 2 & 0 & 1 & 1 \\ 0 & 1 & 1 & 0 & 1 \end{pmatrix} \begin{bmatrix} 0 & 1 & 1 \\ 1 & 0 & 1 \\ 1 & 1 & 0 \\ 0 & 0 & 1 \\ 1 & 1 & 1 \end{bmatrix} \\ \hline \begin{bmatrix} 1 & 2 & 3 & 0 & 1 \end{bmatrix} \begin{bmatrix} 0 & 1 & 1 \\ 1 & 0 & 1 \\ 1 & 1 & 0 \\ 0 & 0 & 1 \\ 1 & 1 & 1 \end{bmatrix} \\ \hline \begin{pmatrix} 1 & 0 & 1 & 1 & 0 \\ 0 & 1 & 0 & 0 & 1 \\ 1 & 1 & 0 & 1 & 0 \\ 0 & 1 & 0 & 0 & 1 \\ 1 & 1 & 1 & 0 & 1 \end{pmatrix} \begin{bmatrix} 0 & 1 & 1 \\ 1 & 0 & 1 \\ 1 & 1 & 0 \\ 0 & 0 & 1 \\ 1 & 1 & 1 \end{bmatrix} \end{bmatrix} =$$

$$\left[\begin{array}{cc|cccc} 2 & 8 & 3 & 2 & 3 & 7 \\ 1 & 3 & 1 & 1 & 1 & 3 \\ \hline 1 & 5 & 2 & 1 & 2 & 4 \\ 3 & 13 & 5 & 3 & 5 & 11 \\ 4 & 18 & 7 & 4 & 7 & 15 \\ 0 & 2 & 1 & 0 & 1 & 1 \end{array}\right] \cup$$

$$\left[\begin{array}{cccc|ccc|c} 5 & 1 & 6 & 4 & 3 & 5 & 4 & 3 \\ 2 & 0 & 2 & 1 & 1 & 2 & 1 & 1 \\ 2 & 1 & 3 & 2 & 2 & 2 & 2 & 2 \\ \hline 3 & 1 & 4 & 3 & 2 & 3 & 3 & 2 \\ 2 & 0 & 2 & 1 & 1 & 2 & 1 & 1 \\ 2 & 1 & 3 & 2 & 2 & 2 & 2 & 2 \end{array}\right] \cup$$

$$\left[\begin{array}{ccc|ccccc|cc} 3 & 6 & 4 & 3 & 3 & 3 & 3 & 2 & 0 & 1 \\ \hline 2 & 5 & 2 & 3 & 3 & 2 & 2 & 2 & 1 & 1 \\ 5 & 11 & 6 & 6 & 6 & 5 & 5 & 4 & 1 & 2 \\ \hline 1 & 4 & 1 & 3 & 3 & 1 & 2 & 2 & 1 & 1 \\ 1 & 3 & 2 & 2 & 2 & 1 & 2 & 1 & 0 & 1 \\ 5 & 8 & 5 & 3 & 3 & 5 & 2 & 2 & 1 & 1 \\ 0 & 1 & 0 & 1 & 1 & 0 & 1 & 1 & 0 & 0 \end{array}\right] \cup$$

$$\left[\begin{array}{cc|ccc|ccc} 5 & 1 & 2 & 0 & 1 & 1 & 3 & 2 & 5 \\ 2 & 0 & 1 & 1 & 0 & 1 & 3 & 2 & 2 \\ \hline 4 & 1 & 2 & 3 & 0 & 1 & 6 & 5 & 4 \\ \hline 2 & 1 & 0 & 1 & 1 & 0 & 1 & 2 & 2 \\ 2 & 0 & 1 & 0 & 0 & 1 & 2 & 1 & 2 \\ 3 & 1 & 1 & 0 & 1 & 0 & 1 & 1 & 3 \\ 2 & 0 & 1 & 0 & 0 & 1 & 2 & 1 & 2 \\ 3 & 1 & 1 & 1 & 0 & 1 & 3 & 3 & 3 \end{array}\right].$$

The resultant is a super 4-matrix which is not a special super row or column n-matrix.

Now we define symmetric super n-matrix and a semi symmetric super n-matrix.

**DEFINITION 4.8:** *Let $T = T_1 \cup T_2 \cup \ldots \cup T_n$ ($n > 3$) be a super n-matrix. If each of the $T_i$ is a symmetric super n-matrix then we call T to be a symmetric super n-matrix: $1 \leq i \leq n$,*

*Example 4.14:* Let $T = T_1 \cup T_2 \cup T_3 \cup T_4 \cup T_5$ be a super 5 matrix where

$$T_1 = \left[\begin{array}{c|c} 3 & 10 \\ \hline 10 & 1 \end{array}\right], \quad T_2 = \left[\begin{array}{c|cc} 3 & 1 & 1 \\ \hline 1 & 0 & 1 \\ 1 & 1 & 8 \end{array}\right],$$

$$T_3 = \left[\begin{array}{c|ccc} 1 & 2 & 0 & 4 \\ \hline 2 & 1 & 5 & 2 \\ 0 & 5 & 1 & 6 \\ 4 & 2 & 6 & 4 \end{array}\right],$$

$$T_4 = \left[\begin{array}{cc|c} 1 & 2 & 3 \\ 2 & 5 & 7 \\ \hline 3 & 7 & 1 \end{array}\right]$$

and

$$T_5 = \left[\begin{array}{cc|ccc|c} 1 & 2 & 3 & 4 & 5 & 6 \\ 2 & 0 & 1 & 1 & 0 & 1 \\ \hline 3 & 1 & 2 & 7 & 1 & 2 \\ 4 & 1 & 7 & 0 & 3 & 5 \\ 5 & 0 & 1 & 3 & 1 & 2 \\ \hline 6 & 1 & 2 & 5 & 2 & 0 \end{array}\right].$$

T is a symmetric super 5-matrix.

Thus a symmetric super n-matrix is either a square super n-matrix or a mixed square super n-matrix. It may so happen that in a super n-matrix $T = T_1 \cup T_2 \cup T_3 \cup \ldots \cup T_n$ ($n > 3$) some of the $T_i$'s are symmetric supermatrices some of them just supermatrices in such case we define T to be a quasi symmetric super n-matrix

*Example 4.15:* Let $T = T_1 \cup T_2 \cup T_3 \cup T_4 \cup T_5 \cup T_6$ be a super 6-matrix where

$$T_1 = \begin{bmatrix} 3 & 1 & 2 & 0 \\ 1 & 1 & 0 & 1 \\ \hline 2 & 0 & 5 & 7 \\ 0 & 1 & 7 & 0 \end{bmatrix},$$

$$T_2 = \begin{bmatrix} 7 & 8 & 1 & 0 & 1 \\ 8 & 1 & 5 & 1 & 3 \\ 1 & 5 & 1 & 0 & 1 \\ 0 & 1 & 0 & 2 & 5 \\ \hline 1 & 3 & 1 & 5 & 0 \end{bmatrix}, T_3 = \begin{bmatrix} 3 & 4 \\ \hline 5 & 7 \end{bmatrix}$$

$$T_4 = \begin{bmatrix} 1 & 2 & 1 & 0 \\ 2 & 1 & 1 & 2 \\ 1 & 1 & 3 & 0 \\ 0 & 2 & 0 & 3 \end{bmatrix},$$

$$T_5 = \begin{bmatrix} 1 & 1 & 1 & 0 & 1 & 0 \\ 1 & 2 & 0 & 1 & 0 & 0 \\ \hline 1 & 0 & 5 & 3 & 1 & 2 \\ 0 & 1 & 3 & 0 & 1 & 7 \\ 1 & 0 & 1 & 1 & 0 & 1 \\ 0 & 0 & 2 & 7 & 1 & 0 \end{bmatrix}$$

and

$$T_6 = \begin{bmatrix} 1 & 0 & | & 1 & 3 & 5 \\ 0 & 1 & | & 2 & 0 & 1 \\ 1 & 2 & | & 7 & 2 & 5 \\ \hline 3 & 0 & | & 2 & 1 & 3 \\ 5 & 1 & | & 5 & 3 & 0 \end{bmatrix}.$$

Clearly $T = T_1 \cup T_2 \cup T_3 \cup T_4 \cup T_5 \cup T_6$ is only a quasi symmetric super 6-matrix.

**DEFINITION 4.9:** *Let $T = T_1 \cup T_2 \cup T_3 \cup \ldots \cup T_n$ ($n > 3$) be an n-matrix in which some of the $T_i$'s are super matrices and some of the $T_j$'s are just matrices $1 \le i, j \le n$. Then we call $T = T_1 \cup T_2 \cup \ldots \cup T_n$ to be a semi super n-matrix.*

*Example 4.16:* Let $T = T_1 \cup T_2 \cup T_3 \cup T_4$ where

$$T_1 = \begin{bmatrix} 3 & 1 & | & 2 & 5 & 6 \\ 1 & 3 & | & 0 & 1 & 1 \\ \hline 2 & 0 & | & 2 & 1 & 0 \\ 5 & 1 & | & 3 & 2 & 1 \\ 6 & 1 & | & 4 & 5 & 6 \end{bmatrix},$$

$$T_2 = \begin{bmatrix} 0 & 1 & 2 & 3 & | & 4 & 5 \\ 6 & 7 & 8 & 9 & | & 0 & 1 \\ \hline 3 & 0 & 1 & 0 & | & 5 & 7 \end{bmatrix},$$

$$T_3 = \begin{bmatrix} 2 & 1 \\ 0 & 3 \\ 1 & 2 \end{bmatrix}$$

and

$$T_4 = \begin{bmatrix} 3 & 1 & 2 & 1 & 5 \\ 1 & 0 & 1 & 3 & 2 \\ \hline 3 & 1 & 3 & 5 & 7 \\ \hline 0 & 2 & 4 & 0 & 6 \end{bmatrix}.$$

T is a semi super 4-matrix

***Example 4.17:*** Let $S = S_1 \cup S_2 \cup S_3 \cup S_4 \cup S_5$ where

$$S_1 = \begin{bmatrix} 3 & 1 & 5 & 1 & 8 \\ 1 & 1 & 6 & 8 & 0 \\ \hline 0 & 2 & 7 & 9 & 2 \\ 1 & 3 & 0 & 1 & 6 \end{bmatrix},$$

$$S_2 = \begin{bmatrix} 3 & 0 & 1 \\ 8 & 5 & 0 \\ 0 & 1 & 2 \end{bmatrix},$$

$$S_3 = \begin{bmatrix} 3 & 6 & 0 \\ 5 & 4 & 3 \\ 2 & 1 & 8 \\ 1 & 0 & 1 \\ 8 & 6 & 7 \end{bmatrix},$$

$$S_4 = \begin{bmatrix} 3 & 6 & 4 & 8 & 1 & 5 \\ 6 & 0 & 6 & 1 & 0 & 3 \\ \hline 4 & 1 & 2 & 4 & 1 & 0 \\ 9 & 0 & 1 & 2 & 4 & 0 \\ 2 & 2 & 6 & 1 & 1 & 2 \end{bmatrix}$$

and

$$S_5 = \left[\begin{array}{cc|cccc} 3 & 1 & 0 & 3 & 1 & 1 \\ 1 & 4 & 6 & 1 & 0 & 2 \\ \hline 1 & 1 & 2 & 2 & 1 & 3 \\ 0 & 1 & 5 & 6 & 8 & 2 \end{array}\right]$$

$S = S_1 \cup S_2 \cup S_3 \cup S_4 \cup S_5$ is a semi super 5-matrix.

***Example 4.18:*** Let $S = S_1 \cup S_2 \cup S_3 \cup S_4$ where

$$S_1 = \left[\begin{array}{cc|c} 3 & 0 & 1 \\ 4 & 5 & 0 \\ \hline 0 & 2 & 6 \end{array}\right],$$

$$S_2 = \left[\begin{array}{cccc} 3 & 6 & 1 & 2 \\ 0 & 1 & 0 & 6 \\ \hline 4 & 1 & 1 & 2 \\ 3 & 0 & 1 & 6 \\ \hline 2 & 1 & 5 & 4 \\ 1 & 0 & 2 & 8 \end{array}\right],$$

$$S_3 = \left[\begin{array}{ccc} 3 & 1 & 2 \\ 1 & 1 & 1 \\ 0 & 1 & 2 \\ 6 & 7 & 0 \\ 0 & 5 & 1 \\ 1 & 2 & 3 \\ 0 & 1 & 5 \end{array}\right]$$

and

$$S_4 = \begin{bmatrix} 3 & 2 & 6 & 4 \\ 2 & 1 & 0 & 1 \\ 6 & 2 & 1 & 1 \\ 4 & 0 & 1 & 2 \end{bmatrix}.$$

$S = S_1 \cup S_2 \cup S_3 \cup S_4$ is a semi super 4-matrix.

**DEFINITION 4.10:** *Let $T = T_1 \cup T_2 \cup ... \cup T_n$ be a semi super n-matrix ($n > 3$). If each $T_i$ is a $m \times p$ matrix with $m < p$ and $1 \leq i \leq n$ with some of the $T_i$'s super matrices and some of the $T_j$'s just matrices then we call T to be a special semi super row n-vector. If $m > p$ then we call $T = T_1 \cup T_2 \cup ... \cup T_n$ to be a special semi super column n-vector.*

We illustrate them by the following examples

***Example 4.19:*** Let $V = V_1 \cup V_2 \cup V_3 \cup V_4 \cup V_5$ be a special semi super row 5-vector where

$$V_1 = \begin{bmatrix} 3 & 1 & 0 & 1 & 3 & 2 & 1 & 1 \\ 0 & 1 & 1 & 0 & 2 & 0 & 2 & 1 \\ 2 & 5 & 1 & 1 & 1 & 1 & 5 & 1 \end{bmatrix},$$

$$V_2 = \begin{bmatrix} 2 & 4 & 6 & 8 & 1 & 0 \\ 1 & 5 & 7 & 9 & 2 & 1 \end{bmatrix},$$

$$V_3 = \begin{bmatrix} 1 & 0 & 0 & 1 & 6 & 1 \\ 2 & 1 & 3 & 4 & 4 & 8 \\ 3 & 2 & 2 & 0 & 5 & 7 \\ 4 & 3 & 1 & 3 & 2 & 6 \end{bmatrix},$$

$$V_4 = \begin{bmatrix} 1 & 0 & 5 & 1 \\ 1 & 1 & 6 & 2 \\ 0 & 1 & 7 & 3 \end{bmatrix}$$

and

$$V_5 = \begin{bmatrix} 3 & 1 & 2 & 3 & 5 & 6 & 7 & 8 \\ 0 & 1 & 9 & 8 & 7 & 6 & 4 & 3 \end{bmatrix}.$$

*Example 4.20:* Let $S = S_1 \cup S_2 \cup S_3 \cup S_4 \cup S_5$ where

$$S_1 = \begin{bmatrix} 3 & 4 & 2 & 1 & 0 & 1 & 6 \\ 1 & 5 & 3 & 1 & 8 & 9 & 4 \end{bmatrix},$$

$$S_2 = \begin{bmatrix} 3 & 0 & 1 & 8 & 1 & 9 & 8 & 1 \\ 1 & 8 & 2 & 6 & 0 & 6 & 5 & 2 \\ \hline 2 & 7 & 0 & 4 & 3 & 7 & 4 & 1 \\ 6 & 6 & 1 & 3 & 8 & 5 & 3 & 1 \end{bmatrix},$$

$$S_3 = \begin{bmatrix} 1 & 4 & 7 & 0 & 9 \\ 2 & 5 & 8 & 1 & 0 \\ 3 & 6 & 9 & 4 & 1 \end{bmatrix},$$

$$S_4 = \begin{bmatrix} 1 & 2 & 3 & 4 & 5 & 6 & 7 & 8 & 9 \\ 0 & 3 & 6 & 7 & 8 & 1 & 0 & 1 & 7 \\ 6 & 1 & 2 & 0 & 1 & 8 & 1 & 0 & 1 \end{bmatrix}$$

and

$$S_5 = \begin{bmatrix} 1 & 0 & 1 \\ 2 & 6 & 2 \end{bmatrix}.$$

$S = S_1 \cup S_2 \cup S_3 \cup S_4 \cup S_5$ is a semi super 5-matrix.

Infact S is a not special semi super row 5-vector as $S_2$ is a super matrix which is divided horizontally also. The row vectors must be partitioned only vertically and never horizontally.

*Example 4.21:* Let T = T$_1$ ∪ T$_2$ ∪ T$_3$ ∪ T$_4$ ∪ T$_5$ where

$$T_1 = \begin{bmatrix} 3 & 1 & 2 \\ 1 & 0 & 1 \\ \hline 1 & 2 & 3 \\ 0 & 1 & 4 \\ 5 & 6 & 7 \\ 8 & 9 & 0 \\ \hline 1 & 2 & 7 \end{bmatrix},$$

$$T_2 = \begin{bmatrix} 1 & 0 & 8 & 1 \\ 3 & 1 & 6 & 8 \\ \hline 1 & 0 & 7 & 6 \\ 6 & 1 & 0 & 2 \\ 1 & 1 & 0 & 1 \\ \hline 3 & 1 & 4 & 6 \\ 7 & 8 & 9 & 0 \\ 1 & 1 & 3 & 4 \end{bmatrix},$$

$$T_3 = \begin{bmatrix} 1 & 1 & 1 & 1 & 1 \\ \hline 0 & 0 & 1 & 0 & 0 \\ \hline 0 & 0 & 0 & 1 & 0 \\ \hline 0 & 0 & 0 & 0 & 1 \\ 1 & 0 & 0 & 0 & 0 \\ 0 & 1 & 0 & 0 & 0 \\ \hline 1 & 0 & 1 & 0 & 0 \\ 0 & 0 & 1 & 0 & 1 \\ 1 & 0 & 1 & 0 & 1 \end{bmatrix}$$

and

$$T_4 = \begin{bmatrix} 1 & 6 & 1 & 4 \\ 1 & 0 & 3 & 4 \\ 4 & 1 & 0 & 1 \\ \hline 1 & 1 & 1 & 1 \\ 0 & 1 & 0 & 1 \\ \hline 3 & 1 & 1 & 1 \\ 1 & 1 & 0 & 1 \\ 2 & 0 & 7 & 1 \\ 1 & 1 & 5 & 3 \end{bmatrix}$$

be a semi super n-matrix (n = 4) but T is not a special semi super column n-matrix as $T_2$ is divided vertically also. If T is to be special semi super column n-vector each matrix must be rectangular m × t matrix with m > t with only horizontal partition of the m × t matrix. If the partition is also vertical even for a single $T_i$ then T is not a special semi super column n-vector.

Now we illustrate major and minor products of super n-matrices before we give the abstract definition.

***Example 4.22:*** Let $T = T_1 \cup T_2 \cup T_3 \cup T_4 \cup T_5$ and $V = V_1 \cup V_2 \cup V_3 \cup V_4 \cup V_5$ be two special super n matrices (n = 5) to find the product TV. Given $T = T_1 \cup T_2 \cup T_3 \cup T_4 \cup T_5$ where

$$T_1 = \begin{bmatrix} 1 & 2 \\ \hline 3 & 0 \\ 1 & 2 \\ 1 & 4 \end{bmatrix},$$

$$T_2 = \begin{bmatrix} 3 & 0 & 1 \\ 1 & 1 & 1 \\ \hline 2 & 3 & 1 \\ 1 & 0 & 1 \\ 1 & 5 & 2 \end{bmatrix}$$

$$T_3 = \begin{bmatrix} 1 & 2 & 3 & 4 \\ 0 & 1 & 1 & 1 \\ 1 & 1 & 0 & 1 \\ \hline 0 & 1 & 1 & 0 \\ \hline 1 & 2 & 0 & 1 \\ 0 & 1 & 1 & 0 \end{bmatrix}$$

$$T_4 = \begin{bmatrix} 1 & 2 & 1 \\ 1 & 1 & 1 \\ 0 & 1 & 0 \\ \hline 1 & 0 & 1 \\ 2 & 1 & 1 \\ 1 & 0 & 1 \\ 1 & 1 & 1 \\ 2 & 0 & 2 \end{bmatrix}$$

and

$$T_5 = \begin{bmatrix} 2 & 1 & 0 \\ 1 & 1 & 1 \\ 0 & 1 & 1 \\ \hline 0 & 1 & 2 \\ 6 & 0 & 0 \\ 1 & 0 & 1 \\ 1 & 1 & 1 \\ 0 & 0 & 1 \\ 0 & 1 & 0 \end{bmatrix}$$

be the special semi column 5-matrix. Given $V = V_1 \cup V_2 \cup V_3 \cup V_4 \cup V_5$ where

$$V_1 = \begin{bmatrix} 1 & 3 & 1 & 0 & 1 & 1 \\ 0 & 1 & 1 & 1 & 0 & 0 \end{bmatrix}$$

$$V_2 = \begin{bmatrix} 1 & 1 & 3 & 1 & 1 & 1 & 1 \\ 3 & 0 & 0 & 0 & 1 & 0 & 1 \\ 0 & 1 & 0 & 2 & 1 & 0 & 1 \end{bmatrix}$$

$$V_3 = \begin{bmatrix} 1 & 0 & 4 & 2 & 0 & 1 & 1 \\ 0 & 1 & 0 & 1 & 0 & 1 & 1 \\ 1 & 0 & 1 & 1 & 1 & 0 & 0 \\ 5 & 1 & 1 & 0 & 0 & 1 & 0 \end{bmatrix}$$

$$V_4 = \begin{bmatrix} 1 & 1 & 0 & 1 & 1 & 0 & 1 & 1 & 0 & 1 & 1 \\ 0 & 1 & 1 & 0 & 0 & 1 & 0 & 1 & 6 & 1 & 0 \\ 0 & 1 & 1 & 1 & 1 & 1 & 0 & 1 & 0 & 0 & 0 \end{bmatrix}$$

and

$$V_5 = \begin{bmatrix} 1 & 1 & 1 & 2 & 2 & 0 & 1 & 1 & 1 \\ 0 & 6 & 1 & 1 & 1 & 1 & 0 & 1 & 0 \\ 1 & 0 & 0 & 0 & 0 & 0 & 1 & 0 & 1 \end{bmatrix}$$

be the special semi row 5-matrix.

$$TV = (T_1 \cup T_2 \cup T_3 \cup T_4 \cup T_5)(V_1 \cup V_2 \cup V_3 \cup V_4 \cup V_5)$$
$$= T_1T_1 \cup T_2V_2 \cup T_3 V_3 \cup V_4T_4 \cup T_5V_5$$

$$= \begin{bmatrix} 1 & 2 \\ \hline 3 & 0 \\ 1 & 2 \\ 1 & 4 \end{bmatrix} \begin{bmatrix} 1 & 3 & 1 & 0 & 1 & 1 \\ 0 & 1 & 1 & 1 & 0 & 0 \end{bmatrix} \cup$$

$$\begin{bmatrix} 3 & 0 & 1 \\ 1 & 1 & 1 \\ \hline 2 & 3 & 1 \\ 1 & 0 & 1 \\ 1 & 5 & 2 \end{bmatrix} \begin{bmatrix} 1 & 1 & 3 & 1 & 1 & 1 & 1 \\ 3 & 0 & 0 & 0 & 1 & 0 & 1 \\ 0 & 1 & 0 & 2 & 1 & 0 & 1 \end{bmatrix} \cup$$

$$\begin{bmatrix} 1 & 2 & 3 & 4 \\ 0 & 1 & 1 & 1 \\ 1 & 1 & 0 & 1 \\ 0 & 1 & 1 & 0 \\ \hline 1 & 2 & 0 & 1 \\ 0 & 1 & 1 & 0 \end{bmatrix} \begin{bmatrix} 1 & 0 & 4 & 2 & 0 & 1 & 1 \\ 0 & 1 & 0 & 1 & 0 & 1 & 1 \\ 1 & 0 & 1 & 1 & 1 & 0 & 0 \\ 5 & 1 & 1 & 0 & 0 & 1 & 0 \end{bmatrix} \cup$$

$$\begin{bmatrix} 1 & 2 & 1 \\ 1 & 1 & 1 \\ 0 & 1 & 0 \\ \hline 1 & 0 & 1 \\ 2 & 1 & 1 \\ 1 & 0 & 1 \\ 1 & 1 & 1 \\ 2 & 0 & 2 \end{bmatrix} \begin{bmatrix} 1 & 1 & 0 & 1 & 1 & 0 & 1 & 1 & 0 & 1 & 1 \\ 0 & 1 & 1 & 0 & 0 & 1 & 0 & 1 & 6 & 1 & 0 \\ 0 & 1 & 1 & 1 & 1 & 1 & 0 & 1 & 0 & 0 & 0 \end{bmatrix} \cup$$

$$\begin{bmatrix} 2 & 1 & 0 \\ \hline 1 & 1 & 1 \\ 0 & 1 & 1 \\ \hline 0 & 1 & 2 \\ 6 & 0 & 0 \\ 1 & 0 & 1 \\ \hline 1 & 1 & 1 \\ 0 & 0 & 1 \\ 0 & 1 & 0 \end{bmatrix} \begin{bmatrix} 1 & 1 & 1 & 2 & 2 & 0 & 1 & 1 & 1 \\ 0 & 6 & 1 & 1 & 1 & 1 & 0 & 1 & 0 \\ 1 & 0 & 0 & 0 & 0 & 0 & 1 & 0 & 1 \end{bmatrix}$$

$$= \left[\begin{array}{c|c} (1\ 2)\begin{bmatrix}1 & 3\\ 0 & 1\end{bmatrix} & (1\ 2)\begin{bmatrix}1 & 0 & 1 & 1\\ 1 & 1 & 0 & 0\end{bmatrix} \\ \hline \begin{bmatrix}3 & 0\\ 1 & 2\\ 1 & 4\end{bmatrix}\begin{bmatrix}1 & 3\\ 0 & 1\end{bmatrix} & \begin{bmatrix}3 & 0\\ 1 & 2\\ 1 & 4\end{bmatrix}\begin{bmatrix}1 & 0 & 1 & 1\\ 1 & 1 & 0 & 0\end{bmatrix} \end{array}\right] \cup$$

$$\left[\begin{array}{c|c} \begin{pmatrix}3 & 0 & 1\\ 1 & 1 & 1\end{pmatrix}\begin{bmatrix}1 & 1 & 3\\ 3 & 0 & 0\\ 0 & 1 & 0\end{bmatrix} & \begin{pmatrix}3 & 0 & 1\\ 1 & 1 & 1\end{pmatrix}\begin{bmatrix}1 & 1 & 1 & 1\\ 0 & 1 & 0 & 1\\ 2 & 1 & 0 & 1\end{bmatrix} \\ \hline \begin{pmatrix}2 & 3 & 1\\ 1 & 0 & 1\\ 1 & 5 & 2\end{pmatrix}\begin{bmatrix}1 & 1 & 3\\ 3 & 0 & 0\\ 0 & 1 & 0\end{bmatrix} & \begin{pmatrix}2 & 3 & 1\\ 1 & 0 & 1\\ 1 & 5 & 2\end{pmatrix}\begin{bmatrix}1 & 1 & 1 & 1\\ 0 & 1 & 0 & 1\\ 2 & 1 & 0 & 1\end{bmatrix} \end{array}\right]$$

$$\cup \left[\begin{array}{c|c} \begin{bmatrix}1 & 2 & 3 & 4\\ 0 & 1 & 1 & 1\\ 1 & 1 & 0 & 1\\ 0 & 1 & 1 & 0\end{bmatrix}\begin{bmatrix}1 & 0 & 4\\ 0 & 1 & 0\\ 1 & 0 & 1\\ 5 & 1 & 1\end{bmatrix} & \begin{bmatrix}1 & 2 & 3 & 4\\ 0 & 1 & 1 & 1\\ 1 & 1 & 0 & 1\\ 0 & 1 & 1 & 0\end{bmatrix}\begin{bmatrix}2 & 0 & 1 & 1\\ 1 & 0 & 1 & 1\\ 1 & 1 & 0 & 0\\ 0 & 0 & 1 & 0\end{bmatrix} \\ \hline \begin{bmatrix}1 & 2 & 0 & 1\\ 0 & 1 & 1 & 0\end{bmatrix}\begin{bmatrix}1 & 0 & 4\\ 0 & 1 & 0\\ 1 & 0 & 1\\ 5 & 1 & 1\end{bmatrix} & \begin{bmatrix}1 & 2 & 0 & 1\\ 0 & 1 & 1 & 0\end{bmatrix}\begin{bmatrix}2 & 0 & 1 & 1\\ 1 & 0 & 1 & 1\\ 1 & 1 & 0 & 0\\ 0 & 0 & 1 & 0\end{bmatrix} \end{array}\right]$$

$$\cup \left[\begin{array}{c|c} \begin{pmatrix}1 & 2 & 1\\ 1 & 1 & 1\\ 0 & 1 & 0\end{pmatrix}\begin{bmatrix}1 & 1 & 0 & 1\\ 0 & 1 & 1 & 0\\ 0 & 1 & 1 & 1\end{bmatrix} & \begin{pmatrix}1 & 2 & 1\\ 1 & 1 & 1\\ 0 & 1 & 0\end{pmatrix}\begin{pmatrix}1 & 0 & 1\\ 0 & 1 & 0\\ 1 & 1 & 0\end{pmatrix} \\ \hline \begin{pmatrix}1 & 0 & 1\\ 2 & 1 & 1\\ 1 & 0 & 1\\ 1 & 1 & 1\\ 2 & 0 & 2\end{pmatrix}\begin{bmatrix}1 & 1 & 0 & 1\\ 0 & 1 & 1 & 0\\ 0 & 1 & 1 & 1\end{bmatrix} & \begin{pmatrix}1 & 0 & 1\\ 2 & 1 & 1\\ 1 & 0 & 1\\ 1 & 1 & 1\\ 2 & 0 & 2\end{pmatrix}\begin{pmatrix}1 & 0 & 1\\ 0 & 1 & 0\\ 1 & 1 & 0\end{pmatrix} \end{array}\right]$$

$$\left[ \frac{\begin{pmatrix} 1 & 2 & 1 \\ 1 & 1 & 1 \\ 0 & 1 & 0 \end{pmatrix} \begin{bmatrix} 1 & 0 & 1 & 1 \\ 1 & 6 & 1 & 0 \\ 1 & 0 & 0 & 0 \end{bmatrix}}{\begin{pmatrix} 1 & 0 & 1 \\ 2 & 1 & 1 \\ 1 & 0 & 1 \\ 1 & 1 & 1 \\ 2 & 0 & 2 \end{pmatrix} \begin{bmatrix} 1 & 0 & 1 & 1 \\ 1 & 6 & 1 & 0 \\ 1 & 0 & 0 & 0 \end{bmatrix}} \right] \cup$$

$$\left[ \begin{array}{c|c|c} [2\ 1\ 0]\begin{bmatrix}1 & 1 & 1\\0 & 6 & 1\\1 & 0 & 0\end{bmatrix} & [2\ 1\ 0]\begin{bmatrix}2 & 2 & 0 & 1\\1 & 1 & 1 & 0\\0 & 0 & 0 & 1\end{bmatrix} & [2\ 1\ 0]\begin{bmatrix}1 & 1\\1 & 0\\0 & 1\end{bmatrix} \\ \hline \begin{bmatrix}0 & 1 & 2\\6 & 0 & 0\\1 & 0 & 1\end{bmatrix}\begin{bmatrix}1 & 1 & 1\\0 & 6 & 1\\1 & 0 & 0\end{bmatrix} & \begin{bmatrix}0 & 1 & 2\\6 & 0 & 0\\1 & 0 & 1\end{bmatrix}\begin{bmatrix}2 & 2 & 0 & 1\\1 & 1 & 1 & 0\\0 & 0 & 0 & 1\end{bmatrix} & \begin{bmatrix}0 & 1 & 2\\6 & 0 & 0\\1 & 0 & 1\end{bmatrix}\begin{bmatrix}1 & 1\\1 & 0\\0 & 1\end{bmatrix} \\ \hline \begin{bmatrix}1 & 1 & 1\\0 & 0 & 1\\0 & 1 & 0\end{bmatrix}\begin{bmatrix}1 & 1 & 1\\0 & 6 & 1\\1 & 0 & 0\end{bmatrix} & \begin{bmatrix}1 & 1 & 1\\0 & 0 & 1\\0 & 1 & 0\end{bmatrix}\begin{bmatrix}2 & 2 & 0 & 1\\1 & 1 & 1 & 0\\0 & 0 & 0 & 1\end{bmatrix} & \begin{bmatrix}1 & 1 & 1\\0 & 0 & 1\\0 & 1 & 0\end{bmatrix}\begin{bmatrix}1 & 1\\1 & 0\\0 & 1\end{bmatrix} \end{array} \right]$$

$$= \left[ \begin{array}{cc|cccc} 1 & 5 & 3 & 2 & 1 & 1 \\ \hline 3 & 9 & 3 & 0 & 3 & 3 \\ 1 & 5 & 3 & 2 & 1 & 1 \\ 1 & 7 & 5 & 4 & 1 & 1 \end{array} \right] \cup \left[ \begin{array}{ccc|ccc} 3 & 4 & 9 & 5 & 4 & 3 & 4 \\ 4 & 2 & 3 & 3 & 3 & 1 & 3 \\ \hline 11 & 3 & 6 & 4 & 6 & 2 & 6 \\ 1 & 2 & 3 & 3 & 2 & 1 & 2 \\ 16 & 3 & 3 & 5 & 8 & 1 & 8 \end{array} \right]$$

$$\begin{bmatrix} 24 & 6 & 11 & 7 & 3 & 7 & 3 \\ 6 & 2 & 2 & 2 & 1 & 2 & 1 \\ 6 & 2 & 5 & 3 & 0 & 3 & 2 \\ 1 & 1 & 1 & 2 & 1 & 1 & 1 \\ \hline 6 & 3 & 6 & 4 & 0 & 4 & 3 \\ 1 & 1 & 1 & 2 & 1 & 1 & 1 \end{bmatrix} \cup$$

$$\begin{bmatrix} 1 & 4 & 3 & 2 & 2 & 3 & 1 & 4 & 12 & 3 & 1 \\ 1 & 3 & 2 & 2 & 2 & 2 & 1 & 3 & 6 & 2 & 1 \\ 0 & 1 & 1 & 0 & 0 & 1 & 0 & 1 & 6 & 1 & 0 \\ \hline 1 & 2 & 1 & 2 & 2 & 1 & 1 & 2 & 6 & 1 & 1 \\ 2 & 4 & 2 & 3 & 3 & 2 & 2 & 4 & 6 & 3 & 2 \\ 1 & 2 & 1 & 2 & 2 & 1 & 1 & 2 & 6 & 1 & 1 \\ 1 & 3 & 2 & 2 & 2 & 2 & 1 & 3 & 6 & 2 & 1 \\ 2 & 4 & 2 & 2 & 4 & 2 & 2 & 4 & 0 & 2 & 2 \end{bmatrix} \cup$$

$$\begin{bmatrix} 2 & 8 & 3 & 5 & 5 & 1 & 2 & 3 & 2 \\ \hline 2 & 6 & 1 & 1 & 1 & 1 & 2 & 1 & 2 \\ 6 & 6 & 6 & 12 & 12 & 0 & 6 & 6 & 6 \\ 2 & 1 & 1 & 2 & 2 & 0 & 2 & 1 & 2 \\ \hline 2 & 7 & 2 & 3 & 3 & 1 & 2 & 2 & 2 \\ 1 & 0 & 0 & 0 & 0 & 0 & 1 & 0 & 1 \\ 0 & 6 & 1 & 1 & 1 & 1 & 0 & 1 & 0 \end{bmatrix}.$$

We see the minor product of T with V yields a super n-matrix (n = 5). It is not a special semi n-matrix or a special semi column n-vector or a special semi row n-vector. It is a super n-matrix (n = 5).

Now we study under what conditions is the product defined and is compatible. We see if T and V are special semi super n-vectors then both T and V must have n elements i.e., n-matrices for the product to be defined. Secondly we need in both $T_iV_i$ is compatible with respect to minor product i.e., $T = T_1 \cup T_2 \cup \ldots \cup T_n$ and $V = V_1 \cup V_2 \cup \ldots \cup V_n$, $1 \leq i \leq n$ we have $T_i$ to be special column super vector and $V_i$ to be special row super vector such that the number of columns in $T_i$ = number of rows of $V_i$; $1 \leq i \leq n$. Thirdly we see the resultant of the minor product yields a super n-matrix and not a semi super n-matrix provided both T and V are special super n-vectors. Now we will find the major product of a special column super n-matrix with its transpose.

*Example 4.23:* Let $T = T_1 \cup T_2 \cup T_3 \cup T_4 \cup T_5$ be a semi super column 5-vector. To find the product $TT^T$. Given $T = T_1 \cup T_2 \cup T_3 \cup T_4 \cup T_5$

$$= \begin{bmatrix} 3 & 1 & 0 \\ 1 & 4 & 1 \\ 0 & 1 & 1 \\ \hline 1 & 1 & 1 \\ 0 & 1 & 0 \\ \hline 1 & 0 & 1 \\ 1 & 1 & 0 \end{bmatrix} \cup \begin{bmatrix} 3 & 4 & 0 & 1 \\ 1 & 2 & 0 & 0 \\ 0 & 1 & 0 & 2 \\ \hline 1 & 1 & 0 & 1 \\ 1 & 1 & 0 & 0 \\ 1 & 0 & 0 & 1 \\ \hline 1 & 0 & 4 & 5 \\ 1 & 1 & 0 & 2 \end{bmatrix} \cup \begin{bmatrix} 3 & 0 & 1 & 1 & 5 \\ 1 & 0 & 5 & 2 & 1 \\ 1 & 1 & 0 & 1 & 4 \\ \hline 2 & 1 & 0 & 1 & 0 \\ 1 & 1 & 0 & 1 & 0 \\ \hline 1 & 0 & 1 & 0 & 0 \\ 0 & 0 & 0 & 1 & 1 \\ \hline 1 & 1 & 0 & 1 & 1 \\ 0 & 0 & 1 & 0 & 0 \\ 0 & 0 & 0 & 1 & 0 \end{bmatrix}$$

$$\cup \begin{bmatrix} 1 & 2 & 1 \\ 1 & 1 & 1 \\ 1 & 0 & 1 \\ \hline 3 & 4 & 0 \\ 0 & 1 & 1 \\ 1 & 1 & 1 \\ 1 & 0 & 1 \\ 2 & 0 & 1 \\ 1 & 2 & 0 \end{bmatrix} \cup \begin{bmatrix} 3 & 1 & 1 & 1 \\ \hline 1 & 0 & 1 & 1 \\ 1 & 1 & 0 & 1 \\ 2 & 1 & 1 & 1 \\ 1 & 0 & 0 & 1 \\ \hline 0 & 1 & 0 & 1 \\ 0 & 1 & 1 & 1 \\ 0 & 1 & 1 & 0 \\ 0 & 1 & 0 & 0 \end{bmatrix}.$$

Now
$$T^T = (T_1 \cup T_2 \cup T_3 \cup T_4 \cup T_5)^T$$
$$= T_1^T \cup T_2^T \cup T_3^T \cup T_4^T \cup T_5^T$$

$$= \begin{bmatrix} 3 & 1 & 0 & 1 & 0 & 1 & 1 \\ 1 & 4 & 1 & 1 & 1 & 0 & 1 \\ 0 & 1 & 1 & 1 & 0 & 1 & 0 \end{bmatrix}$$

$$\cup \begin{bmatrix} 3 & 1 & 0 & 1 & 1 & 1 & 1 \\ 4 & 2 & 1 & 1 & 1 & 0 & 0 & 1 \\ 0 & 0 & 0 & 0 & 0 & 0 & 4 & 0 \\ 1 & 0 & 2 & 1 & 0 & 1 & 5 & 2 \end{bmatrix}$$

$$\cup \begin{bmatrix} 3 & 1 & 1 & 2 & 1 & 1 & 0 & 1 & 0 & 0 \\ 0 & 0 & 1 & 1 & 1 & 0 & 0 & 1 & 0 & 0 \\ 1 & 5 & 0 & 0 & 0 & 1 & 0 & 0 & 1 & 0 \\ 1 & 2 & 1 & 1 & 1 & 0 & 1 & 1 & 0 & 1 \\ 5 & 1 & 4 & 0 & 0 & 0 & 1 & 1 & 0 & 0 \end{bmatrix}$$

$$\cup \begin{bmatrix} 1 & 1 & 1 & 3 & 0 & 1 & 1 & 2 & 1 \\ 2 & 1 & 0 & 4 & 1 & 1 & 0 & 0 & 2 \\ 1 & 1 & 1 & 0 & 1 & 1 & 1 & 1 & 0 \end{bmatrix}$$

$$\cup \begin{bmatrix} 3 & 1 & 1 & 2 & 1 & 0 & 0 & 0 & 0 \\ 1 & 0 & 1 & 1 & 0 & 1 & 1 & 1 & 1 \\ 1 & 1 & 0 & 1 & 0 & 0 & 1 & 1 & 0 \\ 1 & 1 & 1 & 1 & 1 & 1 & 1 & 0 & 0 \end{bmatrix}.$$

$$TT^T = (T_1 \cup T_2 \cup T_3 \cup T_4 \cup T_5)(T_1^T \cup T_2^T \cup \ldots \cup T_5^T)$$
$$= T_1 T_1^T \cup T_2 T_2^T \cup T_3 T_3^T \cup T_4 T_4^T \cup T_5 T_5^T$$

$$= \begin{bmatrix} 3 & 1 & 0 \\ 1 & 4 & 1 \\ 0 & 1 & 1 \\ \hline 1 & 1 & 1 \\ 0 & 1 & 0 \\ 1 & 0 & 1 \\ 1 & 1 & 0 \end{bmatrix} \left[ \begin{array}{c|cc|cccc} 3 & 1 & 0 & 1 & 0 & 1 & 1 \\ 1 & 4 & 1 & 1 & 1 & 0 & 1 \\ 0 & 1 & 1 & 1 & 0 & 1 & 0 \end{array} \right] \cup$$

$$\begin{bmatrix} 3 & 4 & 0 & 1 \\ 1 & 2 & 0 & 0 \\ 0 & 1 & 0 & 2 \\ \hline 1 & 1 & 0 & 1 \\ 1 & 1 & 0 & 0 \\ 1 & 0 & 0 & 1 \\ 1 & 0 & 4 & 5 \\ \hline 1 & 1 & 0 & 2 \end{bmatrix} \left[ \begin{array}{ccc|cccc|c} 3 & 1 & 0 & 1 & 1 & 1 & 1 & 1 \\ 4 & 2 & 1 & 1 & 1 & 0 & 0 & 1 \\ 0 & 0 & 0 & 0 & 0 & 0 & 4 & 0 \\ 1 & 0 & 2 & 1 & 0 & 1 & 5 & 2 \end{array} \right] \cup$$

$$\begin{bmatrix} 3 & 0 & 1 & 1 & 5 \\ 1 & 0 & 5 & 2 & 1 \\ \hline 1 & 1 & 0 & 1 & 4 \\ 2 & 1 & 0 & 1 & 0 \\ 1 & 1 & 0 & 1 & 0 \\ \hline 1 & 0 & 1 & 0 & 0 \\ 0 & 0 & 0 & 1 & 1 \\ 1 & 1 & 0 & 1 & 1 \\ 0 & 0 & 1 & 0 & 0 \\ 0 & 0 & 0 & 1 & 0 \end{bmatrix} \left[ \begin{array}{cc|ccc|ccccc} 3 & 1 & 1 & 2 & 1 & 1 & 0 & 1 & 0 & 0 \\ 0 & 0 & 1 & 1 & 1 & 0 & 0 & 1 & 0 & 0 \\ 1 & 5 & 0 & 0 & 0 & 1 & 0 & 0 & 1 & 0 \\ 1 & 2 & 1 & 1 & 1 & 0 & 1 & 1 & 0 & 1 \\ 5 & 1 & 4 & 0 & 0 & 0 & 1 & 1 & 0 & 0 \end{array} \right]$$

$$\cup \begin{bmatrix} 1 & 2 & 1 \\ 1 & 1 & 1 \\ 1 & 0 & 1 \\ \hline 3 & 4 & 0 \\ 0 & 1 & 1 \\ 1 & 1 & 1 \\ \hline 1 & 0 & 1 \\ 2 & 0 & 1 \\ 1 & 2 & 0 \end{bmatrix} \begin{bmatrix} 1 & 1 & 1 & \vline & 3 & 0 & 1 & 1 & \vline & 2 & 1 \\ 2 & 1 & 0 & \vline & 4 & 1 & 1 & 0 & \vline & 0 & 2 \\ 1 & 1 & 1 & \vline & 0 & 1 & 1 & 1 & \vline & 1 & 0 \end{bmatrix} \cup$$

$$\begin{bmatrix} 3 & 1 & 1 & 1 \\ \hline 1 & 0 & 1 & 1 \\ 1 & 1 & 0 & 1 \\ \hline 2 & 1 & 1 & 1 \\ 1 & 0 & 0 & 1 \\ 0 & 1 & 0 & 1 \\ 0 & 1 & 1 & 1 \\ 0 & 1 & 1 & 0 \\ 0 & 1 & 0 & 0 \end{bmatrix} \begin{bmatrix} 3 & \vline & 1 & 1 & \vline & 2 & 1 & 0 & 0 & 0 & 0 \\ 1 & \vline & 0 & 1 & \vline & 1 & 0 & 1 & 1 & 1 & 1 \\ 1 & \vline & 1 & 0 & \vline & 1 & 0 & 0 & 1 & 1 & 0 \\ 1 & \vline & 1 & 1 & \vline & 1 & 1 & 1 & 1 & 0 & 0 \end{bmatrix} =$$

$$\begin{bmatrix} (3\ 1\ 0)\begin{pmatrix} 3 \\ 1 \\ 0 \end{pmatrix} & \vline & (3\ 1\ 0)\begin{pmatrix} 1 & 0 \\ 4 & 1 \\ 1 & 1 \end{pmatrix} & \vline & (3\ 1\ 0)\begin{bmatrix} 1 & 0 & 1 & 1 \\ 1 & 1 & 0 & 1 \\ 1 & 0 & 1 & 0 \end{bmatrix} \\ \hline \begin{pmatrix} 1 & 4 & 1 \\ 0 & 1 & 1 \end{pmatrix}\begin{pmatrix} 3 \\ 1 \\ 0 \end{pmatrix} & \vline & \begin{pmatrix} 1 & 4 & 1 \\ 0 & 1 & 1 \end{pmatrix}\begin{pmatrix} 1 & 0 \\ 4 & 1 \\ 1 & 1 \end{pmatrix} & \vline & \begin{pmatrix} 1 & 4 & 1 \\ 0 & 1 & 1 \end{pmatrix}\begin{bmatrix} 1 & 0 & 1 & 1 \\ 1 & 1 & 0 & 1 \\ 1 & 0 & 1 & 0 \end{bmatrix} \\ \hline \begin{pmatrix} 1 & 1 & 1 \\ 0 & 1 & 0 \\ 1 & 0 & 1 \\ 1 & 1 & 0 \end{pmatrix}\begin{pmatrix} 3 \\ 1 \\ 0 \end{pmatrix} & \vline & \begin{pmatrix} 1 & 1 & 1 \\ 0 & 1 & 0 \\ 1 & 0 & 1 \\ 1 & 1 & 0 \end{pmatrix}\begin{pmatrix} 1 & 0 \\ 4 & 1 \\ 1 & 1 \end{pmatrix} & \vline & \begin{pmatrix} 1 & 1 & 1 \\ 0 & 1 & 0 \\ 1 & 0 & 1 \\ 1 & 1 & 0 \end{pmatrix}\begin{bmatrix} 1 & 0 & 1 & 1 \\ 1 & 1 & 0 & 1 \\ 1 & 0 & 1 & 0 \end{bmatrix} \end{bmatrix} \cup$$

$$
\left[
\begin{array}{|c|c|}
\hline
\begin{pmatrix} 3 & 4 & 0 & 1 \\ 1 & 2 & 0 & 0 \\ 0 & 1 & 0 & 2 \end{pmatrix}\begin{bmatrix} 3 & 1 & 0 \\ 4 & 2 & 1 \\ 0 & 0 & 0 \\ 1 & 0 & 2 \end{bmatrix} & \begin{pmatrix} 3 & 4 & 0 & 1 \\ 1 & 2 & 0 & 0 \\ 0 & 1 & 0 & 2 \end{pmatrix}\begin{bmatrix} 1 & 1 & 1 & 1 \\ 1 & 1 & 0 & 0 \\ 0 & 0 & 0 & 4 \\ 1 & 0 & 1 & 5 \end{bmatrix} \\
\hline
\begin{pmatrix} 1 & 1 & 0 & 1 \\ 1 & 1 & 0 & 0 \\ 1 & 0 & 0 & 1 \\ 1 & 0 & 4 & 5 \end{pmatrix}\begin{bmatrix} 3 & 1 & 0 \\ 4 & 2 & 1 \\ 0 & 0 & 0 \\ 1 & 0 & 2 \end{bmatrix} & \begin{pmatrix} 1 & 1 & 0 & 1 \\ 1 & 1 & 0 & 0 \\ 1 & 0 & 0 & 1 \\ 1 & 0 & 4 & 5 \end{pmatrix}\begin{bmatrix} 1 & 1 & 1 & 1 \\ 1 & 1 & 0 & 0 \\ 0 & 0 & 0 & 4 \\ 1 & 0 & 1 & 5 \end{bmatrix} \\
\hline
\begin{pmatrix} 1 & 1 & 0 & 2 \end{pmatrix}\begin{bmatrix} 3 & 1 & 0 \\ 4 & 2 & 1 \\ 0 & 0 & 0 \\ 1 & 0 & 2 \end{bmatrix} & \begin{pmatrix} 1 & 1 & 0 & 2 \end{pmatrix}\begin{bmatrix} 1 & 1 & 1 & 1 \\ 1 & 1 & 0 & 0 \\ 0 & 0 & 0 & 4 \\ 1 & 0 & 1 & 5 \end{bmatrix} \\
\hline
\end{array}
\right]
$$

$$
\left[
\begin{array}{c}
\begin{bmatrix} 3 & 4 & 0 & 1 \\ 1 & 2 & 0 & 0 \\ 0 & 1 & 0 & 2 \end{bmatrix}\begin{bmatrix} 1 \\ 1 \\ 0 \\ 2 \end{bmatrix} \\
\hline
\begin{bmatrix} 1 & 1 & 0 & 1 \\ 1 & 1 & 0 & 0 \\ 1 & 0 & 0 & 1 \\ 1 & 0 & 4 & 5 \end{bmatrix}\begin{bmatrix} 1 \\ 1 \\ 0 \\ 2 \end{bmatrix} \\
\hline
\begin{bmatrix} 1 & 1 & 0 & 2 \end{bmatrix}\begin{bmatrix} 1 \\ 1 \\ 0 \\ 2 \end{bmatrix}
\end{array}
\right] \subset
$$

$$\left[\begin{array}{c|c}
\begin{pmatrix} 3 & 0 & 1 & 1 & 5 \\ 1 & 0 & 5 & 2 & 1 \end{pmatrix}\begin{bmatrix} 3 & 1 \\ 0 & 0 \\ 1 & 5 \\ 1 & 2 \\ 5 & 1 \end{bmatrix} & \begin{pmatrix} 3 & 0 & 1 & 1 & 5 \\ 1 & 0 & 5 & 2 & 1 \end{pmatrix}\begin{bmatrix} 1 & 2 & 1 \\ 1 & 1 & 1 \\ 0 & 0 & 0 \\ 1 & 1 & 1 \\ 4 & 0 & 0 \end{bmatrix} \\
\hline
\begin{pmatrix} 1 & 1 & 0 & 1 & 4 \\ 2 & 1 & 0 & 1 & 0 \\ 1 & 1 & 0 & 1 & 0 \end{pmatrix}\begin{bmatrix} 3 & 1 \\ 0 & 0 \\ 1 & 5 \\ 1 & 2 \\ 5 & 1 \end{bmatrix} & \begin{pmatrix} 1 & 1 & 0 & 1 & 4 \\ 2 & 1 & 0 & 1 & 0 \\ 1 & 1 & 0 & 1 & 0 \end{pmatrix}\begin{bmatrix} 1 & 2 & 1 \\ 1 & 1 & 1 \\ 0 & 0 & 0 \\ 1 & 1 & 1 \\ 4 & 0 & 0 \end{bmatrix} \\
\hline
\begin{pmatrix} 1 & 0 & 1 & 0 & 0 \\ 0 & 0 & 0 & 1 & 1 \\ 1 & 1 & 0 & 1 & 1 \\ 0 & 0 & 1 & 0 & 0 \\ 0 & 0 & 0 & 1 & 0 \end{pmatrix}\begin{bmatrix} 3 & 1 \\ 0 & 0 \\ 1 & 5 \\ 1 & 2 \\ 5 & 1 \end{bmatrix} & \begin{pmatrix} 1 & 0 & 1 & 0 & 0 \\ 0 & 0 & 0 & 1 & 1 \\ 1 & 1 & 0 & 1 & 1 \\ 0 & 0 & 1 & 0 & 0 \\ 0 & 0 & 0 & 1 & 0 \end{pmatrix}\begin{bmatrix} 1 & 2 & 1 \\ 1 & 1 & 1 \\ 0 & 0 & 0 \\ 1 & 1 & 1 \\ 4 & 0 & 0 \end{bmatrix}
\end{array}\right]$$

$$\subset \left[\begin{array}{c}
\begin{pmatrix} 3 & 0 & 1 & 1 & 5 \\ 1 & 0 & 5 & 2 & 1 \end{pmatrix}\begin{bmatrix} 1 & 0 & 1 & 0 & 0 \\ 0 & 0 & 1 & 0 & 0 \\ 1 & 0 & 0 & 1 & 0 \\ 0 & 1 & 1 & 0 & 1 \\ 0 & 1 & 1 & 0 & 0 \end{bmatrix} \\
\hline
\begin{pmatrix} 1 & 1 & 0 & 1 & 4 \\ 2 & 1 & 0 & 1 & 0 \\ 1 & 1 & 0 & 1 & 0 \end{pmatrix}\begin{bmatrix} 1 & 0 & 1 & 0 & 0 \\ 0 & 0 & 1 & 0 & 0 \\ 1 & 0 & 0 & 1 & 0 \\ 0 & 1 & 1 & 0 & 1 \\ 0 & 1 & 1 & 0 & 0 \end{bmatrix} \\
\hline
\begin{pmatrix} 1 & 0 & 1 & 0 & 0 \\ 0 & 0 & 0 & 1 & 1 \\ 1 & 1 & 0 & 1 & 1 \\ 0 & 0 & 1 & 0 & 0 \\ 0 & 0 & 0 & 1 & 0 \end{pmatrix}\begin{bmatrix} 1 & 0 & 1 & 0 & 0 \\ 0 & 0 & 1 & 0 & 0 \\ 1 & 0 & 0 & 1 & 0 \\ 0 & 1 & 1 & 0 & 1 \\ 0 & 1 & 1 & 0 & 0 \end{bmatrix}
\end{array}\right]$$

$$\begin{bmatrix} \begin{pmatrix} 1 & 2 & 1 \\ 1 & 1 & 1 \\ 1 & 0 & 1 \end{pmatrix} \begin{bmatrix} 1 & 1 & 1 \\ 2 & 1 & 0 \\ 1 & 1 & 1 \end{bmatrix} & \begin{pmatrix} 1 & 2 & 1 \\ 1 & 1 & 1 \\ 1 & 0 & 1 \end{pmatrix} \begin{bmatrix} 3 & 0 & 1 & 1 \\ 4 & 1 & 1 & 0 \\ 0 & 1 & 1 & 1 \end{bmatrix} & \begin{bmatrix} 1 & 2 & 1 \\ 1 & 1 & 1 \\ 1 & 0 & 1 \end{bmatrix} \begin{bmatrix} 2 & 1 \\ 0 & 2 \\ 1 & 0 \end{bmatrix} \\ \begin{bmatrix} 3 & 4 & 0 \\ 0 & 1 & 1 \\ 1 & 1 & 1 \\ 1 & 0 & 1 \end{bmatrix} \begin{bmatrix} 1 & 1 & 1 \\ 2 & 1 & 0 \\ 1 & 1 & 1 \end{bmatrix} & \begin{bmatrix} 3 & 4 & 0 \\ 0 & 1 & 1 \\ 1 & 1 & 1 \\ 1 & 0 & 1 \end{bmatrix} \begin{bmatrix} 3 & 0 & 1 & 1 \\ 4 & 1 & 1 & 0 \\ 0 & 1 & 1 & 1 \end{bmatrix} & \begin{bmatrix} 3 & 4 & 0 \\ 0 & 1 & 1 \\ 1 & 1 & 1 \\ 1 & 0 & 1 \end{bmatrix} \begin{bmatrix} 2 & 1 \\ 0 & 2 \\ 1 & 0 \end{bmatrix} \\ \begin{bmatrix} 2 & 0 & 1 \\ 1 & 2 & 0 \end{bmatrix} \begin{bmatrix} 1 & 1 & 1 \\ 2 & 1 & 0 \\ 1 & 1 & 1 \end{bmatrix} & \begin{bmatrix} 2 & 0 & 1 \\ 1 & 2 & 0 \end{bmatrix} \begin{bmatrix} 3 & 0 & 1 & 1 \\ 4 & 1 & 1 & 0 \\ 0 & 1 & 1 & 1 \end{bmatrix} & \begin{bmatrix} 2 & 0 & 1 \\ 1 & 2 & 0 \end{bmatrix} \begin{bmatrix} 2 & 1 \\ 0 & 2 \\ 1 & 0 \end{bmatrix} \end{bmatrix}$$

$$\subset \begin{bmatrix} (3\ 1\ 1\ 1)\begin{bmatrix} 3 \\ 1 \\ 1 \\ 1 \end{bmatrix} & (3\ 1\ 1\ 1)\begin{bmatrix} 1 & 1 \\ 0 & 1 \\ 1 & 0 \\ 1 & 1 \end{bmatrix} \\ \begin{pmatrix} 1 & 0 & 1 & 1 \\ 1 & 1 & 0 & 1 \end{pmatrix}\begin{bmatrix} 3 \\ 1 \\ 1 \\ 1 \end{bmatrix} & \begin{pmatrix} 1 & 0 & 1 & 1 \\ 1 & 1 & 0 & 1 \end{pmatrix}\begin{bmatrix} 1 & 1 \\ 0 & 1 \\ 1 & 0 \\ 1 & 1 \end{bmatrix} \\ \begin{bmatrix} 2 & 1 & 1 & 1 \\ 1 & 0 & 0 & 1 \\ 0 & 1 & 0 & 1 \\ 0 & 1 & 1 & 1 \\ 0 & 1 & 1 & 0 \\ 0 & 1 & 0 & 0 \end{bmatrix}\begin{bmatrix} 3 \\ 1 \\ 1 \\ 1 \end{bmatrix} & \begin{bmatrix} 2 & 1 & 1 & 1 \\ 1 & 0 & 0 & 1 \\ 0 & 1 & 0 & 1 \\ 0 & 1 & 1 & 1 \\ 0 & 1 & 1 & 0 \\ 0 & 1 & 0 & 0 \end{bmatrix}\begin{bmatrix} 1 & 1 \\ 0 & 1 \\ 1 & 0 \\ 1 & 1 \end{bmatrix} \end{bmatrix}$$

$$\left| \begin{array}{c} (3\ 1\ 1\ 1) \begin{bmatrix} 2 & 1 & 0 & 0 & 0 & 0 \\ 1 & 0 & 1 & 1 & 1 & 1 \\ 1 & 0 & 0 & 1 & 1 & 0 \\ 1 & 1 & 1 & 1 & 0 & 0 \end{bmatrix} \\ \hline \begin{pmatrix} 1 & 0 & 1 & 1 \\ 1 & 1 & 0 & 1 \end{pmatrix} \begin{bmatrix} 2 & 1 & 0 & 0 & 0 & 0 \\ 1 & 0 & 1 & 1 & 1 & 1 \\ 1 & 0 & 0 & 1 & 1 & 0 \\ 1 & 1 & 1 & 1 & 0 & 0 \end{bmatrix} \\ \hline \begin{bmatrix} 2 & 1 & 1 & 1 \\ 1 & 0 & 0 & 1 \\ 0 & 1 & 0 & 1 \\ 0 & 1 & 1 & 1 \\ 0 & 1 & 1 & 0 \\ 0 & 1 & 0 & 0 \end{bmatrix} \begin{bmatrix} 2 & 1 & 0 & 0 & 0 & 0 \\ 1 & 0 & 1 & 1 & 1 & 1 \\ 1 & 0 & 0 & 1 & 1 & 0 \\ 1 & 1 & 1 & 1 & 0 & 0 \end{bmatrix} \end{array} \right|$$

$$= \left[ \begin{array}{c|cc|ccc} 10 & 7 & 1 & 4 & 1 & 3 & 4 \\ \hline 7 & 18 & 5 & 6 & 4 & 2 & 5 \\ 1 & 5 & 2 & 2 & 1 & 1 & 1 \\ \hline 4 & 6 & 2 & 3 & 1 & 2 & 2 \\ 1 & 4 & 1 & 1 & 1 & 0 & 1 \\ 3 & 2 & 1 & 2 & 0 & 2 & 1 \\ 4 & 5 & 1 & 2 & 1 & 1 & 2 \end{array} \right] \cup$$

$$\begin{bmatrix} 26 & 11 & 6 & 8 & 7 & 4 & 8 & 9 \\ 11 & 5 & 2 & 3 & 3 & 1 & 1 & 3 \\ 6 & 2 & 5 & 3 & 1 & 2 & 10 & 5 \\ \hline 8 & 3 & 3 & 3 & 2 & 2 & 6 & 4 \\ 7 & 3 & 1 & 2 & 2 & 1 & 1 & 2 \\ 4 & 1 & 2 & 2 & 1 & 2 & 6 & 3 \\ 8 & 1 & 10 & 6 & 1 & 6 & 42 & 11 \\ 9 & 3 & 5 & 4 & 2 & 3 & 11 & 6 \end{bmatrix}$$

$$\cup \begin{bmatrix} 36 & 15 & 24 & 7 & 4 & 4 & 6 & 9 & 1 & 1 \\ 15 & 31 & 7 & 4 & 3 & 6 & 3 & 4 & 5 & 2 \\ \hline 24 & 7 & 19 & 4 & 3 & 1 & 5 & 7 & 0 & 1 \\ 7 & 4 & 4 & 6 & 4 & 2 & 1 & 4 & 0 & 1 \\ 4 & 3 & 3 & 4 & 3 & 1 & 1 & 3 & 0 & 1 \\ \hline 4 & 6 & 1 & 2 & 1 & 2 & 0 & 1 & 1 & 0 \\ 6 & 3 & 5 & 1 & 1 & 0 & 2 & 2 & 0 & 1 \\ 9 & 4 & 7 & 4 & 3 & 1 & 2 & 4 & 0 & 1 \\ 1 & 5 & 0 & 0 & 0 & 1 & 0 & 0 & 1 & 0 \\ 1 & 2 & 1 & 1 & 1 & 0 & 1 & 1 & 0 & 1 \end{bmatrix}$$

$$\cup \begin{bmatrix} 6 & 4 & 2 & 11 & 3 & 4 & 2 & 3 & 5 \\ 4 & 3 & 2 & 7 & 2 & 3 & 2 & 3 & 3 \\ 2 & 2 & 2 & 3 & 1 & 2 & 2 & 3 & 1 \\ \hline 11 & 7 & 3 & 25 & 4 & 7 & 3 & 6 & 11 \\ 3 & 2 & 1 & 4 & 2 & 2 & 1 & 1 & 2 \\ 4 & 3 & 2 & 7 & 2 & 3 & 2 & 3 & 3 \\ 2 & 2 & 2 & 3 & 1 & 2 & 2 & 3 & 1 \\ \hline 3 & 3 & 3 & 6 & 1 & 3 & 3 & 5 & 2 \\ 5 & 3 & 1 & 11 & 2 & 3 & 1 & 2 & 5 \end{bmatrix}$$

$$\cup \begin{bmatrix} 12 & 5 & 5 & 9 & 4 & 2 & 3 & 2 & 1 \\ \hline 5 & 3 & 2 & 4 & 2 & 1 & 2 & 1 & 0 \\ 5 & 2 & 3 & 4 & 2 & 2 & 2 & 1 & 1 \\ \hline 9 & 4 & 4 & 7 & 3 & 2 & 3 & 2 & 1 \\ 4 & 2 & 2 & 3 & 2 & 1 & 1 & 0 & 0 \\ 2 & 1 & 2 & 2 & 1 & 2 & 2 & 1 & 1 \\ 3 & 2 & 2 & 3 & 1 & 2 & 3 & 2 & 1 \\ 2 & 1 & 1 & 2 & 0 & 1 & 2 & 2 & 1 \\ 1 & 0 & 1 & 1 & 0 & 1 & 1 & 1 & 1 \end{bmatrix}.$$

The resultant is a symmetric super 5-matrix.

Now we proceed on to define the product of a special semi column super n-matrix with its transpose.

***Example 4.24:*** Let $S = S_1 \cup S_2 \cup S_3 \cup S_4 \cup S_5$ be the given special semi super column 4-vector. To find the minor product of S with $S^T$.
Given

$$S = \begin{bmatrix} 2 & 1 & 1 \\ 0 & 1 & 2 \\ \hline 1 & 0 & 3 \\ 9 & 0 & 1 \\ 1 & 0 & 1 \\ \hline 0 & 1 & 0 \end{bmatrix} \cup \begin{bmatrix} 3 & 0 & 1 & 0 & 1 \\ 4 & 1 & 0 & 1 & 0 \\ 5 & 0 & 0 & 1 & 1 \\ 0 & 1 & 1 & 0 & 0 \\ 1 & 0 & 0 & 0 & 1 \end{bmatrix}$$

$$\cup \begin{bmatrix} 3 & 2 & 5 & 0 \\ \hline 1 & 1 & 0 & 1 \\ 0 & 0 & 1 & 0 \\ 1 & 0 & 0 & 0 \\ 0 & 0 & 0 & 1 \\ \hline 0 & 1 & 0 & 0 \\ 1 & 1 & 0 & 0 \\ 0 & 0 & 1 & 1 \\ \hline 1 & 1 & 0 & 1 \\ 1 & 0 & 0 & 1 \end{bmatrix} \cup \begin{bmatrix} 3 & 0 & 1 & 1 & 2 \\ 1 & 1 & 0 & 1 & 0 \\ 0 & 1 & 1 & 1 & 0 \\ \hline 1 & 1 & 1 & 0 & 0 \\ 0 & 0 & 1 & 1 & 1 \\ \hline 0 & 1 & 0 & 1 & 0 \\ 0 & 1 & 0 & 1 & 1 \\ \hline 0 & 1 & 1 & 0 & 1 \\ 1 & 0 & 0 & 0 & 1 \end{bmatrix}$$

a special semi super column 4-vector.

$$\begin{aligned} S^T &= (S_1 \cup S_2 \cup S_3 \cup S_4)^T \\ &= S_1^T \cup S_2^T \cup S_3^T \cup S_4^T \end{aligned}$$

$$= \begin{bmatrix} 2 & 0 & 1 & 9 & 1 & 0 \\ 1 & 1 & 0 & 0 & 0 & 1 \\ 1 & 2 & 3 & 1 & 1 & 0 \end{bmatrix} \cup$$

$$\begin{bmatrix} 3 & 4 & 5 & 0 & 1 \\ 0 & 1 & 0 & 1 & 0 \\ 1 & 0 & 0 & 1 & 0 \\ 0 & 1 & 1 & 0 & 0 \\ 1 & 0 & 1 & 0 & 1 \end{bmatrix} \cup$$

$$\begin{bmatrix} 3 & 1 & 0 & 1 & 0 & 0 & 1 & 0 & 1 & 1 \\ 2 & 1 & 0 & 0 & 0 & 1 & 1 & 0 & 1 & 0 \\ 5 & 0 & 1 & 0 & 0 & 0 & 0 & 1 & 0 & 0 \\ 0 & 1 & 0 & 0 & 1 & 0 & 0 & 1 & 1 & 1 \end{bmatrix} \cup$$

$$\begin{bmatrix} 3 & 1 & 0 & 1 & 0 & 0 & 0 & 0 & 1 \\ 0 & 1 & 1 & 1 & 0 & 1 & 1 & 1 & 0 \\ 1 & 0 & 1 & 1 & 1 & 0 & 0 & 1 & 0 \\ 1 & 1 & 1 & 0 & 1 & 1 & 1 & 0 & 0 \\ 2 & 0 & 0 & 0 & 1 & 0 & 1 & 1 & 1 \end{bmatrix}.$$

$$\begin{aligned}
SS^T &= (S_1 \cup S_2 \cup S_3 \cup S_4 \cup S_5)(S_1^T \cup S_2^T \cup S_3^T \cup S_4^T) \\
&= S_1 S_1^T \cup S_2 S_2^T \cup S_3 S_3^T \cup S_4 S_4^T
\end{aligned}$$

$$= \begin{bmatrix} 2 & 1 & 1 \\ 0 & 1 & 2 \\ \hline 1 & 0 & 3 \\ 9 & 0 & 1 \\ 1 & 0 & 1 \\ \hline 0 & 1 & 0 \end{bmatrix} \begin{bmatrix} 2 & 0 & 1 & 9 & 1 & 0 \\ 1 & 1 & 0 & 0 & 0 & 1 \\ 1 & 2 & 3 & 1 & 1 & 0 \end{bmatrix} \cup$$

$$\begin{bmatrix} 3 & 0 & 1 & 0 & 1 \\ 4 & 1 & 0 & 1 & 0 \\ 5 & 0 & 0 & 1 & 1 \\ 0 & 1 & 1 & 0 & 0 \\ 1 & 0 & 0 & 0 & 1 \end{bmatrix} \begin{bmatrix} 3 & 4 & 5 & 0 & 1 \\ 0 & 1 & 0 & 1 & 0 \\ 1 & 0 & 0 & 1 & 0 \\ 0 & 1 & 1 & 0 & 0 \\ 1 & 0 & 1 & 0 & 1 \end{bmatrix} \cup$$

$$\begin{bmatrix} \begin{array}{cccc} 3 & 2 & 5 & 0 \\ \hline 1 & 1 & 0 & 1 \\ 0 & 0 & 1 & 0 \\ 1 & 0 & 0 & 0 \\ 0 & 0 & 0 & 1 \\ \hline 0 & 1 & 0 & 0 \\ 1 & 1 & 0 & 0 \\ \hline 0 & 0 & 1 & 1 \\ \hline 1 & 1 & 0 & 1 \\ 1 & 0 & 0 & 1 \end{array} \end{bmatrix} \left[ \begin{array}{c|cccc|ccc|cc} 3 & 1 & 0 & 1 & 0 & 0 & 1 & 0 & 1 & 1 \\ 2 & 1 & 0 & 0 & 0 & 1 & 1 & 0 & 1 & 0 \\ 5 & 0 & 1 & 0 & 0 & 0 & 0 & 1 & 0 & 0 \\ 0 & 1 & 0 & 0 & 1 & 0 & 0 & 1 & 1 & 1 \end{array} \right]$$

$$\cup \begin{bmatrix} \begin{array}{ccccc} 3 & 0 & 1 & 1 & 2 \\ 1 & 1 & 0 & 1 & 0 \\ 0 & 1 & 1 & 1 & 0 \\ 1 & 1 & 1 & 0 & 0 \\ \hline 0 & 0 & 1 & 1 & 1 \\ 0 & 1 & 0 & 1 & 0 \\ 0 & 1 & 0 & 1 & 1 \\ \hline 0 & 1 & 1 & 0 & 1 \\ 1 & 0 & 0 & 0 & 1 \end{array} \end{bmatrix} \left[ \begin{array}{cccc|ccc|cc} 3 & 1 & 0 & 1 & 0 & 0 & 0 & 0 & 1 \\ 0 & 1 & 1 & 1 & 0 & 1 & 1 & 1 & 0 \\ 1 & 0 & 1 & 1 & 1 & 0 & 0 & 1 & 0 \\ 1 & 1 & 1 & 0 & 1 & 1 & 1 & 0 & 0 \\ 2 & 0 & 0 & 0 & 1 & 0 & 1 & 1 & 1 \end{array} \right]$$

$$= \begin{bmatrix} \begin{pmatrix} 2 & 1 & 1 \\ 0 & 1 & 2 \end{pmatrix}\begin{pmatrix} 2 & 0 \\ 1 & 1 \\ 1 & 2 \end{pmatrix} & \begin{pmatrix} 2 & 1 & 1 \\ 0 & 1 & 2 \end{pmatrix}\begin{pmatrix} 1 & 9 & 1 \\ 0 & 0 & 0 \\ 3 & 1 & 1 \end{pmatrix} & \begin{pmatrix} 2 & 1 & 1 \\ 0 & 1 & 2 \end{pmatrix}\begin{pmatrix} 0 \\ 1 \\ 0 \end{pmatrix} \\ \hline \begin{pmatrix} 1 & 0 & 3 \\ 9 & 0 & 1 \\ 1 & 0 & 1 \end{pmatrix}\begin{pmatrix} 2 & 0 \\ 1 & 1 \\ 1 & 2 \end{pmatrix} & \begin{pmatrix} 1 & 0 & 3 \\ 9 & 0 & 1 \\ 1 & 0 & 1 \end{pmatrix}\begin{pmatrix} 1 & 9 & 1 \\ 0 & 0 & 0 \\ 3 & 1 & 1 \end{pmatrix} & \begin{pmatrix} 1 & 0 & 3 \\ 9 & 0 & 1 \\ 1 & 0 & 1 \end{pmatrix}\begin{pmatrix} 0 \\ 1 \\ 0 \end{pmatrix} \\ \hline \begin{pmatrix} 0 & 1 & 0 \end{pmatrix}\begin{pmatrix} 2 & 0 \\ 1 & 1 \\ 1 & 2 \end{pmatrix} & \begin{pmatrix} 0 & 1 & 0 \end{pmatrix}\begin{pmatrix} 1 & 9 & 1 \\ 0 & 0 & 0 \\ 3 & 1 & 1 \end{pmatrix} & \begin{pmatrix} 0 & 1 & 0 \end{pmatrix}\begin{pmatrix} 0 \\ 1 \\ 0 \end{pmatrix} \end{bmatrix}$$

$$\cup \begin{bmatrix} 11 & 12 & 16 & 1 & 4 \\ 12 & 18 & 21 & 1 & 4 \\ 16 & 21 & 27 & 0 & 6 \\ 1 & 1 & 0 & 2 & 0 \\ 4 & 4 & 6 & 0 & 2 \end{bmatrix} \cup$$

$$\left[\begin{array}{c|c}
(3\ 2\ 5\ 0)\begin{bmatrix}3\\2\\5\\0\end{bmatrix} & (3\ 2\ 5\ 0)\begin{bmatrix}1&0&1&0\\1&0&0&0\\0&1&0&0\\1&0&0&1\end{bmatrix} \\
\hline
\begin{pmatrix}1&1&0&1\\0&0&1&0\\1&0&0&0\\0&0&0&1\end{pmatrix}\begin{bmatrix}3\\2\\5\\0\end{bmatrix} & \begin{pmatrix}1&1&0&1\\0&0&1&0\\1&0&0&0\\0&0&0&1\end{pmatrix}\begin{bmatrix}1&0&1&0\\1&0&0&0\\0&1&0&0\\1&0&0&1\end{bmatrix} \\
\hline
\begin{bmatrix}0&1&0&0\\1&1&0&0\\0&0&1&1\end{bmatrix}\begin{bmatrix}3\\2\\5\\0\end{bmatrix} & \begin{bmatrix}0&1&0&0\\1&1&0&0\\0&0&1&1\end{bmatrix}\begin{bmatrix}1&0&1&0\\1&0&0&0\\0&1&0&0\\1&0&0&1\end{bmatrix} \\
\hline
\begin{bmatrix}1&1&0&1\\1&0&0&1\end{bmatrix}\begin{bmatrix}3\\2\\5\\0\end{bmatrix} & \begin{bmatrix}1&1&0&1\\1&0&0&1\end{bmatrix}\begin{bmatrix}1&0&1&0\\1&0&0&0\\0&1&0&0\\1&0&0&1\end{bmatrix}
\end{array}\right.$$

$$\left.\begin{array}{c|c}
(3\ 2\ 5\ 0)\begin{bmatrix}0&1&0\\1&1&0\\0&0&1\\0&0&1\end{bmatrix} & (3\ 2\ 5\ 0)\begin{bmatrix}1&1\\1&0\\0&0\\1&1\end{bmatrix} \\
\hline
\begin{pmatrix}1&1&0&1\\0&0&1&0\\1&0&0&0\\0&0&0&1\end{pmatrix}\begin{bmatrix}0&1&0\\1&1&0\\0&0&1\\0&0&1\end{bmatrix} & \begin{pmatrix}1&1&0&1\\0&0&1&0\\1&0&0&0\\0&0&0&1\end{pmatrix}\begin{bmatrix}1&1\\1&0\\0&0\\1&1\end{bmatrix} \\
\hline
\begin{bmatrix}0&1&0&0\\1&1&0&0\\0&0&1&1\end{bmatrix}\begin{bmatrix}0&1&0\\1&1&0\\0&0&1\\0&0&1\end{bmatrix} & \begin{bmatrix}0&1&0&0\\1&1&0&0\\0&0&1&1\end{bmatrix}\begin{bmatrix}1&1\\1&0\\0&0\\1&1\end{bmatrix} \\
\hline
\begin{bmatrix}1&1&0&1\\1&0&0&1\end{bmatrix}\begin{bmatrix}0&1&0\\1&1&0\\0&0&1\\0&0&1\end{bmatrix} & \begin{bmatrix}1&1&0&1\\1&0&0&1\end{bmatrix}\begin{bmatrix}1&1\\1&0\\0&0\\1&1\end{bmatrix}
\end{array}\right)$$

$$
\left[\begin{array}{ccccc}3 & 0 & 1 & 1 & 2\\1 & 1 & 0 & 1 & 0\\0 & 1 & 1 & 1 & 0\\1 & 1 & 1 & 0 & 0\end{array}\right]\left[\begin{array}{cccc}3 & 1 & 0 & 1\\0 & 1 & 1 & 1\\1 & 0 & 1 & 1\\1 & 1 & 1 & 0\\2 & 0 & 0 & 0\end{array}\right] \quad \left[\begin{array}{ccccc}3 & 0 & 1 & 1 & 2\\1 & 1 & 0 & 1 & 0\\0 & 1 & 1 & 1 & 0\\1 & 1 & 1 & 0 & 0\end{array}\right]\left[\begin{array}{ccc}0 & 0 & 0\\0 & 1 & 1\\1 & 0 & 0\\1 & 1 & 1\\1 & 0 & 1\end{array}\right]
$$

$$
\left(\begin{array}{ccccc}0 & 0 & 1 & 1 & 1\\0 & 1 & 0 & 1 & 0\\0 & 1 & 0 & 1 & 1\end{array}\right)\left[\begin{array}{cccc}3 & 1 & 0 & 1\\0 & 1 & 1 & 1\\1 & 0 & 1 & 1\\1 & 1 & 1 & 0\\2 & 0 & 0 & 0\end{array}\right] \quad \left(\begin{array}{ccccc}0 & 0 & 1 & 1 & 1\\0 & 1 & 0 & 1 & 0\\0 & 1 & 0 & 1 & 1\end{array}\right)\left[\begin{array}{ccc}0 & 0 & 0\\0 & 1 & 1\\1 & 0 & 0\\1 & 1 & 1\\1 & 0 & 1\end{array}\right]
$$

$$
\left(\begin{array}{ccccc}0 & 1 & 1 & 0 & 1\\1 & 0 & 0 & 0 & 1\end{array}\right)\left[\begin{array}{cccc}3 & 1 & 0 & 1\\0 & 1 & 1 & 1\\1 & 0 & 1 & 1\\1 & 1 & 1 & 0\\2 & 0 & 0 & 0\end{array}\right] \quad \left(\begin{array}{ccccc}0 & 1 & 1 & 0 & 1\\1 & 0 & 0 & 0 & 1\end{array}\right)\left[\begin{array}{ccc}0 & 0 & 0\\0 & 1 & 1\\1 & 0 & 0\\1 & 1 & 1\\1 & 0 & 1\end{array}\right]
$$

$$
\left[\begin{array}{ccccc}3 & 0 & 1 & 1 & 2\\1 & 1 & 0 & 1 & 0\\0 & 1 & 1 & 1 & 0\\1 & 1 & 1 & 0 & 0\end{array}\right]\left[\begin{array}{cc}0 & 1\\1 & 0\\1 & 0\\0 & 0\\1 & 1\end{array}\right]
$$

$$
\left(\begin{array}{ccccc}0 & 0 & 1 & 1 & 1\\0 & 1 & 0 & 1 & 0\\0 & 1 & 0 & 1 & 1\end{array}\right)\left[\begin{array}{cc}0 & 1\\1 & 0\\1 & 0\\0 & 0\\1 & 1\end{array}\right]
$$

$$
\left(\begin{array}{ccccc}0 & 1 & 1 & 0 & 1\\1 & 0 & 0 & 0 & 1\end{array}\right)\left[\begin{array}{cc}0 & 1\\1 & 0\\1 & 0\\0 & 0\\1 & 1\end{array}\right]
$$

$$= \begin{bmatrix} 6 & 3 & 5 & 19 & 3 & 1 \\ 3 & 5 & 6 & 2 & 2 & 1 \\ \hline 5 & 6 & 10 & 12 & 4 & 0 \\ 19 & 2 & 12 & 82 & 10 & 0 \\ 3 & 2 & 4 & 10 & 2 & 0 \\ \hline 1 & 1 & 0 & 0 & 0 & 1 \end{bmatrix} \cup \begin{bmatrix} 11 & 12 & 16 & 1 & 4 \\ 12 & 18 & 21 & 0 & 4 \\ 16 & 21 & 27 & 0 & 6 \\ 1 & 1 & 0 & 2 & 0 \\ 4 & 4 & 6 & 0 & 2 \end{bmatrix}$$

$$\cup \begin{bmatrix} 38 & 5 & 5 & 3 & 0 & 2 & 5 & 5 & 5 & 3 \\ \hline 5 & 3 & 0 & 1 & 1 & 1 & 2 & 1 & 3 & 2 \\ 5 & 0 & 1 & 0 & 0 & 0 & 0 & 1 & 0 & 0 \\ 3 & 1 & 0 & 1 & 0 & 0 & 1 & 0 & 1 & 1 \\ 0 & 1 & 0 & 0 & 1 & 0 & 0 & 1 & 1 & 1 \\ \hline 2 & 1 & 0 & 0 & 0 & 1 & 1 & 0 & 1 & 0 \\ 5 & 2 & 0 & 1 & 0 & 1 & 2 & 0 & 2 & 1 \\ 5 & 1 & 1 & 0 & 1 & 0 & 0 & 2 & 1 & 1 \\ \hline 5 & 3 & 0 & 1 & 1 & 1 & 2 & 1 & 3 & 2 \\ 3 & 2 & 0 & 1 & 1 & 0 & 1 & 1 & 2 & 2 \end{bmatrix}$$

$$\cup \begin{bmatrix} 15 & 4 & 2 & 4 & 4 & 1 & 3 & 3 & 5 \\ 4 & 3 & 2 & 2 & 1 & 2 & 2 & 1 & 1 \\ 2 & 2 & 3 & 2 & 2 & 2 & 2 & 2 & 0 \\ 4 & 2 & 2 & 3 & 1 & 1 & 1 & 2 & 1 \\ \hline 4 & 1 & 2 & 1 & 3 & 1 & 2 & 2 & 1 \\ 1 & 2 & 2 & 1 & 1 & 2 & 2 & 1 & 0 \\ 3 & 2 & 2 & 1 & 2 & 2 & 3 & 2 & 1 \\ \hline 3 & 1 & 2 & 2 & 2 & 1 & 2 & 3 & 1 \\ 5 & 1 & 0 & 1 & 1 & 0 & 1 & 1 & 2 \end{bmatrix}.$$

We see the resultant is a symmetric semi super 4-matrix. One can by this product obtain several symmetric semi super 4-matrices.

Thus we can define major product or minor products in case of super n-matrices, $n \geq 4$ as in case of super trimatrices and super bimatrices. The same type of operations are repeated as in case of super trimatrices and superbimatrices. This type of super n-matrices will be helpful in the fuzzy super model applications when we have a multi expert opinion with multi attributes. These matrices will be best suited for data storage.

# FURTHER READING

# INDEX

## A

T

# ABOUT THE AUTHORS

**Dr.W.B.Vasantha Kandasamy** is an Associate Professor in the Department of Mathematics, Indian Institute of Technology Madras, Chennai. In the past decade she has guided 12 Ph.D. scholars in the different fields of non-associative algebras, algebraic coding theory, transportation theory, fuzzy groups, and applications of fuzzy theory of the problems faced in chemical industries and cement industries.

She has to her credit 646 research papers. She has guided over 68 M.Sc. and M.Tech. projects. She has worked in collaboration projects with the Indian Space Research Organization and with the Tamil Nadu State AIDS Control Society. This is her 41$^{st}$ book.

On India's 60th Independence Day, Dr.Vasantha was conferred the Kalpana Chawla Award for Courage and Daring Enterprise by the State Government of Tamil Nadu in recognition of her sustained fight for social justice in the Indian Institute of Technology (IIT) Madras and for her contribution to mathematics. (The award, instituted in the memory of Indian-American astronaut Kalpana Chawla who died aboard Space Shuttle Columbia). The award carried a cash prize of five lakh rupees (the highest prize-money for any Indian award) and a gold medal.
She can be contacted at vasanthakandasamy@gmail.com
You can visit her on the web at: http://mat.iitm.ac.in/~wbv

---

**Dr. Florentin Smarandache** is a Professor of Mathematics and Chair of Math & Sciences Department at the University of New Mexico in USA. He published over 75 books and 150 articles and notes in mathematics, physics, philosophy, psychology, rebus, literature.

In mathematics his research is in number theory, non-Euclidean geometry, synthetic geometry, algebraic structures, statistics, neutrosophic logic and set (generalizations of fuzzy logic and set respectively), neutrosophic probability (generalization of classical and imprecise probability). Also, small contributions to nuclear and particle physics, information fusion, neutrosophy (a generalization of dialectics), law of sensations and stimuli, etc. He can be contacted at smarand@unm.edu